\documentclass[9pt]{extarticle}
\pdfoutput=1
\usepackage{graphicx}
\usepackage{hyperref}
\usepackage[utf8]{inputenc}
\usepackage{array,epsfig}
\usepackage{amsmath}
\usepackage{amsfonts}
\usepackage{amssymb}
\usepackage{dsfont}
\usepackage{amsxtra}
\usepackage{amsthm}
\usepackage{leftidx}
\usepackage{mathrsfs}
\usepackage{color}
\let\originallhook\lhook
\usepackage{MnSymbol}
\let\lhook\originallhook
\usepackage{tikz-cd}
\usepackage{enumerate}
\usepackage{todonotes}
\usepackage[backend=biber,style=alphabetic]{biblatex}
\theoremstyle{definition}
\newtheorem{defn}{Definition}[subsection]
\theoremstyle{definition}

\theoremstyle{definition}
\newtheorem{nota}[defn]{Notation}
\theoremstyle{definition}
\newtheorem{cons}[defn]{Construction}
\theoremstyle{definition}
\newtheorem{incdef}[defn]{Incorrect Definition}
\theoremstyle{definition}
\newtheorem{fact}[defn]{Fact}
\theoremstyle{definition}

\theoremstyle{definition}
\newtheorem{warn}[defn]{Warning}
\theoremstyle{definition}
\newtheorem{prob}[defn]{Problem}
\theoremstyle{plain}
\newtheorem{thm}[defn]{Theorem}
\newtheorem*{thm*}{Theorem}
\theoremstyle{plain}

\theoremstyle{definition}

\theoremstyle{plain}
\newtheorem{cor}[defn]{Corollary}
\theoremstyle{definition}
\newtheorem{rmk}[defn]{Remark}
\theoremstyle{definition}
\newtheorem{vari}[defn]{Variant}
\theoremstyle{plain}
\newtheorem{lem}[defn]{Lemma}
\theoremstyle{plain}

\theoremstyle{definition}
\newtheorem{ex}[defn]{Example}
\theoremstyle{plain}

\theoremstyle{plain}
\newtheorem{prop}[defn]{Proposition}
\theoremstyle{plain}
\newtheorem{conj}[defn]{Conjecture}

\newcommand{\bb}[1]{\mathbb{#1}}
\newcommand{\Z}{\bb{Z}}
\newcommand{\Q}{\bb{Q}}
\newcommand{\R}{\bb{R}}

\newcommand{\C}{\bb{C}}

\newcommand{\coker}{\mathrm{coker}\,}

\newcommand{\F}{\mathcal{F}}
\newcommand{\LS}{\mathcal{L}}
\newcommand{\N}{\bb{N}}

\newcommand{\einfty}{\bb{E}_{\infty}}

\newcommand{\liem}{\mathfrak{m}}

\newcommand{\Of}{\mathcal{O}}

\newcommand{\del}{\partial}

\newcommand{\hooklongrightarrow}{\lhook\joinrel\longrightarrow}

\newcommand{\simp}{\boldsymbol{\Delta}}
\newcommand{\sset}{\mathsf{Set}_{\simp}}
\newcommand{\set}{\mathsf{Set}}
\newcommand{\spa}{\mathcal{S}}
\newcommand{\map}{\mathsf{Map}}
\newcommand{\Hom}{\mathrm{Hom}}
\newcommand{\HOM}{\mathrm{HOM}}

\newcommand{\adj}{\dashv}
\newcommand{\colim}{\mathrm{colim}\,}

\newcommand{\hocolim}{\mathrm{hocolim}\,}
\newcommand{\diff}{\mathcal{T}_{\mathrm{Diff}}}
\newcommand{\diffc}{\mathcal{T}_{\mathrm{Diffc}}}

\newcommand{\spec}{\mathbf{Spec}}
\newcommand{\str}{\mathrm{Str}}
\newcommand{\loc}{\mathrm{loc}}
\newcommand{\lex}{\mathrm{lex}}
\newcommand{\strloc}{\mathrm{Str}^{\mathrm{loc}}}
\newcommand{\specr}{\mathrm{Spec}_{\R}}
\newcommand{\speco}{\mathrm{Spec}}
\newcommand{\mspec}{\mathrm{MSpec}}

\newcommand{\ltop}{{}^{\mathrm{L}}\mathsf{Top}}
\newcommand{\rtop}{{}^{\mathrm{R}}\mathsf{Top}}
\newcommand{\lTop}{{}^{\mathrm{L}}\mathsf{TOP}}

\newcommand{\sring}{sC^{\infty}\mathsf{ring}}
\newcommand{\splring}{sC^{\infty}\mathsf{PLog}}
\newcommand{\slring}{sC^{\infty}\mathsf{Log}}
\newcommand{\scring}{s\mathsf{Cring}}
\newcommand{\Mod}{\mathsf{Mod}}

\newcommand{\cotan}{\mathbb{L}}
\newcommand{\tanc}{\mathbb{T}}
\newcommand{\pregeo}{\mathcal{T}}
\newcommand{\geo}{\mathcal{G}}
\newcommand{\icat}{\mathcal{C}}
\newcommand{\icatd}{\mathcal{D}}
\newcommand{\icate}{\mathcal{E}}

\newcommand{\iop}{\mathcal{O}^{\otimes}}
\newcommand{\iopv}{\mathcal{O}^{\prime\otimes}}
\newcommand{\oinfty}{\otimes^{\infty}}
\newcommand{\cinfty}{C^{\infty}}
\newcommand{\geodiffder}{\mathcal{G}^{\mathrm{der}}_{\mathrm{Diff}}}
\newcommand{\geodiffderc}{\mathcal{G}^{\mathrm{der}}_{\mathrm{Diffc}}}
\newcommand{\geodiffc}{\mathcal{G}_{\mathrm{Diffc}}}
\newcommand{\geodiffmod}{\mathcal{G}^{\mathsf{Mod}}_{\mathrm{Diff}}}

\newcommand{\geodiff}{\mathcal{G}_{\mathrm{Diff}}}

\newcommand{\fib}{\mathrm{fib}}
\newcommand{\cofib}{\mathrm{cofib}}
\newcommand{\spect}{\mathcal{S}\mathrm{p}}
\newcommand{\fun}{\mathrm{Fun}}
\newcommand{\Fun}{\mathrm{FUN}}

\newcommand{\simpop}{\boldsymbol{\Delta}^{op}}
\newcommand{\ner}{\mathbf{N}}
\newcommand{\simpopplus}{\simp^{op}_+}
\newcommand{\sing}{\mathrm{Sing}}
\newcommand{\fin}{\mathsf{Fin}_*}
\newcommand{\xtop}{\mathcal{X}}
\newcommand{\ytop}{\mathcal{Y}}
\newcommand{\ofxtop}{(\xtop,\Of_{\xtop})}
\newcommand{\ofytop}{(\ytop,\Of_{\ytop})}
\newcommand{\ofxtopu}{(\xtop_{/U},\Of_{\xtop}|_U)}
\newcommand{\ofxtopui}{(\xtop_{/U_i},\Of_{\xtop}|_{U_i})}

\newcommand{\sch}{\mathsf{Sch}}

\newcommand{\pshv}{\mathsf{PShv}}
\newcommand{\shv}{\mathsf{Shv}}

\newcommand{\prr}{\mathsf{Pr}^{\mathrm{R}}}
\newcommand{\prl}{\mathsf{Pr}^{\mathrm{L}}}

\newcommand{\infcat}{$\infty$-category }
\newcommand{\infop}{$\infty$-operad }
\newcommand{\infcats}{$\infty$-categories }
\newcommand{\infops}{$\infty$-operads }
\newcommand{\infcatt}{$\infty$-category}

\newcommand{\infcatst}{$\infty$-categories}
\newcommand{\infopst}{$\infty$-operads}
\newcommand{\inftop}{$\infty$-topos }
\newcommand{\inftopt}{$\infty$-topos}
\newcommand{\inftopoi}{$\infty$-topoi }
\newcommand{\inftopoit}{$\infty$-topoi}
\newcommand{\cat}{\mathsf{Cat}}
\newcommand{\catinf}{\mathsf{Cat}_{\infty}}
\newcommand{\catinfh}{\widehat{\mathsf{Cat}}_{\infty}}
\newcommand{\Catinf}{\mathsf{CAT}_{\infty}}

\newcommand{\cartsp}{\mathsf{CartSp}}
\newcommand{\alg}{\mathsf{Alg}}
\newcommand{\calg}{\mathsf{CAlg}}
\newcommand{\ind}{\mathrm{Ind}}
\newcommand{\pro}{\mathrm{Pro}}
\newcommand{\sym}{\mathrm{Sym}}
\newcommand{\tor}{\mathrm{Tor}}
\newcommand{\caff}{C^{\infty}\mathsf{Aff}}
\newcommand{\daff}{\mathsf{d}C^{\infty}\mathsf{Aff}}
\newcommand{\dafffp}{\mathsf{d}C^{\infty}\mathsf{Aff}_{\mathrm{fp}}}

\newcommand{\dstack}{\mathsf{d}C^{\infty}\mathsf{St}}

\newcommand{\ev}{\mathrm{ev}}
\newcommand{\et}{\mathrm{\acute{e}t}}
\newcommand{\lfp}{\mathrm{lfp}}

\newcommand{\fp}{\mathrm{fp}}

\newcommand{\gmt}{\mathrm{gmt}}

\newcommand{\rmalg}{\mathrm{alg}}
\usepackage{subfiles}

\title{Derived $C^{\infty}$-Geometry I: Foundations}
\author{Pelle Steffens\thanks{pelle.steffens@outlook.com} \\ Universit\'{e} Toulouse III-Paul Sabatier}
\date{April 17, 2023}

\begin{document}
\setcounter{page}{1}

\maketitle

\begin{abstract}This work is the first in a series laying the foundations of derived geometry in the $\cinfty$ setting, and providing tools for the construction and study of moduli spaces of solutions of partial differential equations that arise in differential geometry and mathematical physics. To advertise the advantages of such a theory, we start with a detailed introduction to derived $\cinfty$-geometry in the context of symplectic topology and compare and contrast with Kuranishi space theory. In the body of this work, we avail ourselves of Lurie's extensive work on abstract structured spaces \cite{dagv,sag} to define \infcats of \emph{derived $\cinfty$-rings and $\cinfty$-schemes} and \emph{derived $\cinfty$-rings and $\cinfty$-schemes with corners} via a universal property in a suitable $(\infty,2)$-category of \infcats with respect to the ordinary categories of manifolds and manifolds with corners (with morphisms the \emph{$b$-maps} of Melrose \cite{Melrose1} in the latter case), and prove many basic structural features about them. Along the way, we establish some derived flatness results for derived $\cinfty$-rings of independent interest. 
\end{abstract}

\tableofcontents

\section{Introduction}  

The purpose of this work is to lay the foundations of \emph{derived geometry} in the differentiable, that is, $C^{\infty}$-setting for applications in the theory of moduli spaces in differential geometry, symplectic geometry and mathematical physics, using the modern language and powerful tools of higher category theory, higher topos theory, and higher algebra. The corresponding theory of derived \emph{algebraic} geometry has been well established for a number of years due to the seminal works of Lurie (DAG series, \cite{dagv} through \cite{dagxiv}, and \cite{sag}) and To\"{e}n-Vezzosi (Homotopical Algebraic Geometry, \cite{TV1,TV2}). Derived geometry has been established in other contexts as well; there are derived versions of analytic geometry due to Lurie \cite{dagix} and Porta-Yu \cite{P1,PY1}. In fact, a substantial literature on derived differential geometry already exists since the pioneering work of Spivak \cite{spivak}, including a sizable work-in-progress of Joyce \cite{Joy1}, the model categorical efforts of Carchedi-Roytenberg \cite{CR2,CR1}, the work of Borisov-No\"{e}l \cite{BN}, recent work of Behrend-Liao-Xu \cite{BLX} and Amorim-Tu \cite{AmorimTu}, and the thesis of Nuiten \cite{Nui2} (and undoubtedly others that would deserve to be mentioned).\\
Derived geometry is a confluence of classical geometry, homological and homotopical algebra, intersection theory, deformation theory and higher sheaf theory, and may be approached and appreciated from any of these avenues. There are a number of excellent introductory texts available that do the subject justice; let us mention in particular the surveys of To\"{e}n and Anel \cite{Toen,Anel}. We will motivate the theory we wish to develop in this work via an intersection problem, but one quite different from the well-known algebro-geometric account that passes from Serre's intersection formula to Koszul resolutions and derived pushouts of dg-algebras, as in the introduction of \cite{dagv}, for instance. We will be concerned with intersection theory in infinite dimensions.
\begin{rmk}
In the body of this work, we will freely use the language of modern higher category theory. The following section however is aimed at a wider audience interested in homological and homotopical methods in relation to moduli spaces, but not necessarily versed in derived algebraic geometry or intimately familiar with \infcatst. We believe that methods and tools from these fields have the potential to contribute substantively to the global study of moduli spaces in differential geometry and mathematical physics; it is our hope that for those who come new to derived geometry and abstract homotopy theory, the story below can serve as a point of entry. 
\end{rmk}

\subsection{From PDEs to \infcatst: an introduction}
From a sufficiently abstract vantage, in mathematics influenced by Quantum Field Theory such as symplectic geometry and gauge theory, one studies the geometry of moduli spaces of solutions of nonlinear elliptic equations on manifolds -which are usually required to be compact (if not, the function spaces need to satisfy some decay estimates to admit well behaved moduli spaces)- up to the action of a (possibly infinite dimensional) Lie group of symmetries and perhaps suitably compactified. Dispensing with the issues of compactification and symmetries for the moment, we are interested in the following situation:
\begin{enumerate}[$(1)$]
    \item $M$ a compact smooth manifold. 
    \item $V\rightarrow M$ a smooth fibre bundle over $M$.
    \item $F\rightarrow M$ a smooth vector bundle over $M$.
    \item $P:\Gamma(V)\rightarrow\Gamma(F)$ a nonlinear elliptic differential operator acting between smooth sections of $V$ and $F$. 
\end{enumerate}
Let $\mathsf{Sol}(P)=P^{-1}(0)$. Let $x\in \mathsf{Sol}(P)$ and suppose that the linearization $TP_x:\Gamma(x^*TM)\rightarrow\Gamma(F)$, a 2-term Fredholm complex with finite dimensional homology, is surjective. Then $\mathsf{Sol}(P)$ admits the structure of a smooth manifold in a neighbourhood of $x$. If the linearized differential operator is not surjective, we still have the following important principle.
\begin{fact}[Local finite dimensional reduction by Kuranishi models]\label{fact:kur}
Locally, $\mathsf{Sol}(P)$ is given by the zero set of a smooth function $f:\R^n\rightarrow \R^k$ such that at each solution $x$ of $f=0$, the 2-term complex determined by the linearization of $f$ at $x$ is \emph{quasi-isomorphic} to the 2-term complex determined by the linearization of $P$ at $x$. \end{fact}
This follows from an application of the inverse function theorem for Banach manifolds and elliptic bootstrapping methods, after replacing the spaces of smooth sections with Sobolev completions of sufficiently high regularity; we refer to the appendices of \cite{mcduffsalamon} for a textbook account in symplectic topology.\\
Depending on the geometric situation, it may or may not be possible to \emph{perturb} the operator $P$ and obtain a well defined cobordism class of smooth spaces of solutions. When this is not possible (when $\mathsf{Sol}(P)$ is the space of genus 0 pseudo-holomorphic curves on a non semipositive symplectic manifold, for instance), one is forced to make sense of $\mathsf{Sol}(P)$ using the zeroes of local finite dimensional reductions which are not transverse. We have the following two problems:
\begin{enumerate}[$(a)$]
    \item The local finite dimensional reductions are far from unique; only the homology complex induced by $TP_x$ as $x$ varies over $\mathsf{Sol}(P)$ is invariant. 
    \item The space $\mathsf{Sol}$ is a gluing of the zero sets of local finite dimensional reductions, but as these spaces can have arbitrarily badly behaved topology (as subspaces of some Cartesian space) it is not clear how to perform this gluing and obtain some sort of geometric $\cinfty$ structure on $\mathsf{Sol}(P)$.
\end{enumerate}
Broadly speaking, our goal may be stated as follows.
\begin{itemize}
    \item Come up with some sort of `singular' $\cinfty$ structure $\geo$ that a topological space may be endowed with such that a $\geo$-structure on a space $X$ locally turns it into the zero set of a section of a vector bundle on some manifold. The structure $\geo$ should satisfy the following desiderata.
    \begin{enumerate}
        \item[$(G1)$] Moduli spaces of elliptic PDEs on compact manifolds admit such a structure, hopefully in a canonical manner.
        \item[$(G2)$] It is possible to `do' geometry with spaces endowed with the structure $\geo$.
    \end{enumerate}
    $G2$ is a bit vague. Let us try to not be too ambitious for the moment and require something like an `integration theory' for enumerative geometry purposes. Instead we notice that a minimal requirement for being a legitimate geometric space of some sort involves \emph{having a tangent bundle}. So we will ask that
    \begin{enumerate}
        \item[$(G2)$] A space $X$ equipped with the structure $\geo$ has a tangent bundle (in a suitable sense to be determined). 
    \end{enumerate}
\end{itemize}
We will be gluing the local models of Fact \ref{fact:kur} together; to do so, we'll need to talk about \emph{morphisms} of such, that is, we need a category of local models. 
\begin{defn}\label{def:kur}
An \emph{affine Kuranishi model (without isotropy)} is a triple $(X,p:E\rightarrow X,s)$ where $X$ is a smooth manifold, $p:E\rightarrow X$ is a vector bundle on $E$, and $s$ is a section of $p$. We will usually just write $E$ for the bundle $p:E\rightarrow X$. Given two affine Kuranishi models $(X,E,s)$ and $(Y,F,t)$, a morphism $f:(X,E,s)\rightarrow (Y,F,t)$ is a commuting diagram 
\[
\begin{tikzcd}
E\ar[d,"p"] \ar[r,"f_v"] & F \ar[d,"q"] \\
X\ar[r,"f_b"] & Y 
\end{tikzcd}
\]
where $f_v$ is fibrewise linear such that $f_v\circ s=t\circ f_b$. Affine Kuranishi models and morphisms between them form a category, that we denote $\mathsf{AffKur}$. We let $Z(s)$ denote the zero locus of $s$, topologized as a subspace of $X$.
\end{defn}
We notice that these local models have a natural notion of a tangent bundle. If a topological space $Y$ is covered by the zero locus of a single affine Kuranishi model, then there exists a distinguished object (up to isomorphism) in the \emph{derived category} of sheaves of complexes of real vector spaces on $Y$; indeed, for an affine Kuranishi model $\mathbf{K}=(M,E,s)$ we have the complex
\[\tanc_{\mathbf{K}}:= T_M \overset{Ts}{\longrightarrow} s^*T^v_E \in \mathbf{D}(\shv_{\mathrm{Vect}_{\R}}(Z(s))),\]
where $T^v_E$ is the vertical tangent bundle of $E$ with respect to $M$, fitting into a fibre sequence
\[ T^v_E\longrightarrow T_E \longrightarrow p^*T_M  \]
of vector bundles over $E$. We call this object in the derived category of sheaves on $Z(s)$ the \emph{tangent complex} or the \emph{virtual tangent sheaf} of $\mathbf{K}$. If the topological space $Y$ is $\mathsf{Sol}(P)$ for some moduli problem defined by an elliptic equation, then the homology of linearization also determines a well defined object in $\mathbf{D}(\shv_{\mathrm{Vect}_{\R}}(Y))$. \\
In the coming sections, we will go through several candidate structures that we call \emph{Kuranishi atlases}, after \cite{FO} and \cite{fooo1}. Our meta-goal is to explain why higher category theory is necessary to obtain a satisfactory theory.
\begin{warn}
In this introduction, we will be playing fast and loose with well established mathematical terminology, particularly that of a \emph{Kuranishi atlas/structure}. First of all, we do not, for the most part, talk about Kuranishi structures with isotropy, for the sake of pedagogy. More importantly, we will in the sequel introduce several distinct notions of a Kuranishi atlas, most of which are not equivalent to any the definitions found in the literature \cite{fooo1,Joykura,mcduffwehrheim,mcduffwehrheim,Pardon} (nor will we make a serious attempt at a comparison). We give two reasons for our (abuse of) terminology.
\begin{enumerate}[$(1)$]
    \item Our Kuranishi atlases are at least in spirit very similar to what Fukaya and Ono had in mind in their pioneering work and are designed to accomplish the same sort of job.
    \item The `correct' definition upon which we eventually settle is in fact the one forced on us by universal (higher) categorical constructions, as we show in this work and its successors.
\end{enumerate}
\end{warn}
\subsubsection{Attempt 1: atlases up to homotopy}
We wish to define a mathematical structure that glues together affine Kuranishi models, that is, zero sets of $\cinfty$ sections of (finite rank) vector bundles on (finite dimensional) manifolds, in such a manner that the local tangent complexes likewise glue and form a global object. From Fact \ref{fact:kur} we should infer that the tangent complex of an affine Kuranishi model \emph{as an object in the derived category} is the correct substitute for the tangent bundle of a manifold. This introduces a problem: gluing of geometric objects should happen along isomorphisms, but it seems that the category of affine Kuranishi models has too few of those. An isomorphism of affine Kuranishi models $(X,E,s)\rightarrow (Y,F,t)$ induces diffeomorphisms between $X$ and $Y$ and fibrewise linear diffeomorphisms between $E$ and $F$, but there are many more morphisms of affine Kuranishi models that induce a bijection on zero sets and a quasi-isomorphism on tangent complexes; that is to say, the inverse function theorem for affine Kuranishi models does not hold. For instance, if $(X,E,s)$ is an affine Kuranishi model and $W$ is a finite dimensional real vector space, then we have an obvious map  \[ (X\times W,E\times W,s+\mathrm{id}_W)   \longrightarrow  (X,E,s)   \]
which induces a bijection $Z(s+\mathrm{id}_W)\cong Z(s)$ and a quasi-isomorphism $\tanc_{(X\times W,E\times W,s+\mathrm{id}_W) }\simeq \tanc_{(X,E,s) }$ of complexes of sheaves on $Z(s)$. If $(X,E,s)$ is a local model for a moduli space $\mathsf{Sol}(P)$, then so is $(X\times W,E\times W,s+\mathrm{id}_W)$, and there does not seem to be a very compelling reason to prefer one over the other. The most straightforward way to remedy this is to apply the abstract procedure of \emph{localizing a category at a subcategory of weak equivalences} to obtain the `correct' category of affine Kuranishi models. Let $\icat$ be a category and let $W\subset \icat$ be a subcategory of \emph{weak equivalences} that contains all isomorphisms, then define, up to essentially unique equivalence of categories, a new category $\icat[W^{-1}]$ equipped with a functor 
\[L:\icat\longrightarrow \icat[W^{-1}]\]
by declaring that $L$ has the following universal property in the 2-category of categories: for each category $\icatd$, the restriction functor along the functor $L:\icat\rightarrow \icat[W^{-1}]$ induces an equivalence of categories
\[ \fun(\icat[W^{-1}],\icatd)\overset{\simeq}{\longrightarrow} \fun_{W}(\icat,\icatd),\]
where $\fun_{W}(\icat,\icatd)\subset \fun(\icat,\icatd)$ is the full subcategory spanned by functors that carry weak equivalences in $\icat$ to isomorphisms in $\icatd$. The existence of this category follows from the following 
\begin{cons}[Gabriel-Zisman localization]\label{cons:zigzag}
Let $(\icat,W)$ be a category with weak equivalences. We define a new category, denoted $\icat[W^{-1}]$, as follows. 
\begin{enumerate}
    \item[$(O)$] Objects are the objects of $\icat$.
    \item[$(M)$] For $C$ and $D$ two objects of $\icat$, consider the set of \emph{zigzags} between $C$ and $D$, which are diagrams
    \[
    \begin{tikzcd}
    C \ar[r,dash]& E_{1}\ar[r,dash] & E_{2}  \ar[r,dash]  & \ldots \ar[r,dash] & E_{m-1}\ar[r,dash] & D
    \end{tikzcd}
     \]
     in $\icat$ of variable finite length, such that
     \begin{enumerate}[$(1)$]
         \item arrows are allowed to go in either direction.
         \item adjacent arrows go in different directions.
         \item all leftward maps are in $W$.
         \item there are no identity arrows.
     \end{enumerate}
     We let $\Hom_{\icat[W^{-1]}]}(C,D)$ be the set of equivalence classes for the equivalence relation on the set of zigzags between $C$ and $D$ generated by the following operations: if a zigzag contains a piece 
     \[\ldots\longleftarrow E_{i-1} \longrightarrow  E_i\longleftarrow  E_{i+1}\longrightarrow\ldots\] or a piece
     \[\ldots \longrightarrow E_{i-1} \longleftarrow  E_i\longrightarrow E_{i+1} \longleftarrow\ldots\]
and we are given a weak equivalence $f:E_{i}'\rightarrow E_i$, then we may replace these pieces with 
\[\ldots \longleftarrow E_{i-1} \longrightarrow  E_i\overset{f}{\longleftarrow} E_{i}' \overset{f}{\longrightarrow} E_{i} \longleftarrow E_{i+1}\longrightarrow \ldots\]  
respectively with
\[  \ldots \longrightarrow E_{i-1}\longleftarrow E'_i  \longrightarrow E_{i+1}\longleftarrow\ldots \]     
Similarly, if we are given a weak equivalence $g:E_i\rightarrow E''_i$, then we may replace the pieces above with 
\[ \ldots\longleftarrow E_{i-1}\longrightarrow E_i''\longleftarrow E_{i+1}\longrightarrow\ldots \]
respectively with
 \[\ldots\longrightarrow E_{i-1} \longleftarrow  E_i\overset{g}{\longrightarrow} E_{i}'' \overset{g}{\longleftarrow} E_{i} \longrightarrow E_{i+1}\longleftarrow \ldots\]
 The composition is simply the map that concatenates zigzags and composes morphisms and removes identities in the resulting diagram if $(2)$ and $(4)$ are not satisfied.
\end{enumerate}
There is an obvious functor $\icat\rightarrow\icat[W^{-1}]$ by viewing every \emph{nonidentity} arrow as a zigzag of length 1. This functor satisfies the universal property in the 2-category of categories we formulated above.
\end{cons}

\begin{ex}
Let $\mathcal{A}$ be a Grothendieck abelian category and let $\mathrm{Ch}(\mathcal{A})$ be the category of (unbounded, homologically graded) chain complexes in $\mathcal{A}$. The unbounded derived category $\mathbf{D}(\mathcal{A})$ of $\mathcal{A}$ is the category $\mathrm{Ch}(\mathcal{A})[q-iso^{-1}]$ localized at the subcategory whose morphisms are quasi-isomorphisms. Let $X$ be a topological space, then we can identify the category $\shv_{\mathrm{Ch}(\mathcal{A})}(X)$ of sheaves of chain complexes on $X$ with the category of chain complexes in the category $\shv_{\mathcal{A}}(X)$ of $\mathcal{A}$-valued sheaves, and the derived category $\mathbf{D}(\shv_{\mathcal{A}}(X))$ coincides with the localization $\shv_{\mathrm{Ch}(\mathcal{A})}(X)[q-iso^{-1}]$.
\end{ex}
\begin{ex}
Let $\mathsf{CGH}$ be the category of compactly generated Hausdorff topological spaces and let $W\subset \mathsf{CGH}$ be the subcategory whose morphisms are weak homotopy equivalences. Then the category $\mathsf{CGH}[W^{-1}]$ is the classical homotopy category of topological spaces.    
\end{ex}
\begin{nota}
In light of the previous example, we will somewhat abusively also denote the category $\icat[W^{-1}]$ by $h\icat$ and call it the \emph{homotopy category of $\icat$}.    
\end{nota}
We could now localize the category $\mathsf{AffKur}$ at the subcategory of maps that induce a homeomorphism on zero loci and a quasi-isomorphism on tangent complexes, but it will be quite convenient to introduce a different homotopy theory by embedding Kuranishi models in a larger category.
\begin{defn}\label{defn:dgmanifold}
For a manifold $X$, we let $\cinfty_X$ denote its sheaf of $\cinfty$ functions. Also, we will denote smooth vector bundles over $X$ and their corresponding sheaves of $\cinfty_X$-modules of sections by the same symbol. Let $n\in \Z_{\geq 1}$. A \emph{differential graded manifold}, or \emph{dg-manifold of amplitude $n$}, is a triple $(X,\{E_i\}_{1\leq i\leq n},d)$ of a manifold $X$, a collection of (constant finite rank) vector bundles $\{E_i\}_{1\leq i\leq n}$ on $X$ and $d$ is a differential on the sheaf of graded commutative $\cinfty_X$-algebras $\sym_{\cinfty_X}E^{\vee}$ free on the graded vector bundle $E^{\vee}=\oplus_{i=1}^n E_i^{\vee}[i]$, making it a sheaf of differential graded commutative $\cinfty_X$-algebras. By convention, a dg-manifold of amplitude $0$ is simply a manifold.\\
A morphism of dg-manifolds $(X,\{E_i\}_{i=1}^n,d_X)\rightarrow (Y,\{F_i\}_{i=1}^m,d_Y)$ is a $\cinfty$ map $f:X\rightarrow Y$ together with a map $f^*\sym_{\cinfty_Y} F^{\vee}\rightarrow \sym_{\cinfty_X} E^{\vee}$ of sheaves of differential graded $\cinfty_X$-algebras. We thus have a category $\mathbf{dgMan}$ of dg-manifolds (of variable amplitude).
\end{defn}
\begin{rmk}
We can \emph{identify} the category of affine Kuranishi models with the category of dg-manifolds of amplitude $\leq 1$. We have an isomorphism of categories $\mathsf{AffKur}\rightarrow\mathbf{dgMan}_{\leq 1}$ that carries $(X,E,s)$ to the dg-manifold $(X,E,\iota_s)$ where $\iota_s$ is the differential on the exterior algebra bundle $\Lambda^{\bullet}E^{\vee}$ that contracts a section of the dual bundle $E^{\vee}$ with the section $s$.
\end{rmk}
Let $(X,\{E_i\}_{i=1}^n,d)$ be a dg-manifold of amplitude $n$ and write $E=\oplus_{i=1}^n E_i[-i]$. The differential $\sym_{\cinfty_X}E^{\vee}\rightarrow \sym_{\cinfty_X}E^{\vee}$ is determined by a collection of maps $\{d_k:\sym_{\cinfty_X}^k E\rightarrow E\}_{k\geq 0}$ of sheaves of $\cinfty_X$-modules of degree $-1$. The map $d_0$ is simply a section $s$ of $E_1$. We have a functor
\[ Z:\mathbf{dgMan}\longrightarrow \mathsf{Top} \]
carrying a dg-manifold to the zero locus of the section $d_0=s$. Note that $d_1$ consists of a collection of maps $\{E_i\rightarrow E_{i+1}\}_{1\leq i\leq n}$ of vector bundles over $X$.
\begin{defn}\label{def:kurweakeq}
Let $(X,\{E_i\}_{1\leq i\leq n},d)$ be a dg-manifold with $d_0=s$, and let $x\in Z(s)$ be a point. The differential $T_xs$ of $s$ at $x$ is a section $T_xX\rightarrow T_{s(x)}E_1\cong T_xX\oplus (E_1)_x$ of the projection $T_{x(x)}E_1\rightarrow T_xX$ and is thus determined by the linear map $T_xX\rightarrow (E_1)_x$. The \emph{tangent complex of $(X,\{E_i\}_{1\leq i\leq n},d)$ at $x$} is the complex 
\[\begin{tikzcd}
    T_xX\ar[r,"T_{x}s" ]&(E_1)_x  \ar[r,"d_1|_{E_1}"] & (E_2)_x \ar[r,"d_1|_{E_2}"] & \ldots
\end{tikzcd}  \]
Note that if $(X,\{E_i\}_{1\leq i\leq n},d)$ is an affine Kuranishi model, then only the first two terms of this complex are nonzero. Let $(Y,\{F_i\}_{1\leq i\leq m},d)$ be another dg-manifold with $d_0=t$ and suppose we have a morphism $f:(X,\{E_i\}_{i=1}^n,d)\rightarrow (Y,\{F_i\}_{i=1}^m,d)$. Then we have, for each $x\in Z(s)$ an induced map
\[
\begin{tikzcd}
T_xX\ar[r,"T_{x}s"]\ar[d,"T_xf"']&(E_1)_x \ar[d] \ar[r,"d_{1,x}|_{E_1}"] & (E_2)_x \ar[d]\ar[r,"d_{1,x}|_{E_2}"] & \ldots \\
T_{f(x)}Y\ar[r,"T_{f(x)}t"]& (F_1)_{f(x)} \ar[r,"d_{1,x}|_{F_1}"] &  (F_2)_{f(x)} \ar[r,"d_{1,x}|_{F_2}"] & \ldots
\end{tikzcd}
\]
of chain complexes. We say that $f$ is a \emph{weak equivalence at $x$} if the map of chain complexes above is a quasi-isomorphism. We say that $f$ is a \emph{weak equivalence} if $f$ induces a bijection $Z(s)\cong Z(t)$ and $f$ is a weak equivalence at all points of $Z(s)$. Let $W \subset \fun(\Delta^1,\mathbf{dgMan})$ be the full subcategory spanned by the weak equivalences. This full subcategory contains all isomorphisms so the pair $(\mathbf{dgMan},W)$ is a category with weak equivalences.
\end{defn}
\begin{rmk}
It can be shown that if $f$ is a weak equivalence at $p$, then $f$ induces a homeomorphism from a neighbourhood of $p\in Z(s)$ onto a neighbourhood of $f(p)\in Z(t)$; thus, a weak equivalence always induces a homeomorphism on zero sets.
\end{rmk}
\begin{defn}
The \emph{homotopy category} of affine Kuranishi models is the full subcategory of the category $\mathbf{dgMan}[W^{-1}]$ spanned by affine Kuranishi models; that is, it is the image under the composition 
\[ \mathsf{AffKur}\cong \mathbf{dgMan}_{\leq 1} \subset\mathbf{dgMan} \overset{L}{\longrightarrow}\mathbf{dgMan}[W^{-1}].\]
We denote this category by $h\mathsf{AffKur}$.
\end{defn}
\begin{rmk}
The notation $h\mathsf{AffKur}$ is abusive, since this category is not the localization of $\mathsf{AffKur}$ at the weak equivalences of Definition \ref{def:kurweakeq}. We will explain later why it is that the localization $\mathbf{dgMan}[W^{-1}]$ is better behaved (roughly speaking, the category $\mathsf{AffKur}$ does not have enough `path objects'). For now, we will stick with the full subcategory $h\mathsf{AffKur}$, but all constructions and results in this introduction hold for the larger category $\mathbf{dgMan}[W^{-1}]$ as well.
\end{rmk}
\begin{rmk}
There is an obvious forgetful functor $Z:\mathsf{AffKur}\rightarrow\mathsf{Top}$ sending an affine Kuranishi model $\mathbf{K}$ to the topological space $Z(\mathbf{K}):=Z(s)$. It follows from the definition of a weak equivalence that this forgetful functor factors through a functor $h\mathsf{AffKur}\rightarrow \mathsf{Top}$ that we abusively also denote $Z$.
\end{rmk}

Now if we are given a nice topological space $X$, we would like to express the idea that $X$ arises as the gluing of a collection of affine Kuranishi models, and this gluing should be compatible on the overlaps, which means that a cocycle condition should hold, all in the category $h\mathsf{AffKur}$, that is, `up to homotopy'. To express gluing of two Kuranishi models $(X,E,s)$ and $(Y,F,t)$ along a common open subset in $Z(s)$ and $Z(t)$, we have to verify that restricting to open subsets \emph{of the zero locus} is well defined in the homotopy category $h\mathsf{AffKur}$.
\begin{defn}
Let $\mathbf{K}=(X,E,s)$ be an affine Kuranishi model and let $U\subset Z(s)$ be an open subset, then we say that a morphism $\mathbf{H}\rightarrow \mathbf{K}$ in $h\mathsf{AffKur}$ \emph{exhibits $\mathbf{H}$ as a localization of $\mathbf{K}$ with respect to $U$} if for every affine Kuranishi model $\mathbf{J}$, restriction along the map $\mathbf{H}\rightarrow \mathbf{K}$ induces a bijection 
\[ \Hom_{h\mathsf{AffKur}}(\mathbf{J},\mathbf{H}) \overset{\simeq}{\longrightarrow} \Hom_{h\mathsf{AffKur}}^U(\mathbf{J},\mathbf{K}), \]
where $\Hom_{h\mathsf{AffKur}}^U(\mathbf{J},\mathbf{K})\subset \Hom_{h\mathsf{AffKur}}(\mathbf{J},\mathbf{K})$ is the subset of those maps $f:\mathbf{J}\rightarrow\mathbf{K}$ that satisfy the condition that $Z(f):Z(\mathbf{J})\rightarrow Z(\mathbf{K})$ factors through $U$. 
\end{defn}
It follows immediately from the definition that a localization of $\mathbf{K}$ with respect to $U\subset Z(s)$ is unique up to unique isomorphism in the homotopy category $h\mathsf{AffKur}$, provided it exists. As it turns out, a localization at any open subset exists.
\begin{lem}\label{lem:localizationsexist1}
For every open $U\subset Z(s)$ and every open set $V\subset X$ such that $V\cap Z(s)=U$, the morphism of affine Kuranishi models $(V,E|_V,s|_V)\rightarrow (X,E,s)$ exhibits a localization with respect to $U$.
\end{lem}
\begin{proof}[Sketch of proof]
Denote $\mathbf{K}=(X,E,s)$ and $\mathbf{H}=(V,E|_V,s|_V)$ and let $\iota:\mathbf{H}\rightarrow\mathbf{K}$ be the map induced by the inclusion $V\subset X$. Let $\mathbf{J}$ be any affine Kuranishi model and suppose we are given a morphism $f:\mathbf{J}\rightarrow \mathbf{K}$ in the homotopy category such that $Z(\mathbf{J})\rightarrow Z(\mathbf{K})$ factors through $U$. We first show that there is a morphism $\tilde{f}:\mathbf{J}\rightarrow\mathbf{H}$ in the homotopy category such that $f=\iota\circ\tilde{f}$. The morphism $f$ is represented by some zigzag
 \begin{equation}\label{eq:zigzag}
    \begin{tikzcd}
    \mathbf{J} \ar[r,dash]& \mathbf{J}_{1}\ar[r,dash] & \mathbf{J}_{2}  \ar[r,dash]  & \ldots \ar[r,dash] & \mathbf{J}_{n-1}\ar[r,dash] & \mathbf{K}
    \end{tikzcd}
  \end{equation}
in the category $\mathbf{dgMan}$. To define $\tilde{f}$, we distinguish between two cases.
\begin{enumerate}[$(i)$]
\item Suppose that the zigzag \eqref{eq:zigzag} ends with 
\[ \mathbf{J}_{n-1}\longrightarrow \mathbf{K}, \]
then the map $Z(\mathbf{J}_{n-1})\rightarrow Z(\mathbf{K})$ factors through $U$, so writing $\mathbf{J}_{n-1}=(Y,\{F_i\}_i,d)$ and $W=Y\times_XV$, the map $\mathbf{J}_{n-1}|_W=(W,\{F_i|_W\}_i,d|_W)\rightarrow \mathbf{J}_{n-1}$ induced by the open inclusion $W\subset Y$ is a weak equivalence. By definition of the equivalence relation on zigzags, we may replace the map $\mathbf{J}_{n-1}\rightarrow \mathbf{K}$ by the piece 
\[ \mathbf{J}_{n-1}\overset{\simeq}{\longleftarrow} \mathbf{J}_{n-1}|_{W}\overset{\simeq}{\longrightarrow}  \mathbf{J}_{n-1}\longrightarrow \mathbf{K} \]
and compose the last two arrows. Since the map $\mathbf{J}_{n-1}|{_W}\rightarrow \mathbf{K}$ factorizes as $\mathbf{J}_{n-1}|{_W}\rightarrow\mathbf{H}\overset{i}{\rightarrow}\mathbf{K}$, we find the desired map $\tilde{f}$.
 \item Suppose that the zigzag \eqref{eq:zigzag} ends with 
\[ \mathbf{J}_{n-1}\overset{\simeq}{\longleftarrow} \mathbf{K}. \]
Write again $\mathbf{J}_{n-1}=(Y,\{F_i\}_i,d)$, then because the leftward map above is a weak equivalence, we have $Z(d_0)\cong Z(\mathbf{K})$. Choose some open set $W\subset Y$ such that $W\cap Z(d_0)\cong U$, then we set $\mathbf{H}'=(X\times_{Y}W,E|_{X\times_{Y}W},s|_{X\times_{Y}W})$. The maps $\mathbf{H}\rightarrow \mathbf{K}$ and $\mathbf{H}\overset{\simeq}{\rightarrow} \mathbf{J}_{n-1}|_W$ determine for each $\mathbf{J}\in \mathsf{AffKur}$ a commuting diagram 
\[
\begin{tikzcd}
\Hom_{h\mathsf{AffKur}}(\mathbf{J},\mathbf{H}')\ar[r]\ar[d,"\cong"] &\Hom_{h\mathsf{AffKur}}^U(\mathbf{J},\mathbf{K}) \ar[d,"\cong"]\\
\Hom_{h\mathsf{AffKur}}(\mathbf{J},\mathbf{J}_{n-1}|_{W})\ar[r] &\Hom_{h\mathsf{AffKur}}^U(\mathbf{J},\mathbf{J}_{n-1}) 
\end{tikzcd}
\]
of sets, where the vertical maps are bijections. Then we apply $(i)$ to the map $\mathbf{J}_{n-2}\rightarrow \mathbf{J}_{n-1}$ to find some $\tilde{f'}\in \Hom_{h\mathsf{AffKur}}(\mathbf{J},\mathbf{J}_{n-1}|_{W})$ and we take its inverse image along the left vertical map in the commuting diagram above. Since $\mathbf{H}'\cong \mathbf{H}$ in the homotopy category via the zigzag
\[ (X\times_{Y}W,E|_{X\times_{Y}W},s|_{X\times_{Y}W}) \overset{\simeq}{\longleftarrow }(X\times_{Y}W\cap V,E|_{X\times_{Y}W\cap V},s|_{X\times_{Y}W\cap V})\overset{\simeq}{\longrightarrow} (V,E|_V,s|_V),  \]
this determines an element in $\Hom_{h\mathsf{AffKur}}(\mathbf{J},\mathbf{H})$. 
\end{enumerate} 
It is straightforward (but a bit tedious) to verify that the assignment described in $(i)$ and $(ii)$ is well defined on equivalence classes of zigzags and is a left and right inverse of the map that composes with $\iota$.
\end{proof}
For $\mathbf{K}=(X,E,s)$ and $V\subset X$ an open set, we write $\mathbf{K}|_V$ for $(V,E|_V,s|_V)$. For $U\subset Z(\mathbf{K})$ an open set, we will also write, somewhat abusively, $\mathbf{K}|_U$ for the localization of $\mathbf{K}$ with respect to $U$, which we understand as an object in the homotopy category, where it is defined up to unique isomorphism. It follows easily from the preceding lemma that sending a localization of $\mathbf{K}$ to the underlying open subset of $Z(\mathbf{K})$ induces an equivalence of categories between the full subcategory of the slice category $h\mathsf{AffKur}_{/\mathbf{K}}$ spanned by localizations of $\mathbf{K}$ and the lattice of open subsets of $Z(\mathbf{K})$ (viewed as a category). If $f:\mathbf{K}\rightarrow\mathbf{H}$ is a morphism in $h\mathsf{AffKur}$, then for each open $U\subset Z(\mathbf{K})$ we denote by $f|_U$ the composition $\mathbf{K}|_U\rightarrow\mathbf{K}\rightarrow\mathbf{H}$. \\
The following definition is the most naive approach one might be tempted to try.
\begin{incdef}\label{defn:kuratlas}
Let $X$ be a paracompact Hausdorff topological space. A \emph{naive Kuranishi atlas (without isotropy)} on $X$ consists of the following data.
\begin{enumerate}[$(a)$]
    \item An open covering $\{U_i\rightarrow X\}_{i\in I}$ of $X$ (not necessarily finite).
    \item A collection of affine Kuranishi models $\{\mathbf{K}_i\}_{i\in I}$ with zero loci $\{Z(\mathbf{K}_i)\}_{i\in I}$ called \emph{charts}. 
    \item A collection of homeomorphisms $\psi_i:Z(\mathbf{K}_i)\rightarrow U_i$ called \emph{footprint maps} or \emph{chart maps}.
    \item For every pair of indices $i,j\in I$ such that $U_{ij}:=U_i\cap U_j$ is nonempty, an isomorphism  $\phi_{ij}:\mathbf{K}_i|_{\psi_i^{-1}(U_{ij})}\rightarrow \mathbf{K}_j|_{\psi_j^{-1}(U_{ij})}$ in the homotopy category $h\mathsf{AffKur}$.
\end{enumerate}
These data are required to satisfy the following conditions.
\begin{enumerate}[$(1)$]
    \item The transition maps $\phi_{ij}$ are compatible with the footprint maps: for all pairs $i,j\in I$ such that $U_{ij}$ is nonempty, the diagram
    \[
    \begin{tikzcd}
    \psi_i^{-1}(U_{ij})=Z(\mathbf{K}_i|_{\psi_i^{-1}(U_{ij})}) \ar[dr,"\psi_i|_{U_{ij}}"'] \ar[rr,"Z(\phi_{ij})"] && Z(\mathbf{K}_j|_{\psi_j^{-1}(U_{ij})})=\psi_j^{-1}(U_{ij}) \ar[dl,"\psi_j|_{U_{ij}}"]\\
    & U_{ij}
    \end{tikzcd}
    \]
    commutes.
    \item The cocycle condition holds: for every triple $i,j,k\in I$ such that $U_{ijk}:=U_i\cap U_j\cap U_k$ is nonempty, we see that $(1)$ and the universal property of localization imply that the composition 
    \[ \phi_{ij}|_{\psi^{-1}(U_{ijk})}:\mathbf{K}_i|_{\psi_i^{-1}(U_{ijk})} \longrightarrow \mathbf{K}_i|_{\psi_i^{-1}(U_{ij})}\overset{\phi_{ij}}{\longrightarrow} \mathbf{K}_j|_{\psi_j^{-1}(U_{ij})} \]
    where the first map is a localization, factors through $\mathbf{K}|_{\psi_j^{-1}(U_{ijk})}$. Then we can apply $\phi_{jk}$, and we demand that the equality \[\phi_{jk}|_{\psi_j^{-1}(U_{ijk})}\circ \phi_{ij}|_{\psi_i^{-1}(U_{ijk})}=\phi_{ik}|_{\psi_i^{-1}(U_{ijk})}\]
    holds.
\end{enumerate}
\end{incdef}
Note that this description is in almost complete analogy with the notion of an atlas on a manifold; indeed, suppose that for each $i$ in the set $I$ indexing the charts, the section $s_i$ is transverse to the zero section, then a naive Kuranishi atlas in the sense above gives $X$ the structure of a smooth manifold. We will explore the deficiencies of this definition when the sections $s_i$ are not transverse in great detail shortly, but it does have a distinct virtue: for $P$ an elliptic PDE on a compact manifold, the space $\mathsf{Sol}(P)$ can be endowed with a naive Kuranishi atlas in a fairly explicit manner.
\begin{cons}[Rough sketch]\label{cons:sumchart}
Let $x\in \mathsf{Sol}(P)$ be a solution, then Fact \ref{fact:kur} asserts the existence of an affine Kuranishi model $\mathbf{K}_x=(M_x,E_x,s_x)$ and a homeomorphism $Z(s_x)\cong U_x$, where $x\subset U_x\subset \mathsf{Sol}(P)$ is some open neighbourhood. This Kuranishi model is obtained as follows: let $TP_x$ be the differential of $P$ at $x$ and choose a linear splitting $t_x:\coker TP_x\rightarrow \Gamma(F)$ of the projection $\Gamma(F)\rightarrow \coker TP_x$. The $\cinfty$ map (of infinite dimensional manifolds)
\[P+t_x:\Gamma(V)\times \coker TP_x\longrightarrow \Gamma(F)\] has surjective differential at $x$ by construction, so we conclude (after replacing all spaces of $\cinfty$ sections with spaces of Sobolev sections of sufficiently high regularity) that the zero set of $P+t_x$ around $x$, abusively denoted $Z(P+t_x)\subset \Gamma(V)\times \coker TP_x$, is a manifold. It follows that the zero set of $P$ around $x$ coincides with the zero set of the induced map 
\[f_x:Z(P+t_x)\subset \Gamma(V)\times \coker TP_x\longrightarrow  \coker TP_x\]
of \emph{finite dimensional manifolds}, and we thus have a Kuranishi neighbourhood $\mathbf{K}_x=(Z(P+t_x),Z(P+t_x)\times \coker TP_x,f_x+\mathrm{id})$ with trivial bundle. Now if $U_x\cap U_y\neq \emptyset$ for some other affine Kuranishi model $\mathbf{K}_y=(Z(P+t_y),Z(P+t_y)\times \coker TP_y,f_y+\mathrm{id})$ constructed in this manner associated with another solution $y\in \mathsf{Sol}(P)$, we can consider the map 
\[P+t_x+t_y:\Gamma(V)\times\coker TP_x\times\coker TP_y\longrightarrow\Gamma(F)\]
whose differential is surjective on the overlap $U_x\cap U_y$. This yields the \emph{sum chart} $\mathbf{K}_{x,y}=(Z(P+t_x+t_y),(Z(P+t_x+t_y)\times \coker TP_x\times\coker TP_y,f_{x,y}+\mathrm{id})$. The obvious maps
\[\Gamma(V)\times\coker TP_x\times\{0\} \longrightarrow  \Gamma(V)\times\coker TP_x\times\coker TP_y,\]and
\[  \Gamma(V)\times \{0\}\times\coker TP_y \longrightarrow  \Gamma(V)\times\coker TP_x\times\coker TP_y \]
induce maps 
\[ \mathbf{K}_x|_{Z(P+t_x+t_y)} \overset{\simeq}{\longrightarrow} \mathbf{K}_{x,y}\overset{\simeq}{\longleftarrow}  \mathbf{K}_y|_{Z(P+t_x+t_y)} \]
which are weak equivalences (but not isomorphisms!) of affine Kuranishi models by construction. Recalling the definition of morphisms in $h\mathsf{AffKur}$ via zigzags, we see that we have a well defined transition map $\phi_{x,y}:\mathbf{K}_x|_{U_{x}\cap U_y}\rightarrow \mathbf{K}_y|_{U_{x}\cap U_y}$ in the homotopy category. One observes that the cocycle condition holds over triple overlaps $U_x\cap U_y\cap U_z$ by constructing a triple sum chart $\mathbf{K}_{x,y,z}$ and noting that the equivalence relation on zigzags implies that both transition maps $\phi_{y,z}\circ \phi_{x,y}$ and $\phi_{x,z}$ are equivalent to the zigzag $\mathbf{K}_x|_{Z(P+t_x+t_y+t_Z)}  \overset{\simeq}{\rightarrow} \mathbf{K}_{x,y,z}\overset{\simeq}{\leftarrow}  \mathbf{K}_y|_{Z(P+t_x+t_y+t_Z)}$.
\end{cons}
We have satisfied $G1$: moduli spaces of elliptic PDEs can be endowed with naive Kuranishi atlases. For $G2$, we should glue together the various tangent complexes $\tanc_{\mathbf{K}_i}\in \mathbf{D}(\shv_{\mathrm{Vect}_{\R}}(U_i))$. To facilitate this gluing process, it's important to understand how the tangent complex is functorial in the transition maps of the Kuranishi atlas. To this end, it is convenient to recast Incorrect Definition \ref{defn:kuratlas} somewhat differently, which will require us to introduce some terminology from category theory.
\begin{fact}[Grothendieck construction]\label{fact:groth}
Suppose we are given a functor
\[ F:\icat^{op}\longrightarrow \mathsf{Cat} \]
from the opposite of a category $\icat$ to the category of categories\footnote{Since the category of categories is naturally a 2-category, we will allow \emph{pseudofunctors}, that is, `functors' that preserve associativity only up to invertible natural transformation. We will ignore this point for now}. We can repackage this data into a single category, denoted $\int_{\icat}F$, as follows:
\begin{enumerate}
    \item[$(O)$] Objects are pairs $(C,X)$ with $C\in \icat$ and $X\in F(C)$.
    \item[$(M)$] Morphisms between pairs $(C,X)$ and $(D,Y)$ are pairs $(g,h)$ where $g:C\rightarrow D$ is a morphism in $\icat$ and $h:F(g)(Y)\rightarrow X$ is a morphism in $F(X)$.
\end{enumerate}
It's straightforward to figure out how composition of morphisms should work. Note that there is an obvious functor $\int_{\icat}F\rightarrow\icat$ sending $(C,X)$ to $C$. This is the \emph{Grothendieck construction on $F$}. If we are given a natural transformation $\alpha:F\Rightarrow G$ with domain some other functor $G:\icat^{op}\rightarrow\mathsf{Cat}$, that is, we are given a collection of functors $\alpha_C:F(C)\rightarrow G(C)$ such that the obvious squares commute, we can construct a functor $\int_{\icat}F\rightarrow\int_{\icat}G$ that fits into a diagram 
\[
\begin{tikzcd}
\int_{\icat}F \ar[dr] \ar[rr] && \int_{\icat}G \ar[dl] \\
& \mathsf{\icat}
\end{tikzcd}
\]
by simply applying $\alpha$ fibrewise. In this way, the Grothendieck construction determines a functor
\[ \fun(\icat^{op},\mathsf{Cat}) \longrightarrow \mathsf{Cat}_{/\icat}  \]
from the category of (pseudo)functors between $\icat^{op}$ and $\mathsf{Cat}$ to the category of categories \emph{over} $\icat$. This functor is not an equivalence of categories, but we can characterize its essential image precisely: given a functor $p:\icatd\rightarrow\icat$, we say that a map $\tilde{g}:D'\rightarrow D$ is \emph{$p$-Cartesian} if for each diagram of the form 
\[
\begin{tikzcd}
    & D'\ar[dr,"\tilde{g}"] \\
    D''\ar[rr,"\tilde{h}"] && D
\end{tikzcd}
\]
in $\icatd$ and each map $f:p(D'')\rightarrow p(D)$ that makes the diagram 
\[
\begin{tikzcd}
    & p(D')\ar[dr,"p(\tilde{g})"] \\
    p(D'')\ar[ur,"f"]\ar[rr,"p(\tilde{h})"] && p(D)
\end{tikzcd}
\]
in $\icat$ commute, there is a map $\tilde{f}:D'\rightarrow D''$ which satisfies $p(\tilde{f})=f$, a \emph{lift} of $f$, that makes the diagram 
\[
\begin{tikzcd}
    & D'\ar[dr,"\tilde{g}"] \\
    D''\ar[ur,"\tilde{f}"]\ar[rr,"\tilde{h}"] && D
\end{tikzcd}
\]
in $\icatd$ commute. A functor $p:\icatd\rightarrow\icat$ is a \emph{Grothendieck fibration}, or a \emph{Cartesian fibration} if for each morphism $f:C'\rightarrow C$ and each $D\in\icatd$ for which $p(D)=C$, there exists a $p$-Cartesian lift $\tilde{f}:D'\rightarrow D$ for which $p(\tilde{f})=f$. The essential image of the Grothendieck construction for functors $F:\icat^{op}\rightarrow\cat$ consists precisely of the Cartesian fibrations over $\icat$.
\end{fact}
\begin{ex}
The canonical example is the following. Let $\mathsf{CAlg}^0_k$ be the category of commutative $k$-algebras for $k$ a commutative ring\footnote{The superscript 0 indicates that these are `orindary' commutative $k$-algebras instead of the `homotopy coherently commutative higher $k$-algebras' that will concern us later}, then we have for every map $f:A\rightarrow B$ in $\mathsf{CAlg}^0_k$, a functor $f_*:\mathrm{Mod}_B \rightarrow \mathrm{Mod}_A$ by viewing a $B$-module as an $A$-module. These functors assemble and define a functor
\[ (\mathsf{CAlg}^0_k)^{op} \longrightarrow \mathsf{Cat},\quad A\longmapsto \mathrm{Mod}_A.   \]
According to Fact \ref{fact:groth}, we have an associated Grothendieck fibration $q:\mathrm{Mod}\rightarrow \mathsf{CAlg}^0_k$: the category $\mathrm{Mod}$ has as objects pairs $(A,M)$ for $A$ a $k$-algebra and $M$ an $A$-module. A morphism $(A,M)\rightarrow (B,N)$ is $q$-Cartesian if the map of $M\rightarrow N$ of $A$-modules is an isomorphism.
\end{ex}
\begin{ex}
Consider the category $\mathsf{Top}^{\mathsf{open}}$ whose objects are paracompact Hausdorff spaces and whose morphisms are open topological embeddings of such spaces. Similarly, let $h\mathsf{AffKur}^{\mathrm{open}}$ be the subcategory of $h\mathsf{AffKur}$ on the morphisms $f:\mathbf{J}\rightarrow\mathbf{K}$ such that $Z(f)$ is an open topological embedding. We claim that the functor
\[Z:h\mathsf{AffKur}^{\mathrm{open}}\longrightarrow \mathsf{Top}^{\mathsf{open}}  \]
taking the underlying topological space of an affine Kuranishi model is a Grothendieck fibration. To substantiate this claim, we have to provide a sufficient supply of $Z$-Cartesian morphisms. Unwinding the definitions, we observe that a morphism $\mathbf{H}\rightarrow \mathbf{K}$ is $Z$-Cartesian \emph{precisely} if this morphism exhibits $\mathbf{H}$ as a localization of $\mathbf{K}$ with respect to the image of the open embedding $Z(\mathbf{H})\rightarrow Z(\mathbf{K})$. Therefore the assertion that for every open embedding $i:U\rightarrow Y$ and every $\mathbf{K}$ with $Z(\mathbf{K})=Y$, there is a $Z$-Cartesian morphism lying above $i$ and terminating at $\mathbf{K}$ is equivalent to Lemma \ref{lem:localizationsexist1}. According to Fact \ref{fact:groth}, the fibration $Z$ corresponds to the functor $(\mathsf{Top}^{\mathrm{open}})^{op}\rightarrow \cat$ that carries a space $X$ to the category of affine Kuranishi model structures on $X$.
\end{ex}
Having expanded our categorical vocabulary, we can think of a naive Kuranishi atlas as a functor.
\begin{incdef}\label{defn:kuratlas2}
 Let $X\in \mathsf{Top}^{\mathsf{open}}$, then a \emph{naive Kuranishi atlas} on $X$ consists of the following data.
\begin{enumerate}[$(a')$]
    \item A collection of maps $\{V_i\rightarrow X\}_{i\in I}$ in $\mathsf{Top}^{\mathrm{open}}$ of $X$ with images being open sets $\{U_i\subset X\}$ that cover $X$. We can view this data as a functor $\mathfrak{U}:I\rightarrow \mathsf{Top}^{\mathsf{open}}$ from the set $I$ viewed as a category with only identity morphisms. Consider the poset 
    \[P^{\leq 3}_I:=\{J\subset I;\,J\neq \emptyset,\,|J|\leq 3 \} \] 
    of nonempty subsets of $I$ of cardinality at most 3 ordered by reverse inclusion, then the functor $\mathfrak{U}$ induces a functor
    \[  f:P^{\leq 3}_I\longrightarrow \mathsf{Top}^{\mathsf{open}} \]
    which sends $J$ to the limit of the diagram $J\subset I\rightarrow \mathsf{Top}^{\mathrm{open}}_{/X}$.
    \item A dotted lift $\tilde{f}$ of $f$ as follows
    \[
    \begin{tikzcd}
    & h\mathsf{AffKur}^{\mathrm{open}}\ar[d,"Z"] \\
    P^{\leq 3}_I \ar[r,"f"'] \ar[ur,dotted,"\tilde{f}"] & \mathsf{Top}^{\mathsf{open}}
    \end{tikzcd}
    \]
    that makes the diagram of categories (strictly) commute. Moreover, we require that $\tilde{f}$ carries every morphism in $P^{\leq 3}_I$ to a Cartesian morphism with respect to $Z$.
\end{enumerate}
\end{incdef}
The two definitions of a naive Kuranishi atlas are equivalent: given a naive Kuranishi atlas in the sense of Definition \ref{defn:kuratlas2}, the restriction of $\tilde{f}$ to $I$ determines a collection of affine Kuranishi models $\{\mathbf{K}_i\}_{i\in I}$ with homeomorphisms $\psi_{i}:Z(\mathbf{K}_i)=V_i\cong U_i\subset X$. For every nonempty intersection $U_i\cap U_j$ with $i\neq j$, we have maps
\[ \mathbf{K}_{i}|_{\psi_{i}^{-1}(U_i\cap U_j)}\longleftarrow \tilde{f}(\{i,j\}) \longrightarrow \mathbf{K}_{j}|_{\psi_{j}^{-1}(U_i\cap U_j)} \]
induced by the subset inclusions $\{i\}\subset \{i,j\}\supset \{j\}$. These maps are isomorphisms because $\tilde{f}$ carries all morphisms to Cartesian morphisms, giving us the isomorphisms $\{\phi_{ij}\}$. The compatibility conditions $(1)$ and $(2)$ are guaranteed by the fact that the diagram of $(b')$ commutes and that $\tilde{f}$ is a functor. Conversely, from datum $(a)$ of Definition \ref{defn:kuratlas} we can construct a functor as in $(a')$, and given data $(b)$ through $(d)$ satisfying $(1)$ and $(2)$, it is possible to construct a lift $\tilde{f}$, and it can be shown that for a suitable choice of morphisms between naive Kuranishi atlases for both of the definitions we have given, this correspondence determines an equivalence of categories. Using the language of Grothendieck fibrations and Incorrect Definition \ref{defn:kuratlas2}, we can make explicit the functoriality of taking the tangent complex of affine Kuranishi models. 
\begin{cons}\label{cons:dercatfib}
We have a pseudofunctor
\[(\mathsf{Top}^{\mathrm{open}})^{op} \longrightarrow \mathsf{Cat},\quad X\longmapsto \mathbf{D}(\shv_{\mathrm{Vect}_{\R}}(X))^{op} \]
that carries each space $X$ to the \emph{opposite} of the derived category of sheaves of chain complexes of $\R$-vector spaces on $X$. We can apply the Grothendieck construction to this pseudofunctor, obtaining a functor
\[ \int_{\mathsf{Top}^{\mathrm{open}}}\mathbf{D}(\shv_{\mathrm{Vect}_{\R}}(\_))^{op} \longrightarrow \mathsf{Top}^{\mathrm{open}}. \]
Concretely, the category $\int_{\mathsf{Top}^{\mathrm{open}}}\mathbf{D}(\shv_{\mathrm{Vect}_{\R}}(\_))^{op}$ is given as follows. 
\begin{enumerate}[$(1)$]
    \item Objects are pairs $(U,\F)$, where $U\in \mathsf{Top}^{\mathrm{open}}$ and $\F$ is a complex of sheaves on $U$. 
    \item Morphisms $(U,\F)\rightarrow (V,\mathcal{G})$ are maps $i:U\subset V$ together with a map $\F\rightarrow i^*\mathcal{G}$ in $\mathbf{D}(\shv_{\mathrm{Vect}_{\R}}(U))$. 
\end{enumerate}
The \emph{tangent complex functor} $\tanc$ is the functor that carries an object $\mathbf{K}$ to the pair $(Z(\mathbf{K}),\tanc_{\mathbf{K}})$. Note that by definition of a weak equivalence, this assignment is well defined up to equivalence. The tangent complex functor fits into a commuting diagram 
\[
\begin{tikzcd}
h\mathsf{AffKur}^{\mathrm{open}} \ar[rr,"\tanc"] \ar[dr,"Z"']&&\int_{\mathsf{Top}^{\mathrm{open}}}\mathbf{D}(\shv_{\mathrm{Vect}_{\R}}(\_))^{op} \ar[dl,"p"] \\ & \mathsf{Top}^{\mathrm{open}}
\end{tikzcd}
\]
and carries Cartesian morphisms with respect to $Z$ to Cartesian morphisms with respect to $p$.
\end{cons}
More explicitly, the tangent complex assigns to each affine Kuranishi model $\mathbf{K}$ its tangent complex $\tanc_{\mathbf{K}}$, viewed as an object in the derived category, and it assigns to each morphism $\mathbf{J}\rightarrow\mathbf{K}$ such that $i:Z(\mathbf{J})\rightarrow Z(\mathbf{K})$ is an open inclusion, the natural morphism $\tanc_{\mathbf{J}}\rightarrow i^*\tanc_{\mathbf{K}}$ in the derived category. If we have a naive Kuranishi atlas on $X$ for some open cover $\{U_i\subset X\}$, that is, we have a commuting diagram of categories
\[
    \begin{tikzcd}
    & h\mathsf{AffKur}^{\mathrm{open}}\ar[d,"Z"] \\
    P^{\leq 3}_I \ar[r,"f"'] \ar[ur,dotted,"\tilde{f}"] & \mathsf{Top}^{\mathsf{open}}
    \end{tikzcd}
\]
such that $\tilde{f}$ carries every $J\subset J'$ to a localization, then we can compose $\tilde{f}$ with the functor $\tanc$ to get a commuting diagram
\[
    \begin{tikzcd}
    &\int_{\mathsf{Top}^{\mathrm{open}}}\mathbf{D}(\shv_{\mathrm{Vect}_{\R}}(\_))^{op}\ar[d,"p"] \\
    P^{\leq 3}_I \ar[r,"f"'] \ar[ur,dotted,"\tanc\circ\tilde{f}"] & \mathsf{Top}^{\mathsf{open}}
    \end{tikzcd}
\]
such that $\tilde{f}$ carries every $J\subset J'$ to a map $(f(J'),\F)\rightarrow (f(J),\mathcal{G})$ for which $\F\rightarrow i^*\mathcal{G}$ is an isomorphism in the derived category. In exact analogy with the equivalence between Incorrect Definitions \ref{defn:kuratlas} and \ref{defn:kuratlas2}, such a diagram consists precisely of the following data. 
\begin{enumerate}[$(a)$]
    \item The open cover $\{U_i\subset X\}_{i\in I}$.
    \item A collection of objects $\F_i\in \mathbf{D}(\shv_{\mathrm{Vect}_{\R}}(U_i))$.
    \item For each pair of indices $i,j\in I$ such that $U_{ij}=U_i\cap U_j$ is nonempty, an isomorphism $\phi_{ij}:\F_i|_{U_{ij}}\cong \F_j|_{U_{ij}}$ \emph{in the category $\mathbf{D}(\shv_{\mathrm{Vect}_{\R}}(U_{ij}))$}, where $\F_i|_{U_{ij}}=i^*\F_i$ for $i$ the inclusion $U_{ij}\subset U_i$, and similarly for $\F_j|_{U_{ij}}$.
\end{enumerate}
Such that the cocycle conditions holds.
\begin{itemize}
    \item For every triple $i,j,k\in I$ such that $U_{ijk}:=U_i\cap U_j\cap U_k$ is nonempty, the equality 
    \[ \phi_{jk}|_{U_{ijk}}\circ \phi_{ij}|_{U_{ijk}} = \phi_{ik}|_{U_{ijk}}\] 
    \emph{as morphisms of $\mathbf{D}(\shv_{\mathrm{Vect}_{\R}}(U_{ijk}))$}.
\end{itemize}
Cast in this form, the data above should be very familiar; as part of an algebraic geometry or general sheaf theory course, one learns how to \emph{glue} sheaves: given a collection of sheaves $\{\F_i\}_i$ (of $A$-modules for a commutative ring $A$, say) on an open cover $\{U_i\subset X\}_i$ and transition isomorphisms $\phi_{ij}:F_i|_{U_{ij}}\cong F_j|_{U_{ij}}$ that satisfy the cocycle condition, there exists a sheaf $\F$ on $X$ such that $\F|_{U_i}\cong \F_i$. We find it convenient to rephrase this gluing principle as follows.
\begin{fact}[Classical descent]\label{fact:descent}
Let $X$ be a topological space with an open cover $\{U_i\subset X\}_{i\in I}$ determining the diagram 
 \[  f:P^{\leq 3}_I\longrightarrow \mathsf{Top}^{\mathsf{open}} \]
as in Definition \ref{defn:kuratlas2}, and let $\mathcal{A}$ be a Grothendieck abelian category, such as $\mathrm{Mod}_A$ for a $A$ a unital commutative ring. There is a canonical equivalence between the category $\shv_{\mathcal{A}}(X)^{op}$, the opposite of the category of sheaves of objects of $\mathcal{A}$, and the lifts
\[
    \begin{tikzcd}
    & \int_{\mathsf{Top}^{\mathrm{open}}}\shv_{\mathcal{A}}(\_)^{op}\ar[d,"p"] \\
    P^{\leq 3}_I \ar[r,"f"'] \ar[ur,dotted,"\F"] & \mathsf{Top}^{\mathsf{open}}
    \end{tikzcd}
\]
such that $\F$ sends all morphisms to Cartesian morphisms with respect to $p$. This equivalence is implemented by the `coCartesian pushforward', which applies to each $\F(J)$ the functor $i_!$ (the direct image of sheaves) induced by the map $i:f(J)\subset X$, resulting in a diagram  $P^{\leq 3}_I\rightarrow \shv_{\mathcal{A}}(X)^{op}$, and takes the colimit (that is, the limit in $\shv_{\mathcal{A}}(X)$).
\end{fact}
We are not in the above situation. Our transition isomorphisms are not isomorphisms of $\R$-vector spaces, but isomorphisms in the derived category of sheaves of chain complexes of $\R$-vector spaces, and our cocycle condition holds only `up to homotopy'. This is so because the tangent complex determines an honest sheaf of $\R$-vector spaces on the zero locus of $(X,E\rightarrow X,s)$ if and only if $s$ is transverse to the zero section. Unfortunately, this precludes the gluing of sheaves.
\begin{prob}\label{prob:notan}
The statement of Fact \ref{fact:descent} for $\shv_{\mathcal{A}}(\_)$ replaced with $\mathbf{D}(\shv_{\mathcal{A}}(\_))$ is false.
\end{prob}
\begin{rmk}
What goes wrong in our case? While we can `push the diagram forward', the derived category $\mathbf{D}(\shv_{\mathcal{A}}(X))$ has very few limits (while the category $\shv_{\mathcal{A}}(X)$ has all limits and colimits), so we might not be able to perform the gluing construction. Even more seriously, if a limit exists, the resulting object $\tanc_X$ need not have the property that $\tanc_X|_{U_i}$ is isomorphic in the derived category to $\tanc_{\mathbf{K}_i}$, so $\tanc_X$ would not deserve to be called `the gluing' of the complexes of sheaves $\tanc_{\mathbf{K}_i}$.
\end{rmk} 
\subsubsection{Attempt 2: strict atlases}
As the previous remark explains, the lack of well behaved limits in the derived category is the source of our current problem. To glue sheaves up to homotopy, we must somehow specify additional information that allows us to `lift' the diagram $P^{\leq 3}_I\rightarrow \mathbf{D}(\shv_{\mathcal{A}}(X))^{op}$ to the category $\shv_{\mathrm{Ch}(\mathcal{A})}(X)$ and take a suitable limit there. Here is one way to achieve this.
\begin{fact}[Strict descent]\label{fact:descent2}
Let $\mathcal{A}$ be a Grothendieck abelian category. Consider the functor
\[ p:\int_{\mathsf{Top}^{\mathrm{open}}}\shv_{\mathrm{Ch}(\mathcal{A})}(\_)^{op}\longrightarrow \mathsf{Top}^{\mathsf{open}}, \]
the Grothendieck fibration associated to the functor
\[ \mathsf{Top}^{\mathsf{open}}\longrightarrow \mathsf{Cat},\quad\quad X\longmapsto \shv_{\mathrm{Ch}(\mathcal{A})}(X)^{op}. \]
Note that since we work with chain complexes, we now have two different notions of a `Cartesian morphism'.
\begin{enumerate}[$(i)$]
    \item A morphism $(U,\F)\rightarrow (V,\mathcal{G})$ given by a pair $(i:U\subset V,\F\rightarrow i^*\mathcal{G})$ is \emph{Cartesian} if the map $\F\rightarrow i^*\mathcal{G}$ is an isomorphism of chain complexes of sheaves on $U$.
    \item A morphism $(U,\F)\rightarrow (V,\mathcal{G})$ is \emph{homotopy Cartesian} if the map $\F\rightarrow i^*\mathcal{G}$ is a quasi-isomorphism of chain complexes of sheaves on $U$.
\end{enumerate}
Let $X$ be a topological space with an open cover $\{U_i\subset X\}_{i\in I}$ determining the diagram 
 \[  f:P_I\longrightarrow \mathsf{Top}^{\mathsf{open}}, \]
 where now $P_I$ is the partially ordered set of \emph{all} nonempty finite subsets of $I$, so that $f$ carries a subset $J=\{i_1,\ldots,i_n\}\subset I$ to the $n$-fold intersection $U_{i_1}\cap\ldots\cap U_{i_n}$. Consider the category of lifts
\[
    \begin{tikzcd}
    & \int_{\mathsf{Top}^{\mathrm{open}}}\shv_{\mathrm{Ch}(\mathcal{A})}(\_)^{op}\ar[d,"p"] \\
    P_I \ar[r,"f"'] \ar[ur,dotted,"\F"] & \mathsf{Top}^{\mathsf{open}}
    \end{tikzcd}
\]
which we will denote by \[\fun_{\mathsf{Top}^{\mathsf{open}}}(P_I,\int_{\mathsf{Top}^{\mathrm{open}}}\shv_{\mathrm{Ch}(\mathcal{A})}(\_)^{op}),\]
the functor category \emph{over $\mathsf{Top}^{\mathsf{open}}$}. We say that a morphism $\F\rightarrow \mathcal{G}$ in this category, that is, a natural transformation over $\mathsf{Top}^{\mathsf{open}}$, is a \emph{weak equivalence} if for each nonempty finite $J\subset  I$, the map $\F(J)\rightarrow \mathcal{G}(J)$ in $\shv_{\mathrm{Ch}(\mathcal{A})}(f(J))^{op}$ is a quasi-isomorphism. With this collection $W$ of weak equivalences, the pair
\[(\fun_{\mathsf{Top}^{\mathsf{open}}}(P_I,\int_{\mathsf{Top}^{\mathrm{open}}}\shv_{\mathrm{Ch}(\mathcal{A})}(\_)^{op}),W)  \]
is a category with weak equivalences, so we can take its homotopy category 
\[h\fun_{\mathsf{Top}^{\mathsf{open}}}(P_I,\int_{\mathsf{Top}^{\mathrm{open}}}\shv_{\mathrm{Ch}(\mathcal{A})}(\_)^{op}). \]
Then there is a canonical equivalence between
\begin{enumerate}[$(a)$]
    \item The full subcategory of the above homotopy category spanned by lifts $\F$ that carry every morphism $J\subset J'$ to a homotopy Cartesian morphism, in the sense of $(ii)$ defined above (note that the requirement that a lift $\F$ carries every morphism to a `usual' Cartesian morphism in the sense of $(i)$ is not stable under weak equivalence and therefore not well defined in the homotopy category of lifts).
    \item The opposite derived category $\mathbf{D}(\shv_{\mathcal{A}}(X))^{op}$.
\end{enumerate}
This equivalence is implemented by the `coCartesian pushforward', which applies to each $\F(J)$ the functor $i_!$ induced by the map $i:f(J)\subset X$, resulting in a diagram  $P_I\rightarrow \shv_{\mathcal{A}}(X)^{op}$, and takes the \emph{homotopy colimit}.
\end{fact}
\begin{rmk}
Schematically, the fact above says something like this:
\begin{itemize}
    \item For an open cover $\{U_i\subset X\}$, the derived category $\mathbf{D}(\shv_{\mathcal{A}}(X))$ is not the category of compatible collections of sheaves $\{\F_i\}$ in the derived categories $\mathbf{D}(\shv_{\mathcal{A}}(U_i))$, rather it is the derived category of \emph{homotopy} compatible collections of sheaves $\{\F_i\}$ in the categories $\shv_{\mathrm{Ch}(\mathcal{A})}(U_i)$.
\end{itemize}
Apart from some notational overhead in the statement of Fact \ref{fact:descent2}, Facts \ref{fact:descent} and \ref{fact:descent2} are very similar. There are only two important differences, which we highlight.
\begin{enumerate}[$(1)$]
    \item In Fact \ref{fact:descent}, we use the poset $P_{I}^{\leq 3}$, but in Fact \ref{fact:descent2}, we use the larger poset $P_I$.
    \item In Fact \ref{fact:descent}, we use a colimit to glue the local pieces, but in Fact \ref{fact:descent2}, we use a homotopy colimit.
\end{enumerate}
These points are related. Recall that for a diagram $\icate\rightarrow \icat$ to a category $\icat$, the colimit functor is the left adjoint $\colim:\fun(\icate,\icat)\rightarrow \icat$ to the constant diagram functor $\icat\rightarrow \fun(\icate,\icat)$, where the latter functor is defined via its adjoint, the obvious projection $\icat\times\icate\rightarrow\icat$. If $\icat=h\icatd$, the homotopy category of some category with weak equivalences $(\icatd,W)$, then the usual colimit is not the correct notion to use, as Problem \ref{prob:notan} suggests. What's more, the formation of colimits in $\icatd$ is also deficient, as this operation might not be homotopically meaningful: we can make the functor category $\fun(\icate,\icatd)$ into a category with weak equivalences by declaring a natural transformation $F\rightarrow G$ a weak equivalence just in case $F(E)\rightarrow G(E)$ is a weak equivalence for all $E\in \icate$, but the colimit functor $\colim:\fun(\icate,\icatd)\rightarrow\icatd$ will usually not preserve these weak equivalences, so equivalent diagrams may not have equivalent colimits. The \emph{homotopy colimit functor} (if it exists) is instead the left adjoint in the adjunction \[\begin{tikzcd}h\fun(\icate,\icatd)\ar[r,shift left,"\hocolim"]&[3em]  h\icatd \ar[l,shift left,"\mathrm{constant}"],\end{tikzcd}\] where $h\fun(\icate,\icatd)$ is the homotopy category for the `objectwise' weak equivalences we just defined. In the context of Fact \ref{fact:descent2}, we see that taking the homotopy colimit is the only meaningful operation we have available to identify the two homotopy categories in the statement. This also explains the need to shift from $P_I^{\leq 3}$ to $P_I$: gluing up to homotopy in general requires information given on \emph{all} possible overlaps. If we had used $P_I$ in Fact \ref{fact:descent}, it wouldn't have made a difference, but the homotopy colimits over $P_I$ and $P_I^{\leq 3}$ are \emph{not the same}\footnote{Technically, the inclusion $P_I^{\leq 3}\subset P_I$ is cofinal (the poset $P_J^{\leq 3}$ for any set $J$ is connected), but not homotopy cofinal (the poset $P_J^{\leq 3}$ is not weakly contractible in general)}.
\end{rmk}
Since we wish the correct definition of a Kuranishi structure to determine a global tangent complex, a glance at the diagram of Fact \ref{fact:descent2} immediately suggests the following definition.
\begin{incdef}\label{defn:kuratlas4}
Let $X\in \mathsf{Top}^{\mathsf{open}}$, then a \emph{strict Kuranishi atlas} on $X$ consists of the following data.
\begin{enumerate}[$(a)$]
    \item A collection of maps $\{V_i\rightarrow X\}_{i\in I}$ in $\mathsf{Top}^{\mathrm{open}}$ of $X$ with images being open sets $\{U_i\subset X\}$ that cover $X$, which induces a functor
    \[  f:P_I\longrightarrow \mathsf{Top}^{\mathsf{open}}. \]
    \item A dotted lift $\tilde{f}$ of $f$ as follows
    \[
    \begin{tikzcd}
    & \mathsf{AffKur}^{\mathrm{open}}\ar[d,"Z"] \\
    P_I \ar[r,"f"'] \ar[ur,dotted,"\tilde{f}"] & \mathsf{Top}^{\mathsf{open}}
    \end{tikzcd}
    \]
    that makes the diagram of categories commute. Moreover, we require that $\tilde{f}$ carries every morphism in $P_I$ to a homotopy Cartesian morphism with respect to $Z$, that is, the map $\tilde{f}(J')\rightarrow\tilde{f}(J)$ determines a localization in the homotopy category $h\mathsf{AffKur}$ for each $J\subset J'$.
\end{enumerate}
\end{incdef}
We have a tangent complex functor as in the commuting diagram
\[
\begin{tikzcd}
\mathsf{AffKur}^{\mathrm{open}} \ar[rr,"\tanc"] \ar[dr,"Z"']&&\int_{\mathsf{Top}^{\mathrm{open}}}\shv_{\mathrm{Ch}_{\R}}(\_)^{op} \ar[dl,"p"] \\ & \mathsf{Top}^{\mathrm{open}}
\end{tikzcd}
\]
which carries weak equivalences to weak equivalences (the tangent complex functor of Construction \ref{cons:dercatfib} is obtained from this one by passing to homotopy categories). 
Given a strict Kuranishi atlas on a space $X$, we can compose the lift $\tilde{f}:P_I\rightarrow\mathsf{AffKur}^{\mathrm{open}}$ with the tangent complex functor and obtain a diagram
\[
    \begin{tikzcd}
    &\int_{\mathsf{Top}^{\mathrm{open}}}\shv_{\mathrm{Ch}_{\R}}(\_)^{op}\ar[d,"p"] \\
    P_I \ar[r,"f"'] \ar[ur,dotted,"\tanc\circ\tilde{f}"] & \mathsf{Top}^{\mathsf{open}}.
    \end{tikzcd}
\]
By construction $\tanc\circ\tilde{f}$ carries every morphism of $P_I$ to a homotopy Cartesian morphism, so $\tanc\circ\tilde{f}$ determines an object in the full subcategory of the homotopy category $h\fun_{\mathsf{Top}^{\mathsf{open}}}(P_I,\int_{\mathsf{Top}^{\mathrm{open}}}\shv_{\mathrm{Ch}(\mathcal{A})}(\_)^{op})$ specified in point $(a)$ of Fact \ref{fact:descent2} and we obtain a global object $\tanc_{X}$ that glues the local data appropriately. We seem to have made some progress, but we are still a ways away from achieving our goal of satisfying the desiderata $G1$ and $G2$. To see where the problem lies, let us unpack Incorrect Definition \ref{defn:kuratlas4} in more familiar terms, along the lines of our initial attempt. A strict Kuranishi atlas on $X$ consists of 
\begin{enumerate}[$(a)$]
   \item An open covering $\{U_i\rightarrow X\}_{i\in I}$ of $X$.
    \item For every finite $J\subset I$, an affine Kuranishi model $\mathbf{K}_J$ with zero locus $Z(\mathbf{K}_J)$. 
    \item A collection of homeomorphisms $\psi_J:Z(\mathbf{K}_J)\rightarrow U_J$, where $U_J=\cap_{j\in J}U_{j}$.
    \item For every inclusion $J\subset J'$, a morphism $\phi_{JJ'}:\mathbf{K}_{J'}\rightarrow \mathbf{K}_{J}$ of affine Kuranishi models determining a localization in the homotopy category $h\mathsf{AffKur}$.
\end{enumerate}
These data are required to satisfy the following conditions.
\begin{enumerate}[$(1)$]
    \item The transition maps $\phi_{JJ'}$ are compatible with the footprint maps: $ \psi_{J}\circ Z(\phi_{JJ'})=\psi_{J'}$.
    \item The cocycle condition holds: for every triple $J\subset J'\subset J''$ we have 
    \[ \phi_{JJ'}\circ\phi_{J'J''} = \phi_{JJ''}. \]
\end{enumerate}
\begin{prob}\label{prob:nokur}
It is not clear how to endow the space $\mathsf{Sol}(P)$ with a strict Kuranishi atlas.
\end{prob}
\begin{rmk}\label{rmk:nokur}
To see what fails, recall that Construction \ref{cons:sumchart} yields a transition map \emph{in the wrong direction}: from an affine Kuranishi model on the smaller open to the larger one. For Definition \ref{defn:kuratlas}, this did not cause problems because we were working in the homotopy category whose morphisms were equivalence classes of zigzags, but in the context of Definition \ref{defn:kuratlas4} where we require maps in the category $\mathsf{AffKur}$, this issue becomes insurmountable.
\end{rmk}
\subsubsection{A first look at homotopy coherence}
Let us recapitulate: the moduli spaces we wish to study do fall in the framework of our equivalent Incorrect Definitions \ref{defn:kuratlas} and \ref{defn:kuratlas2}, but such a structure is \emph{not strong enough} to admit a global tangent complex. On the other hand, an atlas in the sense of Incorrect Definition \ref{defn:kuratlas4} does admit a tangent complex, but this notion is \emph{too strong} to encompass the examples of interest. Evidently, we need to find a definition that threads the needle. To understand what such an intermediate notion should look like, we return to the problem of gluing sheaves of chain complexes up to homotopy. As we explained in Fact \ref{fact:descent2} and the remark that follows, to glue a diagram 
\begin{equation}\label{eq:hocom}\F:P_I\longrightarrow \int_{\mathsf{Top}^{\mathrm{open}}}\mathbf{D}(\shv_{\mathcal{A}}(\_))^{op}\end{equation}
of objects in the various derived categories, it suffices to lift $\F$ to a diagram in the various ordinary categories
\begin{equation}\label{eq:strict}P_I\longrightarrow \int_{\mathsf{Top}^{\mathrm{open}}}\shv_{\mathrm{Ch}(\mathcal{A})}(\_)^{op}.\end{equation}
The diagram \eqref{eq:hocom} is an example of a \emph{homotopy} descent datum, while \eqref{eq:strict} is a \emph{strict} descent datum. There exists an intermediate notion, that of a \emph{homotopy coherent} descent datum. Very roughly, we expect that in order to produce a glued sheaf we must do more than merely demand the cocycle condition up to homotopy; we should rather \emph{choose} a collection of homotopies that `exhibit' the cocycle condition. These specified homotopies must themselves satisfy an associativity (\emph{coherence}) condition up to homotopy themselves, which involves choosing a collection of higher homotopies exhibiting \emph{those} associativity conditions, which must satisfy an even higher coherence condition up to homotopy, and so on. 
\begin{rmk}\label{rmk:chainhomotopies}
Paths and (higher) homotopies are algebro-topological notions and make sense in reference to topological spaces. What does it mean to have a `homotopy' $\phi_{ik}\circ \phi_{ij}\sim\phi_{ik}$ of maps of sheaves of chain complexes? In homological algebra, we can speak of a \emph{chain homotopy} between maps of complexes $f_0,f_1:\F_{\bullet}\rightarrow\mathcal{G}_{\bullet}$; this is a map $H_{01}:\F_{\bullet}\rightarrow \mathcal{G}_{\bullet +1}$ such that $dH=f_0-f_1$. Suppose we have three maps $f_0,f_1,f_2:\F_{\bullet}\rightarrow\mathcal{G}_{\bullet}$ and three homotopies $H_{01}:f_0\sim f_1$, $H_{02}:f_0\sim f_2$ and $H_{12}:f_1\sim f_2$, then a `2-dimensional' chain homotopy between these three homotopies is a map $K:\F_{\bullet}\rightarrow \mathcal{G}_{\bullet+2}$ such that $dK=H_{01}-H_{02}+H_{12}$ (note that $ddK=0$). It should now be obvious how to generalize to an $n$-dimensional chain homotopy between $n+1$ $(n-1)$-dimensional chain homotopies. For $\F$ and $\mathcal{G}$ sheaves of chain complexes on some topological space $X$, let us construct a CW-complex $|\Hom_{\shv_{\mathrm{Ch}_{\R}}(X)}(\F,\mathcal{G})|$ as follows:
\begin{enumerate}[$(a)$]
    \item $0$-cells are the maps $F_{\bullet}\rightarrow \mathcal{G}_{\bullet}$.
    \item $1$-cells are the chain homotopies $\F_{\bullet}\rightarrow\mathcal{G}_{\bullet+1}$.
    \item $2$-cells are the 2-dimensional chain homotopies $\F_{\bullet}\rightarrow\mathcal{G}_{\bullet+2}$.
    \item etc.
\end{enumerate}
We will give a more precise construction in the next subsection; what matters is that we have a topological space $|\Hom_{\shv_{\mathrm{Ch}_{\R}}(X)}(\F,\mathcal{G})|$ that encodes the homological information of maps $\F\rightarrow\mathcal{G}$. 
\end{rmk} 
What follows is a sketch of a definition, meant to indicate how the combinatorics of the higher coherences rapidly gets out of hand. A \emph{homotopy coherent descent datum} for a chain complex of sheaves consists of the following data.
\begin{enumerate}[$(1)$]
    \item A topological space $X$ with a cover $\{U_i\subset X\}_{i\in I}$. 
    \item For $\{i_1,\ldots,i_n\}\subset I$, write $U_{i_1,\ldots,i_n}$ for the $(n-1)$-fold intersection $U_{i_1}\cap\ldots\cap U_{i_n}$. Then we have for each 
such nonempty subset $\{i_1,\ldots,i_n\}\subset I$, a sheaf of chain complexes $\F_{i_1,\ldots,i_n}$ on $U_{i_1,\ldots,i_n}$.
    \item We have for each inclusion $\{i_1,\ldots,i_n\}\subset \{i_1,\ldots,i_n,i_{n+1}\}$ inducing an inclusion $U_{i_1,\ldots,i_n,i_{n+1}}\subset U_{i_1,\ldots,i_n}$ adding one more intersection, a quasi-isomorphism 
    \[  \phi^{i_1,\ldots,i_n,i_{n+1}}_{i_1,\ldots,i_n}:\F_{i_1,\ldots,i_n,i_{n+1}} \longrightarrow \F_{i_1,\ldots,i_n}|_{U_{i_1,\ldots,i_n,i_{n+1}}}.      \]
    \item For each triple $\{i,j,k\}$ corresponding to an intersection $U_{ijk}=U_i\cap U_j\cap U_k$, we demand the existence of another quasi-isomorphism
    \[ \phi_{i}^{ijk}:\F_{ijk}\longrightarrow\F_i|_{U_{ijk}}  \]
    and two homotopies
    \[\begin{tikzcd}&\phi_{i}^{ijk}\ar[dl] \ar[dr] \\ \phi^{ijk}_{ij}\circ \phi^{ij}_{i}  && \phi^{ijk}_{ik}\circ \phi^{ik}_{i}\end{tikzcd}\]
    between the three possible quasi-isomorphisms $\F_{ijk}\rightarrow \F_i|_{U_{ijk}}$ we are now given; these homotopies exhibit the cocycle condition for triple overlaps. In fact, for any inclusion $ U_{i_1,\ldots,i_n,i_{n+1},i_{n+2}}\subset U_{i_1,\ldots,i_n}$ adding two more intersections, we are given another quasi-isomorphism $\phi_{i_1,\ldots,i_n}^{i_1,\ldots,i_n,i_{n+1},i_{n+2}}$, and two homotopies between the three possible compositions $\F_{i_1,\ldots,i_n,i_{n+1},i_{n+2}}\rightarrow\F_{i_1,\ldots,i_n}|_{i_1,\ldots,i_n,i_{n+1},i_{n+2}}$ in the topological space $|\Hom_{\shv_{\mathrm{Ch}_{\R}}}(\F_{ijk},\F_i|_{U_{ijk}})|$.
    \item Suppose we have a quadruple intersection $U_i\cap U_j\cap U_k\cap U_l$, then the previous data $(1)$-$(4)$ provide twelve maps $\F_{ijkl}\rightarrow\F_i|_{U_{ijkl}}$, and homotopies between them that fit together to form the 1-dimensional polyhedron in $\R^2$ obtained from the boundary of a regular hexagon by taking a barycentric subdivision of each edge, like so
    
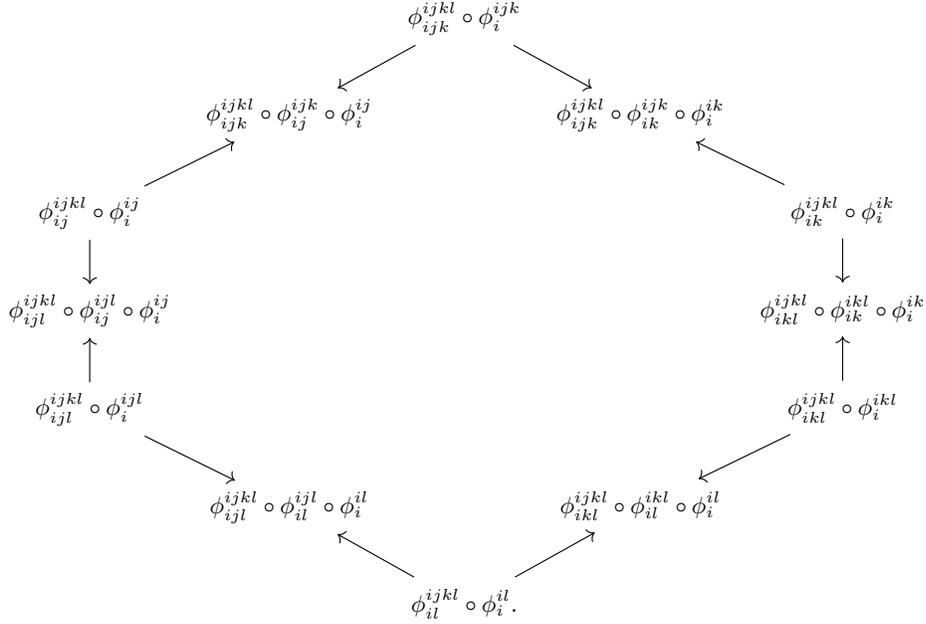
\[
    \begin{tikzcd}[cells={nodes={minimum height=0.7cm}},column sep=tiny]
    && \phi^{ijkl}_{ijk}\circ \phi^{ijk}_i \ar[dl]\ar[dr] \\
    & \phi^{ijkl}_{ijk}\circ \phi^{ijk}_{ij}\circ \phi^{ij}_{i}  &  &  \phi^{ijkl}_{ijk} \circ \phi^{ijk}_{ik}\circ \phi^{ik}_{i} \\
    \phi^{ijkl}_{ij}\circ \phi^{ij}_i \ar[ur]\ar[d]&&&&\phi^{ijkl}_{ik}\circ\phi^{ik}_i\ar[d]\ar[ul]\\
    \phi^{ijkl}_{ijl}\circ \phi^{ijl}_{ij}\circ \phi^{ij}_{i}   & &&& \phi^{ijkl}_{ikl}\circ \phi^{ikl}_{ik}\circ \phi^{ik}_{i}  \\
     \phi^{ijkl}_{ijl}\circ \phi^{ijl}_{i} \ar[u]\ar[dr] &&&& \phi^{ijkl}_{ikl}\circ \phi^{ikl}_i\ar[u]\ar[dl]
    \\
    & \phi^{ijkl}_{ijl}\circ \phi^{ijl}_{il}\circ \phi^{il}_{i} & & \phi^{ijkl}_{ikl}\circ \phi^{ikl}_{il}\circ \phi^{il}_{i} \\
    && \phi^{ijkl}_{il}\circ \phi^{il}_{i}. \ar[ul]\ar[ur]
    \end{tikzcd}
    \]
    Let $\del Q_2\subset\R^2$ denote this 1-dimensional polyhedron, then the diagram above corresponds to a map 
    \[\del Q_2\longrightarrow |\Hom_{\shv_{\mathrm{Ch}_{\R}}(U_{ijkl})}(\F_{ijkl},\F_i|_{U_{ijkl}})|\]
    of CW-complexes. Note that there is a homeomorphism $\del Q_2\cong S^1$, so we see that the previous data $(1)$-$(4)$ yield an element in the fundamental group $\pi_1(|\Hom_{\shv_{\mathrm{Ch}_{\R}}(U_{ijkl})}(\F_{ijkl},\F_i|_{U_{ijkl}})|)$. The cocycle `condition' (we should really talk about cycle \emph{data}) for this quadruple overlap amounts to a choice of a 2-dimensional homotopy trivializing this element. We formulate this combinatorially as follows: we demand the existence of
    \begin{enumerate}[$(i)$]
        \item another quasi-isomorphism
            \[ \phi^{ijkl}_i:\F_{ijkl}\longrightarrow\F_{i}|_{U_{ijkl}}.   \]
            \item homotopies from $\phi^{ijkl}_i$ to all compositions at the vertices of $\del Q_2$.
            \item 2-dimensional homotopies between all possible compositions of the 1-dimensional homotopies.
    \end{enumerate}
    More geometrically, we let $Q_2$ be the polygon obtained by taking the cone on the polyhedron $\del Q_2$, then we demand the existence of a dotted lift 
     \[
    \begin{tikzcd}
    \del Q_2 \ar[d,hook] \ar[r] & {|}\Hom_{\shv_{\mathrm{Ch}_{\R}}(U_{ijkl})}(\F_{ijkl},\F_i|_{U_{ijkl}}){|} \\
    Q_2. \ar[ur,dotted]
    \end{tikzcd}
    \]
    Such a lift corresponds to a diagram 
    
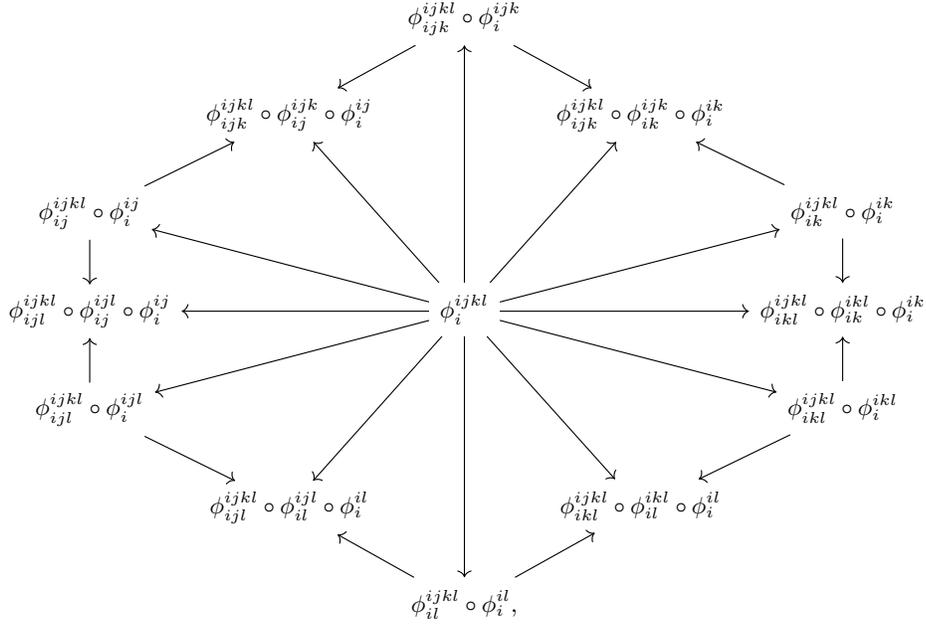
\[
    \begin{tikzcd}[cells={nodes={minimum height=0.7cm}},column sep=tiny]
    && \phi^{ijkl}_{ijk}\circ \phi^{ijk}_i \ar[dl]\ar[dr] \\
    & \phi^{ijkl}_{ijk}\circ \phi^{ijk}_{ij}\circ \phi^{ij}_{i}  &  &  \phi^{ijkl}_{ijk} \circ \phi^{ijk}_{ik}\circ \phi^{ik}_{i} \\
    \phi^{ijkl}_{ij}\circ \phi^{ij}_i \ar[ur]\ar[d]&&&&\phi^{ijkl}_{ik}\circ\phi^{ik}_i\ar[d]\ar[ul]\\
    \phi^{ijkl}_{ijl}\circ \phi^{ijl}_{ij}\circ \phi^{ij}_{i}   & & \phi^{ijkl}_i \ar[uuu]\ar[uul]\ar[uur]\ar[urr]\ar[ull]\ar[ll]\ar[rr]\ar[ddd]\ar[ddl]\ar[ddr]\ar[drr]\ar[dll] && \phi^{ijkl}_{ikl}\circ \phi^{ikl}_{ik}\circ \phi^{ik}_{i}  \\
     \phi^{ijkl}_{ijl}\circ \phi^{ijl}_{i} \ar[u]\ar[dr] &&&& \phi^{ijkl}_{ikl}\circ \phi^{ikl}_i\ar[u]\ar[dl]
    \\
    & \phi^{ijkl}_{ijl}\circ \phi^{ijl}_{il}\circ \phi^{il}_{i} & & \phi^{ijkl}_{ikl}\circ \phi^{ikl}_{il}\circ \phi^{il}_{i} \\
    && \phi^{ijkl}_{il}\circ \phi^{il}_{i}, \ar[ul]\ar[ur]
    \end{tikzcd}
    \]
    in which all triangles are now `filled in' by 2-dimensional homotopies. In fact, for any inclusion $U_{i_1,\ldots,i_n,i_{n+1},i_{n+2},i_{n+3}}\subset U_{i_1,\ldots,i_n}$ adding three more intersections, we have an associated diagram with domain $\del Q_2$, and we demand the existence of a dotted lift as above.
    \item Suppose we have a quintuple intersection $U_i\cap U_j\cap U_k\cap U_l \cap U_m$, then there are 24 quadruple compositions $\F_{ijklm}\rightarrow\F_i$, and the previous data provides 24 squares (arising as various products of homotopies) and eight subdivided hexagons $Q_2$ as above, which fit together to form the boundary of the Archimedean solid known as the \emph{truncated octahedron}. Let $\del Q_3\subset\R^3$ denote this boundary, and let $Q_3$ be the 3-dimensional polyhedron obtained as the cone on $\del Q_3$, then we demand that there exists a lift
    \[
    \begin{tikzcd}
    \del Q_3 \ar[d,hook] \ar[r] & {|}\Hom_{\shv_{\mathrm{Ch}_{\R}}(U_{ijkl})}(\F_{ijklm},\F_i|_{U_{ijklm}}){|} \\
    Q_3. \ar[ur,dotted]
    \end{tikzcd}
    \]
    In fact, for any inclusion $U_{i_1,\ldots,i_n,i_{n+1},i_{n+2},i_{n+3},i_{n+4}}\subset U_{i_1,\ldots,i_n}$ adding four more intersections, we have an associated diagram with domain $\del Q_3$, and we demand the existence of a dotted lift as above.
    \item For every inclusion adding five more intersections, the previous data provides a number of 3-dimensional cubes, cylinders with hexagonal base and truncated octahedra which fit together to form the boundary of a certain 4-dimensional polyhedron $Q_4$, and we demand the existence of a lift of the map from $\del Q_4$ into the appropriate hom space along the boundary inclusion $\del Q_4\subset Q_4$.
    \item[$(8 \rightarrow \infty)$] The previous data fit together to form the boundary of some higher dimensional polyhedron $Q_n$ homeomorphic to the $n$-dimensional disk, and we demand the existence of a lift of the map from $\del Q_n$ into the appropriate hom space along the boundary inclusion $\del Q_n\subset Q_n$.
\end{enumerate}
We claim that if one is given the sheaves $\F_{i_0,\ldots,i_n}$, the quasi-isomorphisms $\phi^{i_1,\ldots,i_n,i_{n+1}}_{i_1,\ldots,i_n}$ and the infinite tower of coherence data sketched above, then one can construct a glued sheaf $\F$ on $X$ that is canonically quasi-isomorphic to $\F_{i}$ on each open $U_i$. To substantiate this claim, we need to do better than vaguely gesture at the higher (i.e. $>3$) coherences: we need to specify the polyhedron $Q_n$ and spell out how its boundary decomposes into products of lower dimensional polyhedra. The most efficient way to achieve this is to repackage a homotopy coherent descent datum into a single functor
\[ \F:Q\longrightarrow  \int_{\mathsf{Top}^{\mathrm{open}}}\shv_{\mathrm{Ch}_{\R}}(\_)^{op} \]
analogous to how a functor $P_I\rightarrow \int_{\mathsf{Top}^{\mathrm{open}}}\shv_{\mathrm{Ch}_{\R}}(\_)^{op}$ encodes a strict descent datum. The functor $\F$
must know about all paths, homotopies between paths, homotopies between homotopies and so on, of all the hom spaces $|\Hom_{\shv_{\mathrm{Ch}_{\R}}}(U_J)(\F_J,\mathcal{G}_J)|$ for $J\subset I$ a nonempty finite subset. As such, $\F$ cannot be an ordinary functor between ordinary categories, but we should treat $\int_{\mathsf{Top}^{\mathrm{open}}}\shv_{\mathrm{Ch}_{\R}}(\_)^{op}$ and $Q$ as \emph{topological categories}, and $\F$ as a \emph{topologically enriched functor}. In fact, we can observe that a homotopical descent datum amounts to the vanishing of certain (inductively constructed) elements in the homotopy groups $\pi_*(|\Hom_{\shv_{\mathrm{Ch}_{\R}}}(U_J)(\F_J,\mathcal{G}_J)|)$, so we only care about the \emph{homotopy type} of the space $|\Hom_{\shv_{\mathrm{Ch}_{\R}}}(U_J)(\F_J,\mathcal{G}_J)|$. To formulate these ideas precisely, we will have to spend some time on a simplicial detour.
\subsubsection{Digression: a bit of abstract homotopy theory}
In abstract homotopy theory, one models homotopy types with the much more combinatorially manageable \emph{simplicial sets} rather than topological spaces. Traditionally, simplicial sets are presented as collections of sets $\{T_n\}_{n\in\Z_{\geq 0}}$ endowed with face maps $\{d_i:T_n\rightarrow T_{n-1}\}$ and degeneracy maps $\{s_j:T_{n-1}\rightarrow T_n\}$ satisfying the simplicial identities. We will instead take the modern perspective and think of simplicial sets as presheaves $X:\simp^{op}\rightarrow \set$ on the simplex category $\simp$ of finite ordinals $\{[n]\}_{n\in\Z_{\geq 0}}$ and order preserving maps among them. Simplicial sets are intimately related to topological spaces and homotopy theory, but also to category theory itself, as we recall in the following examples of simplicial sets.
\begin{ex}
Viewing simplicial sets as presheaves on the simplex category, we see that each finite ordinal $[n]$ determines an object $\Delta^n\in \sset$, the \emph{standard $n$-simplex} via the Yoneda embedding $j:\simp\rightarrow \sset$, so that we have for each simplicial set $T$ a natural bijection $T_n\cong \Hom_{\sset}(\Delta^n,T)$.
\end{ex}
\begin{ex}
Let $\mathsf{CGH}$ be the category of compactly generated Hausdorff spaces. We have a functor $|\_|:\simp\rightarrow \mathsf{CGH}$ carrying a finite ordinal $[n]$ to the \emph{topological $n$-simplex} $|\Delta^n|=\{(x_0,\ldots,x_n)\in \R^{n+1};\, x_i\geq 0\,\forall i,\, \sum_{i=0}^n x_i\leq 1\}$. This functor induces the singular complex functor $\sing:\mathsf{CGH}\rightarrow \sset$, which takes 
a space $Y$ to the simplicial set $\sing(Y)$ whose set of $n$-simplices is the set $\Hom_{\mathsf{CGH}}(|\Delta^n|,Y)$ of continuous maps $|\Delta^n|\rightarrow Y$. The functor $\sing$ admits a left adjoint, the \emph{geometric realization} functor $|\_|$, which carries a simplicial set $T$ to the colimit $\colim_{\Delta^n\rightarrow T}|\Delta^n|$ in the category $\mathsf{CGH}$, where the colimit is taken over all simplices mapping into $T$.
\end{ex}
The category $\mathsf{CGH}$ is the prototypical category with weak equivalences: its weak equivalences are simply the weak homotopy equivalences. A map of simplicial sets $f:T\rightarrow S$ is a \emph{weak equivalence} if the associated map $\sing(X)\rightarrow\sing(Y)$ is a weak homotopy equivalence. Quillen proved that the adjunction $(|\_|\adj \sing)$ induces an equivalence on the associated homotopy categories, so that both categories $\mathsf{CGH}$ and $\sset$ model homotopy types \cite{quillen1967homotopical}. It is well known that the notion of weak homotopy equivalence is not so well behaved with respect to all (compactly generated Hausdorff) spaces. For instance, a weak homotopy equivalence of spaces $f:Y\rightarrow Z$ generally admits a homotopy inverse only if $Y$ and $Z$ are CW complexes. On the simplicial side, the objects for which weak equivalences have inverses are the \emph{Kan complexes}. 
\begin{defn}
For $n\geq 0$ and $0\leq i\leq n$, let $\Lambda^n_i$ be the union of all faces $\Delta^{n-1}\rightarrow \Delta^n$ corresponding to all injective maps $[n-1]\rightarrow [n]$ except the one whose image \emph{does not} contain $i$; this is the $i$'th \emph{horn} of $\Delta^n$. A simplicial set $T$ is a \emph{Kan complex} if for all $n\geq 1$, all $0\leq i\leq n$ and any map $\Lambda^n_i\rightarrow T$, there is a map $\Delta^n\rightarrow T$ fitting into a commuting diagram
\begin{equation*}
\begin{tikzcd}
\Lambda^n_i\ar[d,hook] \ar[r]& T \\
\Delta^n.\ar[ur,dotted]
\end{tikzcd}
\end{equation*}
\end{defn}
For any space $Y$, the singular complex $\sing(Y)$ is a Kan complex. Kan complexes in general are combinatorial gadgets behaving as topological spaces from a homotopical point of view, while general simplicial sets do not admit such an interpretation. For instance, if $T$ is a Kan complex, then the relation of homotopy on vertices (i.e. two vertices $x$ and $y$ are homotopic if there is a simplex $\Delta^1\rightarrow T$ carrying $0$ to $x$ and $1$ to $y$) is an equivalence relation, so it makes sense to speak of the set $\pi_0(T)$ of path components. In fact, there is a completely combinatorial procedure that extracts homotopy groups $\pi_*(T)$ from a Kan complex so that in case $T=\sing(Y)$, there is a canonical isomorphism $\pi_*(T)\cong \pi_*(Y)$. Moreover, every simplicial set may be functorially replaced by a weakly equivalent Kan complex. The unit transformation $T\rightarrow \sing(|T|)$ of the adjunction we just constructed, for instance, provides such a replacement. Motivated by these considerations, we will adhere to the following general principle.
\begin{enumerate}
    \item[$(*)$] A homotopy type \emph{is} a Kan complex. 
\end{enumerate}
\begin{ex}
Let $\cat$ denote the category of categories. We have a fully faithful functor $\simp\hookrightarrow\cat$ simply by viewing a linearly ordered set as a category with at most one morphisms between any two objects. As in the previous example, we have a functor $\cat\rightarrow\sset$, the \emph{nerve}, by carrying a category $\icat$ to the simplicial set $\ner(\icat)$ whose set of $n$-simplices is the set $\Hom_{\cat}([n],\icat)$ of functors from $[n]$ into $\icat$, which we can identify with the set of composable chains of morphisms 
\[ C_0\longrightarrow C_1\longrightarrow \ldots\longrightarrow C_n\]
of length $n$ in $\icat$.
\end{ex}
\begin{defn}
The category $\sset$ of simplicial sets has a symmetric monoidal structure: we define $X\otimes Y$ as $X\times Y$, the product in the category $\sset$, whose set of $n$-simplices is simply the product $X_n\times Y_n$. A \emph{simplicial category} is a category enriched in the symmetric monoidal category of simplicial sets. A \emph{simplicial functor} $f:\icat\rightarrow \icatd$ between simplicial categories is simply a morphism of categories enriched in simplicial sets, that is, a map $f:\mathrm{Ob}(\icat)\rightarrow\mathrm{Ob}(\icatd)$ together with a collection of maps $\{\Hom_{\icat}(C,D)\rightarrow\Hom_{\icatd}(f(C),f(D))\}_{C,D\in \icat}$ of simplicial sets satisfying the usual associativity and unitality conditions. We let $\cat_{\simp}$ denote the category of simplicial categories and simplicial functors among them. \\
Let $\mathsf{Kan}\subset \sset$ be the full subcategory spanned by Kan complexes. This subcategory is stable under limits, so we can view $\mathsf{Kan}\subset \sset$ as a symmetric monoidal subcategory and consider the category of categories enriched in Kan complexes and Kan enriched functors between them as a full subcategory of the category $\cat_{\simp}$. We let $\cat_{\mathsf{Kan}}\subset\cat_{\simp}$ denote the category of \emph{Kan enriched categories}.
\end{defn}
\begin{rmk}
Given a simplicial set $T$, the `Kan replacement functor' $\sing\circ |\_|:\sset\rightarrow \mathsf{Kan}$ preserves products, so we can view this operation as a monoidal functor. It follows that applying $\sing\circ |\_|$ to all simplicial morphism sets of a simplicial category $\icat$ yields a functor $\sing\circ |\_|:\cat_{\simp}\rightarrow \cat_{\mathsf{Kan}}$. 
\end{rmk}
\begin{defn}
Let $\icat$ be a simplicial category. The functor $\pi_0:\mathsf{Kan}\rightarrow \set$ that carries a Kan complex to its set of connected components preserves products, so we can view this operation as a monoidal functor. Applying $\pi_0$ to all simplicial morphism sets of the Kan enriched category $\sing|\icat|$ yields an ordinary category, the \emph{homotopy category} $h\icat$ of $\icat$; this determines a functor $h:\cat_{\simp}\rightarrow\cat$.   
\end{defn}
\begin{ex}
Any ordinary category is a simplicial category (in fact, a Kan enriched category) by viewing the hom sets as constant simplicial sets. This determines a functor $\cat\rightarrow\cat_{\simp}$. We have for any ordinary category $\icat$ viewed as a simplicial category in this manner, an isomorphism $\icat\cong h\icat$ of categories.
\end{ex}
\begin{ex}
The category $\sset$ of simplicial sets is enriched in itself. For $S$ and $T$ simplicial sets, we define a simplicial set $\Hom_{\sset}(S,T)$ as follows: an $n$-simplex of $\Hom_{\sset}(S,T)$ is a map 
\[ S\times\Delta^n\longrightarrow T \]
of simplicial sets, and the simplicial sets $\{\Hom_{\sset}(S,T)\}_{S,T\in \sset}$ make $\sset$ into a simplicial category. If $T$ is a Kan complex, then $\Hom_{\sset}(S,T)$ is again a Kan complex, so the full simplicial subcategory $\mathsf{Kan}\subset \sset$ is Kan enriched. The homotopy category $h\sset$ coincides with the homotopy category $h\mathsf{Kan}$ and is the \emph{classical homotopy category} $\mathcal{H}$ of (compactly generated Hausdorff) spaces.
\end{ex} 
\begin{rmk}
Because the category of simplicial sets has a homotopy theory, the category of simplicial categories `inherits' a homotopy theory: we say that a simplicial functor $f:\icat\rightarrow\icatd$ is a \emph{weak equivalence} if
\begin{enumerate}[$(1)$]
    \item $f$ induces an equivalence of categories $h\icat\rightarrow h\icatd$. 
    \item $f$ is \emph{homotopically fully faithful}: for each pair of objects $C,D\in\icat$, the map $\Hom_{\icat}(C,D)\rightarrow \Hom_{\icatd}(f(C),f(D))$ is a weak homotopy equivalence.
\end{enumerate}
Weak equivalences of simplicial categories are called \emph{Dwyer-Kan equivalences}, after Dwyer and Kan \cite{DK}. Note that we can make the same definition if we replace simplicial sets with compactly generated Hausdorff spaces. This yields the structure of a category with weak equivalences on the category $\cat_{\mathsf{CGH}}$ of categories enriched in compactly generated Hausdorff spaces, and applying the adjunction $(|\_|\adj \sing)$ on all morphism objects induces an adjunction $\begin{tikzcd} \cat_{\simp} \ar[r,shift left] & \cat_{\mathsf{CGH}}\ar[l,shift left] \end{tikzcd}$ which identifies the homotopy categories on both sides. We should think of this homotopy category as the homotopy theory of \emph{categories enriched in homotopy types}, and we have the following extension of our general principle $(*)$.
\begin{enumerate}
    \item[$(**)$] A category enriched in homotopy types is a Kan enriched category. 
\end{enumerate}
Only in a Kan enriched category are the notions of (higher) homotopies of morphisms well behaved.
\end{rmk}
To formulate the notion of a homotopy coherent descent datum, we are now faced with the following tasks.
\begin{enumerate}[$(1)$]
    \item Construct a \emph{Kan enriched} enhancement $\int_{\mathsf{Top}^{\mathrm{open}}}\shv_{\mathrm{Ch}_{\R}}(\_)^{op}_{\simp}$ of the category $\int_{\mathsf{Top}^{\mathrm{open}}}\shv_{\mathrm{Ch}_{\R}}(\_)^{op}$ whose simplicial morphisms sets encode the homological algebra of chain complexes of sheaves, as in Remark \ref{rmk:chainhomotopies}. 
    \item Construct the simplicial category $Q$ which serves a similar function as the poset $P_I$, so that simplicial functors 
    \[ Q\longrightarrow \int_{\mathsf{Top}^{\mathrm{open}}}\shv_{\mathrm{Ch}_{\R}}(\_)^{op}_{\simp}  \]
    are homotopical descent data, rendering precise the sketch we gave in the previous subsection.
\end{enumerate}

We can obtain $Q$ from $P_I$ by suitably `relaxing' the composition rule, following Cordier \cite{cordierhomotopique}.
\begin{cons}[Simplicial thickening of a poset]
Let $P$ be a poset. We define a simplicial category $\mathfrak{C}[P]$ as follows.
\begin{enumerate}
    \item[$(O)$] Objects are the objects of $P$.
    \item[$(M)$] Let $x$ and $y$ be objects of $P$, then we define the simplicial set of morphisms as the nerve of the poset $P_{x,y}$ of linearly ordered subsets $I\subset P$ with least element $x$ and greatest element $y$ (we consider the set $P_{x,y}$ as partially ordered by inclusion).
\end{enumerate}
The composition map $\Hom_{\mathfrak{C}[P]}(x,y)
\times \Hom_{\mathfrak{C}[P]}(y,z)\rightarrow \Hom_{\mathfrak{C}[P]}(x,z)$ is the nerve of the map of posets $P_{x,y}\times P_{y,z}\rightarrow P_{x,z}$ that carries $(I,I')$ to the union $I\cup I'\subset P$, which is again a linearly ordered subset with least element $x$ and greatest element $z$. It is straightforward to check that composition is unital and associative, so $\mathfrak{C}[P]$ is a simplicial category.
\end{cons}
\begin{rmk}
For any two objects $x,y\in P$, the vertices of the simplicial set of morphisms $\Hom_{\mathfrak{C}[P]}(x,y)$ amount to all possible ways to move from $x$ to $y$ in the poset $P$, that is, all nonrepeating chains 
\[x<z_1<z_2<\ldots <z_n <y.\]
If, say, we have a chain $x<v<w<y$, then we can compose the last two morphisms and obtain the chain $x<v<y$, and this composition represents a 1-simplex $\Delta^1\rightarrow \Hom_{\mathfrak{C}[P]}(x,y)$ from $x<v<y$ to $x<v<w<y$. Anytime a chain is obtained from another by composing some number of morphisms in $P$, we add a 1-simplex between these two chains. If we obtain a chain from another by composing in two steps, we have three 1-simplices forming the boundary of a 2-simplex, for instance like so
\[
\begin{tikzcd}
& x<v<y \ar[dr] \\
x<y \ar[rr] \ar[ur]&& x<v<w<y,
\end{tikzcd}
\]
and we add a 2-simplex with this boundary to $\Hom_{\mathfrak{C}[P]}(x,y)$, and so on. 
\end{rmk}
\begin{rmk}
Suppose $P$ is a poset which we view as a simplicial category all of whose morphism simplicial sets consist of a single vertex. There is a canonical simplicial functor of simplicial categories $\mathfrak{C}[P]\rightarrow P$ which collapses every simplicial morphism set to a point.
\end{rmk}
\begin{rmk}
Let $\icat$ be a simplicial category and let $P$ be a poset. We now have three notions of a $P$-diagram in $\icat$:
\begin{enumerate}[$(1)$]
    \item A \emph{(strict) $P$-diagram in $\icat$} is a functor $P\rightarrow \icat$ (which we can think of as either a simplicial functor, viewing $P$ as a simplicial category with no nondegenerate simplices in degrees $\geq 1$, or as an ordinary functor $P\rightarrow \icat_0$).
    \item A \emph{homotopy coherent $P$-diagram in $\icat$} is a simplicial functor $\mathfrak{C}[P]\rightarrow\icat$ from the simplicial thickening of $P$.
    \item A \emph{homotopy $P$-diagram} is a functor $P\rightarrow h\icat$ valued in the homotopy category of $\icat$.
\end{enumerate}
Via the functor $\mathfrak{C}[P]\rightarrow P$, every $P$-diagram in $\icat$ induces a homotopy $P$-diagram for which all higher homotopies are identities. By taking the homotopy category, a homotopy coherent $P$-diagram induces a homotopy $P$-diagram, since the functor $h\mathfrak{C}[P]\rightarrow hP\cong P$ is an equivalence of categories. The latter procedure takes the quotient by the homotopy relation and essentially involves forgetting `how' the diagram $P$ was commutative up to homotopy, which constitutes an undesirable loss of information. In general, there are also obstructions to replacing a homotopy coherent diagram with an equivalent strict one, but from the perspective of the homotopy theory of simplicial categories, the former notion is more natural as will become clear in due course.   
\end{rmk}
\begin{ex}
For the linearly ordered set $[n]$ and any $0\leq i<j\leq n$, the simplicial set of morphisms $\Hom_{\mathfrak{C}[n]}(i,j)$ is the partially ordered set of chains
\[ i=i_0 < i_1 < i_2<\ldots < i_k < j.  \]
We can identify this poset with the power set of $\{i+1,i+2,\ldots,j-1\}\subset [n]$. As a poset, this power set is a product $[1]^{j-i-1}$, so the geometric realization may be identified with an $(j-i-1)$-dimensional cube.
\end{ex}
\begin{ex}\label{ex:simpthickpi}
Let $I$ be a set, like the set indexing a cover $\{U_i\subset X\}_{i\in I}$ of a topological space, and let $J\subset J'\subset I$ be nonempty finite subsets. The simplicial set of morphisms $\Hom_{\mathfrak{C}[P_I]}(J,J')$ may be identified with the nerve of the set of nonrepeating chains of the form 
\[ J=J_0\subset J_1\subset J_2\subset\ldots \subset J_n=J'  \]
partially ordered by inclusion. For $K\subset I$, write $U_K:=\bigcap_{i\in K}U_i$, then such chains correspond to all possible ways of writing the inclusion of intersections $U_{J'}\subset U_J$ as an iterated composition of inclusions of other intersections of the opens $\{U_i\}$. Let $k$ be the cardinality of the set $J\setminus J'$, then the geometric realization of $\Hom_{\mathfrak{C}[P_I]}(J,J')$ is a $(k-1)$-dimensional polyhedron that we can identify with the cone on the double barycentric subdivision of the topological $(k-1)$-simplex, a space we denote by $C(\mathrm{sd}^2(|\del\Delta^{k-1}|))$. This is precisely the object $Q_{k-1}$ that we sketched in the previous subsection and made explicit in the low dimensional example $I=\{ijkl\}$ as $Q_2=|\Hom_{\mathfrak{C}[P_I]}(\{i\},\{ijkl\})|$. The cone point of $C(\mathrm{sd}^2(|\del\Delta^2|))$, corresponding to the smallest chain $\{i\}\subset \{ijkl\}$, was carried to the map $\phi_i^{ijkl}$, and all other vertices were sent to iterated compositions. The boundary $\del Q_{k-1}$ is the subpolyhedron obtained by removing the cone point, so it corresponds to $\mathrm{sd}^2(|\del\Delta^{k-1}|)$, which is the nerve of the poset of nonrepeating chains $\{J\subset \ldots \subset J_i\subset\ldots \subset J'\}$ of length \emph{at least two}. Any such chain is obtained by concatenating subchains $\{J\subset\ldots\subset J_i\}$ and $\{J_i\subset\ldots \subset J'\}$ which precisely corresponds to taking the composition map 
\[  \Hom_{\mathfrak{C}[P_I]}(J,J_i)\times\Hom_{\mathfrak{C}[P_I]}(J_i,J') \longrightarrow \Hom_{\mathfrak{C}[P_I]}(J,J').\]
To sum up, the combinatorics that govern the decomposition of the boundary of $Q_k$ into finite products of $Q_j$s for $j<k$ is fully encoded in the structure of the simplicial category $\mathfrak{C}[P_I]$. 
\end{ex}
With Example \ref{ex:simpthickpi} in hand, we identify the sought after simplicial category $Q$ with the simplicial thickening $\mathfrak{C}[P_I]$ of $P_I$. Our next goal is to formulate a simplicial (Kan enriched) version of the category $\int_{\mathsf{Top}^{\mathrm{open}}}\shv_{\mathrm{Ch}_{\R}}(\_)^{op}$. We turn to the simplicial localization of Dwyer and Kan \cite{DK1}, a much improved version of the localization of categories by zigzags we employed at the start of our discussion. It renders precise the idea that a localization should not involve \emph{identifying} weakly equivalent objects, but instead \emph{adding} data of homotopies between them, that may themselves be homotopic, ad infinitum.  
\begin{cons}[Dwyer-Kan hammock localization]
Let $(\icat,W)$ be a category with weak equivalences. The \emph{hammock localization of $\icat$ with respect to $W$}, denoted $L_H(\icat,W)$, is the simplicial category defined as follows.
\begin{enumerate}
    \item[$(O)$] Objects of $L_H(\icat,W)$ are the objects of $\icat$. 
    \item[$(M)$] Let $C,D\in \icat$ be two objects. The simplicial set of morphisms $\Hom_{L_H(\icat,W)}(C,D)$ is the nerve of the simplicial set of \emph{zigzags} $\mathrm{ZZ}_{\icat}(C,D)$, whose $n$-simplices are diagrams in $\icat$ of the form 
    
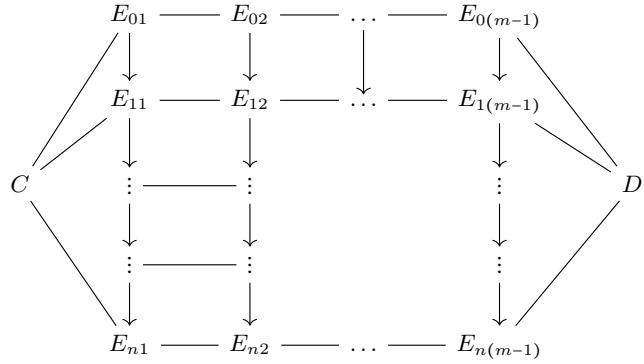
\[
    \begin{tikzcd}
    & E_{01}\ar[d] & E_{02} \ar[d] \ar[l,dash]  \ar[r,dash]& \ldots \ar[r,dash] \ar[d]& E_{0(m-1)} \ar[d] \\
       & E_{11}\ar[d] & E_{12} \ar[d] \ar[l,dash]  \ar[r,dash]& \ldots \ar[r,dash] & E_{1(m-1)} \ar[d] \\
C  \ar[ur,dash] \ar[uur,dash] \ar[ddr,dash]   & \vdots \ar[d]\ar[r,dash]& \vdots\ar[d] &&  \vdots \ar[d] & D \ar[ul,dash] \ar[uul,dash]\ar[ddl,dash]\\
 & \vdots \ar[d] & \vdots \ar[d] \ar[l,dash] & &\vdots\ar[d]  \\
& E_{n1}  & E_{n2} \ar[l,dash]\ar[r,dash]& \ldots \ar[r,dash] & E_{n(m-1)}
    \end{tikzcd}
     \]
     of height $n$ and variable finite length, such that
     \begin{enumerate}[$(1)$]
         \item horizontal arrows are allowed to go in either direction, as long as all arrows in the same column go in the same direction.
         \item arrows in adjacent columns go in different directions.
         \item all leftward maps are in $W$.
         \item all downward maps are in $W$.
         \item no column contains only identities.
     \end{enumerate}
\end{enumerate}
The composition map $\Hom_{L_H(\icat,W)}(C,D)\times\Hom_{L_H(\icat,W)}(E,F)\rightarrow \Hom_{L_H(\icat,W)}(C,F)$ is the map $\mathrm{ZZ}_{\icat}(C,D)\times\mathrm{ZZ}_{\icat}(D,F)\rightarrow \mathrm{ZZ}_{\icat}(C,F)$ that simply concatenates zigzags, and composes columns wherever possible and removes identity columns.\\
We note that viewing an arrow as a 0-simplex of length 1 determines a canonical functor $\icat\rightarrow L_H(\icat,W)$.
\end{cons}
\begin{rmk}\label{rmk:lhfun}
The assignment $(\icat,W)\mapsto L_H(\icat,W)$ is easily seen to be functorial in the following sense: let $\mathsf{RelCat}$ be the category whose objects are pairs $(\icat,W)$ of categories equipped with a class of weak equivalences and whose morphisms $(\icat,W_{\icat})\rightarrow (\icatd,W_{\icatd})$ are functors $\icat\rightarrow \icatd$ that carry $W_{\icat}$ into $W_{\icatd}$. For each such functor, there is an obvious simplicial functor $L_H(\icat,W_{\icat})\rightarrow L_H(\icatd,W_{\icatd})$, so that we have a functor $L_H:\mathsf{RelCat}\rightarrow \cat_{\simp}$.
\end{rmk}
\begin{rmk}\label{rmk:hlh}
The canonical functor $\icat\rightarrow L_H(\icat,W)$ induces an equivalence of categories $h\icat\rightarrow hL_H(\icat,W)$, where the symbol `h' has two different meanings; the first homotopy category is the localization of $\icat$ at $W$, and $hL_H(\icat,W)$ is the homotopy category of the simplicial category $L_H(\icat,W)$. The universal property of localization of categories shows that the functor $\icat\rightarrow h\icat$ can be identified with the composition $\icat\rightarrow L_H(\icat,W)\rightarrow hL_H(\icat,W)$; in particular, if $W$ is the collection of \emph{isomorphisms} in $\icat$, then the functor $\icat\rightarrow hL_H(\icat,W)$ is an equivalence of categories.
\end{rmk}
\begin{defn}
We define a Kan enriched category as follows. Consider the category $\int_{\mathsf{Top}^{\mathrm{open}}}\shv_{\mathrm{Ch}_{\R}}(\_)^{op}$ and let $W$ be the collection of maps $(U,\F)\rightarrow (V,\mathcal{G})$ for which $i:U\hookrightarrow V$ is a homeomorphism and $\F\rightarrow i^*\mathcal{G}$ is a quasi-isomorphism of chain complexes on $U$. This is a class of weak equivalences so we may consider the associated hammock localization. The simplicial category 
\[  L_H(\int_{\mathsf{Top}^{\mathrm{open}}}\shv_{\mathrm{Ch}_{\R}}(\_)^{op},W)\]
is not Kan enriched, so to obtain a category enriched in homotopy types according to our general principle $(*)$ that homotopy types are Kan complexes, we apply the operation $\sing\circ |\_|$ to all simplicial morphism sets. We denote the resulting Kan enriched category by 
\[ \int_{\mathsf{Top}^{\mathrm{open}}}\shv_{\mathrm{Ch}_{\R}}(\_)^{op}_{\simp}. \]
\end{defn}
It follows from Remarks \ref{rmk:lhfun} and \ref{rmk:hlh} that we have a diagram of simplicial categories
\[ \begin{tikzcd}
    \int_{\mathsf{Top}^{\mathrm{open}}}\shv_{\mathrm{Ch}_{\R}}(\_)^{op}\ar[rr] \ar[dr,"p"'] &&  \int_{\mathsf{Top}^{\mathrm{open}}}\shv_{\mathrm{Ch}_{\R}}(\_)^{op}_{\simp}. \ar[dl,"p_{\simp}"] \\
    & \mathsf{Top}^{\mathrm{open}}
\end{tikzcd} \]
which commutes up to natural isomorphism. We have a canonical equivalence of categories
\[h\int_{\mathsf{Top}^{\mathrm{open}}}\shv_{\mathrm{Ch}_{\R}}(\_)^{op}_{\simp} \simeq \int_{\mathsf{Top}^{\mathrm{open}}}\mathbf{D}(\shv_{\mathrm{Vect}_{\R}}(\_)^{op}).  \]
Let us say that a morphism of the simplicial category  $\int_{\mathsf{Top}^{\mathrm{open}}}\shv_{\mathrm{Ch}_{\R}}(\_)^{op}_{\simp}$ is \emph{homotopy $p$-Cartesian} if the associated morphism in its homotopy category is $p$-Cartesian. Finally, we can say what a homotopy coherent descent datum of sheaves of chain complexes for an open cover is.
\begin{defn}
A \emph{homotopy coherent descent datum} for the cover $\{U_i\subset X\}_i$ is a commuting diagram of simplicial categories \[
    \begin{tikzcd}
    & \int_{\mathsf{Top}^{\mathrm{open}}}\shv_{\mathrm{Ch}_{\R}}(\_)^{op}_{\simp}\ar[d,"p_{\simp}"] \\
    \mathfrak{C}[P_I] \ar[r,"f"'] \ar[ur,dotted,"\F"] & \mathsf{Top}^{\mathsf{open}}
    \end{tikzcd}
    \]
    such that $\F$ carries all morphisms to homotopy Cartesian morphisms with respect to $p_{\simp}$.
\end{defn}
\begin{rmk}
The higher homotopies of the simplicial category $\int_{\mathsf{Top}^{\mathrm{open}}}\shv_{\mathrm{Ch}_{\R}}(\_)^{op}_{\simp}$ are constructed via zigzags instead of the chain homotopies as in Remark \ref{rmk:chainhomotopies}, so it is not a priori clear what the Kan complexes of morphisms of $\int_{\mathsf{Top}^{\mathrm{open}}}\shv_{\mathrm{Ch}_{\R}}(\_)^{op}_{\simp}$ have to do with the homological algebra of sheaves of chain complexes. Nevertheless, for $\F$ and $\mathcal{G}$ sheaves of chain complexes on some space $X$, it can be shown that there is a canonical isomorphism in the homotopy category between the CW-complex $|\Hom_{\shv_{\mathrm{Ch}_{\R}}(X)}(\F,\mathcal{G})|$ constructed in Remark \ref{rmk:chainhomotopies} and (the geometric realization of) the Kan complex 
\[\Hom_{\int_{\mathsf{Top}^{\mathrm{open}}}\shv_{\mathrm{Ch}_{\R}}(\_)^{op}_{\simp}}((\mathcal{G},X),(\F,X))\times_{\Hom_{\mathsf{Top}^{\mathrm{open}}}(X,X)}\times\{\mathrm{id}_{X}\}\]
of those maps $(\F,X)\rightarrow(\mathcal{G},X)$ for which the underlying map $X\rightarrow X$ of topological spaces is the identity. 
\end{rmk}
\subsubsection{Attempt 3: homotopy coherent atlases}
With the combinatorial preliminaries taken care of, we can state the analogue of the descent Facts \ref{fact:descent} and \ref{fact:descent2}.
\begin{fact}[Homotopical descent]\label{fact:hodescent}
Fix again a topological space $X$ with an open cover $\{U_i\subset X\}_{i\in I}$ determining a diagram
\[  f:P_I\longrightarrow \mathsf{Top}^{\mathsf{open}}. \]
We define the \emph{homotopy category} of homotopy coherent diagrams, denoted 
\[ h\fun_{\mathrm{Top}^{\mathrm{open}}}^{\mathrm{ho.coh}}(\mathfrak{C}[P_I],\int_{\mathsf{Top}^{\mathrm{open}}}\shv_{\mathrm{Ch}_{\R}}(\_)^{op}_{\simp}) \]
as follows.
\begin{enumerate}
    \item[$(O)$] Objects are commuting diagrams of simplicial categories  
    \[
    \begin{tikzcd}
    & \int_{\mathsf{Top}^{\mathrm{open}}}\shv_{\mathrm{Ch}_{\R}}(\_)^{op}_{\simp}\ar[d,"p_{\simp}"] \\
    \mathfrak{C}[P_I] \ar[r,"f"'] \ar[ur,dotted,"\F"] & \mathsf{Top}^{\mathsf{open}}.
    \end{tikzcd}
    \]
    \item[$(M)$] Let $\F$ and $\mathcal{G}$ be two such objects. Consider the set of homotopy coherent natural transformations
    \[
    \begin{tikzcd}
    & \int_{\mathsf{Top}^{\mathrm{open}}}\shv_{\mathrm{Ch}_{\R}}(\_)^{op}_{\simp}\ar[d,"p_{\simp}"] \\
    \mathfrak{C}[P_I\times \Delta^1] \ar[r] \ar[ur,dotted,"\alpha"] & \mathsf{Top}^{\mathsf{open}}
    \end{tikzcd}
    \]
    that restrict to $\F$ on $\mathfrak{C}[P_I\times \{0\}]$ and restrict to $\mathcal{G}$ on $\mathfrak{C}[P_I\times \{1\}]$ (the lower horizontal map is simply the projection $\mathfrak{C}[P_I\times \Delta^1]\rightarrow \mathfrak{C}[P_I]$ followed by $f$). We define an equivalence relation on this set by declaring $\alpha \sim \beta$ if there exists a commuting diagram
    \[
    \begin{tikzcd}
    & \int_{\mathsf{Top}^{\mathrm{open}}}\shv_{\mathrm{Ch}_{\R}}(\_)^{op}_{\simp}\ar[d,"p_{\simp}"] \\
    \mathfrak{C}[P_I\times \Delta^2] \ar[r] \ar[ur,dotted,"\sigma"] & \mathsf{Top}^{\mathsf{open}}
    \end{tikzcd}
    \]
    such that $\sigma|_{\mathfrak{C}[P_I\times \Delta^{\{0,1\}}]}=\alpha$, $\sigma|_{\mathfrak{C}[P_I\times \Delta^{\{0,2\}}]}=\beta$ and $\sigma|_{\mathfrak{C}[P_I\times \Delta^{\{1,2\}}]}$ is the identity natural transformation on $\mathcal{G}$. One can prove that this prescription does define an equivalence relation, so we define the set of morphisms to be the set of equivalence classes of such homotopy coherent natural transformations.
\end{enumerate}
One can show that composition of homotopy coherent natural transformations is well defined on equivalence classes, so we indeed have a category. Then there is a canonical equivalence between
\begin{enumerate}[$(a)$]
    \item The full subcategory of the above homotopy category spanned by lifts $\F$ that carry every morphism $J\subset J'$ to a homotopy Cartesian morphism.
    \item The opposite derived category $\mathbf{D}(\shv_{\mathcal{A}}(X))^{op}$.
\end{enumerate}
\end{fact}
\begin{rmk}
Note that we do \emph{not} first define a category or simplicial category
\[\fun_{\mathrm{Top}^{\mathrm{open}}}^{\mathrm{ho.coh}}(\mathfrak{C}[P_I],\int_{\mathsf{Top}^{\mathrm{open}}}\shv_{\mathrm{Ch}_{\R}}(\_)^{op}_{\simp})\]
of homotopy coherent diagrams! It is actually quite hard to write down in explicit terms the correct simplicial category of such with morphisms being homotopy coherent natural transformations. The difficulty is that composition of coherent natural transformations is only naturally associative up to coherent homotopy (in the construction above, we quotient out this ambiguity so that we are left with a category). To obtain a genuine simplicial category with a composition law that is associative on the nose, one must apply a `rigidification' procedure, which is fairly inexplicit. All this should be taken as an indication that our current setup for talking about homotopy coherence is suboptimal. Below, we switch from simplicial categories to \emph{\infcatst} which are designed to handle the ambiguity issues we are running into.
\end{rmk}
Facts \ref{fact:descent2} and \ref{fact:hodescent} show that the homotopy categories of strict and homotopical descent data are both equivalent to the derived category of complexes of sheaves on $X$. Here is another (stronger) way to phrase this: since strict diagrams and ordinary natural transformations of such may be viewed as coherent diagrams and coherent natural transformations via the functors $\mathfrak{C}[P_I]\rightarrow P_I$ and $\mathfrak{C}[P_I\times \Delta^1]\rightarrow P_I\times \Delta^1$, there is an obvious functor
\[\fun_{\mathrm{Top}^{\mathrm{open}}}(P_I,\int_{\mathsf{Top}^{\mathrm{open}}}\shv_{\mathrm{Ch}_{\R}}(\_)^{op})\longrightarrow  h\fun_{\mathrm{Top}^{\mathrm{open}}}^{\mathrm{ho.coh}}(\mathfrak{C}[P_I],\int_{\mathsf{Top}^{\mathrm{open}}}\shv_{\mathrm{Ch}_{\R}}(\_)^{op}_{\simp}). \]
With a bit of work, one can show that this functor carries weak equivalences to isomorphisms and thus induces a functor
\[h\fun_{\mathrm{Top}^{\mathrm{open}}}(P_I,\int_{\mathsf{Top}^{\mathrm{open}}}\shv_{\mathrm{Ch}_{\R}}(\_)^{op})\longrightarrow  h\fun_{\mathrm{Top}^{\mathrm{open}}}^{\mathrm{ho.coh}}(\mathfrak{C}[P_I],\int_{\mathsf{Top}^{\mathrm{open}}}\shv_{\mathrm{Ch}_{\R}}(\_)^{op}_{\simp}),\]
which is an equivalence of categories. We stress that the set of objects of the latter category is \emph{much larger} than the set of object of the former, giving us much more freedom to work with homotopy coherent diagrams, at least in principle. On the other hand, it is much harder to specify a homotopical descent datum explicitly than it is to specify a strict one, as one needs to avail oneself of an infinite tower of higher dimensional homotopies that fit together ever more intricately as the number of intersecting opens increases. Upon applying the tangent complex functor, our equivalent Definitions \ref{defn:kuratlas} and \ref{defn:kuratlas2} only provide the first two layers of this tower, namely transition isomorphisms $\phi_{ij}$ on overlaps $U_i\cap U_j$ and homotopies on triple overlaps, which is not enough to construct a glued complex of sheaves. Definition \ref{defn:kuratlas4} provides the entire tower of homotopies, but in a rather trivial manner: they are all identities. If we wish to define an intermediate notion of a Kuranishi atlas that takes advantage of the extra flexibility that Fact \ref{fact:hodescent} affords (so that applying the functor $\tanc$ yields a nontrivial homotopical descent datum) we need to be able to talk about homotopies of maps of Kuranishi atlases, and homotopies of those homotopies, and so on; that is, we need to enhance the morphism sets of $\mathsf{AffKur}$ to Kan complexes. Fortunately, we already know how to do this: we simply replace the homotopy category $h\mathsf{AffKur}$ by the improved hammock localization.
\begin{defn}
Consider the the Kan enriched category obtained by applying the operation $\sing\circ |\_|$ to the hammock localization of the category $\mathbf{dgMan}$ at the weak equivalences of Definition \ref{def:kurweakeq}. The Kan enriched category $\mathsf{AffKur}_{\simp}$ is the full subcategory of this simplicial category spanned by affine Kuranishi models.
\end{defn}
\begin{rmk}\label{rmk:affkurpullback}
The reader might begin to worry that the proliferation of algebro-topological machinery is producing categorical structures that are increasingly resistant to geometric intuition. As reassurance that $\mathsf{AffKur}_{\simp}$ can be justifiable interpreted as a simplicial category of zero sets of sections of vector bundles, we offer the following results: the functor $\mathsf{Man}\rightarrow \mathsf{AffKur}_{\simp}$ is homotopically fully faithful; that is, the map $\cinfty(M,N)\rightarrow \Hom_{\mathsf{AffKur}_{\simp}}(M,N)$ is a homotopy equivalence of Kan complexes. Moreover, this functor carries binary products of manifolds to products of affine Kuranishi models. What's more, every object $\mathbf{K}=(X,E,s)$ is canonically the \emph{homotopy limit} in $\mathsf{AffKur}_{\simp}$ of the diagram 
\[ \begin{tikzcd}
    & X\ar[d,"s"] \\
    X\ar[r,"0"] & E
\end{tikzcd} \]
of manifolds, where the lower horizontal map is the zero section of $E$.
\end{rmk}
Because the hammock localization and the operation $\sing\circ |\_|$ are functorial, we automatically have a tangent complex \emph{simplicial} functor
\begin{equation}\label{eq:derivedtangentcplx}
\begin{tikzcd}
\mathsf{AffKur}_{\simp}^{\mathrm{open}} \ar[rr,"\tanc"] \ar[dr,"Z"']&&\int_{\mathsf{Top}^{\mathrm{open}}}\shv_{\mathrm{Ch}_{\R}}(\_)^{op}_{\simp} \ar[dl,"p_{\simp}"] \\ & \mathsf{Top}^{\mathrm{open}}.
\end{tikzcd}
\end{equation}
Yet another notion of a Kuranishi atlas, stronger than the equivalent Incorrect Definitions \ref{defn:kuratlas} and \ref{defn:kuratlas2}, but weaker than Incorrect Definition \ref{defn:kuratlas4}, now stares us in the face.
\begin{defn}\label{defn:kuratlas3}
Let $X\in \mathsf{Top}^{\mathsf{open}}$, then a \emph{Kuranishi atlas} on $X$ consists of the following data.
\begin{enumerate}[$(a)$]
    \item A collection of maps $\{V_i\rightarrow X\}_{i\in I}$ in $\mathsf{Top}^{\mathrm{open}}$ to $X$ with images being open sets $\{U_i\subset X\}$ that cover $X$, which induces a functor
    \[  f:\mathfrak{C}[P_I]\longrightarrow \mathsf{Top}^{\mathsf{open}} \]
    \item A dotted lift $\tilde{f}$ of $f$ as follows
    \[
    \begin{tikzcd}
    & \mathsf{AffKur}_{\simp}^{\mathrm{open}}\ar[d,"Z"] \\
    \mathfrak{C}[P_I] \ar[r,"f"'] \ar[ur,dotted,"\tilde{f}"] & \mathsf{Top}^{\mathsf{open}}
    \end{tikzcd}
    \]
    that makes the diagram of simplicial categories commute. Moreover, we require that $\tilde{f}$ carries every morphism in $\mathfrak{C}[P_I]$ to a homotopy Cartesian morphism with respect to $Z$.
\end{enumerate}
\end{defn}
We have some reason to believe that this intermediate notion of a `space locally equivalent to a zero set of a smooth function' is the correct one: combining the Diagram \eqref{eq:derivedtangentcplx} with Fact \ref{fact:hodescent}, we deduce that for each space $X$ equipped with a Kuranishi atlas, there exists a global object $\tanc_X$ such that for each open $U_i\subset X$ of the cover, there is a canonical isomorphism $\tanc_{\mathbf{K}_i}\cong \tanc_X|_{U_i}$ in the derived category. Moreover, the simplicial morphism sets of $\mathsf{AffKur}_{\simp}$ are a more sophisticated version of the zigzags of Construction \ref{cons:zigzag} so Problem \ref{prob:nokur} and Remark \ref{rmk:nokur} do not apply (at least not obviously). There are however a number of questions that immediately come to mind.
\begin{enumerate}[$(Q1)$]
    \item Our `correct' definition of a Kuranishi atlas involves a very large amount of data that seems impossible to manipulate: unraveling the definition of an atlas for a cover $\{U_i\subset X\}_{i\in I}$, we have for each nonempty finite $J\subset I$ an affine Kuranishi model $\mathbf{K}_J$ with $Z(\mathbf{K}_J)\cong U_J=\cap_{j\in J}U_j$ and for each inclusion $J\subset J'\subset I$ with $|J'\setminus J|=n$, a continuous map
    \[ C(\mathrm{sd}^2|\del\Delta^{n-1}|) \longrightarrow  |\mathrm{ZZ}(\mathbf{K}_{J},\mathbf{K}_{J'})|  \]
    from the cone on the double barycentric subdivision of the boundary of the topological $(n-1)$-simplex (considered as a subspace of $\R^{n-1}$, say), to the geometric realization of simplicial set of zigzags between $\mathbf{K}_{J}$ and $\mathbf{K}_{J'}$ that is required to take values in the union of those connected components containing localizations. Moreover, the subpolyhedron $\mathrm{sd}^2|\del\Delta^{n-1}|$ is carried to a union of products of previously specified polyhedra in a manner governed by the structure of the simplicial category $\mathfrak{C}[P_I]$. How are we supposed to construct all this data in concrete examples? More to the point: how do we show that for the elliptic PDE $P$, the topological space $\mathsf{Sol}(P)$ admits a Kuranishi atlas in the sense of Definition \ref{defn:kuratlas3}? 
    \item A $\cinfty$-manifold generally admits quite a large set of atlases inducing the same $\cinfty$-structure. For a topological manifold $M$, there is an equivalence relation on the set of all $\cinfty$-atlases by \emph{compatibility}, so that the set of equivalence classes is the set of $\cinfty$-structures on $M$. It seems that a similar story should hold for our Kuranishi atlases; in particular, we would like to define some notion of compatibility which would allow us to talk about \emph{Kuranishi structures} on a topological space $X$ (without having to invoke a specific cover), but it is not so obvious how to proceed starting from Definition \ref{defn:kuratlas3}.
    \item Related to the previous question, one might guess that there is some natural notion of morphism between spaces $X$ and $Y$ equipped with Kuranishi atlases that reduces to the usual notion of a $\cinfty$-map when the Kuranishi atlases give $X$ and $Y$ the structure of $\cinfty$-manifolds, and reduces to a morphism in the simplicial category $\mathsf{AffKur}_{\simp}$ if $X$ and $Y$ are covered by a single chart. In other words, we might ask for a simplicial (Kan enriched) category of \emph{Kuranishi spaces}.
\end{enumerate}
For the reader with practical applications to PDEs in mind, $Q1$ may seem the most pressing concern. It is possible to address this question `by hand', by generalizing the sum chart method of Construction \ref{cons:sumchart} (which is essentially how all constructions of Kuranishi atlases on moduli spaces in the literature appear). Using the notation from Construction \ref{cons:sumchart}, we can choose a collection of points $\{x_i\}_i\subset\mathsf{Sol}(P)$ and an open cover $\{U_i\subset \mathsf{Sol}(P)\}_{i\in I}$ with $x_i\subset U_i$ so that we have affine Kuranishi models $\mathbf{K}_i=(Z(P+t_{x_i}),Z(P+t_{x_i})\times\coker\,TP_{x_i},f_{x_i}+\mathrm{id})$ with $Z(\mathbf{K}_i)\cong U_i$. For each nonempty finite subset $J\subset I$, we may build a `higher' sum chart by considering the zero locus of
\[P+\sum_{j\in J}t_j: \Gamma(V)\times \prod_{j\in J} \coker\,TP_{x_j}\longrightarrow \Gamma(F)  \] 
and setting $\mathbf{K}_{J}:=(Z(t_J),Z(t_J)\times \prod_{j\in J} \coker\,TP_{x_j},f_J\times\mathrm{id})$. For any inclusion $J\subset J'$, we have an obvious weak equivalence $\phi_{JJ'}:\mathbf{K}_{J}|_{Z(t_{J'})}\rightarrow \mathbf{K}_{J'}$, so that for a triple $J\subset J'\subset J''$ the cocycle condition holds, restricted to $Z(t_{J''})$ (recall from Remark \ref{rmk:nokur} that the maps $\phi_{JJ'}$ do \emph{not} constitute a strict Kuranishi atlas; they go in the wrong direction). We can axiomatize the output of this construction.
\begin{defn}\label{defn:kuratlasspecial}
Let $X$ be a paracompact Hausdorff topological space. A \emph{special Kuranishi atlas on $X$} consists of the following data.
\begin{enumerate}[$(a)$]
\item An open cover $\{U_i\subset X\}_i$.
\item For each nonempty finite subset $J\subset I$, an affine Kuranishi model $\mathbf{K}_J$.
\item A collection of homeomorphisms $\phi_J:Z(\mathbf{K}_J)\rightarrow U_J$ with $U_J=\cap_{j\in J}U_j$.
\item For each inclusion $J\subset J'$, there is an open set $V_{J'}$ of the domain of $\mathbf{K}_J$ and a weak equivalence $\phi_{JJ'}:\mathbf{K}_J|_{V_{J'}}\rightarrow \mathbf{K}_{J'}$.
\end{enumerate}
These data are required to satisfy the following conditions.
\begin{enumerate}[$(1)$]
\item The transition maps are compatible with the footprint maps: $\psi_{J'}\circ Z(\phi_{JJ'})=\psi_J$.
\item The cocycle condition holds: for every composition $J\subset J\subset J'$, we have \[\phi_{J'J''}\circ \phi_{JJ'}=\phi_{JJ''}\]
on the open set $\phi_{J'J''}^{-1} (V_{J''})$.
\end{enumerate}
\end{defn}
The definition above is essentially the one that appears in recent literature, compare for instance with McDuff-Wehrheim and Pardon \cite{mcduffwehrheim,Pardon}.
\begin{prop}
Let $\{U_i\subset X\}_{i\in I}$ be a paracompact Hausdorff space equipped with a special Kuranishi atlas $\{\mathbf{K}_J\}_{J\subset I\text{ finite}}$. Then $X$ admits Kuranishi atlas in the sense of Definition \ref{defn:kuratlas3} for the cover $\{U_i\subset X\}_i$ that assigns $\mathbf{K}_J$ to the open set $U_J$. 
\end{prop}
One proves this by inductively constructing zigzags of affine Kuranishi models of increasing height; we leave the details to the motivated reader.\\
The preceding construction of Kuranishi atlases is unsatisfactory in several respects: it is ad-hoc, depends on many choices and is not obviously functorial under chances of input data. To address these concerns, it turns out to be advantageous to deal with $Q2$ and $Q3$ first. To answer \emph{these} questions satisfactorily, it is crucial that we work with Definition \ref{defn:kuratlas3} instead of \ref{defn:kuratlasspecial}; special Kuranishi atlases are still too restrictive. We will see that once we have a good Kan enriched category of Kuranishi spaces, a better solution to $Q1$ will present itself. 
\subsubsection{Gluing Kuranishi atlases}
Recall Fact \ref{fact:hodescent}: giving a homotopical descent datum $\mathfrak{C}[P_I]\rightarrow \int_{\mathsf{Top}^{\mathrm{open}}}\shv_{\mathrm{Ch}_{\R}}(\_)^{op}_{\simp}$ is equivalent to giving a complex of sheaves on $X$, and one passes from the descent datum to the sheaf by taking a homotopy colimit. Similarly, we should think of a Kuranishi atlas $\tilde{f}:\mathfrak{C}[P_I]\rightarrow \mathsf{AffKur}_{\simp}$ as a descent datum for some sort of global object, a \emph{Kuranishi structure} on $X$, which, in analogy with complexes of sheaves, may be defined using \emph{different choices} of covers and associated descent data. The analogy with complexes of sheaves \emph{fails} in that it is in general not possible to take the homotopy colimit of the diagram $\tilde{f}$ in $\mathsf{AffKur}_{\simp}$. To deal with this problem, we apply a well known categorical construction: we formally add all homotopy colimits to $\mathsf{AffKur}_{\simp}$. Recall that for an ordinary category $\icat$, the category $\pshv_{\set}(\icat)=\fun(\icat^{op},\set)$ of \emph{presheaves} on $\icat$ admits all limits and colimits. In fact, the Yoneda embedding $j:\icat\hookrightarrow \pshv_{\set}(\icat)$ exhibits the category of presheaves on $\icat$ as the category \emph{freely generated under colimits} by $\icat$. We expect that a similar result for homotopy colimits holds upon replacing sets with homotopy types.
\begin{defn}
Consider the category $\pshv_{\sset}(\mathsf{AffKur}_{\simp})=\fun(\mathsf{AffKur}_{\simp}^{op},\sset)$ of simplicial presheaves on $\mathsf{AffKur}_{\simp}$. This category is a simplicial category in a natural way: for $K$ a simplicial set, let $\underline{K}:\mathsf{AffKur}_{\simp}^{op}\rightarrow \sset$ be the constant presheaf with value $K$, then an $n$-simplex in the simplicial morphism set $\Hom_{\pshv_{\sset}(\mathsf{AffKur}_{\simp})}(F,G)$ is by definition a natural transformation
\[ F\times\underline{\Delta^n}\longrightarrow  G \]
of simplicial presheaves.
\end{defn}
 Enriched category theory provides a fully faithful simplicial Yoneda embedding
\[ j:\mathsf{AffKur}_{\simp}\hooklongrightarrow \pshv_{\sset}(\mathsf{AffKur}_{\simp}). \]
Let $\pshv_{\mathsf{Kan}}(\mathsf{AffKur}_{\simp})\subset \pshv_{\sset}(\mathsf{AffKur}_{\simp})$ be the full subcategory spanned by presheaves taking values in Kan complexes, then the fact that $\mathsf{AffKur}_{\simp}$ is Kan-enriched is equivalent to the assertion that the Yoneda embedding factors through $\pshv_{\mathsf{Kan}}(\mathsf{AffKur}_{\simp})$. Beware: the category $\pshv_{\mathsf{Kan}}(\mathsf{AffKur}_{\simp})$ is not itself Kan enriched! To obtain Kan enriched category, we must restrict to a certain full subcategory. 
\begin{defn}
Let $\pshv_{\mathsf{Kan}}(\mathsf{AffKur}_{\simp})$ be the category of simplicial functors $\mathsf{AffKur}_{\simp}^{op}\rightarrow \mathsf{Kan}$. A functor $X:\mathsf{AffKur}_{\simp}^{op}\rightarrow\mathsf{Kan}$ is \emph{projectively cofibrant} if for each natural transformation $F\rightarrow G$ such that $F(\mathbf{K})\rightarrow G(\mathbf{K})$ is a trivial Kan fibration for all $\mathbf{K}\in \mathsf{AffKur}_{\simp}$, every map $X\rightarrow G$, i.e. every natural transformation, admits a dotted lift as in the diagram 
\[
\begin{tikzcd}
    & F\ar[d] \\
X\ar[ur,dotted]\ar[r] & G.
\end{tikzcd}
\]
We let $\pshv_{\mathsf{Kan}}(\mathsf{AffKur}_{\simp})^c$ denote the simplicial category of projectively cofibrant Kan valued presheaves.
\end{defn}
\begin{rmk}
The fact that the simplicial category $\pshv_{\mathsf{Kan}}(\mathsf{AffKur}_{\simp})^c$ is Kan enriched is a basic result of \emph{enriched model category theory} and is exactly analogous to the following observation: we can make the category $(\sset,W)$ with the weak homotopy equivalences into a Kan enriched category by taking the hammock localization $L_H(\sset,W)$ and applying the operation $\sing\circ | \_|$ but this is unnecessarily complicated; instead of enlarging $\sset$, we can \emph{restrict} to the full subcategory of Kan complexes, which is automatically Kan enriched and represents the same object in the homotopy theory of simplicial categories, as the canonical simplicial functor $\mathsf{Kan}\rightarrow |\sing(L_H(\sset,W))|$ is a Dwyer-Kan equivalence.
\end{rmk}
For every $\mathbf{K}\in \mathsf{AffKur}_{\simp}$, the Yoneda functor $j(\mathbf{K})$ is a projectively cofibrant Kan valued presheaf. Indeed, by Yoneda's lemma, this is simply the statement that for any objectwise trivial Kan fibration $F\rightarrow G$, the map $F(\mathbf{K})\rightarrow G(\mathbf{K})$ is a surjection in degree 0. We thus have a fully faithful functor
\[ j:\mathsf{AffKur}_{\simp}\hooklongrightarrow \pshv_{\mathsf{Kan}}(\mathsf{AffKur}_{\simp})^c, \]
\emph{the Yoneda embedding}, and it can be shown that $\pshv_{\mathsf{Kan}}(\mathsf{AffKur}_{\simp})^c$ admits all homotopy colimits. It is tempting to now define a \emph{Kuranishi space} as a simplicial presheaf $F$ such that $F$ is the homotopy colimit of a Kuranishi atlas 
\[\mathfrak{C}[P_I]\longrightarrow \mathsf{AffKur}_{\simp}\hooklongrightarrow \pshv_{\mathsf{Kan}}(\mathsf{AffKur}_{\simp})^c\]
for some topological space $X$ with a cover $\{U_i\subset X\}_i$, but this does not work.
\begin{prob}
Recall that we want to be able to express the idea that a single Kuranishi \emph{structure} on $X$ may be presented by many different Kuranishi \emph{atlases}, but in the simplicial category $\pshv_{\mathsf{Kan}}(\mathsf{AffKur}_{\simp})^c$, each atlas determines a \emph{distinct} object because this simplicial category is freely generated under homotopy colimits. Indeed, let $\mathbf{K}$ be an affine Kuranishi model, let $\{U_i\subset Z(\mathbf{K})\}_{i\in I}$ be an open cover and let 
\[\tilde{f}:\mathfrak{C}[P_I]\longrightarrow\mathsf{AffKur}_{\simp}^{\mathrm{open}}\]
be a diagram together with a compatible collection of maps $\{\tilde{f}(J)\rightarrow\mathbf{K}\}_{J\subset I,|J|<\infty}$ that are all localizations, then the homotopy colimit of the diagram $\tilde{f}$ \emph{does} exist in $\mathsf{AffKur}_{\simp}$ and is the object $\mathbf{K}$ itself. However, the Yoneda embedding does not preserve this homotopy colimit: the canonical map 
\[ \hocolim_{\mathfrak{C}[P_I]}j\tilde{f}\longrightarrow j(\hocolim_{\mathfrak{C}[P_I]}\tilde{f}) \simeq  j(\mathbf{K})\]
in the Kan enriched category $\pshv_{\mathsf{Kan}}(\mathsf{AffKur}_{\simp})^{c}$ is not an equivalence.
\end{prob}
 To resolve this issue, we restrict to those presheaves that view the map above as an equivalence.
\begin{defn}
A simplicial presheaf $F$ is a \emph{sheaf} if for all data as above, the canonical map
\[ F(\mathbf{K}) \longrightarrow \Hom_{\pshv_{\mathsf{Kan}}(\mathsf{AffKur}_{\simp})}(\hocolim_{\mathfrak{C}[P_I]}j\tilde{f},F)  \]
is an equivalence. We denote by $\shv_{\mathsf{Kan}}(\mathsf{AffKur}_{\simp})^c$ the Kan enriched full subcategory spanned by the projectively cofibrant Kan valued presheaves on $\mathsf{AffKur}_{\simp}$ that are also sheaves (in the sequel, we will for the sake of brevity just call these objects \emph{simplicial sheaves}).
\end{defn}
\begin{rmk}
The Yoneda embedding factors through $\shv_{\mathsf{Kan}}(\mathsf{AffKur}_{\simp})^c$: every representable presheaf $j(\mathbf{K})$ is a sheaf. The simplicial category $\shv_{\mathsf{Kan}}(\mathsf{AffKur}_{\simp})^c$ again admits all homotopy colimits. In fact, the inclusion $\shv_{\mathsf{Kan}}(\mathsf{AffKur}_{\simp})^c\subset \pshv_{\mathsf{Kan}}(\mathsf{AffKur}_{\simp})^c$ of sheaves into presheaves admits a left adjoint
\[ \begin{tikzcd}  h\pshv_{\mathsf{Kan}}(\mathsf{AffKur}_{\simp})^c \ar[r,"L",shift left] & h\shv_{\mathsf{Kan}}(\mathsf{AffKur}_{\simp})^c\ar[l,hook',shift left]\end{tikzcd}  \]
on the level of homotopy categories, the \emph{sheafification functor}. It follows formally that one computes a homotopy colimit by taking the homotopy colimit of presheaves and applying $L$ to the result.
\end{rmk}
We can now make the following definition.
\begin{defn}\label{defn:kurspace}
The Kan enriched category of \emph{Kuranishi spaces}, denoted $\mathsf{Kur}_{\simp}$, is the full subcategory of $\shv_{\mathsf{Kan}}(\mathsf{AffKur}_{\simp})^c$ spanned by objects that are equivalent to homotopy colimits of Kuranishi atlases. More precisely, an object $F\in \shv_{\mathsf{Kan}}(\mathsf{AffKur}_{\simp})^c$ is a Kuranishi space if we can find a paracompact Hausdorff space $X$, a cover $\{U_i\}$ determining a functor $P_I\rightarrow \mathsf{Top}^{\mathsf{open}}$ and a dotted lift $\tilde{f}$ of $f$
\[
    \begin{tikzcd}
    & \mathsf{AffKur}_{\simp}^{\mathrm{open}}\ar[d,"Z"] \\
    \mathfrak{C}[P_I] \ar[r,"f"'] \ar[ur,dotted,"\tilde{f}"] & \mathsf{Top}^{\mathsf{open}}
    \end{tikzcd}
\]
such that $F$ is isomorphic in the homotopy category to the homotopy colimit of the composition
\[  \mathfrak{C}[P_I] \longrightarrow \mathsf{AffKur}_{\simp}^{\mathrm{open}}\subset \mathsf{AffKur}_{\simp}\overset{j}{\hooklongrightarrow}\shv_{\mathsf{Kan}}(\mathsf{AffKur}_{\simp})^c.\]
\end{defn}
We finally have a Kan enriched category of Kuranishi spaces to each object of which we can assign a tangent complex (although our construction is admittedly quite formal). In fact, one can show that this assignment is functorial -again up to coherent homotopy- so that different atlases presenting the same object in $\mathsf{Kur}_{\simp}$ determine quasi-isomorphic tangent complexes (this is one instance where it is crucially important that we glue Kuranishi atlases in the simplicial category of sheaves, not presheaves).
\begin{rmk}
Let us make another argument for the reasonableness of Definition \ref{defn:kurspace}. Let $\mathsf{CartSp}\subset\mathsf{AffKur}_{\simp}$ be the full subcategory spanned by affine Kuranishi models of the form $(\R^n,\R^n,0)$ that is, Cartesian spaces equipped with the zero vector bundle. Then the full subcategory of $\mathsf{Kur}_{\simp}$ spanned by Kuranishi spaces that admit an atlas $\mathfrak{C}[P_I]\rightarrow \mathsf{AffKur}_{\simp}^{\mathrm{open}}$ that factors through $\mathsf{CartSp}^{\mathrm{open}}$ is Dwyer-Kan equivalent to the ordinary category of $\cinfty$ manifolds.    
\end{rmk}
\begin{rmk}
Just as every affine Kuranishi model has an underlying topological space, every sheaf on affine Kuranishi models has an underlying sheaf on topological spaces. Composing with the functor $Z:\mathsf{AffKur}_{\simp}\rightarrow \mathsf{Top}$ induces a simplicial functor $Z^*:\shv_{\mathsf{Kan}}(\mathsf{Top})^c\rightarrow \shv_{\mathsf{Kan}}(\mathsf{AffKur}_{\simp})^c$ which admits a left adjoint
\[ \begin{tikzcd}  h\shv_{\mathsf{Kan}}(\mathsf{AffKur}_{\simp})^c \ar[r,"Z_!",shift left] & h\shv_{\mathsf{Kan}}(\mathsf{Top})^c\ar[l,"Z^*",shift left]\end{tikzcd}  \]
on the level of homotopy categories such that the diagram 
\[
\begin{tikzcd}
h\mathsf{AffKur}_{\simp} \ar[d,"j"]\ar[r,"Z"] & \mathsf{Top} \ar[d,"j"] \\
h\shv_{\mathsf{Kan}}(\mathsf{AffKur}_{\simp})^c\ar[r,"Z_!"] & h\shv_{\mathsf{Kan}}(\mathsf{Top})^c
\end{tikzcd}
\]
commutes up to natural isomorphism. It follows formally that for a Kuranishi space $F$ determined by a cover $\{U_i\subset X\}_i$ and a descent datum $\mathfrak{C}[P_I]\rightarrow \mathsf{AffKur}_{\simp}$, the underlying sheaf of spaces $Z_!(F)$ is representable and equivalent to $j(X)$.
\end{rmk}
Our careful foundational efforts are paying off: by design, $\mathsf{Kur}_{\simp}$ lies fully faithfully in a much larger simplicial category, which affords many benefits.
For instance, we are now in the position to present the following powerful representability criterion.
\begin{thm}\label{thm:repcritkur}
Let $F$ be simplicial sheaf on $\mathsf{AffKur}_{\simp}$. Suppose that there exists a collection of affine Kuranishi models $\{\mathbf{K}_i\}_{i\in I}$ and a collection of maps $\{j(\mathbf{K}_i)\rightarrow F\}_{i\in I}$ of simplicial sheaves such that for each affine Kuranishi model $\mathbf{J}$ and each map $j(\mathbf{J})\rightarrow F$ of simplicial sheaves, the following conditions are satisfied. 
\begin{enumerate}[$(1)$]
    \item For each $\mathbf{K}_i$, the homotopy pullback 
    \[\begin{tikzcd}
    j(\mathbf{J})\times^h_F j(\mathbf{K}_i) \ar[d] \ar[r] & j(\mathbf{J}) \ar[d] \\
    j(\mathbf{K}_i) \ar[r] & F
    \end{tikzcd}\] 
    is representable by some affine Kuranishi model $\mathbf{H}_i$.
    \item The map $j(\mathbf{J})\times^h_F j(\mathbf{K}_i)\rightarrow j(\mathbf{J})$, which, by $(1)$ and the homotopical fully faithfulness of the Yoneda embedding corresponds to a map $\mathbf{H}_i\rightarrow \mathbf{J}$ in $h\mathsf{AffKur}$, is a localization with respect to an open subset $U_i\subset Z(\mathbf{J})$.
    \item The collection of open subsets $\{U_i\subset Z(\mathbf{J})\}_{i\in I}$ obtained from the localizations $\{\mathbf{H}_i\rightarrow \mathbf{J}\}_{i\in I}$ cover $Z(\mathbf{J})$.
\end{enumerate}
Then the underlying simplicial sheaf on topological spaces $Z_!(F)$ is representable, that is, equivalent to $j(X)$ for some topological space $X$. If $X$ is paracompact Hausdorff, then $F$ is a Kuranishi space.
\end{thm}
The proof of this theorem uses a variety of techniques developed in the body of this paper that would take us too far afield to explain at this point. Instead, we draw attention to the vast simplification the theorem affords:
\begin{itemize}
    \item It is \emph{hard} to endow a concrete topological space $X$ with the structure of a Kuranishi space: one has to find a cover by affine Kuranishi models and \emph{construct} by hand an infinite tower of higher coherences on overlaps. 
    \item It is \emph{easy} to show that a simplicial sheaf $F$ is a Kuranishi space: one has to find a collection of affine Kuranishi models mapping to $F$ and \emph{check} a small number of (in practice reasonably verifiable) conditions.  
\end{itemize}
This juxtaposition highlights the dichotomy between \emph{structure} and \emph{property}: a simplicial sheaf on $\mathsf{AffKur}_{\simp}$ already has all the structure of a `generalized $\cinfty$-space' we could hope to endow it with. What is left is to verify that it lives in a subcategory whose objects satisfy suitable geometricity conditions, like having a tangent complex (compare with the task of showing that a given topological space admits a smooth manifold structure versus the problem of showing it is a topological manifold).
\begin{rmk}
Since the functor
\[ h\fun_{\mathsf{Top}^{\mathrm{open}}}(P_I,\int_{\mathsf{Top}^{\mathrm{open}}}\shv_{\mathrm{Ch}_{\R}}(\_)^{op}) \longrightarrow  h\fun^{\mathrm{ho.coh}}_{\mathsf{Top}^{\mathrm{open}}}(\mathfrak{C}[P_I],\int_{\mathsf{Top}^{\mathrm{open}}}\shv_{\mathrm{Ch}_{\R}}(\_)^{op}_{\simp})  \]
is an equivalence, every homotopical descent datum for complexes of sheaves can be \emph{strictified} to a strict descent datum. This strictification involves performing a simplicially enriched \emph{left Kan extension} along the functor $\mathfrak{C}[P_I]\rightarrow P_I$ \emph{relative to the projection $p_{\simp}$}. This procedure critically requires the existence of homotopy colimits in $\int_{\mathsf{Top}^{\mathrm{open}}}\shv_{\mathrm{Ch}_{\R}}(\_)^{op}_{\simp})$. In contrast, the simplicial category $\mathsf{AffKur}_{\simp}$ does not admit many homotopy colimits, and the analogue of strictification for diagrams in $\mathsf{AffKur}_{\simp}$ is \emph{false}. For instance, Theorem \ref{thm:repcritkur} above crucially requires that we use the more complicated, but also more flexible homotopy coherent Definition \ref{defn:kuratlas3} of a Kuranishi atlas instead of the strict Definition \ref{defn:kuratlas4}. In other words, given a simplicial sheaf $F$ that satisfies the conditions of Theorem \ref{thm:repcritkur}, we can find a simplicial diagram $\mathfrak{C}[P_I]\rightarrow \mathsf{AffKur}_{\simp}$ (carrying all 1-morphisms to localizations) with homotopy colimit $F$ in $\shv_{\mathsf{Kan}}(\mathsf{AffKur}_{\simp})^c$ that \emph{in general cannot be strictified} to a diagram $P_I\rightarrow \mathsf{AffKur}$. Thus, paradoxical though it may appear at the moment, including the seemingly much more involved notion of a homotopy coherent diagram of Kuranishi spaces renders our theory significantly easier to work with. 
\end{rmk}
\begin{rmk}
It is possible to replace $\mathsf{AffKur}_{\simp}$ by a larger Dwyer-Kan equivalent Kan enriched category $\widehat{\mathsf{AffKur}}_{\simp}$ for which every homotopy coherent diagram $\mathfrak{C}[P_I]\rightarrow \widehat{\mathsf{AffKur}}_{\simp}$ \emph{can} be strictified. For instance, let us take $\widehat{\mathsf{AffKur}}_{\simp}$ the (homotopy) essential image of the Yoneda embedding $j:\mathsf{AffKur}_{\simp}\hookrightarrow\shv_{\mathsf{Kan}}(\mathsf{AffKur}_{\simp})^c$, then one can show that for every poset $P$ and every homotopy coherent diagram $\mathcal{J}:\mathfrak{C}[P]\rightarrow \widehat{\mathsf{AffKur}}_{\simp}$, there exists a homotopy coherent natural transformation \[\overline{\mathcal{J}}:\mathfrak{C}[P\times\Delta^1]\longrightarrow \widehat{\mathsf{AffKur}}_{\simp}\]
such that the restriction $\overline{\mathcal{J}}|_{\mathfrak{C}[P\times\{0\}]}$ equals $\mathcal{J}$, the restriction 
$\overline{\mathcal{J}}|_{\mathfrak{C}[P\times\{1\}]}$ factors via $\mathfrak{C}[P]\rightarrow P$ and for each $E\in P$, the map $\overline{\mathcal{J}}|_{\mathfrak{C}[\{E\}\times\Delta^1]}:\Delta^1\rightarrow\widehat{\mathsf{AffKur}}_{\simp}$ is an equivalence. The obvious cost of this maneuver is that most of the objects of the simplicial category $\widehat{\mathsf{AffKur}}_{\simp}$ are simplicial sheaves that do not admit a description as affine Kuranishi models, so that a diagram $P_I\rightarrow \widehat{\mathsf{AffKur}}_{\simp}$ lacks geometric interpretation. 
\end{rmk}
\begin{rmk}\label{rmk:kuriso}
We have not mentioned Kuranishi atlases with isotropy, but we remark here that Theorem \ref{thm:repcritkur} is the starting point for gluing more general Kuranishi-type structures in a manner that entirely avoids spelling out coordinate changes and complicated overlap conditions for atlases. Suppose that $\Gamma$ is a finite group and that $\mathbf{K}=(M,p:E\rightarrow M,s)$ is a $\Gamma$-equivariant Kuranishi atlas (that is, both $M$ and $E$ admit $\Gamma$-actions and $p$ and $s$ are $\Gamma$-equivariant), then we have an action map $m:\Gamma\times \mathbf{K}\rightarrow\mathbf{K}$ and we can form a simplicial diagram $\simp^{op}\rightarrow \mathsf{AffKur}$ like so
\[\begin{tikzcd} \ldots\ar[r,shift left=3] \ar[r,shift left] \ar[r,shift right]\ar[r,shift right=3] & \Gamma\times \Gamma\times\mathbf{K} \ar[r,shift left=2] \ar[r]\ar[r,shift right=2] &\Gamma\times\mathbf{K}\ar[r,shift left,"m"]\ar[r,shift right,"\pi"'] & \mathbf{K} \end{tikzcd}  \] 
where $\pi$ is the projection onto $\mathbf{K}$, the \emph{nerve of the action groupoid}. Composing with the localization functor from $\mathsf{AffKur}$ to $\mathsf{AffKur}_{\simp}$ determines a simplicial diagram $\simp^{op}\rightarrow \mathsf{AffKur}_{\simp}$. Composing with the Yoneda embedding $j:\mathsf{AffKur}_{\simp}\hookrightarrow\shv_{\mathsf{Kan}}(\mathsf{AffKur}_{\simp})^c$, we can form the homotopy colimit $[\mathbf{K}/\Gamma]$ of this object in the simplicial category of simplicial sheaves, that is, the homotopy quotient of $\mathbf{K}$ by $\Gamma$. We call such a simplicial sheaf an \emph{affine Kuranishi model with isotropy}. Suppose that for a simplicial sheaf $F$ we have a collection of maps $\{[\mathbf{K}_i/\Gamma_i]\rightarrow F\}_{i \in I}$ for which $(1)$, $(2)$ and $(3)$ of Theorem \ref{thm:repcritkur} are satisfied, then $F$ is a \emph{Kuranishi space with isotropy}. In general, the underlying simplicial sheaf $Z_!(F)$ on the category of topological spaces will be representable by an \emph{orbispace}. It can be shown that such an object also admits a tangent complex. For each $x\in X$, the tangent complex of $F$ at $x$ has a natural structure of a chain complex of left $\R[\Gamma_x]$-modules, where $\Gamma_x$ is the isotropy group of $F$ at $x$.
\end{rmk}
Let us return to $Q1$: our theory suggests a new two-step strategy for proving that some topological space $X$ admits a Kuranishi atlas.
\begin{description}
    \item[Step 1] Give $X$ the structure of a simplicial sheaf on $\mathsf{AffKur}_{\simp}$, that is, find a simplicial sheaf $F_X$ such that the underlying simplicial sheaf $Z_!(F_X)$ on $\mathsf{Top}$ is representable by the topological space $X$.
    \item[Step 2] Find a collection of maps $\{j(\mathbf{K}_i)\rightarrow F_X\}_{i\in I}$ that satisfies the conditions of Theorem \ref{thm:repcritkur}. 
\end{description}
Supposing that we can perform Steps 1 and 2, Theorem \ref{thm:repcritkur} yields a simplicial diagram $\mathfrak{C}[P_I]\rightarrow \mathsf{AffKur}_{\simp}$ for some index set $I$ whose homotopy colimit in $\shv_{\sset}( \mathsf{AffKur}_{\simp})$ is $F_X$, which yields a Kuranishi structure on $X=Z_!(F_X)$. If we specialize to $X=\mathsf{Sol}(P)$, the space of solutions of an elliptic PDE on a compact manifold, it turns out that there is an easy and conceptually straightforward way to give $\mathsf{Sol}(P)$ the structure of a sheaf on $\mathsf{AffKur}_{\simp}$. As we will see in a moment, this involves no analysis of infinite dimensional function spaces, but only formal categorical constructions. In contrast, to carry out Step 2, we will need to use the kind of arguments that go into proving Fact \ref{fact:kur} which rest on a large body of geometric and functional analysis.
\subsubsection{Representability of moduli of solutions}
To carry out Step 1 for $\mathsf{Sol}(P)$ we use that the simplicial category $\shv_{\mathsf{Kan}}(\mathsf{AffKur}_{\simp})^c$ has very good categorical properties (it is an \emph{\inftopt}), among which are the following.
\begin{enumerate}[$(a)$]
    \item $\shv_{\mathsf{Kan}}(\mathsf{AffKur}_{\simp})^c$ does not only admits all homotopy colimits, but also all homotopy limits which sometimes interact nicely with each other. For instance, the formation of homotopy pullbacks along a given map $F\rightarrow G$ of simplicial sheaves preserves all homotopy colimits: we have for each diagram $\mathcal{J}:K\rightarrow \shv_{\mathsf{Kan}}(\mathsf{AffKur}_{\simp})^c$ a canonical isomorphism 
\[ \hocolim_{k \in K}(\mathcal{J}(i)\times^h_G F) \overset{\cong}{\longrightarrow} (\hocolim_{k\in K}\mathcal{J}(k))\times^h_{G}F \]
in the homotopy category $h\shv_{\mathsf{Kan}}(\mathsf{AffKur}_{\simp})^c$.
    \item A simplicial functor
    \[  T:(\shv_{\mathsf{Kan}}(\mathsf{AffKur}_{\simp})^c)^{op}\longrightarrow \sset \]
    that is, a simplicial presheaf on simplicial sheaves on $\mathsf{AffKur}_{\simp}$, preserves all homotopy limits if and only if there is a simplicial sheaf $G$ such that the Yoneda functor $\Hom(\_,G)$ is equivalent to $T$.  
\end{enumerate}
It follows from $(a)$ that the simplicial functor
\[   (\shv_{\mathsf{Kan}}(\mathsf{AffKur}_{\simp})^c)^{op}\longrightarrow \sset,\quad\quad F\longmapsto \Hom(F\times X,Y)\]
preserves homotopy limits, so $(b)$ implies that this functor is representable by an object that we will denote by $\map(X,Y)$: it is the simplicial sheaf that carries an affine Kuranishi model $\mathbf{K}$ to the Kan complex of \emph{$\mathbf{K}$-parametrized maps from $X$ to $Y$}. If we pick $X$, $Y$ and $\mathbf{K}$ equal to $j(M)$, $j(N)$ and $P$ for three manifolds $M,N,P$ (viewed as affine Kuranishi models), then we have a homotopy equivalence $\map(M,N)(P)\simeq \cinfty(P\times M,N)$, the set of $P$-parametrized $\cinfty$ maps from $M$ to $N$. More generally, if $Y\rightarrow M$ is a fibre bundle over a manifold, then there is a simplicial sheaf $\map_M(M,Y)$ of sections, which is defined so that there is a homotopy equivalence of Kan complexes 
\[ \map_M(M,Y)(\mathbf{K}) \simeq \Hom(\mathbf{K}\times M,Y )\times^h_{\Hom(\mathbf{K}\times M,M )}\{\pi\},  \]
natural in $\mathbf{K}$, where the homotopy pullback denotes the Kan complex of maps $\mathbf{K}\times M\rightarrow Y$ for which the composition with $Y\rightarrow M$ agrees with the obvious projection $\pi:\mathbf{K}\times M\rightarrow M$. Recall that we are given the following data of a compact manifold $M$, a fibre bundle $Y\rightarrow M$, a vector bundle $V\rightarrow M$ and a nonlinear elliptic PDE $P$ acting between sections of $Y$ and $V$. Then one can prove that $P$ naturally extends to a map \[\widehat{P}:\map_M(M,Y)\longrightarrow \map_M(M,V) \]
of simplicial sheaves on $\mathsf{AffKur}_{\simp}$. To define $\mathsf{Sol}(P)$ as a simplicial sheaf, we recall that solutions and zero sets are categorically represented as pullbacks. In fact, Remark \ref{rmk:affkurpullback} tells us that we can recover all finite dimensional affine Kuranishi models as homotopy pullbacks via the derived intersection of the defining section with the zero section. Thus, the only reasonable candidate for the sheafy enhancement of $\mathsf{Sol}(P)$ is the homotopy pullback of $\widehat{P}$ with the zero section.
\begin{thm}[Elliptic representability]\label{thm:ellrepsimp}
Let $\mathsf{Sol}(\widehat{P})$ be the cone in the homotopy pullback diagram 
\[ 
\begin{tikzcd}
\mathsf{Sol}(\widehat{P}) \ar[d] \ar[r] &\map_M(M,Y) \ar[d,"\widehat{P}"] \\
0 \ar[r] & \map_M(M,V)
\end{tikzcd}
\]
of simplicial sheaves on $\mathsf{AffKur}_{\simp}$. Then the following hold true.
\begin{enumerate}[$(1)$]
    \item The object $Z_!(\mathsf{Sol}(\widehat{P}))$ is representable by the topological space $\mathsf{Sol}(P)$.
    \item There is a collection of maps $\{j(\mathbf{K}_i)\rightarrow \mathsf{Sol}(\widehat{P})\}$ that satisfies the representability conditions of Theorem \ref{thm:repcritkur}. Thus, $\mathsf{Sol}(\widehat{P})$ is a Kuranishi space.
\end{enumerate}
\end{thm}
We conclude that the topological space $\mathsf{Sol}(P)$ admits a \emph{canonical} Kuranishi structure (but, of course, not a canonical Kuranishi atlas) whose tangent complex coincides with the sheaf $\{TP_x\}_{x\in \mathsf{Sol}(P)}$. The maps $j(\mathbf{K}_i)\rightarrow \mathsf{Sol}(\widehat{P})$, whose existence is guaranteed by the theorem are of course the local finite dimensional reductions of Fact \ref{fact:kur}. It is a nontrivial task involving somewhat nonstandard methods of functional analysis and elliptic theory to shows that the `familiar' construction of local finite dimensional reductions actually determines a morphism to the simplicial sheaf $\mathsf{Sol}(\widehat{P})$. We take it up in part III of this series \cite{cinftyIII}.
\begin{rmk}
In the situation of Theorem \ref{thm:ellrepsimp}, imposing certain finiteness conditions guarantees that the Kuranishi structure on $\mathsf{Sol}(P)$ admits a very simple description. Suppose that the dimensions of the vector spaces $\ker TP_x$ and $\coker TP_x$ are bounded as $x$ ranges over $\mathsf{Sol}(P)$, then there is an affine Kuranishi model $\mathbf{K}$ and and equivalence $j(\mathbf{K})\simeq \mathsf{Sol}(\widehat{P})$; in other words, there is a triple $(M,E,s)$ such that $\mathsf{Sol}(P)=Z(s)$ which identifies the tangent complexes of $P$ and $s$ at each $x\in\mathsf{Sol}(P)$. We will prove (a generalization of) this in part III as well.
\end{rmk}
\begin{rmk}
We can, of course, define the homotopy pullback diagram in the theorem above for \emph{any} nonlinear PDE, but $(2)$ will usually not hold; in fact, the ellipticity of $P$ is essentially equivalent to the assertion that $\mathsf{Sol}(\widehat{P})$ is a Kuranishi space.
\end{rmk}
\begin{rmk}
In the situation of Theorem \ref{thm:ellrepsimp}, suppose that $M$ admits a $\cinfty$ action by a finite group $\Gamma$ and suppose that the PDE $\widehat{P}$ is equivariant for the induced actions on the mapping sheaves $\map_M(M,Y)$ and $\map_M(M,V)$, then the general yoga of higher topos theory provides a $\Gamma$ action on $\mathsf{Sol}(\widehat{P})$ and the associated quotient object $[\mathsf{Sol(\widehat{P})}/\Gamma]$ is a Kuranishi space with isotropy in the sense of Remark \ref{rmk:kuriso} (note that this works just as well if $\Gamma$ is a Lie group, or if the group acts on $Y$ and $V$). Classically, one obtains Kuranishi structures with isotropy via an equivariant version of Construction \ref{cons:sumchart}. This introduces a host of 
analytic difficulties and smoothness issues, as the group $\Gamma$ does not act in a $\cinfty$ fashion on any Sobolev completion of the mapping space $\map_M(M,Y)$. Our approach does not avoid these difficulties entirely; they are hidden in the proof of Theorem \ref{thm:ellrepsimp}. 
\end{rmk}
\subsubsection{The language of \infcatst}
We have made the case that a theory that can handle singular spaces locally cut out by non transverse $\cinfty$ functions will necessitate a substantial amount of `higher mathematics'. The approach to the theory of derived $\cinfty$-geometry we sketched above is a legitimate one, but it is not without difficulties.
\begin{enumerate}[$(D1)$]
    \item We defined the Kan enriched category of affine Kuranishi models as the natural simplicial generalization of the procedure of localizing a category at a set of weak equivalences, which was the obvious choice in the context of our descent problem for the tangent complex. While this hammock localization (followed by the operation $\sing\circ |\_|$) has good categorical properties, it outputs a simplicial category whose Kan complexes of morphisms are rather mysterious. For instance, recall that we can view $\mathbf{K}$ for $\mathbf{K}=(X,E,s)$ as a dg-manifold with associated dg-algebra of global functions $\Gamma(\Lambda^{\bullet} E^{\vee})$ with differential given by $\iota_s$, the contraction with the section $s$. We should expect that the Kan complex of `functions' between $\mathbf{K}$ and $\R$ in the simplicial category $\mathsf{AffKur}_{\simp}$ recovers this dg-algebra in the following sense: Remark \ref{rmk:affkurpullback} asserts in particular that $\mathsf{Man}\rightarrow \mathsf{AffKur}_{\simp}$ preserves products, so addition and multiplication in $\R$ make the $\Z_{\geq0}$-graded set $\pi_*(\Hom_{\mathsf{AffKur}_{\simp}}(\mathbf{K},\R))$ into a graded $\R$-algebra and there should be a canonical isomorphism 
    \[ \pi_*(\Hom_{\mathsf{AffKur}_{\simp}}(\mathbf{K},\R)) \simeq H_*(\Gamma(\Lambda^{\bullet} E^{\vee}),\iota_s)  \]
    of graded $\R$-algebras, but this is not at all obvious. In fact, it is not even so straightforward to prove that the functor $\mathsf{Man}\rightarrow\mathsf{AffKur}_{\simp}$ is homotopically fully faithful.
    \item In the story we told above we stated several theorems without proof to avoid getting bogged down in technicalities. Currently, these proofs are not quite within reach, partly for the reason explained above that we do not have a solid grasp on the simplicial category $\mathsf{AffKur}_{\simp}$, but also because the language we have been using is not conducive to streamlined proofs. For Theorem \ref{thm:repcritkur} to hold for example, we needed to replace the notion of a diagram in a category with the much more involved notion of a homotopy coherent diagram in a simplicial category. Replacing a poset $P$ with its simplicial thickening $\mathfrak{C}[P]$ is an instance of \emph{cofibrant replacement}: the correct space of maps between two objects $A$ and $B$ in a model category (like the model category of simplicial categories) is obtained by replacing $A$ with an equivalent cofibrant object, and $B$ by an equivalent fibrant object. As we saw, even for a simple object like the poset $P_I$, it is rather cumbersome to work directly with the simplicial category $\mathfrak{C}[P_I]$. If we were to develop our theory using simplicial categories, we would at each step have to ensure that the objects we work with are appropriately fibrant or cofibrant. We would like to place ourselves in a setting of higher category theory where we do not have to worry about the distinction between homotopy coherent and strict diagrams, homotopy Cartesian and Cartesian morphisms and homotopy colimits and ordinary colimits. 
\end{enumerate}
We come back to $D1$ in a moment; first we address $D2$ by moving our discussion into the framework of \emph{\infcatst}. For us, an \infcat is a \emph{quasi-category} in the sense of Joyal or a \emph{weak Kan complex} in the sense of Boardman-Vogt.
\begin{defn}
An \emph{\infcatt} is a simplicial set $X$ such that for each $n\geq 2$, each $0<i<n$ and each \emph{inner} horn $\Lambda^n_i\rightarrow X$, there is a dotted lift as in the diagram 
\begin{equation}\label{eq:innerhorn}
\begin{tikzcd}
\Lambda^n_i\ar[d,hook] \ar[r]& X \\
\Delta^n.\ar[ur,dotted]
\end{tikzcd}
\end{equation}
\end{defn}
Recall the nerve functor $\ner:\mathsf{Cat}\rightarrow \sset$. It is not hard to see that for every category $\icate$ and every inner horn $\Lambda^n_i\rightarrow \ner(\icate)$, there is a \emph{unique} lift 
\[
\begin{tikzcd}
\Lambda^n_i\ar[d,hook] \ar[r]& \ner(\icate) \\
\Delta^n,\ar[ur,dotted,"\exists !"']
\end{tikzcd}
\]
see for instance Proposition 1.1.2.2 of \cite{HTT}. In fact, the nerve functor \emph{identifies} the category $\mathsf{Cat}$ of categories with the full subcategory $\sset$ spanned by those \infcats $X$ for which the diagram \eqref{eq:innerhorn} admits a unique lift. For $n=2$, the lifting problem above expresses that for any two composable maps in an \infcatt, one may choose a composition. For \infcatst, this composition is not unique in general and is not associative on the nose. The solubility of the higher lifting problems expresses that composition is associative up to coherent homotopy. \\
We will freely use the extensive foundational body of work by Lurie \cite{HTT,HA} which gives meaning to all familiar concepts of category theory, like diagram categories, (co)limits, adjoint functors, the Yoneda embedding, (commutative) algebras and modules over them, and many more \emph{internally to \infcatst}. In the setting of simplicial categories, working with homotopy coherent diagrams was arduous due to appearance of the intricate simplicial thickening $\mathfrak{C}[\_]$. In contrast, diagrams in \infcats are entirely straightforward.
\begin{ex}[Functors and diagrams in \infcatst]\label{ex:funcats}
Let $\icat$ be an \infcatt. The category $\sset$ admits the \emph{Joyal model structure}. The fibrant objects are precisely the \infcatst, and \emph{every object for this model structure is cofibrant}. Thus, we can define diagrams without having to resort to cofibrant replacement: let $K$ be a simplicial set, then a \emph{functor} from $K$ to $\icat$ is (simply) a map $f:K\rightarrow \icat$ of simplicial sets. The \emph{\infcat of functors} $\fun(K,\icat)$ is the simplicial set defined as follows: for every simplicial set $K'$, there is a canonical bijection 
\[ \Hom_{\sset}(K',\fun(K,\icat)) \cong \Hom_{\sset}(K\times K',\icat). \]
It is not hard to show that $\fun(K,\icat)$ is again an \infcatt.
\end{ex}
\begin{ex}[Simplicial categories as \infcatst]\label{ex:simpcatsinfcats}
The operation $\mathfrak{C}[\_]$ determines a functor from the category of posets to the category of simplicial categories. Restricting to linearly ordered sets, we have a functor $\mathfrak{C}[\_]:\simp\rightarrow \cat_{\simp}$, which induces a functor $\cat_{\simp}\rightarrow \sset$. This is the \emph{homotopy coherent nerve} and is also denoted by $\ner$; it takes a simplicial category $\icat$ to the simplicial set whose set of $n$-simplices is $\Hom_{\cat_{\simp}}(\mathfrak{C}[n],\icat)$. Observe that if $\icat$ is an ordinary category, then the homotopy coherent nerve coincides with the ordinary nerve, so using the symbol $\ner$ for both functors introduces no ambiguity. The functor $\ner$ admits a left adjoint, also denoted $\mathfrak{C}[\_]$, and it was shown by Lurie that the adjunction $(\mathfrak{C}[\_]\adj \ner)$ identifies the homotopy theory of simplicial categories with the homotopy theory of \infcats \cite{HTT}. In particular, if $\icat$ is Kan enriched, then $\ner(\icat)$ is an \infcatt. To summarize, simplicial categories serve as another model for \infcats in which composition is \emph{strictly} associative. This feature of simplicial categories actually renders them technically less convenient, as most higher categories that appear in practice have associative composition up to coherent homotopy. With the language we now have at our disposal, we can for example write down the \infcat of homotopy coherent diagrams \[\mathfrak{C}[P_I]\longrightarrow \int_{\mathsf{Top}^{\mathrm{open}}}\shv_{\mathrm{Ch}_{\R}}(\_)^{op}_{\simp}\]
for an open cover $\{U_i\subset X\}_i$ as the simplicial set \[\fun(\ner(P_I),\ner(\int_{\mathsf{Top}^{\mathrm{open}}}\shv_{\mathrm{Ch}_{\R}}(\_)^{op}_{\simp}))\times_{\fun(\ner(P_I),\ner(\mathsf{Top}^{\mathrm{open}}))}\times\{f\} \]
where the fibre is taken at the map $f:P_I\rightarrow \mathsf{Top}^{\mathrm{open}}$ specifying the open cover. This simplicial set is an \infcat not isomorphic to the nerve of any simplicial category, reflecting that composition is not strictly associative.
\end{ex}
Working with \infcats may at times be challenging, as the combinatorics of simplicial sets can be intricate. Fortunately, the theory is at this point so well developed that one need not be very well versed in the simplicial technicalities to avail oneself of the power of the formalism. Very often, one can reason `at a high level', employing the usual `1-categorical' semantics with the understanding that every definition or result should be replaced with its $\infty$-categorical counterpart, which one should find somewhere in the works \cite{HTT,HA,sag}. This results in a minimal effort strategy for finding statements of theorems with the correct homotopical content. As a case in point, we offer the following example which, with the hard work we did to define a good category of Kuranishi spaces in mind, should be compelling.
\begin{ex}\label{ex:univpropdg}
The Dwyer-Kan hammock localization $L_H$ is the derived functor of the usual localization, which was defined via a universal property. The localization $L_H$ therefore enjoys a similar, `higher' version of this universal property which, with the notion of a functor \infcat at our disposal, we can now state.
\begin{enumerate}
    \item[$(*)$] Let $\mathbf{A}$ be an ordinary category with a subcategory $W$ of weak equivalences. Consider the functor $f:\mathbf{A}\rightarrow L_H(\mathbf{A},W)\rightarrow \sing|L_H(\mathbf{A},W)|$. The simplicial category $\sing|L_H(\mathbf{A},W)|$ is Kan enriched by design and we may view $\mathbf{A}$ as a Kan enriched category with only nondegenerate $0$-simplices, so applying the homotopy coherent nerve functor $\mathbf{N}$, we obtain a functor
    \[ \ner(f):\ner(A)\longrightarrow\ner(\sing|L_H(\mathbf{A},W)|) \]
    of \infcatst. Then for every \infcat $\icat$, composition with $\ner(f)$ induces an equivalence
    \[ \fun(\ner(\sing|L_H(\mathbf{A},W)|),\icat)\overset{\simeq}{\longrightarrow} \fun_W(\ner(A),\icat)    \]
    of \infcatst, where $\fun_W(\ner(A),\icat)\subset \fun(\ner(A),\icat)$ is the full sub-\infcat spanned by those functors $g:\ner(A)\rightarrow\icat$ that carry every map in $W$ to an equivalence in $\icat$.
\end{enumerate}
Note that modulo the symbols $\ner$ and $\infty$, this is exactly the same statement as we had for the Gabriel-Zisman localization. We also have a similar uniqueness statement: there is, up to equivalence of \infcatst, a unique \infcat $\ner(A)[W^{-1}]$ that comes equipped with a functor $\ner(A)\rightarrow\ner(A)[W^{-1}]$ satisfying the universal property $(*)$ above. It follows that we could have simply \emph{defined} the \infcat $\ner(\mathsf{AffKur}_{\simp})$ as the full subcategory of the $\infty$-categorical localization $\ner(\mathbf{dgMan})[W^{-1}]$ spanned by the objects equivalent to an affine Kuranishi model. Note that here $\ner(\mathbf{dgMan})$ is just the ordinary nerve of the 1-category of dg-manifolds, while $\ner(\mathsf{AffKur}_{\simp})$ is the coherent nerve. Now recall our Definition \ref{defn:kuratlas3} of a Kuranishi structure on a space $X$: it is a cover $\{U_i\subset X\}$ and a lift of the associated functor $\mathfrak{C}[P_I]\rightarrow\mathsf{Top}^{\mathrm{open}}$ to $\mathsf{AffKur}_{\simp}$ carrying every map to a localization of affine Kuranishi models. By definition of the coherent nerve of example \ref{ex:simpcatsinfcats}, we see that a Kuranishi structure on $X$ is \emph{precisely} a cover $\{U_i\subset X\}$ and a lift 
\[
\begin{tikzcd}
& \ner(\mathsf{AffKur}_{\simp})\ar[d,"Z"] \ar[d]\\
\ner(P_I) \ar[ur,dotted,"\tilde{f}"] \ar[r,"f"] & \ner(\mathsf{Top}^{\mathrm{open}}),
\end{tikzcd}
\]
carrying every map to a localization, where again $\ner(P_I)$ and $\ner(\mathsf{Top}^{\mathrm{open}})$ are the ordinary nerves of 1-categories. Combining these observations, we conclude that, if at the start of our discussion, right after Definition \ref{def:kurweakeq} of a weak equivalence of dg-manifolds, we had taken the nerve to identify categories with \infcats having unique horn fillers and replaced all subsequent steps (like the localization at weak equivalences of dg-manifolds) by their $\infty$-categorical counterparts, we would have found the correct definition of a Kuranishi atlas on our first attempt. 
\end{ex}
\subsubsection{Dg-manifolds and derived $\cinfty$-rings}
Let us return to the first difficulty $D1$ of the theory we outlined in the previous subsections. With Example \ref{ex:univpropdg} in mind, we rephrase the definition of $\mathsf{dgMan}_{\simp}$.
\begin{defn}\label{def:dmanloc}
The \infcat of \emph{dg-manifolds}, denoted $\mathsf{dgMan}_{\infty}$, is the $\infty$-categorical localization of the category $\ner(\mathbf{dgMan})$ (viewed as an \infcat with uniquely defined composition) at the weak equivalences $W$ introduced in Definition \ref{def:kurweakeq}.
\end{defn}
Unfortunately, it is quite hard to understand general localizations of \infcatst. To get a better grasp on the \infcat of dg-manifolds, we change perspective; instead of studying the universal \infcat that inverts the weak equivalences of dg-manifolds, we turn our attention to the universal \infcat that adjoins finite limits\footnote{We will also add retracts for idempotents, which is a technical condition the reader may safely ignore for the moment} to the category $\mathsf{Man}$.
\begin{defn}\label{def:dmanfinlim} 
The \infcat of \emph{affine derived manifolds} is the universal \infcat admitting finite limits (and retracts for idempotents) equipped with a functor 
\[  f:\ner(\mathsf{Man}) \longrightarrow \daff \]
that preserves transverse pullback diagrams of manifolds. More precisely, for any \infcat $\icat$ that admits finite limits and retracts for idempotents, we denote by 
\[ \fun^{\pitchfork}(\ner(\mathsf{Man}),\icat)\subset \fun(\ner(\mathsf{Man}),\icat) \]
the full sub-\infcat spanned by functors that preserve transverse pullback diagrams of manifolds. Then for any \infcat $\icat$ that admits finite limits (and retracts for idempotents), composition with $f$ induces an equivalence \[  \fun^{\mathrm{lex}}(\daff,\icat)\overset{\simeq}{\longrightarrow} \fun^{\pitchfork}(\ner(\mathsf{Man}),\icat)\]
of \infcatst, $\fun^{\mathrm{lex}}(\daff,\icat)\subset \fun(\daff,\icat)$ being the full sub-\infcat spanned by functors that preserve all finite limits (the superscript $\mathrm{lex}$ is short for \emph{left exact} which should be taken as synonymous to `finite limit preserving').
\end{defn}
The two universal properties of Definitions \ref{def:dmanloc} and \ref{def:dmanfinlim} may not seem to be obviously related, but one can show that they yield the same \infcatt.
\begin{thm}\label{thm:dgmandaff}
The \infcat $\mathsf{dgMan}_{\infty}$ admits all finite limits (and retracts for idempotents) and the functor $\ner(\mathsf{Man})\rightarrow \mathsf{dgMan}_{\infty}$ lies in the \infcat $\fun^{\pitchfork}(\ner(\mathsf{Man}),\mathsf{dgMan}_{\infty})$. Moreover, the functor $\daff\rightarrow \mathsf{dgMan}_{\infty}$ provided by the universal property of $\daff$ is an equivalence of \infcatst.
\end{thm}
This was proven very recently by Carchedi \cite{carchedidgman}; we will give a give a slightly different proof in part II \cite{cinftyII}. The upshot of this result is that we may now forget about the category of dg-manifolds and focus our attention on the \infcat $\daff$ defined via the universal property of Definition \ref{def:dmanfinlim}. We are thus reduced to studying functors $f:\ner(\mathsf{Man})\rightarrow \icat$ into \infcats admitting finite limits (and retracts for idempotents) that preserve transverse pullback diagrams of manifolds, which as it turns out, is much more manageable than studying functors $\ner(\mathbf{dgMan})\rightarrow \icat$ which carry weak equivalences to equivalences. 
\begin{rmk}
By Whitney's embedding theorem, every manifold $M$ admits an embedding $M\hookrightarrow \R^{n}$ for some $n$. By the tubular neighbourhood theorem, there is an open neighbourhood $M\subset U\subset \R^n$ diffeomorphic to the normal bundle $N_{M/\R^n}$ of $M$ in $\R^n$. It follows that every manifold is a retract of an open subset in some Cartesian space. Every open subset $U\subset\R^n$ has a \emph{characteristic function}: a $\cinfty$ function $\chi:\R^n\rightarrow\R$ that is nonzero precisely on $U$, and $U$ is diffeomorphic to the zero set of the function $\chi y-1:\R^{n+1}\rightarrow\R$, where $y$ denotes the last coordinate. Thus, a functor $\mathsf{Man}\rightarrow \icat$ that carries transverse pullback diagrams of manifolds to pullback diagrams in $\icat$ is completely determined by its values on the model spaces $\R^n$.
\end{rmk}
Here is a more precise version of this statement.
\begin{prop}
Let $\cartsp\subset\mathsf{Man}$ be the full subcategory spanned by \emph{Cartesian spaces} of the form $\R^n$ for some $n\geq 0$. Let $\icat$ be an \infcat that admits finite products, then we denote by $\fun^{\pi}(\ner(\cartsp),\icat)\subset \fun(\ner(\cartsp),\icat)$ the full subcategory spanned by functors preserving finite products. Note that for each \infcat $\icat$ admitting finite limits and idempotents, composition with the full subcategory inclusion $\cartsp\subset\mathsf{Man}$ induces a functor
\[ \fun^{\pitchfork}(\ner(\mathsf{Man}),\icat)\overset{\theta}{\longrightarrow}  \fun^{\pi}(\ner(\cartsp),\icat)\]
(observe that a product is a particular instance of a transverse pullback diagram, namely a pullback over the point). Then $\theta$ is an equivalence of \infcatst.
\end{prop}
Thus, we are further reduced to studying finite product preserving functors $\ner(\cartsp) \rightarrow \icat$. The \infcat of such functors has excellent formal properties akin to those enjoyed by the category of commutative algebras. In fact, if $\mathrm{T}$ is a category that is generated under finite products by a single object $t$ (like how $\cartsp$ is generated by $\R$), we should think of a product preserving functor $\ner(\mathrm{T})\rightarrow \icat$ as some sort of algebra internal to the \infcat $\icat$ whose multiplication operations are determined by the maps in the category $\mathrm{T}$. This point of view is due to Lawvere and classical in the formal theory of algebraic structures: here is a motivating example.
\begin{ex}
Let $\mathsf{Poly}_{\R}\subset\cartsp$ be the subcategory of \emph{polynomial maps}. An ordinary functor $f:\mathsf{Poly}_{\R}\rightarrow \set$ that preserves products is determined by the set $f(\R)$; denote this set by $A$. By functoriality, the multiplication and addition on $\R$ determine maps 
\[ +:A\times A\longrightarrow A,\quad \quad \times:A\times A\longrightarrow A  \]
that are associative and commutative, which make the set $A$ into a commutative $\R$-algebra. It is not hard to see that the assignment $f\mapsto A$ determines an equivalence of categories $\fun^{\pi}(\mathsf{Poly}_{\R},\set)\cong \mathsf{CAlg}^0_{\R}$. Turning to the setting of simplicial categories or \infcatst, we should replace the category $\set$ with the \infcat of Kan complexes, which we will denote by $\spa=\ner(\mathsf{Kan})$. Recall that because we are working with \infcatst, a functor $f:\ner(\mathrm{T})\rightarrow\spa$ (i.e. a map of simplicial sets) preserving products is automatically a homotopy coherent functor preserving homotopy products. Thus, talking about product preserving functors $f:\ner(\mathrm{T})\rightarrow\spa$ is a precise and concise way of describing Kan complexes endowed with an $\R$-module structure up to coherent homotopy and a compatible binary operation that is fully coherently associative and commutative.  
\end{ex}
A very large class of algebraic structures (commutative monoids, abelian groups, Lie algebras,...) admit a description as a category of finite product preserving functors to the category $\set$. 
\begin{defn}
A \emph{$\cinfty$-ring} is a product preserving functor $\cartsp\rightarrow \set$. We let $\cinfty\mathsf{ring}=\fun^{\pi}(\cartsp,\set)$ denote the category of $\cinfty$-rings. A \emph{simplicial $\cinfty$-ring} is a product preserving functor $\ner(\cartsp)\rightarrow\spa$. We let $\sring=\fun^{\pi}(\ner(\cartsp),\spa)$ denote the \infcat spanned by simplicial $\cinfty$-rings.
\end{defn}
\begin{ex}
Let $M$ be a manifold, then the commutative algebra $\cinfty(M)$ admits an obvious structure of a (simplicial) $\cinfty$-ring: for any smooth map $f:\R^n\rightarrow \R^m$, composing with $f$ yields the corresponding map $\cinfty(M)^n=\cinfty(M,\R^n)\rightarrow \cinfty(M,\R^m)$.
\end{ex}
\begin{rmk}
The subcategory inclusion $\mathsf{Poly}_{\R}\subset\cartsp$ obviously preserves finite products and therefore induces a forgetful functor
\[ (\_)^{\rmalg}:\sring\longrightarrow \scring_{\R}  \]
of \infcatst, carrying a simplicial $\cinfty$-ring to its underlying simplicial commutative $\R$-algebra. This functor is \emph{conservative} (if a map $f:A\rightarrow B$ is carried to an equivalence $A^{\rmalg}\overset{\simeq}{\rightarrow} B^{\rmalg}$, then $f$ is an equivalence in the \infcat $\sring$) and it preserves all limits and \emph{sifted} colimits. Many questions about simplicial $\cinfty$-rings can therefore be reduced to questions about simplicial commutative $\R$-algebras, which constitutes a significant simplification.
\end{rmk}
In this work, we will prove the following theorem, which appeared before with Carchedi \cite{univprop}.
\begin{thm}
Let $\sring_{\mathrm{fp}}\subset\sring$ be the full subcategory spanned by \emph{homotopically finitely presented} objects (a simplicial $\cinfty$-ring $A$ lying in this subcategory in particular has the property that the $\cinfty$-ring of connected components $\pi_0(A)$ can be written as $\cinfty(\R^n)/I$ for $I$ a finitely generated ideal). The functor $\cinfty(\_):\mathsf{Man}\rightarrow\sring^{op}$ preserves transverse pullbacks of manifolds and factors through $\sring_{\mathrm{fp}}^{op}$. The induced left exact functor $\daff\rightarrow \sring^{op}_{\mathrm{fp}}$ is an equivalence of \infcatst.  
\end{thm}
Morally, this theorem asserts that just as in (derived) algebraic geometry, we can completely describe derived manifolds dually in terms of their derived algebras of functions, as long as we remember that these `algebras' are simplicial $\cinfty$-rings. 
\begin{rmk}
It follows formally from Yoneda's lemma that if $X$ is a derived manifold corresponding to a simplicial $\cinfty$-ring $A$, that is, a functor $F_A:\ner(\cartsp)\rightarrow\spa$, we have an equivalence
\[\Hom_{\daff}(X,\R)=\Hom_{\sring}(j(\R),F_A)\simeq F_A(\R),\]
which we know is the Kan complex underlying $A$ (and it is also the Kan complex underlying the simplicial commutative algebra $A^{\rmalg}$). Combining this observation with Theorem \ref{thm:dgmandaff}, we see that for an affine Kuranishi model $\mathbf{K}=(X,E,s)$, the space $\Hom_{\mathsf{AffKur}_{\simp}}(\mathbf{K},\R)$, which initially seemed intractably complicated, is actually (up to equivalence) a standard homological algebraic object: it admits the structure of a simplicial $\cinfty$-ring, and the underlying derived $\R$-algebra is equivalent to the differential graded algebra $\Gamma(\Lambda^{\bullet}E^{\vee})$.
\end{rmk}
We have embedded the \infcat of affine Kuranishi models in a much richer algebraic and categorical structure. The goal of this series is to explore this structure, adapting results and constructions from derived algebraic geometry whenever possible and developing novel theory to be put to use in the study of differential geometric moduli spaces. To come full circle, we should explain how one extracts enumerative invariants from derived manifolds via virtual integration, but this would make this already verbose introduction unacceptably lengthy. We will return to this topic ($\cinfty$-motivic homotopy theory and the normal deformation stack) later in this series, as we will other important aspects of the theory that we did not mention yet at all, like derived symplectic geometry and applications to mathematical physics.

\subsection{Overview of contents}
We give a section-by-section overview of the material covered in this article.
\subsubsection*{Section 2}
The theory we will develop is paradigmatically algebro-geometric, that is, the central objects with which we will be concerned are \emph{structured spaces} of the form $(X,\Of_X)$ for a $X$ a topological space and $\Of_X$ a sheaf of algebras of some sort that is local, in a suitable sense, and generalizations of such. Moduli spaces generically have internal symmetries, that is, their geometric points come with stabilizer groups and should be treated as \emph{stacks}. When these stabilizer groups have no geometry themselves -when we are dealing with \emph{Deligne-Mumford} instead of the more general \emph{Artin} stack- we can model stacks as generalized structured spaces $\ofxtop$, where now $\xtop$ is not a topological space but a \emph{(higher) topos}. In the language of topos theory, we can interpret structure sheaves on a topos $\xtop$ as functors
\begin{equation}\label{eq:struct}
 \icat\longrightarrow\xtop 
\end{equation}  
satisfying certain limit and colimit conditions, for $\icat$ the \emph{opposite} of a small, generating subcategory of whatever category of algebraic structures under consideration (see the comprehensive \cite{elephant}, or Chapter 8 of \cite{MLM1} for an excellent introduction). In the $\infty$-categorical setting, the theory of structured higher topoi has been developed by, primarily, Lurie \cite{dagv,sag}. In Section 2, we review parts of this theory and elaborate on others. Subsections 2.1 and 2.2 are concerned with Lurie's theory of \emph{geometries} and \emph{pregeometries} which both serve as domain categories for functors of the form \eqref{eq:struct} and have the structure necessary to speak of \emph{local} rings and \emph{local} morphism. We discuss spectra for (pre)geometries as well as the interaction of structured spaces with the hierarchy of categorical levels and associated truncation functors that appear once one adopts \infcatst. In Subsection 2.3, we interpret the classical theory of $\cinfty$-rings of Lawvere and Dubuc in the framework of geometries. Many results here are contained in the monograph of Moerdijk-Reyes \cite{MR} and the more recent work of Joyce \cite{Joy2}; we have opted to give a different presentation to streamline our discussion of the derived theory in Section 3. We also briefly place real algebraic geometry within the general theory and explain its relation to $\cinfty$-geometry, but offer no substantive further development. Subsection 2.4 treats Lawvere theories in the higher categorical context, which furnish a more basic encoding of algebraic structure than (pre)geometries. This subsection can be regarded as supplemental to Subsection 5.5.8 of \cite{HTT}.
\subsubsection*{Section 3}
In this section, we introduce the \infcat $\sring$ of simplicial $\cinfty$-rings and study its basic properties. The \infcat $\sring$ comes equipped with the forgetful underlying algebra functor
\[ (\_)^{\rmalg}:\sring\longrightarrow\scring_{\R} \]
and we prove using results from the appendix that it behaves well under the formation of pushouts along surjections. In subsection 3.2, we endow the opposite of the subcategory of compact objects of $\sring$ with the structure of a geometry and investigate the associated spectrum functor
\[ \spec:\{\text{simplicial $\cinfty$-rings}\}\longrightarrow \{\text{simplicial $\cinfty$-ringed spaces}\} \]
adjoint to the global sections functor $\Gamma$. We show that the adjunction $(\spec\adj\Gamma)$ is a reflection on the class of \emph{Lindel\"{o}f} simplicial $\cinfty$-rings (that is, the topological space underlying $\spec\,A$ is Lindel\"{o}f, which implies the structure sheaf of $\spec\,A$ admits partitions unity), extending results of Joyce in \cite{Joy2} to the derived setting. This should be contrasted with the algebro-geometric situation, where $\spec$ is fully faithful. We dub the objects in the image of this reflection \emph{geometric}; if $A$ is Lindel\"{o}f, then $A$ is geometric if and only if the unit map \[A\longrightarrow\Gamma(\spec\,A)\]
is an equivalence. The geometric simplicial $\cinfty$-rings form the basic local building blocks of derived $\cinfty$-schemes and derived $\cinfty$-stacks we study in Part II. We construct a functor 
\[ \{\text{simplicial $\cinfty$-ringed \inftopoit}\}\longrightarrow \{\text{presentably symmetric monoidal stable \infcatst}\} \]
carrying a pair $\ofxtop$ to the \infcat $\Mod_{\Of_{\xtop}}^{\otimes}$ of $\Of_{\xtop}$-modules. Also in this subsection we construct a spectrum functor on modules $\mspec_A:\Mod_A\rightarrow \Mod_{\Of_{\speco\,A}}$ for a simplicial $\cinfty$-ring $A$. Again, contrary to the algebro-geometric situation, this functor is not fully faithful, rather it is a reflection onto $\Mod_{\Of_{\speco\,A}}$ with fully faithful right adjoint $\Gamma$ (if $A$ is Lindel\"{o}f). The essential image consists of $A$-modules we call \emph{geometric} and we formulate some criteria to establish geometricity of modules, which will come in handy in later work. In Subsection 3.3, we make contact with the theory of \emph{differential graded} $\cinfty$-rings as developed by Carchedi-Roytenberg \cite{CR2,CR1} and employed by Nuiten \cite{Nui2}. Using this theory, we intepret dg-manifolds as derived $\cinfty$-schemes and take a few steps toward proving this furnishes a fully faithful embedding.
\subsubsection*{Section 4}
This section is an assortment of mostly technical, but for the theory developed later in this work and its successors, crucial results that are quite specific to the $\cinfty$-setting. While $\cinfty$-geometry is in some respects more straightforward than algebraic geometry due to the abundance of $\cinfty$ functions -in particular partitions of unity that render the Grothendieck topologies in $\cinfty$ geometry more classical- many tautological or trivial facts of homological algebra are quite difficult to prove (or just false) for $\cinfty$-rings. We draw on classical work in singularity theory and ideal theory of $\cinfty$ functions due to Malgrange, Thom and others to study the interaction of $\cinfty$ geometry and the categorical structure of $\sring$. In particular, we will prove that the unit of the free simplicial $\cinfty$-ring monad is a flat morphism, and we show that $\cinfty$-rings of Whitney functions on closed sets (like $\cinfty$-rings of the form $\cinfty(\R^n_{\geq 0})$) have excellent formal properties in $\sring$.
\subsubsection*{Section 5}
In the last section, we generalize the theory to include derived manifolds with boundaries and corners. Following Gillam-Molcho and Joyce-Francis-Staite \cite{GilMol,Joyfra}, we declare that the class of \emph{$b$-maps} introduced by Melrose \cite{Melrose1,KottkeMelrose} is the right notion of morphism among manifolds with corners. A major insight due to Kottke-Melrose and Gillam-Molcho is that the combinatorial data of the corner structure and the calculus of $b$-maps is captured by a `positive' variant of \emph{logarithmic geometry}, in the sense of Fontaine-Illusie, Kato and Ogus \cite{Kato,Ogus}. We prove a precise version of this statement: there is a canonical equivalence of \infcatst
\[ \{\text{simplicial $\cinfty$-rings with corners}\} \simeq 
\{\text{simplicial $\cinfty$-rings $A$ with a log structure on $A_{\geq 0}$}\}.  \]
Here, $A_{\geq0}$ is the simplicial commutative monoid of \emph{positive elements} of $A$; to define it, we appeal to the results of Section 4. By virtue of this equivalence, we can place derived $\cinfty$-geometry with corners in the framework of (pre)geometries and we conclude with a discussion of the resulting notions of structured \inftopoi and spectra.
\subsubsection*{Appendix A}
Many of the results of this paper rely on an examination of the properties of the forgetful functor
\[ (\_)^{\rmalg}:\sring\longrightarrow \scring_{\R} \]
that takes a simplicial $\cinfty$-ring to its underlying simplicial commutative $\R$-algebra. This functor is well behaved because $\cinfty$ functions on $\R^n$ have the \emph{Fermat property}: for any $g:\R\rightarrow\R$, there is a \emph{difference quotient} $\delta g:\R^2\rightarrow\R$ satisfying 
\[g(x+y) =g(x) + y\delta g(x,y).\]
A Lawvere theory $\mathrm{T}$ over the theory of commutative $k$-algebras for $k$ some commutative ring satisfying the Fermat property is a \emph{Fermat theory}. As the theory of $\cinfty$-rings is a Fermat theory, one can show that if $I\subset \cinfty(\R^n)$ is any ideal, then the algebra $\cinfty(\R^n)/I$ has a canonical structure of a $\cinfty$-ring; this expresses that $(\_)^{\rmalg}$ preserves certain pushouts. In the appendix, using some ideas from \cite{dagix}, we generalize this observation and give a formula for arbitrary colimits in a Fermat theory $f^*:s\mathrm{T}\alg\rightarrow \scring_k$: if $\mathcal{J}:K\rightarrow s\mathrm{T}\alg$ is any diagram, choose for each $k\in K$ a free $\mathrm{T}$-algebra $F_k$ and a surjection $F_k\rightarrow \mathcal{J}(k)$, then the underlying commutative $k$-algebra of $\colim_K\mathcal{J}$ is a pushout
\[
\begin{tikzcd}
 \coprod_{k\in K}f^*(F_k) \ar[r] \ar[d] & f^*(\coprod_{k\in K} F_k)\ar[d] \\
 \underset{K}{\colim}f^*\mathcal{J}\ar[r] & f^*(\underset{K}{\colim}\mathcal{J}) .
\end{tikzcd}
\]
Thus for many purposes, the analysis of colimits in $s\mathrm{T}\alg$ reduces to the study of properties of the comparison map $ \coprod_{k\in K}f^*(F_k)\rightarrow f^*(\coprod_{k\in K} F_k)$ for coproducts of free $\mathrm{T}$-algebras, which are a priori the only colimits besides sifted ones we understand.
\subsection*{Acknowledgements}
\addcontentsline{toc}{subsection}{\protect\numberline{}Acknowledgements}
This paper grew out of a part of my PhD work under Damien Calaque's supervision. It is hard to express the extent of my gratitude for his guidance, constant encouragement, friendship, and his patience. Parts of this work are concerned with universal properties of derived geometries; I am thankful to Andrew MacPherson and David Carchedi for discussion on the topic, and to Dave again for the joint work we did on the subject. Special thanks to Joost Nuiten for many conversations about the contents of this paper, and derived geometry and homotopy theory in general. I am indebted to Carlos Simpson, Mauro Porta and Corina Keller for reading preliminary versions and providing valuable comments. I would also like to thank Carlos again for suggesting to write a detailed introduction to derived ($\cinfty$-)geometry motivated by moduli of solutions of PDEs.
\newline
\emph{This project has received funding from the European Research Council (ERC) under the European Union’s Horizon 2020 research and innovation programme (grant
agreement No 768679)}.
\subsection*{Notation, conventions and preliminaries}
\addcontentsline{toc}{subsection}{\protect\numberline{}Notation, conventions and preliminaries}
\begin{itemize}
    \item We handle the interplay between small and large categories via the usual device of Grothendieck universes, i.e. we assume Tarski-Grothendieck set theory. For any cardinal $\kappa$, we denote by $\mathcal{U}(\kappa)$ the collection of sets of rank $<\kappa$. We fix once and for all three strongly inaccessible cardinals $\kappa_s<\kappa_l<\kappa_{vl}$; then we call the sets in $\mathcal{U}(\kappa_s)$ \emph{small}, those in $\mathcal{U}(\kappa_l)$ \emph{large}, and those in $\mathcal{U}(\kappa_{vl})$ \emph{very large}.
    \item The ordinary category of (small) sets is denoted as $\set$. The ordinary category of (small) simplicial sets is denoted as $\sset$.
    \item An $\infty$-category or $(\infty,1)$-category is a \emph{weak Kan complex}, also known as a \emph{quasi-category}. Our reference on the foundations of such higher categories is Lurie's book \emph{Higher Topos Theory} \cite{HTT}. The large \infcat of small \infcats is denoted $\catinf$; it is the homotopy coherent nerve of the simplicial category of fibrant cofibrant objects of the simplicial model category of marked simplicial sets. The very large \infcat of large \infcats is denoted $\catinfh$. 
    \item The homotopy category of an $\infty$-category $\icat$ is denoted by $h\icat$.
     \item For $\icat$ an \infcat and $f:K\rightarrow \icat$ a diagram, we let $\icat_{/f}$ and $\icat_{f/}$ denote the slice \infcats defined by the existence of a bijection between the set of maps $S\rightarrow \icat_{/f}$ of simplicial sets and the set of maps $S\star K\rightarrow \icat$, where $\star$ is the join of simplicial sets, and $\icat_{f/}$ is defined similarly. If $K=\Delta^0$ and $f$ classifies an object $C\in \icat$, we have slices $\icat_{/C}$ and $\icat_{C/}$. The joins $S\star\Delta^0$ and $\Delta^0\star S$ are denoted $S^{\rhd}$ and $S^{\lhd}$ respectively.
    \item For $\icat$ an $\infty$-category, the Kan complex of morphisms between two objects $X$ and $Y$ is denoted by $\Hom_{\icat}(X,Y)$. We let $\Hom^{\mathrm{R}}(X,Y)$ denote the fibre $\icat_{/Y}\times_{\icat}\{X\}$ and let similarly $\Hom^{\mathrm{L}}(X,Y)$ be the fibre $\icat_{X/}\times_{\icat}\{Y\}$. Both these Kan complexes are a model for $\Hom_{\icat}(X,Y)$.
    	\item An \infcat $\icat$ is an \emph{$n$-category} if for each pair of objects $X,Y\in\icat$, the space $\Hom_{\icat}(X,Y)$ is $(n-1)$-truncated. We write $\cat_n\subset\catinf$ for the full subcategory spanned by $n$-categories, which is itself a (large) $(n+1)$-category. We will in the sequel no longer distinguish between an ordinary category $\icat$ and its nerve $\ner(\icat)$, that is, we view every category as a 1-category.
    \item For $S$ a simplicial set and $W$ a collection of edges of $S$ containing all degenerate ones, $(\mathbf{A},W)$, we denote by $\mathbf{A}[W^{-1}]$ the $\infty$-category obtained by taking a fibrant replacement of the marked simplicial set $(A,W)$ in the model category $\sset^+$ of marked simplicial sets. 
     \item For $\icat$ a simplicial set (usually an $\infty$-category) and $\icatd$ an $\infty$-category, the simplicial set of morphism from $\icat$ to $\icatd$ is denoted as $\fun(\icat,\icatd)$. It is an $\infty$-category and it is called \emph{the $\infty$-category of functors from $\icat$ to $\icatd$}. When $\icat=\spa$, the \infcat of spaces, we write $\pshv(\icatd)$ for $\fun(\icatd^{op},\spa)$, and call it the \emph{\infcat of presheaves on $\icatd$}.
    
     \item For $\icat$ an \infcatt, we will freely use the straightening-unstraightening equivalences 
     \[ \mathsf{coCart}_{\icat}\simeq\fun(\icat,\catinf),\quad  \mathsf{Cart}_{\icat}\simeq\fun(\icat^{op},\catinf) \]
     and 
      \[ \mathsf{LFib}_{\icat}\simeq\fun(\icat,\spa),\quad  \mathsf{RFib}_{\icat}\simeq\fun(\icat^{op},\spa)=\pshv(\icat) \]
   of Sections 2.3 and 3.2 of \cite{HTT}. Here $\mathsf{coCart}_{\icat}:=\ner((\sset^{+,\mathrm{cocart}})_{/\icat}^{\mathrm{fc}})$, the fibrant-cofibrant objects in the simplicial model categories of marked simplicial sets over $\icat$ with the coCartesian model structure, and $\mathsf{LFib}\subset \mathsf{coCart}_{\icat}$ is the full subcategory spanned by those $p:\icate\rightarrow\icat$ for which every edge of $\icate$ is $p$-coCartesian; $\mathsf{Cart}_{\icat}$ and $\mathsf{RFib}_{\icat}$ are defined similarly. 
    \item  We say that a right fibration $p:\icatd\rightarrow\icat$ is \emph{representable} if $\icatd$ has a final object; this implies that $\icatd\rightarrow \icat$ is equivalent to $\icat_{/p(D)}$, where $D$ is a final object of $\icatd$. It's easy to see that the representable right fibrations are exactly the objects in the essential image of the Yoneda embedding $j:\icat\rightarrow \mathsf{RFib}_{\icat}$.
\begin{enumerate}[$(1)$]
    \item Let $f:\icat\rightarrow \icatd$ be a functor, then we say that an object $\epsilon_{D}\in \icat_{/D}:=\icatd_{/D}\times_{\icatd}\icat$ depicted as a pair $(C,f(C)\rightarrow D)$ is a \emph{counit transformation at $D$} if $\epsilon_D$ is final; that is, if the right fibration $\icat_{/D}\rightarrow\icat$ is representable. We say that $f$ is a \emph{left adjoint} if there is a counit transformation for every $D\in\icatd$.
    \item Let $f:\icat\rightarrow \icatd$ be a functor, then we say that an object $\eta_{D}\in \icat_{D/}:=\icatd_{D/}\times_{\icatd}\icat$ depicted as a pair $(C,D\rightarrow f(C))$ is a \emph{unit transformation at $D$} if $\eta_D$ is initial. We say that $f$ is a \emph{right adjoint} if there is a unit transformation for every $D\in\icatd$.
\end{enumerate}
   An adjunction $L:\icat\leftrightarrows \icatd:R$, with $L$ the left adjoint and $R$ the right adjoint is written as $(L\dashv R)$.
   \item Consider a diagram $\sigma:\Delta^1\times\Delta^1\rightarrow \widehat{\mathsf{Cat}}_{\infty}$ 
  \[
\begin{tikzcd}
\icat \ar[d,"F"] \ar[r,"L"] & \icatd \ar[d,"F'"]  \\
\icat'\ar[r,"L'"] & \icatd'
\end{tikzcd}
\]
commuting up to a specified homotopy $\alpha$. Let $D\in\icatd$ an object, then we say that this diagram is \emph{$L$-right adjointable locally at $D$} (or \emph{horizontally right adjointable}) if $L$ and $L'$ admit right adjoints $U$ and $U'$ respectively, and the \emph{Beck-Chevalley transformation} 
\[F\circ U \longrightarrow   U'\circ L' \circ F\circ U \overset{\alpha}{\simeq} U'\circ F'\circ L\circ U\longrightarrow U'\circ F'   \]
is an equivalence at $D$. A square diagram $\sigma$ as above is \emph{$U$-right adjointable} if it is $U$-right adjointable locally at every $D\in\icatd$. There is an evident dual notion of a (horizontally or vertically) \emph{left adjointable square}. Using elementary manipulations of units and counits, it is easy to see that the diagram is $L$-right adjointable locally at $D$ if $L$ and $L'$ admit right adjoints and the functor $\icat\times_{\icatd}\icatd_{/D}\rightarrow  \icat'\times_{\icatd'}\icatd'_{/F'(D)}$ induced by $F$ and $F'$ preserves final objects; that is, $F$ and $F'$ carry counit transformations at $D$ to counit transformations at $F'(D)$. In particular, to show that a the square $\sigma$ is horizontally right adjointable, we observe that 
\begin{enumerate}[$(a)$]
\item if $L$ and $L'$ are fully faithful, then it suffices to check that $F'$ carries maps exhibiting coreflections in $\icat$ to coreflections in $\icatd'$.
\item if the adjoints $U$ and $U'$ are fully faithful, then it suffices to check that $F$ carries the essential image of $U$ to the essential image of $U'$. 
\end{enumerate}.
We will use these criteria frequently, together with the fact that in case $F$ and $F'$ admit left adjoints, the square $\sigma$ is horizontally right adjointable if and only if it is vertically left adjointable. 
\item We let $\prl\subset\catinf$ denote the subcategory of presentable \infcats and functors admitting right adjoints among them. Similarly, we let $\prr\subset\catinf$ denote the subcategory of presentable \infcats and functors admitting left adjoints between them. For $\kappa$ a regular cardinal, we let $\prl_{\kappa}\subset\prl$ denote the subcategory whose objects are $\kappa$-compactly generated presentable \infcats and whose morphisms are functors admitting $\kappa$-continuous right adjoints. 
\item We let $\catinf^{\mathrm{Ex}}\subset\catinf$ denote the subcategory whose objects are stable \infcats and whose morphisms are exact functors. We write $\mathsf{Pr}^{\mathrm{Ex}}\subset\catinfh$ for the subcategory whose objects are presentable stable \infcats and whose morphisms are exact functors admitting a right adjoint. We denote the \infcat of spectra by $\spect$.
\item For $\iop,\Of^{\otimes\prime}$ , we write $\alg_{\Of}(\Of')$ for the \infcat of $\iop$-algebra objects in $\Of'$. For $\icat^{\otimes}$ a symmetric monoidal \infcat, we let $\mathsf{CAlg}(\icat)$ denote the \infcat of $\einfty$-algebra objects in $\icat$. For $A$ an $\einfty$-algebra in $\icat$, we let $\Mod_A(\icat)$ denote the \infcat of $A$-modules. For $k$ a commutative ring viewed as its Eilenberg-MacLane spectrum, we write $\mathsf{CAlg}_k=\mathsf{CAlg}(\Mod_k(\spect))$ for the \infcat of $k$-modules. We write $\calg^{\geq 0}_k\subset\calg$ for the full subcategory spanned by connective $k$-algebras, and $\calg^0_k$ for the ordinary category of commutative $k$-algebras.
    \item Our grading conventions are \emph{homological}, that is, the differential on a complex \emph{lowers} the degree. Accordingly, a complex of $R$-modules $C \in \mathrm{Ch}(\mathrm{Mod}_R)$ for some commutative ring $R$ is called \emph{connective} if $H_{n}(C)=0$ for all $n\leq -1$. A complex is called \emph{eventually connective} if there exists some $n$ such that $H_k(C)=0$ for all $k<n$.
  \item We record the following \infcatst.
   \begin{enumerate}[$(a)$]
   \item The subcategory $\ltop\subset\catinfh$ whose objects are \inftopoit, and whose morphisms are functors that are left exact and admit a right adjoint. Such morphisms between \inftopoi will be called \emph{algebraic morphisms}. For $\xtop,\ytop\in\ltop$, the full subcategory of $\fun(\xtop,\ytop)$ spanned by algebraic morphisms is denoted $\fun^*(\xtop,\ytop)$.
    \item The subcategory $\rtop\subset\catinfh$ whose objects are \inftopoit, and whose morphisms are functors that admit a left exact left adjoint. Morphisms in $\rtop$ will be called \emph{geometric morphisms}. For $\xtop,\ytop\in\rtop_n$, the full subcategory of $\fun(\ytop,\xtop)$ spanned by geometric morphisms is denoted $\fun_*(\ytop,\xtop)$. The \infcats $\rtop$ and $\ltop$ are canonically antiequivalent, as are the \infcats $\fun^*(\xtop,\ytop)$ and $\fun_*(\ytop,\xtop)$.
 \end{enumerate}
    
Algebraic and geometric morphisms are usually denoted as in the adjoint pair $(f^*,f_*):\begin{tikzcd}\xtop \ar[r,shift left] & \ytop\ar[l,shift left]\end{tikzcd}$.
We will also need to work with \inftopoi relative over a given base: let $\icat$ be an \infcatt. A \emph{topos fibration} over $\icat$ is a presentable fibration $p:\widetilde{\xtop}\rightarrow \icat$ such that for each $C\in \icat$, the fibre $\widetilde{\xtop}_C$ is an \inftop and for each edge $f:\Delta^1\rightarrow \icat$, the coCartesian transformation $f_!:\widetilde{\xtop}_{f(\{0\})}\rightarrow \widetilde{\xtop}_{f(\{1\})}$ is an algebraic morphism. This is equivalent to demanding that the straightening $\mathrm{St}(p):\icat\rightarrow \catinfh$ factors through $\ltop$. We let $\ltop_{\icat}\subset \mathsf{coCart}_{\icat}$ denote the subcategory whose objects are topos fibrations and whose morphisms are commutative diagrams
\[
\begin{tikzcd}
\widetilde{\xtop}\ar[rr] \ar[dr] && \widetilde{\ytop} \ar[dl] \\
& \icat
\end{tikzcd}
\]
such that the horizontal map preserves coCartesian edges and for each $C\in\icat$, the induced map on the fibres is an algebraic morphism of \inftopoit. Let $q:\overline{\ltop}\rightarrow\ltop$ be the coCartesian fibration associated to the subcategory inclusion $\ltop\hookrightarrow\catinfh$. Then $q$ is a topos fibration. In fact $q$ is a \emph{universal} topos fibration, uniquely (up to equivalence) determined by the property that pulling back along $q$ induces, for any \infcat $\icat$, a canonical bijection between equivalence classes of topos fibrations over $\icat$ and functors $\icat\rightarrow\ltop$. 
\item Let $\icat$ be an \infcatt, then we say that \emph{sieve} on $C\in \icat$ is a subobject of $j(C)\in \pshv(\icat)$, where $j:\icat\hookrightarrow\pshv(\icat)$ denotes the Yoneda embedding. A Grothendieck topology on a small \infcat consists of a collection of sieves $\{U\hookrightarrow j(C)\}$ for each object $C\in \icat$, called \emph{covering sieves}, such that 
\begin{enumerate}[$(1)$]
    \item $j(C)\rightarrow j(C)$ is covering.
    \item If $U\rightarrow j(C)$ is covering and $D\rightarrow C$ is any map, then $U\times_{j(C)}j(D)$ is covering on $D$.
    \item If $U\rightarrow j(C)$ is a sieve and $V\rightarrow j(C)$ is a covering sieve, then if for each $j(D)\rightarrow j(C)$ that factors through $V$, $j(D)\times_{j(C)}U$ is a covering sieve on $D$, then $U$ is a covering sieve on $D$.
\end{enumerate}
Let $\tau$ be a Grothendieck topology on $\icat$, then the full subcategory $\mathsf{Shv}(\icat)\subset\pshv(\icat)$ spanned by objects that are $S$-local for $S$ the class of covering sieves $U\rightarrow j(C)$ are \emph{sheaves}. If $\icat$ is small so that $\pshv(\icat)$ is presentable, localizing at the collection of monomorphisms that are covering sieves induces a sheafification functor $L:\pshv(\icat)\rightarrow\shv(\icat)$, which is left exact, so that $\shv(\icat)$ is an \inftopt. Conversely, if $\xtop\subset\pshv(\icat)$ is a localization obtained by inverting a strongly saturated class $S$ of morphisms that is stable under pullbacks and generated by a small set of monomorphisms (so that the class $S$ is \emph{topological} in the sense of \cite{HTT}, Definition 6.1.2.4), then $\xtop$ coincides with the \inftop $\shv(\icat)$ for the Grothendieck topology given by those sieves $i:U\hookrightarrow j(C)$ such that $Li$ is an equivalence. Usually, we will specify a topology by giving a \emph{basis} (also known as a \emph{Grothendieck pretopology}); a basis $\mathcal{B}$ on $\icat$ is the following data. 
\begin{enumerate}
\item[$(*)$] For each object $C\in \icat$, a collection $\mathcal{B}(C)$ of families $\{U_{\alpha}\rightarrow C\}$ of morphisms. Such distinguished families will be called \emph{coverings}.   
\end{enumerate}
These collections are required to satisfy the following conditions.
\begin{enumerate}[$(i)$]
    \item For each $C\in \icat$, the family $\{\mathrm{id}:C\rightarrow C\}$ is a covering. 
    \item For each map $f:C'\rightarrow C$, and each covering $\{U_{\alpha}\rightarrow C\}$ of $C$, the pullbacks $U_{\alpha}\times_CC'$ exist for all $\alpha$ and the family $\{U_{\alpha}\times_{C}C'\rightarrow C'\}$ is a covering of $C'$.
    \item Let $\{U_{\alpha}\rightarrow C\}$ be a covering, and suppose we are given a covering $\{W_{\beta_{\alpha}}\rightarrow U_{\alpha} \}$ for each $\alpha$. Then the induced family $\{W_{\beta_{\alpha}}\rightarrow C\}$ is a covering.
\end{enumerate}
Let $\mathcal{B}$ be a Grothendieck pretopology on an $\infty$-category $\icat$. Consider, for each $C\in \icat$, the collection of those sieves $U\hookrightarrow j(C)$ that contain a sieve generated by some covering in $\mathcal{B}(C)$. Then this collection of sieves specifies a Grothendieck topology on $\icat$.
\item We record the following \infcatst.
\begin{enumerate}[$(a)$]
    \item The subcategory $\ltop_n\subset\catinfh$ whose objects are $n$-topoi, and whose morphisms are functors that are left exact and admit a right adjoint. Such morphisms between \inftopoi will be called \emph{algebraic morphisms}. For $\xtop,\ytop\in\ltop_n$, the full subcategory of $\fun(\xtop,\ytop)$ spanned by algebraic morphisms is denoted $\fun^*(\xtop,\ytop)$.
    \item The subcategory $\rtop_n\subset\catinfh$ whose objects are $n$-topoi, and whose morphisms are functors that admit a left exact left adjoint. Morphisms in $\rtop$ will be called \emph{geometric morphisms}. For $\xtop,\ytop\in\ltop_n$, the full subcategory of $\fun(\ytop,\xtop)$ spanned by geometric morphisms is denoted $\fun_*(\ytop,\xtop)$. The \infcats $\rtop_n$ and $\ltop_n$ are canonically antiequivalent, and the \infcats $\fun^*(\xtop,\ytop)$ and $\fun_*(\ytop,\xtop)$ are canonically equivalent.
\end{enumerate}
\cite{HTT}, theorem 6.4.1.5 in particular implies that taking full subcategories of $(n-1)$-truncated objects induces a functor
\[ \rtop \overset{\tau_{\leq (n-1)}}{\longrightarrow} \rtop_n. \]
This functor admits a fully faithful right adjoint $L_n$ that embeds the \infcat of $n$-topoi into the \infcat of \inftopoit. The essential image of this embedding is characterized as follows: let $n\in\Z_{\geq0}$. An \inftop $\xtop$ is \emph{$n$-localic} if for every \inftop $\ytop$, the canonical map
\[ \fun_*(\ytop,\xtop) \longrightarrow \fun_*(\tau_{\leq (n-1)}\ytop,\tau_{\leq (n-1)}\xtop)  \]
is an equivalence of \infcatst. We let $\ltop_{n-\loc}$ denote the \infcat of $n$-localic \inftopoit. For $\xtop$ an \inftopt, we have for every integer $n\geq -1$ a geometric morphism $\xtop\rightarrow L_n\xtop$, the \emph{$n$-localic reflection}.
\item The \infcat of $n$-localic \inftopoi is equivalent to the \infcat of $n$-topoi; in particular, the \infcat $\ltop_{n-\loc}$ is an $(n+1)$-category. The 1-category of $0$-localic \inftopoi may be identified with the category $\mathsf{Locale}$ of classical \emph{locales} (see \cite{MLM1}, Chapter IX). A locale $\mathcal{U}$ is \emph{spatial} if it has enough points; we have an equivalence of categories $\mathsf{Locale}^{\mathrm{sp}}\simeq\mathsf{Top}$ between spatial locales and \emph{sober} topological spaces; we denote the latter category by $\mathsf{Top}$, which should not cause confusion as all topological spaces we consider are sober. A $0$-localic \inftop is \emph{spatial} if its associated locale is spatial. We let $\ltop^{\mathrm{sp}}_{0-\loc}\subset\ltop_{0-\loc}$ denote the full subcategory spanned by spatial $0$-localic \inftopoit, so that we have an equivalence $\mathsf{Top}\simeq\rtop^{\mathrm{sp}}_{0-\loc}$. 
\item Let $\icat$ be a presentable \infcat and $n\geq 0$ an integer, then an object $C\in\icat$ is \emph{$n$-connective} if the $(n-1)$'th truncation $\tau_{\leq_(n-1)}C$ is a final object of $\icat$. If $\icat$ is $n$-connective for all integers $n\geq 0$, then we say that $C$ is \emph{$\infty$-connective}. A morphism $f:C\rightarrow D$ is $n$-connective if it is $n$-connective as an object of $\icat_{/D}$. An object $X$ in an \inftop $\xtop$ is \emph{hypercomplete} if $X$ is $S$-local for the class of $\infty$-connective morphisms. We denote by $L_{\mathrm{hyp}}:\xtop\rightarrow\xtop^{\mathrm{hyp}}$ the reflection onto the \inftop of hypercomplete objects of $\xtop$. We let $\ltop^{\mathrm{hc}}\subset\ltop$ denote the full subcategory spanned by hypercomplete \inftopoit; that is, those \inftopoi for which all objects are hypercomplete. The construction $\xtop\mapsto \xtop^{\mathrm{hyp}}$ determines a left adjoint to the inclusion $\ltop^{\mathrm{hc}}\subset\ltop$.
\item Let $\icat$ be a presentable \infcatt, then $\icat$ is \emph{Postnikov complete} if the tower
\[ \icat\longrightarrow\ldots\longrightarrow\tau_{\leq n}\icat\longrightarrow\tau_{\leq (n-1)}\icat\longrightarrow\ldots\longrightarrow \tau_{\leq0}\icat \]
is a limit diagram in the \infcat $\catinfh$. The \emph{Postnikov completion} of $\icat$ is the canonical functor $\icat\rightarrow\widehat{\icat}:=\lim_{n\in\Z_{\geq 0}}\tau_{\leq n}\icat$. For each integer $k\geq 0$ the Postnikov completion induces an equivalence $\tau_{\leq k}\icat\simeq \tau_{\leq k}\widehat{\icat}$. If $\xtop$ is an \inftopt, then $\widehat{X}$ is again an \inftop and the functor $\xtop\rightarrow \widehat{X}$ is an algebraic morphism (\cite{sag}, Appendix A.7.2). Let $\ltop^{\mathrm{Pc}}\subset\ltop$ denote the full subcategory spanned by Postnikov complete \inftopoi, then the construction $\xtop\mapsto \widehat{\xtop}$ determines a left adjoint to the inclusion $\ltop^{\mathrm{Pc}}\subset\ltop$.
\item An \inftop $\xtop$ is \emph{bounded} if $\xtop$ arises as a small filtered colimit of $n$-localic \inftopoi in the \infcat $\ltop$ (where $n$ is allowed to vary). Equivalently, $\xtop$ is bounded if $\xtop$ if the diagram
\[ L_0\xtop\longrightarrow\ldots\longrightarrow L_n\xtop\longrightarrow L_{n+1}\xtop\longrightarrow \ldots\longrightarrow \xtop  \]
is a colimit diagram in $\ltop$, where each transition map $L_n\xtop\rightarrow L_{n+1}\xtop$ exhibits an $n$-localic reflection (\cite{sag}, Appendix A.7.1). The colimit $\colim_{n\in \Z_{\geq 0}}L_n\xtop$ is again an \inftop and determines a localization $\xtop\rightarrow \colim_{n\in \Z_{\geq 0}}L_n\xtop$ in the \infcat $\rtop$, the \emph{bounded reflection}. Let $\ltop^\mathrm{b}\subset\ltop$ be the full subcategory spanned by bounded \inftopoit, then the construction $\xtop\mapsto \colim_{n\in \Z_{\geq 0}}L_n\xtop$ determines a right adjoint to the inclusion $\ltop^\mathrm{b}\subset\ltop$. The bounded reflection and Postnikov completion are different ways of approximating an \inftop by its truncations. As such, these constructions carry the same information: the compositions $\ltop^{\mathrm{b}}\subset\ltop\rightarrow\ltop^{\mathrm{Pc}}$ and $\ltop^{\mathrm{Pc}}\subset\ltop \rightarrow \ltop^{\mathrm{b}}$ given by $\xtop\mapsto \widehat{\xtop}$ and $\xtop\mapsto \colim_{n\in \Z_{\geq 0}}L_n\xtop$ are mutually inverse equivalences. 
    \item A \emph{manifold} is a second countable, Hausdorff topological manifold without boundary whose topological dimension is globally bounded, equipped with a maximal $C^{\infty}$-atlas. The category of manifolds is denoted $\mathsf{Man}$. A manifold in our sense may have connected components of differing dimensions, as long as there is not a (countable) sequence of connected components whose dimensions grow to infinity. An \emph{$n$-manifold} is a manifold each connected component of which has dimension $n$.
\end{itemize}
\newpage

\section{Geometries, pregeometries and Lawvere theories}

The \infcats of derived $\cinfty$-geometries, without and with corners, we construct in this paper are defined in terms of a universal property they enjoy, understood in a suitable $(\infty,2)$-categorical sense, with respect to the ordinary categories of $\cinfty$ manifolds without and with corners. Let us outline the general procedure we have in mind.
\begin{enumerate}[$(1)$]
    \item Let $\pregeo$ be a category of \emph{geometric objects}, that we suppose is formally analogous to the category of topological spaces, or smooth manifolds, or the category of nonsingular subschemes of affine $k$-schemes, in the following precise sense.
    \begin{enumerate}[$(a)$]
    \item $\pregeo$ has finite products.
    \item $\pregeo$ comes equipped with a collection of `open embeddings' $\{U_i\rightarrow X\}_i$ satisfying a few natural stability conditions; most importantly, intersections with, that is, pullbacks along, open embeddings exist and are again open embeddings. In the topological or manifold setting, it is obvious what open embeddings are. In the algebraic setting, the open embeddings of $k$-schemes correspond to localizations $A\rightarrow A[1/a]$ of their function algebras. Some variation of the class of open embeddings is possible; in the examples above, for instance, we may replace open embeddings with locally open embeddings (that is, local diffeomorphisms in differential geometry and \'{e}tale morphisms in algebraic geometry) without essentially changing the theory.
    \item $\pregeo$ is endowed with a topology that is generated by the open embeddings. 
    \end{enumerate}
    \item To a category $\pregeo$ as above, we can associate a fully faithful functor $\pregeo\hookrightarrow \geo$ that preserves finite products and pullbacks along morphisms the distinguished  class of open embeddings. The category $\geo$ has all finite limits and should be interpreted as the category of \emph{singular} geometric objects obtained by taking arbitrary intersections of the objects of $\pregeo$. In fact, we have for each natural number $n$, an $(n+1)$-category $\geo_{\leq n}$ equipped with a functor $\pregeo\hookrightarrow\geo_{\leq n}$ that completes $\pregeo$ by finite limits in the $(n+1)$-categorical sense. We view the case $n=\infty$ as the \emph{derived geometry} associated to $\pregeo$. 
    \item Let $\geo$ be an \infcat admitting finite limits, then it makes sense to consider for a topological space $X$ \emph{sheaves of $\geo$-structures} on $X$, which gives us an \infcat $\mathsf{Top}(\geo)$ of pairs $(X,\Of_X)$. Suppose that $\geo$ is obtained from a category $\pregeo$ of geometric objects via the procedure sketched above then, as it turns out, the class of open embeddings on $\pregeo$ provides the means to talk about sheaves of local $\geo$-structures and local morphisms among such (in the same sense that ringed spaces may be \emph{locally} ringed). On general grounds, there exists a spectrum functor 
 \[  \spec:\geo\longrightarrow\mathsf{Top}^{\mathrm{loc}}(\geo)\subset \mathsf{Top}(\geo)  \]
 taking values in the subcategory of spaces equipped with sheaves of local $\geo$-structures and local morphisms among them, which realizes the objects in the abstractly defined \infcat $\geo$ as \emph{structured spaces}, which we declare are bona fide geometric objects. Following algebro-geometric paradigms, we think of the essential image of $\spec$ as \emph{affine $\geo$-schemes}; it is then a straightforward matter to glue affine $\geo$-schemes along equivalences in the \infcat $\mathsf{Top}(\geo)$. Such gluings correspond to pairs $(X,\Of_X)\in \mathsf{Top}(\geo)$ for which $X$ admits an open cover $\{U_i\subset X\}_i$ such that for all $i$, $(U_i,\Of|_{U_i})$ is an affine $\geo$-scheme. This method of gluing affines is equivalent to, but more concrete than, the gluing of affine Kuranishi models as simplicial sheaves in the introductory section.  
\end{enumerate}
\begin{rmk}
In fact, we will work more generally: the spectrum functor is naturally defined on the pro-completion of $\geo$ and takes values in an \infcat not of \emph{spaces} equipped with $\geo$-structures, but of \emph{higher topoi} equipped with $\geo$-structures. The generalization to topoi is not generalizing for the sake of it; rather many moduli spaces of interest (in symplectic topology for instance) are naturally structured topoi instead of structured topological spaces. 
\end{rmk}
Subsections 2.1 and 2.2 contain no original work and are devoted to the recollection Lurie's theory of pregeometries and geometries introduced in \cite{dagv} (parts of which appear in somewhat different form in \cite{sag}), which renders precise the ideas just discussed. Subsection 2.3 treats $\cinfty$-rings in the setting of the general theory outlined below. In the last subsection we prove some results about projectively generated \infcats that we will need later.

\subsection{Basic notions}
The following notion stands at the basis of the theory we will develop.
\begin{defn}[Lurie \cite{dagv}]\label{admissibility}
Let $\mathcal{T}$ be an $\infty$-category. Let $\pregeo^{\mathrm{ad}}\subset\pregeo$ be a subcategory of $\pregeo$, morphisms of which will be called \emph{admissible}. The subcategory $\pregeo^{\mathrm{ad}}$ is an \emph{admissibility structure} on $\pregeo$ if the inclusion $\pregeo^{\mathrm{ad}}\subset\pregeo$ is a categorical fibration, that is, $\pregeo^{\mathrm{ad}}$ is a replete subcategory, and the following conditions are satisfied.
\begin{enumerate}[(a)]
\item For every admissible map $f:U\rightarrow X$ and any map $g:Y\rightarrow X$, there is a pullback 
\begin{center}
\begin{tikzcd}
V \ar[d,"f'"'] \ar[r,"g'"] & U\ar[d,"f"] \\
Y \ar[r,"g"] & X
\end{tikzcd}    
\end{center}
in $\pregeo$ with $f'$ an admissible map.
\item For a commutative diagram 
\begin{center}
\begin{tikzcd}
X \ar[dr,"f"'] \ar[rr,"h"] && Y \ar[dl,"g"]\\
&  Z 
\end{tikzcd}    
\end{center}
with $g$ admissible, $f$ is admissible if and only if $h$ is admissible.   
\item A retract of an admissible map is admissible. 
\end{enumerate}
Suppose that we are also given a Grothendieck topology $\tau$ on $\pregeo$, that is, a specification of covering sieves $\pregeo^{S}_{/X}\subset \pregeo^{S}_{/X}$. We say that the admissibility structure $\pregeo^{\mathrm{ad}}$ \emph{is compatible with the topology $\tau$} (or that the \emph{topology is compatible with the admissibility structure}) if the following condition is satisfied.
\begin{enumerate}
    \item[$(*)$] Every covering sieve $\pregeo^{S}_{/X}\subset \pregeo^{S}_{/X}$ contains a covering sieve generated by admissible morphisms. 
\end{enumerate}
We say that a family $\{U_i\rightarrow X\}_i$ on $X$ in $\pregeo$ that generates a covering sieve and consists of admissible morphism is an \emph{admissible covering}.
\end{defn}
\begin{lem}
Let $\pregeo$ be a pregeometry equipped with a topology $\tau$ and an admissibility structure $\pregeo^{\mathrm{ad}}$. Then the following are equivalent.
\begin{enumerate}[$(1)$]
    \item The admissibility structure is compatible with the topology.
    \item The topology $\tau$ has a basis consisting of admissible coverings. 

\end{enumerate}
\end{lem}
\begin{proof}
Note that $(2)\Rightarrow (1)$ is immediate since the covering sieves for a topology specified by a basis are precisely those sieves that contain a sieve generated by a covering family. For $(1)\Rightarrow (2)$, it suffices to show that the admissible coverings form a basis for a topology on $\pregeo$. Since equivalences are admissible, the families $\{\mathrm{id}:X\rightarrow X\}$ are admissible coverings. Let $\{U_i\rightarrow X\}_i$ be an admissible covering generating a covering sieve $U\rightarrow j(X)$ (viewed as a subobject of the representable sheaf associated to $X$). Then for any $Y\rightarrow X$, the family $\{U_i\times_YX\rightarrow Y\}_i$ exists and consists of admissibles, by definition of an admissibility structure. This family on $Y$ generates the pullback sieve $U\times_{j(X)}j(Y)$, which is therefore a covering sieve. Suppose we are given an admissible covering $\{U_i\rightarrow X\}$ and an admissible covering $\{W_{ij}\rightarrow U_i\}_{ij}$ on each $U_i$. Let $W_i\rightarrow j(U_i)$ be the covering sieve generated by $\{W_{ij}\rightarrow U_i\}_{j}$ and let $W\rightarrow j(X)$ be the sieve generated by the family $\{W_{ij}\rightarrow X\}_{ij}$, then we have a commuting diagram 
\[
\begin{tikzcd}
W_i\ar[d]\ar[r] & j(U_i)\ar[d] \\
W\ar[r] & j(X)
\end{tikzcd}
\]
so the sieve $W\times_{j(X)}j(U_i)$ on $U_i$ is covering as it contains $W_i$. Since the family $\{U_i\rightarrow X\}_i$ is covering, we deduce that the sieve $W$ is covering so that $\{W_{ij}\rightarrow X\}_{ij}$ is an (admissible) covering.
\end{proof}
\begin{rmk}
In view of the preceding lemma, we will usually define a compatible topology on an \infcat $(\pregeo,\pregeo^{\mathrm{ad}})$ equipped with an admissibility structure by giving a basis of families that contain only admissible maps. The lemma also shows that a compatible topology on $\pregeo$ is completely determined by its restriction to $\pregeo^{\mathrm{ad}}$. Conversely, the topologies on $\pregeo^{\mathrm{ad}}$ that extend to a topology on $\pregeo$ compatible with the admissibility structure are those that admit a basis $\mathcal{B}$ satisfying the following condition.
\begin{enumerate}
\item[$(*)$] For any covering family $\{U_i\rightarrow X\}_i$ in $\mathcal{B}$ and any map $Y\rightarrow  X$ in $\pregeo$ \emph{not necessarily admissible}, the family $\{U_i\times_XY\rightarrow Y\}_i$, which consists of admissible morphisms by $(a)$ of Definition \ref{admissibility}, lies in $\mathcal{B}$.
\end{enumerate}
Indeed, a basis $\mathcal{B}$ of a topology on $\pregeo^{\mathrm{ad}}$ viewed as a collection of families of morphisms of $\pregeo$ is a basis for a topology on $\pregeo$ if and only if it is stable under pullbacks in $\pregeo$, that is, if it satisfies the condition $(*)$. 
\end{rmk}
\begin{defn}\label{pregeometry}
\begin{enumerate}[(1)]
    \item A \emph{pregeometry} is a pair $(\mathcal{T},\mathcal{T}^{\mathrm{ad}})$ of an essentially small $\infty$-category $\mathcal{T}$ that admits finite products, together with a topology and an admissibility structure $\mathcal{T}^{\mathrm{ad}}$ on $\mathcal{T}$ compatible with the topology.
    \item A \emph{geometry} is a pair $(\mathcal{G},\mathcal{G}^{\mathrm{ad}})$ of an essentially small idempotent complete $\infty$-category $\mathcal{G}$ that admits finite limits together with an admissibility structure $\mathcal{G}^{\mathrm{ad}}$ on $\mathcal{G}$ compatible with the topology.
\end{enumerate}
\end{defn}
We will usually just write $\mathcal{T}$ (or $\mathcal{G}$) for a pregeometry $(\mathcal{T},\mathcal{T}^{\mathrm{ad}})$ (or a geometry $(\mathcal{G},\mathcal{G}^{\mathrm{ad}})$). 
\begin{ex}[Discrete (pre)geometries]\label{ex:discretepregeo}
Let $\mathcal{T}$ be an $\infty$-category with finite products. We make $\mathcal{T}$ into a pregeometry by declaring that only equivalences are admissible, which generates the trivial Grothendieck topology. Pregeometries for which the subcategory of admissible maps is the subcategory spanned by equivalences are called \emph{discrete pregeometries}. Discrete pregeometries are the same thing as \emph{Lawvere theories}: essentially small \infcats that admit finite products, which we study in more detail in Section \ref{sec:lawvere}. A basic example of a discrete pregeometry is the following: let $k$ be a commutative ring, and let $\pregeo^{\mathrm{disc}}_k:=\mathrm{Poly}_k$, where $\mathrm{Poly}_k$ is the (ordinary) category defined as follows.
\begin{enumerate}
    \item[$(O)$] Objects are the affine $k$-spaces $\mathbb{A}^n_k$ for $n\geq0$.
    \item[$(M)$] Morphisms $f:\mathbb{A}^n\rightarrow \mathbb{A}^m$ are polynomial maps.
\end{enumerate}
Of course, starting from an $\infty$-category $\mathcal{G}$ that has finite limits and is idempotent complete, we have an associated \emph{discrete geometry} with underlying $\infty$-category $\mathcal{G}$.
\end{ex}
\begin{ex}
Let $\geo$ be a category that admits pullbacks, then setting $\geo^{\mathrm{ad}}=\geo$ yields an admissibility structure on $\geo$. Every Grothendieck topology is compatible with this admissibility structure.     
\end{ex}
\begin{ex}
Let $\geo$ be a geometry, then for each object $X\in \geo$, the overcategory $\geo_{/X}$ admits finite limits and is endowed with a natural Grothendieck topology: a sieve $(\geo_{/X})^S_{/Y}$ if and only if $(\geo_{/X})^S_{/Y}\subset \geo_{/Y}$ is a sieve. Declaring a morphism $f:Y\rightarrow Z$ over $X$ to be admissible if and only if $f$ is admissible in $\geo$ determines an admissibility structure on $\geo_{/X}$ compatible with the topology.  
\end{ex}
From Example \ref{ex:discretepregeo}, we see that a discrete pregeometry is the same thing as a Lawvere theory. As we indicated in the introduction, a Lawvere theory T has an associated category of algebras: these are functors $\mathrm{T}\rightarrow \set$ that preserve finite products. For the example $\mathrm{T}=\mathrm{Poly}_k$ above, these algebras are precisely commutative $k$-algebras; To see this, we consider a functor $F:\mathrm{Poly}_k\rightarrow \set$ that preserves finite products, which determines a set $A:=F(\mathbb{A}^1)$. This set has addition and multiplication operations defined by
\[ A^2\cong F(\mathbb{A}\times \mathbb{A}) \overset{F(+)}{\longrightarrow}   F(\mathbb{A}) =A, \quad \quad \quad A^2\cong F(\mathbb{A})\times F(\mathbb{A}) \overset{F(\cdot)}{\longrightarrow}   F(\mathbb{A}) =A,\]
where $+$ and $\cdot$ are the addition and multiplication maps $+,\cdot:\mathbb{A}^2\rightarrow\mathbb{A}$. Multiplication and addition are commutative and associative operations; we have, for instance, commuting diagrams
\[
\begin{tikzcd}
\mathbb{A}^3 \ar[d,"\cdot\times\mathrm{id}"']\ar[r,"\mathrm{id}\times \cdot"] & \mathbb{A}^2\ar[d,"\cdot"] \\
\mathbb{A}^2\ar[r,"\cdot"] & \mathbb{A},
\end{tikzcd}
\quad \quad \begin{tikzcd}
\mathbb{A}^3 \ar[d,"+\times\mathrm{id}"']\ar[r,"\mathrm{id}\times +"] & \mathbb{A}^2\ar[d,"+"] \\
\mathbb{A}^2\ar[r,"+"] & \mathbb{A},
\end{tikzcd}
\]
expressing associativity of multiplication and addition. By functoriality, these diagrams are carried to commuting diagrams
\[
\begin{tikzcd} A^3\ar[d,"F(\cdot)\times\mathrm{id}"']\ar[r,"\mathrm{id}\times F(\cdot)"] & A^2\ar[d,"F(\cdot)"] \\
A^2\ar[r,"F(\cdot)"] & A
\end{tikzcd}
\quad \quad \begin{tikzcd} A^3\ar[d,"F(+)\times\mathrm{id}"']\ar[r,"\mathrm{id}\times F(+)"] & A^2\ar[d,"F(+)"] \\
A^2\ar[r,"F(+)"] & A
\end{tikzcd}
\]
expressing that the multiplication and addition operations on $A$ are associative. The commutativity, unit and distributive law are similarly a consequence of the fact that $F$ is a functor. If we replace $\set$ with the \infcat of spaces $\spa$, we obtain a notion of a simplicial commutative $k$-algebra, which is one version of the notion of $k$-algebra \emph{up to coherent homotopy}. Much more variation is possible: we can also replace $\set$ with the category $\shv_{\set}(X)$ of sheaves (of sets) on a topological space $X$ in which case we obtain the category of commutative $k$-algebras \emph{internal to $\shv_{\set}(X)$}, that is, a sheaf of $k$-algebras on $X$, or with the \infcat $\shv(X)$ of sheaves (of spaces) on $X$, in which case we obtain the \infcat of sheaves of simplicial commutative $k$-algebras on $X$. In full generality, we might be interested in finite product preserving functors
\[ \mathrm{T}\longrightarrow \icat \]
valued in arbitrary \infcatst, which we should think of as $\mathrm{T}$-algebras internal to $\icat$. If we view pregeometries and geometries as a generalization of Lawvere theories, we can also think of functors 
\[\pregeo/\geo \longrightarrow \icat\]
that preserve all the relevant structure as some kind of algebras internal to $\icat$, with possibly very intricate multiplication rules determined by the (pre)geometry.
\begin{defn}\label{pregeomstructure}
Let $\mathcal{T}$ be an \infcat equipped with an admissibility structure $\pregeo^{\mathrm{ad}}$ and let $\icat$ be an $\infty$-category. We denote by $\fun^{\mathrm{ad}}(\mathcal{T},\icat)$ the full subcategory of $\fun(\mathcal{T},\icat)$ spanned by those functors $\Of:\mathcal{T}\rightarrow \icat$ which preserves pullbacks along admissible maps, that is, if 
    \[
    \begin{tikzcd}
    V \ar[d,"f'"'] \ar[r,"g'"] & U\ar[d,"f"] \\
    Y \ar[r,"g"] & X
    \end{tikzcd}
    \] 
    is a pullback square in $\pregeo$ where $f$ and thus $f'$ are admissible, then the diagram 
     \[
    \begin{tikzcd}
    \Of(V) \ar[d,"f'"'] \ar[r,"g'"] & \Of(U)\ar[d,"f"] \\
    \Of(Y) \ar[r,"g"] & \Of(X)
    \end{tikzcd}
    \] 
    is also a pullback square.\\
    Suppose that $\pregeo$ is a pregeometry, then we let $\fun^{\pi\mathrm{ad}}(\pregeo,\icat)\subset\fun^{\mathrm{ad}}(\pregeo,\icat)$ denote the full subcategory spanned by functors that also preserve finite products. 
\end{defn}
In the theory of operads, a map of operads $\iop\rightarrow \iopv$ is an $\iop$-algebra in $\iopv$. For (pre)geometries, the same principle applies.
\begin{defn}
\begin{enumerate}[(1)]
\item Let $\pregeo$ and $\pregeo'$ be pregeometries. A \emph{transformation of pregeometries} is a functor $f\in \fun^{\pi\mathrm{ad}}(\pregeo,\pregeo')$ such that $f(\pregeo^{\mathrm{ad}})\subset (\pregeo')^{\mathrm{ad}}$, and $f$ takes admissible coverings to admissible coverings. 
We let $\fun_{\tau}^{\pi\mathrm{ad}}(\pregeo,\pregeo')\subset \fun^{\pi\mathrm{ad}}(\pregeo,\pregeo')$ denote the full subcategory spanned by transformations of pregeometries.
\item Let $\geo$ and $\geo'$ be geometries. A \emph{transformation of geometries} is a functor $f\in \fun^{\mathrm{lex}}(\geo,\geo')$ such that $f(\geo^{\mathrm{ad}})\subset (\geo')^{\mathrm{ad}}$, and $f$ takes admissible coverings to admissible coverings. We let $\fun^{\mathrm{lex}}_{\tau}(\geo,\geo')\subset \fun^{\mathrm{lex}}(\geo,\geo')$ be the full subcategory spanned by transformations of geometries.
\end{enumerate}
We let $\mathsf{PGeo}$ denote the \infcat whose objects are pregeometries and whose morphisms are transformations of pregeometries. More precisely, we consider the fibrant simplicial category $\mathsf{PGeo}_{\simp}$ defined as follows.
\begin{enumerate}
    \item[$(O)$] Objects are triples $(\pregeo,\tau,\pregeo^{\mathrm{ad}})$ of an \infcat $\pregeo$ that admits finite products, a Grothendieck topology $\tau$ on $\pregeo$ and an admissibility structure on $\pregeo$ compatible with the topology. 
    \item[$(M)$] The Kan complex of morphisms between $\pregeo$ and $\pregeo'$ (we suppress the other data from the notation) is the largest Kan complex of the full subcategory of $\fun^{\pi\mathrm{ad}}(\pregeo,\pregeo')$ spanned by functors $f:\pregeo\rightarrow\pregeo'$ such that $f(\pregeo^{\mathrm{ad}})\subset (\pregeo')^{\mathrm{ad}}$ and that carry admissible coverings to admissible coverings.
\end{enumerate}
We let $\mathsf{PGeo}:=\ner(\mathsf{PGeo}_{\simp})$ be the coherent nerve of the fibrant simplicial category we just defined. Similarly, we have an \infcat $\mathsf{Geo}$ of geometries.
\end{defn}
\begin{rmk}
The \infcats just constructed can be promoted to $(\infty,2)$-categories. We do not intend to use $(\infty,2)$-category theory in a serious way in this work, but we will invoke the formalism sporadically in this section to put some constructions in context. We will in particular perform a variety of \emph{completions}, by adjoining limits, colimits or equivalences in a controlled manner, which are best understood as $(\infty,2)$-functors. For our purposes, an $(\infty,2)$-category is a fibrant $\sset^+$-enriched category, that is, a category enriched in marked simplicial sets such that the marked simplicial morphism sets are \infcats with their equivalences marked. \\
We let $\mathsf{PGEO}$ denote the $(\infty,2)$-category of pregeometries. More precisely, we consider the $\sset^{+}$-enriched category defined as follows.
\begin{enumerate}
    \item[$(O)$] Objects are triples $(\pregeo,\tau,\pregeo^{\mathrm{ad}})$ of an \infcat $\pregeo$ that admits finite products, a Grothendieck topology $\tau$ on $\pregeo$ and an admissibility structure on $\pregeo$ compatible with the topology. 
    \item[$(M)$] The marked \infcat of morphisms between $\pregeo$ and $\pregeo'$ (we suppress the other data from the notation) is the the full subcategory of $\fun^{\pi\mathrm{ad}}(\pregeo,\pregeo')$ spanned by functors $f:\pregeo\rightarrow\pregeo'$ such that $f(\pregeo^{\mathrm{ad}})\subset (\pregeo')^{\mathrm{ad}}$ and that carry admissible coverings to admissible coverings, with equivalences marked.
\end{enumerate}
Similarly, we have an $(\infty,2)$-category $\mathsf{GEO}$ of geometries.
\end{rmk}
\begin{rmk}\label{rmk:geometrypullback}
We have a faithful functor $\mathsf{Forget}:\mathsf{Geo}\rightarrow\catinf^{\mathrm{lex},\mathrm{Idem}}$ to the \infcat of \infcats admitting finite limits and retracts for idempotents and left exact functors among them, forgetting the admissibility structure and the topology. This functor is a Cartesian fibration: if $f:\geo\rightarrow \geo'$ is a left exact functor, $\tau$ a Grothendieck topology and $\geo^{\prime\mathrm{ad}}$ an admissibility structure on $\geo$ compatible with the topology, then 
\begin{enumerate}[$(a)$]
\item declaring the families $\{U_i\rightarrow X\}_i$ in $\geo$ that are carried to families that generate covering sieves for $\tau$ by $f$ to be covering families determines a basis for a Grothendieck topology on $\geo$.
\item the subcategory $f^{-1}(\geo^{\prime})$ of morphisms $f:U\rightarrow X$ that are carried to admissible morphisms in $\geo$ determines an admissibility structure compatible with the the topology defined in $(a)$. 
\end{enumerate}  
For this geometry structure, the functor $f$ is a transformation of geometries. A left exact functor $\geo''\rightarrow \geo$ is a transformation of geometries for the geometry structure on $\geo$ constructed in $(a)$ and $(b)$ if and only if the composition $\geo''\rightarrow\geo'$ is a transformation of geometries, which is a reformulation of the assertion that $\geo\rightarrow \geo'$ is a Cartesian edge. Since the \infcat of geometry structures on $\geo$ has an initial object in the form of the discrete geometry structure, we deduce the existence of a fully faithful left adjoint $\catinf^{\mathrm{lex},\mathrm{Idem}}\rightarrow \mathsf{Geo}$ to $\mathsf{Forget}$ that carries $\geo$ to $\geo_{\mathrm{disc}}$. The construction of the Cartesian edge requires that $f$ at least preserves pullbacks, so the analogous construction for pregeometries does not work.  
\end{rmk}
\begin{rmk}\label{rmk:geometrypushforward}
The functor $\mathsf{Forget}:\mathsf{Geo}\rightarrow \catinf^{\mathrm{lex},\mathrm{Idem}}$ is also a coCartesian fibration: if $f:\geo\rightarrow \geo'$ is a left exact functor, $\tau$ a Grothendieck topology and $\geo^{\mathrm{ad}}$ an admissibility structure on $\geo$ compatible with the topology, then 
\begin{enumerate}[$(a)$]
    \item the intersection of all admissibility structures on $\geo'$ containing $f(\geo^{\mathrm{ad}})$ is an admissibility structure.
    \item the intersection of all Grothendieck topologies on $\geo'$ compatible with the admissibility structure of $(a)$ and containing all sieves generated by the images of admissible coverings under $f$, is a Grothendieck topology compatible with the admissibility structure of $(a)$.
\end{enumerate}
For this geometry structure, the functor $f$ is a transformation of geometries. A left exact functor $\geo'\rightarrow \geo''$ to a geometry $\geo''$ is a transformation of geometries if and only if the composition $\geo\rightarrow \geo''$ is, which is a reformulation of the assertion that $\geo\rightarrow \geo''$ is a coCartesian edge. Since the \infcat of geometry structures on $\geo'$ has a final object in the form of the total geometry structure $\geo^{\mathrm{ad}}=\geo$, equipped with the maximal topology for which all sieves are covering, we deduce the existence of a fully faithful right adjoint $\catinf^{\mathrm{lex},\mathrm{Idem}}\rightarrow \mathsf{Geo}$ to $\mathsf{Forget}$. Again, the analogous assertion is false for pregeometries.
\end{rmk}
\subsubsection{Examples of (pre)geometries}
The first nontrivial examples of pregeometries we give are the motivating ones of this work.
\begin{ex}[Pregeometry of Smooth Manifolds]\label{pregeomsm}
Let $\diff$ be the pregeometry whose underlying $\infty$-category is $\mathsf{Man}$. A morphism $f:U\rightarrow M$ is admissible if it is an injective local diffeomorphism, that is, an open embedding. The usual topology on $\diff$ is compatible with the admissibility structure: let us say that a family of open embeddings $\{U_i\rightarrow M\}_i$ is an \emph{admissible covering} if this family is jointly surjective. The admissible coverings form a basis for the standard \emph{\'{e}tale topology} on $\diff$. Similarly, we let $\diff^{\mathrm{open}}$ be the pregeometry whose underlying $\infty$-category is the nerve of the category of open submanifolds of $\R^n$, for some $n\geq 0$, endowed with the \'{e}tale topology. Admissible morphisms are again open embeddings. There is an obvious transformation of pregeometries $\diff^{\mathrm{open}}\rightarrow \diff$.
\end{ex}
We will also deal with manifolds with corners. 
\begin{ex}[Pregeometry of Smooth Manifolds with Corners]\label{pregeomsmcorners}
Let $\mathsf{Man}_c$ be the category whose objects are manifolds with corners, and whose morphisms are the \emph{$b$-maps} of Melrose \cite{Melrose,Melrose1,Joyfra}. A map $f:M\rightarrow N$ between manifolds with corners is locally of the form 
$\R^n\times\R^k_{\geq 0}\subset U\rightarrow V\subset \R^m\times\R^l_{\geq0}$ for $U$ and $V$ open subsets, so we may replace $M$ and $N$ by these opens. Then $f$ is a \emph{$b$-map} if $f$ is (Seeley) smooth (i.e. there is an extension $\tilde{f}$ of $f$ to some open neighbourhood $U\subset \tilde{U}\subset\R^n\times\R^k$ such that $\tilde{f}$ is smooth) and either of the following two conditions hold.
\begin{enumerate}[$(1)$]
    \item $f$ maps $\R^n\times\R_{\geq 0}^k$ into $\R^m\times\{0\}$
    \item Write $f=(f_1,\ldots,f_n,f_{n+1,\ldots,f_{m+l}})$, then each $f_{n+i}$ decomposes uniquely as $g_{n+i}\prod h_1^{\alpha_1}\ldots h_k^{\alpha_k}$ where $g_{n+i}\in\cinfty(U,\R_{>0})$, the $\alpha_j \in \Z_{\geq0}$ are nonnegative integers and the $\{h_j\}$ form a complete set of boundary defining functions (which we can take to be coordinate functions). 
\end{enumerate}
If a $b$-map $f$ does \emph{not} satisfy $(1)$, we say that $f$ is an \emph{interior $b$-map}. For $M$ a manifold with corners, we let 
\[\cinfty_b(M)\subset\cinfty_{\geq 0}(M):=\cinfty(M,\R_{\geq0})\]
denote the subset of interior $b$-maps. For $p\in U\subset \R^{n}\times\R^k_{\geq 0}$, we say that \emph{$p$ lies in the codimension $i$ corner} for $1\leq i\leq k$ if $i$ out of the last $k$ coordinates of $p$ vanish. This is a property of points of open subsets of $\R^{n}\times\R^k_{\geq 0}$ stable under diffeomorphisms of manifolds with corners. Let $M$ a manifold with corners, then define the subsets
\[\del^i M = \{p\in M;\, \text{in some (hence all) coordinate charts }p\text{ lies in the codimension } i\text{ corner}\}. \]
A \emph{boundary component} of $M$ is the closure of a connected component of $\del^1M$. We let $H^1(M)$ denote the set of boundary components. A manifold with corners $M$ is a \emph{manifold with faces} if either one of the following equivalent conditions are satisfied.
\begin{enumerate}[$(1)$]
    \item For each boundary component $F\in H^1(M)$, the map $F\rightarrow M$ is an injection. 
    \item For each boundary component $F\in H^1(M)$, there exists a boundary defining function $\rho_F\in \cinfty_b(M)$; that is, around each point of $p\in F$ there are local coordinates with $\rho_F$ as last element.
\end{enumerate}
We define a pregeometry $\diffc$ as the nerve of the category of manifolds with faces and interior $b$-maps among them, endowed with the \'{e}tale topology which is generated by jointly surjective families of open embeddings. A map $f:U\rightarrow M$ is admissible if it is an open embedding. We also have a pregeometry $\diffc^{\mathrm{open}}\subset\diffc$, the full subcategory spanned by open subspaces of the model spaces $\{\R^{n}\times\R^{k}_{\geq 0}\}$
\end{ex}
Here is a variant of the pregeometry of $\cinfty$ manifolds with a much coarser topology, that will be of little interest to us. 
\begin{ex}[Smooth Manifolds, finitary \'{e}tale topology]
Let $\diff^{\mathrm{fin}}$ be the pregeometry whose underlying \infcat is $\mathsf{Man}$ equipped with the same admissibility structure as $\diff$. We endow $\diff^{\mathrm{fin}}$ with the following topology compatible with the admissibility structure: a sieve $(\diff^{\mathrm{fin}})^S_{/M}$ on a manifold $M$ is covering if it contains a sieve generated by a finite collection $\{U_i\rightarrow M\}_{i\in I}$ of admissible maps that are jointly surjective. The identity functor on the category of manifolds determines a transformation of pregeometries $\diff^{\mathrm{fin}}\rightarrow\diff$.
\end{ex}

\begin{ex}[Pregeometry of Complex Manifolds]\label{pregeomcplx}
The starting point of \emph{derived analytic geometry} as developed in the final sections of \cite{dagix} and in \cite{P1,PY1} is the \emph{complex analytic pregeometry} $\pregeo_{\mathrm{An}_{\C}}$. Here, the underlying $\infty$-category is the nerve of the category of open submanifolds of $\C^n$ for some $n\in \N$, and a morphism is admissible if it is an injective local biholomorphism.  
\end{ex}

The following pregeometries are among the main players in the passage from classical algebraic geometry to derived algebraic geometry. We will explain later on how these pregeometries `generate' the geometries wherein derived algebraic geometry takes place.
\begin{ex}[Pregeometry of Zariski open subschemes of affine $k$-space]\label{pregeomalgzar}
Let $k$ be a commutative ring, and let $\pregeo_{\mathrm{Zar}}(k):=(\calg^{0\mathrm{Zar}}_k)^{op}$, where $\mathsf{CAlg}^{0\mathrm{Zar}}_k$ is the (ordinary) category of $k$-algebras of the form $k[x_1,\ldots,x_n,(f(x_1,\ldots,x_n))^{-1}]$, where $f$ is a polynomial function on affine $n$-dimensional $k$-space $\mathbb{A}^n_k$. Given an object $A\in \mathsf{CAlg}^{0\mathrm{Zar}}_k$, we denote the corresponding object in $\pregeo_{\mathrm{Zar}}(k)$ by $\mathrm{Spec}\,A$ (for the moment, this is just abstract notation, not meant to indicate that $\mathrm{Spec}\,A$ is a locally ringed space). Recall that, given a commutative $k$-algebra and any element $b\in B$, the \emph{localization} of $B$ by $b$ is the universal object (defined up to isomorphism) $f:B\rightarrow B[1/b]$ such that $f(b)$ is invertible, and should be thought of as the algebra of functions on the open set where $b$ is nonzero. We can give the $\infty$-category $\pregeo_{\mathrm{Zar}}(k)$ the structure of a pregeometry by endowing it with the following admissibility structure and compatible topology:
\begin{enumerate}[(1)]
    \item A morphism $\mathrm{Spec}\,A\rightarrow \mathrm{Spec}\,B$ is admissible if and only if there exists some element $b\in B$ such that the map $B\rightarrow A$ induces an isomorphism $B[1/b]\cong A$. 
    \item A collection of admissible morphism $\{\mathrm{Spec}\,B_i[1/b_i]\rightarrow \mathrm{Spec}\,B\}_i$ is an admissible covering if and only if the elements $\{b_i\}$ generate the unit ideal in $B$.  
\end{enumerate}
\end{ex}

\begin{ex}[Pregeometry of affine schemes \'{e}tale over affine $k$-space]\label{pregeomalget}
Let $k$ be a commutative ring, and let $\pregeo_{\text{\'{e}t}}(k):=(\mathsf{CAlg}^{0sm}_k)^{op}$, where $\mathsf{CAlg}^{0sm}_k$ is the (ordinary) category of $k$-algebras $A$ that admit an \'{e}tale map $f:k[x_1,\ldots,x_n]\rightarrow A$ (that is, $f$ is finitely presented, flat, and the module of relative K\"{a}hler differential $\Omega_f$ vanishes). We can give the $\infty$-category $\pregeo_{\text{\'{e}t}}(k)$ the structure of a pregeometry by endowing it with the following admissibility structure and compatible topology:
\begin{enumerate}[(1)]
    \item A morphism $\mathrm{Spec}\,A\rightarrow \mathrm{Spec}\,B$ is admissible if and only if the map $B\rightarrow A$ is \'{e}tale. 
    \item A collection of admissible morphism $\{\mathrm{Spec}\,B_i\rightarrow \mathrm{Spec}\,B\}_i$ is an admissible covering if and only if there exists a finite set of indices $\{i_j\}_{1\leq j\leq n}$ such that the induced map $g:B\rightarrow \prod_{1\leq j\leq n}B_{i_j}$ is faithfully flat (that is, the base change functor along $g$ preserves and reflects exact sequences of $B$-modules).
\end{enumerate}
\end{ex}

The following two examples of \emph{geometries} describe the arena of \emph{classical} algebraic geometry.

\begin{ex}[Geometry of affine $k$-schemes (Zariski)]\label{geomalgzar}
Let $k$ be a commutative ring, and let $\geo_{\mathrm{Zar}}(k):=(\mathsf{CAlg}^0_k)_{\mathrm{fp}}$, where $(\mathsf{CAlg}^0_k)_{\mathrm{fp}}$ is the (ordinary) category of finitely presented $k$-algebras; that is $k$-algebras of the form $k[x_1,\ldots,x_n]/I$ for some finitely generated ideal $I$. We give $\geo_{\mathrm{Zar}}(k)$ the structure of a geometry by endowing it with the following admissibility structure and compatible topology, which is the obvious extension of the admissibility structure and topology on $\pregeo_{\mathrm{Zar}}(k)\subset \geo_{\mathrm{Zar}}(k)$, using the same notations:
\begin{enumerate}[(1)]
   \item A morphism $\mathrm{Spec}\,A\rightarrow \mathrm{Spec}\,B$ is admissible if and only if there exists some element $b\in B$ such that the map $B\rightarrow A$ induces an isomorphism $B[1/b]\cong A$. \item A collection of admissible morphism $\{\mathrm{Spec}\,B_i[1/b_i]\rightarrow \mathrm{Spec}\,B\}_i$ is an admissible covering if and only if the elements $\{b_i\}$ generate the unit ideal in $B$.  
\end{enumerate}

\end{ex}
\begin{ex}[Geometry of affine $k$-schemes (\'{e}tale)]\label{geomalget}
Continuing the notation of the previous examples, we let $\geo_{\text{\'{e}t}}(k)$ be the geometry that has the same underlying $\infty$-category as $\geo_{\mathrm{Zar}}(k)$, and whose admissibility structure and topology are as follows:
\begin{enumerate}[(1)]
    \item A morphism $\mathrm{Spec}\,A\rightarrow \mathrm{Spec}\,B$ is admissible if and only if the map $B\rightarrow A$ is an \'{e}tale map of commutative $k$-algebras. 
    \item A collection of admissible morphism $\{\mathrm{Spec}\,B_i\rightarrow \mathrm{Spec}\,B\}_i$ is an admissible covering if and only if there exists a finite set of indices $\{i_j\}_{1\leq j\leq n}$ such that the induced map $g:B\rightarrow \prod_{1\leq j\leq n}B_{i_j}$ is faithfully flat.
\end{enumerate}
\end{ex}

\subsubsection{Geometric envelopes}
A pregeometry $\pregeo$ is distillation of the features shared by various theories of \emph{nonsingular affine schemes}, like the examples above. To each such theory, we can associate a geometry $\geo$ corresponding to a theory of singular \emph{derived} affine schemes in a $(\infty,2)$-functorial manner, a procedure which one morally should think of as adjoining to the \infcat $\pregeo$ all pullbacks along maps which are \emph{not} admissible.
\begin{defn}
Let $\pregeo$ be a pregeometry, then a functor $f:\pregeo\rightarrow\geo$ to a geometry $\geo$ \emph{exhibits $\geo$ as geometric envelope of $\pregeo$} if $f$ is a transformation of pregeometries, and for each geometry $\geo$, composition with $f$ induces an equivalence
\[ \fun^{\mathrm{lex}}_{\tau}(\geo,\geo')\overset{\simeq}{\longrightarrow}\fun^{\pi\mathrm{ad}}_{\tau}(\pregeo,\geo). \]
\end{defn}
We will now discuss this notion from another perspective, which we will not use in the rest of this work. Consider the $\sset^+$-enriched functor
\[ \mathsf{PGEO}^{op}\times \mathsf{GEO}\longrightarrow\Catinf   \]
carrying a pair $(\pregeo,\geo)$ to the \infcat $\fun^{\pi\mathrm{ad}}_{\tau}(\pregeo,\geo)$, viewed as a fibrant object of $\sset^+$ with equivalences marked. Here $\Catinf$ denotes the $(\infty,2)$-category of \infcatst, in our model presented by the fibrant marked simplicial category $(\sset^+)^{\mathrm{fc}}$ of fibrant-(cofibrant) objects of $\sset^+$. By adjunction, this functor of $(\infty,2)$-categories is equivalent to a functor
\[ \fun^{\pi\mathrm{ad}}_{\tau}(\_,\_):\mathsf{PGEO}^{op}\longrightarrow \mathsf{FUN}(\mathsf{GEO},\mathsf{CAT}_{\infty}),  \]
where $\mathsf{FUN}$ denotes the internal hom in $(\infty,2)$-categories. 
\begin{prop}[$(\infty,2)$-categorical Yoneda lemma]\label{prop:infty2yoneda}
Let $\icat$ be an $(\infty,2)$-category, then there is a fully faithful functor 
\[j_{\icat}: \icat \hooklongrightarrow \Fun(\icat^{op},\mathsf{CAT}_{\infty}),    \]
the \emph{Yoneda embedding}. Furthermore, for any $F\in \Fun(\icat^{op},\Catinf)$ and $C\in \icat$, there is a natural equivalence 
\[\HOM_{\Fun(\icat^{op},\Catinf)}(j_{\icat}(C),F)\simeq F(C)\] of \infcatst. 
\end{prop}
This is not hard to prove using the model of $(\infty,2)$-categories as fibrant $\sset^+$-enriched categories; essentially the same proof as the the $\infty$-categorical Yoneda lemma given in \cite{HTT}, Section 5.1 applies.  
\begin{defn}
Let $\icat$ be an $(\infty,2)$-category, then an $(\infty,2)$-functor $F:\icat\rightarrow \Catinf$ is \emph{representable} if it lies in the image of the Yoneda embedding $j_{\icat}$.
\end{defn}
With the Yoneda lemma, it becomes easy to interpret operations such as completions under colimits and limits, and localizations as $(\infty,2)$-functors. Here is one example of how this works.
\begin{ex}
Let $\mathsf{SITE}$ be the fibrant $\sset^+$-enriched category whose objects are pairs $(\icat,\tau)$ of a small \infcat $\icat$ admitting finite limits together with a Grothendieck topology $\tau$, and whose simplicial set of morphisms between sites $\icat$ and $\icat'$ is the full subcategory $\fun'(\icat,\icat')\subset\fun(\icat,\icat')$ spanned by functors carrying covering families to covering families. Consider the $\sset^+$-enriched functor 
\[  \mathsf{SITE}^{op}\times \lTop \longrightarrow (\sset^+)^{\mathrm{fc}}  \]
between fibrant $\sset^+$-enriched categories carrying a pair $(\icat,\xtop)$ to the full subcategory $\fun^*(\icat,\xtop)\subset \fun(\icat,\xtop)$ spanned by functors $f$ that satisfy the following conditions.
\begin{enumerate}[$(1)$]
    \item $f$ preserves finite limits.
    \item For every collection of morphism $\{C_i\rightarrow C\}_i$ generating a covering sieve, the morphism $\coprod_i f(C_i)\rightarrow f(C)$ is an effective epimorphism in $\xtop$.
\end{enumerate}
We regard $\fun^*(\icat,\xtop)$ as fibrant marked simplicial set with equivalences marked. For each site $(\icat,\tau)$, consider the composition $\icat\overset{j}{\rightarrow}\pshv(\icat)\overset{L}{\rightarrow}\shv(\icat)$, where $L$ is a sheafification functor for the topology $\tau$. Then it follows from \cite{HTT}, Proposition 6.2.3.20 that composition with $L\circ j$ induces for each \inftop $\xtop$ an equivalence $\fun^*(\shv(\icat),\xtop)\rightarrow\fun^*(\icat,\xtop)$. Thus, composition with $L\circ j$ determines a natural equivalence from $\fun^*(\shv(\icat),\_)$ to $\fun^*(\icat,\_)$. The former functor lies in the image of the Yoneda embedding, so we deduce the existence of an $(\infty,2)$-functor $\mathsf{SITE}\rightarrow\lTop$ carrying $(\icat,\tau)$ to $\shv(\icat)$. On morphism \infcatst, this functor is the composition 
\[\fun(\icat,\icat')\longrightarrow\fun^*(\icat,\shv(\icat')) \overset{\simeq}{\longleftarrow}\fun^* (\shv(\icat),\shv(\icat')).  \]
\end{ex}
Let $\pregeo$ be pregeometry, then a functor $\pregeo\rightarrow\geo$ in $\fun^{\pi\mathrm{ad}}(\pregeo,\geo)$ \emph{exhibits a $\geo$ as a geometric envelope for $\pregeo$} if the corresponding map $j(\geo)\rightarrow \fun^{\pi\mathrm{ad}}_{\tau}(\pregeo,\_)$ in $\mathsf{FUN}(\mathsf{GEO},\Catinf)$ is an equivalence. That is, a pregeometry admits a geometric envelope if and only if the $(\infty,2)$-functor $\fun^{\pi\mathrm{ad}}_{\tau}(\pregeo,\_)$ is representable.\\
Every pregeometry admits a geometric envelope.
\begin{prop}\label{prop:geoenvexist}
Let $\pregeo$ be a pregeometry. Then there is a geometry $\geo$ and a functor $f:\pregeo\rightarrow\geo$ that exhibits a geometric envelope. Moreover, $f$ is fully faithful.
\end{prop}
\begin{proof}
Let $\mathcal{R}$ be the union of the collection of all diagrams $\{S\rightarrow \pregeo^{op}\}$ where $S$ runs over all finite sets viewed as discrete simplicial sets and the collection of all pushout diagrams $(\Lambda^2_0)^{\rhd}\rightarrow\pregeo^{op}$ that carry at least one of the nondegenerate edges of $\Lambda^2_0$ to an admissible morphism, and let $\mathcal{K}$ be the collection of all finite simplicial sets together with the \infcat $\mathrm{Idem}$, the free \infcat on an arrow idempotent up to coherent homotopy. Consider the fully faithful functor $\pregeo^{op}\hookrightarrow\pshv^{\mathcal{K}}_{\mathcal{R}}(\pregeo^{op})$ provided by \cite{HTT}, Proposition 5.3.6.2 and let $\geo:=\pshv^{\mathcal{K}}_{\mathcal{R}}(\pregeo^{op})^{op}$, so that $\geo$ admits finite limits and retracts for idempotents and $f:\pregeo\rightarrow\geo$ is fully faithful. Endow $\geo$ with the intersection of all admissibility structures that contain $f(\pregeo^{\mathrm{ad}})$; this is an admissibility structure. Then endow $\geo$ with the intersection of all Grothendieck topologies compatible with this admissibility structure that contain the images of admissible coverings under $f$, then we claim that the resulting geometry is a geometric envelope for $\pregeo$. Invoking \cite{HTT}, Proposition 5.3.6.2, we only have to verify that if $g:\pregeo\rightarrow\geo'$ carries $\pregeo^{\mathrm{ad}}$ to $\geo^{\prime \mathrm{ad}}$, then the associated left exact functor $G:\geo\rightarrow\geo'$ carries $\geo^{\mathrm{ad}}$ to $\geo^{\prime \mathrm{ad}}$, and that if $g$ carries admissible coverings to admissible coverings, then so does $G$. For the first assertion, we note that $G^{-1}(\geo^{\prime \mathrm{ad}})$ is an admissibility structure on $\geo$ containing $f(\pregeo^{\mathrm{ad}})$, so that $\geo^{\mathrm{ad}}\subset G^{-1}(\geo^{\prime \mathrm{ad}})$ and therefore $G(\geo^{\mathrm{ad}})\subset  \geo^{\prime \mathrm{ad}}$. For the second assertion, we note that the collection of families $\{U_i\rightarrow X\}_i$ of admissible morphisms in $\geo$ that are carried to an admissible covering by $G$ determines a basis for a topology that contains the images of the admissible coverings under $f$.
\end{proof}
\begin{cor}
There is an $(\infty,2)$-functor $\mathsf{PGEO}\rightarrow\mathsf{GEO}$ carrying each pregeometry to a geometric envelope. On morphism \infcatst, this functor is the composition
\[ \fun_{\tau}^{\pi\mathrm{ad}}(\pregeo,\pregeo') \longrightarrow  \fun_{\tau}^{\pi\mathrm{ad}}(\pregeo,\geo')\overset{\simeq}{\longleftarrow} \fun_{\tau}^{\mathrm{lex},\mathrm{ad}}(\geo,\geo').\] 
\end{cor}
\begin{rmk}
For the reader not comfortable with Proposition \ref{prop:infty2yoneda}, we note that the same constructions work after passing to maximal Kan complexes of morphisms in all $(\infty,2)$-categories; that is, we have a functor
\[\fun_{\tau}^{\pi\mathrm{ad}}(\_,\_)^{\simeq}:\mathsf{PGeo}^{op}\longrightarrow \pshv(\mathsf{Geo}^{op}) \]
which Proposition \ref{prop:geoenvexist} guarantees factors through $\mathsf{Geo}^{op}$. The representability of the functor $\fun_{\tau}^{\pi\mathrm{ad}}(\pregeo,\_)^{\simeq}$ is equivalent to the representability of the functor $\fun_{\tau}^{\pi\mathrm{ad}}(\pregeo,\_)$, by Proposition \ref{prop:geoenvexist} again.
\end{rmk}
\begin{ex}[Geometric envelopes of $\pregeo_{\mathrm{Zar}}(k)$]
For $k$ a commutative ring, we may consider the $\infty$-category of \emph{simplicial commutative $k$-algebras}, denoted $\scring_k$, which is defined as the $\infty$-category of algebras for the finite limit theory whose objects are all the affine $k$-spaces $\mathbb{A}^n_k$ and whose morphisms are polynomial maps. The $\infty$-category $\geo^{\mathrm{der}}_{\mathrm{Zar}}(k)$ is defined as the opposite of the $\infty$-category of finitely presented objects in the presentable $\infty$-category $\sring_k$. We have an equivalence $\tau_{\leq 0}\scring_k\simeq \mathsf{CAlg}^0_k$ which remains an equivalence after restricting to finitely presented objects. Again, we may define for each simplicial $k$-algebra $B$ and each $b\in \tau_{\leq0}B$, a localization $B\rightarrow B[1/b]$ defined up to equivalence. We endow $\geo^{\mathrm{der}}_{\mathrm{Zar}}(k)$ with the following admissibility structure:
\begin{enumerate}[(1)]
   \item A morphism $\mathrm{Spec}\,A\rightarrow \mathrm{Spec}\,B$ is admissible if and only if there exists some element $b\in \tau_{\leq0}B$ such that the map $B\rightarrow A$ induces an isomorphism $B[1/b]\cong A$. 
   \item A collection of admissible morphism $\{\mathrm{Spec}\,B_i[1/b_i]\rightarrow \mathrm{Spec}\,B\}_i$ is an admissible covering if and only if the elements $\{b_i\}$ generate the unit ideal in $B$.  
\end{enumerate}
There is an obvious functor $\pregeo_{\mathrm{Zar}}(k)\hookrightarrow \geo^{\mathrm{der}}_{\mathrm{Zar}}(k)$
which, according to \cite{dagv}, Proposition 4.2.3, exhibits $\geo^{\mathrm{der}}_{\mathrm{Zar}}(k)$ as a geometric envelope of $\pregeo_{\mathrm{Zar}}(k)$. 
\end{ex}
The construction of the geometric envelope of the \emph{\'{e}tale} pregeometry follows along similar lines.\\
One may wonder what the geometric envelope of $\diff$ look like. A partial answer is given in the subsection below.
\subsection{Structure sheaves and structured spaces}
Since we asserted that pregeometries and geometries capture salient categorical properties of affine schemes, we should also expect to be able to associate to an object $C$ in a pregeometry or geometry a `spectrum', that is, a pair $(X_C,\Of_{X_C})$ where $X_C$ is a topological space and $\Of_{X_C}$ is some sort of sheaf of local rings on $X_C$. This picture is essentially correct, up to two modifications.
\begin{enumerate}[$(1)$]
    \item The `space' $X_C$ underlying the spectrum will in general not be a topological space, but an \inftopt. If we adopt the very mild condition that our topological spaces be \emph{sober}, letting $X_C$ be an \inftop constitutes a generalization as we can consider the category of sober topological spaces as a full subcategory of the \infcat of \inftopoit, by carrying a space $X$ to its \inftop of sheaves $\shv(X)$. This formulation allows us to interpret the object $\Of_{X}$ in several different ways.
    \begin{enumerate}[$(a)$]
        \item $\Of_X$ is a sheaf of commutative rings on $X$; that is, $\Of_X$ is a functor
        \[ \mathrm{Open}(X)^{op}\longrightarrow \mathsf{CAlg}^0_{\Z}  \]
        satisfying the sheaf condition. By right Kan extension along the Yoneda embedding $\mathrm{Open}(X)^{op}\hookrightarrow\shv(X)^{op}$, such a functor is the same thing as a functor
        \[ \Of_X:\shv(X)^{op}\longrightarrow \mathsf{CAlg}^0_{\Z}  \]
        preserving all limits (\cite{HTT}, Proposition 6.2.3.20). By the adjoint functor theorem, the condition that $\Of_X$ preserves all limits is equivalent to the condition that $\Of_X$ is a right adjoint. 
        \item For formal reasons, we have an equivalence of \infcats \[\fun^{\mathrm{R}}(\shv(X)^{op},\mathsf{CAlg}^0_{\Z})\simeq\fun^{\mathrm{R}}((\mathsf{CAlg}^{0}_{\Z})^{op},\shv(X)),\]
        so we can think of $\Of_X$ as a limit preserving functor $(\mathsf{CAlg}^{0}_{\Z})^{op}\rightarrow\shv(X)$. In terms of the previous descriptions, we can informally describe this functor as carrying a commutative ring $A$ to the sheaf on $X$ that assigns to an open set $U\subset X$ the set $\Hom_{\mathsf{CAlg}^0_{\Z}}(A,\Of_X(U))$.
        \item Because every commutative ring $A$ may be written as a filtered colimit of its finitely presented subrings, $\mathsf{CAlg}^0_{\Z}$ is \emph{compactly generated}: the inclusion of finitely presented commutative rings $(\mathsf{CAlg}^0_{\Z})_{\fp}\subset \mathsf{CAlg}^0_{\Z}$ induces an equivalence of categories $\mathrm{Ind}((\mathsf{CAlg}^0_{\Z})_{\fp})\simeq \mathsf{CAlg}^0_{\Z}$. The inclusion $(\mathsf{CAlg}^0_{\Z})_{\fp}\subset \mathsf{CAlg}^0_{\Z}$ induces an equivalence
        \[   \fun^{\mathrm{R}}((\mathsf{CAlg}^{0}_{\Z})^{op},\shv(X))  \overset{\simeq}{\longrightarrow} \fun^{\mathrm{lex}}((\mathsf{CAlg}^{0}_{\Z})_{\fp}^{op},\shv(X))\]
        with inverse being right Kan extension, so we conclude that $\Of_{X}$ may be regarded as a functor $\Of_X:(\mathsf{CAlg}^{0}_{\Z})_{\fp}^{op}\rightarrow\shv(X)$ that preserves finite limits. Phrased in this manner, the condition of $\Of_X$ being \emph{local} takes on the following form: for every collection $\{a_i\}_i\subset A$ that generates the unit ideal in $A$, the map 
        \[ \coprod_i \Of_X(A[a^{-1}_i]) \longrightarrow \Of_X(A)  \]
        is an effective epimorphism in $\shv(X)$ (we show a version of this statement below in the differential geometric setting). 
    \end{enumerate}
    \item For an object $C$ in some geometry $\geo$, it does not make much sense to ask for an \inftop $\xtop_C$ that comes equipped with a sheaf of commutative rings, since such a datum makes no reference to $\geo$. Only when our geometry describes ordinary algebraic geometry over the ground ring $\Z$ should we expect that a locally ringed space is the correct notion. In fact, we have already encountered that specific geometry: it is $\geo_{\mathrm{Zar}}(\Z)$. Furthermore, we observe that our last formulation of the object $\Of_X$ defines it as a functor from $\geo_{\mathrm{Zar}}(\Z)$ in such a way that the locality condition makes reference only to the admissibility structure. A natural generalization suggests itself: a \emph{$\geo$-structure} on an \inftop $\xtop$ is a functor
    \[ \Of:\geo\longrightarrow\xtop \]
    that is both left exact and takes admissible coverings to effective epimorphisms in $\xtop$. 
\end{enumerate}

\begin{defn}\label{defn:tstruct}
Let $\xtop$ be an $\infty$-topos.
\begin{enumerate}[(1)]
    \item For $\mathcal{T}$ a pregeometry, a \emph{$\pregeo$-structure on $\xtop$} is a functor
    \[ \Of_{\xtop}:\pregeo\longrightarrow\xtop \]
    in $\fun^{\pi\mathrm{ad}}(\mathcal{T},\xtop)$ that satisfies the following condition: for every admissible covering $\{U_i\rightarrow X\}_i$, the induced map 
    \[ \coprod_i  \Of_{\xtop}(U_i) \longrightarrow   \Of_{\xtop}(X)  \]
    is an effective epimorphism in $\xtop$. The full subcategory of $\fun^{\pi\mathrm{ad}}(\mathcal{T},\xtop)$ spanned by $\pregeo$-structures is denoted $\mathrm{Str}_{\mathcal{T}}(\xtop)$.
    \item For $\mathcal{G}$ a geometry, a \emph{$\geo$-structure on $\xtop$} is a left exact functor 
    \[   \Of_{\xtop}:\geo\longrightarrow\xtop \]
   satisfying the same descent condition we formulated above for pregeometries. The full subcategory of $\fun^{\mathrm{lex}}(\geo,\xtop)$ spanned by $\geo$-structures is denoted $\mathrm{Str}_{\geo}(\xtop)$.
\end{enumerate}
\end{defn}
\begin{warn}
Note that while a geometry $\mathcal{G}$ can also be viewed as a pregeometry, a $\mathcal{G}$-structure on an $\infty$-topos $\xtop$ \emph{with $\mathcal{G}$ viewed as a geometry} is \emph{not} the same thing as a $\mathcal{G}$-structure on $\xtop$ \emph{with $\mathcal{G}$ viewed as a pregeometry}. To prevent ambiguity, we will always use the symbol $\mathcal{G}$ to mean a geometry, and $\mathcal{T}$ to mean a pregeometry. 
\end{warn}
\begin{defn}\label{defn:localmorphism}
For $\pregeo$ a pregeometry and $\Of_{\xtop}$, $\Of_{\xtop}'$ functors in $\fun^{\pi\mathrm{ad}}(\pregeo,\xtop)$, then a \emph{local morphism between $\Of_{\xtop}$ and $\Of_{\xtop}'$} is a natural transformation $\alpha:\Of_{\xtop}\rightarrow \Of_{\xtop}'$ such that for all admissible maps $X\rightarrow Y$ of $\pregeo$, the natural commuting diagram
\[
\begin{tikzcd}
\Of_{\xtop}(X) \ar[d] \ar[r] & \Of_{\xtop}'(X) \ar[d] \\
\Of_{\xtop}(Y) \ar[r] & \Of_{\xtop}'(Y)
\end{tikzcd}    
\]
is a pullback square. We denote by $\fun^{\pi\mathrm{ad},\loc}(\pregeo,\xtop)\subset \fun^{\pi\mathrm{ad}}(\pregeo,\xtop)$ the subcategory on the local morphisms. For $\geo$ a geometry, the same definition applied to left exact functors $\geo\rightarrow\xtop$ yields the subcategory $\fun^{\lex,\loc}(\geo,\xtop)\subset\fun^{\lex}(\geo,\xtop)$ the subcategory on the local morphisms.\\
We let $\strloc_{\pregeo}(\xtop)=\str_{\pregeo}(\xtop)\cap\fun^{\pi\mathrm{ad},\loc}(\pregeo,\xtop)$ and $\strloc_{\geo}(\xtop)=\str_{\geo}(\xtop)\cap\fun^{\lex}(\geo,\xtop)$ denote the categories whose objects are $\pregeo/\geo$-structures and whose morphisms are local morphisms.
\end{defn}
\begin{prop}\label{prop:locmorstab}
Let $\geo$ be a geometry and let $\xtop$ be an \inftopt, then the following hold true.
\begin{enumerate}[$(1)$]
\item Consider a commuting diagram
\[
\begin{tikzcd}
& \Of' \ar[dr,"g"] \\
\Of \ar[rr,"h"] \ar[ur,"f"] && \Of''
\end{tikzcd}
\]
of left exact functors $\geo\rightarrow\xtop$. Suppose that $g$ is a local morphism, then $f$ is a local morphism if and only if $h$ is a local morphism.
\item Let $f:\Of\rightarrow\Of'$ be a local morphism of left exact functors. If $\Of'$ is a $\geo$-structure, then $\Of$ is a $\geo$-structure as well.
\end{enumerate}
The same results hold for pregeometries and functors preserving finite products and admissible pullbacks.
\end{prop}
\begin{proof}
The assertion $(1)$ is an immediate consequence of the pasting property of pullback squares. For the second assertion, we wish to show that for any admissible covering $\{U_i\rightarrow X\}_i$, the map $\coprod_i\Of(U_i)\rightarrow\Of(X)$ is an effective epimorphism. Since $\Of'$ is a $\geo$-structure and effective epimorphisms in an \inftop are stable under the formation of pullbacks, it suffices to show that the diagram
\[
\begin{tikzcd}
\coprod_i\Of(U_i)\ar[d]\ar[r] &\Of(X)\ar[d]\\
\coprod_i\Of'(U_i)\ar[r] &\Of'(X)
\end{tikzcd}
\]
induced by $f$ is a pullback. Since colimits are universal in an \inftopt, this follows from the assumption that $f$ is a local morphism. The proof of the analogous results for $\pregeo$-structures is the same.
\end{proof}
We make some observations about the formation of limits and colimits in \infcats of structures.
\begin{prop}\label{prop:struccolimit}
Let $\xtop$ be an \inftopt, then the following hold true.
\begin{enumerate}[$(1)$]
\item Let $\geo$ be a geometry, then all functors in the commuting square 
\[
\begin{tikzcd}
\strloc_{\geo}(\xtop)\ar[d]\ar[r] & \fun^{\lex,\loc}(\geo,\xtop) \ar[d] \\
\str_{\geo}(\xtop) \ar[r] & \fun^{\lex}(\geo,\xtop)
\end{tikzcd}
\]
preserve filtered colimits; in particular, all \infcats in this square admit filtered colimits.
\item Let $\pregeo$ be a pregeometry, then the \infcats $\strloc_{\pregeo}(\xtop)$ and $\fun^{\pi\mathrm{ad},\loc}(\geo,\xtop)$ admit sifted colimits and the inclusions $\strloc_{\pregeo}(\xtop)\subset \fun^{\pi\mathrm{ad},\loc}(\geo,\xtop)$ and $\fun^{\pi\mathrm{ad},\loc}(\geo,\xtop)\subset \fun^{\pi\mathrm{ad}}(\pregeo,\xtop)$ preserve sifted colimits.
\end{enumerate}
\end{prop}
\begin{proof}
For the first assertion, we note that the full subcategory $\fun^{\lex}(\geo,\xtop)\subset\fun(\geo,\xtop)$ is stable under filtered colimits so it suffices to argue the following.
\begin{enumerate}[$(a)$]
\item A filtered colimit of $\geo$-structures taken in $\fun(\geo,\xtop)$ carries admissible coverings to effective epimorphisms.
\item Let $K\rightarrow \fun^{\lex,\loc}(\geo,\xtop)$ be a filtered diagram and let $f:K^{\rhd}\rightarrow\fun(\geo,\xtop)$ be colimit diagram extending this diagram. Then for any $k\in K$, the map $f(k)\rightarrow f(\infty)$ has the property that for any admissible map $U\rightarrow X$, the square 
\[
\begin{tikzcd}
f(k)(U) \ar[d]\ar[r] & f(\infty)(U) \ar[d] \\
f(k)(X)\ar[r]& f(\infty)(X)
\end{tikzcd}
\]
is a pullback diagram.
\end{enumerate}
The claim $(a)$ follows from the fact that the collection of effective epimorphisms is stable under colimits of arrows in any \inftopt. The second follows from the assumption that $f(k\rightarrow k')$ is a local morphism and that filtered colimits are left exact in \inftopoit. To prove the second assertion, it suffices to show that the compositions $\fun^{\pi\mathrm{ad},\loc}(\geo,\xtop)\subset \fun^{\pi\mathrm{ad}}(\pregeo,\xtop)\subset\fun(\pregeo,\xtop)$ and $\strloc_{\pregeo}(\xtop)\subset \fun^{\pi\mathrm{ad},\loc}(\geo,\xtop)\subset\fun(\pregeo,\xtop)$ are stable under sifted colimits. We have to show the following.
\begin{enumerate}[$(i)$]
\item Let $K\rightarrow \fun^{\pi\mathrm{ad},\loc}(\pregeo,\xtop)$ be a sifted diagram and let $f:K^{\rhd}\rightarrow \fun(\pregeo,\xtop)$ be a colimit diagram extending this diagram, then $f(\infty)$ preserves finite limits and pullbacks along admissible maps, and for every $k\in K$, the map $f(k)\rightarrow f(\infty)$ is a local morphism.
\item A sifted colimit of $\pregeo$-structures and local transition morphisms taken in $\fun^{\pi\mathrm{ad},\loc}(\pregeo,\xtop)$ carries admissible coverings to effective epimorphisms.
\end{enumerate}
We prove $(i)$. Since the formation of sifted colimits preserves products, the functor $f(\infty)$ preserves products. We now show that for each $k\in K$ and each admissible map $U\rightarrow X$ in $\pregeo$, the diagram 
\[
\begin{tikzcd}
f(k)(U) \ar[d]\ar[r] & f(\infty)(U) \ar[d] \\
f(k)(X)\ar[r]& f(\infty)(X)
\end{tikzcd}
\]
is a pullback square. This amounts to the assertion that the diagram $K^{\rhd}\times\Delta^1\rightarrow \xtop$ induced by $f$ and the map $\Delta^1 \rightarrow \xtop$ classifying $U\rightarrow X$ is a Cartesian transformation. Since $f$ is a colimit diagram and $f|_{K}$ is a Cartesian transformation by virtue of the assumption that $f|_{K}$ takes values in $\fun^{\lex,\loc}(\geo,\xtop)$, this follows from $(4)$ of \cite{HTT}, Theorem 6.1.3.9. It remains to be shown $f(\infty)$ preserves pullbacks along admissible maps. Let $U\rightarrow X$ be an admissible and $Y\rightarrow X$ any map, then for any $k\in K$, we have pullback diagrams  
\[
\begin{tikzcd}
f(k)(U) \ar[d]\ar[r] & f(\infty)(U) \ar[d] \\
f(k)(X)\ar[r]& f(\infty)(X)
\end{tikzcd} \quad\quad \begin{tikzcd}
f(k)(U\times_XY) \ar[d]\ar[r] & f(\infty)(U\times_XY) \ar[d] \\
f(k)(Y)\ar[r]& f(\infty)(Y)
\end{tikzcd}
\]
indexed by $K$. Since $f(k)$ preserves pullbacks along admissibles, $f(\infty)$ does so as well (this proof is nothing but the assertion that a Cartesian transformation is a realization fibration). As we have just established, we can take a sifted colimit of $\pregeo$-structures and local transition morphisms in $\fun(\pregeo,\xtop)$, so $(ii)$ follows from the fact that the collection of effective epimorphisms is stable under colimits of arrows in any \inftopt.  
\end{proof}
\begin{rmk}\label{gstructurealgebrasheaf}
Let $\icat$ be an $\infty$-category with finite limits, and let $\xtop$ be an $\infty$-topos. By \cite{HTT}, Proposition 5.5.1.9 and 5.3.5.10, the Yoneda embedding $j:\icat\rightarrow \mathrm{Pro}(\icat)$ induces an equivalence
\[\fun^{\mathrm{R}}(\mathrm{Pro}(\icat),\xtop) \overset{\simeq}{\longrightarrow}  \fun^{\mathrm{lex}}(\icat,\xtop),\]
where the left hand side is the $\infty$-category of functors that admit a left adjoint. We have a natural equivalence $\fun^{\mathrm{R}}(\mathrm{Pro}(\icat),\xtop)\simeq \fun^{\mathrm{R}}(\xtop^{op},\mathrm{Ind}(\icat^{op}))$, and because a functor from $\xtop^{op}$ to $\mathrm{Ind}(\icat^{op})$ admits a left adjoint if and only if it preserves small limits (apply the adjoint functor theorem \cite{HTT}, prop 5.5.2.9 to the functor of opposite categories), we find that $\fun^{\mathrm{R}}(\xtop^{op},\mathrm{Ind}(\icat^{op}))= \shv_{\mathrm{Ind}(\icat^{op})}({\xtop})$. We conclude that for a \emph{discrete} geometry $\mathcal{G}$, a $\mathcal{G}$-structure on $\xtop$ is an $\mathrm{Ind}(\mathcal{G}^{op})$-valued sheaf on $\xtop$. To cement this intuition, we encourage the reader to have for $\mathcal{G}^{op}$ the category $(\mathsf{CAlg}^0_k)_{\mathrm{fp}}$ of $k$-algebras of the form $k[x_1,\ldots,x_n]/I$ in mind, where $k$ is a commutative ring and $I$ a finitely generated ideal, whose opposite category is the category of finitely presented affine $k$-schemes. By taking the ind-completion we obtain $\mathrm{Ind}((\mathsf{CAlg}^0_k)_{\mathrm{fp}})\simeq \mathsf{CAlg}^0_k$, the category of all commutative $k$ algebras. Thus, a $(\mathsf{CAlg}^0_k)^{op}_{\mathrm{fp}}$-structure on an $\infty$-topos $\shv(X)$ for $X$ a topological space can be canonically identified with a sheaf of commutative $k$-algebras on the space $X$.   
\end{rmk}
We show that a pregeometry and its geometric envelope describe the same structured spaces.
\begin{prop}\label{prop:geoenvstructures}
Let $\pregeo$ be a pregeometry, and let $f:\pregeo\hookrightarrow\geo$ exhibit $\geo$ as a geometric envelope, then composition with $f$ induces for any \inftop $\xtop$ an equivalence $\str_{\pregeo}(\xtop)\simeq \str_{\geo}(\xtop)$ which restricts to an equivalence $\strloc_{\pregeo}(\xtop)\simeq \strloc_{\geo}(\xtop)$. 
\end{prop}
\begin{proof}
We show the first equivalence. By definition of a geometric envelope, we are reduced to showing the following assertion: a left exact functor $\Of:\geo\rightarrow\xtop$ is a $\geo$-structure if and only if its restriction $\Of|_{\pregeo}$ is a $\pregeo$-structure. The `only if' direction is immediate. For the other direction, we note that the collection of families $\{U_i\rightarrow X\}_i$ of admissible morphisms in $\geo$ that have the property that $\coprod_i\Of(U_i)\rightarrow \Of(X)$ is an effective epimorphism in $\xtop$ determines a Grothendieck pretopology $\tau_{\xtop}$ on $\geo$ compatible with the admissibility structure. By assumption, every family of the form $\{f(V_i)\rightarrow f(Y)\}_i$ for $\{V_i\rightarrow Y\}_i$ a covering family for the pretopology on $\pregeo$ is a covering family for the pretopology $\tau_{\xtop}$. Since the Grothendieck topology $\tau$ on $\geo$ is the coarsest one compatible with the admissibility structure that has this property, we conclude that every covering family of $\tau$ lies in $\tau_{\xtop}$. For the second equivalence, the only nontrivial point is to verify that a morphism $\alpha:\Of\rightarrow \Of'$ is local if its restriction to $\pregeo$ is local. Let $\geo^{\mathrm{ad}}_{\xtop}$ be the subcategory on the morphisms $f:U\rightarrow X$ that have the following property: for every morphism $\Of'\rightarrow \Of$ of $\geo$-structures, the natural diagram 
\[
\begin{tikzcd}
\Of_{\xtop}(U) \ar[d] \ar[r] & \Of_{\xtop}'(U) \ar[d] \\
\Of_{\xtop}(X) \ar[r] & \Of_{\xtop}'(X)
\end{tikzcd}    
\]
is a pullback square, then one readily verifies using pasting of pullback squares that $\geo^{\mathrm{ad}}_{\xtop}$ is an admissibility structure containing $f(\pregeo^{\mathrm{ad}})$, so we conclude that $\geo^{\mathrm{ad}}\subset \geo^{\mathrm{ad}}_{\xtop}$
\end{proof}

If $\pregeo\hookrightarrow\geo$ exhibits a geometric envelope, the theories of $\pregeo$-structures and $\geo$-structures coincide. Similarly, different pregeometries might have the same theory of structured spaces. It will prove useful to investigate the ambiguity.
\begin{defn}\label{defn:moritaequivalencepregeo}
Let $f:\pregeo\rightarrow\pregeo$ be a transformation of pregeometries, then $f$ is a \emph{Morita equivalence of pregeometries} if for any \inftop $\xtop$, composition with $f$ induces an equivalence $\strloc_{\pregeo'}(\xtop)\simeq \strloc_{\pregeo}(\xtop)$.
\end{defn}
In Section 3.2 of \cite{dagv} some criteria for a transformation of pregeometries to be a Morita equivalence are provided. Morita equivalences are typically of the form $\pregeo\subset\pregeo'$ where $\pregeo'$ is generated under colimits by admissibles from $\pregeo$. 
\begin{ex}
The inclusion $\diff^{\mathrm{open}}\subset\diff$ is a Morita equivalence.    
\end{ex}
\begin{ex}
Let $\mathsf{Man}_c$ be the \infcat of manifolds with corners, then the \'{e}tale topology and the admissibility structure by open embeddings give this \infcat the structure of a pregeometry. The inclusions $\diffc^{\mathrm{Open}}\subset\diffc\subset \mathsf{Man}_c$. 
\end{ex}
\subsubsection{Factorization systems on \infcats of structure sheaves}
We have introduced the theory of pregeometries and geometries by appealing to the problem of constructing derived geometric objects from nonsingular ones. Just now, the topology and admissibility structure played a somewhat different role: the topology on a geometry $\geo$ determines the full subcategory of $\fun^{\lex}(\geo,\xtop)$ of \emph{local} objects -the $\geo$-structures- and the admissibility structure determines the subcategory of $\fun^{\lex}(\geo,\xtop)$ of local morphisms. For the case $\geo=\geo_{\mathrm{Zar}}(k)$, we can, as in the preceding remark, identify left exact functors $\Of:\geo_{\mathrm{Zar}}(k)\rightarrow\spa$ with commutative $k$-algebras and under this identification, $\Of$ is a $\geo$-structure if and only if the associated $k$-algebra $A_{\Of}$ is a local ring. Similarly, a morphism $f:\Of\rightarrow \Of$ is local if and only if the associated map of $k$-algebras $f:A_{\Of}\rightarrow A_{\Of'}$ is local, that is, $f$ reflects invertible elements. Replacing $\spa$ by the \inftop of $
\shv(X)$, we find that $\geo$-structures in $\shv(X)$ and local morphisms of left exact functors $\geo\rightarrow\shv(X)$ correspond to the usual notions of a sheaf of commutative $k$-algebras with local stalks and morphisms of sheaves that are local on stalks. A geometry structure on an \infcat $\geo$ allows one to interpret these concepts for $\ind(\geo^{op})$-valued sheaves on arbitrary \inftopoit.\\
There is another feature of the theory of local morphisms of rings we wish to emulate: consider a map $f:A\rightarrow B$ of commutative $k$-algebras and let $S$ be the set of elements in $A$ that are carried to an equivalence in $B$, then $f$ factorizes as
\[ A\overset{f'}{\longrightarrow} A[S^{-1}]\overset{f''}{\longrightarrow} B \] 
where $A\rightarrow A[S^{-1}]$ is the localization of $A$ at $S$ and $f''$ is a local morphism. Moreover, this factorization is functorial in both $A$ and $B$. The following result generalizes this observation to arbitrary geometries and \inftopoit.
\begin{prop}\label{prop:localfactorization}
Let $\geo$ be a geometry and let $\xtop$ be an \inftopt, then the following hold true.
\begin{enumerate}[$(1)$]
\item There is a factorization system (see \cite{HTT}, Definition 5.3.8.8) $(S_L,S_R)$ on the \infcat $\fun^{\lex}(\geo,\xtop)\simeq\shv_{\ind(\geo^{op})}(\xtop)$ for which the collection $S_{R}$ is the class of local morphisms.
\item In case $\xtop=\spa$, so that $\fun^{\lex}(\geo,\xtop)\simeq\ind(\geo^{op})$, the class $S_L$ admits the following description: let us say that morphism $f:A\rightarrow B$ of $\ind(\geo^{op})$ is \emph{ind-admissible} if $f$ is a colimit in $\fun(\Delta^1,\ind(\geo^{op}))$ of a filtered diagram in $\fun(\Delta^1,(\geo^{\mathrm{ad}})^{op})$. Then the collection of ind-admissible morphisms is the collection $S_L$.
\item The factorization system of $(1)$ restricts to a factorization system on $\str_{\geo}(\xtop)$.
\item Let $\fun^{\loc}(\Delta^1,\fun^{\lex}(\geo,\xtop))\subset \fun(\Delta^1,\fun^{\lex}(\geo,\xtop))$ be the full subcategory spanned by local morphisms of left exact functors, then the inclusion admits a left adjoint relative to the projection $\ev_1:\fun(\Delta^1,\fun^{\lex}(\geo,\xtop))\rightarrow \fun^{\lex}(\geo,\xtop)$.
\item Let $\Of_{\xtop}:\geo\rightarrow\xtop$ be a $\geo$-structure then the subcategory inclusion $\strloc_{\geo}(\xtop)_{/\Of_{\xtop}}\subset \fun^{\lex}(\geo,\xtop)_{/\Of_{\xtop}}$ is strongly reflective, that is, the inclusion is full and admits a left adjoint $L_{\Of_{\xtop}}$ which is an $\omega$-accessible localization.
\end{enumerate} 
\end{prop}
\begin{proof}
The assertions $(1)$, $(2)$ and $(3)$ are contained in \cite{sag}, Theorem 20.2.2.5 and Propositions 21.3.1.1. To prove $(4)$, we invoke \cite{HTT}, Lemma 5.2.8.19. Let $\Of_{\xtop}:\geo\rightarrow\xtop$ be a left exact functor, and let $\fun^{\lex}(\geo,\xtop)'_{/\Of_{\xtop}}\subset \fun^{\lex}(\geo,\xtop)_{/\Of_{\xtop}}$ be the full subcategory spanned by local morphisms $\Of\rightarrow\Of_{\xtop}$, then it follows from $(4)$ that this inclusion admits a left adjoint. It follows from Proposition \ref{prop:locmorstab} that $\fun^{\lex}(\geo,\xtop)'_{/\Of_{\xtop}}=\fun^{\lex,\loc}(\geo,\xtop)_{/\Of_{\xtop}}$ and in case $\Of_{\xtop}$ is $\geo$-structure, also $\fun^{\lex}(\geo,\xtop)'_{/\Of_{\xtop}}=\strloc_{\geo}(\xtop)_{/\Of_{\xtop}}$. It remains to be shown that the collection of local morphisms is stable under filtered colimits, which follows from the fact that the formation of finite limits preserves filtered colimits in an \inftopt.
\end{proof}
\subsubsection{$\infty$-Categories of structured higher topoi}
We proceed with some general remarks concerning the functoriality of $\str_{\geo}(\xtop)$ in $\geo$ and $\xtop$. We have $(\infty,2)$-functors
\[ \str_{\_}(\_):\mathsf{PGEO}^{op}\longrightarrow \mathsf{FUN}(\mathsf{LTOP},\Catinf),\quad \quad \str_{\_}(\_):\mathsf{GEO}^{op}\longrightarrow \mathsf{FUN}(\mathsf{LTOP},\Catinf)  \]
where $\mathsf{LTOP}$ is the $(\infty,2)$-category of \inftopoit. We claim that these functors factor through the essential images of the Yoneda embeddings and determine $(\infty,2)$-functors 
\[ \mathsf{PGEO}\longrightarrow \mathsf{LTOP},\quad \quad \mathsf{GEO}\longrightarrow \mathsf{LTOP}.  \]
This is equivalent to the assertion that for any pregeometry $\pregeo$ (any geometry $\geo$) there exists an \inftop $\mathcal{K}_{\pregeo}$ and a $\pregeo$-structure $\Of:\pregeo\rightarrow\mathcal{K}_{\pregeo}$ such that for any \inftop $\xtop$, composition with $\Of$ determines an equivalence
\[   \fun^*(\mathcal{K},\xtop) \longrightarrow \str_{\pregeo}(\xtop), \]
and similarly for geometries. This is the \emph{universal $\pregeo$/$\geo$-structure}. For geometries, this follows immediately from \cite{HTT}, Proposition 6.2.3.20; indeed, the sheafified Yoneda embedding $j:\geo\rightarrow \pshv(\geo)\overset{L}{\rightarrow} \shv(\geo)$ exhibits a universal $\geo$-structure. For pregeometries, one first constructs a geometric envelope $\pregeo\hookrightarrow \geo$, then the composition 
\[ \pregeo \hooklongrightarrow \geo \hooklongrightarrow \pshv(\geo)\overset{L}{\longrightarrow} \shv(\geo) \]
exhibits a universal $\geo$-structure; this follows immediately from Proposition \ref{prop:geoenvstructures} below.
\begin{defn}\label{defn:structuredtopoi}
Taking maximal Kan complexes, the $(\infty,2)$-functor 
\[ \str_{\_}(\_):\mathsf{GEO}^{op}\times\mathsf{LTOP} \longrightarrow \Catinf \]
determines a functor
\[ \mathsf{Geo}^{op}\times\ltop\longrightarrow \catinfh \]
carrying a pair $(\geo,\xtop)$ to the \infcat $\str_{\geo}(\xtop)$. We have a subfunctor $\strloc_{\_}(\_)$ carrying the pair $(\geo,\xtop)$ to the subcategory $\strloc_{\geo}(\xtop)$ on the local morphisms. We let $q_{\mathsf{Geo}}:\overline{\ltop}_{\mathsf{Geo}}\rightarrow \mathsf{Geo}^{op}\times \ltop$ denote a coCartesian fibration associated to the functor $\strloc_{\_}(\_)$. For a fixed geometry $\geo$, we denote by $\ltop(\geo)$ the pullback $\{\geo\}\times_{\mathsf{Geo}^{op}\times\ltop}\overline{\ltop}_{\mathsf{Geo}}$ which comes equipped with a coCartesian fibration \[q_{\geo}:\ltop(\geo)\longrightarrow  \ltop.\] Replacing $\mathsf{Geo}$ with $\mathsf{PGeo}$, we have a coCartesian fibration 
$q_{\mathsf{PGeo}}:\overline{\ltop}_{\mathsf{PGeo}}\rightarrow \mathsf{PGeo}^{op}\times \ltop$, which restricts to coCartesian fibration \[q_{\pregeo}:\ltop(\pregeo)\longrightarrow  \ltop\] over each pregeometry $\pregeo$.
\end{defn}

\begin{rmk}
We have a non-full subcategory inclusion $\ltop\hookrightarrow \widehat{\mathsf{Cat}}_{\infty}$ which classifies a coCartesian fibration over $\ltop$, the \emph{universal topos fibration}, which was denoted as $q:\overline{\ltop}\rightarrow \ltop$. For $\icat\rightarrow\icatd$ a coCartesian fibration associated to a functor $F:\icatd\rightarrow\catinf$, there is a canonical coCartesian equivalence between the unstraightening of the composition
\[ \icatd\overset{F}{\longrightarrow}\catinf\overset{\fun(\icate,\_)}{\longrightarrow} \catinf\]
and the left vertical coCartesian fibration in the pullback diagram 
\[
\begin{tikzcd}
  \fun(\icate,\icat)\times_{\fun(\icate,\icatd)}\icatd\ar[d] \ar[r] & \fun(\icate,\icat)\ar[d] \\
  \icatd \ar[r] & \fun(\icate,\icatd)
\end{tikzcd}
\]
where the lower horizontal map is adjoint to the obvious projection $\icatd\times\icate\rightarrow\icatd$. It follows that for $\mathcal{T}$ a pregeometry, we can identify $\ltop(\mathcal{T})$ with the subcategory of $\fun(\mathcal{T},\overline{\ltop})\times_{\fun(\mathcal{T},\ltop)}\ltop$ defined as follows:
\begin{enumerate}[(1)]
    \item An object of $\fun(\mathcal{T},\overline{\ltop})\times_{\fun(\mathcal{T},\ltop)}\ltop$, which we can identify with a pair $(\xtop,\Of_{\xtop})$, for some $\infty$-topos $\xtop$ and a functor $\Of_{\xtop}:\mathcal{T}\rightarrow \xtop$, lies in $\ltop(\mathcal{T})$ precisely if $\Of_{\xtop}$ is a $\mathcal{T}$-structure on $\xtop$.
    \item For a morphism $\alpha:(\xtop,\mathcal{O}_{\xtop})\rightarrow (\ytop,\mathcal{O}'_{\ytop})$ in $\fun(\mathcal{T},\overline{\ltop})\times_{\fun(\mathcal{T},\ltop)}\ltop$, let $f^*:\xtop\rightarrow \ytop$ be the underlying algebraic morphism. We declare $\alpha$ to be a morphism in $\ltop(\mathcal{T})$ if the induced morphism $f^*\circ \mathcal{O}_{\xtop}\rightarrow \mathcal{O}'_{\ytop}$ is a local morphism of $\mathcal{T}$-structures on $\ytop$. 
\end{enumerate}
\end{rmk}
\begin{defn}
We call the opposite category of $\ltop(\mathcal{T})$ the $\infty$-category of \emph{$\mathcal{T}$-structured spaces}, and denote it $\rtop(\mathcal{T})$. Similarly, we have the category $\rtop(\mathcal{G})$ of $\mathcal{G}$-structured spaces.
\end{defn}
\begin{rmk}\label{rmk:exponentiate}
Let $\geo$ be a discrete geometry. For $f^*:\xtop\rightarrow \ytop$ an algebraic morphism, we can identify the coCartesian pushforward associated to $f^*$ with the functor 
\[ \fun^{\mathrm{lex}}(\geo,\xtop) \longrightarrow  \fun^{\mathrm{lex}}(\geo,\ytop) \]
given by composition with $f^*$. This functor has a right adjoint, given by composition with $f_*$. It follows that $p_{\geo}$ is also a Cartesian fibration. Under the identification $\shv_{\ind(\geo^{op})}(\xtop)\simeq\fun^{\lex}(\geo,\xtop)$, the functor $f_*$ corresponds to the functor $\fun^{\mathrm{R}}(\ytop^{op},\ind(\geo^{op}))\rightarrow\fun^{\mathrm{R}}(\xtop^{op},\ind(\geo^{op}))$ given by composition with $f^*:\xtop^{op}\rightarrow \ytop^{op}$. This identification yields the following alternative construction of the fibration $q_{\geo}:\ltop(\geo)\rightarrow\ltop$, which we will find useful later on. The subcategory inclusion $\rtop\subset\catinf$ of \inftopoi and geometric morphisms between them determines via unstraightening a biCartesian fibration $\overline{\rtop}^{op}\rightarrow\ltop$ whose fibre over $\xtop\in\ltop$ is an \infcat equivalent to $\xtop^{op}$. Define a simplicial set $\widetilde{\ltop}(\ind(\geo^{op}))$ over $\ltop$ by the following universal property: for any map of simplicial sets $K\rightarrow\ltop$, there is a canonical bijection 
\[ \Hom_{(\sset)_{/\ltop}}(K,\widetilde{\ltop}(\ind(\geo^{op})))\cong \Hom_{\sset}(K\times_{\ltop}\overline{\rtop}^{op},\ind(\geo^{op})). \]
It follows from \cite{HTT}, Corollary 3.2.2.12 that the map $q_{\ind(\geo^{op})}:\widetilde{\ltop(\ind(\geo^{op}))}\rightarrow\ltop$ is a biCartesian fibration; in particular, it is a Cartesian fibration associated to the functor \[\left(f^*:\xtop\rightarrow\ytop\right)\longmapsto\left(\_\circ f^*:\fun(\ytop^{op},\ind(\geo^{op}))\rightarrow\fun(\xtop^{op},\ind(\geo^{op}))\right).\]
 Let $\ltop(\ind(\geo^{op}))\subset \widetilde{\ltop(\ind(\geo^{op}))}$ be the full subcategory spanned by functors $\xtop^{op}\rightarrow \ind(\geo^{op})$ that preserve small limits, then $q_{\ind(\geo^{op})}$ restricts to a Cartesian fibration $q_{\ind(\geo^{op})}:\ltop(\ind(\geo^{op}))\rightarrow\ltop$ which is naturally equivalent to $q_{\geo}:\ltop(\geo)\rightarrow\ltop$.  
\end{rmk}

\subsubsection{Spectra and schemes}
We first describe the global sections functor for $\geo$-structured \inftopoit. 
\begin{prop}\label{prop:fibreinclusionadjoint}
Let $p:\icat\rightarrow \icatd$ be a Cartesian fibration, and consider for each object $C\in\icat$ the induced Cartesian fibration $p_C:\icat_{C/}\rightarrow \icatd_{p(C)/}$ (\cite{HTT}, Proposition 2.4.3.2). Then the inclusion of the fibre $(\icat_{p(C)})_{C/}\hookrightarrow \icat_{C/}$ has a right adjoint, and a map $\eta:C''\rightarrow C'$ for which $p(C'')=p(C)$ is a counit transformation if and only if $\eta$ is a $p$-Cartesian lift terminating at $C'$. 
\end{prop}
\begin{proof}
We should show that for each morphism $f:C\rightarrow C'$ of $\icat_{C/}$, the right fibration \[(\icat_{p(C)})_{C/}\times_{\icat_{C/}}\icat_{C//C''}\rightarrow (\icat_{p(C)})_{/C}\]
is representable, that is, the \infcat $(\icat_{p(C)})_{C/}\times_{\icat_{C/}}\icat_{C//C'}$ has an final object. We have a diagram of simplicial sets 
\[
\begin{tikzcd}
\{f\}\times_{\icatd_{p(C)//p(C')}}\icat_{C//C'} \ar[d] \ar[r] & (\icat_{p(C)})_{C/}\times_{\icat_{C/}}\icat_{C//C'} \ar[d] \ar[r] & \icat_{C//C'} \ar[d] \\
\{f\} \ar[r,"i"] & \Hom^{\mathrm{R}}_{\icatd_{p(C)/}}(p(C),p(C')) \ar[r] & \icatd_{p(C)//p(C')}
\end{tikzcd}
\]
where both squares are pullbacks and the indicated map $i$ is a homotopy equivalence of Kan complexes and thus a categorical equivalence. As all objects in the diagram are fibrant and all vertical maps are categorical fibrations, the upper left horizontal map is also a categorical equivalence. But a final object in the \infcat $\{f\}\times_{\icatd_{p(C)//p(C')}}\icat_{C//C'}$ is exactly a $p$-Cartesian lift of the map $p(C)\rightarrow p(C')$ terminating at $C'$.
\end{proof}
\begin{cons}\label{cons:globalsectionsgeo}
Let $\geo$ be a discrete geometry. We can identify, for each \inftop $\xtop$, the \infcat $\strloc_{\geo}(\xtop)$ with the \infcat of right adjoint functors $\fun^{\mathrm{R}}(\mathrm{Pro}(\geo),\xtop)$, so we can identify the fibre at $\spa$ with the ind-completion $\mathrm{Ind}(\geo^{op})$. As $\geo$ is discrete, the coCartesian fibration $q_{\geo}:\ltop(\geo)\rightarrow \ltop$ is also a Cartesian fibration. Let $\emptyset \in \mathrm{Ind}(\geo^{op})$ be an initial object in $\mathrm{Ind}(\geo^{op})$, then the object $(\spa,\emptyset)\in\ltop(\geo)$ is initial by \cite{HTT}, Proposition 4.3.1.10 and 4.3.1.5 and the fact that $\spa$ is initial in $\ltop$, so invoking Proposition \ref{prop:fibreinclusionadjoint}, we have a functor
\[ \Gamma: \ltop(\geo) \longrightarrow \mathrm{Ind}(\geo^{op})\]
characterized by the property that for each \inftop $\xtop$, the functor $\Gamma|_{p_{\geo}^{-1}(\xtop)}:\shv_{\mathrm{Ind}(\geo^{op})}(\xtop)\rightarrow \mathrm{Ind}(\geo^{op})$ may be canonically identified with the global sections functor given by composing with $\spa\rightarrow\xtop$. Let $\geo$ be an arbitrary geometry and consider the canonical transformation of geometries $f:\geo_{\mathrm{disc}}\rightarrow\geo$. We define the \emph{global sections functor} for $\geo$-structured \inftopoi as the composition
\[ \Gamma:\ltop(\geo)\overset{f^*}{\longrightarrow}\ltop(\geo_{\mathrm{disc}})\overset{\Gamma}{\longrightarrow}  \mathrm{Ind}(\geo^{op}) \]
so that for each \inftop $\xtop$, we may canonically identify $\Gamma|_{p_{\geo}^{-1}(\xtop)}$ with the functor
\[ \strloc_{\geo}(\xtop) \subset \strloc_{\geo_{\mathrm{disc}}}(\xtop) \simeq \shv_{\mathrm{Ind}(\geo^{op})} \overset{\Gamma}{\longrightarrow} \mathrm{Ind}(\geo^{op}).  \]
\end{cons}
The global sections functor admits a left adjoint, the \emph{$\geo$-spectrum} functor. In fact, as shown by Lurie, this spectrum functor is a consequence of a much more general construction, the \emph{relative spectrum} functor associated to an arbitrary transformation of geometries. The relative spectrum is an indispensable tool for comparing different geometries, as we will do later in this work.
\begin{thm}[Relative spectrum]\label{thm:relativespectrum}
Let $f:\geo\rightarrow\geo'$ be transformation of geometries, then the functor $f^*:\ltop(\geo')\rightarrow \ltop(\geo)$ admits a left adjoint $\spec^{\geo'}_{\geo}$.   \end{thm}

The preceding theorem may be phrased somewhat more globally. Let $\overline{\mathsf{Geo}}\rightarrow \mathsf{Geo}$ be the coCartesian fibration associated to $\mathsf{Geo}\hookrightarrow\catinf$, and define a simplicial set $\overline{\ltop}\mathsf{Geo}$ over $\mathsf{Geo}$ by the universal property that for any map of simplicial sets $K\rightarrow \mathsf{Geo}$, there is a canonical bijection 
\[ \Hom_{(\sset)_{/\mathsf{Geo}}}(K,\overline{\ltop}\mathsf{Geo}) \cong \Hom_{\sset}(K\times_{\mathsf{Geo}}\overline{\mathsf{Geo}},\overline{\ltop}) . \]
It follows from \cite{HTT}, Cor. 3.2.2.12 that the canonical map $\overline{\ltop}\mathsf{Geo}\rightarrow \mathsf{Geo}$ is a Cartesian fibration. In a similar fashion, we have Cartesian fibration $\ltop\mathsf{Geo}\rightarrow \mathsf{Geo}$ which admits a canonical map from $\ltop\times\mathsf{Geo}$; we define a Cartesian fibration over $\mathsf{Geo}$ as the pullback of Cartesian fibrations
\[
\begin{tikzcd}
\ltop\mathsf{Geo} \ar[d]\ar[r] &     \overline{\ltop}\mathsf{Geo}\ar[d] \\
\ltop\times\mathsf{Geo} \ar[r] & \ltop\mathsf{Geo}.
\end{tikzcd}
\]
Note that the fibre over a geometry $\geo$ is equivalent to the \infcat $\ltop(\geo)$, so we will denote objects of $\ltop\mathsf{Geo}$ by triples $(\geo,\xtop,\Of_{\xtop})$. Unwinding the definitions, the unstraightening of the Cartesian fibration $q_{\mathsf{Geo}}:\ltop\mathsf{Geo}\rightarrow\mathsf{Geo}$ is the functor $\mathsf{Geo}^{op}\rightarrow\catinfh$ that carries a transformation $f:\geo\rightarrow \geo'$ of geometries to the functor $f^*:\ltop(\geo')\rightarrow\ltop(\geo)$. Theorem \ref{thm:relativespectrum} therefore implies the following.
\begin{cor}
The Cartesian fibration $q_{\mathsf{Geo}}:\ltop\mathsf{Geo}\rightarrow\mathsf{Geo}$ is also a coCartesian fibration.
\end{cor}
\begin{rmk}
Given a functor $p:\icat\rightarrow \icatd\times\icate$, the following are equivalent. 
\begin{enumerate}[$(1)$]
    \item The composition $p_{\icatd}:\icat\rightarrow \icatd$ is a Cartesian fibration, $p$ carries $p_{\icatd}$-Cartesian morphisms to equivalences in $\icate$ and for each $D\in \icatd$, the restricted functor $p_{\icate,D}:\icat_{D}\rightarrow \icate\times\{D\}$ is a coCartesian fibration. 
    \item The composition $p_{\icat}:\icat\rightarrow\icate$ is a coCartesian fibration, $p$ carries $p_{\icate}$-Cartesian morphisms to equivalences in $\icatd$ and for each $E\in \icate$, the restricted functor $p_{E,\icatd}:\icat_E\rightarrow \icatd\times\{E\}$ is a Cartesian fibration.
\end{enumerate}
From the description of the $q_{\mathsf{Geo}}$-Cartesian morphism, we see that $(1)$ is satisfied for the functor $\ltop\mathsf{Geo}\rightarrow \ltop\times\mathsf{Geo}$. It follows that the functor $q_{\ltop}:\ltop\mathsf{Geo}\rightarrow\ltop$ is a coCartesian fibration. 
\end{rmk}
Here is a feature of the relative spectrum we will use on a few occasions. Let $f:\geo\rightarrow\geo'$ be a transformation of geometries and let $\ofxtop\in \ltop(\geo)$ and $(\xtop',\Of_{\xtop'}')\in \ltop(\geo')$. We say that a map $\ofxtop\rightarrow f^*(\xtop',\Of_{\xtop'}')=(\xtop',\Of_{\xtop'}'\circ f)$ \emph{exhibits $(\xtop',\Of_{\xtop'}')$ as a relative spectrum of $\ofxtop$} if this map exhibits a unit transformation at $\ofxtop$. In other words, if the corresponding map $\ofxtop \rightarrow (\xtop',\Of_{\xtop'}')$ in the \infcat $\ltop\mathsf{Geo}$ is $q_{\mathsf{Geo}}$-coCartesian. 
\begin{prop}\label{prop:relspecpushout}
Let $f:\geo\rightarrow\geo'$ be a transformation of geometries and let $\ofxtop\rightarrow f^*(\xtop',\Of_{\xtop'}'\circ f)$ be a map that exhibits a relative spectrum. Suppose that the underlying map on \inftopoi fits into pushout diagram 
\[
\begin{tikzcd}
\xtop \ar[r] \ar[d,"g^*"] & \xtop'\ar[d,"{g'}^{*}"] \\
\ytop \ar[r] & \ytop'
\end{tikzcd}
\]
in $\ltop$. Then the canonical map $(\ytop,g^*\circ \Of_{\xtop})\rightarrow (\ytop',{g'}^{*}\circ \Of'_{\xtop'}\circ f)$ of $\geo$-structured \inftopoi exhibits a relative spectrum.
\end{prop}
\begin{proof}
The data in the statement of the proposition correspond to a diagram 
\[
\begin{tikzcd}
(\geo,\xtop,\Of_{\xtop}) \ar[r] \ar[d] & (\geo',\xtop',\Of_{\xtop'}') \ar[d] \\
(\geo,\ytop,\Of_{\ytop}) \ar[r] & (\geo',\ytop',\Of_{\ytop'}')
\end{tikzcd}
\]
in $\ltop\mathsf{Geo}$ for which the upper horizontal map is $q_{\mathsf{Geo}}$-coCartesian, and both vertical maps are coCartesian for the fibration $q_{\ltop}:\ltop\mathsf{Geo}\rightarrow \ltop$. It follows that the diagram is a $q_{\ltop}$-pushout. Since the underlying square of \inftopoi is a pushout, the diagram is in fact a pushout diagram in $\ltop\mathsf{Geo}$. Since the underlying diagram of geometries is also a pushout as the vertical maps induce identity maps on the underlying geometries, the diagram is also a $q_{\mathsf{Geo}}$-pushout. Since the upper horizontal map is $q_{\mathsf{Geo}}$-coCartesian, the lower horizontal map is $q_{\mathsf{Geo}}$-coCartesian as well. 
\end{proof}
\begin{rmk}
In fact, one proves Theorem \ref{thm:relativespectrum} via the preceding proposition which allows one to replace an arbitrary $\geo$-structured \inftop by the classifying \inftop equipped with its universal $\geo$-structure, for which the relative spectrum can be explicitly constructed.
\end{rmk}
For $\geo$ a geometry, we let $\spec^{\geo}$ denote the composition
\[ \ind(\geo) \longrightarrow \ltop(\geo_{\mathrm{disc}}) \overset{\spec^{\geo}_{\geo{\mathrm{disc}}}}{\longrightarrow} \ltop(\geo), \]
where the first map is the inclusion of the fibre over the initial \inftopt. By construction, $\spec^{\geo}$ is a left adjoint to $\Gamma$; we call it the \emph{$\geo$-spectrum functor}. As shown by Lurie, this functor admits a much more concrete description.
\begin{cons}[$\geo$-spectrum]\label{spectrumfunctor}
Let $\geo$ be a geometry. We say that a morphism $f:U\rightarrow X$ in $\mathrm{Pro}(\geo)$ is \emph{admissible} if $f$ fits into a pullback diagram
\[
\begin{tikzcd}
U\ar[r]\ar[d,"f"] & j(U')\ar[d,"j(f')"] \\
X \ar[r] & j(X')
\end{tikzcd}
\]
where $f':U'\rightarrow X'$ is an admissible morphism in $\geo$ and $j:\geo\rightarrow \mathrm{Pro}(\geo)$ is the Yoneda embedding. Lemma 2.2.4 of \cite{dagv} tells us that 
\begin{enumerate}[(1)]
    \item Every equivalence in $\mathrm{Pro}(\geo)$ is admissible.
    \item The collection of admissible morphism is closed under taking pullbacks along any morphism in $\mathrm{Pro}(\geo)$.
    \item For a commutative diagram 
\begin{center}
\begin{tikzcd}
A \ar[dr,"f"'] \ar[rr,"h"] && C \ar[dl,"g"]\\
&  B 
\end{tikzcd}    
\end{center}
with $g$ admissible, $f$ is admissible if and only if $h$ is admissible. 
\end{enumerate}
Let $\mathrm{Pro}(\geo)^{\mathrm{ad}}_{/X}$ be the full subcategory of $\mathrm{Pro}(\geo)_{/X}$ spanned by admissible morphisms to $X$. This $\infty$-category is essentially small, as all admissible maps are pullbacks of admissible maps in $\geo$. We make $\mathrm{Pro}(\geo)^{\mathrm{ad}}_{/X}$ a site by endowing it with the Grothendieck topology generated by the covering families of the form $\{j'(V_{\alpha})\times_{j(U')}U\rightarrow U\}_{\alpha \in I}$ for $U\rightarrow j(U')$ some morphism and $\{V_{\alpha}\rightarrow U'\}_{\alpha \in I}$ an admissible covering of $U'$ in $\geo$.\\
Let $\mathrm{Spec}\,X$ denote the $\infty$-topos $\shv(\mathrm{Pro}(\geo)^{\mathrm{ad}}_{/X})$. We define a $\geo$-structure on $\mathrm{Spec}\,X$ by assigning to an object $Y\in \geo$ the sheafification of the presheaf sending an admissible map $U\rightarrow X$ to $\colim_{i\in I}\, \Hom_{\geo}(U_i,Y)$, where $\{U_i\}_{i\in I}$ is a filtered diagram with colimit $U$ in $\geo^{op}$. More formally, we have a functor
\[\geo \times \mathrm{Pro}(\geo)^{\mathrm{ad}}_{/X} \longrightarrow  \geo \times \mathrm{Pro}(\geo) = \geo \times \fun^{\mathrm{lex}}(\geo,\mathcal{S})\overset{\mathrm{ev}}{\longrightarrow} \mathcal{S}  \]
where the last functor is the evaluation pairing. By adjunction we obtain a functor $\widetilde{\Of_{\speco\,X}}:\geo\rightarrow \mathsf{PShv}(\mathrm{Pro}(\geo)^{\mathrm{ad}}_{/X})$. Let $\Of_{\mathrm{Spec}\,X}$ be the composition
\[\Of_{\mathrm{Spec}\,X}: \geo\overset{\rho}{\longrightarrow} \mathsf{PShv}(\mathrm{Pro}(\geo)^{\mathrm{ad}}_{/X}) \overset{L}{\longrightarrow} \shv(\mathrm{Pro}(\geo)^{\mathrm{ad}}_{/X}),  \]
where $L$ is a sheafification functor. The functor $\Of_{\mathrm{Spec}\,X}$ is indeed a $\geo$-structure on $\mathrm{Spec}\,X$: $\widetilde{\Of_{\speco\,X}}$ is manifestly left exact and $L$ is a left exact localization, so $\Of_{\mathrm{Spec}\,X}$ is also left exact. It's easy to see that $\Of_{\mathrm{Spec}\,X}$ is local, that is, for an admissible covering $\{V_{\beta}\rightarrow Y\}_{\beta\in J}$ the map $\coprod_{\beta}\Of_{\mathrm{Spec}\,X}(V_{\beta})\rightarrow \Of_{\mathrm{Spec}\,X}(Y)$ is an effective epimorphism. The natural transformation $\mathrm{id}\rightarrow L$ from the identity on $\mathsf{PShv}(\mathrm{Pro}(\geo)^{\mathrm{ad}}_{/X})$ to the sheafification induces a morphism of $\geo_{\mathrm{disc}}$-structured \inftopoi 
\[ (\mathsf{PShv}(\mathrm{Pro}(\geo)^{\mathrm{ad}}_{/X}),\widetilde{\Of_{\speco\,X}})\longrightarrow (\mathsf{Shv}(\mathrm{Pro}(\geo)^{\mathrm{ad}}_{/X}),\Of_{\speco\,X}).  \]
We can identify the global sections of $\Gamma(\widetilde{\Of_{\speco\,X}})$ with $X$, so that we have a morphism
\[ X\longrightarrow \Gamma(\Of_{\speco\,X})  \]
in $\ind(\geo^{op})$.
\end{cons}
The following is Theorem 2.2.12 of \cite{dagv}.
\begin{thm}[Absolute spectrum]\label{thm:spectrum}
Let $\geo$ be a geometry and let $X\in \mathrm{Pro}(\geo)=\ind(\geo^{op})^{op}$, then the morphism  
\[ X\longrightarrow \Gamma(\Of_{\speco\,X})  \]
exhibits a unit transformation for the global sections functor.
\end{thm}
\begin{cor}\label{cor:0truncspec}
Let $X\in \mathrm{Pro}(\geo)$ and suppose that for every object $U\rightarrow X$ in $\mathrm{Pro}(\geo)^{\mathrm{ad}}_{/X}$, the object $U$ is $n$-truncated as an object of the opposite category $\mathrm{Ind}(\geo^{op})$, then the structure sheaf $\Of_{\speco\,X}$ takes $n$-truncated values in $\speco\,X$.
\end{cor}
\begin{proof}
The functor $\widetilde{\Of_{\speco\,X}}$ carries $Y\in \geo$ to the presheaf that assigns to $U\in \mathrm{Pro}(\geo)^{\mathrm{ad}}_{/X}$ the space $\Hom_{\mathrm{Pro}(\geo)}(U,Y)=\Hom_{\mathrm{Ind}(\geo^{op})}(Y,U)$, where we view $Y$ as an object of $\mathrm{Pro}(\geo)$. This space is $n$-truncated by assumption. Since sheafification carries $n$-truncated objects to $n$-truncated objects, we conclude.  
\end{proof}
\begin{nota}
Let $X\in \mathrm{Pro}(\geo)$, then we denote the underlying \inftop $\mathsf{Shv}(\mathrm{Pro}(\geo)^{\mathrm{ad}}_{/X})$ of $\spec^{\geo}X$ by $\speco\,X$.
\end{nota}
\begin{ex}
Let $A$ be a commutative $k$-algebra, viewed as an object of $\mathrm{Pro}(\geo_{\mathrm{Zar}}(k))$, then for $\spec^{\geo_{\mathrm{Zar}}(k)}=(\xtop,\Of_{\xtop})$, the \inftop $\shv(\ofxtop)$ is the $0$-localic \inftop of sheaves on the prime spectrum $\speco'\,A$ of $A$, and $\Of_{\xtop}$ is the tautological structure sheaf carrying a basic open $D_a\subset \speco'\,A$ associated to a localization $A\rightarrow A[a^{-1}]$ to $A[a^{-1}]$.   
\end{ex}

\begin{ex}
Let $A$ be a commutative $k$-algebra, viewed as an object of $\mathrm{Pro}(\geo_{\mathrm{Zar}}(k))$, then $\mathrm{Pro}(\geo_{\mathrm{Zar}}(k))^{\mathrm{ad}}_{/A}$ can be identified with the nerve of the partially ordered set of localizations $A\rightarrow A[a^{-1}]$. The spectrum $\spec^{\geo_{\mathrm{Zar}}(k)}$ is the pair $(\xtop,\Of_{\xtop})$ with $\xtop$ the \inftop of sheaves on the prime spectrum $\speco'\,A$ of $A$, and $\Of_{\xtop}$ is the tautological structure sheaf carrying a basic open $D_a\subset \speco'\,A$ associated to a localization $A\rightarrow A[a^{-1}]$ to $A[a^{-1}]$. To see this, it is enough to show that $\shv(\mathrm{Pro}(\geo_{\mathrm{Zar}}(k))^{\mathrm{ad}}_{/A})$ is a spatial $0$-localic \inftopt, since this will imply that the spectrum functor $\spec^{\geo}$ takes values in topological spaces, which in turn guarantees it must coincide with the familiar Zariski spectrum, as both are left adjoint to the global sections functor. Since $\mathrm{Pro}(\geo_{\mathrm{Zar}}(k))^{\mathrm{ad}}_{/A}$ is a $(-1)$-category, the sheaf  \inftop is $0$-localic, and since this category has finite limits and the Zariski topology is finitary, $\tau_{\leq 0}\shv(\mathrm{Pro}(\geo_{\mathrm{Zar}}(k))^{\mathrm{ad}}_{/A})$ is coherent and thus has enough points by Deligne's completeness theorem.  
\end{ex}
\begin{rmk}
If $\geo$ is a geometry for which $\mathrm{Pro}(\geo)^{\mathrm{ad}}_{/X}$ is a $(-1)$-category equipped with a topology that is not finitary, Hochster duality fails and the $0$-localic \inftop $\shv(\mathrm{Pro}(\geo)^{\mathrm{ad}}_{/X})$ need not be spatial. This in fact occurs in our main examples of interest.    
\end{rmk}
\begin{ex}
Let $A$ be a commutative $k$-algebra, viewed as an object of $\mathrm{Pro}(\geo_{\et}(k))$, then the the underlying \inftop of the spectrum $\spec^{\geo_{\et}(k)}\,A$ is the (\inftop associated to the) \'{e}tale topos of $A$ and admits a description as $\shv(\mathsf{CAlg}^{0}_{k})^{op,\et}_{/A}$, the sheaves on the category of affine $k$-schemes \'{e}tale over $A$. Cast in this form, the structure sheaf is simply the composition \[(\mathsf{CAlg}^{0}_{k,\mathrm{fp}})^{op}\subset (\mathsf{CAlg}_k^{0})^{op}\longrightarrow\shv(\mathsf{CAlg}^{0}_{k})^{op,\et}\longrightarrow \shv(\mathsf{CAlg}^{0}_{k})^{op,\et}_{/A}\]   
where the second functor is the Yoneda embedding.
\end{ex}

\begin{defn}
A $\geo$-structured \inftop $\ofxtop$ is an \emph{affine $\geo$-scheme} if $\ofxtop$ is equivalent to $\spec^{\geo}X$ for some $X\in \ind(\geo^{op})$. We say that $\ofxtop$ is a \emph{finitely presented affine $\geo$-scheme} if $\ofxtop$ is equivalent to $\spec^{\geo}X$ for some $X\in \geo^{op}\subset \ind(\geo^{op})$   
\end{defn}
The spectrum functor allows us to discuss \emph{schemes} for a geometry. 
\begin{defn}
Let $\geo$ be a geometry. A morphism $h:\ofxtop\rightarrow\ofytop$ is an \emph{\'{e}tale morphism} of $\geo$-structured \inftopoi if the underlying morphism $\xtop\rightarrow \ytop$ is an \'{e}tale geometric morphism of \inftopoi and $h$ is $p_{\geo}$-coCartesian. We let $\ltop^{\et}(\geo)\subset\ltop(\geo)$ denote the subcategory on the \'{e}tale morphisms. Similarly, we have the subcategory $\rtop^{\et}(\geo)=(\ltop^{et}(\geo))^{op}\subset\rtop(\geo)$ on the \'{e}tale morphisms.  \end{defn}
\begin{nota}
Let $\ofxtop$ be a $\geo$-structured \inftop for some geometry $\geo$. The codomain fibration $\fun(\Delta^1,\xtop)\rightarrow\xtop$ induces via unstraightening an equivalence $\xtop\simeq \rtop^{\et}_{/\xtop}$ carrying $U\in \xtop$ to the geometric morphism $\xtop_{/U}\rightarrow\xtop$ forgetting the map to $U$. We have an associated \'{e}tale algebraic morphism $\pi^*_U:\xtop\rightarrow \xtop_{/U}$ (given by taking products with $U$). We denote the corresponding $\geo$-structure 
\[\geo\overset{\Of_X}{\longrightarrow}\xtop\overset{\pi_U^*}{\longrightarrow }\xtop_{/U}\]
by $\Of_{\xtop}|_U$. Of course, every \'{e}tale morphism of structured topoi is of the form $\ofxtop\rightarrow\ofxtopu$.
\end{nota}
For the reader's convenience, we repeat the statement of Proposition 2.3.5 of \cite{dagv}. The proof follows from basic properties of \'{e}tale algebraic morphisms of \inftopoi and coCartesian morphisms.
\begin{prop}\label{prop:etalemorphismprop}
Let $\geo$ be a geometry. Then the following hold true.
\begin{enumerate}[$(1)$]
    \item In a commuting triangle
    \[
    \begin{tikzcd}
     & \ofytop \ar[dr,"g"] \\
     \ofxtop \ar[rr,"h"] \ar[ur,"f"] && (\mathcal{Z},\Of_{\mathcal{Z}})
    \end{tikzcd}
    \]
    in $\ltop(\geo)$ for which $f$ is \'{e}tale, $g$ is \'{e}tale if and only if $h$ is \'{e}tale.
    \item The \infcat $\rtop^{\et}(\geo)$ admits small colimits which are preserved by the subcategory inclusion $\rtop^{\et}(\geo)\subset\rtop(\geo)$.
    \item The functor $\rtop^{\et}(\geo)_{/\ofxtop}\rightarrow \rtop^{\et}_{/\xtop}$ is a trivial Kan fibration. In particular, there is a canonical equivalence $\xtop\simeq \rtop^{\et}(\geo)_{/\ofxtop}$ of \infcatst. 
\end{enumerate}
\end{prop}
\begin{defn}
Let $\geo$ be a geometry. An object $\ofxtop\in\ltop$ is a \emph{$\geo$-scheme} if there exists a collection of objects $\{U_i\}_{i\in I}\subset\xtop$ such that $\coprod_iU_i\rightarrow 1_{\xtop}$ is an effective epimorphism and for each $i\in I$, the object $\ofxtopui$ is an affine $\geo$-scheme. A $\geo$-scheme is \emph{locally of finite presentation} if it is possible to choose the cover $\{U_i\}_i$ so that $\ofxtopui$ is an affine $\geo$-scheme of finite presentation. We have \infcats $\sch_{\lfp}(\geo)\subset \sch(\geo)\subset \rtop(\geo)$ of $\geo$-schemes and $\geo$-schemes locally of finite presentation.
\end{defn}
In Part II, we will develop some axiomatics around certain full subcategories $\mathcal{L}\subset \rtop(\geo)$ of \emph{scheme theories} that encapsulates some of the formal properties enjoyed by the category of affine $k$-schemes of finite presentation. For any such $\LS$, there is an associated \infcat $\sch(\geo,\LS)$ of schemes which has excellent formal properties; in particular, we will develop some tools that will allow for efficient comparison of different scheme theories. The flexibility our approach affords will pay off during the course of this work and its successors, since we will construct various theories of affine schemes for different ($\cinfty$-)geometries (with and without corners) and associated theories of higher geometric stacks that we wish to compare to one another.
\subsubsection{Truncation of structure sheaves}
Let $\geo$ be a geometry, then the \infcat $\str_{\geo}(\xtop)$ is a full subcategory of the \infcat $\str_{\geo_{\mathrm{disc}}}(\xtop)\simeq\shv_{\mathrm{Ind}(\geo^{op})}(\xtop)$. We have \emph{two} notions of what it means to have an $n$-truncated $\geo$-structure on $\xtop$, but these fortunatey coincide.
\begin{lem}\label{lem:truncatedsheaves}
Let $n\geq-1$ an integer. Let $\xtop$ be an \inftopt, let $\geo$ be an \infcat that admits finite limits and let $\geo\rightarrow\xtop$ be a left exact functor, then the following are equivalent.
\begin{enumerate}[$(1)$]
    \item The functor $\Of_{\xtop}$ takes $n$-truncated values in $\xtop$.
    \item The right adjoint functor $F:\xtop^{op}\rightarrow \mathrm{Ind}(\geo^{op})$ associated to $\Of_{\xtop}$ takes $n$-truncated values.
    \item Let $\geo_{\leq n}$ be the opposite of the $(n+1)$-category of compact objects in the $(n+1)$-category $\tau_{\leq n}\mathrm{Ind}(\geo^{op})$. Since $\tau_{\leq n}$ preserves compact objects, the composition   
\[ \geo^{op}\hooklongrightarrow\mathrm{Ind}(\geo^{op})\overset{\tau_{\leq n}}{\longrightarrow}  \tau_{\leq n}\mathrm{Ind}(\geo^{op})\]
 factors via $\geo_{\leq n}^{op}$ as a left exact functor $f:\geo\rightarrow \geo_{\leq n}$. Then the functor $\Of_{\xtop}$ factors as $\geo\overset{f}{\rightarrow}\geo_{\leq n}\overset{\Of_{\xtop}'}{\rightarrow }\xtop$ for a left exact functor $\Of_{\xtop}'$. 
\end{enumerate}
\end{lem}
\begin{proof}
Suppose that $\Of_{\xtop}$ factors as $\geo\rightarrow \tau_{\leq n}\xtop\subset\xtop$, then the pro-completed functor $\pro(\geo)\rightarrow\xtop$ also factors via $\tau_{\leq n}\xtop\subset\xtop$ as a small limit preserving functor since the full subcategory spanned by $n$-truncated objects is stable under limits. Invoking the adjoint functor theorem, we deduce that the left adjoint $G:\xtop\rightarrow \pro(\geo)$ to $\Of_{\xtop}$ factors as $\xtop\overset{\tau_{\leq n}}{\rightarrow}\tau_{\leq n}\xtop\overset{G'}{\rightarrow}\pro(\geo)$. Since the opposite functor $G':(\tau_{\leq n}\xtop)^{op}\rightarrow\ind(\geo^{op})$ preserves limits, its essential image lies in the full subcategory spanned by $n$-truncated objects, by \cite{HTT}, Proposition 5.5.6.16, and we conclude $(1)\Rightarrow (2)$. If $\xtop^{op}\rightarrow \ind(\geo^{op})$ takes values in $\tau_{\leq n}\ind(\geo^{op})$ then the adjoint functor $\pro(\geo)\rightarrow\xtop$ factors as a composition of limit preserving functors $\pro(\geo)\rightarrow\pro(\geo_{\leq n})\rightarrow \xtop$, so $\Of_{\xtop}$ factors via $\geo_{\leq n}$, which proves $(2)\Rightarrow (3)$. The assertion $(3)\Rightarrow (1)$ follows from \cite{HTT}, Proposition 5.5.6.16 again.
\end{proof}
\begin{defn}\label{defn:ntruncatedstr}
Let $\pregeo$ be a pregeometry and let $n\in\Z_{\geq -1}$. A functor $\Of:\pregeo\rightarrow \xtop$ preserving finite products and pullbacks along admissibles is \emph{$n$-truncated} if for each $X\in \pregeo$, the object $\Of(X)$ is $n$-truncated in $\xtop$. We let $\fun^{\pi\mathrm{ad},\leq n}(\pregeo,\xtop)\subset \fun^{\pi\mathrm{ad}}(\pregeo,\xtop)$ denote the full subcategory spanned by $n$-truncated objects. Similarly, we have full subcategories 
\[\str^{\leq n}_{\pregeo}(\xtop)\subset\str_{\pregeo}(\xtop),\quad \fun^{\pi\mathrm{ad},\loc,\leq n}(\pregeo,\xtop)\subset \fun^{\pi\mathrm{ad},\loc}(\pregeo,\xtop),\quad \str^{\loc,\leq n}_{\pregeo}(\xtop)\subset\strloc_{\pregeo}(\xtop)\]
of $n$-truncated objects. Similarly for a geometry $\geo$, we have a full subcategory $\fun^{\lex,\leq n}(\geo,\xtop)\subset\fun^{\lex}(\geo,\xtop)$ of $n$-truncated objects, and full subcategories
\[\str^{\leq n}_{\geo}(\xtop)\subset\str_{\geo}(\xtop),\quad \fun^{\lex,\loc,\leq n}(\geo,\xtop)\subset \fun^{\lex,\loc}(\geo,\xtop),\quad \str^{\loc,\leq n}_{\geo}(\xtop)\subset\strloc_{\geo}(\xtop).\]
\end{defn}
\begin{defn}\label{def:nstub}
For any geometry $\geo$, consider the left exact functor $f:\geo\rightarrow\geo_{\leq n}$ of Lemma \ref{lem:truncatedsheaves}. According to Remark \ref{rmk:geometrypushforward}, we can choose an $\mathsf{Forget}$-coCartesian lift of $f$, with $\mathsf{Forget}$ the forgetful functor $\mathsf{Geo}\rightarrow\catinf^{\lex,\mathrm{Idem}}$. The resulting geometry that we also denote $\geo_{\leq n}$ is the \emph{$n$-stub} of $\geo$.
\end{defn}
It follows immediately from the definition that composition with the transformation $f:\geo\rightarrow\geo_{\leq n}$ induces equivalences
\[\str^{\leq n}_{\geo}(\xtop)\simeq\str_{\geo_{\leq n}}(\xtop),\quad \fun^{\lex,\loc,\leq n}(\geo,\xtop)\simeq \fun^{\lex,\loc}(\geo_{\leq n},\xtop),\quad \str^{\loc,\leq n}_{\geo}(\xtop)\simeq \strloc_{\geo_{\leq n}}(\xtop).\]
The $n$-stub has further pleasant properties.
\begin{defn}
Let $\pregeo$ be a pregeometry, then a transformation of pregeometries $f:\pregeo\rightarrow \geo$ where $\geo$ is a geometry the underlying \infcat of which is $(n+1)$-truncated \emph{exhibits $\geo$ as an $n$-truncated geometric envelope} if for every $n$-truncated geometry $\geo'$, composition with $f$ induces an equivalence \
\[ \fun^{\mathrm{lex}}_{\tau}(\geo,\geo')\overset{\simeq}{\longrightarrow}\fun^{\pi\mathrm{ad}}_{\tau}(\pregeo,\geo). \]
\end{defn}
\begin{prop}\label{prop:truncatedgeoenv}
Let $\pregeo$ be a pregeometry, let $\geo$ be a geometric envelope for $\pregeo$ and let $\geo_{\leq n}$ be an $n$-stub of $\geo$ for $n\geq -1$, then the composition
\[\pregeo\hooklongrightarrow \geo\longrightarrow \geo_{\leq n} \]
exhibits an $(n+1)$-truncated geometric envelope.
\end{prop}
For any geometry $\geo$, we have a transformation of geometries $f:\geo\rightarrow\geo_{\leq n}$ and thus a relative spectrum $\spec^{\geo_{\leq n}}_{\geo}$. Beware that in general, this functor \emph{does not} coincide with the operation 
\[  \ofxtop \longmapsto (\xtop,\tau_{\leq n}\Of_{\xtop}).\]
This formula does not make sense for general (pre)geometries, since the functor $\tau_{\leq n}:\xtop\rightarrow\xtop$ does not preserve finite limits and thus may not take $\geo$-structures to $\geo$-structures; the relative spectrum $\spec^{\geo_{\leq n}}_{\geo}$ need not even be the identity on the underlying \inftopoit. For pregeometries however, composition with $\tau_{\leq n}:\xtop\rightarrow\xtop$ carries functors in $\fun^{\pi\mathrm{ad}}(\pregeo,\xtop)$ to itself if a mild condition is satisfied.
\begin{defn}
Let $\pregeo$ be a pregeometry and $n\geq -1$ an integer, then $\pregeo$ is \emph{compatible with $n$-truncations} if for any \inftop $\xtop$ and any functor $\Of:\pregeo\rightarrow \xtop$ preserving finite products and pullbacks along admissible maps, the following conditions are satisfied.
\begin{enumerate}[$(1)$]
\item The functor 
\[\pregeo\overset{\Of_{\xtop}}{\longrightarrow}\xtop\overset{\tau_{\leq n}}{\longrightarrow}\Of_{\xtop}\]
preserves finite products and pullbacks along admissible maps. 
\item The transformation $\Of_{\xtop}\rightarrow\tau_{\leq n}\circ\Of_{\xtop}$ in $\fun^{\pi\mathrm{ad}}(\pregeo,\xtop)$ is a local morphism. 
\end{enumerate} 
\end{defn}
\begin{rmk}
Let $\pregeo$ be a pregeometry, then $\pregeo$ is compatible with $n$-truncations if and only if for each \inftop $\xtop$, each $\pregeo$-structure $\Of\rightarrow \xtop$ and each admissible map $U\rightarrow X$ in $\pregeo$, the diagram
\[
\begin{tikzcd}
\Of(U)\ar[d]\ar[r] & \tau_{\leq n}\Of(U)\ar[d] \\
\Of(X) \ar[r] &\tau_{\leq n}\Of(X) 
\end{tikzcd}
\] 
is a pullback square in $\xtop$; this is the content of Proposition 3.3.3 of \cite{dagv}. This condition is satisfied if every admissible map in $\pregeo$ is $(n-1)$-truncated; this is the content of Proposition 3.3.5 of \cite{dagv}. We will use this criterion in the next section.
\end{rmk}

\begin{prop}\label{prop:pgeotruncfunctor}
Let $\pregeo$ be pregeometry that is compatible with $n$-truncations, then the following hold true.
\begin{enumerate}[$(1)$]
\item Let $\xtop$ be an \inftopt, then composition with $\tau_{\leq n}$ on $\xtop$ induces left adjoints  
\[\str^{\leq n}_{\pregeo}(\xtop)\rightarrow \str_{\pregeo}(\xtop),\quad \fun^{\pi\mathrm{ad},\loc,\leq n}(\pregeo,\xtop)\rightarrow \fun^{\pi\mathrm{ad},\loc}(\pregeo,\xtop),\quad \str^{\loc,\leq n}_{\pregeo}(\xtop)\rightarrow \strloc_{\pregeo}(\xtop)\]
to the inclusions of Definition \ref{defn:ntruncatedstr}.
\item  Let $\xtop$ be an \inftop and $\Of:\pregeo\rightarrow\xtop$ a functor preserving finite products and pullbacks along admissible maps, then $\Of$ is a $\pregeo$-structure if and only if $\tau_{\leq n}\Of$ is a $\pregeo$-structure. If $\Of':\geodiffder\rightarrow\xtop$ is another functor preserving finite products and pullbacks along admissible maps, then a morphism $f:\Of\rightarrow \Of'$ is local if and only if $\tau_{\leq n}\Of\rightarrow\tau_{\leq n}\Of'$ is local. That is, both commuting squares of \infcats in the diagram 
\[
\begin{tikzcd}
\str_{\pregeo}(\xtop) \ar[d,"\tau_{\leq n}"]\ar[r] & \fun^{\pi\mathrm{ad}}(\pregeo,\xtop) \ar[d,"\tau_{\leq n}"] & \fun^{\pi\mathrm{ad},\mathrm{loc}}(\pregeo,\xtop)\ar[d,"\tau_{\leq n}"]\ar[l]\\
\str^{\leq n}_{\pregeo}(\xtop)\ar[r] & \fun^{\pi\mathrm{ad},\leq n}(\pregeo,\xtop) & \fun^{\pi\mathrm{ad},\mathrm{loc},\leq n}(\pregeo,\xtop)\ar[l]
\end{tikzcd}
\]
are pullback squares.
\end{enumerate}
\end{prop}
\begin{proof}
The first assertion follows from \cite{dagv}, Proposition 3.3.3 and the second follows immediately from the assumption that $\Of\rightarrow\tau_{\leq n}\Of$ is a local morphism and Proposition \ref{prop:locmorstab}.
\end{proof}
\begin{prop}\label{prop:geotruncfunctor}
Let $\pregeo$ be a pregeometry compatible with $n$-truncations and let $f:\pregeo\hookrightarrow\geo$ be a geometric envelope for $\pregeo$, then the transformation $(\xtop,\Of_{\xtop})\rightarrow (\xtop,\tau_{\leq n}\Of_{\xtop})$ of \emph{$\pregeo$-structures} exhibits a unit transformation for the fully faithful embedding $\ltop(\geo_{\leq n})\hookrightarrow\ltop(\geo)$. In other words, the relative spectrum $\spec^{\geo_{\leq n}}_{\geo}$ can be identified with the composition of $\pregeo$-structures with the truncation functor $\tau_{\leq n}$ on the underlying \inftopoit.
\end{prop}
\begin{proof}
The truncation functor $\geo\rightarrow\geo_{\leq n}$ induces a functor $\iota:\ltop(\ind(\geo_{\leq n}^{op}))\rightarrow\ltop(\ind(\geo^{op}))$ which carries $q_{\ind(\geo_{\leq n}^{op})}$-Cartesian edges to $q_{\ind(\geo^{op})}$-Cartesian edges, by unstraightening the natural transformation 
\[ \fun^{\mathrm{R}}((\_)^{op},\ind(\geo_{\leq n}^{op}))\longrightarrow \fun^{\mathrm{R}}((\_)^{op},\ind(\geo^{op})) \]
of functors $\ltop^{op}\rightarrow\catinf$ given by composition with the inclusion $\ind(\geo_{\leq n}^{op})\simeq \tau_{\leq n}\ind(\geo^{op})\subset\ind(\geo^{op})$. In view of \cite{HA}, Proposition 7.3.2.6, it suffices to argue that $\iota$ admits a left adjoint on each fibre, which is an immediate consequence of Proposition \ref{prop:pgeotruncfunctor}. If $\ofxtop$ is a $\geo$-structure, then applying $\tau_{\leq n}:\xtop\rightarrow\xtop$ to $\Of_{\xtop}$ viewed as a $\pregeo$-structure yields a $\pregeo$-structure $\tau_{\leq n}\Of_{\xtop}$ and thus an $n$-truncated $\geo$-structure, that is, a $\geo_{\leq n}$-structure. For each $\geo_{\leq n}$-structured \inftop $\ofytop$ and each map $\ofxtop\rightarrow\ofytop$, both maps in the induced composition 
\[ \ofxtop\longrightarrow (\xtop,\tau_{\leq 0}\Of_{\xtop})\longrightarrow \ofytop  \]
are local morphisms by Proposition \ref{prop:pgeotruncfunctor}, so the left adjoint to $\iota$ restricts to a left adjoint to the functor $\ltop(\geo_{\leq n})\rightarrow\ltop(\geo)$.  
\end{proof}

\subsection{The geometry of finitely presented $C^{\infty}$-rings}\label{sectiongeometryofsmoothrings}
We have introduced a pregeometry $\diff$ of smooth manifolds. The goal in this section is to construct a \emph{$0$-truncated} geometric envelope for this pregeometry (although we will prove that it is one a bit later) in the guise of Lawvere's \emph{$\cinfty$-rings}. We take some time to study $C^{\infty}$-rings and $C^{\infty}$-schemes within the framework of geometries and structured spaces. While not strictly necessary, some familiarity with ordinary $C^{\infty}$-rings and synthetic differential geometry, as exposed, for instance, in the monograph \cite{MR} or the more recent \cite{Joy2}, will be advantageous to the reader.
\begin{defn}
Let $\mathsf{CartSp}\subset \diff$ be the full subcategory spanned by objects of the form $\R^n$ for $n\geq 0$. A \emph{$C^{\infty}$-ring} is an algebra for the Lawvere theory $\mathsf{CartSp}$, that is, a finite product preserving functor $\mathsf{CartSp}\rightarrow \mathsf{Set}$. The full subcategory of $\fun(\mathsf{CartSp},\mathsf{Set})$ spanned by $C^{\infty}$-rings is denoted $C^{\infty}\mathsf{ring}$. This is a strongly reflective, and thus presentable, subcategory of $\fun(\mathsf{CartSp},\mathsf{Set})$. 
\end{defn}
We will discuss Lawvere theories and their $\infty$-categories of space-valued algebras in more detail in Subsection \ref{sec:lawvere}. Unwinding the definition, a $C^{\infty}$-ring $A$ consists of a set $A(\R)$ equipped with a functional calculus for all smooth functions; that is, we have functorial operations 
\[ f_*:A(\R)^n\longrightarrow A(\R)^m\]
for each smooth map $f:\R^n\rightarrow \R^m$. We will often abuse notation and write `an element $a\in A$' for $a\in A(\R)$. Let $\pregeo^{disc}_{\R}$ be the category that has the same objects as $\mathsf{CartSp}$, but only polynomial maps; this is a Lawvere theory, and its algebras are precisely commutative $\R$-algebras. The transformation of Lawvere theories $\pregeo^{disc}_{\R}\rightarrow \mathsf{CartSp}$ induces an `underlying commutative $\R$-algebra' functor $(\_)^{\rmalg}:C^{\infty}\mathsf{ring}\rightarrow \mathsf{CAlg}^0_{\R}$. Many $C^{\infty}$-rings of interest are subsumed by the following examples.
\begin{enumerate}[(1)]
    \item The forgetful functor $(\_)^{\rmalg}:C^{\infty}\mathsf{ring}\rightarrow \mathsf{CAlg}^0_{\R}$ preserves limits and sifted colimits, and thus admits a left adjoint, the \emph{free $C^{\infty}$-ring functor} $F^{\cinfty}$ which takes the polynomial algebra $\R[x_1,\ldots,x_n]$ to $C^{\infty}(\R^n)$.
    \item Let $A$ be a $C^{\infty}$-ring and let $I\subset A$ be an ideal of the underlying $\R$-algebra. Then $A/I$ is a $C^{\infty}$-ring and the map $A\rightarrow A/I$ is regular epimorphism, i.e. it is the coequalizer in $C^{\infty}\mathsf{ring}$ of the equivalence relation determined by $I$. This is a consequence of Hadamard's lemma. 
    \item For a subset $X\subset \R^n$, the algebra of smooth functions 
    \[C^{\infty}(X):= \{f:X\rightarrow \R;\,\forall x\in X\text{ there exists } x\in U\subset \R^n\text{ open and } \tilde{f}\in C^{\infty}(U)\text{ such that }\tilde{f}|_{X\cap U}=f|_U\}\]
    is a $C^{\infty}$-ring by composition. If $X$ is a closed subset, then an application of Tietze's extension theorem shows that the natural map $C^{\infty}(\R^n)\rightarrow C^{\infty}(X)$ induces an isomorphism $C^{\infty}(\R^n)/\mathfrak{m}^{0}_X\rightarrow C^{\infty}(X)$, where $\mathfrak{m}^{0}_X$ is the ideal of functions that vanish on $X$. In particular (by the Whitney embedding theorem), the algebra of smooth functions on a manifold $M$ is a $C^{\infty}$-ring of the form $C^{\infty}(\R^n)/\mathfrak{m}^0_M$.   
    \item For $x\in \R^n$, the local algebra of germs of smooth functions at $x$ is a $C^{\infty}$-ring, given by $C^{\infty}(\R^n)_x:=C^{\infty}(\R^n)/{\mathfrak{m}_x}$ with $\mathfrak{m}^g_x$ the ideal of smooth functions vanishing in some neighbourhood of $x$.
    \item Every local Artin $\R$-algebra $W=\R\oplus \mathfrak{m}$ is a $C^{\infty}$-ring, whose $C^{\infty}$-ring structure is uniquely determined by the underlying algebra. Such $C^{\infty}$-rings are also called \emph{Weil algebras}.
    \item Let $\mathfrak{m}$ be the maximal ideal of the $C^{\infty}$-ring of germs $C^{\infty}(\R^n)_0$. The $C^{\infty}$-ring $J^k_n:=C^{\infty}(\R^n)_0/\mathfrak{m}^k$ of \emph{$k$'th order jets at $0$} is a Weil algebra.  
    \item Let $R$ be a complete local Noetherian $\R$-algebra with residue class field $\R$, then $R$ is of the form $R\cong\R[[x_1,\ldots,x_n]]/I$ for some ideal $I$, by Cohen's structure theorem. By Borel's lemma on formal power series, there is an equivalence $C^{\infty}(\R^n)/\mathfrak{m}^{\infty}_0\cong \R[[x_1,\ldots,x_n]]$, where $\mathfrak{m}^{\infty}_0$ is the ideal of functions that are flat at $0$ (all partial derivatives vanish at $0$). Thus, $R$ can be written as a quotient by $\mathfrak{m}^{\infty}_0$ of $C^{\infty}(\R^n)/\tilde{I}$, where $ \tilde{I}$ is a finitely generated (because $\R[[x_1,\ldots,x_n]]$ is Noetherian) ideal lifting $I$, so we conclude that $R$ is a $C^{\infty}$-ring. It's easy to see that all algebra morphism between complete local Noetherian $\R$-algebras are morphisms of $C^{\infty}$-rings, so the $C^{\infty}$-ring structure of $R$ is also uniquely determined by the underlying algebra, as in the case of Weil algebras (which this example subsumes). 
\end{enumerate}
\begin{rmk}
The essential image of the free $C^{\infty}$-ring functor $F^{\cinfty}:\mathsf{CAlg}^0_{\R}\rightarrow C^{\infty}\mathsf{ring}$ already contains many interesting objects. For instance, $F^{\cinfty}(\mathsf{CAlg}^0_{\R})$ contains all $C^{\infty}$-rings of smooth functions on \emph{compact} manifolds. This is an immediate consequence of the Nash-Tognoli theorem \cite{Nash,Tognoli}, which extends an older result of Whitney that all compact submanifolds of Euclidean space are diffeomorphic to zero loci of systems of real analytic equations.  
\end{rmk}
\begin{rmk}
Clearly, functions on manifolds that have less regularity also form $C^{\infty}$-rings: let $M$ be a manifold, then there are $C^{\infty}$-rings $C^k(M)$ of $k$-times differentiable functions and $\mathrm{Lip}^k(M)$ of $k$-times differentiable functions with locally Lipschitz derivatives. Let $M$ be an $n$-dimensional manifold and let $k\in \Q_{\geq0}$ and $p\in \Z_{>0}$ such that $kp>n$, then the space $W_{\mathrm{loc}}^{k,p}(M)$ of Sobolev functions of class $(k,p)$ is also a $\cinfty$-ring by an extension of the Sobolev multiplication theorems, which can be deduced from the Gagliardo-Nirenberg interpolation estimates (see, for instance, \cite{Melrose2}).
\end{rmk}
\begin{nota}
The functor $(\_)^{\rmalg}$ does not preserve pushouts nor coproducts in general. We reserve the symbol $\otimes^{\infty}$ for the pushout of $C^{\infty}$-rings.
\end{nota}
\begin{defn}
A $C^{\infty}$-ring $A$ is \emph{finitely generated} if $A\simeq C^{\infty}(\R^n)/I$ for some $n<\infty$. A is \emph{finitely presented} if $A\simeq C^{\infty}(\R^n)/I$ for some $n<\infty$ and $I$ a finitely generated ideal.
\end{defn}
\begin{rmk}
A $C^{\infty}$-ring $A$ is finitely presented if and only if the functor corepresented by $A$ (on the category of $C^{\infty}$-rings) preserves filtered colimits. $A$ is finitely generated if and only if the functor corepresented by $A$ preserves filtered colimits of diagrams consisting only of monomorphisms. See \cite{AR1}, chapter 3 for proofs of these facts. As the category of $C^{\infty}$-rings is presentable, we see that the full subcategories spanned by finitely generated and finitely presented $C^{\infty}$-rings have finite colimits. 
\end{rmk}
\begin{rmk}\label{rmk:transverse}
Let $f:N\rightarrow M$ and $g:P\rightarrow M$ be smooth maps of manifolds. We say that the pullback $N\times_M P$ is \emph{transverse} if for each $x_1\in N$, $x_2\in P$ such that $f(x_1)=x=p(x_2)$, the induced map $T_{x_1}f\oplus T_{x_2}g:T_{x_1}N\oplus T_{x_2}P\rightarrow T_xM$ is a surjection. An elementary but crucial result in synthetic differential geometry is the following: the functor $C^{\infty}:\diff\rightarrow C^{\infty}\mathsf{ring}^{op}$ is fully faithful, takes values in finitely presented objects, and preserves finite products and transverse pullbacks. For a proof, see \cite{MR}, chapter 1, theorem 2.8. The next chapter shall be concerned with proving a derived version of this result. 
\end{rmk}
\begin{rmk}
For $f:M\rightarrow \R^n$ a function on a manifold, we call the set $\mathrm{Carr}(f):=f^{-1}(\R^n\setminus\{0\})$ the \emph{carrier} of $f$, and we call the set $\mathrm{Supp}(f):=\overline{\mathrm{Carr}(f)}$ the \emph{support} of $f$. We will use frequently that any open set $U\rightarrow M$ in a manifold has a \emph{characteristic function} $\chi_U:M\rightarrow \R$, a function on $M$ such that $\mathrm{Carr}(\chi_U)=U$. We will also use that any function $f\in C^{\infty}(U)$ defined on an open set $U\subset M$ of a manifold $M$ is divisible by some function $g|_U$ where $g$ is defined on all of $M$ and nonzero on $U$. 
\end{rmk}
\begin{rmk}
Another cornerstone of our theory is the elementary lemma of Hadamard: for any smooth function $f:\R^n\rightarrow\R$ and any $p=(p_1,\ldots,p_n)\in \R^n$, there exists a collection of $n$ smooth function $\{g_i\}$ on $\R^n$ such that $f(\mathbf{x})-f(p)=\sum_{i=1}^n g_i(\mathbf{x})(x_i-p_i)$. All our tools for computing colimits in the 1-category $\cinfty\mathsf{ring}$ and the \infcat $\sring$ depend on it. 
\end{rmk}
For any $n>0$, the $C^{\infty}$-rings $C^{\infty}(\R^n)$ do not satisfy the the conclusion of the Nullstellensatz for arbitrary ideals; instead, we single out three classes of ideals for which the weak version of the Nullstellensatz does hold. Let $M$ be a smooth manifold of dimension $n>0$ and let $I$ be an ideal of the commutative algebra $C^{\infty}(M)$. For $x\in M$, we have the ideals
\begin{enumerate}[$(1)$]
    \item $\mathfrak{m}^0_x$ of functions that vanish at $x$, and the quotient map $C^{\infty}(M)\rightarrow C^{\infty}(M)/\mathfrak{m}^0_x\cong\R$ is the map $\mathrm{ev}_{x}$ evaluating at $x$.
    \item $\mathfrak{m}^{\infty}_x$ of functions that are flat at $x$, and choosing coordinates centered at $x$, the quotient map $C^{\infty}(M)\rightarrow C^{\infty}(M)/\mathfrak{m}^{\infty}_x\cong\R[[x_1,\ldots,x_n]]$ is the map $j^{\infty}_x$ taking the infinite jet at $x$.
    \item $\mathfrak{m}^{g}_x$ of functions that vanish in a neighbourhood of $x$, and choosing coordinates centered at $x$, the quotient map $C^{\infty}(M)\rightarrow C^{\infty}(M)/\mathfrak{m}^{g}_x\cong C^{\infty}(\R^n)_0$ is the map taking the germ at $x$.
\end{enumerate}
Since surjections of ring maps carry ideals to ideals, it makes sense to ask whether a function $f\in C^{\infty}(M)$ is pointwise, formally, or locally contained in an ideal $I$. 
\begin{defn}
Let $M$ be a smooth manifold of dimension $n>0$, and let $I\subset C^{\infty}(M)$ be an ideal. Write $Z(I)$ for the common zero locus of the functions in $I$.
\begin{enumerate}[(1)]
    \item $I$ is \emph{point determined} iff for all $f\in C^{\infty}(\R^n)$, $f\in I$ iff $f(x)=0$ for all $x\in Z(I)$.
    \item $I$ is \emph{jet determined} or \emph{closed} iff for all $f\in C^{\infty}(\R^n)$, $f\in I$ iff $j^{\infty}_x(f) \in j^{\infty}_x(I)$ for all $x\in Z(I)$, where $j^{\infty}_x:C^{\infty}(\R^n)\rightarrow \R[[x_1,\ldots,x_n]]$ carries a function to its formal power series at $x$.
    \item $I$ is \emph{locally determined} or \emph{germ determined} iff for all $f\in C^{\infty}(\R^n)$, $f\in I$ iff $f_x\in I_x$ for all $x\in Z(I)$, where $f_x$ and $I_x$ are the germ of $f$ at $x$ and the ideal of $C^{\infty}(M)_x$ of germs at $x$ of functions in $I$.
\end{enumerate}
\end{defn} 
\begin{rmk}
Here are some properties of the classes of ideals just defined.
\begin{enumerate}[(1)]\label{idealproperties}
    \item Point determined implies jet determined implies germ determined. None of these implications can be reversed in general. For instance, the ideal $I\subset\cinfty(\R)$ of functions whose jet at $0$ vanishes is jet determined but not point determined. An ideal of $\cinfty(\R)$ generated by a compactly supported function is germ determined but not jet determined. Finally, for an ideal that satisfies none of the conditions above -and for which the Nullstellensatz fails completely- consider the ideal of compactly supported functions in $\cinfty(\R)$. 
    \item A collection of functions $\{f_1,\ldots,f_m\}$ on $M$ generates a point determined ideal if the functions $\{f_1,\ldots,f_m\}$ are \emph{independent}, that is, the zero locus of $(f_1,\ldots,f_m):M\rightarrow \R^m$ consists of regular points. 
    \item Recall that given a collection of functions $\{f_{\alpha}\}\subset C^{\infty}(M)$ such that their carriers furnish a locally finite collection of opens on $M$, the pointwise sum $\sum_{\alpha}f_{\alpha}$ exists in $C^{\infty}(M)$ and is called a locally finite sum. An ideal $I\subset C^{\infty}(M)$ is germ determined if and only if $I$ is closed under taking locally finite sums. It follows easily that finitely generated ideals are germ determined. 
    \item Let $I\subset \R[[x_1,\ldots,x_n]]$ be an ideal, then in order to conclude that $h\in I$, it suffices to show that for all $k\in \Z_{>0}$, $h\in I+\mathfrak{m}^k$, where $\mathfrak{m}=(x_1,\ldots,x_n)$, the unique maximal ideal. Indeed, it suffices to show that $I=\cap_{k\geq 1} (I+\mathfrak{m}^k)$. The inclusion $I\subset \cap_{k\geq 1} (I+\mathfrak{m}^k)$ is obvious. For the other inclusion, it suffices to show that $p(\cap_{k\geq 1} (I+\mathfrak{m}^k))=0$, where $p:\R[[x_1,\ldots,x_n]]\rightarrow \R[[x_1,\ldots,x_n]]/I$ is the projection, but we clearly have
    \[ p(\bigcap_{k>1} (I+\mathfrak{m}^k))\subset \bigcap_{k\geq1}p(I+\mathfrak{m}^k)=\bigcap_{k\geq 1}p(\mathfrak{m})^k.  \]
    Since $p$ is a local morphism, $p(\mathfrak{m})$ is the maximal ideal of $\R[[x_1,\ldots,x_n]]/I$ so Krull's intersection theorem yields $\cap_{k\geq 1}p(\mathfrak{m})^k=0$ as $\R[[x_1,\ldots,x_n]]/I$ is local and Noetherian. 
    \item We say that a finitely generated $C^{\infty}$-ring $A=C^{\infty}(\R^n)/I$ is \emph{point determined} (\emph{closed}, \emph{germ determined}) if $I$ is a point determined (closed, germ determined) ideal. This does \emph{not} depend on the presentation of $A$. Thus, if $C^{\infty}(\R^n)/I\cong C^{\infty}(\R^m)/J$ and $I$ is point determined (closed, germ determined), then $J$ is point determined (closed, germ determined) as well. As an application, let $M$ be a manifold and note that as $M$ lies in some $\R^n$ as a closed submanifold, $C^{\infty}(M)$ is a point determined $C^{\infty}$-ring. As $C^{\infty}(M)$ is finitely presented, this shows that we have a presentation $C^{\infty}(M)\cong C^{\infty}(\R^n)/I$ where $I$ is a finitely generated and point determined ideal. 
\end{enumerate}
\end{rmk}
\begin{rmk}
Let $X\subset \R^n$ be a subset. We define the following ideals of $C^{\infty}(\R^n)$ associated to $X$:
\begin{align*}
    \mathfrak{m}^0_X:&=\{f\in C^{\infty}(\R^n)|\,f(p)=0\,\forall p\in X\},\\
    \mathfrak{m}^{\infty}_X:&=\{f\in C^{\infty}(\R^n);\,D_{\alpha}f(p)=0\,\forall p\in X\},\\
    \mathfrak{m}^g_X:&=\{f\in C^{\infty}(\R^n);\,\exists U\supset X \text{ open, }f|_U=0\}.
\end{align*}
In the second line, $D_{\alpha}$ denotes the differential operator $\del^{\alpha_1}_{x_1}\ldots \del^{\alpha_n}_{x_n}$ for $\alpha$ a multi-index $(\alpha_1,\ldots,\alpha_n)\in (\Z_{\geq 0})^n$. $\mathfrak{m}_X^0$ is point determined, $\mathfrak{m}_X^{\infty}$ is closed and $\mathfrak{m}_X^{g}$ is germ determined. If $X\subset \overline{X^{\circ}}$, then $\mathfrak{m}^0_X=\mathfrak{m}^{\infty}_X$. 
\end{rmk} 
As we have seen, Tietze's extension theorem shows that for $X\subset \R^n$ closed subset, we have $C^{\infty}(\R^n)/\mathfrak{m}_X^0\cong C^{\infty}(X)$. There is a similar characterization of $C^{\infty}$-rings of the form $C^{\infty}(\R^n)/\mathfrak{m}_X^{\infty}$ that uses Whitney's extension theorem. 
\begin{defn}\label{def:whitneyfunctions}
Let $X\subset U$ be a closed subset of an open subset in $\R^n$, and let $F=(f^k)_{k\in \Z_{\geq 0}^n}$ be a collection of continuous functions for each multi-index $k$. $F$ is a \emph{Whitney function} if for each $m\geq 0$, we have for $\mathbf{x}, \mathbf{y}\in X$
\[ f^k(\mathbf{x}) = \sum_{|l|=m-|k|} \frac{f^{l+k}(\mathbf{y})}{l!}(\mathbf{x}-\mathbf{y})^l + R^k(\mathbf{x},\mathbf{y}),  \]
where $R^k(\mathbf{x},\mathbf{y})$ is a term that goes to $0$ as $|\mathbf{x}-\mathbf{y}|\rightarrow 0$ faster than $|\mathbf{x}-\mathbf{y}|^{m-|k|}$. We let $C^{\infty}(X;U)$ denote the algebra of Whitney functions on $X\subset U$; it is a $\cinfty$-ring in a canonical fashion.
\end{defn}
The following easy lemma shows that if $X\subset\R^n$ is a closed quadrant, then the Whitney functions on $X$ coincide with the functions that have infinitely many derivatives up to the boundary.
\begin{lem}
Let $X\subset\R^n$ be a closed convex subset with nonempty interior. Then restriction to $X^{\circ}$ induces an equivalence between $\cinfty(X;\R^n)$ and the space
\[ \{f\in\cinfty(X^{\circ});\,D_{\alpha}f\text{ is bounded on }X^{\circ} \forall \alpha\in \Z_{\geq 0}^n\}. \]
\end{lem}
\begin{prop}[Whitney's Extension Theorem \cite{Whitney}]\label{prop:whitneyextension}
Let $X\subset U$ be a closed subset of an open subset in $\R^n$, then taking the infinite jet prolongation and restricting to $X$ yields an isomorphism $C^{\infty}(U)/\mathfrak{m}_X^{\infty}\cong C^{\infty}(X;U)$.
\end{prop}
A proof can be found in \cite{Malgrange}. We record the following pleasant property of flat ideals, i.e. ideals of the form $\mathfrak{m}^{\infty}_X$ for $X\subset \R^n$ closed.
\begin{thm}[Reyes-van Qu\^{e} \cite{reyesvanque}]\label{productflatideals}
Let $X\subset \R^n$ and $Y\subset\R^m$ be closed, then as ideals of $C^{\infty}(\R^{n+m})$ we have the equality $(\mathfrak{m}^{\infty}_X,\mathfrak{m}^{\infty}_Y)=\mathfrak{m}^{\infty}_{X\times Y}$.
\end{thm}
\begin{cor}
Let $X\subset\R^n$ and $Y\subset \R^m$ be closed subsets, then the canonical map 
\[ C^{\infty}(X;\R^n)\otimes^{\infty}C^{\infty}(Y;\R^m)\longrightarrow C^{\infty}(X\times Y;\R^{n+m}) \]
is an equivalence.
\end{cor}
\begin{rmk}
We will also prove a derived version of the result above (Theorem \ref{thm:reyesvanquederived}), which has important consequences for the deformation theory of closed $\cinfty$-rings. The result also shows that the local models for manifolds with corners behave well under the derived tensor product of $C^{\infty}$-rings, which is the starting point for derived $\cinfty$-geometry with corners.
\end{rmk}
\subsubsection{Locality of $\cinfty$-rings}
In this subsection, we define the admissible maps for a geometry with underlying category $C^{\infty}\mathsf{ring}^{op}_{\mathrm{fp}}$. These maps will correspond to open inclusions of \emph{$C^{\infty}$-schemes}. First, we state for convenience a special case of Remark \ref{rmk:transverse}.
\begin{lem}\label{lem:localizeopens}
Let $U\subset M$ be an open subset of a manifold $M$ and let $f:N\rightarrow M$ be a $\cinfty$ map of manifolds. Then the diagram 
\[
\begin{tikzcd}
 f^{-1}(U)\ar[d]\ar[r,hook] & N\ar[d,"f"]   \\
 U\ar[r,hook] & M
\end{tikzcd}
\]
induces a pushout diagram 
\[
\begin{tikzcd}
\cinfty(M)\ar[d]\ar[r] & \cinfty(U) \ar[d]\\
\cinfty(N) \ar[r] & \cinfty(f^{-1}(U)).
\end{tikzcd}
\]
of $\cinfty$-rings.
\end{lem} 
\begin{defn}
Let $A$ be a $C^{\infty}$-ring and let $a\in A$. A map $f:A\rightarrow B$ such that $f(a)$ is invertible \emph{exhibits $B$ as a localization of $A$} if for each $C^{\infty}$-ring $C$, composition with $f$ induces a bijection 
\[ \Hom_{C^{\infty}\mathsf{ring}}(B,C)\overset{\simeq}{\longrightarrow} \Hom^0_{C^{\infty}\mathsf{ring}}(A,C)\]
where $\Hom^0_{C^{\infty}\mathsf{ring}}(A,C)$ is the subset of maps that send $a$ to an invertible element.
\end{defn}
Localizations are unique up to unique isomorphism; given some $a\in A$, we denote the localization at this element by $A\rightarrow A[a^{-1}]$. The existence of localizations is guaranteed by the following result.
\begin{lem}\label{discretelocalization}
The following hold true.
\begin{enumerate}[$(1)$]
    \item The restriction map $\cinfty(\R)\rightarrow \cinfty(\R\setminus\{0\})$ exhibits $\cinfty(\R\setminus\{0\})$ as a localization of $\cinfty(\R)$ with respect to $\mathrm{id}_{\R}$.
    \item Let $f:A\rightarrow B$ be a morphism of $\cinfty$-rings and let $a\in A$ an element. Let $h:B\rightarrow B'$ be a map for which the composition $A\rightarrow B\rightarrow B'$ carries $a$ to an invertible element. Then $h$ exhibits $B'$ as a localization with respect to $f(a)$ if and only if the diagram 
   \[\begin{tikzcd}
   A\ar[d]\ar[r,"f"] & B\ar[d] \\
    A[a^{-1}] \ar[r] & B'
    \end{tikzcd}\]
    induced by the universal property of $A[a^{-1}]$ is a pushout.
    \item Let $A$ be a $\cinfty$-ring and $a\in A$ an element, and let $\tilde{a}:\cinfty(\R)\rightarrow A$ be the map classifying $a$, that is, $\tilde{a}(\mathrm{id}_{\R})=a$. Let $h:A\rightarrow A'$ be a map for which the composition $\cinfty(\R)\rightarrow A\rightarrow A'$ carries $\mathrm{id}_{\R}$ to an invertible element. Then $h$ exhibits $A'$ as a localization of $A$ with respect to $A$ if and only if the diagram 
    \[    \begin{tikzcd}
    \cinfty(\R)\ar[d]\ar[r] & A\ar[d] \\
    \cinfty(\R\setminus\{0\}) \ar[r] & A'
    \end{tikzcd}\]
    induced by the universal property of $\cinfty(\R\setminus\{0\})$ from $(1)$ is a pushout. 
    \item Let $M$ be a manifold and $f\in \cinfty(M)$ a function, then the map $\cinfty(M)\rightarrow\cinfty(f^{-1}(\R\setminus\{0\}))$ exhibits a localization with respect to $f$.
    \item Let $A=\cinfty(M)/I$ for $M$ a manifold and let $p:\cinfty(M)\rightarrow A$ denote the projection. Let $f\in A$ and let $\tilde{f}\in \cinfty(M)$ be a lift of $f$ so that $p(\tilde{f})=p$, then the map $A\rightarrow \cinfty(\tilde{f}^{-1}(\setminus\{0\}))/I|_U$, where $I|_U$ is the ideal generated by the image of $I$ under $p$, exhibits a localization with respect to $f$.
\end{enumerate} 
\end{lem}
\begin{proof}
Notice that $(2)$ is obvious from the definitions, $(3)$ follows from $(1)$ and $(2)$, $(4)$ follows at once from $(3)$ and Lemma \ref{lem:localizeopens} and $5$ is an immediate consequence of $(2)$, $(4)$ and the fact that taking quotients commutes with pushouts. It suffices to argue $(1)$. The $\cinfty$ function $\R\setminus\{0\}\rightarrow \R^2$ carrying $x$ to $(x,1/x)$ is a closed embedding onto the zero set of the submersion $xy-1:\R^2\rightarrow \R$, so the map $\cinfty(\R)\rightarrow \cinfty(\R\setminus\{0\})$ may be identified with the map $\cinfty(\R)\rightarrow \cinfty(\R^2)\rightarrow\cinfty(\R^2)/(xy-1)$. This map is obtained by applying the free $\cinfty$-ring functor $F^{\cinfty}$ to the map $\R[x]\rightarrow \R[x,y]\rightarrow \R[x,y]/(xy-1)$, which exhibits the $\R$-algebraic localization $\R[x]\rightarrow \R[x,x^{-1}]$. Now we conclude by observing that $F^{\cinfty}$ carries $\R$-algebraic localizations to localization of $\cinfty$-rings.
\end{proof} 

\begin{rmk}
The analysis of the previous lemma shows that in many cases, the $C^{\infty}$-ring localization is very different from the $\R$-algebraic localization. Indeed, inverting the identity in $C^{\infty}(\R)$ yields only those smooth functions $f(x)$ on $\R\setminus\{0\}$ that approach infinity at most polynomially fast as $x\rightarrow 0$.
\end{rmk}
\begin{cor}
Let $A$ be a finitely generated (presented) $\cinfty$-ring, then for any $a\in A$, the localization $A\rightarrow A[a^{-1}]$ is again finitely generated (presented).    
\end{cor}
If $A$ is finitely presented, then $A$ is germ determined, so all localizations of a finitely presented $\cinfty$-ring are germ determined. If $A$ is finitely generated and germ determined, then a localization of $A$ need not be germ determined. Consider the $\cinfty$-ring $\cinfty(\R)_0$ of germs of functions at the origin; this $\cinfty$-ring is obviously germ determined, but the localization $\cinfty(\R)_0[x^{-1}]\cong\cinfty(\R\setminus\{0\})/\mathfrak{m}_g^0|_{\R\setminus\{0\}}$ at the germ of the identity function is not: since $Z(\mathfrak{m}_g^0|_{\R\setminus\{0\}})=\emptyset$, the ideal $\mathfrak{m}_g^0|_{\R\setminus\{0\}}$ being germ determined would imply that $1\in \mathfrak{m}_g^0|_{\R\setminus\{0\}}$ which is not the case. The following result asserts that a localization of a germ determined $\cinfty$-ring is well behaved if it respects the zero locus.  
\begin{prop}\label{prop:localizegermdet}
Let $M$ be a manifold and let $I\subset\cinfty(M)$ be a germ determined ideal. Let $Z(I)\subset U\subset M$ be an open set, then the map $\cinfty(M)/I\rightarrow \cinfty(U)/I|_U$ restricting to $U$ is an isomorphism. 
\end{prop}
\begin{proof}
Write $p:\cinfty(M)\rightarrow \cinfty(M)/I$ for the projection. Invoking Lemma \ref{discretelocalization}, the map $\cinfty(M)/I\rightarrow \cinfty(U)/I|_U$ is the localization with respect to the image under $p$ of a characteristic function $\chi_U$ of $U$. Thus, it suffices to show that $p(\chi_U)$ is invertible. Since $M\setminus U$ and $Z(I)$ are disjoint closed sets, we may invoke normality of $M$ (twice) to find open sets $V,W$ such that 
\[ Z(I)\subset V\subset \overline{V}\subset W\subset \overline{W}\subset U.\]
Choose a partition of unity $\{\phi_W,\phi_{M\setminus\overline{V}}\}$ subordinate to the cover $\{W,M\setminus\overline{V}\}$, then $\phi_W$ equals $1$ on $\overline{V}$ and vanishes on $M\setminus W$, so $\chi_U^{-1}\phi_W$ is well defined and $\chi_U\chi_U^{-1}\phi_W=\phi_W$. Since $1-\phi_W$ vanishes in a neighbourhood of $Z(I)$ and $I$ is germ determined, we conclude. 
\end{proof}
\begin{nota} We will denote $\geodiff$ for the opposite category of the category of finitely presented $C^{\infty}$-rings. To notationally distinguish a finitely presented $C^{\infty}$-ring $A$ from $A$ as an object of $\geodiff$, we use the notation $\mathrm{Spec}\,A$ in the latter case. We also say that an ideal $J$ of a finitely presented object $C^{\infty}(\R^n)/I$ is \emph{germ determined} if the pullback of $J$ along the quotient map $C^{\infty}(\R^n)\rightarrow C^{\infty}(\R^n)/I$ is germ determined (equivalently, the $C^{\infty}$-ring $(C^{\infty}(\R^n)/I)/J$ is germ determined). 
\end{nota}

We endow $\geodiff$ with the structure of a geometry according to the following prescription: 
\begin{enumerate}[(1)]
    \item A map $f:\mathrm{Spec}\, A\rightarrow \mathrm{Spec}\, B$ in $\geodiff$ is admissible if and only if there exists an element $b\in B$ such that the image of $b$ under $f$ is invertible in $A$ and the induced map $B[1/b]\rightarrow A$ is an equivalence. 
    \item A collection $\{\mathrm{Spec}\,B[1/b_{\alpha}]\rightarrow \mathrm{Spec}\,B\}_{\alpha \in J}$ generates a covering sieve if and only if the germ determined ideal generated by the elements $b_{\alpha}$ in $B$ contains the unit.
\end{enumerate}
The definition of the topology is motivated by the following observation.
\begin{lem}\label{lem:germdetidealcovering}
Let $M$ be a manifold and let $\{f_i\}$ be a collection of functions in $\cinfty(M)$. Then the germ determined ideal of $\cinfty(M)$ generated by the collection $\{f_i\}_i$ contains $1$ if and only if the collection of open subsets $\{f_i^{-1}(\R\setminus\{0\})\}_i$ covers $M$.  
\end{lem}
\begin{proof}
Since the smallest germ determined ideal of $\cinfty(M)$ generated by an ideal $I$ is the ideal $J$ of those functions $f$ such that for each $x\in Z(I)$, there is a neighbourhood $x\in U$ such that $f|_U\subset I|_U$, we immediately see that $Z(I)=Z(J)$. Now the collection of open subsets $\{f_i^{-1}(\R\setminus\{0\})\}_i$ covers $M$ if and only if $Z((f_i)_i)=Z(J)=\emptyset$ for $J$ the germ determined ideal generated by $(f_i)_i$, but $J$ germ determined implies that $1\in Z(J)$ if and only if $Z(J)=\emptyset$.
\end{proof}
\begin{lem}\label{lem:cartspgeneratescovers}
Let $A$ be a finitely presented $\cinfty$-ring of the form $\cinfty(\R^n)/(f_1,\ldots,f_n)$ and let $p:\cinfty(\R^n)\rightarrow A$ denote the projection map. Let $\{a_i\}_{i\in I}$ be a collection of elements of $A$. The following are equivalent.
\begin{enumerate}[$(1)$]
    \item The collection $\{a_i\}_{i\in I}$ determines an admissible covering $\{A\rightarrow A[a^{-1}_i]\}$.
    \item For any collection $\{b_i\}_{i\in I}$ of elements in $\cinfty(\R^n)$ such that $p(b_i)=a_i$, the open sets $\{b_i^{-1}(\R\setminus\{0\})\}_{i\in I}$ cover $Z(f_1,\ldots,f_n)$.
    \item There exists a collection $\{b_i\}_{i\in I}$ of elements in $\cinfty(\R^n)$ such that $p(b_i)=a_i$ and the open sets $\{b_i^{-1}(\R\setminus\{0\})\}_{i\in I}$ cover $Z(f_1,\ldots,f_n)$.
    \item There exists an injection $I\subset J$ and a collection $\{b_j\}_{j\in J}$ of elements in $\cinfty(\R^n)$ that determines an admissible covering such that $p(b_j)=a_j$ for $j\in I\subset J$ and $p(b_j)=0$ for $j\in J\setminus I$. 
    \item For any collection $\{b_i\}_{i\in I}$ of elements in $\cinfty(\R^n)$ such that $p(b_i)=a_i$, let $U_i$ denote the open set $b_i^{-1}(\R\setminus\{0\})$. Then the map 
    \[ \coprod_i Z(f_1|_{U_i},\ldots,f_n|_{U_i})\longrightarrow Z(f_1,\ldots,f_n) \]
    is a surjection.
\end{enumerate}
\end{lem}
\begin{proof}
$(1)$ amounts to the following assertion: let $\{b_i\}_{i\in I}$ be any collection in $\cinfty(\R^n)$ lifting $\{a_i\}_{i\in I}$, then the germ determined ideal of $\cinfty(\R^n)$ generated by $\{b_i\}_{i\in I}\cup \{f_1,\ldots,f_n\}$ contains the unit. This is in turn equivalent to the statement that $Z((b_i)_{i},f_1,\ldots,f_n)=\emptyset$, that is, the the open sets $\{b_i^{-1}(\R\setminus\{0\})\}_{i\in I}$ cover $Z(f_1,\ldots,f_n)$. This shows that $(1)$ is equivalent to $(2)$ and $(3)$. To show $(3)\Rightarrow (4)$, we choose, using that disjoint closed sets can be separated in $\R^n$ by opens, an open set $V$ satisfying 
\[\R^n\setminus \overline{\cup_i b_i^{-1}(\R\setminus\{0\})}\subset \R^n\setminus \cup_i b_i^{-1}(\R\setminus\{0\})\subset V\subset \overline{V}\subset \R^n\setminus Z(f_1,\ldots,f_n),\]
then for a characteristic function $\chi_V$ of $V$, the collection $\{\chi_V^{-1}(\R\setminus\{0\})\}\bigcup \{b_i^{-1}(\R\setminus\{0\})\}_{i\in I}$ covers $\R^n$ and $p(\chi_V)=0$. $(4)\Rightarrow (3)$ is obvious. Note that $(5)$ is a reformulation of $(2)$ using that $Z(f_1|_{U_i},\ldots,f_n|_{U_i})=Z(f_1,\ldots,f_n)\cap U_i$. 
\end{proof}

As a consequence of Lemmas \ref{lem:germdetidealcovering} and \ref{lem:cartspgeneratescovers}, the Grothendieck topology on $\geodiff$ is generated by the open cover topology on $\mathsf{CartSp}$, in the sense that every covering family in $\geodiff$ is refined by a covering pulled back from a covering family in $\mathsf{CartSp}$. 
\begin{prop}\label{proofgeometrydisc}
 The collection of admissible morphisms and admissible coverings described above endow $\geodiff$ with an admissibility structure and a compatible topology.   
\end{prop}
\begin{proof}
We check that the admissible maps are stable under pullbacks, retracts and that, if $g$ is admissible and $h$ another map with codomain being the domain of $g$, then $h$ is admissible if and only if $g\circ h$ is admissible. The stability under pullbacks follows at once from Lemma \ref{discretelocalization}. For stability under retracts, consider a localization $f:A\rightarrow A[1/a]$ and a retraction diagram
\[
\begin{tikzcd}
A'\ar[r]\ar[d]& A \ar[r,"h"] \ar[d,"f"] &A'\ar[d]\\
B \ar[r] & A[1/a]\ar[r] & B
\end{tikzcd}
\]
Now $B$ is the localization $A'[1/h(a)]$. To see this, note that for a map $A'\rightarrow C$ that inverts $h(a)$, we get a unique map $q:A[1/a]\rightarrow C$ as in the commuting diagram
\[
\begin{tikzcd}
A'\ar[r]\ar[d]& A \ar[r,"h"] \ar[d,"f"] &A'\ar[d] \ar[dr,bend left=30]\\
B \ar[r,"g"] & A[1/a]\ar[r]\ar[rr,bend right=30,"q"] & B & C
\end{tikzcd}
\]
so we have map $B\rightarrow C$ as $q \circ g$. Note that this map is unique: suppose we have $p$ and $p'$ as in the commuting diagram
 \[
\begin{tikzcd}
A'\ar[r]\ar[d]& A \ar[r,"h"] \ar[d,"f"] &A'\ar[d] \ar[dr,bend left=30]\\
B \ar[r,"g"] & A[1/a]\ar[r,"k"] & B \ar[r,"p",shift left=1]\ar[r,"p'"',shift right=1] & C
\end{tikzcd}
\]
then by uniqueness $k$ equalizes $p$ and $p'$, and because the diagram is a retraction we have $p=p\circ k\circ g=p'\circ k \circ g=p'$. We have the claims about compositions of admissibles left to check. Choose a regular epimorphism $\cinfty(\R^n)\rightarrow A$ and a localization $A\rightarrow [a^{-1}]$, then Lemma \ref{discretelocalization} provides a pushout diagram  
\[
\begin{tikzcd}
   \cinfty(\R^n)\ar[d]\ar[r] & \cinfty(U)\ar[d] \\
   A\ar[r] & A[a^{-1}] 
\end{tikzcd}
\]
with $U=\tilde{a}^{-1}(\R\setminus\{0\})$ where $\tilde{a}$ is a lift of $a$ to $\cinfty(\R^n)$. It follows from Corollary \ref{cor:algeffepi} that the right vertical map is also a regular epimorphism, so we can lift any $b\in A[a^{-1}]$ to some $\tilde{b}$ in $\cinfty(U)$ and form the right pushout diagram
\[
\begin{tikzcd}
   \cinfty(\R^n)\ar[d]\ar[r] & \cinfty(U)\ar[d] \ar[r] & \cinfty(V)\ar[d]\\
   A\ar[r] & A[a^{-1}] \ar[r] & A[a^{-1}][b^{-1}]
\end{tikzcd}
\]
with $V=\tilde{b}^{-1}(\R\setminus\{0\})$. Then the outer square is a pushout diagram, which implies that $A\rightarrow  A[a^{-1}][b^{-1}]$ is a localization by Lemma \ref{discretelocalization} again. Conversely, if we are given a map $ A[a^{-1}]\rightarrow A'$ such that the composition $A\rightarrow A[a^{-1}]\rightarrow A'$ is a localization, we consider the diagram
\[
\begin{tikzcd}
   A\ar[d]\ar[r] & A[a^{-1}] \ar[d,equal] \ar[r] &A'\ar[d,equal]\\
   A[a^{-1}] \ar[r,equal] & A[a^{-1}] \ar[r] & A'.
\end{tikzcd}
\]
In view of Lemma \ref{discretelocalization}, it suffices to show that both squares are pushouts. The right square is obviously a pushout. For the left square, Lemma \ref{discretelocalization} guarantees that it suffices to show that the right vertical map exhibits a localization with respect to the image of $a$ in $A[a^{-1}]$, which is clearly the case as this image is already invertible. It remains to be shown that the collection of admissible coverings specify a Grothendieck pretopology. This is clear from the characterization provided by Lemma \ref{lem:cartspgeneratescovers}.
\end{proof}
\begin{rmk}
The relation between the pregeometry $\diff$ and the geometry $\geodiff$ is as follows: by Corollary \ref{smoothringtruncated}, $\geodiff$ is a $0$-truncated geometric envelope of $\diff$. 
\end{rmk}

\begin{vari}\label{var:finitarygeodiff}
Let $\geodiff^{\mathrm{fin}}$ denote the geometry that has the same underlying \infcat and admissibility structure as $\geodiff$, but is endowed with the following coarser topology compatible with the admissibility structure: a sieve on $\speco\,A\in \geodiff^{\mathrm{fin}}$ is covering if and only it contains a finite collection $\{\speco\,A[a^{-1}_i]\rightarrow\speco\,A\}_{i\in I}$ such that the elements $\{a_i\}$ generate the unit ideal in $A$. The identity functor $\geodiff^{\mathrm{fin}}\rightarrow \geodiff$ is a transformation of geometries which fits into a commuting diagram 
\[
\begin{tikzcd}
 \diff^{\mathrm{fin}}\ar[d,"\cinfty(\_)"]\ar[r] & \diff\ar[d,"\cinfty(\_)"]\\
 \geodiff^{\mathrm{fin}} \ar[r] & \geodiff
\end{tikzcd}
\]
where the vertical functors, which carry a manifold to its $\cinfty$-ring of $\cinfty$-functions, exhibit $0$-truncated geometric envelopes.
\end{vari}

Let $\Of:\geodiff\rightarrow \mathcal{S}$ be a discrete $\geodiff$-structure in spaces, which can be identified with a $C^{\infty}$-ring by the equivalences $\mathrm{Str}_{\geodiff}(\mathcal{S})\simeq \mathrm{Ind}(\geodiff^{op})\simeq C^{\infty}\mathsf{ring}$; the corresponding $C^{\infty}$-ring $A_{\mathcal{O}}$ is up to unique isomorphism determined by $\Hom_{C^{\infty}\mathsf{ring}}(B,A_{\Of})=\Of(B)$ for $B$ a finitely presented $C^{\infty}$-ring. We'd like to give a characterization of what it means to be a (non-discrete) $\geodiff$-structure on $\mathcal{S}$ in terms of the corresponding $C^{\infty}$-ring. We need the following lemma, due to Bunge, Dubuc and Joyal \cite{BD1}.
\begin{lem}\label{opencoveringlemma}
Any open covering on $\R^n$ is generated under pullbacks, composition and refinement by coverings on $\R$.
\end{lem}
\begin{proof}
Fix an open covering $\{U_{\alpha}\rightarrow \R^n\}_{a\in J}$. This covering is a composition of the coverings $W_{I}=\coprod_{i\in I}U_{\alpha_i}$ for the finite subsets $I\subset J$. To see that such finite coverings are pulled back from coverings on $\R$, we first reduce any finite covering to a covering $U_1\coprod U_2=W$ by induction. Choosing characteristic functions $\chi_{U_1}$ and $\chi_{U_2}$, we may replace them by $\frac{\chi^2_{U_1}}{\chi_{U_1}^2+\chi_{U_1}^2}$ and  $\frac{\chi^2_{U_2}}{\chi_{U_1}^2+\chi_{U_1}^2}$, so that $\chi_{U_1}+\chi_{U_2}=1$. Now $U_1=\chi_1(\R\setminus\{0\})$ and $U_2=\chi_1(\R\setminus\{1\})$. \\
Now we show that the covering $\{W_{I}\rightarrow \R^n\}_{I\subset J,|I|\leq \infty}$ is refined by a covering pulled back from $\R$. Choose some covering $\{Y_k\rightarrow \R\}$ by bounded open sets, and a proper smooth function $\varphi:\R^n\rightarrow \R$ (for instance, the square length function $(x_1,\ldots,x_n)\rightarrow x_1^2+\ldots+ x_n^2$). The opens $\varphi^{-1}(Y_k)$ cover $\R^n$ and are bounded, so each such open is covered by a finite collection $U_{\alpha_{j_1}},\ldots,U_{\alpha_{j_n}}$ and thus $\varphi^{-1}(Y_k)\subset W_S$ for some finite index set $S$. Consequently, there is a refinement map $\coprod_k \varphi^{-1}(Y_k)\rightarrow \coprod_{I\subset J,|I|<\infty} W_I$. \end{proof}
\begin{prop}\label{prop:archimedeanlocal}
$\Of$ is a $\geodiff$-structure on $\mathcal{S}$ if and only if $A_{\Of}$ is local as a commutative ring and $A_{\Of}$ is \emph{pointed} in the sense that there is an $\R$-algebra map $A_{\Of}\rightarrow \R$ (which is necessarily surjective so that this map is the projection onto the residue field).
\end{prop}
\begin{proof}
We should check that $A_{\Of}$ is a local ring with residue field $\R$ if and only if for each finitely presented $C^{\infty}$-ring $B$ and each admissible covering $\{B\rightarrow B[1/b_{\alpha}]\}_{\alpha \in J}$, the map 
\[\coprod_{\alpha}\Hom_{C^{\infty}\mathsf{ring}}(B[1/b_{\alpha}],A_{\Of})\longrightarrow \Hom_{C^{\infty}\mathsf{ring}}(B,A_{\Of})\]
is an epimorphism. By Lemma \ref{lem:cartspgeneratescovers}, any admissible covering on $B$ is pulled back from a covering on a free $C^{\infty}$-ring, so, because epimorphisms are stable under pullbacks in $\mathsf{Set}$, we note that it is enough to prove the claim for the collection of free $C^{\infty}$-rings. By Lemma \ref{lem:germdetidealcovering}, an admissible covering of $\mathrm{Spec}\,C^{\infty}(\R^n)$ is the same thing as an open covering of $\R^n$. Thus, we should check that $A_{\Of}$ is local with residue field $\R$ if only if for each open cover $\{U_{\alpha}\rightarrow \R^n\}_{\alpha \in J}$, the map $\coprod_{\alpha}\Of(U_{\alpha})\rightarrow \Of(\R^n)$ is an epimorphism. Because epimorphisms are stable under pullback, composition and refinement in $\mathsf{Set}$, we reduce further to having to check the statement only for covering families of $\R$. In one direction, consider the open covering $\R\setminus \{0\}\cap \R\setminus \{1\}\rightarrow \R$, and note that we have transverse pullback diagrams
\[
\begin{tikzcd}
\R\setminus \{0\} \ar[d]\ar[r] & \R^2\ar[d,"{(a,b)\mapsto ab}"]\\
* \ar[r,"1"] & \R
\end{tikzcd}
 \qquad
\begin{tikzcd}
\R\setminus \{1\} \ar[d]\ar[r] & \R^2\ar[d,"{(a,b)\mapsto (a-1)b}"]\\
* \ar[r,"1"] & \R
\end{tikzcd}\]
showing that $\Of (\R\setminus \{0\})=A_{\Of}^2\times_{A_{\Of}} \{1\}$ is the set $A_{\Of}^{\times}$ of invertible elements of $A_{\Of}$, and similarly $\Of(\R\setminus \{1\})$ is the set $1-A^{\times}_{\Of}$ of elements $a\in A_{\Of}$ such that $1-a$ is invertible. Clearly, $A_{\Of}^{\times} \coprod (1-A^{\times}_{\Of})\rightarrow A$ is an epimorphism if and only if $A_{\Of}$ is local as a commutative ring. Now we show that there is a map $p:A_{\Of}\rightarrow \R$ of $C^{\infty}$-rings which is nonzero if $A_{\Of}$ is nonzero: $\Of$ gives a functor 
\[\mathrm{Open}(\R)\longrightarrow \mathrm{Sub}(A_{\Of}(\R)),\quad U\longmapsto \Of(U), \]
(note that since $\Of$ preserves all pullbacks in $\mathrm{Open}(\R)$, $\Of(U)$ is a subobject of $A_{\Of}(\R)$), which is a map of locales as it is left exact and sends coverings to epimorphism. Since the underlying topological spaces of these locales are sober, we get a map $p:A_{\Of}(\R)\rightarrow \R$ of sets. This map is a morphism of $\cinfty$-rings and thus a morphism of $\R$-algebras: it suffices to show that for a smooth map $f:\R^n\rightarrow \R^m$, the diagram of locales 
\[
\begin{tikzcd}
\mathrm{Open}(\R^m)\ar[d,"f^{-1}"] \ar[r,"\Of"] & \mathrm{Sub}(A_{\Of}(\R)^m) \ar[d,"f^*"]\\
\mathrm{Open}(\R^n)\ar[r,"\Of"] & \mathrm{Sub}(A_{\Of}(\R)^n) 
\end{tikzcd}
\]
commutes, where the right vertical map $f^*$ sends a subobject $X$ of $A_{\Of}(\R)^m$ to the pullback $X\times_{A_{\Of}(\R)^m}A_{\Of}(\R)^n$ along the map $f_*:A_{\Of}(\R)^n\rightarrow A_{\Of}(\R)^m$. Concretely, for any smooth map $f:\R^n\rightarrow \R^m$ and any open in $U$, we ask that $\Of(f^{-1}(U))\simeq \Of(U)\times_{\Of(\R^m)}\Of(\R^n)$. This is clearly the case since $\Of$ preserves pullbacks along open inclusions. The kernel of the map $p:A_{\Of}\rightarrow \R$ just constructed is a maximal ideal, so by locality of $A_{\Of}$, $p$ must be the projection onto the residue field. \\
For the converse, suppose that $A_{\Of}$ is local and comes equipped with a map $A_{\Of}\rightarrow\R$ of commutative $\R$-algebras. We want to show that for any open covering $\coprod_{\alpha} U_{\alpha}\rightarrow \R$, the induced map $\coprod_{\alpha}\Of(U_{\alpha})\rightarrow A_{\Of}(\R)$ is an epimorphism. Points in $A_{\Of}$ corresponds by Yoneda to maps $q:C^{\infty}(\R)\rightarrow A_{\Of}$, so it suffices to show that each such map factors as $C^{\infty}(\R)\rightarrow C^{\infty}(U_{\alpha})\rightarrow A_{\Of}$ for some index $\alpha$. Consider the composition 
\[\varphi:C^{\infty}(\R)\overset{q}{\longrightarrow} A_{\Of}\longrightarrow \R,\]
which is a morphism of commutative $\R$-algebras. Denote $p=\varphi(x)\in\R$, where $x$ denotes the identity on $\R$.  Suppose that $f(p)$ vanishes, then Hadamard's lemma implies that $f(x)=(x-p)h(x)$ for some $h\in\cinfty(\R)$, so $\varphi(f)=(\varphi(x)-\varphi(p))\varphi(h)=(\varphi(x)-p)\varphi(h)=0$ as $\varphi$ is an $\R$-algebra morphism. It follows that the maximal ideal $\ker(\ev_p)$ is contained in $\ker(\varphi)$ which establishes that $\ker(\ev_p)=\ker(\varphi)$ as both are maximal ideals. We conclude that $\varphi$ is given by evaluation $\mathrm{ev}_p$ at some $p\in  U_{\alpha} \subset \R$ for some index $\alpha$; let $\chi_{U_{\alpha}}$ be a characteristic function for $U_{\alpha}$, then $\mathrm{ev}_p(\chi_{U_{\alpha}})\neq 0$, implying that the image of $\chi_{U_{\alpha}}$ under $q$ is not in the maximal ideal $\ker(q)$ of $A_{\Of}$, so $q(\chi_{U_{\alpha}})$ is invertible in $A_{\Of}$ by locality. Now $q$ factors through the localization $C^{\infty}(\R)\rightarrow C^{\infty}(U_{\alpha})$ of $\chi_{U_{\alpha}}$ so we are done.
\end{proof}
\begin{rmk}\label{rmk:evaluationmorphism}
The proof above shows that morphisms $\cinfty(\R)\rightarrow \R$ of $\R$-algebras are of the form $\ev_p$ for some $p\in \R$ and are therefore $\cinfty$-ring morphisms. Suppose we have an $\R$-algebra map $\varphi:\cinfty(\R^n)\rightarrow \R$, then we have real numbers $p_i=\varphi(x_i)$ where $x_i$ is the $i'th$ coordinate function $\R^n\rightarrow \R$. Suppose $f:\R^n\rightarrow\R$ vanishes at $p=(p_i)_i$, then Hadamard's lemma implies that $f$ may be written as $f=\sum_i g_i(x_i-p_i)$ so that $\varphi(f)=0$. It follows that $\varphi$ coincides with the map $\ev_p$ and is therefore a $\cinfty$-ring morphism. Now consider a finitely generated $\cinfty$-ring $\cinfty(\R^n)/I$, then the regular epimorphism $\cinfty(\R^n)\rightarrow\cinfty(\R^n)/I$ induces an injection
\[ \Hom_{\cinfty\mathsf{ring}}(A,\R)\longrightarrow \Hom_{\cinfty\mathsf{ring}}(\cinfty(\R^n),\R).\]
It follows right away that any $\R$-algebra morphism $\cinfty(\R^n)/I\rightarrow \R$ is a $\cinfty$-ring morphism of the form $\ev_p$ for some $p\in Z(I)$ and that this correspondence determines a bijection $\Hom_{\cinfty\mathsf{ring}}(A,\R)\cong Z(I)$. Since any $\cinfty$-ring is a filtered colimit of its finitely generated sub $\cinfty$-rings, we deduce a bijection $\Hom_{\cinfty\mathsf{ring}}(A,\R)\cong \Hom_{\mathsf{CAlg}^0_{\R}}(A^{\rmalg},\R)$.
\end{rmk}
The preceding remark has the following consequence.
\begin{prop}
Let $A$ be a finitely presented $\cinfty$-ring, then a collection $\{A\rightarrow A[a_i^{-1}]\}_i$ of localization morphisms determines an admissible covering if and only if the map 
\[ \coprod_i\Hom_{\cinfty\mathsf{ring}}(A[a^{-1}],\R)\longrightarrow  \Hom_{\cinfty\mathsf{ring}}(A,\R) \]
is a surjection.
\end{prop}
\begin{proof}
This is an immediate consequence of point $(5)$ of Lemma \ref{lem:cartspgeneratescovers} and the bijection $\Hom_{\cinfty\mathsf{ring}}(A,\R)\cong Z(I)$ of Remark \ref{rmk:evaluationmorphism}.
\end{proof}
\begin{prop}\label{localmorphisms}
Let $\alpha:\Of\rightarrow \Of'$ be a morphism of left exact functors $\geodiff\rightarrow\mathcal{S}$. Then $\alpha$ is a local morphism if and only if the corresponding morphism $f_{\alpha}: A_{\Of}\rightarrow A_{\Of'}$ is local as a map of commutative rings.
\end{prop}
\begin{proof}
The map $f_{\alpha}$ is local as a map of commutative rings if and only if $f_{\alpha}$ reflects invertibility, which is true if and only if $A^{\times}_{\Of}\simeq A_{\Of}\times_{A_{\Of'}}A^{\times}_{\Of'}$. If $\alpha$ is local, this obviously holds. In the other direction, we want to show that for each localization $B\rightarrow B[1/b]$ of finitely presented $C^{\infty}$-rings, the naturality square induced by $\alpha$ gives an equivalence $\Of(B[1/b])\simeq \Of(B)\times_{\Of'(B)}\Of'(B[1/b])$. Because $B$ is finitely presented, we have a pushout $B[1/b]\simeq C^{\infty}(\R\setminus\{0\})\otimes^{\infty}_{C^{\infty}(\R)}B$ so we get a commuting cube 
\[
\begin{tikzcd}
& A^{\times}_{\Of} \arrow[rr] \arrow[dd] & & A^{\times}_{\Of' } \arrow[dd] \\
\Of(B[1/b])  \arrow[ur ]\arrow[rr, crossing over] \arrow[dd] & & \Of'(B[1/b])\arrow[ur]  \\
& A_{\Of}  \arrow[rr] & & A_{\Of'}\\
\Of(B)\arrow[rr] \arrow[ur]& & \Of'(B)\arrow[ur] \arrow[from=uu, crossing over]\\
\end{tikzcd}
\]
Because the side faces are pullbacks and the back face is a pullback by assumption, the front face is a pullback as well. 
\end{proof}
\begin{cor}\label{cor:evermorphismlocal}
Let $\alpha:\Of\rightarrow \Of'$ be a morphism of $\geodiff$-structures in spaces. Then $\alpha$ is a local morphism.
\end{cor}
\begin{proof}
As the residue field of both $A_{\Of}$ and $A_{\Of'}$ is $\R$, any morphism of commutative $\R$-algebras between them is local. 
\end{proof}
In summary, the 1-category $\strloc_{\spa}(\geodiff)\subset\cartsp$ is the \emph{full} subcategory spanned by \emph{Archimedean local} $C^{\infty}$-rings, that is, $C^{\infty}$-rings $A$ for which the underlying commutative $\R$-algebra of $A$ is a local ring equipped with a $\cinfty$-ring morphism (equivalently, an $\R$-algebra morphism) $A\rightarrow\R$. Whenever we talk about local $C^{\infty}$-rings in the sequel, we mean this notion. Note that in particular, local $\cinfty$-rings are pointed: since there is at most one $\cinfty$-ring morphism $A\rightarrow\R$ for $A$ local, $\R$ is a final object of $\strloc_{\spa}(\geodiff)$ and we have a factorization 
\[ \strloc_{\spa}(\geodiff)\hooklongrightarrow \cinfty\mathsf{ring}_{/\R}\longrightarrow \cartsp\]
where the first functor is also fully faithful.
\begin{rmk}\label{localenoughpoints}
Let $(\xtop,\Of_{\xtop})$ be a $(\geodiff)_{\mathrm{disc}}$-structured $\infty$-topos for which the underling 1-topos $\tau_{\leq 0}\xtop$ has enough points. Then $\Of_{\xtop}$ is a $\geodiff$-structure if and only if for each point $p_*:\spa\rightarrow\xtop$, the induced $C^{\infty}$-ring is local. Moreover, every morphism $\alpha:\Of_{\xtop}\rightarrow \Of'_{\xtop}$ of $\geodiff$-structures is local. 
\end{rmk}
\begin{rmk}\label{rmk:classtoposlocal}
It follows from the preceding results and discussion that the classifying \inftop of $\geodiff$-structures $\shv(\geodiff)$ is a classifying \inftop for (Archimedean) local $\cinfty$-rings. Let $\geodiff^{\mathrm{fin}}$ be the geometry from Variant \ref{var:finitarygeodiff}, then one deduces from the arguments of Proposition \ref{prop:archimedeanlocal} and \ref{localmorphisms} that a $\geodiff^{\mathrm{fin}}$-structure in spaces is simply a $\cinfty$-ring that is local as a commutative ring, and that 
a local morphism of $\geodiff^{\mathrm{fin}}$-structures in spaces corresponds to a local morphism of local rings. It follows that $\shv(\geodiff^{\mathrm{fin}})$ is a classifying \inftop for $\cinfty$-rings that are local as commutative algebras.  
\end{rmk}
It follows from Proposition \ref{prop:localfactorization} that for each local $\cinfty$-ring $B$, the full subcategory inclusion $\strloc_{\geodiff}(\spa)_{/B}\subset \cinfty\mathsf{ring}_{/B}$ admits an accessible left adjoint $L_B$, which carries a map $A\rightarrow B$ to the second map in the factorization 
\[A\longrightarrow \colim_{a\in A,f(a)\in B^{\times}}A[a^{-1}]\longrightarrow B,  \]
in particular, $\strloc_{\geodiff}(\spa)$ is a presentable 1-category.  
\begin{lem}\label{lem:localizationrn}
Let $\cinfty(\R^n)/I\rightarrow \R$ be an $\R$-algebra morphism, corresponding to an evaluation map $\ev_p$ for some $p\in Z(I)$ by Remark \ref{rmk:evaluationmorphism}, then the quotient $q_p:\cinfty(\R^n)/I\rightarrow \cinfty(\R^n)/(I,\mathfrak{m}^g_p)$ taking germs at $p$ exhibits a unit transformation for the left adjoint $L_{\R}$ at $\cinfty(\R^n)/I$.
\end{lem}
\begin{proof}
We are required to show that for any local $\cinfty$-ring $A$, composition with $q_p$ induces a bijection 
\[\Hom_{\cinfty\mathsf{ring}_{/\R}}(\cinfty(\R^n)/(I,\mathfrak{m}^g_p),A)\overset{\cong}{\longrightarrow} \Hom_{\cinfty\mathsf{ring}_{/\R}}(\cinfty(\R^n)/I,A). \]
Since $q_p$ is a regular epimorphism, this map is an injection, so it suffices to show it is surjective. For any map $f:\cinfty(\R^n)/I\rightarrow A$ the composition $\cinfty(\R^n)/I\rightarrow \R$ coincides with $\ev_p$. Suppose that $g\in \cinfty(\R^n)/I$ comes from a function $\tilde{g}$ that vanishes in a neighbourhood $U$ of $p$, then there exists a function $\tilde{h}:\R^n\rightarrow\R$ that equals $1$ on a smaller neighbourhood $p\in V\subset U$ and vanishes outside $p$, which induces an element $h\in\cinfty(\R^n)/I$. We have $\ev_p(h)=1$, so $f(h)$ is invertible in $A$, but $gh=0$ so $f(g)f(h)=0$ which then implies $f(g)=0$. It follows that $f$ factors through $\cinfty(\R^n)/I\rightarrow \cinfty(\R^n)/(I,\mathfrak{m}^g_p)$. 
\end{proof}
The functor $L_B$ is an accessible localization, but we can say a bit more.
\begin{prop}\label{prop:localizationsiftedcolim}
Let $B$ be a local $\cinfty$-ring, then $\strloc_{\geodiff}(\spa)_{/B}\rightarrow \cinfty\mathsf{ring}_{/B}$ preserves sifted colimits.  \end{prop}
\begin{proof}
Since the right fibration $\icat_{/C}\rightarrow \icat$ detects colimits for all $C\in \icat$, it suffices to consider the case of $B=\R$. We first show the following.
\begin{enumerate}
    \item[$(*)$] Let $A$ be a $\cinfty$-ring and write $\Of_A$ for the associated $(\geodiff)_{\mathrm{disc}}$-structure, then $A$ is local, that is, $\Of_A$ is a $\geodiff$-structure, if and only if for every \emph{good} open cover $\{U_i\subset\R\}$ (a cover on $\R$ for which any nonempty finite intersection is diffeomorphic to $\R$), the induced map $\coprod_i\Of_A(U_i)\cong \coprod_iA\rightarrow A$ is an effective epimorphism.
\end{enumerate}
To see this, we recall that Proposition \ref{prop:archimedeanlocal} guarantees that it suffices to consider open coverings of $\R$. Now the assertion $(*)$ follows easily from the fact that every open cover on $\R$ may be refined by a good open cover. Let $\mathcal{J}:K\rightarrow \strloc_{\geodiff}(\spa)$ be a colimit diagram with $K$ a small sifted simplicial set, then using $(*)$, we see that it suffices to show that
\begin{enumerate}[$(a)$]
    \item The colimit of the composition $K\rightarrow \strloc_{\geodiff}(\spa) \subset \cinfty\mathsf{ring}\rightarrow \fun(\mathsf{CartSp},\set)$ is a $\cinfty$-ring.
    \item The colimit $\cinfty$-ring carries good open covers on $\R$ to effective epimorphisms in $\set$.
\end{enumerate}
The assertion $(a)$ follows immediately from the fact that the formation of sifted colimits preserves products in $\set$. Let $A=\colim_{\mathcal{J}}K$, then we wish to show that for any good open cover $\{U_i\subset\R\}_i$, the map $\coprod_iA\rightarrow A$ is an effective epimorphism, but this map is a colimit of the diagram 
\[\mathcal{J}_{\Delta^1}: K\longrightarrow \fun(\Delta^1\times\cartsp,\set) \]
given by $\{\coprod_i\mathcal{J}(k)\rightarrow \mathcal{J}(k)\}_{k\in K}$ and each of these maps is an effective epimorphism as all the objects $\mathcal{J}(k)$ are assumed to be local.
\end{proof}
\begin{cor}\label{cor:locepi}
For any $\cinfty$-ring $A$ and any $\R$-algebra morphism $p:A\rightarrow\R$, the map $A\rightarrow A_p$ is a regular epimorphism.    
\end{cor}
\begin{proof}
Since the class of regular epimorphisms is stable under colimits and every $\cinfty$-ring is a filtered colimit of its finitely generated subrings, it suffices to argue that the construction 
\[ (p:A\rightarrow \R) \longmapsto (A\rightarrow A_p) \]
preserves filtered colimits, in view of Lemma \ref{lem:localizationrn}, but this construction preserves sifted colimits (in fact, it preserves all colimits as we show in the proof of Proposition \ref{prop:lawveretheoryofgerms} below).
\end{proof}
More generally, one can show that for any regular epimorphism $A\rightarrow B$ of $\cinfty$-rings (that is, a surjective map on the underlying sets) with $B$ local, the map $A\rightarrow L_BA$ is again a regular epimorphism. \\
Lemma \ref{lem:localizationrn} and Proposition \ref{prop:localizationsiftedcolim} imply that the 1-category of local $\cinfty$-rings is in fact the category of algebras for a Lawvere theory of germs of smooth functions.
\begin{prop}\label{prop:lawveretheoryofgerms}
The following hold true.
\begin{enumerate}[$(1)$]
    \item The full subcategory inclusion $\strloc_{\spa}(\geodiff)\subset\strloc_{\spa}((\geodiff)_{\mathrm{disc}})=\cinfty\mathsf{ring}$ preserves colimits.
    \item Let $\mathsf{CartSp}_{g}$ be the category defined as follows.
    \begin{enumerate}
        \item[$(O)$] Objects are pairs $(\R^n,p)$ of a Cartesian spaces with a point $p\in\R^n$.
        \item[$(M)$] A morphism $(\R^n,p)\rightarrow (\R^m,q)$ is a germ at $p$ of a smooth function that takes $p$ to $q$.
    \end{enumerate}
    Then there is a canonical equivalence between the 1-category of algebras for this Lawvere theory and the 1-category $\strloc_{\spa}(\geodiff)$ of local $\cinfty$-rings.
\end{enumerate}
 \end{prop}
 In the proof, we use some language from Subsection \ref{sec:lawvere}.
\begin{proof}
 We first argue that $\strloc_{\spa}(\geodiff)$ is the 1-category of algebras for \emph{some} Lawvere theory. It follows from Proposition \ref{prop:lawverestability} $(2)$ that the 1-category $\cinfty\mathsf{ring}_{/\R}$ is the category of algebras for the Lawvere theory whose objects are pairs $(\R^n,p)$ with $p\in\R^n$ a point and whose morphisms $(\R^n,p)\rightarrow (\R^m,q)$ are smooth maps $f:\R^n\rightarrow\R^m$ such that $f(p)=q$. Proposition \ref{prop:localizationsiftedcolim} tells us that the inclusion $\strloc_{\spa}(\geodiff)\subset \cinfty\mathsf{ring}_{/\R}$ preserves sifted colimits; it is a consequence of \cite{HA}, Proposition 7.1.4.12 that $\strloc_{\spa}(\geodiff)$ is the 1-category of algebras for the \emph{opposite} of the essential image of the composition
\[ \cartsp^{op}_{/\R}\subset \cinfty\mathsf{ring}_{/\R} \overset{L_{\R}}{\longrightarrow}  \strloc_{\spa}(\geodiff), \]
where the second functor is the left adjoint taking $g:A\rightarrow \R$ to the filtered colimit $\colim_{a\in A,\,g(a)\neq 0}A[a^{-1}]$. Lemma \ref{lem:localizationrn} shows that $L_{\R}$ carries $\ev_p:\cinfty(\R^n)\rightarrow \R$ to the local $\cinfty$-ring $\cinfty(\R^n)/{\mathfrak{m}^g_p}$. The opposite of the category of these $\cinfty$-rings is precisely the category we described in the statement of the proposition. We are left to show $(1)$. Since $\strloc_{\spa}(\geodiff)$ is projectively generated by the objects $\cinfty(\R^n)/\mathfrak{m}^g_p$ for $p\in \R^n$ a point, it suffices to show that the inclusion $\strloc_{\spa}(\geodiff)\subset\cinfty\mathsf{ring}$ preserves binary coproducts of this class of objects in view of \cite{HTT}, Proposition 5.5.8.15, so we should show that the canonical map 
\[ \cinfty(\R^n)/\mathfrak{m}^g_p \oinfty \cinfty(\R^m)/\mathfrak{m}^g_q \longrightarrow \cinfty(\R^{n+m})/\mathfrak{m}^g_{(p,q)} \]
is an isomorphism for any $p:*\subset\R^n$ and $q:*\subset\R^m$. Unwinding the definitions, this amounts to the assertion that the ideal $\mathfrak{m}^g_{(p,q)}$ is generated by the images of the ideals $\mathfrak{m}^g_p$ and $\mathfrak{m}^g_q$. To see this, we let $\pi_1:\R^{n+m}\rightarrow\R^n$ and $\pi_2:\R^{n+m}\rightarrow\R^m$ denote the projections and consider a function $f$ vanishing on some open $(p,q)\in U$. Then $U$ contains a product open neighbourhood $(p-\epsilon,p+\epsilon)\times (q-\delta,q+\delta)$, so choosing a function $h_1$ on $\R^n$ that vanishes on $(p-\epsilon,p+\epsilon)$ and equals $1$ on $\R^n\setminus \pi_1(U)$ and a function $h_2$ on $\R^m$ that vanishes on $(q-\delta,q+\delta)$ and equals $1$ on $\R^m\setminus \pi_2(U)$, we conclude that $f$ is divisible by $h_1+h_2$. 
\end{proof}
\begin{rmk}
$\geodiff^{\mathrm{fin}}$-structures in spaces, that is, $\cinfty$-rings that are merely local as a commutative ring are much more poorly behaved: the full subcategory of such in $\cinfty\mathsf{ring}$ is not stable under binary coproducts.
\end{rmk}
For a general geometry $\geo$, the category $\sch_{\fp}(\geo)$ of $\geo$-schemes of finite presentation admits limits and the full subcategory inclusion $\sch_{\fp}(\geo)\subset \rtop(\geo)$ preserves these limits, but $\rtop(\geo)$ need not admit many limits in general. For $\geo=\geodiff$ however, the previous results imply that there is in fact a rather large supply of limits of structured spaces. 
\begin{cor}\label{cor:geodiffstableunderlimits}
Let $\rtop^{\mathrm{sp}}_{{0-\mathrm{loc}}}\subset \rtop$ be the full subcategory spanned by spatial $0$-localic \inftopoit, and consider the pullbacks
\[\rtop^{\mathrm{sp}}_{{0-\mathrm{loc}}}(\geodiff)=\rtop(\geodiff)\times_{\rtop}\rtop^{\mathrm{sp}}_{{0-\mathrm{loc}}},\quad\quad \rtop^{\mathrm{sp}}_{{0-\mathrm{loc}}}((\geodiff)_{\mathrm{disc}})=\rtop((\geodiff)_{\mathrm{disc}})\times_{\rtop}\rtop^{\mathrm{sp}}_{{0-\mathrm{loc}}}.\]
Then $\rtop^{\mathrm{sp}}_{{0-\mathrm{loc}}}(\geodiff)$ admits all small limits and the subcategory inclusion $\rtop^{\mathrm{sp}}_{{0-\mathrm{loc}}}(\geodiff)\subset \rtop^{\mathrm{sp}}_{{0-\mathrm{loc}}}((\geodiff)_{\mathrm{disc}})$ is stable under small limits.
\end{cor}
\begin{proof}
Since the \infcat of $\rtop^{\mathrm{sp}}_{0-\mathrm{loc}}$ of spatial $0$-localic \inftopoi is equivalent to the nerve of the category of sober topological spaces, it has small limits.  It follows from a routine application of the theory of relative colimits (\cite{HTT}, Lemma 4.3.1.5 and Propositions 4.3.1.9 and 4.3.1.10) that it suffices to argue that for any $\xtop\in \rtop^{\mathrm{sp}}_{0-\mathrm{loc}}$, the full subcategory $\strloc_{\spa}(\geodiff)\subset\strloc_{\spa}((\geodiff)_{\mathrm{disc}})=\fun^{\mathrm{lex}}(\geodiff,\xtop)$ is stable under small colimits. Since every $(\geodiff)_{\mathrm{disc}}$-structure takes $0$-truncated values, we may test this on the underlying $1$-topoi. Now the condition that a $(\geodiff)_{\mathrm{disc}}$-structure on $\tau_{\leq 0}\geodiff$, (that is, a sheaf of $\cinfty$-rings on $\tau_{\leq 0}\xtop$) is a $\geodiff$-structure may be tested on points, since $\tau_{\leq 0}\xtop$ has enough points, so we reduce to the case $\xtop=\spa$, which follows from Proposition \ref{prop:lawveretheoryofgerms}.
\end{proof}
\begin{rmk}
Many of the results above hold in much greater generality. For instance, if $\pregeo$ is any pregeometry, and $\Of$ a $\pregeo$-structure, then one can show that the \infcat $\strloc_{\pregeo}(\spa)_{/\Of}$ is the \infcat of algebras for a Lawvere theory. In particular, the \infcat $\strloc_{\pregeo}(\spa)$ of $\pregeo$-structures and local morphisms in spaces is an \infcat of algebras for a Lawvere theory if and only if it admits a final object, in which case we have an equality $\strloc_{\pregeo}(\spa)=\str_{\pregeo}(\spa)$.
\end{rmk}
\subsubsection{Locality in real algebraic geometry}
We now turn to some (classical) results for $\cinfty$-rings that have local underlying commutative algebras, that is, $\geodiff^{\mathrm{fin}}$-structures, giving alternative proofs of some of the results of \cite{JOYAL1986271}, \cite{kockrealclosed} and \cite{MR} (this material will not be used in the rest of this paper). We have seen in Remark \ref{rmk:classtoposlocal} that $\geodiff^{\mathrm{fin}}$-structures and local morphism between them are simply $\cinfty$-rings that are local as commutative $\R$-algebras and local morphism between them. With the language of geometries at our disposal, we can put this observation in context. Recall that the underlying algebra functor $(\_)^{\rmalg}$ admits a left adjoint $F^{\cinfty}$, the free $\cinfty$-ring functor, which carries finitely presented objects to finitely presented objects (since $(\_)^{\rmalg}$ preserves sifted colimits), that is, $F^{\cinfty}$ restricts to a right exact functor \[F^{\cinfty}:(\mathsf{CAlg}^0_{\R})_{\fp}\longrightarrow\cinfty\mathsf{ring}_{\fp}.\] Now the observation of Remark \ref{rmk:classtoposlocal} can be reformulated as follows.
\begin{prop}
 The functor $F^{\cinfty}$ determines a transformation of geometries $F^{\cinfty}:\geo_{\mathrm{Zar}}(\R)\rightarrow\geodiff^{\mathrm{fin}}$. Moreover, the functor $F^{\cinfty}$ generates the geometry structure on $\cinfty\mathsf{ring}_{\fp}^{op}$ in the sense of Remark \ref{rmk:geometrypushforward}.  
\end{prop}
To prove this, we note it suffices to show that
\begin{enumerate}[$(1)$]
    \item $F^{\cinfty}$ carries $\R$-algebraic localizations to localization of $\cinfty$-rings.
    \item any localization of $\cinfty$-rings is pulled back from the map $F^{\cinfty}(\R[x]\rightarrow \R[x,x^{-1}])$.
    \item a (finite) admissible covering $\{A\rightarrow A[a_i^{-1}]\}_{i\in I}$ induces a finite admissible covering $\{F^{\cinfty}(A)\rightarrow F^{\cinfty}(A[a_i^{-1}])\}_{i\in I}$ 
    \item any admissible covering for the finitary \'{e}tale topology is obtained by pullback, composition and refinement from the covering $\{\R\setminus\{0\}\hookrightarrow \R,\R\setminus\{1\}\hookrightarrow \R\}$
\end{enumerate}
Over the course of this subsection, only $(3)$ was not explicitly verified, but note that it follows from the fact that if $\R[x_1,\ldots,x_n]\rightarrow A$ is a regular epimorphism, then $\cinfty(\R^n)\rightarrow F^{\cinfty}(A)$ is a regular epimorphism. Naturally, we may ask what changes if we replace $\geo_{\mathrm{Zar}}(\R)$ with the geometry $\geo_{\et}(\R)$ controlling algebraic geometry over $\R$ for the \'{e}tale topology. To treat this question, we first need to introduce the appropriate notion of an \'{e}tale morphism of finitely presented $\cinfty$-rings.
\begin{defn}\label{defn:etalemorphismgeodiff}
 Let $f:A\rightarrow B$ be a morphism of finitely presented $\cinfty$-rings, then $f$ is \emph{\'{e}tale} if there exists an admissible covering $\{B\rightarrow B[b_i^{-1}]\}_{i\in I}$ such that each composition $A\rightarrow B[b_i^{-1}]$ exhibits a localization.    
\end{defn}
\begin{rmk}
As in algebraic geometry, one can show that to be an \'{e}tale map of finitely presented $\cinfty$-rings is equivalent to being a \emph{formally \'{e}tale map}; that is, having the left lifting property with respect to maps of $\cinfty$-rings with nilpotent kernel.    
\end{rmk}
One readily verifies that the subcategory on the \'{e}tale morphisms determines an admissibility structure on $\cinfty\mathsf{ring}_{\fp}$.  
\begin{vari}
 The 1-category $\cinfty\mathsf{ring}_{\fp}$ admits the following admissibility structure with compatible topology.
 \begin{enumerate}[$(1)$]
     \item A morphism $A\rightarrow B$ is admissible if and only if it is \'{e}tale.
     \item A collection $\{A\rightarrow B_i\}_i$ of \'{e}tale morphism generates a covering sieve if and only if it generates a covering sieve for the \'{e}tale topology on $\geodiff$ and $I$ is finite. 
 \end{enumerate}
 We let $\geo_{\mathrm{Diff},\et}^{\mathrm{fin}}$ denote this geometry.
\end{vari}

The identity functor on $\cinfty\mathsf{ring}_{\fp}$ induces a transformation of geometries $\geo_{\mathrm{Diff}}^{\mathrm{fin}}\rightarrow \geo_{\mathrm{Diff},\et}^{\mathrm{fin}}$ since localization morphisms are obviously \'{e}tale. The following proposition makes precise the sense in which the Zariski and \'{e}tale topology coincide in differential geometry.
\begin{prop}
The transformation of geometries $\geo_{\mathrm{Diff}}^{\mathrm{fin}}\rightarrow \geo_{\mathrm{Diff},\et}^{\mathrm{fin}}$ induces for any \inftop $\xtop$ equivalences of \infcats $\str_{\geo_{\mathrm{Diff}}^{\mathrm{fin}}}(\xtop)\simeq \str_{\geo_{\mathrm{Diff}}^{\mathrm{fin}}}(\xtop)$ and $\strloc_{\geo_{\mathrm{Diff}}^{\mathrm{fin}}}(\xtop)\simeq \strloc_{\geo_{\mathrm{Diff}}^{\mathrm{fin}}}(\xtop)$.   
\end{prop}
\begin{rmk}
Removing the condition that a covering sieve contains a finite covering family of \'{e}tale morphisms, we obtain a geometry $\geo_{\mathrm{Diff},\et}$ which the functor $\pregeo_{\mathrm{Diff},\et}\rightarrow \geo_{\mathrm{Diff},\et}$ exhibits as a geometric envelope. It follows that the identity functor $\geodiff\rightarrow\geo_{\mathrm{Diff},\et}$ determines a transformation of geometries which induces equivalences $\str_{\geo_{\mathrm{Diff},\et}}(\xtop)\simeq \str_{\geodiff}(\xtop)$ and $\strloc_{\geo_{\mathrm{Diff},\et}}(\xtop)\simeq \strloc_{\geodiff}(\xtop)$
\end{rmk}
It may seem reasonable to suggest that the functor $F^{\cinfty}$ induces a transformation of geometries $\geo_{\et}(\R)\rightarrow\geo_{\mathrm{Diff},\et}^{\mathrm{fin}}$, but this is not the case: the separable field extension $\R\rightarrow \C$ is obviously an \'{e}tale covering of commutative $\R$-algebras, but $F^{\cinfty}$ carries carries this map to the map $\R\rightarrow 0$. Since $\shv(\geo_{\et}(\R))$ is a classifying \inftop for strictly Henselian (local) commutative $\R$-algebras we deduce the following.
\begin{cor}
 No $\cinfty$-ring with local underlying commutative algebra is strictly Henselian.   
\end{cor}
As it turns out, the lack of algebraic closure is the only obstruction to the functor $F^{\cinfty}$ being a transformation of geometries. To make this precise, we need to systematically remove the \'{e}tale coverings $A\rightarrow \prod_iA_i$ such that for any $\R$-point $A\rightarrow \R$, the codomain of the \'{e}tale map $\R\rightarrow \R\otimes_{A}\prod_iA_i$, which can be identified with a finite product $\prod_j R$ of finite separable extensions of $\R$ (that is, $\prod_j R$ is a finite product of $\C$'s and $\R$'s), contains no factor of $\R$. We proceed by introducing a version of \emph{real algebraic geometry}. Recall the following notion.
\begin{defn}
 A field $F$ is \emph{formally real} if $-1\in F$ is not a finite sum of squares. Let $A$ be a commutative ring, then we say that a closed point $\speco\,F\rightarrow \speco\,A$ is \emph{formally real} if $F$ is a formally real field.  
\end{defn}
We define another geometry structure on $(\mathsf{CAlg}^0_{\R,\fp})^{op}$ replacing faithful flatness by the condition that a covering family covers all formally real closed points.
\begin{vari}\label{vari:formrealgeom}
The 1-category $(\mathsf{CAlg}^0_{\R,\fp})^{op}$ admits the following admissibility structure with compatible topology.
 \begin{enumerate}[$(1)$]
     \item A morphism $A\rightarrow B$ is admissible if and only if it is an \'{e}tale morphism of commutative $\R$-algebras.
     \item A collection $\{A\rightarrow B_i\}_i$ of \'{e}tale morphism generates a covering sieve if and only if the induced map 
     \[ f:\coprod_i\speco\,B_i \longrightarrow\speco\,A \]
     on prime spectra has the property that for every formally real closed point $\mathfrak{p}\subset \speco\,A$, there is an index $i$ and a formally real closed point $\mathfrak{q}$ of $B_i$ such that $f(\mathfrak{q})=\mathfrak{p}$.  
 \end{enumerate}
 We call the coverings of $(2)$ \emph{formally real \'{e}tale coverings} and the Grothendieck topology they generate on $(\mathsf{CAlg}^0_{\R,\fp})^{op}$ the \emph{formally real \'{e}tale topology}.
 We let $\geo_{\et,\mathrm{f}\R}(\R)$ denote the resulting geometry.
\end{vari}
\begin{rmk}
Given a family $\{A\rightarrow B_i\}_{i\in I}$ of finitely presented commutative $\R$-algebras, the condition in point $(2)$ of Variant \ref{vari:formrealgeom} is equivalent to the following:
\begin{enumerate}
    \item[$(2')$] For every maximal ideal $\mathfrak{p}$ of $A$ determining a formally real field $A/\mathfrak{p}$, there exists an index $i$ such that the tensor product $A/\mathfrak{p}\otimes_{A}B_i$ is nonzero and admits a maximal ideal $\mathfrak{q}$ such that $(A/\mathfrak{p}\otimes_{A}B_i)/\mathfrak{q}$ is formally real. 
\end{enumerate}
From this condition, it is obvious that a pullback of a formally real \'{e}tale covering is again a formally real \'{e}tale covering. It follows easily that the formally real \'{e}tale coverings indeed define a Grothendieck topology compatible with the admissibility structure. 
\end{rmk}
\begin{rmk}
The \'{e}tale topology and the formally real \'{e}tale topology on $(\mathsf{CAlg}^0_{\R,\fp})^{op}$ are incomparable: neither refines the other.   
\end{rmk}
\begin{prop}\label{prop:finettrans}
The functor $F^{\cinfty}$ determines a transformation of geometries $\geo_{\et,\mathrm{f}\R}(\R)\rightarrow \geo_{\mathrm{Diff},\et}^{\mathrm{fin}}$.
\end{prop}
\begin{proof}
Clearly $F^{\cinfty}$ is left exact. Let $f:A\rightarrow B$ be an \'{e}tale map of finitely presented commutative $\R$-algebras, then we wish to show that $F^{\cinfty}$ carries this map to an \'{e}tale map of finitely presented $\cinfty$-rings. The structure theory of \'{e}tale maps guarantees the existence of a pushout diagram
\[
\begin{tikzcd}
 \R[x_1,\ldots,x_n]\ar[d,"\phi"] \ar[r]& A\ar[d,"f"] \\
 \R[y_1,\ldots,y_n][\mathrm{Det}^{-1}] \ar[r] & B
\end{tikzcd}
\]
where $\phi(x_i)=g_i(y_1,\ldots,y_n)$ are polynomials and $\mathrm{Det}$ is the determinant of the matrix $\{\frac{\del g_i}{\del y_j}\}_{1\leq i,j\leq n}$, so it suffices to consider \'{e}tale maps that are of the form as the left vertical one in the diagram above. Using that $F^{\cinfty}$ carries localizations of commutative rings to localizations of $\cinfty$-rings, we observe that we have to argue that composite $\varphi:\cinfty(\R^n)\overset{g^*}{\rightarrow} \cinfty(\R^n)\rightarrow  \cinfty(U)$ is \'{e}tale, where $g:\R^n\rightarrow\R^n$ is some polynomial map and $U\subset\R^n$ is the open set where the determinant of $\{\frac{\del g_i}{\del y_j}\}_{1\leq i,j\leq n}$ does not vanish. It follows immediately from the standard implicit function theorem for $\cinfty$ functions and Lemma \ref{discretelocalization} that $\varphi$ is an \'{e}tale map of $\cinfty$-rings of finite presentation. It remains to be shown that $F^{\cinfty}$ carries admissible covers to admissible covers. Suppose $\{A\rightarrow B_i\}_i$ is a formally real \'{e}tale covering, then we wish to show that the map \[\coprod_i\Hom_{\cinfty\mathsf{ring}}(F^{\cinfty}(B_i),\R)\longrightarrow \Hom_{\cinfty\mathsf{ring}}(F^{\cinfty}(A),\R)\] is a surjection. By adjointness, it suffices to show that for each $\R$-point $p:\speco\,\R\rightarrow \speco\,A$, there is an index $i$ and an $\R$-point $q:\speco\,\R\rightarrow \speco\,B_i$ such that the composition $\speco\,\R\overset{q}{\rightarrow} \speco\,B_i\rightarrow\speco\,A$ equals $p$. By assumption, there exists an index $i$ such that the tensor product $\R\otimes_{A}B_i$ is not empty and contains a maximal ideal corresponding to a formally real field. Since $A\rightarrow B_i$ is \'{e}tale, the map $\R\rightarrow \R\otimes_{A}B_i$ is a finite product $\prod_j R_j$ of finite separable extensions of $\R$. Since the only separable extensions of $\R$ are $\R$ and $\C$ and $\C$ is not formally real, we conclude that at least one of the factors $R_j$ equals $\R$, which gives the desired $\R$-point $q:B_i\rightarrow \R\otimes_{A}B_i\rightarrow\R$.
\end{proof}
It follows that for $A$ a $\geo^{\mathrm{fin}}_{\mathrm{Diff},\et}$-structure, that is, a $\cinfty$-ring whose underlying commutative ring is local, the commutative ring $A^{\rmalg}$ is a $\geo_{\et,\mathrm{f}\R}(\R)$-structure. To understand what it means to be a $\geo_{\et,\mathrm{f}\R}(\R)$-structure in terms of classical commutative algebra, we need a formally real analogue of the correspondence between strictly Henselian rings and $\geo_{\et}(k)$-structures.
\begin{defn}
Let $(A,\mathfrak{m})$ be a local ring, then we say that $A$ is \emph{separably real closed} if 
\begin{enumerate}[$(1)$]
    \item $A$ is a Henselian local ring.
    \item the residue field $A/\mathfrak{m}$ is a \emph{real closed} field; that is, $A/\mathfrak{m}$ is formally real and does not admit nontrivial formally real algebraic extensions.
\end{enumerate}
\end{defn}
The following is an analogue of Corollary 7.19 of \cite{dagvii}. We leave its proof, which one adapts readily from the strictly Henselian case, to the motivated reader.
\begin{prop}
A commutative $\R$-algebra $A$ is a $\geo_{\et,\mathrm{f}\R}(\R)$-structure in spaces if and only if $A$ is a separably real closed local ring. Moreover, a morphism of such local rings is local as a morphism of $\geo_{\et,\mathrm{f}\R}(\R)$-structures if and only if it is a local $\R$-algebra homomorphism.
\end{prop}
In other words, $\shv(\geo_{\et,\mathrm{f}\R}(\R))$ is a classifying \inftop for separably real closed local commutative $\R$-algebras.
\begin{cor}[Kock,Joyal-Reyes]
 Let $A$ be a $\cinfty$-ring with local underlying commutative algebra. Then $A^{\rmalg}$ is separably real closed.   
\end{cor}
\begin{rmk}
While the \'{e}tale topology and the formally real \'{e}tale topology on $(\calg^0_{\R,\fp})^{op}$ are incomparable, they both refine the Nisnevich topology. Using Proposition \ref{prop:finettrans} we conclude that composition with the functor $F^{\cinfty}$ carries sheaves for the finitary \'{e}tale topology on the category of finitely presented $\cinfty$-rings to Nisnevich sheaves. Since the points for the Nisnevich topology are Henselian local rings, this corresponds to the fact that $\cinfty$-rings that are local as commutative rings are already Henselian.
\end{rmk}
\begin{rmk}
The formally real spectrum $\spec^{\geo_{\et,\mathrm{f}\R}(\R)}\,A$ of a commutative ring $A$ for the geometry $\geo_{\et,\mathrm{f}\R}(\R)$ is quite different from the \'{e}tale spectrum $\spec^{\geo_{\et}(\R)}\,A$: while the latter is an affine Deligne-Mumford stack and has in particular a 1-localic (but not 0-localic) underlying \inftopt, the underlying \inftop of $\spec^{\geo_{\et,\mathrm{f}\R}(\R)}\,A$ is 0-localic. Since it is coherent, Hochster duality identifies the underlying \inftop of $\spec^{\geo_{\mathrm{f}\R}(\R)}\,A$ with the sheaf topos on a coherent topological space. The points of this space correspond to maps $f:A\rightarrow A'$ where $A'$ is separably real closed local $\R$-algebra and $f$ is a filtered colimit of \'{e}tale morphisms. The structure sheaf of $\spec^{\geo_{\et,\mathrm{f}\R}(\R)}\,A$ is a familiar object in real algebraic geometry: it is the sheaf of \emph{Nash functions}. For a textbook account of real algebraic geometry and a detailed discussion of the formally real spectrum, we refer to \cite{coste2018real}. 
\end{rmk}
\subsubsection{Archimedean spectra and locally $\cinfty$-ringed spaces}
The general theory of geometries provides us with an \infcat of locally $\cinfty$-ringed \inftopoi and a spectrum functor $\geodiff$ left adjoint to the global sections functor. The $\geodiff$-spectrum is well behaved on the class of finitely presented $\cinfty$-rings, but in general it has the undesirable feature that the underlying \inftop of $\spec^{\geodiff}\,A$ need not be spatial. We now show how to recover the \emph{Archimedean spectrum} of a $\cinfty$-ring, first defined by Dubuc and investigated by Moerdijk-Reyes-van Qu\^{e} and (extensively) by Joyce, from the $\geodiff$-spectrum. Recall that a $0$-localic $\infty$-topos $\xtop$ that is spatial arises as the $\infty$-category of sheaves on a sober topological space, which determines an equivalence of $\infty$-categories between $\mathsf{Top}$, the $\infty$-category of sober topological spaces and continuous maps, and $\rtop^{\mathrm{sp}}_{0-\mathrm{loc}}$, the full subcategory of $\rtop$ spanned by spatial $0$-localic $\infty$-topoi. We can lift this equivalence to one on the level of structured spaces for the geometry $\geodiff$.
\begin{defn}
Let $\mathsf{Top}(\cinfty\mathsf{ring})$ be the category of \emph{$\cinfty$-ringed spaces} defined as follows.
\begin{enumerate}
    \item[$(O)$] Objects are pairs $(X,\Of_X)$ of a sober topological space equipped with  sheaves of $C^{\infty}$-rings.
    \item[$(M)$] Morphisms are pairs $(f,\alpha):(X,\Of_X)\rightarrow (Y,\Of_Y)$ where $f:X\rightarrow Y$ is a continuous map and $\alpha:f^*\Of_Y\rightarrow \Of_X$ a morphisms of sheaves of $\cinfty$-rings.
\end{enumerate}
Let $\mathsf{Top}^{\loc}(\cinfty\mathsf{ring})\subset \mathsf{Top}(\cinfty\mathsf{ring})$ be the subcategory of \emph{locally $\cinfty$-ringed spaces} defined as follows.
\begin{enumerate}
    \item[$(O)$] Objects are pairs $(X,\Of_X)$ such that for each point $p:*\rightarrow X$, the $\cinfty$-ring $p^*\Of_X$ is (Archimedean) local.  
    \item[$(M)$] Morphisms are pairs $(f,\alpha):(X,\Of_X)\rightarrow (Y,\Of_Y)$ where $f:X\rightarrow Y$ is a continuous map and $\alpha:f^*\Of_Y\rightarrow \Of_X$ induces for each $p:*\rightarrow X$ a local morphism $(p\circ f)^*\Of_Y\rightarrow p^*\Of_X$ of local $\cinfty$-rings.
\end{enumerate}
\end{defn}

\begin{rmk}
 It follows from the observation in the proof of Corollary \ref{cor:evermorphismlocal} that the condition of being a local morphism is superfluous, so that the subcategory of locally $\cinfty$-ringed spaces is a full subcategory of the category of $\cinfty$-ringed spaces. The argument of Corollary \ref{cor:geodiffstableunderlimits} shows that this full subcategory is stable under arbitrary small limits.
\end{rmk}
Combining Propositions \ref{prop:archimedeanlocal}, \ref{localmorphisms} and Remark \ref{localenoughpoints}, we deduce the following.
\begin{prop}\label{prop:comparisonringedspaces}
Let $\rtop^{\mathrm{sp}}_{0-\mathrm{loc}}((\geodiff)_{\mathrm{disc}})$ be the $\infty$-category of pairs $(\xtop,\Of_{\xtop})$ where $\xtop$ is a spatial $0$-localic $\infty$-topos. There is a canonical equivalence of $\infty$-categories
\[\zeta: \rtop^{\mathrm{sp}}_{0-\mathrm{loc}}((\geodiff)_{\mathrm{disc}})\overset{\simeq}{\longrightarrow} \mathsf{Top}(\cinfty\mathsf{ring})  \]
that associates to a pair $(\shv(X),\Of_{X})$ the pair $(X,\Of_X)$. Moreover, $\zeta$ restricts to an equivalence 
\[\zeta:\rtop^{\mathrm{sp}}_{0-\mathrm{loc}}(\geodiff)\overset{\simeq}{\longrightarrow}  \mathsf{Top}^{\loc}(\cinfty\mathsf{ring}),  \]
where $\rtop^{\mathrm{sp}}_{0-\mathrm{loc}}(\geodiff)$ is the $\infty$-category of spatial 0-localic $\geodiff$-structured $\infty$-topoi.
\end{prop}
We now relate the $\geodiff$-spectrum to the \emph{Archimedean spectrum} of $\cinfty$-rings, first introduced by Dubuc in \cite{dubucschemes}. Let $A$ be a $\cinfty$-ring. Note that a map $A\rightarrow B$ of $\cinfty$-rings is admissible as a morphism of $\mathrm{Pro}(\geodiff)$ if and only if it $A\rightarrow B$ exhibits a localization at some element of $A$; this follows at once from Lemma \ref{discretelocalization}. Let $\cinfty\mathsf{ring}_{A}^{\mathrm{ad}}\subset \cinfty\mathsf{ring}_{A}$ be the full subcategory spanned by maps exhibiting a localization with respect to some element of $A$ so that we have an equivalence of $1$-categories $\cinfty\mathsf{ring}_{/A}^{op,\mathrm{ad}}\simeq \mathrm{Pro}(\geo)^{\mathrm{ad}}_{/A}$. It is an immediate consequence of the definition of a localization that the 1-category $\cinfty\mathsf{ring}_{/A}^{op,\mathrm{ad}}$ is in fact a 0-category, that is, the nerve of a poset. Since we have an equivalence $\fun^{\mathrm{lex}}(\icat,\xtop)\simeq \shv_{\pro(\icat)^{op}}(\xtop)$ for any \inftop $\xtop$ and any small \infcat admitting finite limits $\icat$, we may think of the presheaf
\[\widetilde{\Of_{\speco\,A}}:\geodiff \longrightarrow \pshv(\cinfty\mathsf{ring}_{/A}^{op,\mathrm{ad}})  \]
as the tautological $\cinfty\mathsf{ring}$-valued presheaf
\[ \widetilde{\Of_{\speco\,A}} :\cinfty\mathsf{ring}_{A/}^{\mathrm{ad}} \longrightarrow \cinfty\mathsf{ring},\quad \quad (A\rightarrow A[a^{-1}])\longmapsto A[a^{-1}]. \]
Unwinding the definitions, the $\geodiff$-spectrum $\spec^{\geodiff}\,A$ is the pair $(\shv(\cinfty\mathsf{ring}_{/A}^{op,\mathrm{ad}}),\Of_{\speco\,A})$ where $\Of_{\speco\,A}$ is the sheafification of the presheaf of $\cinfty$-rings $\widetilde{\Of_{\speco\,A}}$. 
\begin{cons}[Archimedean spectrum]\label{cons:realspectrum}
Let $A$ be a $C^{\infty}$-ring. Consider the functor $\Hom_{C^{\infty}\mathsf{ring}}(\_,\R):\cinfty\mathsf{ring}^{op}\rightarrow \set$ represented by $\R$, which induces a functor
\[ \Hom_{C^{\infty}\mathsf{ring}}(\_,\R):\cinfty\mathsf{ring}_{/A}^{op,\mathrm{ad}} \longrightarrow \mathrm{Sub}(\Hom_{C^{\infty}\mathsf{ring}}(A,\R)). \]
Denote by $\mathcal{B}(A)$ the essential image of this functor, which consists of the sets
\[\{U_{a}\}_{a\in A},\quad U_a:=\mathrm{ev}_a^{-1}(\R\setminus \{0\}) \]
where
\[\mathrm{ev}_a:\Hom_{C^{\infty}\mathsf{ring}}(A,\R) \longrightarrow \R,\quad \mathrm{ev}_a(f)=f(a).\]
Note that the subposet $\mathcal{B}(A)\subset\mathrm{Sub}(\Hom_{C^{\infty}\mathsf{ring}}(A,\R))$ is stable under pullbacks pullbacks: the meet of $U_a$ and $U_b$ equals $U_{ab}$. In particular, $\mathcal{B}(A)$ is a basis for a topology $\mathcal{B}(A)\subset\mathrm{Open}(\Hom_{C^{\infty}\mathsf{ring}}(A,\R))$. Let $\specr\,A$ be the topological space whose underlying set is $\Hom_{_{C^{\infty}\mathsf{ring}}}(A,\R)$, equipped with the topology generated by the basis $\mathcal{B}(A)$. Since $\mathcal{B}(A)\subset \mathrm{Open}(\Hom_{C^{\infty}\mathsf{ring}}(A,\R))$ is stable under finite limits, restricting sheaves along this inclusion induces an equivalence $\shv_{\icat}(\specr\,A)\simeq \shv_{\icat}(\mathcal{B}(A))$ for any presentable \infcat $\icat$. Since the functor $\Hom_{\cinfty\mathsf{ring}}(\_,\R)$ preserves finite limits, the composition 
\[ \cinfty\mathsf{ring}^{op,\mathrm{ad}}_{/A}\longrightarrow  \mathcal{B}(A)\overset{j}{\hooklongrightarrow}\shv( \mathcal{B}(A)) \]
induces an algebraic morphism $\varphi_A^*:\shv(\cinfty\mathsf{ring}^{op,\mathrm{ad}}_{/A})\rightarrow\shv(\mathcal{B}(A))$ by \cite{HTT}, Proposition 6.2.3.20. Let $\Of_{\specr\,A}$ be the composition 
\[\geodiff\overset{\Of_{\speco\,A}}{\longrightarrow} \shv(\cinfty\mathsf{ring}^{op,\mathrm{ad}}_{/A}) \overset{ \varphi_A^*}{\longrightarrow} \shv(\mathcal{B}(A))\simeq \shv(\specr\,A),  \]
then $\Of_{\specr\,A}$ is a $\geodiff$-structure on the spatial $0$-localic \inftop $\shv(\specr\,A)$, which we can identify with the sheafification of the diagonal left Kan extension in the (non commuting) diagram
\[
\begin{tikzcd}
\cinfty\mathsf{ring}^{\mathrm{ad}}_{/A}\ar[d] \ar[r,"\widetilde{\Of_{\speco\,A}}"] &[2em] \cinfty\mathsf{ring} \\
\mathcal{B}(A). \ar[ur,"\widetilde{\Of_{\specr\,A}}"']
\end{tikzcd}
\]
By Proposition \ref{prop:comparisonringedspaces}, the $\geodiff$-structure $\Of_{\specr\,A}$ corresponds to a locally $\cinfty$-ringed space $(\specr\,A,\Of_{\specr\,A})$, which we call the \emph{Archimedean spectrum} of $A$.
\end{cons}
\begin{defn}
 We call a locally $\cinfty$-ringed space equivalent to one of the form $(\specr\,A,\Of_{\specr\,A})$ an \emph{affine $\cinfty$-scheme}, and denote by $\caff\subset \mathsf{Top}^{\loc}(\cinfty\mathsf{ring})$ the full subcategory spanned by affine $\cinfty$-schemes. A locally $C^{\infty}$-ringed space $(X,\Of_X)$ is a \emph{$C^{\infty}$-scheme} if there is a covering $\{U_i\rightarrow X\}$ such that $(U,\Of_X|_{U_i})$ is equivalent to the Archimedean spectrum of some $C^{\infty}$-ring. We denote the full subcategory of $\mathsf{Top}^{\loc}(\cinfty\mathsf{ring})$ spanned by $C^{\infty}$-schemes by $C^{\infty}\sch$.    
\end{defn}
\begin{ex}
Let $M$ be a manifold, then the collection of $\cinfty$-ring maps $\cinfty(M)\rightarrow \R$ is in bijection with the points of $M$ by Remark \ref{rmk:evaluationmorphism}. It follows from Construction \ref{cons:realspectrum} and the fact that any open set of $M$ admits a characteristic function that the poset $\mathcal{B}(\cinfty(M))$ coincides with $\mathrm{Open}(M)$ and the presheaf $\widetilde{\Of_{\speco\,\cinfty(M)}}$ coincides with the sheaf of $\cinfty$ functions on $M$. Moreover, the functor $\cinfty\mathsf{ring}^{op,\mathrm{ad}}_{/\cinfty(M)}\rightarrow\mathcal{B}(\cinfty(M))$ is an isomorphism of posets. To see this, it suffices to show that if for $a,b\in \cinfty(M)$, we have an inclusion $U_a\subset U_b$, then $a$ is invertible in $\cinfty(M)[b^{-1}]$, but this is clear from the observation that the latter assertion holds if and only if $a:M\rightarrow\R$ does not vanish on $b^{-1}(\R\setminus\{0\})$.   
\end{ex}
More generally, we have the following result.
\begin{prop}\label{prop:germdetbasis}
Let $A$ be a finitely generated $\cinfty$-ring of the form $\cinfty(\R^n)/I$ and let $p:\cinfty(\R)\rightarrow A$ denote the projection.
\begin{enumerate}[$(1)$]
    \item The topological space $\specr\,A$ coincides with the subspace $Z(I)\subset \R^n$.
    \item Suppose that for all $a\in A$, the localization $A[a^{-1}]$ is germ determined. Then the map of posets $\Hom_{\cinfty\mathsf{ring}}(\_,\R):\cinfty\mathsf{ring}^{op,\mathrm{ad}}_{/A}\rightarrow \mathcal{B}(A)$ is an equivalence. 
\end{enumerate}
\end{prop}
\begin{proof}
For $(1)$, we note that the epimorphism $p$ induces an injection $\Hom_{\cinfty\mathsf{ring}}(A,\R)\subset \Hom_{\cinfty\mathsf{ring}}(\cinfty(\R^n),\R)$ and identifies the former set with the collection of evaluation maps at points $x\in Z(I)$, by Remark \ref{rmk:evaluationmorphism}. Let $a\in A$ and choose a lift $\tilde{a}\in \cinfty(\R^n)$ such that $p(\tilde{a})=a$, then it is obvious that the set $\ev_{a}^{-1}(\R\setminus\{0\})$ coincides with the set $\tilde{a}^{-1}(\R\setminus \{0\})\cap Z(I)$ so the basis $\mathcal{B}(A)$ is a topology and coincides with the subspace topology on $Z(I)$. For $(2)$, we first suppose that $Z(I)=\emptyset$. Since $I$ is germ determined, $I$ contains the unit and $\cinfty(\R^n)/I=0$, in which case we obviously have the desired equivalence of posets. Suppose $Z(I)$ is not empty. The map $\Hom_{\cinfty\mathsf{ring}}(\_,\R)$ is clearly essentially surjective. To see it is fully faithful, suppose that for $a,b\in A$, we have $U_a\subset U_b$, that is, we can find lifts $\tilde{a}$ and $\tilde{b}$ such that $\tilde{a}^{-1}(\R\setminus\{0\})\cap Z(I)\subset \tilde{b}^{-1}(\R\setminus\{0\})\cap Z(I)$, then we wish to show that $b$ is invertible in $A[a^{-1}]$. Write $U_{\tilde{a}}=\tilde{a}^{-1}(\R\setminus\{0\})$ and $U_{\tilde{b}}=\tilde{b}^{-1}(\R\setminus\{0\})$ so that $A[a^{-1}]\cong \cinfty(U_{\tilde{a}})/I|_{U_{\tilde{a}}}$ and $A[b^{-1}]\cong \cinfty(U_{\tilde{b}})/I|_{U_{\tilde{b}}}$ and form a pushout diagram
\[
\begin{tikzcd}
 A \ar[d] \ar[r] & A[a^{-1}] \ar[d] \\
 A[b^{-1}] \ar[r] & A',
\end{tikzcd}
\]
then we may identify the map $A[a^{-1}]\rightarrow A'$ with the localization $\cinfty(U_{\tilde{a}})/I|_{U_{\tilde{a}}}\rightarrow \cinfty(U_{\tilde{a}}\cap U_{\tilde{b}})/I|_{U_{\tilde{a}}\cap U_{\tilde{b}}}$. It suffices to show that this map is an isomorphism, but since $Z(I|_{U_{\tilde{a}}})=Z(I_{U_{\tilde{a}}\cap U_{\tilde{b}}})$ and $A[a^{-1}]$ is assumed germ determined, this follows from Proposition \ref{prop:localizegermdet}. 
\end{proof}
In general, the map of posets $\Hom_{\cinfty\mathsf{ring}}(\_,\R):\cinfty\mathsf{ring}^{op,\mathrm{ad}}_{/A}\rightarrow \mathcal{B}(A)$ is not an equivalence.
\begin{ex}
Let $A$ be the $\cinfty$-ring $\cinfty(\R)/I$ where $I$ is the ideal of compactly supported functions. Then $\cinfty\mathsf{ring}^{op,\mathrm{ad}}_{/A}$ can be identified with the set $\mathcal{B}(A)_{\infty}$ of equivalence classes of open subsets of $\R$ by the equivalence relation that identifies $U$ and $V$ if there is a closed bounded set $K\in \R$ for which $U\cap \R\setminus K=V\cap \R\setminus K$. This set has an induced poset structure that admits finite meets for which the map of sets $\cinfty\mathsf{ring}^{op,\mathrm{ad}}_{/A}\rightarrow \mathcal{B}(A)_{\infty}$ becomes a functor preserving finite meets. We therefore have a Grothendieck topology on $\mathcal{B}(A)_{\infty}$. By construction, this topology must contain the covering pulled back from some covering $\{U_i\rightarrow\R\}$ of $\R$ by compactly supported opens, but since $\cinfty(U_i)/I|_{U_i}=0$ for all $i$, we deduce that the family $\{\cinfty(\R)_I\rightarrow 0\}$ is a covering. We conclude that the topology on $\mathcal{B}(A)_{\infty}$ is the dense one, that is, every nonempty sieve is covering. In particular, the associated sheaf topos is not trivial. On the other hand, the set $\Hom_{\cinfty\mathsf{ring}}(A,\R)$ is empty so the associated sheaf topos is trivial.
\end{ex}
Even when the functor $\cinfty\mathsf{ring}^{op,\mathrm{ad}}_{/A}\rightarrow \mathcal{B}(A)$ is an equivalence, it need not identify the Grothendieck topology on both sides. 
\begin{ex}
Let $A=\cinfty(\R^{\mathbb{N}})=\colim_{n\in\N}\cinfty(\R^n)$ where the maps in the cotower are induced by the maps $\R^{n}\rightarrow\R^{n-1}$ forgetting the last coordinate, then $\specr\,A=\R^{\mathbb{N}}$, and for each $n\in \N$, we have a projection $\pi_{n}:\R^{\mathbb{N}}\rightarrow\R^n$; the basis $\mathcal{B}(A)\subset \mathrm{Open}(\specr\,A)$ is spanned by the objects of the form $\pi_n^{-1}(U)$ for $U\subset\R^n$ any open. Using that all the maps $\pi_n$ are surjective, one readily verifies that for $a,b\in A$, the inclusion $U_a\subset U_b$ implies that $b$ is invertible in $A[a^{-1}]$, which in turn implies that the functor $\cinfty\mathsf{ring}^{op,\mathrm{ad}}_{/A}\rightarrow \mathcal{B}(A)$ is an equivalence. Consider the open set $V=\pi_1^{-1}(V')$ where $V'\subset\R$ is the open set $\cup_{k\in \N} (k,k+1/2)$, then $\pi_1^{-1}(V)=\cup_{k\in \N} \pi_1^{-1}(k,k+1/2)$. For each $k\in \N$, we may choose a cover $\{W_{j_k}\rightarrow \pi_1^{-1}(k,k+1/2)\}$ where each $W_{j_k}$ is of the form $\pi_k^{-1}(C)$ for $C\subset\R^n$ a translation of an open cube $(0,1/2)^k$. Then the collection $\{W_{j_k}\rightarrow V\}_{j_k}$ in $\mathcal{B}(A)$ is manifestly a covering for the topology induced from $\mathrm{Open}(\specr\,A)$, but does not contain a covering family of the form $\{g^{-1}(U_i)\rightarrow V\}$ for $g:V\rightarrow \specr\,B$ a map induced by a morphism $\hat{g}:B\rightarrow \cinfty(V)$ with $B$ finitely presented and $\{U_i\subset\specr\,B\}$ a cover. Indeed, $B$ is finitely presented so $\hat{g}$ must factor through some $\cinfty(Q)$ for $Q\subset \R^l$ the inverse image of $\cup_{k\in \N} (k,k+1/2)$ by the projection $\R^l\rightarrow \R$. But none of the opens $W_{j_k}$ for $k>l$ can contain an open set of the form $\pi^{-1}_l(U)$, $U\subset\R^l$.     
\end{ex}
These examples notwithstanding, the algebraic morphism $\shv(\cinfty\mathsf{ring}^{op,\mathrm{ad}}_{/A})\rightarrow \shv(\mathcal{B}(A))$ is not far from an equivalence, in the sense specified below.
\begin{defn}\label{defn:spatialreflection}
 Let $\xtop,\ytop$ be 0-localic \inftopoi with $\ytop$ spatial, then an algebraic morphism $f^*:\xtop\rightarrow\ytop$ \emph{exhibits $\ytop$ as the spatial reflection of $\xtop$} if for any spatial $0$-localic \inftop $\mathcal{Z}$, composition with $f^*$ induces an equivalence
 \[ \Hom_{\ltop}(\ytop,\mathcal{Z})\overset{\simeq}{\longrightarrow} \Hom_{\ltop}(\xtop,\mathcal{Z})  \]
 of (discrete) Kan complexes.
\end{defn}
Let $\xtop\in \rtop_{0-\mathrm{loc}}$ be a $0$-localic \inftop of with associated locale $\mathcal{U}=\tau_{\leq -1}\xtop$. To this locale, we may associate in turn a sober topological space $X_{\mathcal{U}}$: its points are the morphisms of locales $p:1\rightarrow \mathcal{U}$, that is, the frame homomorphisms $p^{-1}:\mathcal{U}\rightarrow \{0,1\}$, equipped with the topology that declares a set $U\subset X_{\mathcal{U}}$ open precisely if there is some $V\in \mathcal{U}$ such that for all $p\in U$ we have $p^{-1}(V)=1$. We obtain a map of locales 
\[ \mathrm{Open}(X_{\mathcal{U}}) \longrightarrow \mathcal{U}  \]
which induces an algebraic morphism 
\[ f^*:\xtop\longrightarrow \shv(X_{\mathcal{U}}), \]
where $\shv(X_{\mathcal{U}})$ is spatial by construction. This map exhibits a unit transformation for the inclusion $\ltop^{\mathrm{sp}}_{0-\mathrm{loc}}\subset\ltop_{0-\mathrm{loc}}$ and we let $L_{\mathrm{sp}}$ denote the resulting left adjoint. We wish to lift this adjunction to the level of structured spaces. To achieve this, we will now detail a certain construction of left adjointable squares, which will prove quite useful in a variety of contexts. 
\begin{prop}\label{prop:pullbackcocartleftadj}
Let $p:\icat\rightarrow\icatd$ be a categorical fibration of \infcats and $g:\icatd'\rightarrow \icatd$ a functor admitting a left adjoint $f:\icatd\rightarrow\icatd'$. Suppose that the following condition is satisfied.
\begin{enumerate}
\item[$(*)$] For each object $C\in\icat$ and each unit transformation $(D',p(C)\overset{e}{\rightarrow} g(D'))\in \icatd'\times_{\icatd}\icatd_{p(C)/}$ at $p(C)\in \icatd$, there is a $p$-coCartesian lift $\overline{e}:C\rightarrow C'$ of $e$ starting at $C$.  
\end{enumerate}
 Let $\icat'=\icatd'\times_{\icatd}\icat$, then the pullback diagram
\[
\begin{tikzcd}
\icat'\ar[r,"h"]\ar[d,"p'"] & \icat\ar[d,"p"] \\
\icatd'\ar[r,"g"] & \icatd
\end{tikzcd}
\]
is horizontally left adjointable. Moreover, a pair $(C',C\rightarrow h(C'))\in \icat'\times_{\icat}\icat_{C/}$ is a unit transformation at an object $C\in \icat$ if and only if the following conditions are satisfied.
\begin{enumerate}[$(1)$]
\item The pair $(p'(C'),p(C)\rightarrow gp'(C'))\in \icatd'\times_{\icatd}\icatd_{p(C)/}$ is a unit transformation.
\item The morphism $C\rightarrow h(C')$ in $\icat$ is a $p$-coCartesian lift of $p(C)\rightarrow gp'(C')$.
\end{enumerate}
\end{prop}
\begin{proof}
The assumption that $g$ admits a left adjoint is equivalent to the assertion that for any $D\in\icatd$, the \infcat $\icatd'_{D/}=\icatd'\times_{\icatd}\icatd_{D/}$ has an initial object. To prove the proposition, we are required to show that for any object $C\in \icat$ the following holds.
\begin{enumerate}
    \item[$(\bullet)$] An object $(C',C\rightarrow h(C')) \in\icat'_{C/}=\icat'\times_{\icat}\icat_{C/}$ is initial if and only if $(p'(C'),p(C)\rightarrow gp'(C'))\in \icatd'_{p(C)/} $ is initial and $C\rightarrow h(C')$ is a $p$-coCartesian lift of $p(C)\rightarrow gp'(C')$.
\end{enumerate}
We have a pullback diagram
\[
\begin{tikzcd}
\icat'_{C/}\ar[r]\ar[d,"p'_C"] & \icat_{C/}\ar[d,"p_C"] \\
\icatd'_{p(C)/}\ar[r] & \icatd_{p(C)/}
\end{tikzcd}
\]
where $p_C$ is a categorical fibration. It follows from the assumption $(*)$ and \cite{HTT}, Proposition 2.4.3.1 that the lower horizontal functor carries each initial object $(D',p(C)\rightarrow g(D'))$ of $\icatd'_{p(C)/}$ to an object $p(C)\rightarrow g(D')\in \icatd_{p(C)/}$ for which the morphism 
\[
\begin{tikzcd}
& p(C) \ar[dl,"\mathrm{id}"'] \ar[dr] \\
p(C) \ar[rr] && g(D')
\end{tikzcd}
\]
from the initial object in $\icatd_{p(C)/}$ admits a $p_C$-coCartesian lift starting at the initial object $\mathrm{id}:C\rightarrow C\in \icat_{C/}$. The `if' direction of assertion $(\bullet)$ will therefore follow from the following one.
\begin{enumerate}
    \item[$(\bullet\bullet)$] Let 
     \[
\begin{tikzcd}
\icat'\ar[r,"h"]\ar[d,"p'"] & \icat\ar[d,"p"] \\
\icatd'\ar[r,"g"] & \icatd
\end{tikzcd}
\]
be a pullback diagram of \infcats where $p$ is a categorical fibration, $\icat$ and $\icatd$ have initial objects $C$ and $D$ that are preserved by $p$ and $\icatd'$ has an initial object $D'$. Then the pair $(D',C')\in \icat'$, where $C'\in\icat$ is the codomain of a coCartesian lift of the arrow $D\rightarrow g(D')$ with domain $C$, is an initial object of $\icat'$. 
\end{enumerate}
To prove this, we note that we have a diagram 
\[
\begin{tikzcd}
 \icat'_{(C',D')/}\ar[d] \ar[r]& \icat_{C'/} \ar[d]\ar[r] & \icat_{C/} \ar[d] \\
 \icatd'_{D'/} \ar[r] &\icatd_{g(D')/} \ar[r] &\icatd_{D/}
\end{tikzcd}
\]
where both squares are homotopy pullbacks for the Joyal model structure. Since the functors $\icat_{C/}\rightarrow \icat$, $\icatd_{D/}\rightarrow\icatd$ and $\icatd'_{D'/}\rightarrow\icatd'$ are trivial fibrations, the functor $\icat'_{(C',D')/}\rightarrow \icat'$ is a trivial fibration as well. The `only if' direction of $(\bullet)$ follows immediately from essential uniqueness of initial objects.
\end{proof}

\begin{defn}\label{defn:correctspectrum}
We let $L_{\mathrm{sp}}(\geodiff)$ denote a left adjoint to the inclusion $\ltop^{\mathrm{sp}}_{0-\mathrm{loc}}(\geodiff)\subset\ltop_{0-\mathrm{loc}}(\geodiff)$ provided by Proposition \ref{prop:pullbackcocartleftadj}. We define a functor as the composition
\[ \cinfty\mathsf{ring} \overset{\spec^{\geodiff}}{\longrightarrow} \ltop_{0-\mathrm{loc}}(\geodiff) \overset{L_{\mathrm{sp}}(\geodiff)}{\longrightarrow} \ltop^{\mathrm{sp}}_{0-\mathrm{loc}}(\geodiff).  \]
and denote it $\spec$.
\end{defn}
We will show that on objects, the functor $\spec$ coincides with the Archimedean spectrum of Construction \ref{cons:realspectrum}. First, we record the following feature of the functor $\spec$, which is particular to $\cinfty$ geometry.
\begin{prop}\label{prop:specpreservelimsp}
The composition
\[\cinfty\mathsf{ring}^{op}\overset{\spec}{\longrightarrow} \rtop^{\mathrm{sp}}_{0-\mathrm{loc}}(\geodiff)\longrightarrow \rtop^{\mathrm{sp}}_{0-\mathrm{loc}}\]
preserves small limits. 
\end{prop}
\begin{proof}
The composition factors via $\rtop^{\mathrm{sp}}_{0-\mathrm{loc}}(\geodiff)\rightarrow \rtop^{\mathrm{sp}}_{0-\mathrm{loc}}((\geodiff)_{\mathrm{disc}})\rightarrow \rtop^{\mathrm{sp}}_{0-\mathrm{loc}}$. The first functor preserves limits by Corollary \ref{cor:geodiffstableunderlimits} and the functor $p_{(\geodiff)_{\mathrm{disc}}}$ preserves limits and colimits, so we conclude as $\spec$ is a right adjoint.
\end{proof}
\begin{prop}\label{prop:spatialreflection}
 Let $A$ be a $\cinfty$-ring, then the algebraic morphism $\varphi_A^*$ of Construction \ref{cons:realspectrum} exhibits $\shv(\specr\,A)$ as the spatial reflection of $\shv(\cinfty\mathsf{ring}_{/A}^{op,\mathrm{ad}})$. 
\end{prop}
\begin{proof}
It follows from Construction \ref{cons:realspectrum} that the assignment $A\mapsto \varphi_A^*$ is functorial in $A$, so we have a natural transformation 
\[
\cinfty\mathsf{ring}\times\Delta^1\longrightarrow \ltop_{0-\mathrm{loc}},\quad\quad A\longmapsto (\shv(\cinfty\mathsf{ring}_{/A}^{op,\mathrm{ad}})\rightarrow \shv(\specr\,A)),
\] 
and since $\shv(\specr\,A)$ is spatial, the functor above induces a functor 
\[ f:\cinfty\mathsf{ring}\times\Delta^1\longrightarrow \ltop^{\mathrm{sp}}_{0-\mathrm{loc}}\]
carrying $(A,0)$ to the spatial reflection of $\shv(\cinfty\mathsf{ring}_{/A}^{op,\mathrm{ad}})$. We need to show that $f|_{A\times\Delta^1}$ is an equivalence for all $\cinfty$-rings $A$. If $A$ is finitely presented, then $f|_{A\times\Delta^1}$ is an equivalence; this follows from Proposition \ref{prop:germdetbasis} and the fact that finitely generated ideals of rings of $\cinfty$ functions on manifolds are germ determined. Therefore, it suffices to show that $f|_{\cinfty\mathsf{ring}\times\{0\}}$ and $f|_{\cinfty\mathsf{ring}\times\{1\}}$ preserve filtered colimits. We will show that both functors preserve all small colimits. The functor $f|_{\cinfty\mathsf{ring}\times\{0\}}$ is the composition 
\[ \cinfty\mathsf{ring} \overset{\spec^{\geodiff}}{\longrightarrow} \ltop_{0-\mathrm{loc}}(\geodiff)\overset{p_{\geodiff}}{\longrightarrow} \ltop_{0-\mathrm{loc}}\overset{L_{\mathrm{sp}}}{\longrightarrow}\ltop^{\mathrm{sp}}_{0-\mathrm{loc}} \]
which Proposition \ref{prop:pullbackcocartleftadj} identifies with the composition
\[ \cinfty\mathsf{ring} \overset{\spec}{\longrightarrow} \ltop^{\mathrm{sp}}_{0-\mathrm{loc}}(\geodiff){\longrightarrow} \ltop^{\mathrm{sp}}_{0-\mathrm{loc}} \]
which preserves small colimits by Proposition \ref{prop:specpreservelimsp}. For the functor $f|_{\cinfty\mathsf{ring}\times\{1\}}$, it suffices to show that the functor $A\mapsto \specr\,A$ sends colimits of $\cinfty$-rings to limits of topological spaces. From the definition of $\specr\,A$, this is clear on the underlying sets, and it is a straightforward verification that $\specr\,\colim_{k\in K} A_k$ carries the initial topology with respect to the maps $\{\specr\,\colim_{k\in K} A_k\rightarrow \specr\,A_k\}_{k\in K}$ for any small diagram $K\rightarrow \cinfty\mathsf{ring}$. 
\end{proof}

\begin{prop}\label{comparingspectra}
Identifying the 1-categories $\rtop^{\mathrm{sp}}_{0-\mathrm{loc}}(\geodiff)$ and $\mathsf{Top}^{\loc}(\cinfty\mathsf{ring})$ as in Proposition \ref{prop:comparisonringedspaces}, the object $\spec\,A$ is equivalent to the Archimedean spectrum of Construction \ref{cons:realspectrum} of $A$. The functor $L_{\mathrm{sp}}(\geodiff)$ carries the 1-category $\sch_{0}(\geodiff)$ of $0$-localic $\geodiff$-schemes onto the 1-category $\cinfty\sch$ of $\cinfty$-schemes and restricts to an equivalence on the full subcategory spanned by those $0$-localic $\geodiff$-schemes that are spatial.
\end{prop}
\begin{proof}
This is an immediate consequence of Construction \ref{cons:realspectrum}, Proposition \ref{prop:pullbackcocartleftadj} and Proposition \ref{prop:spatialreflection}.
\end{proof}

We now study some basic features of the Archimedean spectrum, starting with some topological properties of the space $\specr\,A$.
\begin{lem}
 Let $A$ be a $\cinfty$-ring and consider for each $a\in A$ the map $\ev_a:\specr\,A\rightarrow \R$. Then $\specr\,A$ carries the initial topology with respect to the map
 \[ \prod_{a\in A}\ev_a:\specr\,A\longrightarrow \prod_{a\in A}\R. \]
\end{lem}
\begin{proof}
The topology on $\specr\,A$ is the coarsest one which contains the basis open sets $\ev_a^{-1}(\R\setminus\{0\})$, so it is coarser than the initial topology for the map $\prod_{a\in A}\ev_a$. To see that it also finer, we note that for any open $U\in \R$ and any $a\in A$, the set $\ev_a^{-1}(U)$ coincides with the set $\ev_{(\chi_U)_*(a)}(\R\setminus\{0\})$ where $\chi_U$ is a characteristic function for $U$ and $(\chi_U)_*:A\rightarrow A$ is the $\cinfty$ structure map.
\end{proof}
\begin{cor}
For any $\cinfty$-ring $A$, the topological space $\specr\,A$ is completely regular Hausdorff.    
\end{cor}
\begin{proof}
The map $\prod_{a\in A}\ev_a$ is injective and the properties of being Hausdorff and completely regular are hereditary (that is, stable under initial topologies).    
\end{proof}
We proceed by investigating the structure sheaf $\Of_{\specr\,A}$. 
\begin{lem}
Let $A$ be a $\cinfty$-ring and let $a\in A$ be an element with associated open set $U_a=\Hom_{\cinfty\mathsf{ring}}(A[a^{-1}],\R)\subset\specr\,A$. Let $\mathcal{M}_{U_a}$ be the category $\cinfty\mathsf{ring}^{op,\mathrm{ad}}_{/A}\times_{\mathcal{B}(A)}\{U_a\}$ of localizations $A\rightarrow A[b^{-1}]$ that determine the same subset $U_b=U_a$. Then $\mathcal{M}_{U_a}$ is filtered and there is an equivalence of $\cinfty$-rings 
\[ \widetilde{\Of_{\specr\,A}}(U_a) \simeq \underset{A\rightarrow A[b^{-1}]\in \mathcal{M}^{op}_{U_a}}{\colim} A[b^{-1}]\]
natural in $A[a^{-1}]$.
\end{lem}
\begin{proof}
The category $\mathcal{M}_{U_a}^{op}$ is filtered because $\mathcal{M}_{U_a}$ admits finite limits. Let $\mathcal{M}^{op}$ be the categorical mapping cylinder associated to the functor $\Hom_{\cinfty\mathsf{ring}}(\_,\R):\cinfty\mathsf{ring}^{\mathrm{ad}}_{/A}\rightarrow \mathcal{B}(A)^{op}$. Concretely, consider the category $\mathcal{M}$ defined as follows.
\begin{enumerate}
    \item[$(O)$] An object of $\mathcal{M}$ is either a localization $A\rightarrow A[a^{-1}]$ or a basis open $U_b\subset\specr\,A$.
    \item[$(M)$] The set $\Hom_{\mathcal{M}^{op}}(A[a^{-1}],A[b^{-1}])$ contains one element if $a$ is invertible in $A[b^{-1}]$ and is empty otherwise. The set $\Hom_{\mathcal{M}^{op}}(U_a,U_b)$ contains one element if $U_b\subset U_a$ and is empty otherwise. The set $\Hom_{\mathcal{M}^{op}}(A[a^{-1}],U_b)$ contains one element if $U_b\subset \Hom_{\cinfty\mathsf{ring}}(A[a^{-1}],\R)$ and is empty otherwise. The set $\Hom_{\mathcal{M}^{op}}(U_A,A[b^{-1}])$ is always empty.
\end{enumerate}
Then the functor $\Hom_{\cinfty\mathsf{ring}}(\_,\R)$ factors as $i_0:\cinfty\mathsf{ring}^{\mathrm{ad}}_{/A}\subset \mathcal{M}^{op}\rightarrow \mathcal{B}(A)^{op}$ where the second map is a retraction of the full subcategory inclusion $i_1:\mathcal{B}(A)^{op}\subset\mathcal{M}^{op}$ and the presheaf $\widetilde{\Of_{\specr\,A}}$ is the restriction to the full subcategory $\mathcal{B}(A)^{op}\subset\mathcal{M}^{op}$ of a left Kan extension as in the diagram 
\[
\begin{tikzcd}
\cinfty\mathsf{ring}^{\mathrm{ad}}_{A/} \ar[d,hook,"i_0"]\ar[r,"\widetilde{\Of_{\speco\,A}}"] &[2em] \cinfty\mathsf{ring}\\
\mathcal{M}^{op} \ar[ur,"i_{0!}\widetilde{\Of_{\spec\,A}}"']
\end{tikzcd}
\]
Now it suffices to observe that for any $U_a\in\mathcal{B}(A)$, the category $\mathcal{M}_{U_a}^{op}$ is a left cofinal subcategory of $\cinfty\mathsf{ring}^{\mathrm{ad}}_{/A}\times_{\mathcal{M}^{op}} \mathcal{M}^{op}_{/U_a}$.
\end{proof}
\begin{lem}\label{lem:stalklocalization}
Let $A$ be a $\cinfty$-ring and let $A\rightarrow\R$ be an $\R$-algebra morphism corresponding to an $\R$-point $p:*\rightarrow \specr\,A$, which induces a map of local $\cinfty$-rings $p^*\Of_{\specr\,A}\rightarrow \R$. Then the composition 
\[A\longrightarrow \Gamma(\Of_{\specr\,A})\longrightarrow p^*\Of_{\specr\,A}\longrightarrow \R\]
exhibits $p^*\Of_{\specr\,A}\rightarrow \R$ as a unit transformation for the localization $L_{\R}$ at $p:A\rightarrow \R$.
\end{lem}
\begin{proof}
By definition of the topology on $\specr\,A$ and the presheaf $\widetilde{\Of_{\specr\,A}}$, the map $A\rightarrow p^*\Of_{\specr\,A}$ is a filtered colimit $A\rightarrow \colim_{a\in A,p(a)\neq 0}A[a^{-1}]$ over admissible morphisms and thus ind-admissible. The map $p^*\Of_{\specr\,A}\rightarrow \R$ is a (local) morphism between local $\cinfty$-rings, so we conclude.    
\end{proof}
The product of stalks 
\[\prod_{p\in\specr\,A}p^*:\shv_{\cinfty\mathsf{ring}}(\specr\,A)\longrightarrow \prod_{p\in\specr\,A} \cinfty\mathsf{ring}\] 
is conservative and preserves colimits and split equalizers; it follows from Lemma \ref{lem:stalklocalization} and Beck's monadicity theorem that the sheafification of the presheaf $\widetilde{\Of_{\specr\,A}}$ is the map to the equalizer in the diagram
\[ \begin{tikzcd}
   \widetilde{\Of_{\specr\,A}}\ar[r] & \mathrm{Eq}\ar[r]&  \prod_{p\in \specr\,A}p_*A_p \ar[r,shift right]\ar[r,shift left] & \prod_{p\in \specr\,A}p_*\left(\prod_{p'\in \specr\,A}p^{\prime *}A_{p'}\right)
\end{tikzcd}  \]
Since the equalizer is subsheaf of $\prod_{p\in \specr\,A}p_*A_p$, we deduce the following.
\begin{prop}[Godement resolution]\label{prop:godementsheaf}
Consider the sheaf $\prod_{p\in \specr\,A}p_*A_p$ which admits a canonical map from $\widetilde{\Of_{\specr\,A}}$ carrying a section to its various stalks. Then the sheafification 
\[  \widetilde{\Of_{\specr\,A}}\longrightarrow \Of_{\specr\,A}\longrightarrow \prod_{p\in \specr\,A}p_*A_p\]
identifies the sheaf $\Of_{\specr\,A}$ with the smallest subsheaf of $\prod_{p\in \specr\,A}p_*A_p$ that contains the image of the map $\widetilde{\Of_{\specr\,A}}\rightarrow \prod_{p\in \specr\,A}p_*A_p$. Concretely, for any open $U_a\subset \specr\,A$, the domain of the map $\widetilde{\Of_{\specr\,A}}(U_a)\rightarrow \Of_{\specr\,A}(U_a)$ is identified with
\[\{(s_p)_{p\in U_a};\,\text{for all }p\in U_a\text{ there exist }b\in A,\,t\in A[b^{-1}]\text{ such that }U_b\subset U_a\text{ and }{s_p'}={t_p'}\text{ for all }p'\in U_b\}\subset \prod_{p\in U_a}A_{p}. \]
\end{prop}
\begin{cor}\label{cor:presheafalreadysheaf}
Let $A$ be a finitely generated $\cinfty$-ring of the form $\cinfty(\R^n)/I$.
\begin{enumerate}[$(1)$]
    \item The map $A\rightarrow \Gamma(\Of_{\specr\,A})$ coincides with projection $\cinfty(\R^n)/I\rightarrow \cinfty(\R^n)/\tilde{I}$, where $\tilde{I}$ is the smallest germ determined ideal containing $I$.
    \item If $A$ is finitely presented, then the presheaf $\widetilde{\Of_{\specr\,A}}$ is already a sheaf.
\end{enumerate}
\end{cor}
\begin{proof}
The assertion $(1)$ is an immediate consequence of the following assertions.
\begin{enumerate}[$(a)$]
\item The map $A\rightarrow \Gamma(\Of_{\specr\,A})$ is a surjection.
\item For any $p\in \specr\,A$ corresponding to an evaluation map $\ev_p:\cinfty(\R^n)/I\rightarrow \R$, the localization $\cinfty(\R^n)/I\rightarrow (\cinfty(\R^n)/I)_p$ factors through $\cinfty(\R^n)/\tilde{I}$ and induces an isomorphism $(\cinfty(\R^n)/I)_p\cong (\cinfty(\R^n)/\tilde{I})_p$.
\item If $I$ is germ determined, then the map $A\rightarrow \prod_{p\in \specr\,A}A_p$ is an injection.
\end{enumerate}
To show $(a)$, we use Proposition \ref{prop:godementsheaf} and Lemma \ref{lem:localizationrn} to deduce that an element of $\Gamma(\Of_{\specr\,A})$ is a tuple $(f_p)_{p\in Z(I)}$ for which each $p$ has a neighbourhood $U\subset Z(I)$ such that $f_{p'}=g_{Up'}$ for some $g_U\in \cinfty(U)/I|_U$ for all $p'\in U$. We can find an element $g\in\cinfty(\R^n)/I$ for which $f_p=g_p$ for all $p\in Z(I)$ via a partition of unity. The claims $(b)$ and $(c)$ follow immediately from Lemma \ref{lem:localizationrn} and the definition of being germ determined. To show $(2)$, we note that if $I$ is finitely generated, $I|_U$ is germ determined for any open $U\subset\R^n$, and the preceding arguments show that the map $\cinfty(U)/I\rightarrow \Of_{\specr\,A}(U)$ is an isomorphism.
\end{proof}
The following result is a consequence of the $\cinfty$ regularity of an arbitrary affine $\cinfty$-scheme.
\begin{prop}\label{prop:gammastalk}
Let $A$ be a $\cinfty$-ring, and let $p:\Gamma(\Of_{\specr\,A})\rightarrow\R$ be an $\R$-algebra morphism. Let $p':A\rightarrow\R$ be the composition of the unit transformation at $A$ with $p$, then the diagram 
\[
\begin{tikzcd}
 A\ar[r] \ar[d] &    \Gamma(\Of_{\specr\,A})\ar[d] \\
 A_{p'}\ar[r] & \Gamma(\Of_{\specr\,A})_p
\end{tikzcd}
\]
is a pushout of $\cinfty$-rings and the lower horizontal map is an isomorphism. Moreover, the map $p$ equals the map $\Gamma(\Of_{\specr\,A})\rightarrow A_{p'}\rightarrow \R$.
\end{prop}
\begin{proof}
The map $A\rightarrow A_{p'}$ is ind-admissible so the map $\Gamma(\Of_{\specr\,A})\rightarrow A_{p'}\oinfty_A\Gamma(\Of_{\specr\,A})$ is ind-admissible as well. Thus, to see that the latter map exhibits a unit transformation for the localization $L_{\R}$, it suffices to argue that the ring $A_{p'}\oinfty_A\Gamma(\Of_{\specr\,A})$ is local, since it comes equipped with a map to $\R$. Since the left vertical map is a regular epimorphism by Corollary \ref{cor:locepi}, it follows from Corollary \ref{cor:algeffepi} that it suffices to check that the pushout $A_{p'}\otimes_A\Gamma(\Of_{\specr\,A})$ is a local ring. We may identify this ring with the quotient of $\Gamma(\Of_{\specr\,A})$ by the following equivalence relation, identifying global sections with families of stalks using the Godement resolution.
\begin{enumerate}
\item[$(*)$]Write $(s_x)_x$ with $x\in \specr\,A$ for global sections in $\Gamma(\Of_{\specr\,A})$, then $(s_x)_x\sim (t_x)_x$ if and only if $a_x(s_x-t_x)_x=0$ for all $x\in\specr\,A$, for some $a\in A$ for which $a_p\neq 0$, where $a_x$ denotes the stalk of $a$ at $x$.
\end{enumerate}
It suffices to show that the composition $A\rightarrow\Gamma(\Of_{\specr\,A})\rightarrow  A_{p'}\otimes_A\Gamma(\Of_{\specr\,A})$ is a surjection; then $A_{p'}\rightarrow A_{p'}\otimes_A\Gamma(\Of_{\specr\,A})$ is also a surjection which forces $A_{p'}\otimes_A\Gamma(\Of_{\specr\,A})$ to be a local ring. Let $(s_x)_x$ be a global section, then there is an $a\in A$ with $p'(a)\neq 0$ and an element $b\in A[a^{-1}]$ such that $s_x=b_x$ for all $x\in U_a$. Since $A$ is a filtered colimit of its finitely generated sub $\cinfty$-rings, we can find
\begin{enumerate}[$(1)$]
    \item A monomorphism $\psi:\cinfty(\R^n)/I\hookrightarrow A$ whose image contains $a$,
    \item An open set $U\subset Z(I)$ and a characteristic function $\chi_U$ for $U$ in $\cinfty(\R^n)/I$ such that $\psi(\chi)=a$, inducing a pushout diagram 
    \[ 
    \begin{tikzcd}
    \cinfty(\R^n)/I \ar[d,hook] \ar[r,"\psi"] & A\ar[d] \\
    \cinfty(U)/I \ar[r,"\phi"] & A[a^{-1}]
    \end{tikzcd}
    \]
    \item An element $g\in \cinfty(U)/I$ such that $\phi(g)=b$.
\end{enumerate}
The map 
\[\cinfty(U)/I\longrightarrow A[a^{-1}]\overset{p'}{\longrightarrow} \R\]
corresponds to evaluation at some point $p'\in Z(I)$. Now we choose open sets $W,V,Y$ such that $p'\subset W\subset \overline{W}\subset V\subset \overline{V} \subset Y \subset \overline{Y}\subset U\subset  Z(p)$ and a bump function $\lambda$ on $\R^n$ such that
\[\lambda|_{\overline{V}} =1,\quad\quad \lambda|_{\R^n\setminus Y}=0.\] 
Abusing notation, we pass to the quotient by $I$ and view $\lambda$ as an element of $\cinfty(\R^n)/I$. Let $\chi_W$ be a characteristic function for $W$, then by construction, $\psi(\chi_W)_x(s_x-\psi(\lambda)_xs_x)=0$ for all $x\in\specr\,A$, so $(s_x)_x\sim (\psi(\lambda)_xs_x)_x$. The function $\lambda g$ on $U$ can be extended by 0 to an element $\lambda \tilde{g}\in \cinfty(\R^n)/I$. Since $\psi(\lambda)_x=0$ for all $x\notin U_a$, we have $\psi(\lambda)_xs_x=0=\psi(\lambda \tilde{g})_x=\psi(\lambda)_x\psi(\tilde{g})_x$ for all $x\notin U_a$, and for all $x\in U_a$ we have $\psi(\lambda)_xs_x=\psi(\lambda)_xb_x=\psi(\lambda)_x\psi(\tilde{g})_x=\psi(\lambda\tilde{g})_x$. Thus, $(s_x)_x=(\psi(\lambda\tilde{g})_x)_x$ and we conclude that the map $A\rightarrow A_{p'}\otimes_A\Gamma(\Of_{\specr\,A})$ is a surjection. From the fact that the diagram in the lemma is now a pushout we deduce that $p$ factors through $\Gamma(\Of_{\specr\,A})\rightarrow A_{p'}$. It follows that we have a retraction
\[ A_{p'}\longrightarrow \Gamma(\Of_{\specr\,A})_p \longrightarrow A_{p'},\]
but we have just argued that the first map is a surjection so we conclude.
\end{proof}
\begin{prop}[Joyce \cite{Joy2}]\label{prop:joyce1}
Let $A$ be a $\cinfty$-ring. Then $\spec$ carries the unit map $A\rightarrow \Gamma(\Of_{\specr\,A})$ to an equivalence.    
\end{prop}
\begin{proof}
Proposition \ref{prop:gammastalk} guarantees that the map $\specr\,\Gamma(\Of_{\specr\,A})\rightarrow\specr\,A$ is a bijection. 
It suffices to show that the map on presheaves
\[ \widetilde{\Of_{\specr\,A}}\longrightarrow \widetilde{\Of_{\Gamma(\Of_{\specr\,A})}}\]
induces an equivalence on stalks. This follows from Lemma \ref{lem:stalklocalization} and Proposition \ref{prop:gammastalk}.
\end{proof}
Joyce introduces the following terminology.
\begin{defn}\label{defn:completerings}
A $\cinfty$-ring $A$ is \emph{complete} if the unit map $A\rightarrow \Gamma(\Of_{\specr\,A})$ is an isomorphism. The full subcategory of $\cinfty\mathsf{ring}$ spanned by complete $\cinfty$-rings is denoted $\cinfty\mathsf{ring}_{\mathrm{cplt}}$    
\end{defn}
\begin{cor}
The functor $\Gamma:\caff\rightarrow\cinfty\mathsf{ring}$ is a fully faithful embedding onto the full subcategory $\cinfty\mathsf{ring}_{\mathrm{cplt}}$. The full subcategory $\cinfty\mathsf{ring}_{\mathrm{cplt}}$ is reflective and for each $A\in \cinfty\mathsf{ring}$, the unit map $A\rightarrow \Gamma(\Of_{\specr\,A})$ exhibits $\Gamma(\Of_{\specr\,A})$ as the completion of $A$.
\end{cor}
We will shortly deal with derived versions of $\cinfty$-rings and their Archimedean spectra. In this setting, we need to impose a stronger condition than completeness to ensure that passage from derived $\cinfty$-rings to their spectra is well behaved. 
\begin{defn}\label{defn:geometricring}
Let $A$ be a $\cinfty$-ring, then we say that $A$ is \emph{Lindel\"{o}f} if $\specr\,A$ is a Lindel\"{o}f topological space. We say that $A$ is \emph{geometric} if $A$ is complete, that is, the unit map $A\rightarrow\Gamma(\Of_{\specr\,A})$ is an equivalence, and $A$ is Lindel\"{o}f.
\end{defn}
The following is a basic topological result. For a proof, see Theorem 4.40 of \cite{Joy2}. A very general result of this kind is proven by Kriegl and Michor as Theorem 16.10 of \cite{krieglmichor}.
\begin{prop}
Let $(X,\Of_{X})$ be a $\cinfty$-scheme such that $X$ is Lindel\"{o}f and Hausdorff. Then $X$ is paracompact and the sheaf $\Of_{X}$ is fine, that is, it admits partitions of unity for arbitrary covers.  
\end{prop}
\begin{rmk}
As a consequence of the previous proposition, sheaves of discrete $\Of_{\specr\,A}$-modules are also fine, hence acyclic. It follows that for left bounded sheaves of $\Of_{\specr\,A}$-modules, the openwise homotopy groups are already sheaves (see Proposition \ref{prop:parahausdorffpresheaf}).
\end{rmk}

Below we state without proof a useful result due to Joyce, which provides a large supply of affine $\cinfty$-schemes. 
\begin{defn}
A locally $\cinfty$-ringed space $(X,\Of_X)$ is \emph{$\cinfty$-regular} if $X$ carries the initial topology with respect to the canonical map $X\rightarrow \specr\,\Gamma(\Of_X)$.     
\end{defn}
\begin{thm}[\cite{Joy2}]\label{thm:joyce}
Let $(X,\Of_X)$ be a locally $\cinfty$-ringed space for which the following conditions are satisfied.
\begin{enumerate}[$(1)$]
    \item $X$ is Hausdorff and Lindel\"{o}f.
    \item $(X,\Of_X)$ is $\cinfty$-regular.
\end{enumerate}
Then the unit map $(X,\Of_{X})\rightarrow (\specr\,\Gamma(\Of_{X}),\Of_{\specr\,\Gamma(\Of_X)})$ is an isomorphism.
\end{thm}
\begin{rmk}
In \cite{moerdijkreyesvanque}, Moerdijk, van Qu\^{e} and Reyes consider the spectrum functor for the geometry $\geodiff^{\mathrm{fin}}$. Since the Grothendieck topology on $\mathrm{Pro}(\geodiff^{\mathrm{fin}})^{\mathrm{ad}}_{/A}$ is finitary by stipulation, the sheaf topos $\shv(\mathrm{Pro}(\geodiff^{\mathrm{fin}})^{\mathrm{ad}}_{/A})$ is spatial by Deligne's completeness theorem and thus corresponds to a coherent topological space. This space can be identified with the \emph{$\cinfty$-radical prime ideal spectrum} of $A$: consider for an ideal $I\subset A$, the ideal $\sqrt[\infty]{I}$ of those elements of $a$ that are carried to $0$ by every $p\in \specr\,A$ that lies in $Z(I)$. Then $I$ is \emph{$\cinfty$-radical} if $I=\sqrt[\infty]{I}$. The spectrum $\spec^{\geodiff^{\mathrm{fin}}}$ is the collection of $\cinfty$-radical prime ideal endowed with the Zariski topology and equipped with a sheaf of $\cinfty$-rings whose stalks are local as commutative rings. It is not hard to see that for any $\cinfty$-ring $A$, the presheaf $\widetilde{\Of_{\speco^{\mathrm{fin}}\,A}}$ is already a sheaf. As a consequence, the spectrum functor $\spec^{\geodiff^{\mathrm{fin}}}$ is fully faithful. In particular, the ordinary category of smooth manifolds embeds fully faithfully in the category of coherent topological spaces equipped with sheaves of $\cinfty$-rings with local stalks. By the general theory of geometries, the functor $\spec^{\geodiff}$ is equivalent to the composition $\spec^{\geodiff}_{\geodiff^{\mathrm{fin}}}\circ\spec^{\geodiff^{\mathrm{fin}}}$, where $\spec^{\geodiff}_{\geodiff^{\mathrm{fin}}}$ is the relative spectrum associated to the transformation $\geodiff^{\mathrm{fin}}\rightarrow\geodiff$. 
\end{rmk}

\subsubsection{Modules of $\cinfty$-rings and locally $\cinfty$-ringed spaces}
We now discuss (sheaves of) modules for $\cinfty$-rings, which are simply (sheaves of) modules for their underlying algebras. 
\begin{defn}
Let $\mathrm{Mod}^{\rmalg}$ be the category defined as follows.
\begin{enumerate}
    \item[$(O)$] Objects are pairs $(A,M)$ of a commutative $\R$-algebra together with a (left and right) $A$-module $M$.
    \item[$(M)$] Morphisms $(A,M)\rightarrow (B,N)$ are pairs of a map $A\rightarrow B$ and a map $B\otimes_AM\rightarrow N$ of $B$-modules.
\end{enumerate}
The obvious functor $p:\mathrm{Mod}^{\rmalg}\rightarrow\mathsf{CAlg}^0_{\R}$ is a presentable fibration: a map $f:(A,M)\rightarrow (B,N)$ is $p_{\mathrm{Mod}}$-coCartesian if the map $B\otimes_AM\rightarrow N$ is an isomorphism, and the map is $p_{\mathrm{Mod}}$-Cartesian if the adjoint map $M\rightarrow N$ is an isomorphism. We define a category $\mathrm{Mod}$ as the cone in the pullback diagram 
\[
\begin{tikzcd}
\mathrm{Mod}\ar[d,"p_{\mathrm{Mod}}"]\ar[r] & \mathrm{Mod}^{\rmalg}\ar[d,"p"] \\
\cinfty\mathsf{ring} \ar[r,"(\_)^{\rmalg}"] & \mathsf{CAlg}^0_{\R}.
\end{tikzcd}
\]
Its objects are pairs $(A,M)$ of a $\cinfty$-ring $A$ with an $A^{\rmalg}$-module. For $A$ a $\cinfty$-ring, we let $\mathrm{Mod}_A$ denote the pullback $\{A\}\times_{\cinfty\mathsf{ring}}\mathrm{Mod}$, the category of (left and right) $A^{\rmalg}$-modules. Define a category as follows.
\begin{enumerate}
    \item[$(O)$] Objects are triples $(X,\Of_X,\F)$ of a sober topological space together with a an object $(\Of_X,\F)\in \shv_{\mathrm{Mod}}(X)$ such that the underlying sheaf of $\cinfty$-rings $\Of_X$ determines a locally $\cinfty$-ringed space $(X,\Of_X)$.
    \item[$(M)$] Morphisms are triples $(f,\alpha,\beta):(X,\Of_X,\F)\rightarrow (Y,\Of_Y,\mathcal{G})$ of a continuous map $f:X\rightarrow Y$, a morphism $\alpha:f^*\Of_Y\rightarrow\Of_X$ of sheaves of $\cinfty$-rings (which is automatically a local morphism on stalks), and a morphism $\beta:f^*\F\otimes_{f^*\Of_Y}\Of_X\rightarrow \mathcal{G}$ of sheaves of $\Of_X$-modules. 
\end{enumerate}
We denote this category by $\mathsf{Top}^{\loc}(\mathrm{Mod})$. It comes equipped with a Cartesian fibration \[p_{\mathsf{ModTop}}:\mathsf{Top}^{\loc}(\mathrm{Mod})\longrightarrow \mathsf{Top}^{\loc}(\cinfty\mathsf{ring}).\]    
A map $(X,\Of_X,\F)\rightarrow (Y,\Of_Y,\mathcal{G})$ is $p_{\mathsf{ModTop}}$-Cartesian if and only if the induced map $f^*\F\otimes_{f^*\Of_Y}\Of_X\rightarrow \mathcal{G}$ is an isomorphism. $p_{\mathsf{ModTop}}$ is also a coCartesian fibration: a map $(X,\Of_X,\F)\rightarrow (Y,\Of_Y,\mathcal{G})$ is $p_{\mathsf{ModTop}}$-coCartesian precisely if the adjoint map $f_*\F\rightarrow \mathcal{G}$ is an isomorphism of sheaves of $\Of_Y$-modules. For $(X,\Of_X)$ a locally $\cinfty$-ringed space, we denote by $\mathrm{Mod}_{\Of_X}$ the fibre at $(X,\Of_X)$ of the functor 
\[p_{\mathsf{ModTop}}:\mathsf{Top}^{\loc}(\mathrm{Mod})^{op}\longrightarrow \mathsf{Top}^{\loc}(\cinfty\mathsf{ring})^{op}\]
on \emph{opposite} categories.
\end{defn}
There is an obvious \emph{global sections functor} 
\[\Gamma:\mathsf{Top}^{\loc}(\mathrm{Mod})^{op}\longrightarrow \mathrm{Mod}\]
which fits into a commuting diagram
\begin{equation}\label{eq:globsec}
     \begin{tikzcd}
     \mathrm{Mod}\ar[d,"p_{\mathrm{Mod}}"] &    \mathsf{Top}^{\loc}(\mathrm{Mod})^{op}\ar[d,"p_{\mathsf{ModTop}}"]\ar[l,"\Gamma"'] \\
    \cinfty\mathsf{ring} & \mathsf{Top}^{\loc}(\cinfty\mathsf{ring})^{op}. \ar[l,"\Gamma"']   
    \end{tikzcd}
\end{equation}
\begin{thm}\label{thm:modulespectrum}
The global sections functor on modules admits a left adjoint $\mathrm{M}\spec$ with the following properties.
\begin{enumerate}[$(1)$]
    \item There is a canonical natural isomorphism making the diagram 
    \[\begin{tikzcd}
    \mathrm{Mod}\ar[d,"p_{\mathrm{Mod}}"]\ar[r,"\mathrm{M}\spec"] &[3em]    \mathsf{Top}^{\loc}(\mathrm{Mod})^{op}\ar[d,"p_{\mathsf{ModTop}}"] \\
    \cinfty\mathsf{ring}\ar[r,"\spec"] &[3em] \mathsf{Top}^{\loc}(\cinfty\mathsf{ring})^{op}
    \end{tikzcd}
    \]
    commute.
    \item Let $\mspec_A:\mathrm{Mod}_A\rightarrow\mathrm{Mod}_{\Of_{\specr\,A}}$ be the induced functor on the fibre over $A$ and $\Of_{\specr\,A}$, then $\mspec_A$ is naturally equivalent to the composition 
    \[ \mathrm{Mod}_A \longrightarrow \mathrm{Mod}_{\widetilde{\Of_{\specr\,A}}}\overset{L}{\longrightarrow}   \mathrm{Mod}_{\Of_{\specr\,A}}\]
    where the second functor is a sheafification and the first functor carries $M$ to the presheaf of $\widetilde{\Of_{\specr\,A}}$-modules given by 
    \[ U_a\longmapsto M\otimes_A\widetilde{\Of_{\specr\,A}}(U_a). \]
    \item For any map $A\rightarrow B$ of $\cinfty$-rings inducing a map $f:\specr\,B\rightarrow\specr,A$, there is a canonical natural isomorphism making the diagram 
    \[
    \begin{tikzcd}
    \mathrm{Mod}_A\ar[d] \ar[r,"\_\otimes_AB"] &[6em]    \mathrm{Mod}_{B}\ar[d]\\
\mathrm{Mod}_{\Of_{\specr\,A}} \ar[r,"f^*\_\otimes_{f^*\Of_{\specr\,A}}\Of_{\speco\,B}"] &[6em] \mathrm{Mod}_{\Of_{\speco\,B}},  
\end{tikzcd}
    \]
    commute.
\end{enumerate}
\end{thm}
\begin{proof}
We claim that it suffices to provide a collection of functors $\mspec\,A:\mathrm{Mod}_A\rightarrow\mathrm{Mod}_{\Of_{\specr\,A}}$ that have the property specified in point $(2)$; these functors will then assemble to yield the desired functor $\mathrm{M}\spec$ which has the properties $(1)$ and $(3)$. To see this, we consider the categorical mapping cylinders associated to the global sections functors. Concretely, define a category $\mathcal{M}$ as follows.
\begin{enumerate}
    \item[$(O)$] An object is either an object of $\mathsf{Top}^{\loc}(\cinfty\mathsf{ring})^{op}$ or of $\cinfty\mathsf{ring}$.
    \item[$(M)$] Let $A,B\in \cinfty\mathsf{ring}$ and $(X,\Of_X),(Y,\Of_Y)\in\mathsf{Top}^{\loc}(\cinfty\mathsf{ring})$. Morphism sets are given by
    \begin{align*}
        \Hom_{\mathcal{M}}(A,B)&=\Hom_{\cinfty\mathsf{ring}}(A,B),&\Hom_{\mathcal{M}}(A,(X,\Of_X))&=\Hom_{\cinfty\mathsf{ring}}(A,\Gamma(\Of_X)),\\ 
        \Hom_{\mathcal{M}}((X,\Of_X),(Y,\Of_Y))&=\Hom_{\mathsf{Top}^{\loc}(\cinfty\mathsf{ring})}((Y,\Of_Y),(X,\Of_X)),&\Hom_{\mathcal{M}}((X,\Of_X),A)&=\emptyset
    \end{align*} 
\end{enumerate}
By construction, there is a functor $q:\mathcal{M}\rightarrow\Delta^1$ which is a Cartesian fibration associated to $\Gamma$, so that we can identify the fibres of $q$ over $0$ and $1$ with $\cinfty\mathsf{ring}$ and $\mathsf{Top}^{\loc}(\cinfty\mathsf{ring})^{op}$ respectively. Similarly, we have a Cartesian fibration $q_{\mathrm{Mod}}:\mathcal{M}_{\mathrm{Mod}}\rightarrow\Delta^1$ associated to the global sections functor on modules. The commuting diagram \eqref{eq:globsec} determines the horizontal functor in the diagram
\[
\begin{tikzcd}
\mathcal{M}_{\mathrm{Mod}} \ar[rr,"\Phi"]\ar[dr,"q_{\mathrm{Mod}}"']&&\mathcal{M} \ar[dl,"q"]\\
& \Delta^1,
\end{tikzcd}
\]
which carries $q_{\mathrm{Mod}}$-Cartesian edges to $q$-Cartesian edges. We first claim that $\Phi$ is also a Cartesian fibration. According to \cite{luriegoodwillie}, Lemma 1.4.14\footnote{Note that the statement of this lemma is missing the condition that $\Phi$ preserves Cartesian edges}, this amounts to the following two assertions.
\begin{enumerate}[$(a)$]
    \item The functors $p_{\mathrm{Mod}}:\mathrm{Mod}\rightarrow\cinfty\mathsf{ring}$ and $p_{\mathsf{ModTop}}:\mathsf{Top}^{\loc}(\mathrm{Mod})^{op}\longrightarrow \mathsf{Top}^{\loc}(\cinfty\mathsf{ring})^{op}$ are Cartesian fibrations. 
    \item In the commuting diagram 
    \[
    \begin{tikzcd}
     \mathrm{Mod}\ar[d,"p_{\mathrm{Mod}}"] &    \mathsf{Top}^{\loc}(\mathrm{Mod})^{op}\ar[d,"p_{\mathsf{ModTop}}"]\ar[l,"\Gamma"'] \\
    \cinfty\mathsf{ring} & \mathsf{Top}^{\loc}(\cinfty\mathsf{ring})^{op}, \ar[l,"\Gamma"']   
    \end{tikzcd}
    \]
    the global sections functors carry $p_{\mathsf{ModTop}}$-Cartesian morphisms to $p_{\mathrm{Mod}}$-Cartesian morphisms.
\end{enumerate}
Both $(a)$ and $(b)$ are immediate from the construction of $p_{\mathrm{Mod}}$ and $p_{\mathsf{ModTop}}$. Since $q$ is a Cartesian fibration associated to $\Gamma$ which admits the left adjoint $\spec$, the functor $q$ is also a coCartesian fibration. Invoking \cite{luriegoodwillie}, Lemma 1.4.14 again, we deduce that the following are equivalent.
\begin{enumerate}[$(i)$]
    \item $\Phi$ is a coCartesian fibration, which, in view of the fact that $\Phi$ is a Cartesian fibration, amounts to the assertion that for each map $X\rightarrow Y$ in $\mathcal{M}$, the Cartesian pullback functor $\Phi^{-1}(Y)\rightarrow\Phi^{-1}(X)$ admits a left adjoint.
    \item The following three conditions are satisfied.
    \begin{enumerate}[$(a')$]
        \item The functor $q_{\mathrm{Mod}}$ is a coCartesian fibration, which, in view of the fact that $q_{\mathrm{Mod}}$ is a Cartesian fibration, amounts to the assertion that the global sections functor $\Gamma$ on modules admits a left adjoint, that we denote $\mathrm{M\spec}$.
        \item The functor $\Phi$ carries $q_{\mathrm{Mod}}$-coCartesian edges to $q$-coCartesian edges.
        \item In the induced diagram
        \[\begin{tikzcd}
\mathrm{Mod}\ar[d,"p_{\mathrm{Mod}}"]\ar[r,"\mathrm{M}\spec"] &[3em]    \mathsf{Top}^{\loc}(\mathrm{Mod})^{op}\ar[d,"p_{\mathsf{ModTop}}"] \\
    \cinfty\mathsf{ring}\ar[r,"\spec"] &[3em] \mathsf{Top}^{\loc}(\cinfty\mathsf{ring})^{op},
    \end{tikzcd}
    \]
    which commutes up to natural isomorphism by $(b')$, the functor $\mathrm{M}\spec$ carries $p_{\mathrm{Mod}}$-coCartesian morphisms to $p_{\mathsf{ModTop}}$-coCartesian morphisms.
    \end{enumerate}
\end{enumerate}
If $(ii)$ holds, we have the functor $\mathrm{M\spec}$ and conditions $(b')$ and $(c')$ imply $(1)$ and $(3)$. It suffices to verify $(i)$. There are three types of maps in $\mathcal{M}$ for which we need to check that the Cartesian pullback functor associated to $\Phi$ admits a left adjoint. For maps of the form $A\rightarrow B$ and $(X,\Of_X)\rightarrow (Y,\Of_Y)$, the pullbacks are the functors $\mathrm{Mod}_B\rightarrow\mathrm{Mod}_A$ and $\mathrm{Mod}_{\Of_X}\rightarrow\mathrm{Mod}_{\Of_Y}$, which clearly have left adjoints. We are left to deal with the case of maps of the form $A\rightarrow (X,\Of_X)$. Using that $q$ is a Cartesian and coCartesian fibration, such a map factorizes as 
\[  A\longrightarrow \Gamma(\Of_{\specr\,A})\longrightarrow (\specr\,A,\Of_{\specr\,A})\longrightarrow (X,\Of_X).\]
The last map is one we have already dealt with, so it suffices to give a left adjoint to the coCartesian pullback for the composition of the first two maps. Unwinding the definitions, we need to show that the functor 
\[ \mathrm{Mod}_{\Of_{\specr\,A}} \overset{\Gamma}{\longrightarrow} \mathrm{Mod}_{\Gamma(\Of_{\specr\,A})} \longrightarrow \mathrm{Mod}_A \]
admits a left adjoint. This functor is equal to the composition
\[ \mathrm{Mod}_{\Of_{\specr\,A}} {\longrightarrow} \mathrm{Mod}_{\widetilde{\Of_{\specr\,A}}} \overset{\Gamma}{\longrightarrow} \mathrm{Mod}_A  \]
where the first functor is the inclusion of sheaves into the category of presheaves of $\widetilde{\Of_{\specr\,A}}$-modules on $\specr\,A$. This admits a left adjoint by sheafification, and we leave to the reader the easy verification that the construction 
\[ M\longmapsto (U_a\mapsto M\otimes_A\widetilde{\Of_{\specr\,A}}(U_a))\]
is a left adjoint to the second functor.
\end{proof}

\begin{cor}\label{cor:modulestalk}
Let $A$ be a $\cinfty$-ring and let $M$ be an $A$-module, then for every $p\in \specr\,A$, the composition $M\rightarrow \Gamma(\mspec_A\,M)\rightarrow p^*\mspec_A\,M$ is naturally isomorphic to the quotient map $M\rightarrow M\otimes_AA_p$.
\end{cor}
\begin{proof}
The isomorphism is provided by the up to isomorphism commutative diagram
\[
\begin{tikzcd}
\mathrm{Mod}_A\ar[d] \ar[r] &    \mathrm{Mod}_{A_p}\ar[d,"\simeq"]\\
\mathrm{Mod}_{\Of_{\specr\,A}} \ar[r] & \mathrm{Mod}_{\Of_{\speco\,A_p}},
\end{tikzcd}
\]
where the vertical maps are the module spectra functors.
\end{proof}
As a consequence, we have the following analogue of Proposition \ref{prop:godementsheaf}.
\begin{prop}[Godement resolution]\label{prop:godementsheafmod}
Let $M$ be a module of a $\cinfty$-ring $A$. Consider the sheaf $\prod_{p\in \specr\,A}p_*(M\otimes_AA_p)$ which admits a canonical map from $\widetilde{\Of_{\specr\,A}}$ carrying a section to its various stalks. Then the sheafification 
\[  \widetilde{\mspec_A\,M}\longrightarrow \mspec_A\,M\longrightarrow \prod_{p\in \specr\,A}p_*(M\otimes_AA_p)\]
identifies the sheaf $\mspec_A\,M$ with the smallest subsheaf of $\prod_{p\in \specr\,A}p_*(M\otimes_AA_p)$ that contains the image of the map $\widetilde{\mspec_A\,M}\rightarrow \prod_{p\in \specr\,A}p_*(M\otimes_AA_p)$. Concretely, for any open $U_a\subset \specr\,A$, the codomain of the map $\widetilde{\mspec_A\,M}(U_a)\rightarrow \mspec_A\,M(U_a)$ is identified with
\[ \{(s_p)_{p\in U_a};\,\text{for all }p\in U_a\text{ there exist }b\in A,\,t\in  M\otimes_AA[b^{-1}]\text{ such that }\,U_b\subset U_a \text{ and }{s_p'}={t_p'}\text{ for all }p'\in U_b\}\subset \prod_{p\in\specr\,A}M\otimes_AA_p\]    
\end{prop}
Spectra of projective modules over finitely presented $\cinfty$-rings are particularly easy to describe.
\begin{prop}\label{prop:modulealreadysheaf}
 Let $A$ be a finitely presented $\cinfty$-ring and let $M$ be a finitely generated projective $A$-module, then the presheaf 
 \[ U_a\longmapsto A[a^{-1}]\otimes_AM\]
 is already a sheaf.
\end{prop}
\begin{proof}
Since the category of sheaves is stable under retracts, we may suppose that $M$ is a finitely generated and free $A$-module. In this case, the presheaf
\[ U_a\longmapsto A[a^{-1}]\otimes_AM\]
is a finite product of the presheaf 
\[ U_a\longmapsto A[a^{-1}]\]
which is a sheaf by Corollary \ref{cor:presheafalreadysheaf}.
\end{proof}
\begin{prop}[Joyce \cite{Joy2}]\label{prop:moduladjunction}
Let $(\specr\,A,\Of_{\specr\,A})$ be an affine $\cinfty$-scheme and let $\F$ be an $\Of_{\specr\,A}$-module. Then the counit map $ \mathrm{MSpec}_A \Gamma(\F)\rightarrow \F$ is an isomorphism.  
\end{prop}
\begin{proof}
From Theorem \ref{thm:modulespectrum}, it is straightforward to see that the counit at $\F$ is induced by the map of presheaves $\widetilde{\mspec_A\,\Gamma(\F)}\rightarrow\F$ that carries $U_a$ to the map $\F(\specr\,A)\otimes_{A}A[a^{-1}]\rightarrow\F(U_a)$. We show that this map of presheaves induces an equivalence on all stalks. Let $p\in \specr\,A$ and $t \in \F_p$ a stalk of section. Choose some $U_a$ and an $s\in \F(U_a)$ such that $s_p=t$. Choosing a pushout diagram 
\[ 
    \begin{tikzcd}
    \cinfty(\R^n)/I \ar[d] \ar[r,"\psi",hook] & A\ar[d] \\
    \cinfty(U)/I \ar[r,"\phi"] & A[a^{-1}],
    \end{tikzcd}
\]
we can find an element $b\in A$, an open set $V\subset \overline{V}\subset U_a$ such that $b_p=1$ and $b_x=0$ for all $x\notin V$. Let $b|_{U_a}$ denote the image of $b$ under the localization map $A\rightarrow A[a^{-1}]$, then we have $(b|_{U_a}s)_p=b_ps_p=t$, and since $b|_{U_a}s\in\F(U)$ vanishes away from $V$, we can use the assumption that $\F$ is sheaf to find a global section $\tilde{s}$ that restricts to $b|_{U_a}s$ on $U_a$. We conclude that the counit map is surjective on stalks. To see it is injective, we assume that a global section $s$ vanishes in some neighbourhood $U_a$ of $p$, then we can once again choose a smaller neighbourhood $V\subset U_a$ and element $\chi_V\in A$ such that $(\chi_V)_p=1$ and $(\chi_V)_x=0$ for all $x\notin V$. Then $(\chi_Vs)_p=s_p$ but using again that $\F$ is a sheaf, $\chi_vs$ vanishes globally.  
\end{proof}
It follows from the preceding proposition that the global sections functor $\Gamma:\mathrm{Mod}_{\Of_{\specr\,A}}\rightarrow \mathrm{Mod}_A$ is a fully faithful right adjoint to a localization. 
\begin{defn}\label{defn:geometricmodules}
Let $A$ be a $\cinfty$-ring. An $A$-module $M$ is \emph{geometric} if $M$ lies in the essential image of the global sections functor $\Gamma:\mathrm{Mod}_{\Of_{\specr\,A}}\rightarrow \mathrm{Mod}_A$. We let $\mathrm{Mod}^{\gmt}_A\subset\mathrm{Mod}_A$ be the full subcategory spanned by geometric $A$-modules. Let $f:M\rightarrow N$ be a map of $A$-modules with $N$ geometric, then we say that \emph{$f$ exhibits $N$ as a geometrization of $M$} if for every geometric $A$-module $K$, composition with $f$ induces a bijection
\[ \Hom_{\mathrm{Mod}_A}(K,M)\overset{\cong}{\longrightarrow} \Hom_{\mathrm{Mod}_A}(K,N), \]
equivalently, if the map $\mspec_A\,M\rightarrow \mspec_A\,N$ is an isomorphism of sheaves of $\Of_{\specr\,A}$-modules.
\end{defn}
\begin{rmk}
What we call geometric modules, Joyce calls \emph{complete modules} \cite{Joy2}.
\end{rmk}
\begin{rmk}
The category of geometric modules is presentable, since it is equivalent to the category of modules of a commutative algebra object in a topos.    
\end{rmk}
Here is a criterion for detecting geometrizations.
\begin{lem}\label{lem:geometrizationcriterion}
Let $A$ be a $\cinfty$-ring and let $f:M\rightarrow N$ be a map of $A$-modules with $N$ geometric. Then $f$ exhibits $N$ as a geometrization of $M$ if and only if for each $p\in \specr\,A$, the associated map $M\otimes_AA_p\rightarrow N\otimes_AA_p$ is an isomorphism of $A_p$-modules.     
\end{lem}
\begin{proof}
The map $f$ exhibits a geometrization if and only if for each $p\in\specr\,A$, the functor $p^*:\mathrm{Mod}_{\Of_{\specr\,A}}\rightarrow \mathrm{Mod}_{\Of_{\speco\,A_p}}$ caries the map $\mspec_A\,f$ to an equivalence. It follows from Corollary \ref{cor:modulestalk} that $p^*(\mspec_A\,f)$ can be identified with the map $M\otimes_AA_p\rightarrow N\otimes_AA_p$. 
\end{proof}
We now have two notions of `locally determined' for ideals of finitely generated $\cinfty$-rings: being geometric as an $A$-module and being germ determined. We show that these notions coincide. We will need the following easy lemma, that will be crucial when we discuss localizations of derived $\cinfty$-rings.
\begin{lem}\label{lem:openfunctionisflat}
Let $U\hookrightarrow M$ be an open embedding of manifolds, then the induced map $C^{\infty}(M)\rightarrow C^{\infty}(U)$ is a flat map of commutative $\R$-algebras.
\end{lem}
\begin{proof}
Take a finite linear combination of 0 as $\sum_{i=1}^n g_if_i|_U=0$ with $g_i\in C^{\infty}(U)$ and $f_i\in C^{\infty}(M)$, then we should show that there exists a finite set of elements $\{h_j\}_j\subset C^{\infty}(U)$ and linear combinations $g_i=\sum_j h_jb_{ij}|_U$ with $b_{ij}\in C^{\infty}(M)$ such that $\sum_{i}f_ib_{ij}=0$ for all $j$. We can write each $g_i$ as a quotient $g'_i/\gamma_i$, where $g'_i$ and $\gamma_i$ are defined on $M$ such that $\gamma_i$ does not vanish on $U$. Now pick a characteristic function $\chi_U$ for $U$ and set $h_i=1/(\gamma_i\chi_U)$, $b_{ij}=0$ if $i\neq j$ and $b_{ii}=\gamma_i g_i \chi_U$. 
\end{proof}
\begin{cor}
Let $A$ be a $\cinfty$-ring and let $a\in A$ an element, then the map $A\rightarrow A[a^{-1}]$ is a flat map of commutative $\R$-algebras.    
\end{cor}
\begin{proof}
Writing $A$ as a filtered colimit of its finitely generated subrings that contain $a$, we may suppose that $A=\cinfty(\R^n)/I$ is finitely generated, using that the class of flat maps is stable under filtered colimits. The localization at $a$ corresponds to a map $\cinfty(\R^n)/I\rightarrow \cinfty(U)/I|_U$, which is a pushout of $\R$-algebras of the map $\cinfty(\R^n)\rightarrow \cinfty(U)$, which is flat.    
\end{proof}
\begin{rmk}
Note that this lemma uses the existence of $\cinfty$ bump functions. Its counterpart in $\R$- or $\C$-analytic geometry is false: an open immersion of analytic spaces need not induce a flat map on the corresponding analytic algebras.
\end{rmk}

\begin{prop}
Let $A=\cinfty(\R^n)/I$ be a finitely generated and germ determined $\cinfty$-ring and $J\subset \cinfty(\R^n)/I$ an ideal, then $J$ is germ determined if and only $J$ is a geometric $A$-module.   
\end{prop}
\begin{proof}
We have a commuting diagram 
\[
\begin{tikzcd}
J\ar[d,hook] \ar[r] & \prod_{p\in \specr\,A} J\otimes_AA_p  \ar[d]\\
A\ar[r] & \prod_{p\in \specr\,A} A_p
\end{tikzcd}
\]
Since $A$ is germ determined, the lower horizontal map is an injection by the proof of Corollary \ref{cor:presheafalreadysheaf}, so the upper horizontal map is injective as well (in fact, the right vertical map is also injective as $A_p$ is flat over $A$). To see the surjectivity of the map $J\rightarrow \Gamma(\mathrm{MSpec}_A\,J)$, we note that Lemma \ref{lem:openfunctionisflat} guarantees that the map $J\otimes_{A}A[a^{-1}]\rightarrow A[a^{-1}]$ is injective and identifies the $J\otimes_{A}A[a^{-1}]$ with the ideal generated by the image of $J$ under the localization $A\rightarrow A[a^{-1}]$. We then see using the Godement resolution that the surjectivity of the map $J\rightarrow \Gamma(\mathrm{MSpec}_A\,J)$ is simply a reformulation of the definition of being germ determined.
\end{proof}
\begin{rmk}
We will prove in the next section that the geometrization functor for $A$-modules is an equivalence if and only if the underlying topological space $\specr\,A$ is compact. If this is the case and $A$ is finitely generated, then in particular every ideal of $A$ is germ determined.    
\end{rmk}
Let $f:A\rightarrow B$ be a map of $\cinfty$-rings. The pushforward $\mathrm{Mod}_A\rightarrow\mathrm{Mod}_B$ by tensoring $\_\otimes_AB$ generally does not preserve geometric modules, while its right adjoint always does. We now give a useful result asserting that in certain cases, this right adjoint \emph{commutes} with geometrization.
\begin{prop}\label{prop:surjectiveleftadj}
Let $f:A\rightarrow B$ be a regular epimorphism of $\cinfty$-rings. Then the diagram of right adjoints
\[
\begin{tikzcd}
 \mathrm{Mod}_A & \mathrm{Mod}_{B}\ar[l,"f^*"'] \\
  \mathrm{Mod}_{\Of_{\specr\,A}} \ar[u,"\Gamma"] &  \mathrm{Mod}_{\Of_{\speco\,B}} \ar[u,"\Gamma"] \ar[l]
\end{tikzcd}
\]
is vertically left adjointable.
\end{prop}
\begin{proof}
Since the vertical maps are fully faithful right adjoints of localizations, it suffices to show that if $M\rightarrow N$ exhibits a geometrization of $B$-modules, then it exhibits a geometrization of $A$-modules. Since $f^*$ preserves geometric modules, it follows from Lemma \ref{lem:geometrizationcriterion} that it suffices to verify the following two assertions.
\begin{enumerate}[$(a)$]
    \item For each $p\in \specr\,B$, if the map $M\otimes_BB_p\rightarrow N\otimes_BB_p$ is an isomorphism, then the map $M\otimes_AA_{p'}\rightarrow N\otimes_AA_{p'}$ where $p'$ is the map $p':A\rightarrow B\overset{p}{\rightarrow}\R$, is an isomorphism.
    \item For each $p'\in\specr\,A$ which does not factor through $\specr\,B$, the map $M\otimes_AA_{p'}\rightarrow N\otimes_AA_{p'}$ is the zero map.
\end{enumerate}
Let $p:B\rightarrow\R$ a point, with $p':A\rightarrow \R$ the composition with $f$, then we claim that the diagram 
\[
\begin{tikzcd}
A\ar[d]\ar[r] & B\ar[d] \\
A_{p'} \ar[r]  & B_p
\end{tikzcd}
\]
is a pushout of commutative $\R$-algebras. Since $A\rightarrow A_{p'}$ (and $A\rightarrow B$) is a  surjection, the comparison map $A_{p'}\otimes_{A}B\rightarrow (A_{p'}\oinfty_{A}B)^{\rmalg}$ is an isomorphism, by Corollary \ref{cor:algeffepi}. In particular, the map $A_{p'}\rightarrow A_{p'}\oinfty_{A}B$ is surjective, so by elementary algebra $A_{p'}\oinfty_{A}B$ is a local ring and by construction, the map $p:B\rightarrow \R$ factors through $A_{p'}\oinfty_{A}B$, so $A_{p'}\oinfty_{A}B$ is a local $\cinfty$-ring. The map $B\rightarrow A_{p'}\oinfty_{A}B$ is a pushout of the ind-admissible map $A\rightarrow A_{p'}$ and thus ind-admissible, so we conclude that $B\rightarrow A_{p'}\oinfty_{A}B$ is the localization at $p$. Since the diagram above is a pushout of $\R$-algebras, the diagram
\[
\begin{tikzcd}
\mathrm{Mod}_A\ar[d,"\_\otimes_AA_{p'}"']\ar[r] & \mathrm{Mod}_B\ar[d,"\_\otimes_BB_p"] \\
\mathrm{Mod}_{A_{p'}} \ar[r]  & \mathrm{Mod}_{B_p}
\end{tikzcd}
\]
is horizontally right adjointable, which implies $(a)$. For $(b)$, we consider a point $p':A\rightarrow\R$ that does not factor through $f$. If the pushout $A_{p'}\oinfty_AB$ is not the zero $\cinfty$-ring, it must be a local ring and $A_{p'}\rightarrow B_p$ a surjective map of local rings, as we argued above. Taking the quotient by the maximal ideals, we conclude that $B_p/\mathfrak{m}$ is an $\R$-algebra that admits a surjective map from $\R$ and therefore must be $\R$ itself, which contradicts the assumption that $p'$ did not factor through $f$. Using the right adjointable diagram above with $0$ in place of $B_p$, we deduce $(b)$.
\end{proof}
\begin{cor}\label{cor:epigeometric}
Let $f:A\rightarrow B$ be a regular epimorphism of $\cinfty$-rings, then a $B$-module $M$ is a geometric if and only if $M$ is geometric as an $A$-module.    
\end{cor}
\begin{proof}
The functor $\mathrm{Mod}_B\rightarrow\mathrm{Mod}_A$ is conservative and carries geometrizations to geometrizations by Proposition \ref{prop:surjectiveleftadj}.    
\end{proof}
\begin{cor}
Let $A$ be a geometric $\cinfty$-ring and let $I$ be an ideal of $A$, then the following are equivalent.
\begin{enumerate}[$(1)$]
    \item $I$ is a geometric $A$-module.
    \item The $\cinfty$-ring $A/I$ is a geometric $A$-module.
    \item The $\cinfty$-ring $A/I$ is geometric.
\end{enumerate}
\end{cor}
\begin{proof}
We have a fibre sequence
\[ 0\longrightarrow I\longrightarrow A\longrightarrow A/I\longrightarrow 0  \]
of $A$-modules, so we deduce that $I$ is geometric if and only if $A/I$ is geometric, since the category of geometric modules is stable under kernels and cokernels in $\mathrm{Mod}_A$ and $A$ itself is a geometric $A$-module. We have a regular epimorphism $A\rightarrow A/I$ so it follows from Corollary \ref{cor:epigeometric} that $A/I$ is geometric as an $A$-module if and only if $A/I$ is geometric as an $A/I$-module, which is in turn equivalent to the $\cinfty$-ring $A/I$ being complete. The map $\specr\,A/I\rightarrow \specr\,A$ exhibits a closed subspace, so $\specr\,A/I$ is Lindel\"{o}f.
\end{proof}

\begin{warn}
We say that a property $P$ on modules of $C^{\infty}$-rings is \emph{local} if the following holds for every pair $(A,M)$ where $A$ is a $C^{\infty}$-ring and $M$ is an $A$-module.
\begin{enumerate}[$(1)$]
    \item If $A\rightarrow B$ is admissible, then $M\otimes_AB$ has the property $P$.
    \item Suppose there exists an admissible covering $\{A\rightarrow A[1/a_i]\}_{i\in I}$ such that for each $i\in I$, the module $M\otimes_AA[1/a_i]\in \mathrm{Mod}_{A[1/a_i]}$ has the property $P$, then $M$ has the property $P$.
\end{enumerate}

In algebraic geometry, many natural properties, such as being a finitely presented, finitely generated, or being a finite rank vector bundle, are local (for the Zarsiki/\'{e}tale/fppf topology). In contrast, these same properties are \emph{not} local in $C^{\infty}$-geometry because the \'{e}tale topology is not finitary. For instance, let $U:=\coprod_{i=0}^{\infty}B_i\subset\R^n$ be a countable disjoint union of open balls in $B_i\subset\R^n$, then it is easy to construct $C^{\infty}(U)$-modules that are locally finitely generated, but not globally, such as a vector bundle whose rank is $i$ on each open ball $B_i$. Similarly, it easy to construct geometric $C^{\infty}$-rings which are locally finitely presented, but not globally.  
\end{warn}

\subsection{Lawvere theories}\label{sec:lawvere}
This subsection may be regarded as an elaboration on section 5.5.8 of \cite{HTT}.
\begin{defn}\label{multisortedlawvere}
A \emph{Lawvere theory} is a small $\infty$-category $\mathrm{T}$ with finite products. A \emph{transformation of Lawvere theories} is a functor $f:\mathrm{T}\rightarrow \mathrm{T'}$ that preserves finite products. We let $\mathsf{LawThy}\subset\catinf$ denote the subcategory whose objects are Lawvere theories and whose morphisms are transformations of Lawvere theories.\\
Let $\mathrm{T}$ be a Lawvere theory. A \emph{set of sorts for $\mathrm{T}$} is a (small) set $S$ together with an injective function $i:S\hookrightarrow \mathrm{Ob}_{h\mathrm{T}}$, such that every object of $\mathrm{T}$ is equivalent to a product of objects in the image of $i$. A Lawvere theory with a specified set of sorts $S$ is an \emph{$S$-sorted Lawvere theory}. If the set $S$ is the subset $\{1,\ldots,n\}\subset \N$, we call an $S$-sorted Lawvere theory an $n$-sorted Lawvere theory.
\end{defn}

\begin{defn}
Let $\xtop$ be an $\infty$-topos and let $\mathrm{T}$ be Lawvere theory. A \emph{$\mathrm{T}$-algebra in $\xtop$} is a product preserving functor $F:\mathrm{T}\rightarrow \xtop$. The full subcategory of $\fun(\mathrm{T},\xtop)$ spanned by T-algebras in $\xtop$ is denoted $\mathrm{T}\mathsf{Alg}(\xtop)$. T-algebras in the $\infty$-topos of spaces are called \emph{simplicial $\mathrm{T}$-algebras} and the $\infty$-category thereof is denoted $s\mathrm{T}\mathsf{Alg}$.
\end{defn}
Obviously, the $\infty$-category $\mathrm{T}\mathsf{Alg}(\xtop)$ has all limits and sifted colimits (which are computed objectwise in $\xtop$). \cite{HTT}, Proposition 5.5.8.10 shows that $\mathrm{T}\mathsf{Alg}(\xtop)$ a compactly generated presentable $\infty$-category; there exists an accessible localization $L:\fun(\mathrm{T},\xtop)\rightarrow \mathrm{T}\mathsf{Alg}(\xtop)$ that carries compact objects to compact objects. For any $\xtop$, we have a canonical equivalence \[\mathrm{T}\alg(\xtop)\simeq \shv_{s\mathrm{T}\alg}(\xtop) \simeq s\mathrm{T}\alg\otimes\xtop. \]
Recall that an \infcat $\icat$ is \emph{projectively generated} if it is generated under colimits by its compact projective objects, that is, objects $C\in\icat$ that have the property that the functor $\Hom_{\icat}(C,_)$ corepresented by $\icat$ preserves sifted colimits. For $\xtop=\spa$, \cite{HTT}, Lemma 5.5.8.14 shows that $s\mathrm{T}\mathsf{Alg}$ is projectively generated by the essential image of the Yoneda embedding $\mathrm{T}^{op}\hookrightarrow s\mathrm{T}\alg$. Clearly, a transformation of Lawvere theories $f:\mathrm{T}\rightarrow \mathrm{T'}$ induces for each $\infty$-topos $\xtop$ a functor $f^*:\mathrm{T'}\mathsf{Alg}(\xtop)\rightarrow \mathrm{T}\mathsf{Alg}(\xtop)$ preserving small limits and small sifted colimits. For $\xtop=\spa$, the relationship between Lawvere theories and projectively generated presentable \infcats can be made very precise. For the following proposition, we note that a left adjoint $f:s\mathrm{T}\alg\rightarrow s\mathrm{T}'\alg$ admits a right adjoint $g$ that preserves sifted colimits if and only if $f$ carries compact projective objects to compact projective objects. The only `if direction' is an immediate check and for the other direction, it suffices to show that for each $t\in \mathrm{T}$ the composition 
\[ s\mathrm{T}'\alg\longrightarrow s\mathrm{T}\alg\subset\pshv(\mathrm{T}^{op})\overset{\ev_{t}}{\longrightarrow}\spa\] 
preserves sifted colimits, but this functor is corepresented by $f(t)$ which is compact projective by assumption. 
\begin{prop}\label{prop:lawvereprojgen}
Let $\prl_{\mathrm{Proj}}\subset\prl$ be the subcategory whose objects are projectively generated presentable \infcats and whose morphisms are functors admitting a right adjoint that preserves sifted colimits. Let $\mathsf{LawThy}^{\mathrm{Idem}}\subset\mathsf{LawThy}$ denote the full subcategory spanned by idempotent complete Lawvere theories. 
\begin{enumerate}[$(1)$]
    \item The construction 
\[ \mathrm{T}  \longmapsto s\mathrm{T}\alg \]
extends to an equivalence of \infcats $\mathsf{LawThy}^{\mathrm{Idem}}\simeq \prl_{\mathrm{Proj}}$.
    \item The \infcat $\prl_{\mathrm{Proj}}$ is presentable.
    \item The \infcat $\prl_{\mathrm{Proj}}$ is semiadditive.
    \item The subcategory inclusion $\prl_{\mathrm{Proj}}\subset \prl$ preserves colimits.
\end{enumerate}
\end{prop}
\begin{proof}
We first construct the functor $\mathsf{LawThy}^{\mathrm{Idem}}\rightarrow \mathsf{Pr}^{\mathrm{L}}$. For $\mathcal{K}$ a class of small simplicial sets, let $\catinf(\mathcal{K})$ ($\catinfh(\mathcal{K})$) denote the (very) large \infcat of small (large) \infcats that admit colimits indexed by simplicial sets in $\mathcal{K}$ and functors that preserve colimits indexed by elements in $\mathcal{K}$. Let $\mathcal{P}$ be the collection of finite discrete simplicial sets together with the \infcat $\mathrm{Idem}$ and let $\mathcal{P}'$ be the collection of all small simplicial sets. It follows from \cite{HTT}, Proposition 5.3.6.2 that the inclusion $i:\catinfh(\mathcal{P}')\subset\catinfh(\mathcal{P})$ admits a left adjoint $L$. If $\icat$ is a \emph{small} \infcat that admits finite coproducts, then \cite{HTT}, Proposition 5.5.8.15 asserts that the Yoneda embedding 
\[  \icat\hooklongrightarrow \fun^{\pi}(\icat^{op},\spa) =s\icat^{op}\alg \]
exhibits a unit transformation for the adjunction $(L\adj i)$ so we conclude that the restriction of $L$ to the full subcategory $\mathsf{LawThy}^{\mathrm{Idem}}\cong\catinf(\mathcal{P})\subset \catinfh(\mathcal{P})$ takes values in the full subcategory of $\prl\subset \catinfh(\mathcal{P}')$ spanned by projectively generated presentable \infcatst. To see that a transformation of Lawvere theories is carried to a morphism in $\prl_{\mathrm{Proj}}$, we note that for a coproduct preserving functor $f:\icat\rightarrow\icatd$ between the opposite categories of two Lawvere theories, we have a commuting diagram 
\[
\begin{tikzcd}
\icat \ar[d,hook,"j_{\icat}"] \ar[r,"f"] & \icatd \ar[d,hook,"j_{\icatd}"] \\
s\icat^{op}\alg \ar[r,"L(f)"] & s\icatd^{op}\alg
\end{tikzcd}
\]
of coproduct preserving functor where moreover $L(f)$ preserves all colimits. The right adjoint to $L(f)$ preserves sifted colimits if and only if $L(f)$ carries compact projective objects to compact projective objects, so since \cite{HTT}, Proposition 5.5.8.25 implies every compact projective object in $s\icat^{op}\alg$ is a retract of one in the image of $j_{\icat}$ we conclude using the diagram above and the stability of compact projectives under retracts. This concludes the construction of the desired functor. \\    
Now $(2)$, $(3)$ and $(4)$ follow from $(1)$ and the following assertions.
\begin{enumerate}[$(2')$]
    \item[$(2')$] The \infcat $\catinf(\mathcal{P})$ is presentable.
    \item[$(3')$] The \infcat $\catinf(\mathcal{P})$ is semiadditive.
    \item[$(4')$] The functor $L|_{\catinf(\mathcal{P})}:\catinf(\mathcal{P})\rightarrow \catinfh(\mathcal{P'})$ preserves small colimits.
\end{enumerate}
Assertion $(2')$ follows from \cite{HA}, Lemma 4.8.4.2. To prove the semiadditivity of $\catinf(\mathcal{P})$, we note that the assignments $C\mapsto (C,\emptyset_{\icatd})$ and $D\mapsto (\emptyset_{\icat},D)$ where $\emptyset_{\icat}$ and $\emptyset_{\icatd}$ are initial objects $\icat$ and $\icatd$ determine fully faithful inclusions $\icat\hookrightarrow \icat\times\icatd$ and $\icatd\hookrightarrow\icat\times\icatd$ left adjoint to the projections. Let $\mathcal{A}$ be an \infcat admitting finite coproducts and suppose that we are given coproduct preserving functors $f:\icat\rightarrow\mathcal{A}$ and $g:\icatd\rightarrow\mathcal{A}$. Since the inclusion $\icat\subset\icat\times\icatd$ is a left adjoint, the \infcat $\icat_{/(s,t)}$ admits a final object $(s,(s,\emptyset)\rightarrow (s,t))$ for any $(s,t)\in \icat\times\icatd$, so the functor $f$ admits a left Kan extension $F:\icat\times\icatd\rightarrow\mathcal{A}$. Similarly, $g$ admits a left Kan extension $G:\icat\times\icatd\rightarrow\mathcal{A}$. Composing the inclusions $\icat\hookrightarrow\icat\times\icatd$ and $\icatd\hookrightarrow\icat\times\icatd$ with $F\coprod G$ yields the functors $f$ and $g$, and given any other functor $H:\icat\times\icatd\rightarrow\mathcal{A}$ compatible with $f$ and $g$, we have a natural transformation $F\coprod G\rightarrow H$. This natural transformation is an equivalence whenever $H$ preserves binary coproducts since $\icat$ and $\icatd$ generate $\icat\times\icatd$ under binary coproducts. We conclude that $\icat\times\icatd$ is a coproduct of $\icat$ and $\icatd$ in the homotopy category $h\catinf(\mathcal{P})$. For $(4')$, we only have to show that the inclusion $\catinf(\mathcal{P})\subset\catinfh(\mathcal{P})$ preserves small colimits since $L$ preserves colimits, but this is obvious.\\
We are left to prove $(1)$. For any projectively generated presentably \infcat $\icat$, the full subcategory $\icat_0$ spanned by compact projective objects is idempotent complete and admits finite coproducts, and \cite{HTT}, Proposition 5.5.8.25 asserts that $\icat\simeq s\icat_0^{op}\alg$, so the functor $L|_{\catinf(\mathcal{P})}$ is essentially surjective. Now let $\icat$ and $\icatd$ be idempotent complete \infcats admitting finite coproducts. Let $\fun'(\icat,\icatd)$ denote the full subcategory spanned by finite coproduct preserving functors, $ \fun'(s\icat^{op}\alg,s\icatd^{op}\alg)$ the full subcategory spanned by colimit preserving functors whose right adjoint preserves sifted colimits, and $\fun'(\icat,s\icatd^{op}\alg)$ the full subcategory spanned by functors preserving finite coproducts and taking values in compact projective objects. We have a commuting diagram 
\[ \begin{tikzcd} \fun'(\icat,\icatd)\ar[rr,"\theta"] \ar[dr,"\theta'"] &&  \fun'(s\icat^{op}\alg,s\icatd^{op}\alg) \ar[dl,"\theta''"]  \\ 
& \fun'(\icat,s\icatd^{op}\alg)
\end{tikzcd} \]
where the diagonal functors are induced by the Yoneda embeddings for $\icat$ and $\icatd$. It suffices to show that the diagonal functors are equivalences. \cite{HTT}, Proposition 5.5.8.15 implies that $\theta''$ is an equivalence and since \cite{HTT}, Proposition 5.5.8.25 asserts that the Yoneda embedding $\icatd\hookrightarrow s\icat^{op}\alg$ is an equivalence on the full subcategory spanned by compact projective objects by virtue of the assumption that $\icatd$ is idempotent complete, we deduce that the functor $\theta'$ is also an equivalence.
\end{proof}
We will refer to the functor $\mathrm{T}\mapsto s\mathrm{T}\alg$ as the \emph{sifted colimit completion}.
\begin{rmk}
It follows from Proposition \ref{prop:lawvereprojgen} that the sifted colimit completion carries a product $\mathrm{T}\times\mathrm{T}'$ to the coproduct $s\mathrm{T}\alg\coprod s\mathrm{T}\alg$ in $\prl$, which is the product in $(\prl)^{op}\simeq \prr$. Since $\prr\subset\catinfh$ preserves limits, the functor $\mathrm{T}\mapsto s\mathrm{T}\alg$ carries the product $\mathrm{T}\times\mathrm{T}'$ to the product $s\mathrm{T}\alg\times s\mathrm{T}'\alg$ so that $s\mathrm{T}\alg\times s\mathrm{T}'\alg$ is generated under sifted colimits by the coproduct preserving functor 
\[ \mathrm{T}^{op} \times\mathrm{T}^{\prime op} \overset{j\times j}{\longrightarrow} s\mathrm{T}\alg\times s\mathrm{T}'\alg. \]
Note that this may fail for other limits in $\mathsf{LawThy}^{\mathrm{Idem}}$; for instance, the fibre product $s\mathrm{T}'\alg\times_{s\mathrm{T}\alg}s\mathrm{T}''\alg$ in $\prl_{\mathrm{Proj}}$ is the sifted colimit completion of the Lawvere theory $\mathrm{T}'\times_{\mathrm{T}}\mathrm{T}''$ (the fibre product in $\catinf$), which need not coincide with the fibre product in $\catinfh$. \\
In general, Proposition \ref{prop:lawvereprojgen} shows that a limit of $K\rightarrow\mathsf{LawThy}^{\mathrm{Idem}}$ is obtained by taking the limit of the composition $K\rightarrow\mathsf{LawThy}^{\mathrm{Idem}}\subset \catinf$, while the colimit is obtained by taking the limit of the diagram $K^{op}\rightarrow (\prl_{\mathrm{Proj}})^{op} \subset \prr \subset \catinfh$ and extracting the compact projective objects in the resulting \infcatt.
\end{rmk}
\begin{rmk}
Let $(\prl)^{\otimes}$ be the symmetric monoidal \infcat of presentable \infcats equipped with the Lurie tensor product, then the full subcategory spanned by pairs $(\langle n\rangle,\icat_1\oplus\ldots\oplus\icat_n)$ where each $\icat_i$ is a projectively generated presentable \infcat is a symmetric monoidal subcategory. To see this, it suffices to show that the \infcat $\fun^{\mathrm{R}}(\icat^{op},\icatd)$ is projectively generated if both $\icat$ and $\icatd$ are, which follows from the equivalences
\[ \fun^{\mathrm{R}}(\icat^{op},\icatd) \simeq \fun^{\mathrm{L}}(\icatd,\icat^{op})^{op}\simeq \fun'(\icatd_0,\icat^{op})^{op}\simeq \fun^{\pi}(\icatd_0^{op},s\icat^{op}_0\alg)\simeq \fun^{\pi}(\icat_0^{op}\times\icatd_0^{op},\spa). \]
Now we show that this symmetric monoidal structure determines a symmetric monoidal structure on $\prl_{\mathrm{Proj}}$. It suffices to show that if $f:\icat\rightarrow\icatd$ preserves colimits and compact projective objects, then $f\otimes\mathrm{id}:\icat\otimes\icate\rightarrow\icatd\otimes\icate$ preserves compact projective objects as well. We can identify the functor $f\otimes\mathrm{id}$ with the functor
\[ \fun^{\pi}(\icate_0^{op},\icat)\longrightarrow \fun^{\pi}(\icate_0^{op},\icatd) \]
that composes with $f$. Since a functor $f:\icate_0^{op}\rightarrow\icat$ is compact projective if and only if $f(E)$ is compact projective for each $E\in \icate_0$ (this follows immediately from the equivalence $\fun^{\pi}(\icate_0^{op},s\icat^{op}_0\alg)\simeq \fun^{\pi}(\icat_0^{op}\times\icate_0^{op},\spa)$) and $f$ preserves compact projectives, we conclude. From the equivalence $\fun^{\mathrm{R}}(\icat^{op},\icatd) \simeq \fun^{\pi}(\icat_0^{op}\times\icatd_0^{op},\spa)$ we see that the symmetric monoidal structure is in fact the direct sum, that is, it is Cartesian and coCartesian.
 \end{rmk}
\begin{rmk}
Let $\mathsf{LawThy}'\subset \catinfh$ be the subcategory whose objects are \emph{essentially} small \infcats that admit finite products and whose morphisms are functors preserving finite products. Since every (essentially small) \infcat admits a (small) minimal model, the fully faithful inclusion $\mathsf{LawThy}\subset \mathsf{LawThy}'$ is essentially surjective. Thus, in the sequel we may work with essentially small Lawvere theories and the \infcats of algebras they generate without difficulty. 
\end{rmk}
\begin{rmk}
The sifted colimit completion may also be obtained by using the self-enrichment of $\catinfh$: the functor $\fun^{\pi}(\_,\spa)$ determines a functor $\mathsf{LawThy}^{\mathrm{Idem}}\rightarrow(\prr)^{op}$ which coincides with the sifted colimit completion after passage to adjoints. We will not prove this rigorously, but give a few hints on how to proceed. First, one can repeat the proof of Proposition \ref{prop:lawvereprojgen} (minus the semiadditivity result) to obtain an equivalence between the \infcat $\catinf^{\mathrm{Idem}}$ of small idempotent complete \infcats and the \infcat $\prl_{\mathrm{cc}}$ whose objects are presentable \infcats admitting a small set of \emph{completely compact} objects and whose morphisms are left adjoints that admit a right adjoint that admits a further left adjoint. This equivalence is implemented by the (small) colimit completion functor $L$ which carries $\icat$ to $\pshv(\icat)$. The construction $\fun((\_)^{op},\spa)$ determines another colimit preserving functor from $\catinf$ to $\prl_{\mathrm{cc}}$. Composing $\fun((\_)^{op},\spa)$ with the inverse of $L$, we obtain a colimit preserving functor $\catinf\rightarrow\catinf^{\mathrm{Idem}}$. It is not hard to see that this functor carries the full subcategory $\simp$ to itself, which implies that $L$ and $\fun((\_)^{op},\spa)$ are in fact equivalent (both are the functor $\icat\mapsto \mathrm{Idem}(\icat)$). Restricting $L$ (or $\fun((\_)^{op},\spa)$) to $\catinf(\mathcal{P})$ yields a functor $\catinf(\mathcal{P})\rightarrow\catinfh$. Let $\mathcal{Q}\rightarrow\catinf(\mathcal{P})$ be a coCartesian fibration associated to this functor, then one readily verifies that the sifted colimit completion and the functor $\fun^{\pi}((\_)^{op},\spa)$ determine the same full subcategory of $\mathcal{Q}$.  
\end{rmk}
\begin{rmk}\label{rmk:lawverethytruncation}
It is observed in \cite{HTT}. Remark 5.5.8.26 that the $n$'th truncation $\tau_{\leq n}s\mathrm{T}\mathsf{Alg}$ is precisely the full subcategory spanned by functors $\mathrm{T}\rightarrow \mathcal{S}$ taking values in $n$-truncated objects. Since we have an equivalence $\fun(\mathrm{T},\tau_{\leq0}\spa)\simeq \fun(h\mathrm{T},\set)$ and the functor $\mathrm{T}\rightarrow h\mathrm{T}$ preserves and reflects finite products, the $1$-category $\tau_{\leq 0}s\mathrm{T}\mathsf{Alg}$ can be identified with $h\mathrm{T}\mathsf{Alg}$, the nerve of the category of $h\mathrm{T}$-algebras and we have a fully faithful inclusion $h\mathrm{T}\mathsf{Alg}\hookrightarrow s\mathrm{T}\mathsf{Alg}$. In turn, this inclusion determines a morphism $sh\mathrm{T}\alg\rightarrow s\mathrm{T}\alg$ in $\prl_{\mathrm{Proj}}$. Quite often, this functor is not an equivalence, as the following example shows.
\end{rmk}
The Lawvere theories below are the basic ones we deal with in this work.
\begin{ex}
Let $\icat^{\otimes}$ be a symmetric monoidal projectively generated presentable \infcat such that the tensor product commutes with colimits separately in each variable, then the \infcat of $\calg(\icat)$ is projectively generated. To see this, we note that forgetful functor $\calg(\icat)\rightarrow\icat$ preserves limits and sifted colimits by \cite{HA}, Corollary 3.2.2.3 and 3.2.3.2 and admits a left adjoint, the free commutative algebra functor $\mathrm{Sym}^{\bullet}$. Let $\icat_0\subset\icat$ be a full subcategory spanned by a collection of compact projective generators stable under coproducts, and let $F(\icat_0)\subset \calg(\icat)$ be the essential image of $\icat_0$ under $\mathrm{Sym}^{\bullet}$, then it follows from \cite{HA}, Proposition 7.1.4.12 that the inclusion $F(\icat_0)\subset \calg(\icat)$ induces an equivalence $sF(\icat_0)^{op}\alg\simeq  \calg(\icat)$ (all this actually holds for algebras for an arbitrary \infop in the symmetric monoidal \infcat $\icat^{\otimes}$). Now suppose that 
\begin{enumerate}
    \item[$(1)$] The full subcategory $\tau_{\leq 0}\icat$ is stable under the tensor product,
\end{enumerate} 
then $\tau_{\leq 0}\icat\subset\icat$ is symmetric monoidal, and the fully faithful inclusion $\calg(\tau_{\leq 0}\icat)\subset \calg(\icat)$ can be identified with the nerve of the category $hF(\icat_0)\alg$ as the full subcategory spanned by $0$-truncated objects, since the forgetful functor $\calg(\icat)\rightarrow\icat$ preserves and detects truncations, which follows from Remark \ref{rmk:lawverethytruncation}. We have a (strictly) commuting diagram 
\[
\begin{tikzcd}
\calg(\icat)  \ar[r] & \icat  \\
\calg(\tau_{\leq 0}\icat) \ar[u,hook] \ar[r] & \tau_{\leq 0}\icat\ar[u,hook]
\end{tikzcd}
\]
of right adjoints. If we also suppose that 
\begin{enumerate}
    \item[$(2)$] For each $X\in\tau_{\leq 0}$, the $\Sigma_n$-coinvariants of $X^{\otimes n}$ are $0$-truncated,
\end{enumerate}
 then this diagram is horizontally left adjointable. If we moreover assume that 
 \begin{enumerate}
     \item[$(3)$] The \infcat $\icat_0$ lies in the full subcategory $\tau_{\leq 0}\icat\subset\icat$, 
 \end{enumerate}
then $\calg(\icat)$ is generated under sifted colimits by the image of the functor
\[  \icat_0 \overset{\mathrm{Sym}^{\bullet}}{\longrightarrow} \calg(\tau_{\leq 0}\icat)\hooklongrightarrow\calg(\icat).  \]
If $k$ is a field of containing $\Q$, then $(1)$, $(2)$ and $(3)$ hold for $\icat=\Mod_k^{\geq 0}$, the \infcat of connective $k$-modules, so $\calg^{\geq0}_k$ can be identified with $s\mathrm{T}\alg$, where $\mathrm{T}$ is the opposite of the category of discrete $k$-algebras of the form $k[x_1,\ldots,x_n]$. We may also view $\mathrm{T}$ as the \emph{discrete pregeometry} $\pregeo^{\mathrm{disc}}_k$ whose objects are affine $k$-spaces $\mathbb{A}_k^n$ and whose morphisms are polynomial maps among them. \\
When $\icat^{\otimes}=\spa^{\times}$, only $(1)$ and $(3)$ hold. $\spa$ is projectively generated by its full subcategory $\mathsf{Fin}$ of finite discrete spaces, and we can characterize its essential image under $\sym^{\bullet}$ as a certain $(2,1)$-category $\mathcal{F}\subset \calg(\spa)=\mathsf{Mon}_{\einfty}$ whose objects are parametrized by $\Z_{\geq 0}$, and whose morphisms are disjoint unions of classifying spaces of symmetric groups. For instance, we have $\Hom_{\mathcal{F}}(0,0)\simeq \coprod_n B\Sigma_n$. Then $s\mathcal{F}^{op}\alg\simeq\mathsf{Mon}_{\einfty}$, and using Remark \ref{rmk:lawverethytruncation} and the diagram above, we deduce that the homotopy category of $\mathcal{F}$ must coincide with the Lawvere theory $\mathsf{FCMon}$ of free commutative monoids. We will let $s\mathsf{CMon}$ denote the \infcat of \emph{simplicial commutative monoids}, the algebras for the 1-sorted Lawvere theory $\mathsf{FCMon}$. The transformation of Lawvere theories $\mathcal{F}\rightarrow \mathsf{FCMon}$ induces a functor $s\mathsf{CMon}\rightarrow \mathsf{Mon}_{\einfty}$. This functor is not an equivalence, but it is conservative; in fact, it is both monadic and comonadic.  
\end{ex}
\begin{ex}
The category $\cartsp$ whose set of objects is $\{\R^k;\, k \in \Z_{\geq 0}\}$ and whose morphisms are smooth maps is a Lawvere theory, generated under finite products by $\R$. A $\cartsp$-algebra in an $\infty$-topos $\xtop$ is called a \emph{$C^{\infty}$-ring in $\xtop$}.
\end{ex}
\begin{ex}
The category $\cartsp_c$ whose set of objects is $\{\R^k\times \R_{\geq 0}^m;\, k,m \in \Z_{\geq 0}\}$ and whose morphisms are \emph{interior $b$-maps} is a 2-sorted Lawvere theory, generated under finite products by $\R$ and $\R_{\geq 0}$. A $\cartsp_c$-algebra in an $\infty$-topos $\xtop$ is called a \emph{$C^{\infty}$-ring with pre-corners in $\xtop$} (the prefix `pre' will be explained in section \ref{sec:cornerlog}).
\end{ex}
The obvious functors $\mathcal{T}^{\mathrm{disc}}_{\R}\hookrightarrow \cartsp\hookrightarrow \cartsp_c$ show that every $C^{\infty}$-ring with pre-corners in $\xtop$ has an underlying $C^{\infty}$-ring, and every $C^{\infty}$-ring has an underlying commutative $\R$-algebra.\\
Anticipating the results in the next section, we develop here the theory of simplicial T-algebras to some extent. The \infcat $s\mathrm{T}\mathsf{Alg}$ is far from being an $\infty$-topos (colimits are not universal), but it does have a few redeeming features: as limits and sifted colimits are computed in the $\infty$-topos $\mathsf{PShv}(\mathrm{T}^{op})$ and geometric realizations are sifted colimits, we see that the $\infty$-category $s\mathrm{T}\mathsf{Alg}$ inherits the following properties from $\mathsf{PShv}(\mathrm{T}^{op})$.
\begin{prop}
Let $\mathrm{T}$ be a Lawvere theory.
\begin{enumerate}[$(1)$]
    \item Sifted colimits are universal in $s\mathrm{T}\mathsf{Alg}$.
    \item Every groupoid object in $s\mathrm{T}\mathsf{Alg}$ is effective.
    \item $s\mathrm{T}\mathsf{Alg}$ has an \emph{epi-mono factorization system}: there exists a factorization system $(S_L,S_R)$ on $s\mathrm{T}\mathsf{Alg}$ such that $S_L$ consists of effective epimorphisms and $S_R$ consists of monomorphisms.
    \item For each small \emph{sifted} simplicial set $K$ and each natural transformation $\overline{\alpha}:p\rightarrow q$ between functors $p,q:K^{\rhd}\rightarrow s\mathrm{T}\mathsf{Alg}$ the following holds: if $q$ is a colimit diagram and $\overline{\alpha}|_{K}$ is a Cartesian transformation, then $p$ is a colimit diagram if and only if $\alpha$ is a Cartesian transformation.  
\end{enumerate}
\end{prop}

\begin{rmk}
Let T be an $S$-sorted Lawvere theory, then we may associate to any simplicial T-algebra $A$ a collection of homotopy sets as follows: for each object $s$ of $\mathrm{T}$ in the image of $i:S\rightarrow \mathrm{Ob}_T$, there is a functor $\theta_s:s\mathrm{T}\mathsf{Alg}\rightarrow \mathcal{S}$ given by evaluating at $s$. It is customary to identify a simplicial $T$-algebra with the $S$-tuple of spaces $(\theta_s(A))_{s\in i(S)}$: we will usually denote the $C^{\infty}$-ring of smooth functions on a manifold $M$ (possibly with corners) as $C^{\infty}(M)$ and the $C^{\infty}$-ring with pre-corners as $(C^{\infty}(M),C^{\infty}_{b}(M))$. For each $n\geq 0$ and $s \in S$, we denote by $\pi_n(A)_s$ the $n$'th homotopy set of $\theta_s(A)$ which is an abelian group for $n\geq 1$. By the previous remark, the homotopy sets $\pi_0(A)_s$ can be identified with $\tau_{\leq 0}A(s)$ and if T is a 1-category, the $S$-tuple $\pi_0(A):=(\pi_0(A)_s)_{s\in i(S)}$ carries the structure of an ordinary T-algebra; we will use both notations in the sequel. 
\end{rmk}
\begin{rmk}\label{effepilawvere}
From the generating properties of the objects $s \in i(S)$ we deduce immediately that the functor 
\[\theta:s\mathrm{T}\mathsf{Alg}\overset{\prod_{s\in i(S)}\theta_s}{\longrightarrow} \prod_{s\in i(S)}\mathcal{S}  \]
is conservative. Combining this observation with the fact that each $\theta_s$ preserves geometric realizations and \cite{HTT}, Corollary 7.1.2.15, we see that a morphism $A\rightarrow B$ of simplicial T-algebras is an effective epimorphism if and only if the induced map $\pi_0(A)\rightarrow \pi_0(B)$ (with the notation from the previous remark) on sets is surjective.
\end{rmk}
\begin{rmk}\label{yonedalawvere}
By Yoneda, the space $\Hom_{s\mathrm{T}\mathsf{Alg}}(j(s),A)$ coincides with $\theta_s(A)$, where $j:\mathrm{T}^{op}\rightarrow s\mathrm{T}\mathsf{Alg}$ is the Yoneda embedding; it follows that there is a bijection of sets $\Hom_{hs\mathrm{T}\mathsf{Alg}}(j(s),A)\simeq \pi_0(A)_s$.
\end{rmk}
\begin{defn}
Let T be a Lawvere theory. A simplicial T-algebra $A$ is 
\begin{enumerate}[(1)]
    \item \emph{finitely generated to order $n$} if the functor $\Hom_{s\mathrm{T}\mathsf{Alg}}(A,\_): s\mathrm{T}\mathsf{Alg}\rightarrow \mathcal{S}$ corepresented by $A$ preserves colimits of small filtered diagrams $f:K\rightarrow s\mathrm{T}\alg$ that have the following properties.
    \begin{enumerate}[$(a)$]
    \item For each $k\in K$, the object $f(k)$ is $n$-truncated.
    \item For each edge $k\rightarrow k'$, the induced map of multisorted sets $\pi_n(f(k))\rightarrow \pi_n(f(k'))$ is a monomorphism.
    \end{enumerate}
    \item \emph{finitely generated} if $A$ is finitely generated to order $0$. 
    \item \emph{almost finitely presented} if $A$ is finitely generated to order $n$ for all $n\in \Z_{\geq 0}$.
    \item \emph{finitely presented} or \emph{compact} if the functor $\Hom_{s\mathrm{T}\mathsf{Alg}}(A,\_): s\mathrm{T}\mathsf{Alg}\rightarrow \mathcal{S}$ corepresented by $A$ preserves small filtered colimits, that is, if $A$ is a compact object.
    \item \emph{finitely presented and projective} or \emph{compact projective} if the functor $\Hom_{s\mathrm{T}\mathsf{Alg}}(A,\_): s\mathrm{T}\mathsf{Alg}\rightarrow \mathcal{S}$ corepresented by $A$ preserves small sifted colimits.
\end{enumerate}
\begin{rmk}
The following hold true.
\begin{enumerate}[$(1)$]
\item A simplicial $\mathrm{T}$-algebra $A$ is of finite generation to order $n$ if and only if $\tau_{\leq n}A$ is of finite generation to order $n$. 
\item If $A$ is of finite generation to order $n$, then $\tau_{\leq (n-1)}A$ is compact in $\tau_{\leq (n-1)}s\mathrm{T}\alg$. Conversely, if $A$ is $(n-1)$-truncated, then $A$ is compact in $\tau_{\leq (n-1)}s\mathrm{T}\alg$ if and only if $A$ is of finite generation to order $n$.
\end{enumerate}
For $\mathrm{T}'$ a 1-categorical Lawvere theory, an object $A$ in the ordinary category $\mathrm{T}'\alg$ is \emph{finitely presented (generated)} if $\Hom_{\mathrm{T}'\alg}(A,\_)$ preserves filtered diagrams (whose transition maps are monomorphisms). We see that an object $A\in s\mathrm{T}\alg$ is finitely generated to order $0$ if and only if $\pi_0(A)$ is finitely generated as $h\mathrm{T}$-algebra, and if $A$ is finitely generated to order $1$, then $\pi_0(A)$ is a finitely presented $h\mathrm{T}$-algebra. 
\end{rmk}

\end{defn}
In ordinary commutative algebra, an algebra $A$ is finitely generated if there is some free algebra $F$ on finitely many generators and a quotient map $F\rightarrow A$. The following proposition shows that the same principle can be applied to finitely generated T-algebras, with the caveat that an effective equivalence relation must be replaced by an effective groupoid.
\begin{prop}\label{fgcpt}
Let $\mathrm{T}$ be an $S$-sorted Lawvere theory, and let $A$ be a simplicial $\mathrm{T}$-algebra. The following are equivalent.
\begin{enumerate}[(1)]
     \item $A$ is finitely generated. 
     \item There exists an object $t$ of $\mathrm{T}$ and an effective epimorphism $q:j(t)\rightarrow A$, where $j:\mathrm{T}^{op}\rightarrow s\mathrm{T}\mathsf{Alg}$ is the Yoneda embedding. 
\end{enumerate}
\end{prop}
\begin{proof}
Consider the poset of subobjects $C\subset \pi_0(A)$ that admit an effective epimorphism $\pi_0(j(t))\rightarrow C$.   Since we can factor any map $j(t)\rightarrow A$ as an effective epimorphism followed by a monomorphism, this poset is filtered and its colimit is $\pi_0(A)$. Then the identity on $\pi_0(A)$ factors through some $C$ so that $C\simeq \pi_0(A)$. Since effective epimorphisms are detected on connected components, we deduce the existence of an effective epimorphism $j(t)\rightarrow A$. Now we show that $(2)\Rightarrow (1)$. Let $Y=\colim_{i\in I}Y_i$ be a colimit of a filtered diagram of $0$-truncated  objects whose transition maps are monomorphisms. A map $A\rightarrow Y$ induces a map $j(t)\rightarrow Y$ which must factor through one of the $Y_i$'s as $j(t)$ is a compact projective object in $s\mathrm{T}\mathsf{Alg}$. Because $Y_i\rightarrow Y$ is a monomorphism and the class of effective epimorphisms is left orthogonal to the class of monomorphisms, we can find the dotted arrow that makes the diagram
\begin{equation*}
\begin{tikzcd}
j(t)\ar[d]\ar[r] & Y_i\ar[d]\\
A\ar[r]\ar[ur,dotted] & Y
\end{tikzcd}    
\end{equation*}
commute, which proves that $A$ is finitely generated.
\end{proof}
\begin{rmk}
Below, we will extend the result above to objects of finite generation to order $n$ for the example $\mathrm{T}=\cartsp$, which yields presentations of simplicial $\cinfty$-rings in terms of generators and relations.
\end{rmk}
The following results give alternative characterizations of the full subcategory spanned by finitely presented T-algebras.
\begin{lem}\label{finpresenvelope}
Let $\mathrm{T}$ be a Lawvere theory. The full subcategory spanned by finitely presented $\mathrm{T}$-algebras is the smallest full subcategory of $s\mathrm{T}\mathsf{Alg}$ that contains the essential image of the embedding $\mathrm{T}^{op}\hookrightarrow s\mathrm{T}\mathsf{Alg}$ and is stable under finite colimits and retracts.
\end{lem}
\begin{proof}
Let $\icat$ be the smallest full subcategory of $s\mathrm{T}\mathsf{Alg}$ that contains the essential image of the fully faithful embedding $j:\mathrm{T}^{op}\hookrightarrow s\mathrm{T}\mathsf{Alg}$ and is stable under finite colimits and retracts. Since $s\mathrm{T}\mathsf{Alg}_{\mathrm{fp}}$ is stable under finite colimits and retracts and contains the objects of $\mathrm{T}^{op}$ as a set of compact projective generators of $s\mathrm{T}\mathsf{Alg}$, we have $\icat\subset s\mathrm{T}\mathsf{Alg}_{\mathrm{fp}}$. To establish the other inclusion, we show that every finitely presented simplicial $\mathrm{T}$-algebra is a retract of a finite colimit of objects in $j(\mathrm{T}^{op})$. Any simplicial T-algebra is a small colimit of free T-algebras, the objects in $j(\mathrm{T}^{op})$. By decomposing the index simplicial set $K$ of a small colimit into the partially ordered set of finite simplicial subsets of $K$, we may write the colimit of $K$ as a filtered colimit of finite colimits (\cite{HTT}, Corollary 4.2.3.11). Applying this to a finitely presented simplicial T-algebra $A$, we have a filtered colimit $A=\colim_{i\in \mathcal{J}}A_i$, where each $A_i$ is a finite colimit of free simplicial T-algebras. Because $A$ is finitely presented, the identity map $A\rightarrow A$ factors trough some $A_i\rightarrow A$ which shows that the desired retraction exists. 
\end{proof}
\begin{prop}\label{finlimittheory}
Let $\mathrm{T}$ be a Lawvere theory. For each idempotent complete $\infty$-category $\icat$ that admits finite limits, the restriction map 
\[\theta:\fun^{\mathrm{lex}}((s\mathrm{T}\mathsf{Alg}_{\mathrm{fp}})^{op},\icat) \longrightarrow \fun^{\pi}(\mathrm{T},\icat)   \]
induced by the fully faithful embedding $j:\mathrm{T}\rightarrow (s\mathrm{T}\mathsf{Alg}_{\mathrm{fp}})^{op}$ is an equivalence with inverse given by a functor taking right Kan extensions along $j$.
\end{prop}
\begin{proof} 
The Yoneda embedding $\icat\hookrightarrow\pshv(\icat)$ induces a commuting diagram
\[
\begin{tikzcd}
 \fun'(s\mathrm{T}\alg^{op},\icat) \ar[d,hook]\ar[r] &\fun^{\mathrm{lex}}(s\mathrm{T}\alg^{op}_{\fp},\icat) \ar[d,hook] \ar[r] & \fun^{\pi}(\mathrm{T},\icat)  \ar[d,hook] \\
 \fun'(s\mathrm{T}\alg^{op},\pshv(\icat)) \ar[r]&
\fun^{\mathrm{lex}}(s\mathrm{T}\alg^{op}_{\fp},\pshv(\icat)) \ar[r,"j^*"] & \fun^{\pi}(\mathrm{T},\pshv(\icat))
\end{tikzcd}
\]
where $\fun'(s\mathrm{T}\alg^{op},\icat)$ and $\fun'(s\mathrm{T}\alg^{op},\pshv(\icat))$ denote full subcategories of functors preserving small limits. As $\pshv(\icat)$ admits small limits and the \infcat $s\mathrm{T}\alg$ is compactly generated, left Kan extension induces an equivalence 
\[ \fun^{\omega-\mathrm{cont}}(s\mathrm{T}\alg,\pshv(\icat)^{op}) {\longrightarrow} \fun(s\mathrm{T}\alg_{\fp},\pshv(\icat)^{op}).\]
We first show that this functor restricts to an equivalence between functors $F:s\mathrm{T}\alg\rightarrow\pshv(\icat)^{op}$ preserving all colimits and right exact functors $f:s\mathrm{T}\alg_{\fp}\rightarrow\pshv(\icat)^{op}$, which in turn shows that the left lower horizontal map is an equivalence. It is clear that if $F$ preserves colimits, then the restriction $F|_{s\mathrm{T}\alg_{\fp}}$ is right exact. Now suppose that $f:s\mathrm{T}\alg_{\fp}\rightarrow \pshv(\icat)^{op}$ is right exact and let $F:s\mathrm{T}\alg\rightarrow \pshv(\icat)^{op}$ be a left Kan extension of $f$ obtained by applying the inverse of the equivalence above. Let 
\[ F':\pshv(s\mathrm{T}\alg_{\fp}) \longrightarrow \pshv(\icat)^{op} \]
be a left Kan extension of $f$ along the Yoneda embedding $j:s\mathrm{T}\alg_{\fp}\hookrightarrow\pshv(s\mathrm{T}\alg_{\fp})$. We may identify $s\mathrm{T}\alg$ with the full subcategory $\mathrm{Ind}(s\mathrm{T}\alg_{\fp})\subset \pshv(s\mathrm{T}\alg_{\fp})$ spanned by right exact functors (\cite{HTT}, Corollary 5.3.5.4), so that $F'|_{\mathrm{Ind}(s\mathrm{T}\alg_{\fp})}$ is identified with $F$. It follows from \cite{HTT}, Proposition 5.5.2.9 and remark 5.5.2.10 that $F'$ admits a right adjoint $G$. It suffices to prove that $G$ factors through $\mathrm{Ind}(s\mathrm{T}\alg_{\fp})$, then $G$ is also a right adjoint to $F'|_{\mathrm{Ind}(s\mathrm{T}\alg_{\fp})}$ which implies that $F$ preserves colimits. But the value of $G$ on some $X\in\pshv(\icat)^{op}$ is the presheaf
\[ s\mathrm{T}\alg_{\fp}^{op} \overset{f^{op}}{\longrightarrow}\pshv(\icat)\overset{\Hom(X,\_)}{\longrightarrow}\spa, \]
which is left exact as $f^{op}$ is left exact. It follows that the lower horizontal functor 
\[\fun'(s\mathrm{T}\alg^{op},\pshv(\icat))\longrightarrow  \fun^{\mathrm{lex}}(s\mathrm{T}\alg^{op}_{\fp},\pshv(\icat)) \]
is an equivalence, and the lower horizontal composition $\fun'(s\mathrm{T}\alg^{op},\pshv(\icat))\rightarrow \fun^{\pi}(\mathrm{T},\pshv(\icat))$ is an equivalence with inverse given by right Kan extension by \cite{HTT}, Proposition 5.5.8.13. To prove that the functor $\theta$ is an equivalence, it now suffices to show that for any left exact functor $f:s\mathrm{T}\alg^{op}_{\fp}\rightarrow\pshv(\icat)$ such that $f|_{j(\mathrm{T})}$ takes values in the image $\icat'$ of the Yoneda embedding $j:\icat\rightarrow\pshv(\icat)$, then $f$ also takes values in $\icat'$. This follows because Lemma \ref{finpresenvelope} shows that every $A\in s\mathrm{T}\alg_{\fp}^{op}$ is a retract of a finite limit of objects in $j(\mathrm{T})$ and $\icat'$ is stable under finite limits and retracts in $\pshv(\icat)$ by assumption. We have established that $\theta$ is an equivalence; to see that the inverse is given by the functor taking right Kan extensions, we may replace $\icat$ with $\icat'\subset\pshv(\icat)$, the essential image of the Yoneda embedding. Unraveling the proof above, the inverse of $\theta$ carries a product preserving functor $f:\mathrm{T}\rightarrow\icat'$ to the restriction to $s\mathrm{T}\alg^{op}_{\fp}$ of a right Kan extension of $\mathrm{T}\overset{f}{\rightarrow}\icat'\overset{\iota}{\subset}\pshv(\icat)$ along the embedding $j:\mathrm{T}\hookrightarrow\mathrm{T}\alg^{op}$. Since the inclusion $\mathrm{T}\hookrightarrow s\mathrm{T}\alg^{op}_{\fp}$ is fully faithful, the restriction $j_*(\iota f)|_{s\mathrm{T}\alg^{op}_{\fp}}$ is also a right Kan extension of $\iota f$. Since $j_*(\iota f)|_{s\mathrm{T}\alg^{op}_{\fp}}$ factors through $\icat'$, it is a right Kan extension of $f$. 
\end{proof}
\begin{cor}
The functor $\catinf^{\mathrm{lex},\mathrm{Idem}}\rightarrow \mathsf{LawThy}$ carrying an idempotent and finitely complete \infcat $\icat$ to $\icat$ viewed as a Lawvere theory admits a left adjoint given on objects by the assignment $\mathrm{T}\mapsto s\mathrm{T}\alg_{\fp}^{op}$. 
\end{cor}

\begin{rmk}\label{rmk:lexfunctors}
The argument in the proof above can be used to show the following: let $\icat$ be a $\kappa$-compactly generated \infcat and let $\icat_{\kappa}$ be the full subcategory spanned by $\kappa$-compact objects. Let $\icatd$ be a (not necessarily presentable) \infcat that admits all small colimits. Then a functor $F:\icat\rightarrow\icatd$ preserves colimits if and only if it admits a right adjoint, and restriction along the inclusion $\icat_{\kappa}\subset\icat$ induces an equivalence
\[ \fun^{\mathrm{L}}(\icat,\icatd) \longrightarrow\fun^{\kappa-\mathrm{rex}}(\icat_{\kappa},\icatd) \]
where $\fun^{\kappa-\mathrm{rex}}(\icat_{\kappa},\icatd)$ denotes the full subcategory spanned by functors preserving $\kappa$-small colimits. 
\end{rmk}
We turn to some stability properties of the \infcat of projectively generated \infcatst.
\begin{prop}\label{prop:lawverestability}
Let $\mathrm{T}$ be a Lawvere theory and let $s\mathrm{T}\alg$ be the associated projectively generated presentable \infcatt. 
\begin{enumerate}[$(1)$]
    \item For $K$ a small simplicial set, $\fun(K,s\mathrm{T}\alg)$ lies in $\prl_{\mathrm{Proj}}$.
    \item For $A\in s\mathrm{T}\alg$, the overcategory $s\mathrm{T}\alg_{/A}$ is equivalent to \infcat of algebras for the Lawvere theory $(\mathrm{T}^{op}_{/A})^{op}:=\mathrm{T}^{op}\times_{s\mathrm{T}\alg^{op}}(s\mathrm{T}\alg_{/A})^{op}$.
    \item For $A\in s\mathrm{T}\alg$, the undercategory $s\mathrm{T}\alg_{A/}$ is equivalent to \infcat of algebras for the Lawvere theory $\mathrm{T}_A$ defined as the full subcategory of $s\mathrm{T}\alg_{A/}$ spanned by objects of the form $A\coprod j(t)$.
\end{enumerate}
\end{prop}

\begin{rmk}
It is easy to see (using \cite{HTT}, Lemma 5.4.5.5 for instance) that $(\mathrm{T}^{op}_{/A})^{op}$ admits finite products and both the functors $(\mathrm{T}^{op}_{/A})^{op}\rightarrow \mathrm{T}$ and $(\mathrm{T}^{op}_{/A})^{op}\rightarrow s\mathrm{T}\alg^{op}$ preserve finite products. Also note that because $\mathrm{T}$ is small and $s\mathrm{T}\alg$ is locally small, the \infcat $(\mathrm{T}^{op}_{/A})^{op}$ is essentially small. The \infcat $\mathrm{T}_A$ may be described as the essential image of the finite product preserving functor 
\[\mathrm{T}\longrightarrow s\mathrm{T}\alg^{op}\overset{\_\coprod A}{\longrightarrow} (s\mathrm{T}\alg_{A/})^{op}\]
and therefore also admits finite products.  
\end{rmk}
\begin{proof}
\begin{enumerate}[$(1)$]
    \item It follows from Proposition \ref{prop:lawvereprojgen} that $(\prl_{\mathrm{Proj}})^{op}\subset \catinfh$ is stable under small products, so we deduce that for $\icat\in \prl_{\mathrm{Proj}}$ and each small simplicial set $K$, the \infcat $\fun(K,\icat)$ is projectively generated, as the evaluation functor $\fun(K,\icat)\rightarrow\prod_{k\in K}\icat$ is conservative and preserves limits and colimits.
    \item In view of \cite{HTT}, Proposition 5.5.8.22, it is sufficient to show that the fully faithful functor $\mathrm{T}^{op}_{/A}\hookrightarrow s\mathrm{T}\alg_{/A}$ takes values in compact projective objects of $s\mathrm{T}\alg_{/A}$ and that the essential image generates $s\mathrm{T}\alg_{/A}$ under sifted colimits. Let $B\rightarrow A$ a morphism in $s\mathrm{T}\alg$, then according to \cite{HTT}, Lemma 5.5.8.13, we may choose a sifted diagram $\mathcal{J}:K\rightarrow \mathrm{T}^{op}\rightarrow s\mathrm{T}\alg$ with colimit $B$, determining a colimit diagram $\overline{\mathcal{J}}:K^{\rhd}\rightarrow  s\mathrm{T}\alg$. To lift this diagram to $s\mathrm{T}\alg_{/A}$ is to lift the object $A\in s\mathrm{T}\alg_{\mathcal{J}/}$ to $s\mathrm{T}\alg_{\overline{\mathcal{J}}/}$ which is possible since the projection $s\mathrm{T}\alg_{\overline{\mathcal{J}}/}\rightarrow s\mathrm{T}\alg_{\mathcal{J}/}$ is a trivial fibration. The resulting diagram $K^{\rhd}\rightarrow\mathrm{T}^{op}_{/A}\rightarrow s\mathrm{T}\alg_{/A}$ is also a colimit diagram since the right fibration $p$ preserves and reflects colimits. It remains to be shown that each object of the form $j(t)\rightarrow A$ is compact projective. Let $\mathcal{J}:K\rightarrow s\mathrm{T}\alg_{/A}$ be a sifted diagram, then we have a fibre sequence
\[ 
\begin{tikzcd}
\Hom_{s\mathrm{T}\alg_{/A}}(j(t),\colim \mathcal{J}) \ar[d] \ar[r] & \Hom_{s\mathrm{T}\alg}(j(t),\colim \mathcal{J}) \ar[d]\\
* \ar[r] & \Hom_{s\mathrm{T}\alg}(j(t),A)
\end{tikzcd}
\]
where the fibre is taken at $j(t)\rightarrow A$, using that the functor $p$ preserves and reflects colimits. We conclude using that $j(t)$ is compact projective in $s\mathrm{T}\alg$ and the fact that colimits are universal in $\spa$.
\item Since the \infcat $s\mathrm{T}\alg_{A/}$ is presentable and the projection $p:s\mathrm{T}\alg_{A/}\rightarrow s\mathrm{T}\alg$ is conservative and admits a left adjoint carrying $B$ to $A\coprod B$, it suffices to show that the functor $p$ preserves sifted colimits, in view of \cite{HA}, Proposition 7.1.4.12. Since sifted simplicial sets are weakly contractible, this follows from \cite{HTT}, Proposition 4.4.2.9. 
\end{enumerate}
\end{proof}
\begin{rmk}
In view of the previous result, we will say that a map $f:A\rightarrow B$ of simplicial $\mathrm{T}$-algebras \emph{exhibits $B$ as finitely generated to order $n$ over $A$} (almost finitely presented over $A$, finitely presented over $A$) if $B$ is finitely generated to order $n$ in the \infcat $s\mathrm{T}_A\alg$ (almost finitely presented, finitely presented in $s\mathrm{T}_A\alg$). 
\end{rmk}
\begin{rmk}
Using that sifted colimits are universal in \infcats of algebras for Lawvere theories, it can be shown along the lines of the proof of Rezk descent for \inftopoi (\cite{HTT}, section 6.1.3) that sifted colimits in $s\mathrm{T}\alg$ are \emph{van Kampen}: the functor $s\mathrm{T}\alg^{op}\rightarrow \prl\subset\catinfh$ associated to the Cartesian fibration $\ev_0:\fun(\Delta^1,s\mathrm{T}\alg)\rightarrow s\mathrm{T}\alg$ preserves (co)sifted limits. As the \infcat $\mathrm{T}^{op}_{/A}$ is sifted for each simplicial T algebra $A$, it follows that the functor 
\[ \left(\mathrm{T}^{op}_{/A}\right)^{op} \longrightarrow \catinfh, \quad (j(t)\rightarrow A) \longmapsto s\mathrm{T}_{t/}\alg  \]
has limit $s\mathrm{T}\alg_{/A}$.
\end{rmk}

\newpage

\section{Derived $\cinfty$-geometry}
In the previous section, we introduced a \emph{pregeometry} of smooth manifolds $\diff$ and variants thereof. According to the general yoga discussed in Subsections 2.1 and 2.2, we may associate to $\diff$ a geometric envelope
\[ \diff\hooklongrightarrow \geodiffder  \]
which comes with its own theory of structured spaces. In this section, we construct this geometry explicitly by endowing the opposite \infcat of the full subcategory of compact objects of $\sring:=\fun^{\pi}(\cartsp,\spa)$ with the structure of a geometry. Subsection 3.1 is concerned with basic properties of the \infcat $\sring$ relative to the underlying algebra functor $(\_)^{\rmalg}:\sring\rightarrow\scring_{\R}$. In Subsection 3.2, we define what we believe is the correct derived version of the Archimedean spectrum and we prove that it is fully faithful on a large class of simplicial $\cinfty$-rings. The proof of this result will require us to consider modules of simplicial $\cinfty$-rings and spectra thereof; we will revisit this theory in much more detail in Part II. In Subsection 3.3, we consider a differential graded variant due to Carchedi-Roytenberg \cite{CR1}, which we show is equivalent to the simplicial theory, a result due to Nuiten \cite{Nui2}, of which we give a slightly different proof.

\subsection{Simplicial $C^{\infty}$-rings}

Using the results from the appendix and Subsection 2.4, we will show that the theory of simplicial $C^{\infty}$-rings is controlled in large part by the underlying algebraic model; in this case given by the transformation of Lawvere theories $\mathsf{Poly}_{\R}\rightarrow \cartsp$. We write $(\_)^{\rmalg}$ for the functor induced by this transformation; it takes values in $s\mathsf{Cring}_{\R}$, the $\infty$-category of simplicial commutative $\R$-algebras, and is clearly conservative. 
\begin{nota}
We reserve the symbol $\otimes^{\infty}$ for the pushout of simplicial $C^{\infty}$-rings to distinguish it from the pushout of simplicial commutative $\R$-algebras.
\end{nota}
\begin{rmk}
If we are given a diagram of $\cinfty$-rings $B\leftarrow A\rightarrow C$, we may form the pushout $B\oinfty_AC$ of $\cinfty$-rings, as in the previous section. Let $i:\cinfty\mathsf{ring}\subset\sring$ be the inclusion of $0$-truncated objects, then we may also form the pushout $i(B)\oinfty_{i(A)}i(C)$, which usually does not coincide with $i(B\oinfty_AC)$ since $i$ does not preserve colimits. We will suppress the functor $i$ from the notation and identify $\cinfty$-rings with $0$-truncated simplicial $\cinfty$-rings, so that the pushout $B\oinfty_AC$ will henceforth denote the \emph{derived} pushout, even if $A$, $B$ and $C$ are 0-truncated. The ordinary pushout of $\cinfty$-rings is instead given by $\tau_{\leq 0}(B\oinfty_AC)$.
\end{rmk}
We will make no notational distinction between simplicial commutative $\R$-algebras and connective $\mathbb{E}_{\infty}$-algebras over $\R$. Also, for $M$ a manifold, we will usually avoid writing $C^{\infty}(M)^{\rmalg}$, to avoid cluttering up notation; hopefully context will suffice to indicate whether we think of $C^{\infty}(M)$ as a simplicial $\cinfty$-ring or as an $\R$-algebra.
\begin{rmk}\label{torsionspectralsequence}
Recall that for a pushout diagram
\begin{equation*}
\begin{tikzcd}
A\ar[d]\ar[r] & B\ar[d] \\
C\ar[r] &D
\end{tikzcd}    
\end{equation*}
of simplicial commutative algebras (over any ring), there is a convergent spectral sequence
\begin{align}\label{torsionspecseq}
E_2^{p,q} = \tor^{\pi_*(A)}_p(\pi_*B,\pi_*C)_q\Rightarrow \pi_{p+q}(D),
\end{align}
see for instance, \cite{HA}, Proposition 7.2.1.19 and \cite{dagv}, Corollary 4.1.14. 
\end{rmk}
\begin{rmk}\label{rmk:hypercohomologyss}
We will also use the following version of the hypercohomology spectral sequence. Let $\xtop$ be a 1-localic \inftop and let $k$ be a commutative ring. Let $\Of_{\xtop}$ be a sheaf of simplicial commutative $k$-algebras on $\xtop$. Let $\F\in \Mod_{\Of_{\xtop}}$ be a left bounded sheaf of $\Of_X$-modules, then there is a convergent cohomological spectral sequence
\[ E^{p,q}_2=H^p(\pi_q(\F),\xtop) \Rightarrow \pi_{q-p}(\Gamma(\F))  \]
valued in the abelian 1-category $\Mod_{\Gamma(\Of_{\xtop})}^{\heartsuit}\simeq \mathrm{Mod}_{\pi_0(\Gamma(\Of_{\xtop}))}$. Here, given a sheaf $F$ of $\pi_0(\Of_{\xtop})$-modules, $H^p(F,\xtop)$ denotes the $p$'th abelian sheaf cohomology of $F$, which we can identify with the composition
\[ \mathrm{Mod}_{\pi_0(\Of_{\xtop})} \simeq \Mod_{\Of_{\xtop}}^{\heartsuit}\subset \Mod_{\Of_{\xtop}}^{\leq 0} \longrightarrow  \Mod_{\Gamma(\Of_{\xtop})}^{\leq 0}\overset{\pi_{-p}}{\longrightarrow}\Mod_{\Gamma(\Of_{\xtop})}^{\heartsuit},\]
using that $\Gamma$ is left t-exact (see the subsection on sheaves of modules below).
\end{rmk}
We start our discussion of the \infcat $\sring$ by developing several tools to aid the computation of the derived $\cinfty$ tensor product. For the next proposition, which is a general observation about sheaves on paracompact spaces with vanishing cohomology sheaves, recall that an \emph{$F_{\sigma}$-subset} of a space $X$ is a countable union of closed sets. An $F_{\sigma}$ subset of a paracompact space is paracompact, and the collection of open $F_{\sigma}$-sets of a paracompact space forms a basis for the topology that is closed under finite intersections.
\begin{prop}\label{prop:parahausdorffpresheaf}
Let $X$ be a paracompact Hausdorff space, and let $\Of_{X}$ be sheaf of simplicial commutative $k$-algebras for $k$ a commutative ring $k$ (which we can view as a connective $\mathbb{E}_{\infty}$-algebra object in $\shv_{\Mod_k}(\xtop)$ or as a $\calg_k^{\geq 0}$-valued sheaf on $X$). Suppose that $\pi_0(\Of_{\xtop})$ is a fine sheaf on $X$. Then 
\begin{enumerate}[$(1)$]
    \item For each left bounded sheaf $\F\in \Mod_{\Of_X}$ of $\Of_X$-modules, the map $\Gamma(\F)\rightarrow \Gamma(\tau_{\leq n}\F)$ induced by the unit of the truncation functor $\F\rightarrow \tau_{\leq n}\F$ exhibits $\Gamma(\tau_{\leq n}\F)$ as a $\tau_{\leq n}$-localization of $\Gamma(\F)$.
    \item Let $\mathcal{B}$ be the basis of open $F_{\sigma}$-sets of $X$, so that restriction induces an equivalence $\shv(X)\simeq \shv(\mathcal{B})$. Then for each left bounded sheaf $\F\in \Mod_{\Of_X}$, the presheaf $\widetilde{\tau_{\leq n}\F}\in \pshv_{\Mod_k}(\mathcal{B})$ given by applying the functor $\tau_{\leq n}$ objectwise is already a sheaf.
    \item For each left bounded sheaf $\F\in \Mod_{\Of_X}$, the presheaf $\widetilde{\pi_n(\F)}$ on $\mathcal{B}$ given by applying the $n$'th homotopy group functor objectwise is already a sheaf.
\end{enumerate}
\end{prop}
\begin{proof}
\begin{enumerate}[$(1)$]
    \item  Because $\F$ is left bounded, all the objects in the fibre sequence
\[ \F'\longrightarrow \F\longrightarrow \tau_{\leq n}\F \]
are left bounded. Since $\tau_{\leq n}\F$ is $n$-truncated, $\Gamma(\tau_{\leq n}\F)$ is also $n$-truncated. For any sheaf $\mathcal{G}$ of $\Of_{X}$-modules, the homotopy sheaves of $\mathcal{G}$ are sheaves of $\pi_0(\Of_{X})$-modules, which are fine sheaves on a paracompact Hausdorff space and therefore acyclic, as $\pi_0(\Of_{X})$ is fine. It follows that the higher abelian sheaf cohomology of $\pi_k(\mathcal{G})$ vanishes for all $k\in \Z$, so in case $\mathcal{G}$ is left bounded, the hypercohomology spectral sequence of Remark \ref{rmk:hypercohomologyss} collapses at the second page. Applying this to $\F'$, we see that because this sheaf of modules is $(n+1)$-connective and left bounded, $\Gamma(\F')$ is also $(n+1)$-connective. Since $\Gamma$ preserves fibre sequences, the result follows.\\

\item We have to show that the sheafification map $\widetilde{\tau_{\leq n}\F}\rightarrow \tau_{\leq n}\F$ is an equivalence, but as truncation is preserved by passing to slice topoi, the map $\widetilde{\tau_{\leq n}\F}(U)\rightarrow \tau_{\leq n}\F(U)$ is identified with the global sections of the map $\widetilde{\tau_{\leq n}\F|_U}\rightarrow \tau_{\leq n}(\F|_U)$ for each open set $U\subset X$. Then letting $U$ range over the basis $\mathcal{B}$, we see that $(1)$ applies because each $U\in\mathcal{B}$ is paracompact Hausdorff and left boundedness is preserved by passing to slice topoi, which implies that $\Gamma(\widetilde{\tau_{\leq n}\F|_U})\rightarrow \Gamma(\tau_{\leq n}(\F|_U))$ is an equivalence.
\item The homotopy groups of sheaves are given by a composition of $\Omega^n$ for some integer $n$, $\tau_{\leq 0}$ and $\tau_{\geq 0}$. The functor $\Omega^n$ is clearly defined objectwise because the functor $\shv_{\Mod_k}(\xtop)\hookrightarrow\fun(\mathcal{B}^{op},\Mod_k)$ is exact, and by $(2)$, the functor $\tau_{\leq 0}$ is defined objectwise on left bounded sheaves. We wish to show that $\tau_{\geq 0}$ is also defined objectwise on a left bounded sheaf $\F$ of $\Mod_{\Of_{X}}$. We have a morphism of fibre sequences 
\[ 
\begin{tikzcd}
\widetilde{\tau_{\geq 0}\F} \ar[r] \ar[d] &\F \ar[d,equal] \ar[r] & \widetilde{\tau_{\leq -1}} \F \ar[d] \\
\tau_{\geq 0}\F \ar[r] & \F  \ar[r] & \tau_{\leq -1} \F.
\end{tikzcd}
\]
in $\fun(\mathcal{B}^{op},\Mod_k)$. Since the right vertical map is an equivalence by $(2)$, the left vertical map is one as well.
\end{enumerate}
\end{proof}

The next lemma is a derived analogue of the fact that ideals of independent functions are point determined (see Remark \ref{idealproperties}).
\begin{lem}\label{projresolutiontransverse}
Let $M$ be an $m$-dimensional manifold and let $\{f_1,\ldots,f_n\}$, $n\leq m$, be an independent collection of functions on $M$, so that topological subspace $Z(f_1,\ldots,f_n)\subset M$ has the structure of a closed submanifold of $M$. Consider the Koszul algebra $C^{\infty}(M)[y_1,\ldots,y_n]$ with $|y_i| =1$ for $1\leq i\leq n$ and $\del y_i=f_i$ as a differential graded $\cinfty(M)$-algebra. Then the map $C^{\infty}(M)[y_1,\ldots,y_n]\rightarrow \cinfty(Z(f_1,\ldots,f_n))$ of dg algebras induced by the inclusion $Z(f_1,\ldots,f_n)\subset M$ exhibits $C^{\infty}(M)[y_1,\ldots,y_n]$ as a projective resolution of $C^{\infty}(Z(f_1,\ldots,f_n))$ in the category of dg $C^{\infty}(M)$-modules. In particular, the map $\cinfty(M)\rightarrow\cinfty(Z(f_1,\ldots,f_n))$ induces an isomorphism $\cinfty(M)/(f_1,\ldots,f_n)\cong \cinfty(Z(f_1,\ldots,f_n))$.
\end{lem}
\begin{proof}
Clearly, the Koszul complex is a complex of projective $C^{\infty}(M)$-modules, so we should show that the complex is a resolution. Let $\cinfty_{M}$ denote the sheaf of $\cinfty$ functions on $M$. Consider the sheaf of bounded differential graded $\cinfty_{M}$-modules on $M$ given by 
\[ \mathcal{F}:U\longmapsto C^{\infty}(U)[y_1,\ldots,y_n],\,\quad \del y_i=f_i|_U,\,1\leq i\leq n,\]
whose complex of global sections is the Koszul algebra $C^{\infty}(M)[y_1,\ldots,y_n]$. The closed embedding $Z(f_1,\ldots,f_n)\subset M$ induces a map of sheaves $\F\rightarrow \cinfty_{Z(f_1,\ldots,f_n)}$ where the latter sheaf carries $U\subset M$ to $\cinfty(Z(f_1|_U,\ldots,f_n|_U))$. Proposition \ref{prop:parahausdorffpresheaf} implies that the homology presheaves of $\F$ are already sheaves, so in order to show that the map $C^{\infty}(M)[y_1,\ldots,y_n]\rightarrow \cinfty(Z(f_1,\ldots,f_n))$ exhibits $C^{\infty}(M)[y_1,\ldots,y_n]$ as a projective resolution of $C^{\infty}(Z(f_1,\ldots,f_n))$, it suffices to give for each point $x\in M$ a neighbourhood basis $\{V_{\beta}\}$ of $x$ in $M$ such that $C^{\infty}(V_{\beta})[y_1,\ldots,y_n]\rightarrow \cinfty(Z(f_1|_{V_{\beta}},\ldots,f_n|_{V_{\beta}}))$ is a quasi-isomorphism. The function $(f_1,\ldots,f_n):M\rightarrow \R^n$ has full rank at $Z(f_1,\ldots,f_n)$, so it has full rank in some open neighbourhood $Z(f_1,\ldots,f_n)\subset V$. By the constant rank theorem, there is an open cover $\{U_{\alpha}\}$ of $V$ such that $U_{\alpha}\cong \R^{m}$ and in these coordinates, the function $(f_1,\ldots,f_n)$ is the projection $(x_1,\ldots,x_n):\R^m\rightarrow \R^n$ onto the first $n$ coordinates. We have a cover $\{U_{\alpha}\}\coprod \{M\setminus Z(f_1,\ldots,f_n)\}$ of $M$ so each point in $M$ has a neighbourhood basis on which $\mathcal{F}$ evaluates as either a complex of the form $C^{\infty}(V)[y_1,\ldots,y_n]$, $V\subset M\setminus Z(f_1,\ldots,f_n)$, which is acyclic because all $f_i|_{M\setminus Z(f_1,\ldots,f_n)}$ are invertible, or we have $C^{\infty}(U)[y_1,\ldots,y_n]$, where $U\subset \R^m$ is an open subset and $\del y_i=x_i$, the projection onto the $i$'th coordinate. Applying Hadamard's Lemma repeatedly, one finds that the map $C^{\infty}(U)/(x_1,\ldots,x_i)\rightarrow C^{\infty}(Z(x_1|_U,\ldots,x_i|_{U}))$ is an isomorphism for $1\leq i\leq n$. Since $Z(x_1|_U,\ldots,x_i|_{U})=U\cap \{0\}\times\R^{n-i}$ the zero locus of the function $x_{i+1}$ has measure zero, so $x_{i+1}$ is a nonzerodivisor of $C^{\infty}(U)/(x_1,\ldots,x_i)$. Thus, the sequence $(x_1,\ldots,x_n)$ is a regular sequence on $C^{\infty}(U)$ showing that the map $C^{\infty}(U)[y_1,\ldots,y_n]\rightarrow\cinfty(Z(x_1|_U,\ldots,x_n|_U))$ is a quasi-isomorphism and thus exhibits the $\cinfty(U)$-module $C^{\infty}(U)[y_1,\ldots,y_n]$ as a projective resolution of $\cinfty(Z(x_1|_U,\ldots,x_n|_U))$.
\end{proof}
In the following result, we make reference to our study of colimits in Lawvere theories exposed in Appendix A.
\begin{lem}\label{diffdiscunramified}
The transformation of Lawvere theories $\mathsf{Poly}_{\R}\rightarrow \cartsp$ satisfies the pushout axiom $(P)$ of Theorem \ref{thm:unramified2}.
\end{lem}
\begin{proof}
We will show that for graph inclusions $\R^p\rightarrow \R^{p+q}$ and $\R^k\rightarrow\R^{k+l}$ the commuting diagram
\[
    \begin{tikzcd}
    C^{\infty}(\R^{p})\otimes\cinfty(\R^q)\otimes C^{\infty}(\R^{k})\otimes\cinfty(\R^l)\ar[r]\ar[d]& C^{\infty}(\R^{p})\otimes C^{\infty}(\R^{k})\ar[d]\\
    C^{\infty}(\R^{p+q+k+l})\ar[r] & C^{\infty}(\R^{p+k})
    \end{tikzcd} 
    \]
is a pushout of simplicial commutative $\R$-algebras. By pasting of pushout squares, it suffices to argue that the diagram 
\[
    \begin{tikzcd}
    C^{\infty}(\R^{p+q})\otimes C^{\infty}(\R^{k+l})\ar[r]\ar[d]& C^{\infty}(\R^{p})\otimes C^{\infty}(\R^{k})\ar[d]\\
    C^{\infty}(\R^{p+q+k+l})\ar[r] & C^{\infty}(\R^{p+k})
    \end{tikzcd} 
    \]
is a pushout, and that for any graph inclusion $\R^n\rightarrow\R^{n+m}$, the diagram
\[
\begin{tikzcd}
\cinfty(\R^n)\otimes\cinfty(\R^m) \ar[d]\ar[r] & \cinfty(\R^{n+m}) \ar[d] \\
\cinfty(\R^n) \ar[r,equal] & \cinfty(\R^n) 
\end{tikzcd}
\]
is a pushout. We treat the upper diagram by induction on $q$ and $l$. First suppose that that $q=l=1$. The upper horizontal map is then induced by taking graphs of functions $f:\R^p\rightarrow \R$ and $g:\R^k\rightarrow \R$. As we work with discrete objects, the torsion spectral sequence collapses at the second page, so we should show that
\[\mathrm{Tor}_p^{C^{\infty}(\R^{p+1})\otimes C^{\infty}(\R^{k+1})}(C^{\infty}(\R^{p+1+k+1}),C^{\infty}(\R^{p})\otimes C^{\infty}(\R^{k})) =0,\quad p\geq 1,\]
and that 
\[\mathrm{Tor}_0^{C^{\infty}(\R^{p+1})\otimes C^{\infty}(\R^{k+1})}(C^{\infty}(\R^{p+1+k+1}),C^{\infty}(\R^{p})\otimes C^{\infty}(\R^{k})) \cong \cinfty(\R^{p+k}).\]
Denote the first n coordinates on $\R^{p+1}$
collectively by $\mathbf{x}$ and the last coordinate by $y$. The function $y-f(\mathbf{x})$ is a
submersion and its zero locus is $\mathrm{Graph}(f) \cong \R^p$, so applying Lemma \ref{projresolutiontransverse}, the ring $C^{\infty}(\R^p)$ has a resolution $C^{\infty}(\R^{p+1})[z_1]$ as a $C^{\infty}(\R^{p+1})$-module with $|z_1|=1$ and $\del z_1=y-f(\mathbf{x})$. Similarly $C^{\infty}(\R^k)$ has a resolution $C^{\infty}(\R^{k+1})[z_2]$ as a $C^{\infty}(\R^{k+1})$-module. Computing the torsion groups of the pushout using these resolution shows that the homotopy groups of the pushout are given by the homology of the complex $C^{\infty}(\R^{p+k+1+1})[z_1,z_2]$, which is by Lemma \ref{projresolutiontransverse} concentrated in degree 0 and isomorphic to $\cinfty(\R^{p+k})$ as the algebra of functions on the graph of the function $(f, g):\R^{p+k}\rightarrow \R^{1+1}$. Now suppose that for $q\leq m$, the statement is true for $l=1$ and all $p$ and $k$. Consider the diagram 
\[
    \begin{tikzcd}
    C^{\infty}(\R^{p+m+1})\otimes C^{\infty}(\R^{k+l})\ar[r]\ar[d]&  C^{\infty}(\R^{p+1})\otimes C^{\infty}(\R^{k+l})\ar[d]\ar[r]& C^{\infty}(\R^{p})\otimes C^{\infty}(\R^{k})\ar[d]\\
    C^{\infty}(\R^{p+m+1+k+l})\ar[r] & C^{\infty}(\R^{p+1+k+l})\ar[r]&  C^{\infty}(\R^{p+k}).
    \end{tikzcd} 
\]
The left square is a pushout by the inductive hypothesis applied to $p+1$, $k$ and $l=1$, and the right square is the case $q=1$ just treated; this completes the inductive step for $q$. The same argument applies to the inductive step for $l$.\\ 
We treat the other diagram. By an inductive argument similar to the one we just employed, we may assume that $m=1$ so that the graph inclusion is induced by a map $f:\R^n\rightarrow \R$. The element $1\otimes x_{n+1}-f(\mathbf{x})\otimes 1$ is a nonzerodivisor and the quotient by ideal generated by this element is easily seen to be $\cinfty(\R^n)$. Using the torsion spectral sequence, the homotopy groups of the relevant pushout are given by the homology groups of the complex $\cinfty(\R^{n+1})[\epsilon]$ with $\del\epsilon=x_{n+1}-f(\mathbf{x})$ and we conclude using Lemma \ref{projresolutiontransverse} again.
\end{proof}
\begin{cor}\label{cor:algeffepi}
$(\_)^{\rmalg}$ preserves and reflects pushouts along effective epimorphisms.
\end{cor}
\begin{proof}
Apply Theorem \ref{thm:unramified2} and Corollary \ref{cor:pushouteffepi} to the transformation $\mathsf{Poly}_{\R}\rightarrow \cartsp$ and use the fact that $(\_)^{\rmalg}$ is conservative.
\end{proof} 
The previous corollary unlocks the powerful techniques available in the setting of connective $\einfty$-algebras. We will appeal to it frequently.
\begin{rmk}\label{rmk:fermatthy}
Let $k$ be a commutative ring and let $f:\mathsf{Poly}_k\rightarrow\mathrm{T}$ be a transformation of 1-categorical Lawvere theories and let $F_{\mathrm{T}}:\mathsf{CAlg}^{0}_{k}\rightarrow \mathrm{T}\alg$ be the free $\mathrm{T}$ algebra functor left adjoint to $f^*$. We say that $f$ \emph{exhibits $\mathrm{T}$ as a Fermat theory over $k$} if the following conditions are satisfied.
\begin{enumerate}[$(1)$]
\item The functor $f$ carries $\mathbb{A}_k^1$ to a generator of $\mathrm{T}$.
\item For every $g\in F_{\mathrm{T}}(k[x,z_1,\ldots,z_n])$, there exists a unique $\delta g\in F_{\mathrm{T}}(k[x,y,z_1,\ldots,z_n])$ such that $g(x+y,\mathbf{z})=g(x,\mathbf{z})+y \delta g(x,y,\mathbf{z})$.
\end{enumerate}
It can be shown that in case $f$ exhibits $\mathrm{T}$ as a Fermat theory over $k$, then $f$ satisfies the pushout axiom $(P)$ of Theorem \ref{thm:unramified2}, so that colimits in $s\mathrm{T}\alg$ are computed by taking pushouts along the comparison maps for coproducts of free $\mathrm{T}$-algebras.
\end{rmk}
\begin{lem}\label{lem:zerolocus}
Let $M$ be a manifold and let $\{f_1,\ldots,f_n\}$ be an independent collection of $\cinfty$ functions on $M$. Let $f:M\rightarrow \R^n$ denote the induced function, then the functor $\cinfty(\_):\diff^{op}\rightarrow\sring$ carries the pullback diagram 
\[
\begin{tikzcd}
Z(f_1,\ldots,f_n)\ar[d] \ar[r] & M\ar[d,"f"] \\
* \ar[r,"0"] & \R^n
\end{tikzcd}
\]
to a pushout diagram
\[
\begin{tikzcd}
\cinfty(\R^n)\ar[d,"\ev_0"] \ar[r] & \cinfty(M)\ar[d] \\
\R \ar[r] & \cinfty(Z(f_1,\ldots,f_n))
\end{tikzcd}
\]
of simplicial $\cinfty$-rings.
\end{lem}
\begin{proof}
Since the map $\ev_0$ is an effective epimorphism, it suffices to argue that the underlying diagram of simplicial commutative algebras is a pushout, by Corollary \ref{cor:algeffepi}. According to Lemma \ref{projresolutiontransverse}, the map $\cinfty(M)\rightarrow \cinfty(Z(f_1,\ldots,f_n))$ identifies $\cinfty(Z(f_1,\ldots,f_n))$ with the quotient $\cinfty(M)/(f_1,\ldots,f_n)$, so it will suffice to show that the higher homotopy groups of the pushout $\R\otimes_{\cinfty(\R^n)}\cinfty(M)$ vanish. Since the map $\cinfty(\R^n)[\epsilon_1,\ldots,\epsilon_n]\rightarrow\R$ given by evaluation at $0$ exhibits $\cinfty(\R^n)[\epsilon_1,\ldots,\epsilon_n]$ with $|\epsilon_i|=i$ and $\del\epsilon_i=x_i$, as a projective resolution of $\R$ as a dg $\cinfty(\R^n)$-module by Lemma \ref{projresolutiontransverse}, the torsion spectral sequence shows that the higher homotopy groups are given by the higher homology of the complex $\cinfty(M)[\epsilon_1,\ldots,\epsilon_n]$ with differential given by $\del\epsilon_i=f_i$. It follows once more from Lemma \ref{projresolutiontransverse} that this homology vanishes.  
\end{proof}
\begin{lem}\label{lem:openset}
Let $U\subset M$ be an open subset of a manifold $M$ and let $\chi_U$ be a characteristic function for $U$, so that the map $U\subset M$ factors as $U\rightarrow M\times \R\rightarrow M$ where the second map is the obvious projection and the first map carries $p\in U$ to $(p,1/\chi_{U}(p))$ and determines a diffeomorphism $U\cong Z(f)\subset M\times\R$ where $f(p,y):M\times\R\rightarrow\R$ is the function $f(p,y)=\chi_U(p)y-1:M\times\R\rightarrow\R$. Then the diagram 
\[
\begin{tikzcd}
 \cinfty(\R) \ar[d,"\ev_0"]\ar[r,"f^*"] & \cinfty(M\times\R) \ar[d] \\
\R \ar[r] & \cinfty(U)
\end{tikzcd}
\]
is a pushout diagram of simplicial $\cinfty$-rings.
\end{lem}
\begin{proof}
This follows immediately from Lemma \ref{lem:zerolocus}.
\end{proof}

\begin{lem}\label{preservecoproduct}
Let $M,N$ be manifolds, then the natural map $C^{\infty}(M)\otimes^{\infty}C^{\infty}(N)\rightarrow C^{\infty}(M\times N)$ is an equivalence.
\end{lem}
\begin{proof}
First we suppose that both $M$ and $N$ are open submanifolds of Euclidean spaces $\R^m$ and $\R^n$ respectively. Choose characteristic functions $\chi_M$ and $\chi_N$ respectively, then we have a pullback diagram 
\[
\begin{tikzcd}
M\times N\ar[d] \ar[r] & \R^{m+n}\times\R^2\ar[d,"(f{,}g)"] \\
* \ar[r,"0"] & \R^2
\end{tikzcd}
\]
of manifolds, where the upper horizontal map carries $(p,q)\in M\times N$ to $(p,q,1/\chi_M(p),1/\chi_M(q))$ and $f$ and $g$ are defined by $f(p,y)=\chi_M(p)y-1$ and $g(q,z)=\chi_M(q)z-1$. It follows from Lemma \ref{lem:zerolocus} that the diagram
\[
\begin{tikzcd}
 \cinfty(\R^2) \ar[d,"\ev_0"]\ar[r,"(f{,}g)^*"] & \cinfty(\R^{m+n}\times\R^2) \ar[d] \\
\R \ar[r] & \cinfty(M\times N)
\end{tikzcd}
\]
is a pushout. We can identify this diagram with the diagram
\[
\begin{tikzcd}
 \cinfty(\R)\oinfty  \cinfty(\R)\ar[d]\ar[r] & \cinfty(\R^{m+1}) \oinfty \cinfty(\R^{n+1})\ar[d] \\
\R \ar[r] & \cinfty(M\times N)
\end{tikzcd}
\]
and we conclude using Lemma \ref{lem:openset}. For general manifolds $M,N$, we use that $\diff\simeq \mathrm{Idem}(\diff^{\mathrm{open}})$ to realize $M$ and $N$ as retracts of some $U$ respectively $V$ in $\diff^{\mathrm{open}}$. Then $M\times N$ is a retract of $U\times V$. $C^{\infty}(M)\otimes^{\infty}C^{\infty}(N)$ is a retract of $C^{\infty}(U)\otimes^{\infty}C^{\infty}(V)$ and $C^{\infty}(M\times N)$ is a retract of $C^{\infty}(U\times V)$. But as the natural map $C^{\infty}(U)\otimes^{\infty}C^{\infty}(V)\rightarrow C^{\infty}(U\times V)$ is an equivalence, $C^{\infty}(M)\otimes^{\infty}C^{\infty}(N)$ and $C^{\infty}(M\times N)$ split equivalent idempotents, so the natural map $C^{\infty}(M)\otimes^{\infty}C^{\infty}(N)\rightarrow C^{\infty}(M\times N)$ must be an equivalence. 
\end{proof}
\begin{rmk}\label{mfdscptobj}
Notice that the proof of Lemma \ref{lem:openset} shows that the essential image of the functor $C^{\infty}(\_):\diff\rightarrow sC^{\infty}\mathsf{ring}$ consists of retracts of pushouts of compact projective objects of $sC^{\infty}\mathsf{ring}$ which are thus compact.
\end{rmk}
\begin{lem}\label{lem:opentransverse}
The functor $C^{\infty}(\_):\diff\rightarrow sC^{\infty}\mathsf{ring}^{op}$ sending a manifold $M$ to the discrete simplicial $C^{\infty}$-ring of smooth functions on $M$ preserves transverse pullbacks of the form
\[
\begin{tikzcd}
N\times_U M\ar[r]\ar[d] & N\ar[d]\\    
M\ar[r,] & U 
\end{tikzcd}
\]
where $U$ is an open submanifold of $\R^n$ for some $n\geq 1$.
\end{lem}
\begin{proof}
We note that the pullback $N\times_U M$ is equivalent to the pullback 
\[
\begin{tikzcd}
(M\times N)\times_{U\times U}U\ar[d]\ar[r] & N\times M\ar[d,"g "]\\
U \ar[r] & U\times U
\end{tikzcd}
\] 
so, as the functor $C^{\infty}(\_):\diff\rightarrow sC^{\infty}\mathsf{ring}^{op}$ preserves binary products by Lemma \ref{preservecoproduct}, we only have to deal with pullback diagrams of the form above. Since $U\subset \R^n$ is open, the diagonal embedding $U\rightarrow U\times U$ is cut out by an independent collection of functions $\{f_1,\ldots,f_n\}$ determining a map $f:U\times U\rightarrow \R^n$. Consider the diagram
\[
\begin{tikzcd}
(M\times N)\times_{U\times U}U\ar[d]\ar[r] & N\times M\ar[d,"g "]\\
U \ar[r] \ar[d]& U\times U\ar[d,"f"] \\
*\ar[r] & \R^{n}
\end{tikzcd}
\] 
in which both squares are pullbacks. Because $g:N\times M\rightarrow U\times U$ is transverse to $U\rightarrow U\times U$, the functions $f_i\circ g$ form an independent collection; Lemma \ref{lem:zerolocus} now guarantees that the upper square and the outer rectangle in the diagram
\[
\begin{tikzcd}
\R\ar[r]\ar[d] & \cinfty(\R^n)\ar[d]\\
 \cinfty(U)\ar[d]\ar[r] & \cinfty(U\times U)\ar[d] \\
\cinfty(N\times_UM) \ar[r] & \cinfty(M\times N)
\end{tikzcd}
\]
are pushouts, so we conclude. 
\end{proof}
We will momentarily show that the functor $C^{\infty}(\_)$ preserves \emph{all} transverse pullbacks.
\subsubsection{Good cell objects}
If $f:C\rightarrow C'$ is a morphism in a presentable category $\icat$ and $S\subset\fun(\Delta^1,\icat)$ is a small set of morphisms in $\icat$, then an application of the \emph{small object argument} allows us to factor $f$ as a composition $C\rightarrow C''\rightarrow C'$ where the second map has the right lifting property with respect to the morphisms in $S$ and the first map is a sequential colimit of maps that are pushouts of coproducts of maps in $S$. Our goal in this subsection is to prove a similar, but more refined result for maps of simplicial $\cinfty$-rings. Roughly speaking, we will prove that any map $f:A\rightarrow B$ of simplicial $\cinfty$-rings can be obtained as a countable sequence of `cell attachments', which yields a description of simplicial $\cinfty$-rings in terms of `generators and relations'. This procedure interacts well with the hierarchy of finiteness conditions we introduced above; suppose that $f:A\rightarrow B$ is a map of ordinary algebras for a 1-sorted Lawvere theory, then $f$ exhibits $B$ as finitely generated over $A$ if there is some finite set $E$ and a factorization
\[ A\longrightarrow A\coprod j(t^E) \longrightarrow B \]
where the second map is a regular epimorphims. If $f$ is also finitely presented, then the set of relations $A\coprod j(t^E)\times_B A\coprod j(t^E)$ is finite. We will prove a version of this result for simplicial $\cinfty$-rings: a morphism $f:A\rightarrow B$ is finitely generated to order $n$ precisely if one can find a cell decomposition for which all cells up to degree $n$ are finite, in a suitable sense.
\begin{rmk}
The results in this subsection are stated only for simplicial $\cinfty$-rings, but it is easy to see that they hold for Fermat theory $\mathsf{Poly}_k\rightarrow\mathrm{T}$ over a commutative ring $k$.
\end{rmk}

\begin{nota}
Let $V$ be a real vector space, possibly of infinite dimension. We write 
\[C^{\infty}(V^{\vee}):=\underset{V'\subset V\,\mathrm{dim}V'<\infty}{\colim}C^{\infty}((V')^{\vee})\] 
for the free simplicial $C^{\infty}$-ring on $V$. Evaluation at $0\in V^{\vee}$ yields a map $C^{\infty}(V^{\vee})\rightarrow \R$ of simplicial $C^{\infty}$-rings, so $C^{\infty}(V^{\vee})$ is augmented over the initial object in $\sring$, and we may consider the $n$-fold suspension $\Sigma^nC^{\infty}(V^{\vee})$ with respect to the augmentation. 
\end{nota}
\begin{defn}\label{defn:cellobject}
Let $A\rightarrow B$ be a morphism of simplicial $C^{\infty}$-rings. We say that $B$ is a \emph{good $A$-cell object} if $B$ is a colimit of a sequential diagram
\[A=A_{-1}\overset{\phi_{-1}}{\longrightarrow} A_0\overset{\phi_0}{\longrightarrow} A_1 \longrightarrow \ldots,\]
where $\phi_{-1}$ is a pushout along a map of the form $A\rightarrow A\otimes^{\infty}C^{\infty}(V_{-1})$ for $V_{-1}$ a possibly infinite dimensional vector space, and $\phi_n$ for $n\geq0$ is a pushouts along a map of the form $ A\otimes^{\infty}\Sigma^nC^{\infty}(V_n)\rightarrow A$ for $V_n$ a possibly infinite dimensional vector space. A good $A$-cell object is
\begin{enumerate}[$(1)$]
	\item \emph{$n$-finite} if the dimension of the vector space $V_k$ is finite for all $k< n$. 
    \item \emph{almost finite} if it is $n$-finite for all $n\in \Z_{\geq 0}$.
    \item \emph{finite} if it is almost finite and the directed colimit over $\Z_{\geq -1}$ in the definition may be replaced by a finite directed subset $\{n\}_{0\leq n\leq j}$.
\end{enumerate} \end{defn}
\begin{prop}\label{prop:afpcell}
Let $f:A\rightarrow B$ be a morphism of simplicial $C^{\infty}$-rings, then the following hold. 
\begin{enumerate}[(1)]
    \item $B$ is equivalent to a good $A$-cell object
    \[A=A_{-1}\longrightarrow A_0\longrightarrow A_1\longrightarrow\ldots\longrightarrow B. \]
    If $f$ is $m$-connective, then may assume that $A=A_0=A_1=\ldots=A_m$. 
    \item $B$ is finitely generated to order $n$ over $A$ if and only if $B$ is equivalent to an $n$-finite good $A$ cell object. 
    	\item $B$ is almost finitely presented over $A$ if and only if $B$ is equivalent to an almost finite good $A$-cell object.
    \item $B$ is finitely presented over $A$ if and only if $B$ is equivalent to a retract of a finite good $A$-cell object. 
\end{enumerate}
\end{prop}
The proof of Proposition \ref{prop:afpcell} will require  a few preliminaries. The following lemmas are adapted from \cite{dagix}, Lemmas 12.18 and 12.19.
\begin{rmk}
The free $C^{\infty}$-ring functor $F^{\cinfty}$ preserves colimits, so we have $\Sigma^nC^{\infty}(V^{\vee})\simeq F^{\cinfty}(\mathrm{Sym}^{\bullet}(V[n]))$ for each $\R$-module $V$. The forgetful-free adjunction between $\mathbb{E}_{\infty}$-algebras and simplicial $C^{\infty}$-rings now establishes the equivalence
\[ \Hom_{\sring}(\Sigma^nC^{\infty}(V^{\vee}),A)\simeq \Hom_{\mathsf{Mod}_{\R}}(V[n],A^{\rmalg})  \]
for all $A\in \sring$.
\end{rmk}

\begin{lem}\label{looping1}
Let $V$ be a real vector space. The map $V[n]\rightarrow \Sigma^nC^{\infty}(V^{\vee})^{\rmalg}$ corresponding to the identity $\Sigma^nC^{\infty}(V^{\vee})\rightarrow \Sigma^nC^{\infty}(V^{\vee})$ via the equivalences above induces an equivalence $\mathrm{Sym}^{\bullet}(V[n])\rightarrow \Sigma^nC^{\infty}(V^{\vee})^{\rmalg}$ of $\mathbb{E}_{\infty}$-algebras over $\R$ for $n>0$.
\end{lem}
\begin{proof}
Since all forgetful and free functors involved commute with filtered colimits, we may write $V=\colim_{V'\subset V,\,\mathrm{dim}\,V'<\infty}$ and suppose that $V$ is finite dimensional. We work by induction on $n$. For $n=1$, we are asked to prove that the natural map \[\R\otimes_{\mathrm{Sym}^{\bullet}(V)}\R\rightarrow \R\otimes^{\infty}_{C^{\infty}(V^{\vee})}\R\simeq \R\otimes_{C^{\infty}(V^{\vee})^{\rmalg}}\R\]
is an equivalence (the last equivalence follows by Corollary \ref{cor:algeffepi}). Suppose that $V$ is 1-dimensional, then $\mathrm{Sym}^{\bullet}(V)\simeq \R[x]$ and we have a map of projective resolutions
\[
\begin{tikzcd}
0\ar[r] & \R[x] \ar[r,"x"] \ar[d] &  \R[x] \ar[r] \ar[d]& \R\ar[d,"\mathrm{id}"] \\
0\ar[r] & C^{\infty}(\R)\ar[r,"x"]  & C^{\infty}(\R) \ar[r] & \R
\end{tikzcd}
\]
where $x$ denotes multiplication by the function $x\mapsto x$ on $\R$
which shows that $\tor_i^{\R[x]}(\R,\R)\cong \tor_i^{C^{\infty}(\R)}(\R,\R)$ for all $i\geq0$, so we are done for $n=1$ and $\mathrm{dim}\,V=1$. For $V$ $k$-dimensional, the map $ \mathrm{Sym}^{\bullet}(V[1])\rightarrow \Sigma C^{\infty}(V^{\vee})^{\rmalg}$ is simply the $k$-fold tensor product of the equivalence we have just established. The induction step for $n\geq 1$ follows at once from Corollary \ref{cor:algeffepi}.
\end{proof}

\begin{lem}\label{looping2}
Let $A$ be a simplicial $C^{\infty}$-ring and let $V$ be a vector space. Let $n>0$ and $V[n]\rightarrow A^{\rmalg}$ be a map of $\R$-modules adjoint to a map $\varphi:V\otimes_{\R} A^{\rmalg}[n]\rightarrow A^{\rmalg}$ of $A^{\rmalg}$-modules. By taking the symmetric algebra and the free simplicial $C^{\infty}$-ring, $V[n]\rightarrow A^{\rmalg}$ is adjoint to a map $\Sigma^nC^{\infty}(V^{\vee})\rightarrow A$. Consider the pushout diagram
\[
\begin{tikzcd}
\Sigma^nC^{\infty}(V^{\vee})\ar[d]\ar[r]& A\ar[d]\\
\R\ar[r] & B
\end{tikzcd}
\]
Then there is a natural map $\mathrm{cofib}(\varphi)\rightarrow B^{\rmalg}$ of $A^{\rmalg}$-modules which has (2n+2)-connective cofibre. 
\end{lem}
\begin{proof}
By Corollary \ref{cor:algeffepi} and Lemma \ref{looping1}, we have $B^{\rmalg}\simeq \R\otimes_{\mathrm{Sym}^{\bullet}(V[n])}A^{\rmalg}$. The composition
\[ V\otimes_{\R} A^{\rmalg}[n] \overset{\varphi}{\longrightarrow} A^{\rmalg} \longrightarrow B \]
of morphisms of $A^{\rmalg}$-modules is homotopic to the composition 
\[ V\otimes_{\R} A^{\rmalg}[n] {\longrightarrow} \mathrm{Sym}^{\bullet}(V[n]) \longrightarrow B  \]
which is nullhomotopic by construction, yielding the desired map $\cofib(\varphi)\rightarrow B$. Since taking cofibres commutes with tensor products, we have an equivalence 
\[ \cofib(V[n]\otimes_{\R}\mathrm{Sym}^{\bullet}(V[n])\rightarrow \mathrm{Sym}^{\bullet}(V[n]))\otimes_{\mathrm{Sym}^{\bullet}(V[n])}A^{\rmalg}\simeq \cofib(V\otimes_{\R} A^{\rmalg}[n]\rightarrow A^{\rmalg})=\cofib(\varphi). \]
One readily verifies that  $\cofib(V[n]\otimes_{\R}\mathrm{Sym}^{\bullet}(V[n])\rightarrow \mathrm{Sym}^{\bullet}(V[n]))$ has vanishing homotopy groups in degrees $0<i\leq 2n$, so the map $\cofib(V[n]\otimes_{\R}\mathrm{Sym}^{\bullet}(V[n])\rightarrow \mathrm{Sym}^{\bullet}(V[n]))\rightarrow \R$ has (2n+2)-connective cofibre, showing that the map
\[\cofib(V[n]\otimes_{\R}\mathrm{Sym}^{\bullet}(V[n])\rightarrow \mathrm{Sym}^{\bullet}(V[n]))\otimes_{\mathrm{Sym}^{\bullet}(V[n])}A^{\rmalg}\simeq \cofib(\varphi) \longrightarrow  B^{\rmalg}\simeq \R\otimes_{\mathrm{Sym}^{\bullet}(V[n])}A^{\rmalg} \]
has $(2n+2)$-connective cofibre as well. 
\end{proof}

\begin{proof}[Proof of Proposition \ref{prop:afpcell}]
Let $A\rightarrow B$ be a simplicial $C^{\infty}$-ring, then we prove $(1)$. We will inductively define a sequence of $n$-connective maps $\psi_n:A_n\rightarrow B$ formed by pushouts as in Definition \ref{defn:cellobject}. If $f$ is not $0$-connective, choose an effective epimorphism $A\otimes^{\infty}C^{\infty}(\R^{J_0})\rightarrow B$; for instance, $J_0$ may be the set of those generators of $\pi_0(B)$ that are not in the image of $\pi_0(A)\rightarrow \pi_0(B)$. Suppose that $n>0$. Assuming we have constructed an $(n-1)$-connective map $\psi_{n-1}:A_{n-1}\rightarrow B$ (which may be the map $A\rightarrow B$ if this map is already $(n-1)$-connective), we construct $\psi_n$. We have $\pi_j(A_{n-1})\simeq \pi_j(B)$ for $j<(n-1)$. The algebraic fibre $\mathrm{fib}(\psi_{n-1}^{\rmalg})$ of the map $\psi_{n-1}^{\rmalg}:A_{n-1}^{\rmalg}\rightarrow B^{\rmalg}$ of connective $\mathbb{E}_{\infty}$-algebras over $\R$ fits into a long exact sequence
\[\ldots\longrightarrow \pi_{n}(A^{\rmalg}_{n-1})\longrightarrow \pi_{n}(B^{\rmalg})\longrightarrow \pi_{n-1}(\mathrm{fib}(\psi_{n-1}^{\rmalg}))\longrightarrow \pi_{n-1}(A^{\rmalg}_{n-1})\longrightarrow \pi_{n-1}(B^{\rmalg})\longrightarrow 0\longrightarrow\ldots \]
Choose a set $J_n$ and a map $\R^{J_n}\otimes_{\R}A_{n-1}^{\rmalg}[n-1]\rightarrow \mathrm{fib}(\psi_{n-1}^{\rmalg})$ of $A_{n-1}^{\rmalg}$-modules that induces a surjective map $\R^{J_n}[n-1]\otimes_{\R}\pi_0(A_{n-1}^{\rmalg})\rightarrow \pi_{n-1}(\mathrm{fib}(\psi_{n-1}^{\rmalg}))$. The composition \[\varphi:\R^{J_n}\otimes_{\R}A_{n-1}^{\rmalg}[n-1]\longrightarrow \mathrm{fib}(\psi_{n-1}^{\rmalg})\longrightarrow A_{n-1}^{\rmalg}\]
in the $\infty$-category of $A_{n-1}^{\rmalg}$-modules is adjoint to a map 
\[ \R^{J_n}[n-1]\longrightarrow   A_{n-1}^{\rmalg}\]
of $\R$-modules. This map yields a map $\mathrm{Sym}^{\bullet}(\R^{J_n}[n-1])\rightarrow A_{n-1}^{\rmalg}$ in $s\mathsf{Cring}_{\R}$, which is in turn adjoint to a map $\Sigma^{n-1} C^{\infty}((\R^{J_n})^{\vee})\rightarrow A_{n-1}$ of simplicial $C^{\infty}$-rings, with $\Sigma^{n-1} C^{\infty}((\R^{J_n})^{\vee})$ the $(n-1)$'th suspension of $C^{\infty}((\R^{J_n})^{\vee})$ at the basepoint $0\in (\R^{J_n})^{\vee}$. Now we define $A_{n}$ as the right pushout square in the diagram
\begin{equation*}
\begin{tikzcd}
\Sigma^{n-1} C^{\infty}((\R^{J_n})^{\vee})\ar[r]\ar[d] & A\otimes^{\infty}\Sigma^{n-1} C^{\infty}((\R^{J_n})^{\vee})\ar[r,"f"]\ar[d]  &A_{n-1}\ar[d]\\
\R \ar[r] & A\ar[r] & A_{n}
\end{tikzcd}    
\end{equation*}
where the left square and the outer rectangle are pushouts as well. The canonical nullhomotopy of the map $\R^{J_n}\otimes_{\R}A_{n-1}^{\rmalg}[n-1]\rightarrow \mathrm{fib}(\psi_{n-1}^{\rmalg})\rightarrow B^{\rmalg}$ yields a homotopy between $\psi_{n-1}\circ f$ and $\Sigma^{n-1} C^{\infty}((\R^{J_n})^{\vee})\rightarrow A\rightarrow B$, so we get a map $\psi_n:A_n\rightarrow B$. We check that $\psi_n$ is $n$-connective: notice that the left and middle vertical maps in the diagram above induce surjections on connected components, so by Corollary \ref{cor:algeffepi}, we have an equivalence $A_n^{\rmalg}\simeq \R\otimes_{\Sigma^{n-1}C^{\infty}((\R^{J_n})^{\vee})^{\rmalg}}A_{n-1}^{\rmalg}$. For $n=1$, we observe that $\pi_0(A_1)\simeq \pi_0(C^{\infty}((\R^{J_0})^{\vee})/\pi_0(\fib(\psi_0^{\rmalg})))\simeq \pi_0(B)$. For $n>1$, Lemma \ref{looping2} provides us with a map $\cofib(\varphi)\rightarrow A_n$ with $(2n)$-connective cofibre. Comparing the $\pi_{n-1}$-terms in the long exact sequence associated with the fibre sequence $\mathrm{fib}(\psi^{\rmalg}_{n})\rightarrow A_{n}^{\rmalg}\rightarrow B$ with those of the long exact sequence associated to the cofibre sequence of $\varphi$ yields the desired connectivity estimate.\\
We now prove $(2)$ and $(3)$. We leave to the reader the straightforward verification that colimits of $n$-finite good cell objects over $A$ are finitely generated to order $n$ over $A$. Note that an easy inductive argument reduces the reverse implication to the following one.
\begin{enumerate}
\item[$(*)$] Let $A\rightarrow B$ be an $n$-finite good $A$-cell object $\{A_i\}_{i\in \Z_{\geq 0}}$, and suppose that $B$ is finitely generated to order $n+1$ over $A$, then $A\rightarrow B$ is equivalent to an $(n+1)$-finite good $A$-cell object $\{A'_i\}_{i\in \Z_{\geq 0}}$ that coincides with $\{A_i\}_{i\in \Z_{\geq 0}}$ up to the $n$'th cell. 
\end{enumerate}
Write the cell decomposition as
\[ A=A_{-1}\longrightarrow A_0\longrightarrow A_1\longrightarrow \ldots\longrightarrow B.\]
Note that we have $\tau_{\leq n}B\simeq \tau_{\leq n}A_k$ for all $k\geq n+1$. We may write $A_{n+1}$ as a colimit 
\[ A_{n+1} = \R\otimes^{\infty}_{\Sigma^{n}C^{\infty}((\R^{J_{n+1}})^{\vee})} A_{n}\simeq \underset{S\subset J_{n+1},\,|S|<\infty}{\colim}\R\otimes^{\infty}_{\Sigma^{n}C^{\infty}(\R^S)} A_{n}. \]
Since the full subcategory of $\sring$ spanned by $n$-truncated objects is stable under filtered colimits, we have an equivalence $\tau_{\leq n}A_{n+1}\simeq \underset{S\subset J_{n+1},\,|S|<\infty}{\colim}\tau_{\leq n}(\R\otimes^{\infty}_{\Sigma^{n}C^{\infty}(\R^S)} A_{n})$. Denote $C_S=\tau_{\leq n}(\R\otimes^{\infty}_{\Sigma^{n}C^{\infty}(\R^S)} A_{n})$, so that $\tau_{\leq n}A_{n+1}$ is a colimit of the filtered diagram $\{C_S\}_{S\subset J_{n+1},|S|<\infty}$. Since $B$ is finitely generated to order $(n+1)$ over $A$, the map $B\rightarrow \tau_{\leq n}B\simeq \tau_{\leq n}A_{n+1}$ factors through some $C_S$ so that $\tau_{\leq n}B$ is a retract of this $C_S$. We now show that there is some $S\subset S'$ for which $\tau_{\leq n}A_{n+1}\simeq C_{S'}$. Since $\R\otimes^{\infty}_{\Sigma^{n}C^{\infty}(\R^S)} A_{n}$ is finitely generated to order $n+1$ over $A$, its $n$-truncation $C_S$ is compact in $\tau_{\leq n}\sring_{/A}$, so that the composition $r:C_S\rightarrow \tau_{\leq n}A_{n+1}\rightarrow  C_S$ is a compact object in the \infcat $\fun(\Delta^1,\tau_{\leq n}\sring_{A/})$, by \cite{HTT}, Proposition 5.3.4.13. We have a commuting square
\[
\begin{tikzcd}
C_S\ar[d]\ar[r,"r"] & C_S\ar[d] \\
\tau_{\leq n}A_{n+1}\ar[r,"\mathrm{id}_A"] & \tau_{\leq n}A_{n+1}
\end{tikzcd}
\]  
Since $\mathrm{id}_A:A\rightarrow A$ is a filtered colimit  of the diagram of arrows $\{C_{S'}\overset{\mathrm{id}}{\rightarrow} C_{S'}\}_{S\subset S'\subset J_{n+1},|S'|<\infty}$, we deduce the existence of a factorization
\[
\begin{tikzcd}
C_S\ar[d]\ar[r,"r"] & C_S\ar[d] \\
C_{S'} \ar[r,"\mathrm{id}"] \ar[d]& C_{S'}\ar[d]\\
\tau_{\leq n}A_{n+1}\ar[r,"\mathrm{id}"] & \tau_{\leq n}A_{n+1}
\end{tikzcd}
\]
for some $S\subset S'$. The map $C_{S'}\rightarrow \tau_{\leq n}A_{n+1}$ is an equivalence on all homotopy groups in degrees $<n$. Because $\tau_{\leq n}A_{n+1}$ is a retract of $C_{S'}$, the map $\pi_n(A_{n+1})\rightarrow\pi_n(C_{S'})$ is a monomorphism. From the factorization above, we deduce that the composition
\[ C_S \longrightarrow C_{S'}\longrightarrow \tau_{\leq n}A_{n+1} \longrightarrow C_S\longrightarrow C_{S'}  \]
is equivalent to the map $C_S\rightarrow C_{S'}$ induced by the map $S\rightarrow S'$ of subsets of $J_{n+1}$. Since this map is an inclusion, the map on homotopy groups $\pi_n(C_S)\rightarrow \pi_n(C_{S'})$ is an epimorphism, so $\pi_n(A_{n+1})\rightarrow \pi_n(C_{S'})$ is an epimorphism as well. We conclude that the map $C_{S'}\rightarrow \tau_{\leq n}A_{n+1}$ is an equivalence on all homotopy groups in degrees $\leq n$. Define 
\[ A_{n+1}':= \R\otimes^{\infty}_{\Sigma^{n}C^{\infty}(\R^{S'})} A_{n}, \]
then $A_{n+1}'$ is finitely presented over $A$ and the map $A'_{n+1}\rightarrow B$ induces an equivalence $\tau_{\leq n}A'_{n+1}\simeq \tau_{\leq n}B$ (but need not be $(n+1)$-connective. Choose a real vector space $V$ and a factorization 
\[A_{n+1}'\longrightarrow A_{n+1}'\oinfty\Sigma^{n+1}\cinfty(V)\longrightarrow B\]
 where the second map is $(n+1)$-connective; for instance, we can choose a surjective map $V\rightarrow \pi_{n+1}(B)$ of $\R$-vector spaces. Since $\sring$ is the \infcat of algebras for a Lawvere theory, it comes endowed with an ($(n+1)$-connective, $n$-truncated) factorization system which yields a left adjoint $L:\sring_{/B}\rightarrow\tau_{\leq n}\sring_{/B}$. Considering the finite subspaces of $V$, we have a filtered diagram $\{A_{n+1}'\oinfty\Sigma^{n+1}\cinfty(V')\}_{V'\subset V,\mathrm{dim}(V)<\infty}$ with colimit $A_{n+1}'\oinfty\Sigma^{n+1}\cinfty(V)$. Since the canonical map $L(A_{n+1}'\oinfty\Sigma^{n+1}\cinfty(V))\rightarrow B$ is $n$-truncated and $(n+1)$-connective, it is an equivalence. It follows that the filtered diagram 
 \[\{\tau_{\leq (n+1)}L(A_{n+1}'\oinfty\Sigma^{n+1}\cinfty(V'))\}_{V'\subset V,\mathrm{dim}(V)<\infty}\]
  has the following properties:
\begin{enumerate}[$(i)$]
\item The $(n+1)$'th truncation $\tau_{\leq (n+1)}B$ of $B$ is a colimit of this diagram.
\item The transition maps $\pi_{n+1}(L(A_{n+1}'\oinfty\Sigma^{n+1}\cinfty(V')))\rightarrow \pi_{n+1}(L(A_{n+1}'\oinfty\Sigma^{n+1}\cinfty(V'')))$ are monomorphisms. 
\end{enumerate}
Since $B$ is finitely generated to order $n+1$ over $A$, the truncation map $B\rightarrow \tau_{\leq (n+1)}B$ factors through some $\tau_{\leq (n+1)}L(A_{n+1}'\oinfty\Sigma^{n+1}\cinfty(V'))$. It follows that the map $A_{n+1}'\oinfty\Sigma^{n+1}\cinfty(V')\rightarrow B$ is $(n+1)$-connective. We can write $A_{n+1}'\oinfty\Sigma^{n+1}\cinfty(V')$ as a pushout 
\[A_{n+1}''=\R\otimes_{\Sigma^n\cinfty(\R^{S'}\times V')}A_n, \]
so invoking $(1)$, we find an $(n+1)$-finite good $A$-cell decomposition of $B$. Note that $(4)$ follows from $(3)$ by choosing an almost finite good cell decomposition for a finitely presented object $A$ and applying the compactness of $A$ to the identity $\mathrm{id}:A\rightarrow A$.
\end{proof}

\subsubsection{Locality of simplicial $\cinfty$-rings}
We show in this section $\sring^{op}_{\mathrm{fp}}$ has a natural structure of a geometry. 
\begin{defn}
Let $A$ be a simplicial $C^{\infty}$-ring and let $a\in \pi_0(A)$. We say that a map $f:A\rightarrow B$ such that $f(a)\in \pi_0(B)$ is invertible is a \emph{localization of $A$ with respect to $a$} if for each $C\in sC^{\infty}\mathsf{ring}$, the map $\Hom_{sC^{\infty}\mathsf{ring}}(B,C)\rightarrow \Hom_{sC^{\infty}\mathsf{ring}}(A,C)$ given by composition with $f$ induces a homotopy equivalence of Kan complexes
\[ \Hom_{sC^{\infty}\mathsf{ring}}(B,C)\overset{\simeq}{\longrightarrow} \Hom^0_{sC^{\infty}\mathsf{ring}}(A,C),\]
where $\Hom^0_{sC^{\infty}\mathsf{ring}}(A,C)$ is the union of those connected components of $\Hom_{sC^{\infty}\mathsf{ring}}(A,C)$ spanned by those maps $g$ such that $g(a)$ is invertible in the commutative $\R$-algebra $\pi_0(C)^{\rmalg}$. 
\end{defn}
In the case of an ordinary $C^{\infty}$-ring $A$ and some $a\in A$, the above definition reduces to the usual $C^{\infty}$ localization $A[1/a]$ given up to equivalence by the pushout
\begin{equation*}
\begin{tikzcd}
C^{\infty}(\R)\ar[r,"q_a"]\ar[d]& A\ar[d]\\
C^{\infty}(\R\setminus \{0\})\ar[r] & A[1/a]
\end{tikzcd}    
\end{equation*}
of $C^{\infty}$-rings. The localization of a simplicial $C^{\infty}$-ring admits a similar characterization, for which we recall some standard terminology.
\begin{defn}\label{strongmap}
Let $f:A\rightarrow B$ be a morphism of simplicial $\cinfty$-rings, then $f$ is \emph{strong} if the map $f^{\rmalg}:A^{\rmalg}\rightarrow B^{\rmalg}$ of simplicial commutative rings is \emph{strong} (in the sense of To\"{e}n-Vezzosi \cite{TV2}, Definition 2.2.2.1.), that is, the natural map
\[ \tor_{0}^{\pi_0(A)}(\pi_n(A),\pi_0(B))\longrightarrow \pi_n(B)\]
is an isomorphism for all $n\geq 0$. 
\end{defn}
\begin{lem}\label{stabilitystrongmorphism}
The following hold true.
\begin{enumerate}[(1)]
    \item A retract of a strong morphism is strong.
    \item In a diagram 
    \[\begin{tikzcd} & C \ar[dl,"h"'] \ar[dr,"g"]\\
    B\ar[rr,"f"]&& A\end{tikzcd} 
    \] where $h$ is strong, $f$ is strong if and only if $g$ is strong.
\end{enumerate}
\end{lem}
\begin{proof}
\begin{enumerate}[(1)]
    \item A retraction diagram 
    \[ 
    \begin{tikzcd}
    A\ar[d,"g"]\ar[r] &B\ar[d,"f"]\ar[r]& A\ar[d,"g"]\\
    A'\ar[r] &B'\ar[r]& A'\\
    \end{tikzcd}
    \]
    where $f$ is a strong morphism induces for each $n\geq 0$ a diagram  
        \[ 
    \begin{tikzcd}
    \pi_n(A)\ar[d]\ar[r] &\pi_n(B)\ar[d]\ar[r]& \pi_n(A)\ar[d]\\
    \tor^{\pi_0(A)}_0(\pi_n(A),\pi_0(A'))\ar[d]\ar[r] & \tor^{\pi_0(B)}_0(\pi_n(B),\pi_0(B'))\ar[d,"\cong"]\ar[r]&  \tor^{\pi_0(A)}_0(\pi_n(A),\pi_0(A'))\ar[d]\\
    \pi_n(A')\ar[r] &\pi_n(B')\ar[r]& \pi_n(A')\\
    \end{tikzcd}
    \]
    where both horizontal rectangles are retraction diagrams. The inverse of the indicated isomorphism yields a map $\pi_n(A')\rightarrow  \tor^{\pi_0(A)}_0(\pi_n(A),\pi_0(A'))$ which is the inverse of $ \tor^{\pi_0(A)}_0(\pi_n(A),\pi_0(A'))\rightarrow \pi_n(A')$.
    \item We are asked to prove that the map $ \tor^{\pi_0(B)}_0(\pi_n(B),\pi_0(A))\rightarrow \pi_n(A)$ is an isomorphism if and only if the induced map $ \tor^{\pi_0(C)}_0(\pi_n(C),\pi_0(A))\rightarrow \pi_n(A)$ is an isomorphism, given that $ \tor^{\pi_0(C)}_0(\pi_n(C),\pi_0(B))\rightarrow \pi_n(B)$ is an isomorphism. The last isomorphism implies that the map $ \tor^{\pi_0(C)}_0(\pi_n(C),\pi_0(A))\rightarrow  \tor^{\pi_0(B)}_0(\pi_n(B),\pi_0(A))$ is also an isomorphism, so the desired statement reduces to the 2-out-of-3 property for isomorphisms in the commuting diagram
    \[
    \begin{tikzcd}
    \tor^{\pi_0(C)}_0(\pi_n(C),\pi_0(A)) \ar[dr] \ar[rr,"\cong"] &&\tor^{\pi_0(B)}_0(\pi_n(B),\pi_0(A)) \ar[dl] \\
    & \pi_n(A)
    \end{tikzcd}
    \]
\end{enumerate}
\end{proof}

We now recall the notion of a \emph{flat map} of simplicial commutative algebras. 
\begin{defn}\label{def:flatmorphism}
A morphism $f:A\rightarrow B$ of simplicial commutative $\R$-algebras is \emph{flat} (respectively \emph{faithfully flat}) if $\pi_0(f)$ is a flat (respectively faithfully flat) morphism of commutative $\R$-algebras and $f$ is strong. A morphism $f:A\rightarrow B$ of simplicial $\cinfty$-rings is \emph{flat} (respectively \emph{faithfully flat}) if $f^{\rmalg}$ is flat (respectively faithfully flat).
\end{defn}
We record the following permanence properties of flat morphisms. 
\begin{prop}\label{prop:flatproperties}
Let $\mathrm{Flat}\subset\fun(\Delta^1,\scring_{\R})$ be the full subcategory spanned by flat morphisms. Then the following hold.
\begin{enumerate}[$(1)$]
    \item The full subcategory $\mathrm{Flat}\subset\fun(\Delta^1,\scring_{\R})$ is stable under composition.
    \item If $f:A\rightarrow B$ is flat, and $g:A\rightarrow C$ is any morphism in $\scring_{\R}$ then the base change $B\rightarrow A\otimes_B C$ is flat.
    \item $\mathrm{Flat}\subset\fun(\Delta^1,\scring_{\R})$ is stable under filtered colimits.
    \item $\mathrm{Flat}\subset\fun(\Delta^1,\scring_{\R})$ is stable under retracts.
    \item $\mathrm{Flat}\subset\fun(\Delta^1,\scring_{\R})$ is stable under finite products.
    \item If $A$ is a coherent simplicial $\R$-algebra, then the \infcat $\mathrm{Flat}_{A}$ is stable under arbitrary small products. 
\end{enumerate}
\end{prop}
\begin{proof}
All of these assertions are standard. For instance, in the case of $(3)$, we let $f:K\rightarrow \fun(\Delta^1,\scring_{\R})$ be a filtered diagram of flat morphisms. Let $\overline{f}:K^{\rhd}\rightarrow \scring_{\R}$ be a colimit diagram of the composition
\[ K\overset{f}{\longrightarrow} \fun(\Delta^1,\scring_{\R}) \overset{\ev_0}{\longrightarrow} \scring_{\R}, \]
then we have a commuting square 
\[
\begin{tikzcd}
K\ar[d] \ar[r,"f"] & \fun(\Delta^1,\scring_{\R}) \ar[d,"{\ev_0}"] \\
K^{\rhd} \ar[r,"{\overline{f}}"] & \scring_{\R} 
\end{tikzcd}
\]
of \infcats and it is easy to see that an ${\ev_0}$-colimit of this diagram is a colimit of $f$. All the fibres of the coCartesian fibration $\ev_0$ admit colimits, and for each map $g:A\rightarrow B$, the coCartesian pushforward functor $g_!$ can be identified with the base change functor along $g$ which preserves colimits, so $f$ admits an $\ev_0$-colimit. We may compute this colimit by taking a coCartesian transformation $F:K\times\Delta^1 \rightarrow \fun(\Delta^1,\scring_{\R})$ such that $F|_{K\times \{0\}}=f$ and $\ev_0\circ F|_{K\times\{1\}}$ is constant on $\overline{f}(-\infty)$, and taking the colimit of $F|_{K\times \{1\}}$ in $(\scring_{\R})_{\overline{f}(-\infty)/}$. Since base change preserves flatness, $F|_{K\times \{1\}}$ is a filtered diagram of flat morphisms, whose colimit may be computed in $\Mod_{\overline{f}(-\infty)}$; it follows that this colimit is flat. Now $(4)$ follows from $(3)$ since the \infcat $\mathrm{Idem}$ classifying idempotents is filtered.
\end{proof}
\begin{lem}\label{lem:freecinftypreserveloc}
Let $f:A\rightarrow B$ be a morphism of simplicial commutative $\R$ algebras that exhibits $B$ as a localization of $A$ with respect to some $a\in\pi_0(A)$. Then $F^{\cinfty}(A)\rightarrow F^{\cinfty}(B)$ exhibits $F^{\cinfty}(B)$ as a localization of $F^{\cinfty}(A)$ with respect to the image of the element $a$ under the unit map $A\rightarrow F^{\cinfty}(A)^{\rmalg}$.    
\end{lem}
\begin{proof}
For any simplicial $\cinfty$-ring $C$, we have a commuting diagram 
\[
\begin{tikzcd}
\Hom_{\scring_{\R}}(B,C^{\rmalg}) \ar[d,"\simeq"] \ar[r] &   \Hom_{\scring_{\R}}(A,C^{\rmalg}) \ar[d,"\simeq"]\\
\Hom_{\sring}(F^{\cinfty}(B),C) \ar[r] &   \Hom_{\scring}(F^{\cinfty}(A),C)
\end{tikzcd}
\]
of spaces in which the vertical maps are equivalences. The upper horizontal map is an inclusion of the connected components containing the maps $A\rightarrow C^{\rmalg}$ that carry $a$ to an invertible element. Given a map $g:A\rightarrow C^{\rmalg}$, the map on underlying algebras of the adjoint map $F^{\cinfty}(A)\rightarrow C$ of simplicial $\cinfty$-rings is naturally equivalent to the map induced by the universal property of the unit $A\rightarrow F^{\cinfty}(A)^{\rmalg}$, so $g$ carries $a$ to an invertible element in $\pi_0(C^{\rmalg})$ if and only if the map $F^{\cinfty}(A)\rightarrow C$ carries the image of $a$ in $\pi_0(F^{\cinfty}(A))\cong F^{\cinfty}(\pi_0(A))$ to an invertible element in $\pi_0(C)$. 
\end{proof}
\begin{prop}\label{prop:localization}
The following hold true.
\begin{enumerate}[$(1)$]
    \item The restriction map $\cinfty(\R)\rightarrow \cinfty(\R\setminus\{0\})$ exhibits $\cinfty(\R\setminus\{0\})$ as a localization of $\cinfty(\R)$ with respect to $\mathrm{id}_{\R}$.
    \item Let $f:A\rightarrow B$ be a morphism of simplicial $\cinfty$-rings and let $a\in \pi_0(A)$ an element. Let $h:B\rightarrow B'$ be a map for which the composition $A\rightarrow B\rightarrow B'$ carries $a$ to an invertible element. Then $h$ exhibits $B'$ as a localization with respect to $f(a)$ if and only if the diagram 
   \[\begin{tikzcd}
   A\ar[d]\ar[r,"f"] & B\ar[d] \\
    A[a^{-1}] \ar[r] & B'
    \end{tikzcd}\]
    induced by the universal property of $A[a^{-1}]$ is a pushout.
    \item Let $A$ be a simplicial $\cinfty$-ring and $a\in A$ an element, and let $\tilde{a}:\cinfty(\R)\rightarrow A$ be the map classifying $a$, that is, $\tilde{a}(\mathrm{id}_{\R})=a$. Let $h:A\rightarrow A'$ be a map for which the composition $\cinfty(\R)\rightarrow A\rightarrow A'$ carries $\mathrm{id}_{\R}$ to an invertible element. Then $h$ exhibits $A'$ as a localization of $A$ with respect to $A$ if and only if the diagram 
    \[    \begin{tikzcd}
    \cinfty(\R)\ar[d]\ar[r] & A\ar[d] \\
    \cinfty(\R\setminus\{0\}) \ar[r] & A'
    \end{tikzcd}\]
    induced by the universal property of $\cinfty(\R\setminus\{0\})$ from $(1)$ is a pushout. 
    \item Let $M$ be a manifold and $f\in \cinfty(M)$ a function, then the map $\cinfty(M)\rightarrow\cinfty(f^{-1}(\R\setminus\{0\}))$ exhibits $\cinfty(f^{-1}(\R\setminus\{0\}))$ as a localization of $\cinfty(M)$ with respect to $f$.
    \item Let $f:A\rightarrow A'$ be a morphism of simplicial $\cinfty$-rings and let $a\in \pi_0(A)$ an element. Then the following are equivalent. 
    \begin{enumerate}[$(a)$]
        \item The map $f$ exhibits $A'$ as a localization with respect to $a$ 
        \item The following conditions are satisfied.
        \begin{enumerate}[$(i)$]
        \item The map $\pi_0(f):\pi_0(A)\rightarrow \pi_0(A')$ exhibits $\pi_0(A')$ as a localization of $\pi_0(A)$ with respect to $a$.
        \item The map $f^{\rmalg}:A^{\rmalg}\rightarrow A^{\prime\rmalg}$ is a flat map of simplicial commutative $\R$-algebras.
            \end{enumerate}
        \item The following conditions are satisfied.
        \begin{enumerate}
        \item[$(i)$] The map $\pi_0(f):\pi_0(A)\rightarrow \pi_0(A')$ exhibits $\pi_0(A')$ as a localization of $\pi_0(A)$ with respect to $a$.
        \item[$(ii')$] The map $f^{\rmalg}:A^{\rmalg}\rightarrow A^{\prime\rmalg}$ is a strong map of simplicial commutative $\R$-algebras.

    \end{enumerate}
    \end{enumerate}
\end{enumerate} 
\end{prop}
\begin{proof}
The map $\R[x]\rightarrow \R[x,x^{-1}]$ inverting $x$ exhibits $\R[x,x]$ as a localization at the identity map. It follows from Lemma \ref{lem:freecinftypreserveloc} that the map $\cinfty(\R)\rightarrow F^{\cinfty}(\R[x,x^{-1}])$ exhibits a localization with respect to $\mathrm{id}_{\R}$. The map $\R[x]\rightarrow \R[x,x^{-1}]$ can be identified with the composition $\R[x]\rightarrow\R[x,y]\rightarrow \R[x,y]/(xy-1)$. It follows from Lemma \ref{lem:openset} that the free simplicial $\cinfty$-ring functor applied to this map can be identified with the map $\cinfty(\R)\rightarrow\cinfty(\R\setminus\{0\})$, which concludes the proof of $(1)$. The diagram in point $(2)$ is a pushout if and only if for every simplicial $\cinfty$-ring $C$, the diagram 
\[
\begin{tikzcd}
 \Hom_{\sring}(B',C) \ar[d] \ar[r] &    \Hom_{\sring}(B,C)\ar[d] \\
 \Hom_{\sring}(A[a^{-1}],C)  \ar[r] &    \Hom_{\sring}(A,C)
\end{tikzcd}
\]
is a pullback square. But this diagram is a pullback square precisely if the upper horizontal map induces an equivalence on those connected components of $\Hom_{\sring}(B',C)$ containing the maps $B'\rightarrow C$ for which the composition $A\rightarrow B\rightarrow C$ carries $a$ to an invertible element, which is the case if and only if $B\rightarrow C$ carries $f(a)$ to an invertible element. Now $(3)$ is an immediate consequence of $(1)$ and $(2)$ and $(4)$ follows from $(3)$ and Lemma \ref{lem:opentransverse}. We complete the proof by showing $(5)$. It follows from $(2)$, Lemma \ref{lem:localizationsexist1} and the fact that $\pi_0:\sring\rightarrow\cinfty\mathsf{ring}$ preserves colimits that we may suppose that $(i)$ is satisfied. In particular, the element $a$ is invertible in $\pi_0(A')$. Consider the $0$-category $\tau_{\leq -1}\sring_{\cinfty(\R)//A}$ of factorizations $\cinfty(\R)\rightarrow B\overset{h}{\rightarrow} A$ where $h$ is $(-1)$-truncated and the composition classifies $a$ and let $\mathrm{Sub}_a(A)$ be the poset of equivalence classes of objects of $\tau_{\leq -1}\sring_{\cinfty(\R)//A}$, that is, the poset of sub simplicial $\cinfty$-rings that contain $a$. Let $\mathrm{Sub}'_a(A)\subset\mathrm{Sub}_a(A)$ be the subposet spanned by finitely generated subrings, and consider the diagram 
\[\mathrm{Sub}'_a(A)\subset \mathrm{Sub}_a(A) \simeq \tau_{\leq -1}\sring_{\cinfty(\R)//A}\hooklongrightarrow \sring_{\cinfty(\R)//A}, \]
then it follows from the argument in the proof of Proposition \ref{fgcpt} that the diagram $F:\mathrm{Sub}'_a(A)^{\lhd\rhd}\rightarrow \sring$ is a colimit diagram. Let $-\infty,\infty\in \mathrm{Sub}'_a(A)^{\lhd\rhd}$ denote the cone respectively the cocone and consider the diagram 
\[ \mathcal{J}_0:\mathrm{Sub}'_a(A)^{\lhd\rhd}\coprod_{\Delta^{\{-\infty,\infty\}}}  \Delta^{\{-\infty,\infty\}}\times\Delta^1\longrightarrow\sring\]
defined by the amalgamation of the functor $\Delta^{\{-\infty,\infty\}}\times\Delta^1\rightarrow \sring$ classifying the diagram 
 \[    \begin{tikzcd}
    \cinfty(\R)\ar[d]\ar[r] & A\ar[d] \\
    \cinfty(\R\setminus\{0\}) \ar[r] & A'
    \end{tikzcd}\]
induced by the universal property of $\cinfty(\R\setminus\{0\})$ and the diagram $F$. Let \[\mathrm{Sub}'_a(A)^{\lhd\rhd}\coprod_{\Delta^{\{-\infty,\infty\}}}  \Delta^{\{-\infty,\infty\}}\times\Delta^1\overset{i}{\longrightarrow} \icat_0\subset \mathrm{Sub}'_a(A)^{\lhd\rhd}\times\Delta^1\] be the map onto the essential image of the map $\mathrm{Sub}'_a(A)^{\lhd\rhd}\coprod_{\Delta^{\{-\infty,\infty\}}}  \Delta^{\{-\infty,\infty\}}\times\Delta^1\rightarrow \mathrm{Sub}'_a(A)^{\lhd\rhd}\times\Delta^1$, then $i$ is a categorical equivalence, so we may replace $\mathcal{J}_0$ by a functor $\mathcal{J}_0:\icat_0\rightarrow\sring$. Let $\mathcal{J}$ be a left Kan extension of $\mathcal{J}_0$ along the inclusion $\icat_0\subset \mathrm{Sub}'_a(A)^{\lhd\rhd}\times\Delta^1$, then $\mathcal{J}|_{\{A'\}\times\Delta^1}$ exhibits a localization for each $B\in \mathrm{Sub}'_a(A)$ with respect to an element $b$ that the map $B\rightarrow A$ carries to $A$. Using $(3)$ and \cite{HTT}, Propositions 4.3.2.8 and 4.2.3.9, we deduce that if the map $\mathcal{J}|_{\{\infty\}\times\Delta^1}$, that is, the map $A\rightarrow A'$ exhibits a localization with respect to $a$, then the filtered diagram $\mathcal{J}|_{\mathrm{Sub}'_a(A)^{\rhd}\times\{1\}}$ is a colimit diagram. To show $(a)\Rightarrow (b)$, we note that since flat maps are stable under filtered colimits and $(\_)^{\rmalg}$ preserves filtered colimits, it suffices to argue that $A\rightarrow A[a^{-1}]$ is a flat map for $A$ a finitely generated simplicial $\cinfty$-ring. Choose using Proposition \ref{fgcpt} an effective epimorphism $\cinfty(\R^n)\rightarrow A$ and lift the element $a$ to some $\tilde{a}\in \cinfty(\R^n)$, then according to $(2)$ and $(3)$, we have a pushout diagram 
\[
\begin{tikzcd}
\cinfty(\R^n)\ar[d] \ar[r] & \cinfty(\tilde{a}^{-1}(\R\setminus\{0\})) \ar[d] \\
A\ar[r] & A'.
\end{tikzcd}
\]
It follows from Corollary \ref{cor:algeffepi} that the underlying diagram of simplicial commutative $\R$-algebras is a pushout diagram. Now we invoke Lemma \ref{lem:openfunctionisflat} and the stability of flat maps under base change. Flat maps are strong by definition, so $(b)\Rightarrow (c)$. Now we suppose that $A\rightarrow A'$ is strong. The universal property of $A[a^{-1}]$ provides a diagram
\[
\begin{tikzcd}
A\ar[d,equal] \ar[r] & A[a^{-1}] \ar[d] \\
A \ar[r]& A'.
\end{tikzcd}
\]
Invoking Lemma \ref{stabilitystrongmorphism} and the flatness of $A\rightarrow A[a^{-1}]$ just established, we see that $A[a^{-1}]\rightarrow A'$ is strong. Since this map induces an isomorphism on connected components, we conclude.  
\end{proof}

\begin{rmk}
Combining Remark \ref{mfdscptobj} with Proposition \ref{prop:localization} shows that the localization of an (almost) finitely presented simplicial $C^{\infty}$-ring with respect to any $a\in \pi_0(A)$ is again (almost) finitely presented.
\end{rmk}
\begin{cor}\label{cor:truncatedsite}
Let $A$ be a simplicial $C^{\infty}$-ring and let $\sring_{/A}^{op,\mathrm{ad}}$ denote the full subcategory of $\sring^{op}_{/A}$ spanned by localization morphisms. The functor $\tau_{\leq 0}:\sring^{op,\mathrm{ad}}_{/A}\rightarrow (C^{\infty}\mathsf{ring}^{op})_{/\pi_0{A}}^{\mathrm{ad}}$ is an equivalence of $\infty$-categories.
\end{cor}
\begin{proof}
We show that $\tau_{\leq 0}$ is fully faithful and essentially surjective. For essential surjectivity, let $\pi_0(A)\rightarrow B$ be a localization morphism in $C^{\infty}\mathsf{ring}$ determined by some $a\in \pi_0(A)$, then Proposition \ref{prop:localization} immediately shows that $B$ is isomorphic to the image under $\tau_{\leq 0}$ of the morphism $A\rightarrow A[a^{-1}]$. For fully faithfulness, note that the functor $\tau_{\leq 0}$ sends the Hom-spaces in $\sring^{op,\mathrm{ad}}_{/A}$ to their zero'th truncation. Thus, to show that $\tau_{\leq 0}$ is fully faithful, it suffices to show that the Hom-spaces of $\sring^{op,\mathrm{ad}}_{/A}$ are already discrete. Let $A\rightarrow B$ and $A\rightarrow C$ be localization morphisms. The space $\mathrm{Hom}_{\sring^{op,\mathrm{ad}}_{/A}}(C,B)$ is equivalent to the space $\mathrm{Hom}_{\sring_{A/}}(B,C)$, but by \cite{HTT}, Lemma 5.5.5.12, this space fits into a homotopy fibre sequence
\[
\Hom_{\sring_{A/}}(B,C) \longrightarrow  \Hom_{\sring}(B,C) \longrightarrow \Hom_{\sring}(A,C) 
\]
where the fibre is taken over the chosen morphism $A\rightarrow C$. Because $A\rightarrow B$ is a localization, the second map in the fibre sequence is an inclusion of connected components, so $\Hom_{\sring_{A/}}(B,C) $ is empty or weakly contractible.
\end{proof}
\begin{nota} We will denote $\geodiffder$ for the opposite category of the $\infty$-category of compact objects in $sC^{\infty}\mathsf{ring}$. To notationally distinguish a finitely presented simplicial $C^{\infty}$-ring $A$ from $A$ as an object of $\geodiffder$, we use the notation $\mathrm{Spec}\,A$ in the latter case (the next subsection will provide motivation for this notation).
\end{nota}
We endow $\geodiffder$ with the structure of a geometry according to the following prescription: 
\begin{enumerate}[(1)]
    \item A map $f:\mathrm{Spec}\, A\rightarrow \mathrm{Spec}\, B$ in $\geodiffder$ is admissible if and only if there exists an element $b\in \pi_0(B)$ such that the image of $b$ under $f$ is invertible in $\pi_0(A)$ and the induced map $B[1/b]\rightarrow A$ is an equivalence. 
    \item A collection $\{\mathrm{Spec}\,B[1/b_{\alpha}]\rightarrow \mathrm{Spec}\,B\}_{\alpha \in J}$ generates a covering sieve if and only if the germ determined ideal generated by the elements $b_{\alpha}$ in $\pi_0(B)$ contains the unit.
\end{enumerate}

\begin{prop}
The admissible morphisms and admissible coverings determine an admissibility structure with compatible topology on $\geodiffder$ such that the truncation functor $\tau_{\leq 0}:\geodiffder\rightarrow\geodiff$ is a transformation of geometries and exhibits a $0$-stub.
\end{prop}
\begin{proof}
 To see that we have indeed defined an admissibility structure on $\geodiffder$ and a compatible topology, we only have to check that admissible maps are stable under retracts and that in a diagram
\[ 
\begin{tikzcd}
\mathrm{Spec}\,A\ar[rr,"f"]\ar[dr,"g"']&& \mathrm{Spec}\,B\ar[dl,"h"]\\
&\mathrm{Spec}\,C
\end{tikzcd}
\]
where $h$ is admissible, $f$ is admissible if and only if $g$ is admissible. Using characterization $(c)$ of localization morphisms of point $(5)$ Proposition \ref{prop:localization}, and Lemma \ref{stabilitystrongmorphism}, we reduce to the discrete case, which is handled in Proposition \ref{proofgeometrydisc}. Similarly, it follows from Proposition \ref{prop:localization} that a collection $\{A\rightarrow B_i\}_i$ determines an admissible covering if and only if the collection $\{\pi_0(A)\rightarrow \pi_0(B_i)\}_i$ determines an admissible covering in $\geodiff$ so the assertion that these coverings define a topology follows from Proposition \ref{proofgeometrydisc}.  
\end{proof}
It follows from the preceding proposition that composition with $\tau_{\leq 0}$ induces equivalences
\[\str_{\geodiff}(\xtop)\cap \str_{\geodiffder}^{\leq0}(\xtop),\quad   \strloc_{\geodiff}(\xtop)\overset{\simeq}{\longrightarrow} \str^{\leq 0}_{\geodiffder}(\xtop)\cap \strloc(\geo)  \] 
for any \inftop $\xtop$. We wish to define \emph{truncation functors} on the \infcats of $\geodiffder$-structured \inftopoit, left adjoints to the subcategory inclusions above. 

\begin{prop}\label{prop:geoenvopen}
The functor $\cinfty(\_):\diff^{\mathrm{open}}\hookrightarrow \geodiffder$ lies in $\fun^{\pi,\mathrm{ad}}_{\tau}(\diff^{\mathrm{open}},\geodiffder)$ and exhibits $\geodiffder$ as a geometric envelope of $\diff^{\mathrm{open}}$.
\end{prop}
We need an easy lemma.
\begin{lem}\label{lem:cinftysheaf}
Let $M$ be a manifold and let $\{U_i\hookrightarrow M\}_{i}$ be a countable admissible covering $\diff$, so that the map $h:\coprod_i U_i\rightarrow M$ admits a \v{C}ech nerve $\check{C}(h)_{\bullet}$ in $\diff$ which is a colimit diagram. Then the coaugmened cosimplicial diagram diagram $\cinfty(\check{C}(h)_{\bullet})$ is a limit diagram in $\sring$.   
\end{lem}
\begin{proof}
We can identify the diagram $\cinfty(\check{C}(h)_{\bullet})$ with the functor $\cartsp\rightarrow\spa$ carrying $\R^n$ to $\Hom_{\diff}(\check{C}(h)_{\bullet},\R^n)$. Since the evaluation functor $\ev_{\R}:\sring\rightarrow\spa$ is conservative and preserves limits, it suffices to shows that the diagram of spaces $\Hom_{\diff}(\check{C}(h)_{\bullet},\R^n)$ is a limit diagram, which follows immediately from the fact that $\check{C}(h)_{\bullet}$ is a colimit diagram.    
\end{proof}
\begin{proof}[Proof of Proposition \ref{prop:geoenvopen}]
It follows from Proposition \ref{prop:localization} and Lemma \ref{lem:opentransverse} that $\cinfty(\_)$ lies in $\fun^{\pi,\mathrm{ad}}_{\tau}(\diff^{\mathrm{open}},\geodiffder)$. We now show that the functor $\cinfty(\_):\diff^{\mathrm{open}}\rightarrow\sring^{op}$ is fully faithful. We wish to show that for $U,V\in\diff^{\mathrm{open}}$, the canonical map
\[ \Hom_{\diff^{\mathrm{open}}}(U,V)\longrightarrow\Hom_{\sring}(\cinfty(V),\cinfty(U)) \]
is an equivalence. As in Lemma \ref{lem:openset}, we may write $V$ as a transverse pullback of Cartesian spaces so upon invoking Lemma \ref{lem:opentransverse}, we may assume that $V$ is a Cartesian space. Choose a good open cover $\{U_i\subset U\}$, then invoking Lemma \ref{lem:cinftysheaf}, we may assume that $U$ is a Cartesian space, in which case the result is obvious. We now have fully faithful functors
\[ \cartsp\hooklongrightarrow\diff^{\mathrm{open}}\overset{\cinfty(\_)}{\hooklongrightarrow} \geodiffder, \]
which induce a commuting diagram of \infcatst  
\[
\begin{tikzcd}
& \fun^{\pi\mathrm{ad}}(\diff^{\mathrm{open}},\icat)\ar[dr,"\theta''"] \\
\fun^{\mathrm{lex}}(\geodiffder,\icat)\ar[rr,"\theta'"] \ar[ur,"\theta"]&& \fun^{\pi}(\cartsp,\icat) 
\end{tikzcd}.
\]
Note that the functor $\theta$ indeed carries left exact functors to funtors lying in $\fun^{\pi\mathrm{ad}}(\diff^{\mathrm{open}},\icat)$ as $\cinfty(\_)$ preserves finite products and transverse pullbacks in $\diff^{\mathrm{open}}$. It follows from Proposition \ref{finlimittheory} that $\theta'$ is an equivalence so it suffices to argue that $\theta''$ is an equivalence. We first show that $\theta''$ has a right adjoint given by right Kan extension along $\iota:\cartsp\hookrightarrow\diff^{\mathrm{open}}$. Let $F:\cartsp\rightarrow\icat$ be a functor preserving finite products, then Proposition \ref{finlimittheory} asserts that $F$ admits a right Kan extension $\widehat{F}$ along the Yoneda embedding $\cartsp\hookrightarrow \geodiffder$. Since this embedding factors as $ \cartsp\overset{\iota}{\rightarrow}\diff^{\mathrm{open}}\hookrightarrow \geodiffder$ where the second functor is fully faithful, the composition $\iota_*F:\diff^\mathrm{open}\rightarrow\geodiffder\overset{\widehat{F}}{\rightarrow}\icat$ is a right Kan extension of $F$ along $\iota$ which lies in $\fun^{\pi,\mathrm{ad}}(\diff^{\mathrm{open}},\icat)$ by Lemma \ref{lem:opentransverse} again (in other words, the functor taking right Kan extensions right adjoint to $\theta''$ is the compostion of $\theta$ with the inverse of $\theta'$). By construction, the counit map $(\iota_*F)|_{\cartsp}\rightarrow F$ is an equivalence. Let $G:\diff^{\mathrm{open}}\rightarrow \icat$ be a functor preserving finite products and pullbacks along admissibles. The unit map $G\rightarrow \iota_*(G|_{\cartsp})$ is an equivalence when restricted to $\cartsp$. Since both $G$ and $\iota_*G$ preserve pullbacks along admissible maps, it suffices to show that any open submanifold of a Cartesian space can be obtained as a pullback along an open inclusion involving only Cartesian spaces. By the existence of charactersistic functions for open sets, it suffices to show that $\R\setminus\{0\}$ is diffeomorphic to a pullback along an open inclusion involving only Cartesian spaces, but this is easy to arrange: choose a smooth bump function $\psi(x):\R\rightarrow [0,1]$ whose value is equal to $-1/2$ on $(-1,1)$ and equal to $1/2$ on $(-\infty,-2)\cup(2,\infty)$ without local minima or maxima on $(-2,-1)\cup(1,2)$, then $\R\setminus\{0\}$ is diffeomorphic to the intersection of the graph of $\psi$ with the open set $\R^2\cong \R\times\R_{>0}\subset \R^2$. \\
We are left to show that the topology and admissiblity structure on $\geodiffder$ are the coarsest ones for which the functor $\diff^{\mathrm{open}}\hookrightarrow \geodiffder$ is a transformation of pregeometries. It follows from Proposition \ref{prop:localization} that every admissible map of $\geodiffder$ is pulled back from the image of an admissible map in $\diff^{\mathrm{open}}$. Now consider an admissible cover $\{A\rightarrow A[a_i^{-1}]\}_i$. It follows from Lemma \ref{lem:cartspgeneratescovers} that we can find an open cover $\{U_j\subset \R^n\}_j$ and a regular epimorphism $f:\cinfty(\R^n)\rightarrow A$ so that the admissible cover $\{\pi_0(A)\rightarrow \pi_0(A)[a_i^{-1}]\}_i$ is a refinement of the cover $\{\pi_0(A)\rightarrow\pi_0(A)\oinfty_{\cinfty(\R^n)}\cinfty(U_j)\}_j$. We can lift $f$ to an effective epimorphism $f':\cinfty(\R^n)\rightarrow A$ and it follows immediately from Proposition \ref{prop:localization} that the original cover on $A$ is a refinement of the cover $\{A\rightarrow A\oinfty_{\cinfty(\R^n)}\cinfty(U_j)\}_j$. 
\end{proof}
\begin{cor}\label{cor:geoenvcons}
Let $\xtop$ be an \inftopt, then composition with $\cinfty(\_)$ induces an equivalence $\str_{\geodiffder}(\xtop)\simeq  \str_{\diff^{\mathrm{open}}}(\xtop)$ that restricts to an equivalence $\strloc_{\geodiffder}(\xtop)\simeq \strloc_{\diff^{\mathrm{open}}}(\xtop)$.
\end{cor}
\begin{cor}\label{cor:geodiffdersifted}
Let $\xtop$ be an \inftopt, then the subcategory inclusions $\strloc_{\geodiffder}(\xtop)\subset \fun^{\lex,\loc}(\geodiffder,\xtop)$ and $\fun^{\lex,\loc}(\geodiffder,\xtop)\subset \fun^{\lex}(\geodiffder,\xtop)$ preserve sifted colimits.
\end{cor}
\begin{proof}
In view of Corollary \ref{cor:geoenvcons}, this follows from Proposition \ref{prop:pgeotruncfunctor}.
\end{proof}
\begin{cor}\label{cor:localpullback}
Let $\xtop$ be an \inftop and $\Of:\geodiffder\rightarrow\xtop$ a left exact functor, that is, a sheaf of simplicial $\cinfty$-rings on $\xtop$. Then $\Of_{\xtop}$ is a $\geodiffder$-structure if and only if $\pi_0(\Of_{\xtop})$ is a $\geodiffder$-structure. If $\Of':\geodiffder\rightarrow\xtop$ is another left exact functor, then a morphism $f:\Of\rightarrow \Of'$ is local if and only if $\pi_0(\Of)\rightarrow\pi_0(\Of')$ is local. That is, both commuting squares of \infcats in the diagram 
\[
\begin{tikzcd}
\str_{\geodiffder}(\xtop) \ar[d,"\tau_{\leq 0}"]\ar[r] & \fun^{\mathrm{lex}}(\geodiffder,\xtop) \ar[d,"\tau_{\leq 0}"] & \fun^{\mathrm{lex},\mathrm{loc}}(\geodiffder,\xtop)\ar[d,"\tau_{\leq 0}"]\ar[l]\\
\str_{\geodiff}(\xtop) \ar[r] & \fun^{\mathrm{lex}}(\geodiff,\xtop)  & \fun^{\mathrm{lex},\mathrm{loc}}(\geodiff,\xtop)\ar[l]
\end{tikzcd}
\]
are pullback squares.
\end{cor}
\begin{proof}
In view of Corollary \ref{cor:geoenvcons}, this follows immediately from Proposition \ref{prop:pgeotruncfunctor}.
\end{proof}
\begin{prop}\label{prop:truncationrelspec}
Let $\ofxtop$ be a $\geodiffder$-structured \inftopt, then the transformation $(\xtop,\Of_{\xtop})\rightarrow (\xtop,\tau_{\leq 0}\Of_{\xtop})$ of \emph{$\diff^{\mathrm{open}}$-structures} exhibits a unit transformation for the fully faithful embedding $\ltop(\geodiff)\hookrightarrow\ltop(\geodiffder)$. In other words, the relative spectrum $\spec^{\geodiff}_{\geodiffder}$ can be identified with the composition of $\diff^{\mathrm{open}}$-structures with the truncation functor $\tau_{\leq 0}$ on the underlying \inftopoit. 
\end{prop}
\begin{proof}
This follows immediately from Proposition \ref{prop:geotruncfunctor}.
\end{proof}
In the previous section, we studied various notions of locality for $\cinfty$-rings, phrased both in terms of the underlying algebra and in terms of $\geodiff$-structures in spaces. Here, we extend these results to simplicial $\cinfty$-rings.

\begin{prop}\label{prop:localringproperties}
The following hold true.
\begin{enumerate}[$(1)$]
    \item $A$ is a $\geodiffder$-structure in spaces if and only $\pi_0(A)$ is a $\geodiff$-structure, that is, an (Archimedean) local $\cinfty$-ring. Every morphism of $\geodiffder$-structures is a local morphism, that is, the subcategory inclusion $\strloc_{\geodiffder}(\spa)\subset \sring$ is fully faithful.
    \item The object $\R\in \strloc_{\geodiffder}(\spa)$ is final. 
    \item The functor $L_{\R}:\sring_{/\R}\rightarrow\strloc_{\geodiffder}(\spa)_{/\R}\simeq \strloc_{\geodiffder}(\spa)$ carries $p:\cinfty(\R^n)\rightarrow \R$ to the quotient map $\cinfty(\R^n)\rightarrow \cinfty(\R^n)_p$. 
    \item The functor $\strloc_{\geodiffder}(\spa)\rightarrow\sring$ preserves all colimits.
    \item Let $\mathsf{CartSp}_g$ be the Lawvere theory of pointed Cartesian spaces and germs of $\cinfty$-functions between them introduced in Proposition \ref{prop:lawveretheoryofgerms}. Then $\strloc_{\geodiffder}(\spa)$ is canonically equivalent to the \infcat of algebras for $\mathsf{CartSp}_g$. 
\end{enumerate}
\end{prop}
\begin{proof}
The first assertion follows directly from Corollary \ref{cor:localpullback} and the fact that the inclusion $\strloc_{\geodiff}(\spa)\subset\cinfty\mathsf{ring}$ is fully faithful. To prove $(2)$ we note that $(1)$ implies that the inclusion $\strloc_{\geodiffder}(\spa)\subset\sring$ commutes with the truncation functor $\tau_{\leq 0}$. Since $\R$ is a final object in $\strloc_{\geodiff}(\spa)$, it is a final object in $\strloc_{\geodiffder}(\spa)$ as well. To prove $(3)$, we note that the commuting diagram 
\[
\begin{tikzcd}
\strloc_{\geodiffder}(\spa) \ar[d,"\tau_{\leq 0}"] \ar[r] & \sring_{/\R}\ar[d,"\tau_{\leq 0}"] \\
\strloc_{\geodiff}(\spa)\ar[r] & \cinfty\mathsf{ring}_{/\R}
\end{tikzcd}
\]
is vertically right adjointable. Since the horizontal functors admit left adjoint localizations $L_{\R}$, this square is also horizontally left adjointable. Thus, we can identify the composition $\cinfty(\R^n)\rightarrow L_{\R}(\cinfty(\R^n))_p\rightarrow \tau_{\leq0}(L_{\R}(\cinfty(\R^n))_p)$ with the localization at $p$ of $\cinfty(\R^n)$ as an ordinary $\cinfty$-ring. It follows from Lemma \ref{lem:localizationrn} that the latter localization is the quotient map $\cinfty(\R^n)\rightarrow\cinfty(\R^n)/\mathfrak{m}^0_p$, so it suffices to show that $L_{\R}(\cinfty(\R^n))_p$ is $0$-truncated. Since the map $\cinfty(\R^n)\rightarrow L_{\R}(\cinfty(\R^n))_p$ is ind-admissible, it suffices to argue the following.
 \begin{enumerate}
     \item[$(*)$] Let $f:A\rightarrow B$ be an ind-admissible morphism of simplicial $\cinfty$-rings (for the geometry structure on $\geodiffder$). If $A$ is $n$-truncated, then so is $B$.  
 \end{enumerate}
To show $(*)$, we note that since the \infcat of $n$-truncated objects in a presentable \infcat is stable under filtered colimits, it suffices to show that if $f:A\rightarrow B$ exhibits a localization of $A$ with respect to some $a\in A$ and $A$ is $n$-truncated, then $B$ is $n$-truncated. This follows immediately from Proposition \ref{prop:localization} which asserts that localizations of simplicial $\cinfty$-rings are flat maps. We now show $(5)$: first note that as we have an equivalence $\strloc_{\geodiffder}(\spa)_{/\R}\simeq \strloc_{\geodiffder}(\spa)$ and $\strloc_{\geodiffder}(\spa)_{/\R}$ is an accessible localization of $\sring_{\R}$ with localization functor $L_{\R}$, the \infcat $\strloc_{\geodiffder}(\spa)_{/\R}$ is presentable. The functor $\strloc_{\geodiffder}(\spa)\rightarrow \sring_{/\R}$ also preserves sifted colimits by Corollary \ref{cor:geodiffdersifted}, so invoking \cite{HA}, Proposition 7.1.4.12 and Proposition \ref{prop:lawverestability}, we conclude that $\strloc_{\geodiffder}(\spa)$ is projectively generated by the essential image of the objects $\cinfty(\R^n)\rightarrow\R\in \sring_{/\R}$ under $L_{\R}$. It follows from $(3)$ that this essential image can be identified with the opposite of the 1-category $\mathsf{CartSp}_g$. To prove $(4)$, we observe that it suffices to argue that the map
 \[ \cinfty(\R^n)_p\oinfty\cinfty(\R^m)_q\longrightarrow\cinfty(\R^{n+m})_{(p,q)} \]
 is an equivalence. It follows from the left adjointable diagram above and Proposition \ref{prop:lawveretheoryofgerms} that this map becomes an equivalence after applying the truncation functor $\tau_{\leq 0}$, so we are reduced to showing that the object $\cinfty(\R^n)_p\oinfty\cinfty(\R^m)_q$ is $0$-truncated. In view of the assertion $(*)$, it will suffice to show that the map $\cinfty(\R^{n+m})\rightarrow \cinfty(\R^n)_p\oinfty\cinfty(\R^m)_q$ is ind-admissible. This follows immediately from the fact that the class of ind-admissible maps is stable under finite coproducts.
\end{proof}

\subsection{Sheaves of simplicial $\cinfty$-rings and modules}
Our main goal in this subsection is to introduce a spectrum functor 
\[\spec:\{\text{simplicial }\text{$\cinfty$-rings}\}\longrightarrow \{\text{locally simplicial }\text{$\cinfty$-ringed sober spaces}\}  \]
that naturally generalizes the spectrum functor
\[ \spec:\{\text{$\cinfty$-rings}\}\longrightarrow \{\text{locally }\text{$\cinfty$-ringed sober spaces}\}. \]
More precisely, we wish to assign to a simplicial $\cinfty$-ring $A$ a pair $(X,\Of_X)$ where $X$ is a topological space and $\Of_X$ is a sheaf valued in $\sring$ on $X$ that is \emph{local} in the sense determined by the geometry $\geodiffder$. Moreover, we require the following:
\begin{enumerate}[$(G1)$]
\item The locally $\cinfty$-ringed space $(X,\pi_0(X))$ coincides with the Archimedean spectrum of $\pi_0(A)$. In particular, the topological space $X$ coincides with the space $\specr\,\pi_0(A)$.
\item Recall that the global sections and the Archimedean spectrum determine a localization onto the category of complete $\cinfty$-rings which induces an equivalence between the category of geometric $\cinfty$-rings and the Hausdorff, Lindel\"{o}f and $\cinfty$-regular locally $\cinfty$-ringed spaces. We wish to establish a similar equivalence between a suitable subcategory of geometric simplicial $\cinfty$-rings and Hausdorff, Lindel\"{o}f and $\cinfty$-regular locally simplicial $\cinfty$-ringed spaces.    
\end{enumerate}
There are in fact several incarnations of the \infcat of locally simplicial $\cinfty$-ringed spaces; this is a consequence of the fact that there are several \emph{different} fully faithful embeddings 
\[\{\text{sober topological spaces}\}\hooklongrightarrow \{\text{$\infty$-topoi}\},  \]
which we investigate below.
\subsubsection{Digression: localic and complicial \inftopoit}
Let $\rtop_{0}\simeq \mathsf{Locale}$ be the 1-category of $0$-topoi and geometric morphisms between them, that is, the 1-category of locales. This category contains the 1-category of sober topological spaces as a coreflective full subcategory via the assignment $X\mapsto\mathrm{Open}(X)$. 
\begin{enumerate}[$(1)$]
\item The truncation functor $\tau_{\leq -1}: \rtop\rightarrow \rtop_{0}\simeq \mathsf{Locale}$ admits a fully faithful right adjoint, which carries a locale $\mathcal{U}$ to its \inftop of sheaves $\shv(\mathcal{U})$, that is, the Bousfield localization of the presheaf \infcat $\pshv(\mathcal{U})$ localized at the \v{C}ech covers. The assignment $\mathcal{U}\mapsto \shv(\mathcal{U})$ determines a fully faithful embedding of $\mathsf{Locale}$ into $\rtop$ which is right adjoint to the truncation functor $\tau_{\leq -1}$.
\item We may stipulate that equivalences of higher sheaves are detected on points: we can assign to a locale $\mathcal{U}$ the underlying \infcat of the simplicial model category $\pshv_{\sset}(\mathcal{U})$ of simplicial presheaves on $\mathcal{U}$ equipped with Jardine's local model structure, in which a map $F\rightarrow G$ is a weak equivalence if for each point $x:\{0,1\}\rightarrow\mathcal{U}$, the map $F_x\rightarrow G_x$ is weak homotopy equivalence. The resulting \inftop may be identified with the hypercompletion $\shv(\mathcal{U})
\rightarrow \shv(\mathcal{U})^{\mathrm{hyp}}$. The assignment $\mathcal{U}\mapsto \shv(\mathcal{U})^{\mathrm{hyp}}$ determines another fully faithful embedding $\mathsf{Locale}\rightarrow\rtop^{\mathrm{hc}}$.
\item We may stipulate that a higher sheaf is well approximated by its truncations: for each $n\geq 0$, we can assign to a locale $\mathcal{U}$ the $n$-topos $\shv_{\tau_{\leq (n-1)}\spa}(\mathcal{U})$ and define an \infcat as the limit of the tower 
\[ \longrightarrow \shv_{\tau_n\spa}(\mathcal{U})\overset{\tau_{\leq(n-1)}}{\longrightarrow} \shv_{\tau_{(n-1)\spa}}(\mathcal{U}) \overset{\tau_{\leq(n-2)}}{\longrightarrow} \ldots\longrightarrow \shv_{\tau_{\leq 0}\spa}(\mathcal{U})\overset{\tau_{\leq(-1)}}{\longrightarrow} \mathcal{U}.  \]
The limit $\{\shv_{\tau_{\leq n}\spa}(\mathcal{U})\}_{n\in\Z_{\geq0}}$ is again an \inftop and coincides with the Postnikov completion $\shv(\mathcal{U})\rightarrow \widehat{\shv}(\mathcal{U})$. The assignment $\mathcal{U}\mapsto \widehat{\shv}(\mathcal{U})$ determines another fully faithful embedding $\mathsf{Locale}\rightarrow \rtop^{\mathrm{Pc}}$. 
\end{enumerate}
By definition, the essential image of the first embedding $\mathsf{Locale}\hookrightarrow\rtop$ consists of \emph{$0$-localic} \inftopoit. The three constructions above do not coincide; in paricular, the hypercompletion or Postnikov completion of an $n$-localic \inftop for $n\geq 0$ is usually not $n$-localic. We wish to provide an intrinsic characterization of the class of hyper/Postnikov complete \inftopoi that arise as the hyper/Posntikov completion of an $n$-localic \inftopt. 
\begin{defn}\label{defn:complicial}
Let $\xtop$ be an \inftop and let $n\geq -1$ be an integer. We say that $\xtop$ is 
\begin{enumerate}[$(1)$]
\item \emph{$n$-complicial} if $\xtop$ is generated under colimits by $n$-truncated objects. We let $\ltop_{n-\mathrm{com}}\subset\ltop$ denote the full subcategory spanned by $n$-complicial \inftopoit.
\item \emph{weakly $n$-complicial} if for each integer $k\geq -1$, the $(k+1)$-topos $\tau_{\leq k}\xtop$ is generated under colimits by $n$-truncated objects (note that this condition is vacuous for $k\leq n$). We let $\ltop_{n-\mathrm{wcom}}\subset\ltop$ denote the full subcategory spanned by weakly $n$-complicial \inftopoit.
\end{enumerate}
\end{defn}
\begin{rmk}\label{rmk:complicial}
As the condition of being weakly $n$-complicial depends only on the truncations $\{\tau_{\leq k}\xtop\}_{k}$, the following are equivalent.
\begin{enumerate}[$(a)$]
\item The \inftop $\xtop$ is weakly $n$-complicial.
\item The bounded reflection of $\xtop$ is weakly $n$-complicial.
\item The hypercompletion of $\xtop$ is weakly $n$-complicial.
\item The Postnikov completion of $\xtop$ is weakly $n$-complicial.
\end{enumerate}
\end{rmk}
For bounded \inftopoi (that is, \inftopoi that arise as filtered colimits of localic \inftopoit) the condition of being (weakly) $n$-complicial is equivalent to being $(n+1)$-localic.
\begin{lem}\label{lem:boundedcomplicial}
Let $\xtop$ be a bounded \inftopt, then the following are equivalent. 
\begin{enumerate}[$(1)$]
\item $\xtop$ is $n$-localic.
\item $\xtop$ is $(n-1)$-complicial.
\item $\xtop$ is weakly $(n-1)$-complicial.
\end{enumerate}
\end{lem}
\begin{proof}
Since $(1)$ asserts that $\xtop\simeq\shv(\icat)$ for $\icat$ an $n$-category, the implications $(1)\Rightarrow (2)\Rightarrow (3)$ are obvious. Suppose that $\xtop$ is bounded and weakly $(n-1)$-complicial. The bounded reflection of $\xtop$ is the colimit $\colim_{k}L_k\xtop$ in $\ltop$ of the cotower of $n$-localic reflections, so it suffices to show that for $k\geq n$, the functor $L_k\xtop\rightarrow L_{k+1}\xtop$ is an equivalence. Since $\tau_{\leq k}\xtop$ is generated by $(n-1)$-truncated objects for all $k$, this amounts to the following assertion: if $\xtop$ is $m$-localic for some $m\geq n$ and generated under small colimits by $(n-1)$-truncated objects, then $\xtop$ is $n$-localic. To see this, we use \cite{HTT}, Lemma 6.4.3.5 and choose an uncountable regular cardinal $\kappa$ such that $\xtop$ is $\kappa$-compactly generated and the collection of $\kappa$-compact objects in $\xtop$ is stable under finite limits and subobjects in $\xtop$. Then the collection $\tau_{\leq (n-1)}\xtop^{\kappa}$ of $\kappa$-compact objects of $\tau_{\leq (n-1)}\xtop$ has the same properties (with respect to $\xtop$). It follows from \cite{HTT}, Corollary 5.5.7.4 that $\tau_{\leq (n-1)}\xtop^{\kappa}$ generates the collection of $(n-1)$-truncated objects of $\xtop$ under $\kappa$-filtered colimits, so $\tau_{\leq (n-1)}\xtop^{\kappa}$ generates all of $\xtop$ under small colimits. It follows from the proof of \cite{HTT}, Proposition 6.4.3.6 that the algebraic morphism of $m$-topoi $\shv_{\leq (m-1)\spa}(\tau_{\leq (n-1)}\xtop^{\kappa})\rightarrow\tau_{\leq (m-1)}\xtop$ is an equivalence. Since $\xtop$ is $m$-localic, we have an equivalence $\shv(\tau_{\leq (n-1)}\xtop^{\kappa})\simeq \xtop$ so that $\xtop$ is $n$-localic.
\end{proof}
\begin{prop}\label{prop:ncompiffhypncomp}
Let $\xtop$ be an \inftop and $n\geq -1$ an integer, then the following hold true.
\begin{enumerate}[$(1)$]
\item Suppose that $\xtop$ is hypercomplete, then $\xtop$ is $n$-complicial if and only if there is an algebraic morphism $\shv(\icat)\rightarrow\xtop$ from an $(n+1)$-localic \inftop that induces an equivalence $\shv^{\mathrm{hyp}}(\icat)\simeq\xtop$. 
\item The \inftop $\xtop$ is weakly $n$-complicial if and only if the bounded reflection of $\xtop$ is $(n+1)$-localic.
\end{enumerate}
\end{prop}
\begin{proof}
For $(1)$, we note that the hypercompletion of an $(n+1)$-localic \inftop is clearly $n$-complicial. For the converse, suppose that $\xtop$ is $n$-complicial. The $(n+1)$-topos $\tau_{\leq n}\xtop$ is of the form $\shv_{\tau_{\leq n}\spa}(\icat)$ for $\icat$ a small $(n+1)$-category with finite limits and equipped with some Grothendieck topology. Since $\xtop$ is hypercomplete and $n$-complicial, the induced functor $\icat\hookrightarrow\xtop$ induces an equivalence $\shv(\icat)^{\mathrm{hyp}}\simeq \xtop$. Note that $(2)$ follows immediately from Remark \ref{rmk:complicial} and Lemma \ref{lem:boundedcomplicial}.
\end{proof}
\begin{cor}
For each integer $n\geq -1$, the equivalence $\ltop^{\mathrm{b}}\simeq \ltop^{\mathrm{Pc}}$ induced by Postnikov completion and bounded reflection restricts to an equivalence $\ltop^{\mathrm{Pc}}_{n-\mathrm{wcom}}\simeq \ltop_{(n+1)-\loc}$ between Postnikov complete weakly $n$-complicial \inftopoi and $(n+1)$-localic \inftopoit.
\end{cor}
\begin{rmk}
It follows from the preceding corollary that the \infcat of Postnikov complete weakly $(-1)$-complicial \inftopoi is equivalent to the 1-category of locales, and this equivalence is implemented by the truncation functor $\tau_{\leq 0}:\rtop^{\mathrm{Pc}}_{(-1)-\mathrm{wcom}}\rightarrow\mathsf{Locale}$ and the assignment $\mathcal{U}\mapsto \widehat{\shv}(\mathcal{U})$. We say that a Postnikov complete weakly $(-1)$-complicial \inftop is \emph{spatial} if its associated locale is spatial. Let $\rtop^{\mathrm{Pc},\mathrm{sp}}_{(-1)-\mathrm{wcom}}$ be the 1-category of spatial Postnikov complete weakly $(-1)$-complicial \inftopoi, then we have an equivalence $ \rtop^{\mathrm{Pc},\mathrm{sp}}_{(-1)-\mathrm{wcom}}\simeq\mathsf{Top}$ given by sending a sober topological space to $\widehat{\shv}(X)$.
\end{rmk}

\subsubsection{Archimedean spectra of simplicial $\cinfty$-rings}
Let $\mathsf{Top}(\sring)$, $\mathsf{Top}^{\mathrm{hc}}(\sring)$ and $\mathsf{Top}^{\mathrm{Pc}}(\sring)$ be the \infcats of pairs $(X,\Of_X)$ of a sober topological space equipped with a sheaf of simplicial $\cinfty$-rings, a hypercomplete sheaf of simplicial $\cinfty$-rings and a Postnikov complete sheaf of $\cinfty$-rings respectively (we will define these \infcats rigorously below), then we have equivalences 
\[ \tau_{\leq 0}(\mathsf{Top}(\sring))\simeq \tau_{\leq 0}\mathsf{Top}^{\mathrm{hc}}(\sring)\simeq \tau_{\leq 0}\mathsf{Top}^{\mathrm{Pc}}(\sring)\simeq \mathsf{Top}(\cinfty\mathsf{ring}).\]
We have a commuting diagram 
\[
\begin{tikzcd}
\mathsf{Top}(\cinfty\mathsf{ring})^{op} \ar[r,"\Gamma"]\ar[d,hook] & \cinfty\mathsf{ring}\ar[d,hook] \\
\mathsf{Top}(\sring)^{op}\ar[r,"\Gamma"] & \sring,
\end{tikzcd}
\]
and similar commuting diagrams for $\mathsf{Top}^{\mathrm{hc}}(\sring)$ and $\mathsf{Top}^{\mathrm{Pc}}(\sring)$, so if we can establish a left adjoint to the lower horizontal map restricted to \emph{locally} simplicial $\cinfty$-ringed spaces, we automatically satisfy $G1$. However, to satisfy $G2$, we will need to restrict our attention to the Postnikov complete version $\mathsf{Top}^{\mathrm{Pc}}(\sring)$.
\begin{rmk}
We incarnate locales as Postnikov complete weakly $(-1)$-complicial \inftopoi because contrary to the situation in algebraic geometry, we cannot guarantee that the structure sheaf of an arbitrary affine derived $\cinfty$-scheme is Postnikov complete; passing to the Postnikov completion leads to a somewhat more satisfactory theory. In practice, we will rarely be concerned with simplicial $\cinfty$-rings for which $\shv(\specr\,A)$ is not already Postnikov complete (if $A$ is finitely generated for instance, then $\shv(\specr\,A)$ is Postnikov complete) so the reader who is not familiar with our nonstandard choice of embedding 
\[\{\text{sober topological spaces}\}\hooklongrightarrow \{\text{$\infty$-topoi}\}  \]
is encouraged to assume throughout that $\shv(\specr\,A)$ is Postnikov complete and disregard the additional layer of formalism.
\end{rmk}
From the geometry $\geodiffder$ we deduce the existence of the spectrum functor 
\[\mathbf{Spec}^{\geodiffder}:\mathrm{Pro}(\geodiffder)\simeq sC^{\infty}\mathsf{ring}\longrightarrow \ltop(\geodiffder)\]
left adjoint to the global sections functor on the \infcat of $\geodiffder$-structured \inftopoit.
\begin{defn}
Define the \infcats 
\[\ltop^{\mathrm{sp}}_{0-\mathrm{loc}}(\geodiffder) := \ltop_{0-\mathrm{loc}}(\geodiffder)\times_{\ltop_{0-\mathrm{loc}}}\ltop^{\mathrm{sp}}_{0-\mathrm{loc}},\quad\quad \ltop^{\mathrm{Pc}}(\geodiffder):= \ltop(\geodiffder)\times_{\ltop}\ltop^{\mathrm{Pc}}, \]
and define left adjoints 
\[L_{\mathrm{sp}}(\geodiffder):\ltop_{0-\mathrm{loc}}(\geodiffder)\longrightarrow \ltop_{0-\mathrm{loc}}^{\mathrm{sp}}(\geodiffder),\quad\quad L_{\mathrm{Pc}}(\geodiffder):\ltop(\geodiffder)\longrightarrow \ltop^{\mathrm{Pc}}(\geodiffder) \]
to the inclusions 
\[\ltop_{0-\mathrm{loc}}^{\mathrm{sp}}(\geodiffder)\subset \ltop_{0-\mathrm{loc}}(\geodiffder),\quad \quad \ltop^{\mathrm{Pc}}(\geodiffder)\subset \ltop(\geodiffder)\]
using Proposition \ref{prop:pullbackcocartleftadj}. We define a functor as the composition
\[ \sring\overset{\spec^{\geodiffder}}{\longrightarrow }\ltop_{0-\mathrm{loc}}(\geodiffder)\overset{L_{\mathrm{sp}}(\geodiffder)}{\longrightarrow}\ltop_{0-\mathrm{loc}}^{\mathrm{sp}}(\geodiffder)\subset \ltop(\geodiffder)\overset{L_{\mathrm{Pc}}(\geodiffder)}{\longrightarrow} \ltop^{\mathrm{Pc}}(\geodiffder)\]
and denote it $\spec$. Since the composition of the first two functors takes values in the \infcat of spatial $0$-localic \inftopoit, the total composition takes values in spatial weakly -1-complicial Postnikov complete \inftopoit; we thus have a functor 
\[ \spec: \sring\longrightarrow \ltop^{\mathrm{Pc},\mathrm{sp}}_{(-1)-\mathrm{wcom}}(\geodiffder),\]
where the codomain is the pullback $\ltop^{\mathrm{Pc}}(\geodiffder)\times_{\ltop^{\mathrm{Pc}}}\ltop^{\mathrm{Pc},\mathrm{sp}}_{(-1)-\mathrm{wcom}}$ and $\ltop^{\mathrm{Pc},\mathrm{sp}}_{(-1)-\mathrm{wcom}}\subset \ltop$ is the full subcategory spanned by weakly $(-1)$-complicial Postnikov complete \inftopoit.  
\end{defn}
\begin{rmk}
The construction $X\mapsto \widehat{\shv}(X)$ with $\widehat{\shv}(X)$ the Postnikov completion of $\shv(X)$ determines an equivalence of $\mathsf{Top}^{op}\simeq \ltop^{\mathrm{Pc},\mathrm{sp}}_{(-1)-\mathrm{wcom}}$ and we can identify $\ltop^{\mathrm{Pc},\mathrm{sp}}_{(-1)-\mathrm{wcom}}(\geodiffder)$  with the coCartesian fibration associated to the functor $\mathsf{Top}^{op}\rightarrow \catinfh$ that takes a sober topological space $X$ to the \infcat $\strloc_{\geodiffder}(\widehat{\shv}(X))$. 
\end{rmk}
\begin{nota}
For convenience down the line and to ease the notational burden, we will shift our perspective on structured topoi, in accordance with Remark \ref{rmk:exponentiate}. 
\begin{enumerate}[$(1)$]
\item We will denote $\ltop((\geodiffder)_{\mathrm{disc}})$ by $\ltop(\sring)$. Its objects are \emph{simplicial $\cinfty$-ringed \inftopoit.}
\item The subcategory $\ltop(\geodiffder)\subset \ltop((\geodiffder)_{\mathrm{disc}})$ will be denoted $\ltop^{\loc}(\sring)$. Its objects are \emph{locally simplicial $\cinfty$-ringed \inftopoit.}
\end{enumerate}
In other words, we prefer to think of an object $\ofxtop$ as an \inftop equipped with a sheaf $\Of_{\xtop}:\xtop^{op}\rightarrow\sring$, rather than as being equipped with left exact functor $\geodiffder\rightarrow\xtop$.
Similarly, we introduce the following notation.
\begin{enumerate}[$(1)$]
\item We will denote $\rtop^{\mathrm{Pc},\mathrm{sp}}_{(-1)-\mathrm{wcom}}((\geodiffder)_{\mathrm{disc}})$ by $\mathsf{Top}^{\mathrm{Pc}}(\sring)$. In light of the previous remark, its objects are \emph{(Postnikov) complete) simplicial $\cinfty$-ringed spaces}; pairs $(X,\Of_X)$ of a topological space $X$ together with a Postnikov complete sheaf $\Of_X:\widehat{\shv}(X)\rightarrow\sring$ of simplicial $\cinfty$-rings. 
\item The subcategory $\rtop^{\mathrm{Pc},\mathrm{sp}}_{(-1)-\mathrm{wcom}}(\geodiffder)$ will be denoted $\mathsf{Top}^{\mathrm{Pc},\loc}(\sring)$. Its objects are \emph{(Postnikov complete) locally simplicial $\cinfty$-ringed spaces.} 
\end{enumerate}
\end{nota}
We show that limits of locally simplicial $\cinfty$-ringed spaces are taken in simplicial $\cinfty$-ringed spaces.
\begin{prop}\label{prop:structurespacelimits}
The subcategory inclusion $\mathsf{Top}^{\mathrm{Pc},\loc}(\sring)\subset \mathsf{Top}^{\mathrm{Pc}}(\sring)$ is fully faithful and is stable under small limits. In particular, $\mathsf{Top}^{\mathrm{Pc},\loc}(\sring)$ admits all small limits.    
\end{prop}
\begin{proof}
The same proof as the one of Corollary \ref{cor:geodiffstableunderlimits} applies, using Proposition \ref{prop:localringproperties} in place of Proposition \ref{prop:lawveretheoryofgerms}.    
\end{proof}
\begin{rmk}
We have formulated the previous result for Postnikov complete structure sheaves but it also clearly holds for the plain and hypercomplete variants.
\end{rmk}
\begin{rmk}
By construction $\spec$ is a left adjoint to the global sections functor 
\[\mathsf{Top}^{\mathrm{Pc},\loc}(\sring)^{op}\subset \ltop(\sring)^{op}\overset{\Gamma}{\longrightarrow}\sring.\]
The functor $\spec$ performs the following steps.
\begin{enumerate}[$(1)$]
\item First, we carry the a simplicial $\cinfty$-ring $A$ to the $\geodiffder$-spectrum, whose underlying \inftop is the \infcat of sheaves on $(\sring_{A/})^{\mathrm{ad},op}$. This \inftop is $0$-localic but need not be spatial, as we saw in the previous section.
\item As in the previous section, we apply the spatial reflection to obtain a pair $(\specr\,A,\Of_{\specr\,A})$, a derived version of the Archimedian spectrum.
\item Finally, we take the Postnikov completion to obtain a pair $(\specr\,A,\Of_{\specr^{\wedge}\,A})$. 
\end{enumerate}
\end{rmk} 
\begin{rmk}
By Proposition \ref{prop:comparisonringedspaces}, the full subcategory of the \infcat $\ltop^{\mathrm{sp}}_{0-\loc}(\sring)$ spanned by pairs $(\xtop,\Of_{\xtop})$ for which $\Of_{\xtop}$ is 0-truncated (that is, it either takes $0$-truncated values as a functor $\geodiffder\rightarrow\xtop$ or as a sheaf $\xtop^{op}\rightarrow \sring$; see Lemma \ref{lem:truncatedsheaves}) may be identified with the 1-category $\mathsf{Top}(\cinfty\mathsf{ring})^{op}$, and this identification agrees with the subcategories of locally $\cinfty$-ringed spaces and local morphisms. The operation of Postnikov completion yields for each \inftop $\xtop$ a unit map $\xtop\rightarrow \widehat{\xtop}$ which induces an equivalence on $n$-truncated objects for all $n$. It follows that the full subcategory of the \infcat $\mathsf{Top}^{\mathrm{Pc}}(\sring)$ spanned by pairs $(X,\Of_X)$ for which $\Of_X$ is $0$-truncated is also equivalent the 1-category $\mathsf{Top}(\cinfty\mathsf{ring})$, and this equivalence persists on locally ringed objects and local morphisms. It follows from Proposition \ref{prop:truncationrelspec} that the assignment
\[ (X,\Of_{X})\longmapsto (X,\pi_0(\Of_X)) \]
where we view $\Of_X$ as a $\diff^{\mathrm{open}}$-structure (or simply as a functor $\cartsp\rightarrow\widehat{\shv}(X)$) determines a right adjoint to the inclusion $\mathsf{Top}(\cinfty\mathsf{ring})\subset \mathsf{Top}^{\mathrm{Pc}}(\sring)$ which fits as the right vertical map into a commuting diagram 
\[
\begin{tikzcd}
\sring\ar[d,"\pi_0"] \ar[r,"\spec"] & \mathsf{Top}^{\mathrm{Pc},\loc}(\sring)^{op}\ar[d,"\pi_0"] \\
\cinfty\mathsf{ring}\ar[r,"\spec"] & \mathsf{Top}^{\loc}(\cinfty\mathsf{ring})^{op}
\end{tikzcd}
\]
of left adjoints. In particular, for $A$ a simplicial $\cinfty$-ring, there is a canonical isomorphism between $(\specr\,A,\pi_0(\Of_{\specr^{\wedge}\,A}))$ and the Archimidean specrum $(\specr\,\pi_0(A),\Of_{\specr\,\pi_0(A)})$ of $\pi_0(A)$.
\end{rmk}
Below we define various subcategories of derived $\cinfty$-schemes, based one some full subcategories of $\sring$ for which the operations $\spec$ and $\Gamma$ interact well with the homotopy theory of $\sring$.
\begin{defn}[Topological properties]
Let $A$ be a simplicial $\cinfty$-ring. Let $P$ be a property of topological spaces. We say that a $A$ \emph{has the property $P$} if the topological space $\specr\,\pi_0(A)$ has the property $P$ (for instance, we say that $A$ is Lindel\"{o}f, paracompact, locally compact,... if the topological space $\specr\,A$ is Linedel\"{o}f, paracompact, locally compact...)
\end{defn}
From the previous section, recall the notions of \emph{geometric} $\cinfty$-rings and \emph{geometric} modules.
\begin{defn}[Geometricity]
Let $A$ be a simplicial $\cinfty$-ring and $n\geq 0$ an integer. Then $A$ is \emph{geometric to order $n$} if the ordinary $\cinfty$-ring $\pi_0(A)$ is geometric in the sense of Definition \ref{defn:geometricring} and for every $k\leq n$, the $\pi_0(A)$-module $\pi_k(A)$ is geometric. The full subcategory spanned by objects geometric to order $n$ is denoted $\sring_{n-\mathrm{gmt}}\subset\sring$.
\end{defn} 
If $A$ is an ordinary $\cinfty$-ring, then $A$ being finitely presented implies that $A$ is geometric. We will prove a similar statement in the derived setting, for which we need to introduce the following hierarchy of finiteness conditions.
\begin{defn}[Finiteness properties]
Let $A$ be a simplicial $\cinfty$-ring.
\begin{enumerate}[$(1)$]
 \item $A$ is \emph{of finite generation to order $n$} if $A$ is compact to order $n$ in $\sring$. The full subcategory spanned by objects of finite generation to order $n$ is denoted $\sring_{n-\mathrm{fg}}\subset\sring$.
    \item $A$ is \emph{almost of finite presentation} if for every $n\geq 0$, $A$ is of finite generation to order $n$. 
    \item $A$ is \emph{of finite presentation} if $A$ is compact in $\sring$. We let $\sring_{\fp}$ denote the full subcategory spanned by objects of finite presentation.
\end{enumerate}
\end{defn}

\begin{defn}\label{affinederivedmfds}
Let $\ofxtop$ be a locally simplicial $\cinfty$-ringed \inftopt.
\begin{enumerate}[$(1)$]
    \item Suppose that there exists a simplicial $\cinfty$-ring $A$ and an equivalence $\spec\,A\simeq \ofxtop$, then $\ofxtop$ is an \emph{affine derived $\cinfty$-scheme}. We denote the full subcategory spanned affine derived $\cinfty$-schemes by $\daff\subset\rtop^{\loc}(\sring)$.
    \item Suppose that there exists a Lindel\"{o}f simplicial $\cinfty$-ring $A$ and an equivalence $\spec\,A\simeq \ofxtop$, then $\ofxtop$ is an \emph{affine derived geometric $\cinfty$-scheme}. We denote the full subcategory spanned affine derived geometric $\cinfty$-schemes by $\daff_{\gmt}\subset\rtop^{\loc}(\sring)$.
    \item Suppose that there exists a finitely generated simplicial $\cinfty$-ring $A$ and an equivalence $\spec\,A\simeq \ofxtop$, then $\ofxtop$ is an \emph{affine derived $\cinfty$-scheme of finite generation to order $n$}. We denote the full subcategory spanned affine derived $\cinfty$-schemes of finite generation to order $n$ by $\daff_{n-\mathrm{fg}}\subset\rtop^{\loc}(\sring)$. The full subcategory spanned by objects finitely generated (to order $0$) is also denoted $\daff_{\mathrm{fg}}$.
    \item Suppose that there exists a finitely presented simplicial $\cinfty$-ring $A$ and an equivalence $\spec\,A\simeq \ofxtop$, then $\ofxtop$ is an \emph{affine derived $\cinfty$-scheme of finite presentation}. We denote the full subcategory spanned affine derived $\cinfty$-schemes of finite presentation by $\daff_{\fp}\subset\rtop^{\loc}(\sring)$.
\end{enumerate}
\end{defn}
\begin{rmk}
As we show below, an affine derived $\cinfty$-scheme $\ofxtop$ is geometric if and only if there exists a geometric simplicial $\cinfty$-ring $A$ and an equivalence $\spec\,A\simeq \ofxtop$.  
\end{rmk}
We can now state the main result of this section.
\begin{thm}\label{thm:spectrumglobalsections}
The following hold true.
\begin{enumerate}[$(1)$]
\item Let $A$ be a Lindel\"{o}f simplicial $\cinfty$-ring, then $A$ is geometric if and only if the unit map $A\rightarrow\Gamma(\spec\,A)$ is an equivalence, and the functor $\Gamma\circ \spec$ exhibits a reflection of $\sring_{\mathrm{Lin}}$ onto $\sring_{\gmt}$.
\item Let $A$ be a Lindel\"{o}f simplicial $\cinfty$-ring, then the map 
\[ \spec\,A\longrightarrow \spec\,\Gamma(\spec\,A) \]
induced by the unit is an equivalence, that is, the adjunction $(\spec\adj\Gamma)$ restricts to an idempotent adjunction between the \infcat of Lindel\"{o}f simplicial $\cinfty$-rings and $\daff_{\gmt}$.
\item For $\ofxtop$ a locally simplicial $\cinfty$-ringed \inftopt, there is a Lindel\"{o}f simplicial $\cinfty$-ring $A$ and an equivalence $\spec\,A\simeq\ofxtop$ if and only if the following conditions are satisfied.
\begin{enumerate}[$(a)$]
\item There is a Lindel\"{o}f Hausdorff topological space $X$ and an equivalence $\widehat{\shv}(X)\simeq\xtop$. 
\item The locally $\cinfty$-ringed space $(X,\pi_0(\Of_\xtop))$ is $\cinfty$-regular.
\end{enumerate}
\end{enumerate} 
\end{thm}
\begin{cor}
The adjunction $(\spec\adj\Gamma)$ induces an equivalence between the \infcat of geometric simplicial $\cinfty$-rings and the full subcategory of $\ltop^{\loc}(\sring)$ spanned by the pairs $\ofxtop$ satisfying conditions $(a)$ and $(b)$ of Theorem \ref{thm:spectrumglobalsections}.
\end{cor}
We will postpone the proof of Theorem \ref{thm:spectrumglobalsections} until after we have discussed modules of $\cinfty$-ringed \inftopoi and spectra of modules. Instead, we apply Theorem \ref{thm:spectrumglobalsections} to prove that the functor $\cinfty(\_):\diff\rightarrow\sring$ is also a geometric envelope. We rely on the following result, the proof of which we also postpone until the next subsection.
\begin{prop}\label{afpisfair}
Let $f:A\rightarrow B$ be a morphism of Lindel\"{o}f simplicial $\cinfty$-rings. Suppose that the following conditions are satisfied.
\begin{enumerate}[$(1)$]
    \item $A$ is geometric to order $m$.
    \item $f$ exhibits $B$ as finitely generated to order $m+1$ over $A$.
    \item $f$ is an effective epimorphism.
\end{enumerate}
Then $B$ is geometric to order $m$. 
\end{prop}
\begin{cor}\label{cor:afpgeometric}
Let $A$ be geometric and suppose that $f:A\rightarrow B$ is an effective epimorphism that exhibits $B$ as almost finitely presented over $A$. Then $B$ is geometric.
\end{cor}
\begin{cor}
If $A$ is an almost finitely presented simplicial $\cinfty$-ring, then $A$ is geometric.
\end{cor}
\begin{prop}\label{manifoldssmoothring}
The functor $C^{\infty}(\_):\diff\rightarrow sC^{\infty}\mathsf{ring}^{op}$ carrying a manifold $M$ to the discrete simplicial $C^{\infty}$-ring of smooth functions on $M$ is fully faithful, and preserves finite products and transverse pullbacks. 
\end{prop}
\begin{proof}
We may identify the 1-category $\diff$ with the full subcategory of $\rtop^{\loc}(\sring)$ spanned by pairs $(\shv(X),\Of_{X})$ for $X$ a paracompact Hausdorff topological space that admits a cover $\{U_i\subset X\}_i$ such that for each $i$ there is an equivalence $(\shv(X)_{/U_i},\Of_X|_{U_i})\simeq (\shv(\R^n),\cinfty_{\R^n})$. It follows from Theorem \ref{thm:spectrumglobalsections} that the functor $\Gamma$ restricted to this full subcategory is fully faithful. The claim about finite products is proven in Lemma \ref{preservecoproduct}. Consider a transverse pullback diagram
\begin{equation*}
\begin{tikzcd}
Y\times_Z X\ar[r]\ar[d] & X\ar[d]\\    
Y\ar[r] & Z 
\end{tikzcd}
\end{equation*}
of manifolds. Since the pushout $\cinfty(X)\oinfty_{\cinfty(Z)}\cinfty(Z)$ is geometric by Corollary \ref{cor:afpgeometric}, it suffices to show that the diagram 
\begin{equation*}
\begin{tikzcd}
\cinfty(Z)\ar[r]\ar[d] & \cinfty(X)\ar[d]\\    
\cinfty(Y)\ar[r] & \cinfty(Y\times_Z X) 
\end{tikzcd}
\end{equation*}
is a pushout in the \infcat of geometric simplicial $\cinfty$-rings. Invoking Theorem \ref{thm:spectrumglobalsections} again, we conclude that this is a local question on $X$, $Y$ and $Z$, so we may suppose that $Z$ is a Cartesian space. Now the result follows from Lemma \ref{lem:opentransverse}.
\end{proof}
We proved the following with David Carchedi \cite{univprop}.
\begin{thm}\label{thm:geoenv}
The functor $\cinfty(\_):\diff\rightarrow\geodiffder$ is a transformation of pregeometries and exhibits a $\geodiffder$ as a geometric envelope of $\diff$. 
\end{thm}
\begin{proof}
The same proof as the one of Proposition \ref{prop:geoenvopen} applies, using Proposition \ref{manifoldssmoothring} instead of Lemma \ref{lem:opentransverse}.
\end{proof}
\begin{cor}\label{smoothringtruncated}
Let $(\geodiffder)_{\leq n}$ be the opposite category of the $(n+1)$-category of compact objects in $\tau_{\leq n}sC^{\infty}\mathsf{ring}$ for $n\geq 0$. The inclusion $\diff\hookrightarrow (\geodiffder)_{\leq n}$ exhibits $(\geodiffder)_{\leq n}$ as an $n$-truncated geometric envelope of $\diff$. In particular, the inclusion $\diff\hookrightarrow C^{\infty}\mathsf{ring}_{\mathrm{fp}}^{op}$ exhibits $C^{\infty}\mathsf{ring}_{\mathrm{fp}}^{op}$ as a 0-truncated geometric envelope of $\diff$.
\end{cor}
\begin{proof}
Easy consequence of Theorem \ref{thm:geoenv} and Proposition \ref{prop:truncatedgeoenv}.
\end{proof}

\subsubsection{Sheaves of modules}
Let $\ofxtop$ be a simplicial $\cinfty$-ringed \inftopt, then $\Of_{\xtop}$ has an underlying sheaf of commutative $\R$-algebras which we can think of as a commutative algebra object in the symmetric monoidal \infcat of sheaves of $\R$-modules on $\xtop$. Given a map $f:\ofxtop\rightarrow\ofytop$ of simplicial $\cinfty$-ringed \inftopoi with underlying geometric morphism $f_*:\xtop\rightarrow \ytop$, we have functor 
\[\Mod_{\Of_{\xtop}}\longrightarrow\Mod_{f_*(\Of_{\xtop})}\longrightarrow\Mod_{\Of_{\ytop}}.\]
For $\ytop=\spa$, these functor assemble into a global sections functor $\Gamma$ on the \infcat of triples $(\xtop,\Of_{\xtop},\F)$ with $\F$ a sheaf of $\Of_{\xtop}^{\rmalg}$-modules. Our goals in this section are to
\begin{enumerate}[$(1)$]
\item Construct a left adjoint $\mathrm{M}\spec$ to $\Gamma$ and study the properties of this adjunction for well behaved (i.e. Lindel\"{o}f) simplicial $\cinfty$-rings.
\item Endow $\Mod_{\ofxtop}$ with a symmetric monoidal structure that is compatible with small colimits in such a manner that the assignment $\ofxtop\mapsto \Mod_{\ofxtop}^{\otimes}$ determines a sheaf of presentably symmetric monoidal \infcats on the \infcat of locally simplicial $\cinfty$-ringed \inftopoit.
\end{enumerate}
\begin{defn}
Let $\mathrm{MComm}^{\otimes}$ be the 2-coloured \infop governing $\einfty$-algebras and modules over them. Concretely, define a simplicial operad $\mathbf{MComm}^{\otimes}\rightarrow \mathsf{Fin}_*$ as follows.
\begin{enumerate}
\item[$(O)$] Objects of $\mathbf{MComm}^{\otimes}$ are pairs $(\langle n\rangle,T)$ of a pointed finite set $\langle n\rangle \in \mathsf{Fin}_*$ together with a subset $T\subset \langle n\rangle^{\circ}$, where $\langle n\rangle^{\circ}=\langle n\rangle\setminus\{*\}$, the complement of the base point.
\item[$(M)$] A morphism between two pairs $(\langle n\rangle,T)$ and $(\langle m\rangle,S)$ is a map $f:\langle n\rangle\rightarrow\langle m\rangle$ preserving the base point such that
\begin{enumerate}[$(i)$]
\item $f(T\cup \{*\})\subset S\cup\{*\}$.
\item For each $s\in  S$, $f^{-1}(s)$ contains exactly one element of $t$.
\end{enumerate} 
\end{enumerate}
Then we let $\mathrm{MComm}^{\otimes}$ be the operadic nerve of this simplicial operad. Let $\Mod^{\otimes}_{\R}=\Mod_{\R}(\spect)$ be the symmetric monoidal \infcat of $\R$-modules spectra. The forgetful functor $\phi:\Mod_{\R}\rightarrow\spect$ induces an accessible t-structure $(\phi^{-1}(\spect^{\leq0}),\phi^{-1}(\spect^{\geq0}))$ on $\Mod_{\R}$ with the following properties (see \cite{HA}, Sections 7.1.1, 7.1.2 and 7.1.3).
\begin{enumerate}[$(1)$]
\item The t-structure is left and right complete.
\item The t-structure is compatible with filtered colimits. 
\item The t-structure is compatible with the symmetric monoidal structure.
\item The functor $\pi_0:\Mod_{\R}\rightarrow \mathrm{Mod}_{\R}$ determines an equivalence $\Mod^{\heartsuit}_{\R}\simeq \mathrm{Mod}_{\R}$. This equivalence induces a symmetric monoidal t-exact equivalence $\mathbf{D}^{\otimes}(\R)\simeq \Mod^{\otimes}_{\R}$ where $\mathbf{D}^{\otimes}(\R)$ is the derived \infcat of chain complexes of $\R$-modules.
\end{enumerate}
The map of \infops $\mathrm{Comm}^{\otimes}\rightarrow\mathrm{MComm}^{\otimes}$ carrying $\langle n\rangle$ to $(\langle n\rangle,\emptyset)$ determines a functor $\alg_{\mathrm{MComm}}(\Mod_{\R})\rightarrow \calg_{\R}=\alg_{\mathrm{Comm}}(\Mod_{\R})$. We define an \infcat $\Mod$ as the cone in the pullback diagram 
\[
\begin{tikzcd}
\Mod\ar[d,"p_{\Mod}"]\ar[r] & \alg_{\mathrm{MComm}}(\Mod_{\R})\ar[d] \\
\sring\ar[r] & \calg_{\R}
\end{tikzcd}
\] 
among \infcatst, where the lower horizontal functor is the composition $\sring\overset{(\_)^{\rmalg}}{\rightarrow} \scring_{\R}\simeq \mathsf{CAlg}_{\R}^{\geq 0}\subset  \mathsf{CAlg}_{\R}$. Its objects may be identified with pairs $(A,M)$ of a simplicial $\cinfty$-ring $A$ together with an $A^{\rmalg}$-module $M$. For $A$ a simplicial $\cinfty$-ring, we let $\Mod_A$ denote the fibre $p_{\Mod}^{-1}(A)$ and call it the \emph{\infcat of $A$-modules}. Replacing $\Mod_{\R}$ with $\Mod_{\R}^{\geq 0}$, we also have an \infcat $\Mod^{\geq 0}$.
\end{defn}
\begin{prop}
The following hold true.
\begin{enumerate}[$(1)$]
\item The \infcat $\Mod$ is compactly generated. Let $A$ be a simplicial $\cinfty$-ring, then we say that an $A$-module $M$ is \emph{perfect} if $M$ is a finitely presented object of $\Mod_A$. An object $(A,M)$ is compact if and only if $A$ is finitely presented and $M$ is a perfect $A$-module.
\item The \infcat $\Mod^{\geq 0}$ is projectively generated. An object $(A,M)$ is compact projective if and only if $A\simeq\cinfty(\R^n)$ for some $n\geq0$ and $M$ is a finitely generated and free $A$-module. 
\end{enumerate}
\end{prop}
\begin{proof}
For $(2)$, we note that the functor $\Theta:\Mod^{\geq0}\rightarrow\sring\times \Mod_{\R}^{\geq 0}$ is conservative and preserves limits and sifted colimits, so we may invoke \cite{HA}, Proposition 7.1.4.12 and Proposition \ref{prop:lawvereprojgen} to deduce that $\Mod^{\geq 0}$ is projectively generated by the essential image of the composition \[\cartsp^{op}\times \mathrm{Vect}_{\R,\mathrm{fd}}\hooklongrightarrow \sring\times\Mod^{\geq 0}_{\R}\longrightarrow\Mod^{\geq0},\]
where the second functor is left adjoint to $\Theta$, which we can identify with the coproduct of the functors $\sring\rightarrow\Mod^{\geq0}$ carrying $A$ to $(A,0)$ and $\Mod_{\R}^{\geq 0}\rightarrow\Mod^{\geq 0}$ carrying $M$ to $(\R,M)$. It follows that a compact projective object of $\Mod^{\geq 0}$ is a retract of an object of the form $(\cinfty(\R^n),M)$ with $M$ a finitely generated and free $\cinfty(\R^n)$-module. We conclude that the category of Cartesian spaces and smooth maps among them is idempotent complete, as is the category of trivial vector bundles over a Cartesian space.
\end{proof}
\begin{nota}
We let $\mathsf{VBCartSp}$ denote the category of vector bundles on Cartesian spaces, whose nerve we may identify with the opposite of the 1-category of compact projective generators of $\Mod^{\geq 0}$.
\end{nota}

\begin{cons}\label{cons:relativetensor}
Let $\mathrm{CTens}^{\otimes}$ be the category defined as follows.
\begin{enumerate}
\item[$(O)$] Objects are triples $(\langle n\rangle,\langle m\rangle,\{T_i\}_{i\in \langle m\rangle^{\circ}})$ where $\{T_i\}_{i\in M}$ is a collection of pairwise disjoint subsets $T_i\subset \langle n\rangle^{\circ}$ indexed by the elements of $\langle m\rangle^{\circ}$.
\item[$(M)$] Morphisms between triples $(\langle n\rangle,\langle m\rangle,\{T_i\}_{i\in \langle m\rangle^{\circ}})$ and $(\langle k\rangle,\langle l\rangle,\{S_j\}_{j\in \langle l\rangle^{\circ}})$ are pairs of morphisms $f:\langle n\rangle\rightarrow\langle k\rangle$ and $g:\langle m\rangle\rightarrow\langle l\rangle$ of pointed finite sets such that the following conditions are satisfied.
\begin{enumerate}[$(i)$]
\item $f(T_i\cup \{*\})\subset S_{g(j)}\cup \{*\}$ for all $i\in\langle n\rangle^{\circ}$. If $g(j)=*$, then we set $S_{g(j)}=\emptyset$.
\item For each $j\in \langle m\rangle^{\circ}$ and each $s\in S_j$, we have
\begin{equation*}
|f^{-1}(s)\cap T_i|=\begin{cases}
1& \text{if }g(j)=i \\
0&\text{otherwise} 
\end{cases}
\end{equation*}
\end{enumerate}
\end{enumerate}
Consider the functors $p:\mathrm{CTens}^{\otimes}\rightarrow\fin$ and $q:\mathrm{CTens}^{\otimes}\rightarrow \fin$ defined by 
\[p(\langle n\rangle,\langle m\rangle,\{T_i\}_{i\in \langle m\rangle^{\circ}})=\langle n\rangle, \quad q(\langle n\rangle,\langle m\rangle,\{T_i\}_{i\in \langle m\rangle^{\circ}})=\langle m\rangle.\]
Let us notationally distinguish the two copies of $\fin$ that serve as the codomain of the functors above by $(\fin)_p$ and $(\fin)_q$, then the map $p\times q:\mathrm{CTens}^{\otimes}\rightarrow(\fin)_p\times(\fin)_q$ is a $(\fin)_q$-family of \infopst. The fibre over $\langle m\rangle$ of  the projection $q$ is an \infop that we denote by $\mathrm{M}^m\mathrm{Comm}^{\otimes}$; an algebra for this operad can be thought of as commutative algebra $A$ together with $m$ $A$-modules. For $\icat^{\otimes}$ a symmetric monoidal \infcatt, define a simplicial set $\widetilde{\Mod}(\icat)^{\otimes}$ with a map to $\fin$ by the following universal property: for any map of simplicial sets $K\rightarrow\fin$, there is a canonical bijection
\[\Hom_{(\sset)_{/(\fin)_q}}(K,\widetilde{\Mod}(\icat)^{\otimes}) =\Hom_{(\sset)_{/(\fin)_p}}(K\times_{(\fin)_q}\mathrm{CTens}^{\otimes},\icat^{\otimes}). \]
It is not hard to see that the projection $\mathrm{CTens}^{\otimes}\rightarrow(\fin)_q$ is a flat categorical fibration, so the simplicial set $\widetilde{\Mod}(\icat)^{\otimes}$ is in particular an \infcatt. Let $\Mod(\icat)^{\otimes}\subset\widetilde{\Mod}(\icat)^{\otimes}$ be the full subcategory spanned by \infop maps $\mathrm{M}^{m}\mathrm{Comm}^{\otimes}\rightarrow\icat^{\otimes}$, so that objects of $\Mod(\icat)^{\otimes}$ can be identified with pairs $(A,\{M_i\}_{i\in\langle n\rangle^{\circ}})$ with $A$ an $\einfty$-algebra in $\icat$, $\langle n\rangle\in\fin$ and $\{M_i\}_{i\in\langle n\rangle^{\circ}}$ an $n$-tuple of $A$-modules.\\
Let $\mathrm{CTens}^{\otimes}_0\subset\mathrm{CTens}^{\otimes}$ be the full subcategory spanned by triples $(\langle n\rangle,\langle m\rangle,\{\emptyset\}_{i\in\langle m\rangle^{\circ}})$, then the projection $\mathrm{CTens}^{\otimes}_0\rightarrow (\fin)_p\times(\fin)_q$ is an isomorphism. Let $\widetilde{\Mod}(\icat)^{\otimes}_0$ be defined by the existence of a canonical bijection
\[\Hom_{(\sset)_{/(\fin)_q}}(K,\widetilde{\Mod}(\icat)^{\otimes}_0 )=\Hom_{(\sset)_{/(\fin)_p}}(K\times_{(\fin)_q}\mathrm{CTens}_0^{\otimes},\icat^{\otimes}), \]
for each map $K\rightarrow (\fin)_q$, then we have an isomorphism $\widetilde{\Mod}(\icat)^{\otimes}_0 \cong \fun_{(\sset)_{/(\fin)_p}}((\fin)_p,\icat^{\otimes})\times(\fin)_q$ so that the inclusion $\mathrm{CTens}_0^{\otimes}\rightarrow \mathrm{CTens}^{\otimes}$ induces a functor $\Mod(\icat)^{\otimes}\rightarrow\mathsf{CAlg}(\icat)\times\fin$.
\end{cons}
The following result is a straightforward adaptation of some of the material in \cite{HA}, Sections 4.4 and 4.5.
\begin{prop}
Let $\icat^{\otimes}$ be a symmetric monoidal \infcat that admits $\simpop$-indexed colimits separately in each variable, then functor $\Mod(\icat)^{\otimes}\rightarrow \mathsf{CAlg}(\icat)\times\fin$ is a coCartesian $\mathsf{CAlg}(\icat)$-family of symmetric monoidal \infcatst. For each $A\in\mathsf{CAlg}(\icat)$, the symmetric monoidal structure is given by the relative tensor product
\[ \Mod_A\times \Mod_A\longrightarrow \Mod_A,\quad (M,N)\longmapsto M\otimes_AN, \]
and for each map $f:A\rightarrow B$, the coCartesian pushforward $\Mod_A\rightarrow\Mod_B$ is given by the relative tensor product $M\mapsto B\otimes_AM$. The projection $p'_{\Mod}:\Mod(\icat)^{\otimes}\rightarrow\fin$ is a Cartesian fibration and an arrow $(A,\{M_i\}_{i\in\langle n\rangle^{\circ}})\rightarrow (B,\{N_j\}_{j\in \langle m\rangle }^{\circ})$ if $m=n$ and for each $i\in n$, the map underlying map $M_i\rightarrow N_i$ in $\icat$ is an equivalence.
\end{prop}

\begin{defn}\label{defn:tensormodules}
Let $\xtop$ be an \inftopt. Applying Construction \ref{cons:relativetensor} to the symmetric monoidal \infcat $\shv_{\xtop}(\Mod_{\R})^{\otimes}$ yields a coCartesian $\shv_{\xtop}(\mathsf{CAlg}_{\R})$-family of symmetric monoidal \infcatst. For each $\Of_{\xtop}\in \shv_{\xtop}(\mathsf{CAlg}_{\R})$, there is a forgetful functor $\theta:\Mod_{\Of_{\xtop}}\rightarrow \shv_{\xtop}(\Mod_{\R})$ which determins a symmetric monoidal equivalence $\Mod_{\Of_{\xtop}}^{\otimes}\rightarrow\shv_{\xtop}(\Mod_{\R})^{\otimes}$ if and only if $\Of_{\xtop}$ is trivial. We define a coCartesian $\shv_{\sring}(\xtop)$-family of symmetric monoidal \infcats as the cone in the pullback diagram 
\[
\begin{tikzcd}
\shv_{\Mod}(\xtop)^{\otimes}\ar[d] \ar[r]& \Mod(\shv_{\Mod_{\R}}(\xtop))^{\otimes}  \ar[d ] \\
\shv_{\sring}(\xtop)\times\fin \ar[r] & \mathsf{CAlg}(\shv_{\Mod_{\R}}(\xtop))\times\fin
\end{tikzcd}
\]
among \infcatst.
\end{defn}
The following result is proven in Section 2.1.1 of \cite{sag}.
\begin{prop}\label{prop:tstructureshmodproperties}
Let $\xtop$ be an \inftop and let $\Of_{\xtop}$ be a sheaf of simplicial $\cinfty$-rings on $\xtop$.
\begin{enumerate}[$(1)$]
\item The \infcat $\Mod_{\Of_{\xtop}}$ is stable, presentable and the tensor product preserves small colimits separately in each variable.
\item The pair of full subcategories $(\theta^{-1}(\shv_{\Mod_{\R}(\xtop)^{\leq 0}}),\theta^{-1}(\shv_{\Mod_{\R}(\xtop)^{\geq 0}}))$ determines an accessible $\mathrm{t}$-structure on $\Mod_{\Of_{\xtop}}$. This $\mathrm{t}$-structure is right complete, compatible with filtered colimits and compatible with the symmetric monoidal structure.
\item For each algebraic morphism $f^*:\xtop\rightarrow\ytop$, the induced functor $\Mod_{\Of_{\xtop}}\rightarrow \Mod_{f^*(\Of_{\xtop})}$ is $\mathrm{t}$-exact.
\item For each geometric morphism $f_*:\xtop\rightarrow\ytop$, the induced functor $\Mod_{\Of_{\xtop}}\rightarrow \Mod_{f_*(\Of_{\xtop})}$ is right $\mathrm{t}$-exact.
\item Let $\xtop\rightarrow \xtop^{\mathrm{hyp}}$ be an algebraic morphism exhibiting a hypercompletion. Then $\Mod_{\Of_{\xtop}}\rightarrow \Mod_{\Of_{\xtop^{\mathrm{hyp}}}}$ exhibits a left separated quotient.
\item Let $\xtop\rightarrow\widehat{\xtop}$ be an algebraic morphism exhibiting a Postnikov completion. Then $\Mod_{\Of_{\xtop}}\rightarrow \Mod_{\Of_{\widehat{\xtop}}}$ exhibits a left completion.
\end{enumerate}
\end{prop}

The following observations about the abelian hearts of the $\mathrm{t}$-structures constructed above are proven in Section 2.1.2 of \cite{sag}.
\begin{prop}
Let $\xtop$ be an \inftop and let $\Of_{\xtop}$ be a sheaf of simplicial $\cinfty$-rings on $\xtop$. 
\begin{enumerate}[$(1)$]
\item The inclusion $\tau_{\leq 0}\xtop\subset\xtop$ induces a functor $\mathrm{Mod}_{\pi_0(\Of_{\xtop})}\rightarrow \Mod_{\Of_{\xtop}}$ which identifies the former 1-category with the heart $\Mod^{\heartsuit}_{\ofxtop}$.
\item Suppose that $\xtop$ is generated under colimits by subobjects of the final object (that is, $\xtop$ is $(-1)$-complicial). Then $\Of_{\xtop}$ is $0$-truncated if and only if the $\mathrm{t}$-exact functor $\mathbf{D}^+(\mathrm{Mod}_{\pi_0(\Of_{\xtop})})\rightarrow \Mod_{\Of_{\xtop}}$ extending the inclusion $\mathrm{Mod}_{\pi_0(\Of_{\xtop})}\rightarrow \Mod_{\Of_{\xtop}}$ is fully faithful and determines an equivalence $\mathbf{D}^+(\mathrm{Mod}_{\pi_0(\Of_{\xtop})})\simeq \Mod_{\Of_{\xtop}}^+$. If either of these equivalent conditions are satisfied and $\xtop$ is moreover hypercomplete, then the equivalence above extends to an equivalence $\mathbf{D}(\mathrm{Mod}_{\pi_0(\Of_{\xtop})})\simeq \Mod_{\Of_{\xtop}}$. 
\end{enumerate} 
\end{prop}

\begin{cons}
The subcategory inclusion $\rtop\rightarrow\catinfh$ of \inftopoi and geometric morphisms between them determines a biCartesian fibration $\overline{\rtop}^{op}\rightarrow\ltop$. Consider the coCartesian fibration $\Mod^{\otimes}\rightarrow\sring\times\mathsf{Fin}_*$ and define a simplicial set $\widetilde{\ltop(\Mod)}^{\otimes}$ over $\ltop\times\fin$ by the following universal property: for any pair of maps of simplicial sets $K\rightarrow\ltop$ and $K\rightarrow\fin$, there is a canonical bijection between the set $\Hom_{(\sset)_{/\ltop\times\fin}}(K,\widetilde{\ltop(\Mod^{\otimes})})$ and the set of commuting diagrams
\[
\begin{tikzcd}
K\times_{\ltop}\overline{\rtop}^{op} \ar[d]\ar[r] & \Mod^{\otimes}\ar[d] \\
K\ar[r] & \fin.
\end{tikzcd}
\]
Recall the \infcat $\widetilde{\ltop(\sring)}$ defined by the existence of a canonical bijection
\[ \Hom_{(\sset)_{/\ltop}}(K,\widetilde{\ltop(\sring)})\cong \Hom_{\sset}(K\times_{\ltop}\overline{\rtop}^{op},\sring) . \]
It is an immediate consequence of \cite{HTT}, Corollary 3.2.2.12 that the functor $\widetilde{\ltop(\Mod)}^{\otimes}\rightarrow\widetilde{\ltop(\sring)}\times\fin$ is a coCartesian fibration, so that we have a commuting diagram 
\[
\begin{tikzcd}
\widetilde{\ltop(\Mod)^{\otimes}} \ar[dr,"q^{\otimes}_{\mathsf{ModTop}}"']\ar[rr,"p^{\otimes}_{\mathsf{ModTop}}"] &&\widetilde{\ltop(\sring)}\times\fin\ar[dl,"q_{\sring}\times \mathrm{id}"] \\
& \ltop\times\fin
\end{tikzcd}
\]
where all three functors are coCartesian fibrations and $p^{\otimes}_{\mathsf{ModTop}}$ carries $q^{\otimes}_{\Mod}$-coCartesian edges to $q_{\sring}\times\mathrm{id}$-coCartesian edges. Since the functor $p'_{\Mod}:\Mod^{\otimes}\rightarrow\sring$ is also a Cartesian fibration, projecting to $\widetilde{\ltop(\sring)}$ determines a commuting triangle 
\[
\begin{tikzcd}
\widetilde{\ltop(\Mod)^{\otimes}} \ar[dr,"q'_{\mathsf{Mod}}"']\ar[rr,"p'_{\mathsf{ModTop}}"] &&\widetilde{\ltop(\sring)}\ar[dl,"q_{\sring}"] \\
& \ltop
\end{tikzcd}
\]
of biCartesian fibrations, where $p'_{\mathsf{ModTop}}$ carries $q'_{\Mod}$-Cartesian edges to $q_{\sring}$-Cartesian edges.\\
 For each pair $(\xtop,\langle n\rangle)\in \ltop\times\fin$, the fibre of $p^{\otimes}_{\mathsf{ModTop}}$ is identified with the \infcat $\fun(\xtop^{op},\Mod^{\otimes}_{\langle n\rangle})$ and the fibre of $p_{\sring}\times \mathrm{id}$ is identified with \infcat $\fun(\xtop^{op},\sring)\times \fin$. We let $\ltop(\mathsf{Mod})^{\otimes}\subset \widetilde{\ltop(\Mod)^{\otimes}}$ and $\ltop(\sring)\times\fin\subset \widetilde{\ltop(\sring)}\times\fin$ be the full subcategories spanned by functors $\xtop^{op}\rightarrow\Mod^{\otimes}_{\langle n\rangle}$ and $\xtop^{op}\rightarrow \sring$ that preserve small limits, then $p^{\otimes}_{\mathsf{ModTop}}$ restricts to a functor $\ltop(\Mod)^{\otimes}\rightarrow\ltop(\sring)\times\fin$ that we also denote $p^{\otimes}_{\mathsf{ModTop}}$. Projecting to $\ltop$ determines a Cartesian fibration $p'_{\mathsf{ModTop}}:\ltop(\Mod)^{\otimes}\rightarrow\ltop(\sring)$.  
\end{cons}
The following result shows that a map $\ofxtop\rightarrow\ofytop$ in $\ltop(\sring)$ determines a symmetric monoidal functor $\Mod^{\otimes}_{\Of_{\xtop}}\rightarrow\Mod^{\otimes}_{\Of_{\ytop}}$. 
\begin{prop}\label{prop:cocartfamily}
The functor $\ltop(\Mod)^{\otimes}\rightarrow\ltop(\sring)\times\fin$ exhibits a coCartesian $\ltop(\sring)$-family of symmetric monoidal \infcats with fibre $\Mod^{\otimes}_{\xtop}$, that is, for each morphism $\ofxtop\rightarrow\ofytop$ of simplicial $\cinfty$-ringed \inftopoit, the functor $\Mod_{\Of_{\xtop}}\rightarrow\Mod_{\Of_{\ytop}}$ has a canonical structure of a monoidal functor.
\end{prop}
\begin{rmk}
In light of this result, we will write objects of $\ltop(\Mod)^{\otimes}$ as triples $(\xtop,\Of_{\xtop},\{\F_i\}_{i\in\langle n\rangle^{\circ}})$ with $\ofxtop$ a simplicial $\cinfty$-ringed \inftop and $\{\F_i\}_{i\in\langle n\rangle^{\circ}}$ a collection of $n$ sheaves of $\Of_{\xtop}$-modules.
\end{rmk}
The proof of the preceding proposition will require a couple of lemmas, that will be useful further down the line as well.
\begin{lem}\label{lem:ajointablepreserve}
Consider a commuting square of \infcats 
\[
\begin{tikzcd}
\icat\ar[d,"p"]\ar[r,"F"] & \icat'\ar[d,"p'"] \\
\icatd\ar[r,"L"] & \icatd'
\end{tikzcd}
\]
where $p$ and $p'$ are inner fibrations. Suppose that the diagram is horizontally right adjointable, then $F$ caries $p$-coCartesian edges to $p'$-coCartesian edges. 
\end{lem}
\begin{proof}
Viewing $p$ and $p'$ as a natural transformation of functors $\Delta^1\rightarrow\catinf$, taking the relative nerve over $\Delta^1$ yields a diagram 
\[
\begin{tikzcd}
\mathcal{M}_{\icat}\ar[rr,"\overline{p}"] \ar[dr,"\overline{F}"']&&\mathcal{M}_{\icatd}\ar[dl,"\overline{L}"] \\
& \Delta^1
\end{tikzcd}
\]
where $\overline{F}$ and $\overline{L}$ are coCartesian fibrations and $\overline{p}$ preserves coCartesian morphisms. The assumption that the diagram is horizontally right adjointable is equivalent to the assumption that $\overline{F}$ and $\overline{L}$ are Cartesian fibrations and $\overline{p}$ preserves Cartesian morphisms. The assertion to be proven is equivalent to the following.
\begin{enumerate}
\item[$(*)$]
Let $f:X\rightarrow Y$ be a $p$-coCartesian morphism in $\overline{F}^{-1}(0)\simeq\icat$ and choose $\overline{F}$-coCartesian lifts $e_X:X\rightarrow X'$ and $e_Y:Y\rightarrow Y'$ of the nondegenerate 1-simplex of $\Delta^1$ so that we have a commuting square 
\begin{equation}\label{eq:lemadj1}
\begin{tikzcd}
X\ar[r,"f"] \ar[d,"e_X"] & Y\ar[d,"e_Y"] \\
X'\ar[r,"f'"] & Y',
\end{tikzcd}
\end{equation}
then $f'$ is $p'$-coCartesian.
\end{enumerate} 
We wish to show that we can solve lifting problems
\begin{equation}\label{eq:lemadj2}
\begin{tikzcd}
\Delta^{\{0,1\}}\ar[d,hook]\ar[dr,"f'"] & \\
\Lambda^n_0\ar[d,hook]\ar[r,"h'"] &  \icat'\ar[d,"p'"] \\
\Delta^n\ar[r,"g'"] \ar[ur,dotted]& \icatd'.
\end{tikzcd}
\end{equation}
Amalgamating the square \eqref{eq:lemadj1} with the diagram $\Lambda^n_0\overset{h'}{\rightarrow} \icat'\subset \mathcal{M}_{\icat}$ determines a map $\Lambda^n_0\times\{1\}\coprod_{\Delta^{\{0,1\}}\times\{1\}}\Delta^{\{0,1\}}\times\Delta^1\rightarrow\mathcal{M}_{\icat}$. Let $\left(\Delta^n\times\Delta^1\right)^{\sharp}$ denote the marked simplicial set with the edges $\{i\}\times\Delta^1$ for $i>1$ marked (and all degenerate ones), then the inclusion $\left(\Lambda^n_0\times\{1\}\coprod_{\Delta^{\{0,1\}}\times\{1\}}\Delta^{\{0,1\}}\times\Delta^1\right)^{\flat}\subset \left(\Delta^n\times\Delta^1\right)^{\sharp}$ is Cartesian marked anodyne and we can solve the lifting problem
\[
\begin{tikzcd}
\left(\Lambda^n_0\times\{1\}\coprod_{\Delta^{\{0,1\}}\times\{1\}}\Delta^{\{0,1\}}\times\Delta^1\right)^{\flat}\ar[d,hook] \ar[r] & \mathcal{M}_{\icat}\ar[d,"\overline{F}"]\\
\left(\Lambda^n_0\times\Delta^1\right)^{\sharp}\ar[r,"\pi_{\Delta^1}"]\ar[ur,dotted,"H"]& \Delta^1
\end{tikzcd}
\]
so that $H$ carries the edge $\{i\}\times\Delta^1$ to an $\overline{F}$-Cartesian edge for $i>1$. Since $\overline{p}$ preserves Cartesian arrows, the composition $\overline{p}\circ H$ carries the edge $\{i\}\times\Delta^1$ to an $\overline{L}$-Cartesian edge for every $i>1$. Amalgamating $\overline{p}\circ H$ with the map $g':\Delta^n\rightarrow\icatd'\subset\mathcal{M}_{\icatd}$, we have a map $\Delta^n\times\{1\}\coprod_{\Lambda^n_0\times\{1\}}\Lambda^n_0\times\Delta^1\rightarrow \mathcal{M}_{\icatd}$. Let $\left(\Delta^n\times\{1\}\coprod_{\Lambda^n_0\times\{1\}}\Lambda^n_0\times\Delta^1\right)^{\sharp}$ and $\left(\Delta^n\times\Delta^1\right)^{\sharp}$ be the marked simplicial sets with the edges $\{i\}\times\Delta^1$ for $i>1$ marked (and all degenerate ones), then the inclusion $\left(\Delta^n\times\{1\}\coprod_{\Lambda^n_0\times\{1\}}\Lambda^n_0\times\Delta^1\right)^{\sharp}\subset \left(\Delta^n\times\Delta^1\right)^{\sharp}$ is Cartesian marked anodyne and we solve the lifting problem
\[
\begin{tikzcd}
\left(\Delta^n\times\{1\}\coprod_{\Lambda^n_0\times\{1\}}\Lambda^n_0\times\Delta^1\right)^{\sharp}\ar[d,hook] \ar[r] & \mathcal{M}_{\icatd}\ar[d,"\overline{L}"]\\
\left(\Delta^n\times\Delta^1\right)^{\sharp}\ar[r,"\pi_{\Delta^1}"]\ar[ur,dotted,"G"]& \Delta^1
\end{tikzcd}
\]
to find a homotopy $G$. It suffices to solve the lifting problem
\[
\begin{tikzcd}
\Lambda^n_0\times\Delta^1\ar[d,hook]\ar[r,"H"] & \mathcal{M}_{\icat}\ar[d]\\
\Delta^n\times\Delta^1\ar[r,"G"]\ar[ur,dotted,"K"] & \mathcal{M}_{\icatd},
\end{tikzcd}
\]
then $K|_{\Delta^n\times\{1\}}$ solves the original lifting problem. Since $H|_{\Delta^{\{0,1\}}\times\{0\}}=f$ is $p$-coCartesian and $H|_{\Lambda^{n}_0\times\{0\}}$ and $G|_{\Delta^n\times\{0\}}$ take values in $\icat$ and $\icatd$, we can solve the lifting problem 
\[
\begin{tikzcd}
\Delta^{\{0,1\}}\ar[d,hook]\ar[dr,"f"]  \\
\Lambda^n_0\ar[d,hook]\ar[r,"H|_{\Lambda^{n}_0\times\{0\}}"] &[3em]  \icat'\ar[d,"p'"] \\
\Delta^n\ar[r,"G|_{\Delta^n\times\{0\}}"] \ar[ur,dotted]&[3em] \icatd'.
\end{tikzcd}
\]
Amalgamating the solution of this lifting problem with $H$ determines a map $H':\Lambda^n_0\times\Delta^1\coprod_{\Lambda^n_0\times\{0\}}\Delta^n\times\{0\}\rightarrow\mathcal{M}_{\icat}$, so it suffices to solve the lifting problem
\[
\begin{tikzcd}
\Lambda^n_0\times\Delta^1\coprod_{\Lambda^n_0\times\{0\}}\Delta^n\times\{0\}\ar[d,hook]\ar[r,"H'"] & \mathcal{M}_{\icat}\ar[d]\\
\Delta^n\times\Delta^1\ar[r,"G"]\ar[ur,dotted] & \mathcal{M}_{\icatd}.
\end{tikzcd}
\]
 The arrows $H'|_{\{0\}\times\Delta^1}$ and $H'|_{\{1\}\times\Delta^1}$ are the maps $e_X$ and $e_Y$, which are $\overline{F}$-coCartesian. Since $\overline{p}$ preserves coCartesian arrows over $\Delta^1$, the arrows $e_X$ and $e_Y$ are also $\overline{r}$-coCartesian (\cite{HTT}, Proposition 2.4.1.3). Let $\left(\Lambda^n_0\times\Delta^1\coprod_{\Lambda^n_0\times\{0\}}\Delta^n\times\{0\}\right)^{\sharp}$ and $\left(\Delta^n\times\Delta^1\right)^{\sharp}$ be marked simplicial sets with the edges $\{i\}\times\Delta^1$ marked for $i\leq 1$ (and all degenerate ones), then we conclude the proof by observing that the inclusion $\left(\Lambda^n_0\times\Delta^1\coprod_{\Lambda^n_0\times\{0\}}\Delta^n\times\{0\}\right)^{\sharp}\subset\left(\Delta^n\times\Delta^1\right)^{\sharp}$ is coCartesian marked anodyne. 
\end{proof}

\begin{lem}\label{lem:rcocart}
Consider a commuting diagram of simplicial sets
\[
\begin{tikzcd}
\icat\ar[dr,"p"] \ar[rr,"r"] && \icatd\ar[dl,"q"] \\
& \icate
\end{tikzcd}
\]
such that the following conditions are satisfied.
\begin{enumerate}[$(a)$]
\item The maps $p$ and $q$ are biCartesian fibrations and $r$ is an inner fibration. 
\item For each $E\in \icate$, the induced functor $r_E:\icat_E\rightarrow\icatd_E$ on the fibres over $E$ is a coCartesian fibration.
\item The map $r$ carries $p$-Cartesian edges to $q$-Cartesian edges. 
\end{enumerate}
Then the following are equivalent.
\begin{enumerate}[$(1)$]
\item The map $r$ carries $p$-coCartesian edges to $q$-coCartesian edges.
\item The map $r$ is a coCartesian fibration.
\end{enumerate}
\end{lem}
\begin{proof}
Suppose $(1)$ is satisfied. It follows from \cite{HTT}, Proposition 2.4.2.11 that $r$ is a locally coCartesian fibration and that an edge $f:X\rightarrow X'$ of $\icat$ is locally $r$-coCartesian if and only if $f$ admits a factorization $X\overset{e''}{\rightarrow}X''\overset{e'}{\rightarrow} X'$ where $e''$ is $p$-coCartesian and $e'$ is $r_{p(X'')}$-coCartesian, where $r_{E}:\icat_E\rightarrow\icatd_E$ denotes the induced functor on the fibre over $E\in \icate$. It suffices to show that the composition of two locally coCartesian edges is locally coCartesian. Since the composition of two $p$-coCartesian edges or two $r_E$-coCartesian edges for some $E\in\icate$ is $p$-coCartesian respectively $r_E$-coCartesian, it suffices to show that if $f_E:X\rightarrow Y$ is $r_E$-coCartesian and $g:Y\rightarrow Y'$ is $p$-coCartesian, then $g\circ f_E$ is locally coCartesian. Denote $E'=p(Y')$ and choose a $p$-coCartesian lift of $p(g)$ starting at $X$, then we have a commuting diagram
\[
\begin{tikzcd}
X\ar[d,"g_X"]\ar[r,"f_E"] & Y\ar[d,"g_Y"] \\
X'\ar[r,"f_{E'}"] & Y'
\end{tikzcd}
\]
where the vertical maps are $p$-coCartesian. It suffices to show that $f_{E'}$ is $r_{E'}$-coCartesian. Invoking $(a)$, $(c)$ and (the proof of) Lemma \ref{lem:ajointablepreserve}, we conclude. For $(2)\Rightarrow (1)$, we note that an $r$-coCartesian lift of a $q$-coCartesian lift of an edge in $\icate$ is a $p$-coCartesian lift of that edge, by \cite{HTT}, Proposition 2.4.1.3 so the result follows from the essential uniqueness of coCartesian lifts.
\end{proof}

\begin{proof}[Proof of Proposition \ref{prop:cocartfamily}]
It suffices to show the following two assertions
\begin{enumerate}[$(1)$] 
\item For each $\ofxtop\in \ltop(\sring)$, the functor $\ltop(\Mod)^{\otimes}\times_{\ltop(\sring)}\{(\xtop,\Of_{\xtop})\}\rightarrow\fin$ is equivalent to the symmetric monoidal \infcat $\Mod^{\otimes}_{\Of_{\xtop}}\rightarrow \fin$.
\item The functor $\ltop(\Mod)^{\otimes}\rightarrow \ltop(\sring)\times\fin$ is a coCartesian fibration.
\end{enumerate}
Consider the pullback diagram 
\[
\begin{tikzcd}
\fun(\xtop^{op},\Mod)^{\otimes}\ar[d]\ar[r] & \fun(\xtop^{op},\Mod^{\otimes})\ar[d]\\
\fun(\xtop^{op},\sring)\times\fin\ar[r] & \fun(\xtop^{op},\sring\times\fin).
\end{tikzcd}
\]
Restricting on the left hand side to functors $\xtop^{op}\rightarrow\Mod^{\otimes}_{\langle n\rangle}$ and $\xtop^{op}\rightarrow\sring$ that preserve small limits, we obtain the functor $\shv_{\xtop}(\Mod)^{\otimes}\rightarrow\shv(\sring)\times\fin$ of Definition \ref{defn:tensormodules}, which is a coCartesian $\shv_{\xtop}(\sring)$-family of symmetric monoidal \infcatst. This map fits into a pullback diagram 
\[
\begin{tikzcd}
\shv_{\xtop}(\Mod)^{\otimes} \ar[d] \ar[r] & \ltop(\Mod)^{\otimes}\ar[d] \\
\shv_{\xtop}(\sring)\times\fin \ar[r] & \ltop(\sring)\times\fin
\end{tikzcd}
\]
of simplicial sets, so $(1)$ is satisfied. We now show that the functor $p_{\mathsf{ModTop}}^{\otimes}:\ltop(\Mod)^{\otimes}\rightarrow\ltop(\sring)\times\fin$ is a coCartesian fibration. Consider the commuting diagram 
\[
\begin{tikzcd}
\ltop(\Mod)^{\otimes}\ar[rr,"p_{\mathsf{ModTop}}^{\otimes}"] \ar[dr,"q_{\mathsf{Mod}}^{\otimes}"']&& \ltop(\sring)\times\fin\ar[dl,"\pi"]\\
& \ltop
\end{tikzcd}
\]
then $q_{\mathsf{Mod}}^{\otimes}$ is a Cartesian fibration, $\pi$ is a biCartesian fibration and $p^{\otimes}_{\mathsf{ModTop}}$ carries $q^{\otimes}_{\mathsf{Mod}}$-Cartesian edges to $\pi$-Cartesian edges. Invoking Lemma \ref{lem:rcocart}, it suffices to show that $q^{\otimes}_{\mathsf{Mod}}$ is a coCartesian fibration and that $p^{\otimes}_{\mathsf{ModTop}}$-carries $q^{\otimes}_{\mathsf{Mod}}$-coCartesian edges to $\pi$-coCartesian edges. Since $q^{\otimes}_{\mathsf{Mod}}$ is a Cartesian fibration and $p^{\otimes}_{\mathsf{ModTop}}$ preserves Cartesian morphisms over $\ltop$, it suffices to argue that for any algebraic morphism $f^*:\xtop\rightarrow\ytop$ of $\ltop$, the diagram
\[
\begin{tikzcd}
\shv_{\ytop}(\Mod)^{\otimes}\ar[d,"p_{\Mod}^{\ytop\otimes}"]\ar[r,"f_*^{\otimes}"] & \shv_{\xtop}(\Mod)^{\otimes}\ar[d,"p_{\Mod}^{\xtop\otimes}"] \\
\shv_{\ytop}(\sring)\times\fin\ar[r] & \shv_{\xtop}(\sring)\times\fin
\end{tikzcd}
\] 
is horizontally left adjointable, where the horizontal functor are given by composing with the geometric morphism $f_*$. The algebraic morphism $f^*$ determines a symmetric monoidal functor $\shv_{\Mod_{\R}}(\xtop)^{\otimes}\rightarrow\shv_{\Mod_{\R}}(\ytop)^{\otimes}$ which induces a functor $\Mod(\shv_{\Mod_{\R}}(\xtop))^{\otimes}\rightarrow \Mod(\shv_{\Mod_{\R}}(\ytop))^{\otimes}$ which in turn yields a left adjoint to the functor $f_*^{\otimes}$. Moreover, the unit and counit transformations of this adjunction are carried to equivalences by the projections $ \shv_{\ytop}(\Mod)^{\otimes}\rightarrow\fin$ and $\shv_{\xtop}(\Mod)^{\otimes}\rightarrow\fin$. Restrict the diagram above to the base point $\langle *\rangle$, then the vertical functors are equivalences, so we are reduced to proving the diagram is horizontally left adjointable at every triple $(\ofxtop,\{\F_i\}_{i\in\langle n\rangle^{\circ}})$ for $n\geq 1$. Let $\fin^{\geq 1}\subset\fin$ be the full subcategory spanned by all objects except the base point $\langle *\rangle$ and let $\shv_{\ytop}(\Mod)^{\otimes\geq1}$ be the fibre of $p^{\ytop\otimes}_{\Mod}$ over $\shv_{\sring}(\ytop)\times \fin^{\geq 1}$ and define $\shv_{\xtop}(\Mod)^{\otimes\geq1}$ similarly, then it suffices to the diagram
\[
\begin{tikzcd}
\shv_{\ytop}(\Mod)^{\otimes\geq1}\ar[d,"p_{\Mod}^{\ytop\otimes}"]\ar[r,"f_*^{\otimes}"] & \shv_{\xtop}(\Mod)^{\otimes\geq 1}\ar[d,"p_{\Mod}^{\xtop\otimes}"] \\
\shv_{\ytop}(\sring)\times\fin^{\geq 1}\ar[r] & \shv_{\xtop}(\sring)\times\fin^{\geq 1}
\end{tikzcd}
\]
is vertically right adjointable. The vertical functors admit fully faithful right adjoints given by the sections that carry a triple $(\ofxtop,\langle n\rangle)$ to $(\ofxtop,\{0\}_{i\in\langle n\rangle^{\circ}})$ and similarly for $\ytop$, so we conclude by remarking that $f^{\otimes}_*$ preserves zero objects.
\end{proof}
\begin{cor}\label{cor:sheafofmodulecats}
The functor $\ltop(\sring)\rightarrow\prl$ carrying a pair $(\xtop,\Of_{\xtop})$ to $\Mod_{\Of_{\xtop}}$ canonically extends to a functor $\Mod^{\otimes}_{\_}: \ltop(\sring)\rightarrow \mathsf{CAlg}(\prl)$. 
\end{cor}
\begin{rmk}
Since we have for a symmetric monoidal \infcat $\icat^{\otimes}$ and an $\einfty$-algebra object $A$ of $\icat$ a canonical equivalence $\mathsf{CAlg}(\icat)_{A/}\simeq \mathsf{CAlg}(\Mod_A(\icat))$ for $\Mod_A^{\otimes}(\icat)$ the symmetric monoidal \infcat of $A$-module objects in $\icat$, we deduce that the functor $\Mod^{\otimes}_{\_}$ of Corollary \ref{cor:sheafofmodulecats} determines for each $\ofxtop\in\ltop$ a functor $\ltop(\sring)_{\ofxtop/}\rightarrow \mathsf{CAlg}(\Mod_{\Mod_{\Of_{\xtop}}}(\prl))$. Since $(\spa,\R)$ is an initial object of $\ltop(\sring)$, we deduce that $\Mod^{\otimes}_{\_}$ determines a functor
\[ \ltop(\sring)\longrightarrow \mathsf{CAlg}(\mathsf{LinCat}_{\R}) \]
to the \infcat of symmetric monoidal $\R$-linear \infcatst.  
\end{rmk}
Restricting the fibrations $q^{\otimes}_{\Mod}$ and $p^{\otimes}_{\mathsf{ModTop}}$ to $\langle 1\rangle\in \fin$, we obtain a commuting diagram 
\[
\begin{tikzcd}
\ltop(\Mod) \ar[dr,"q_{\mathsf{Mod}}"']\ar[rr,"p_{\mathsf{ModTop}}"] &&\ltop(\sring)\ar[dl,"q_{\sring}"] \\
& \ltop
\end{tikzcd}
\]
of presentable fibrations, by Proposition \ref{prop:cocartfamily}. It follows from Proposition \ref{prop:fibreinclusionadjoint} that the full subcategory inclusion $\Mod\subset \ltop(\Mod)$ of the fibre over $\spa$ admits a right adjoint $\Gamma$. This global sections functor is compatible with the global sections functor for simplicial $\cinfty$-ringed \inftopoit.
\begin{prop}
The strictly commuting diagram 
\[
\begin{tikzcd}
\Mod\ar[d] \ar[r,hook] & \ltop(\Mod)\ar[d]\\
\sring\ar[r,hook] & \ltop(\sring)
\end{tikzcd}
\]
is horizontally right adjointable. In the resulting commuting diagram 
\[
\begin{tikzcd}
\Mod\ar[d,"p_{\Mod}"]  & \ltop(\Mod)\ar[d,"p_{\mathsf{ModTop}}"]\ar[l,"\Gamma"']\\
\sring & \ltop(\sring)\ar[l,"\Gamma"']
\end{tikzcd}
\]
$\Gamma$ carries $p_{\mathsf{ModTop}}$-Cartesian edges to $p_{\Mod}$-Cartesian edges.
\end{prop}
\begin{proof}
Since the horizontal functors are fully faithful, it suffices to argue that the vertical functors carry counit transformations to counit transformations. A map $f:(\xtop,\Of_{\xtop},\F)\rightarrow(\ytop,\Of_{\ytop},\mathcal{G})$ exhibits a counit if and only if $\xtop\simeq \spa$ and $f$ is $p_{\Mod}$-Cartesian, so it suffices to show that $p_{\mathsf{ModTop}}$ preserves Cartesian edges, which is the case by construction. The second assertion follows from Lemma \ref{lem:ajointablepreserve}.
\end{proof}
\begin{defn}
We define a subcategory $\ltop^{\mathrm{loc}}(\Mod)\subset \ltop(\Mod)$ via the pullback diagram
\[
\begin{tikzcd}
\ltop^{\mathrm{loc}}(\Mod) \ar[d] \ar[r] &\ltop(\Mod)\ar[d] \\
\ltop^{\mathrm{loc}}(\sring)\ar[r] & \ltop(\sring).
\end{tikzcd}
\]
\end{defn}
\begin{prop}\label{prop:modulespectrumladj}
The diagram above is horizontally left adjointable. If $(\xtop,\Of_{\xtop},\F)$ is a simplicial $\cinfty$-ringed \inftop equipped with a sheaf of modules, then a map $f:(\xtop,\Of_{\xtop},\F)\rightarrow (\ytop,\Of_{\ytop},\mathcal{G})$ in $\ltop(\Mod)$ is a unit transformation if and only if the following conditions are satisfied.
\begin{enumerate}[$(1)$]
\item $\ofytop$ is a locally simplicial $\cinfty$-ringed \inftopt.
\item The map $\ofxtop\rightarrow\ofytop$ exhibits a relative spectrum.
\item The map $f$ is $p_{\mathsf{ModTop}}$-coCartesian; that is, $f$ induces an equivalence $f^*\F\otimes_{f^*\Of_{\xtop}}\Of_{\ytop}\rightarrow \mathcal{G}$ of $\Of_{\ytop}$-modules. 
\end{enumerate}
\end{prop}
\begin{proof}
This follows immediately from Proposition \ref{prop:pullbackcocartleftadj}.
\end{proof}
\begin{rmk}
Let $\geodiffmod$ be the opposite of the full subcategory of $\Mod$ spanned by compact objects. The \infcat $\geodiffmod$ admits an admissibility structure with a compatible topology defined as follows.
\begin{enumerate}[$(1)$]
\item A map $f:(A,M)\rightarrow (B,N)$ is admissible if the underlying map $A\rightarrow B$ exhibits $B$ as a localization with respect to some $a\in A$ and $f$ is a $p_{\Mod}$-coCartesian morphism.
\item A collection of admissibles $\{(A,M)\rightarrow (B_i,N_i)\}$ determines a covering if and only if the collection $\{A\rightarrow B_i\}$ determines a covering.
\end{enumerate}
Then an object $(\xtop,\Of_{\xtop},\F)\in \ltop(\Mod)$ lies in $\ltop(\geodiffmod)$ if and only if $\ofxtop$ lies in $\ltop^{\mathrm{loc}}(\sring)$ and a map $(\xtop,\Of_{\xtop},\F)\rightarrow (\ytop,\Of_{\ytop},\mathcal{G})$ lies in $\ltop(\geodiffmod)$ if and only if the map $\ofxtop\rightarrow\ofytop$ lies in $\ltop^{\mathrm{loc}}(\sring)$. It follows that we can identify $\ltop(\geodiffmod)$ with the pullback $\ltop^{\mathrm{loc}}(\Mod)$ and the relative spectrum with the left adjoint that appears in the proposition above.
\end{rmk}
According to the preceding propositions, we have a horizontally left adjointable square
\[
\begin{tikzcd}
\mathsf{Top}^{\mathrm{loc},\mathrm{Pc}}(\Mod)^{op}\ar[d]\ar[r,"\Gamma"] &  \Mod\ar[d] \\
\mathsf{Top}^{\mathrm{loc},\mathrm{Pc}}(\sring)^{op}  \ar[r,"\Gamma"]& \sring.
\end{tikzcd}
\]
We let $\mathrm{M}\spec$ denote the left adjoint to the upper horizontal map. For $A$ a simplicial $\cinfty$-ring, we let $\mspec^{\wedge}_A$ denote the induced functor on the fibre $\Mod_A\rightarrow\Mod_{\Of_{\specr^{\wedge}\,A}}$.
\begin{prop}
The following hold true. 
\begin{enumerate}[$(1)$]
\item Let $A$ be a simplicial $\cinfty$-ring, then the functor $\mspec^{\wedge}_A:\Mod_A\rightarrow \Mod_{\Of_{\specr^{\wedge}\,A}}$ can be identified with the composition 
\[ \Mod_A\longrightarrow \Mod_{\underline{A}} \longrightarrow \Mod_{\widetilde{\Of_{\speco\,A}}} \longrightarrow \Mod_{\Of_{\specr^{\wedge}\,A}},  \]
where the first map is induced by the algebraic morphism $i^*:\spa\rightarrow \pshv((\sring^{\mathrm{ad}}_{A/})^{op})$ from the initial \inftop and $\underline{A}=i^*A$ is the constant presheaf on $A$, the second functor is the relative tensor product $\_\otimes_{\underline{A}}\widetilde{\Of_{\speco\,A}}$ and the third functor is induced by the geometric morphism $\pshv((\sring^{\mathrm{ad}}_{/A})^{op})\rightarrow\widehat{\shv}(\specr\,A)$ that sheafifies, takes the spatial reflection and the Postnikov completion.  
\item The functor $\mathrm{M}\spec$ carries $p_{\Mod}$-coCartesian edges to $p_{\mathsf{ModTop}}$-coCartesian edges. In particular, for any map $A\rightarrow B$ of simplicial $\cinfty$-rings, we have homotopy rendering the diagram 
\[
\begin{tikzcd}
\Mod_A\ar[d,"\mspec^{\wedge}_A"]\ar[r,"\_\otimes_AB"] &[8em] \Mod_B\ar[d,"\mspec^{\wedge}_B"] \\
\Mod_{\Of_{\specr^{\wedge}\,A}}\ar[r,"f^*\_\otimes_{f^*\Of_{\specr^{\wedge}\,A}}\Of_{\specr^{\wedge}\,B}"] &[8em] \Mod_{\Of_{\specr^{\wedge}\,B}}
\end{tikzcd}
\]
commutative, where $f^*$ is the algebraic morphism on Postnikov complete \inftopoi associated to the map $\specr\,B\rightarrow\specr\,A$.
\end{enumerate}
\end{prop}
\begin{proof}
The assertion $(1)$ is immediate from the definition of $\mathrm{M}\spec$ and $(2)$ follows from Lemma \ref{lem:ajointablepreserve}.
\end{proof}
\begin{cor}
The spectrum functor $\mathrm{M}\speco^{\wedge}_A:\Mod_A\rightarrow \Mod_{\Of_{\specr^{\wedge}\,A}}$ is $\mathrm{t}$-exact.
\end{cor}
\begin{proof}
Invoking $(3)$ of Proposition \ref{prop:tstructureshmodproperties}, it suffices to show that $\_\otimes_{\underline{A}}\widetilde{\Of_{\speco\,A}}$ is t-exact. For this, it suffices to show that $\underline{A}\rightarrow \widetilde{\Of_{\speco\,A}}$ determines a flat map of commutative algebra objects in $\pshv_{\Mod_{\R}}((\sring^{\mathrm{ad}}_{A/})^{op})$. It suffices to show that for any $a\in A$, the map $\underline{A}(A[a^{-1}])\rightarrow \widetilde{\Of_{\speco\,A}}(A[a^{-1}])$ is a flat map of commutative $\R$-algebras. By construction, the latter map can be identified with the localization $A\rightarrow A[a^{-1}]$, which is flat by Proposition \ref{prop:localization}
\end{proof}
\subsubsection{Geometric modules}
For an arbitrary simplicial $\cinfty$-ring $A$, we cannot say much about the spectrum-global sections adjunction on modules. If we impose that $A$ is Lindel\"{o}f, it turns out there is a good theory of geometric modules which naturally extends the same notion for ordinary $\cinfty$-rings.
\begin{thm}\label{thm:geometricmodulesderived}
Let $(\xtop,\Of_{\xtop})$ be a simplicial $\cinfty$-ringed \inftop and suppose that $\xtop\simeq\widehat{\shv}(X)$ for $X$ a paracompact Hausdorff topological space and $\pi_0(\Of_{\xtop})^{\rmalg}$ is a fine sheaf of commutative $\R$-algebras on $X$.
\begin{enumerate}[$(1)$]
\item The functor $\Gamma:\Mod_{\Of_{X}}\rightarrow\Mod_{\Gamma(\Of_{X})}$ is $\mathrm{t}$-exact.
\end{enumerate}
Let $A$ be a Lindel\"{o}f simplicial $\cinfty$-ring. Then the following hold.
\begin{enumerate}
\item[$(2)$] For every object $\F\in \Mod_{\Of_{\specr^{\wedge}\,A}}$, the counit map $\mspec^{\wedge}_{A}\,\F\rightarrow \F$ is an equivalence. Consequently, the functor $\Gamma:\Mod_{\Of_{\specr^{\wedge}\,A}}\rightarrow\Mod_A$ is fully faithful. 
\item[$(3)$] For any $n\in \Z$, the diagram 
\[
\begin{tikzcd}
\Mod_{\Of_{\specr^{\wedge}\,A}}\ar[d,"\pi_n"] \ar[r,"\Gamma"] & \Mod_A\ar[d,"\pi_n"]\\
\Mod_{\Of_{\specr\,A}}^{\heartsuit} \ar[r,"\Gamma"] & \Mod_A^{\heartsuit}
\end{tikzcd}
\]
which commutes by virtue of $(1)$ (see Lemma \ref{lem:texactness} below), is horizontally left adjointable.
\item[$(4)$] An $A$-module $M$ lies in the essential image of the functor $\Gamma$ if and only if $\pi_n(M)$ is a geometric $\pi_0(A)$-module in the sense of Definition \ref{defn:geometricmodules} for all integers $n$. 
\end{enumerate}
\end{thm}
Before the proof, we make a basic observation about t-structures whose proof we leave to the reader.
\begin{lem}\label{lem:texactness}
Let $f:\icat\rightarrow\icatd$ be an exact functor among stable \infcat equipped with $\mathrm{t}$-structures $(\icat^{\leq 0},\icat^{\geq 0})$ and $(\icatd^{\leq 0},\icatd^{\geq 0})$. Suppose that $f$ is right $\mathrm{t}$-exact. Then $f$ is left $\mathrm{t}$-exact if and only if for each $C\in \icat$, the map $f(C)\rightarrow f(\tau_{\leq 0})$ exhibits a $\tau_{\leq 0}$-localization.
\end{lem}

\begin{proof}[Proof of Theorem \ref{thm:geometricmodulesderived}]
We prove $(1)$. The global sections functor is automatically right t-exact. To see it is t-exact, it suffices to show that for each $\F\in \Mod_{\Of_{\xtop}}$, the map $\Gamma(\F)\rightarrow\Gamma(\tau_{\leq 0}\F)$ exhibits a $\tau_{\leq 0}$-localization, in view of Lemma \ref{lem:texactness}. Write $\F$ as the limit $\F=\lim_{n\in \Z}\tau_{\leq n}\F$ of its left bounded truncations. To show that the fibre of $\theta_n:\Gamma(\F)\rightarrow \Gamma(\tau_{\leq n}\F)$ is $(n+1)$-connective, consider the factorization
\[\Gamma(\F)\overset{\theta_{n+1}}{\longrightarrow}\Gamma(\tau_{\leq (n+1)}\F) \overset{\theta_{n+1,n}}{\longrightarrow}\Gamma(\tau_{\leq n}\F)\]
then the octahedral axiom provides a fibre sequence
\[ \fib(\theta_{n+1})\longrightarrow \fib(\theta_n) \longrightarrow \fib(\theta_{n+1,n}) \]
whose long exact sequence implies that it suffices to show that $\fib(\theta_{n+1})$ and $\fib(\theta_{n+1,n})$ are $(n+1)$-connective. Using Proposition \ref{prop:parahausdorffpresheaf}, we see that $\theta_{n+1,n}$ exhibits an $n$-truncation so that $\fib(\theta_{n+1,n})$ is indeed $(n+1)$-connective. As we have the equivalence $\Gamma(\F)\simeq\lim_{k\geq n+1}\Gamma(\tau_{\leq k}\F)$, we also have an equivalence $\fib(\theta_{n+1})\simeq\lim_{k\geq n+2}\fib(\theta_{k,n+1})$, but by Proposition \ref{prop:parahausdorffpresheaf}, the map $\theta_{k,n+1}$ exhibits an $(n+1)$-truncation so $\fib(\theta_{k,n+1})$ is $(n+2)$-connective for all $k\geq n+2$. Since the limit of a tower of $(n+2)$-connective objects in $\Mod_{\Gamma(\Of_{\xtop})}$ is $(n+1)$-connective, we conclude. We first prove $(3)$. Note that since suspension and looping are equivalences in a stable \infcatt, we only have to check the left adjointability of the diagram for the case $n=0$. The functor $\Gamma:\Mod_{\Of_{\specr^{\wedge}\,A}}\longrightarrow \Mod_A $ admits a factorization as 
\[ \Mod_{\Of_{\specr^{\wedge}\,A}}\longrightarrow \Mod_{\Gamma(\Of_{\specr^{\wedge}\,A})}\longrightarrow \Mod_A \]
where the second functor is induced by the unit map $A\rightarrow \Gamma(\Of_{\specr^{\wedge}\,A})$, which is t-exact, so invoking $(1)$, the composition is t-exact. From Lemma \ref{lem:texactness}, we deduce a commuting diagram
\[
\begin{tikzcd}
\Mod_{\Of_{\specr^{\wedge}\,A}}\ar[d,"\tau_{\geq 0}"] \ar[r,"\Gamma"] & \Mod_A\ar[d,"\tau_{\geq 0}"]\\
\Mod^{,\geq 0}_{\Of_{\specr^{\wedge}\,A}}\ar[d,"\tau_{\leq 0}"] \ar[r,"\Gamma"] & \Mod_A^{\geq 0}\ar[d,"\tau_{\leq 0}"]\\
\Mod_{\Of_{\specr\,A}}^{\heartsuit} \ar[r,"\Gamma"] & \Mod_A^{\heartsuit}
\end{tikzcd}
\]
so it suffices to argue that both squares are horizontally left adjointable. In the lower square, the vertical maps are left adjoint localizations to fully faithful inclusions, so the horizontal left adjointability of the lower square amounts to the assertion that $\Gamma$ carries $0$-truncated sheaves of $\Of_{\specr^{\wedge}\,A}$-modules to $0$-truncated $A$-modules. The upper square is a square of right adjoints, so it is horizontally left adjointable if and only if the associated square
\[
\begin{tikzcd}
\Mod_{\Of_{\specr^{\wedge}\,A}}  &[2em] \Mod_A\ar[l,"\mspec^{\wedge}_A"'] \\
\Mod^{\geq 0}_{\Of_{\specr^\wedge\,A}}\ar[u] &[2em] \Mod_A^{\geq 0}\ar[u]\ar[l,"\mspec^{\wedge}_{A}"']
\end{tikzcd}
\]
of left adjoints is vertically right adjointable. Since the vertical maps are fully faithful, this amounts to the verification that if a map $f:M\rightarrow N$ of $A$-modules exhibits $M$ as a connective cover of $N$, then the map $\mspec^{\wedge}_A\,M\rightarrow\mspec^{\wedge}_A\,N$ exhibits a connective cover. Invoking Lemma \ref{lem:texactness} again, this is simply the assertion that $\mspec^{\wedge}_A$ is t-exact. We now prove $(2)$. It is a consequence of $(1)$ and $(3)$ that the functor $\pi_n$ carries the counit map $g:\mspec^{\wedge}_A\Gamma(\F)\rightarrow\F$ to the counit map $\mspec_{\pi_0(A)}\pi_n(\F)\rightarrow \pi_n(\F)$, which is an equivalence by virtue of Proposition \ref{prop:moduladjunction}. Since the $\mathrm{t}$-structure on $\Mod_{\Of_{\specr^{\wedge}\,A}}$ is left and right complete, we conclude that $g$ is an equivalence. To show $(4)$, we first note that the `only if' direction is immediately supplied by $(1)$. In the other direction, it suffices to show that the unit map $f:M\rightarrow \Gamma(\mspec^{\wedge}_A\,M)$ is an equivalence. It follows from $(3)$ that for every $n\in \Z$, the map $\pi_n(M)\rightarrow \pi_n(\Gamma(\mspec^{\wedge}_A\,M))$ exhibits a geometrization in the sense of Definition \ref{defn:geometricmodules} so that $f$ induces an equivalence on all homotopy groups if the modules $\{\pi_n(M)\}_n$ are already geometric.  
\end{proof}

\begin{cor}\label{cor:thicksubcat}
Let $A$ be a Lindel\"{o}f $\cinfty$-ring, then the category of geometric $A$-modules is a thick abelian subcategory of the category of $A$-modules: it is closed under finite direct sums and if any two objects in a short exact sequence
\[ 0\longrightarrow M\longrightarrow M'\longrightarrow M''\longrightarrow 0 \]
are geometric, then so is the third.
\end{cor}
\begin{proof}
This follows at once from $(2)$ and $(4)$ of Theorem \ref{thm:geometricmodulesderived}, since $\mathrm{t}$-exact fully faithful functor $f:\icat\rightarrow\icatd$ of stable \infcats determines an equivalence of $\icat^{\heartsuit}$ onto a thick subcategory of $\icatd^{\heartsuit}$.
\end{proof}
\begin{defn}\label{defn:geometricmodulesderived}
Let $A$ be a Lindel\"{o}f simplicial $\cinfty$-ring. An $A$-module is \emph{geometric} if $M$ lies in the essential image of the global sections functor. We denote $\Mod_A^{\gmt}\subset\Mod_A$ the full subcategory spanned by geometric $A$-modules. Let $f:M\rightarrow N$ be a morphism of $A$-modules with $N$ geometric, then we say that $f$ \emph{exhibits $N$ as a geometrization of $M$} if for every geometric $A$-module $K$, composition with $f$ induces an equivalence of Kan complexes
\[\Hom_{\Mod_A}(N,K)\overset{\simeq}{\longrightarrow}\Hom_{\Mod_A}(M,K),\]
equivalently, if the map $\mspec^{\wedge}_AM\rightarrow\mspec^{\wedge}_AN$ is an equivalence in the \infcat $\Mod_{\Of_{\specr^{\wedge}\,A}}$.
\end{defn}
\begin{rmk}
As a consequence of Theorem \ref{thm:geometricmodulesderived}, the geometrizations assemble into a functor $L_{\gmt}$ that can be identified with $\Gamma\circ\mspec^{\wedge}_A$. Since the \infcat $\Mod_{A}^{\gmt}$ is presentable, the functor $L_{\gmt}$ is an accessible localization. The t-structure induced by the equivalence $\Mod^{\gmt}\simeq \Mod_{\Of_{\specr^{\wedge}\,A}}$ is left and right complete and coincides with the t-structure $(\Mod^{\gmt}\cap \Mod_A^{\leq0},\Mod_A^{\gmt}\cap\Mod_A^{\geq 0})$. 
\end{rmk}
Here is a simple consequence of Theorem \ref{thm:geometricmodulesderived}.
\begin{lem}\label{lem:globalsecgeometric}
Let $(X,\Of_X)$ be a Postnikov complete simplicial $\cinfty$-ringed space and suppose that $\Gamma(\Of_{X})$ is a Lindel\"{o}f simplicial $\cinfty$-ring. Then for any $\F\in \Mod_{\Of_{X}}$, the global sections $\Gamma(\F)\in \Mod_{\Gamma(\Of_{X})}$ are geometric. 
\end{lem}
\begin{proof}
Since the map $(\spa,\Gamma(\Of_{X}))\rightarrow (X,\Of_{X})$ of simplicial $\cinfty$-ringed spaces factorizes as $\Gamma(\Of_{X})\rightarrow \spec\,\Gamma(\Of_{X})\rightarrow \ofxtop$,
the global sections functor $\Gamma:\Mod_{\Of_{X}}\rightarrow \Mod_{\Gamma(\Of_{X})}$ factorizes as $\Mod_{\Of_{X}}\rightarrow\Mod_{\Of_{\specr^{\wedge}\,\Gamma(\Of_{X})}} \rightarrow \Mod_{\Gamma(\Of_{X})}$. 
\end{proof}
We wish to establish some criteria for the detection of geometricity, starting with finiteness properties. In order to work inductively, we will make the following definition.
\begin{defn}
Let $A$ be a Lindel\"{o}f simplicial $\cinfty$-ring and let $M$ be an $A$-module. Let $n\in \Z$ be an integer, then $M$ is \emph{geometric to order $n$} if for all $k\leq n$, the module $\pi_k(M)$ is a geometric $\pi_0(A)$-module. 
\end{defn}
It follows from Theorem \ref{thm:geometricmodulesderived} that $M$ is geometric in the sense of Definition \ref{defn:geometricmodulesderived} if and only if $M$ is geometric to order $n$ for all integers $n$.
\begin{defn}
Let $A$ be a simplicial $\cinfty$-ring and let $M$ be an $A$-module. Let $n\in \Z$ be an integer, then we say that $M$ is \emph{perfect to order $n$} if $\Hom_{\Mod_A}(M,\_)\rightarrow \spect$ corepresented by $M$ preserves colimits of small filtered diagrams $f:K\rightarrow \Mod_A$ that have the following properties.
\begin{enumerate}[$(a)$]
\item For each $k\in K$, the object $f(k)$ is $n$-truncated.
\item For each edge $k\rightarrow k'$ of $K$, the induced map $\pi_n(f(k))\rightarrow \pi_n(f(k'))$ of $\pi_0(A)$-modules is a monomorphism.
\end{enumerate} 
An $A$ module $M$ is \emph{almost perfect} if $M$ is perfect to order $n$ for all $n\in \Z$.
\end{defn}
\begin{rmk}
Let $A$ be a simplicial $\cinfty$-ring and let $M$ be an $A$-module, then the following hold true.
\begin{enumerate}[$(a)$]
\item $M$ is perfect to order $n$ if and only if $\tau_{\leq n}M$ is perfect to order $n$.
\item If $M$ is perfect to order $n$ and $f:M\rightarrow M'$ induces a surjection $\pi_n(M)\rightarrow\pi_n(M')$ and a bijection $\pi_i(M)\rightarrow \pi_i(M')$ for $i<n$, then $M'$ is perfect to order $n$.
\item If $M$ is perfect to order $n$, then $\tau_{\leq (n-1)}M$ is a compact object in $\Mod^{\leq (n-1)}_A$; in particular, $M$ is eventually connective. Conversely, an $(n-1)$-truncated $A$-module is compact in $\Mod^{\leq (n-1)}_A$ if and only if it is perfect to order $n$. 
\end{enumerate}
\end{rmk}
We now show that perfection implies geometricity.
\begin{prop}\label{prop:perfectisgeom}
Let $A$ be a Lindel\"{o}f simplicial $\cinfty$-ring and $n\geq 0$ an integer. Suppose that $M$ is perfect to order $(n+1)$ and that $A$ is geometric to order $n$, then $M$ is geometric to order $n$.  
\end{prop}
\begin{proof}
Since $M$ is eventually connective, we may suppose that $M$ is connective. We use the following result, which is an easier version of the argument we provided for the existence of $n$-finite cell decompositions of objects finitely generated to order $n$ (see also \cite{dagviii}, Proposition 2.6.12).
\begin{enumerate}
\item[$(*)$] Suppose that $M$ is a connective $A$-module, then 
\begin{enumerate}[$(a)$]
\item $M$ is perfect to order $0$ if and only if $\pi_0(M)$ is a finitely generated $\pi_0(A)$-module. 
\item $M$ is perfect to order $n>0$ if and only if for each map $\phi:A^k\rightarrow M$ of $A$-modules that induces a surjection on connected components (which exists for some finite $k$ by $(a)$), the fibre $\fib(\phi)$ is perfect to order $n-1$.
\end{enumerate}
\end{enumerate}
 We proceed by induction on $n\geq 0$. For $n=0$, we are required to show that $\pi_0(M)$ is a geometric $A$-module. Since $M$ is perfect to order $1$, $\pi_0(M)$ is a finitely presented $\pi_0(A)$-module, so we conclude by invoking Corollary \ref{cor:thicksubcat} and the assumption that $\pi_0(A)$ is geometric. For the inductive step, assume that $M$ is perfect to order $n+1$, then $M$ fits into a fibre sequence
\[ \fib(f)\longrightarrow A^k\overset{f}{\longrightarrow} M  \]
where $k$ is a positive integer chosen so that $f$ is a surjection on connected components and $\fib(f)$ is perfect to order $n$. The result follows from the long exact sequence associated to the fibre sequence above and Corollary \ref{cor:thicksubcat}.
\end{proof}
The following result is a generalization of Corollary \ref{cor:epigeometric}.
\begin{prop}\label{prop:eipigeometricderived}
Let $f:A\rightarrow B$ be an effective epimorphism of Lindel\"{o}f simplicial $\cinfty$-rings. Let $M$ be a $B$-module, then $M$ is geometric to order $n$ as a $B$-module if and only if $M$ is geometric to order $n$ as an $A$-module.
\end{prop}
\begin{proof}
This follows immediately from Corollary \ref{cor:epigeometric} and the t-exactness of the functor $\Mod_B\rightarrow\Mod_A$.    
\end{proof}
Now we can prove Proposition \ref{afpisfair} on the geometricity for finitely presented simplicial $\cinfty$-rings.

\begin{proof}[Proof of Proposition \ref{afpisfair}]
We are given an effective epimorphism $f:A\rightarrow B$ which exhibits $B$ as finitely generated to order $(m+1)$ over $A$, which is geometric to order $m$. Form the pushout diagram 
\[
\begin{tikzcd}
A \ar[d]\ar[r,"f"] & B \ar[d] \\
\tau_{\leq m}A\ar[r,"f'"] & C
\end{tikzcd}
\]
of simplicial $\cinfty$-rings, then $f'$ exhibits $C$ as finitely generated to order $(m+1)$ over $\tau_{\leq m}A$. Since the truncation functors preserve colimits, the right vertical map induces an equivalence $\pi_k(B)\cong\pi_k(C)$ for all $k\leq m$, so we may replace $A$ by $\tau_{\leq m}A$ and $B$ by $C$ and assume that $A$ is geometric. Using Proposition \ref{prop:afpcell}, we may assume that $B$ is an $(m+1)$-finite good $A$-cell object. Proceeding inductively, we are reduced to proving the following special case of the proposition.
\begin{enumerate}
    \item[$(*)$] Let $A$ be a geometric simplicial $\cinfty$-ring and let $V$ be a finite dimensional real vector space. Suppose we are given a morphism $\Sigma^n\cinfty(V)\rightarrow A$ for some $0\leq n$. Consider the pushout diagram
    \[
    \begin{tikzcd}
    \Sigma^n\cinfty(V) \ar[d]\ar[r]& A\ar[d]\\
    \R\ar[r] & B.
    \end{tikzcd}
    \]
    Then $B$ is geometric.
\end{enumerate}
It follows from Proposition \ref{prop:eipigeometricderived} that it suffices to argue that the underlying $A$-module of the simplicial $\cinfty$ ring $B$ is geometric. It follows from Corollary \ref{cor:algeffepi} that the canonical comparison map $\R\otimes_{\Sigma^n\cinfty(V)^{\rmalg}}A^{\rmalg}\rightarrow B^{\rmalg}$ is an isomorphism. We have a convergent homological spectral sequence
\[ E^{p,q}_2= \mathrm{Tor}_p^{\pi_*(\Sigma^n\cinfty(V))}(\R,\pi_*(A))_q \Rightarrow \pi_{p+q}(B),\]
of $\pi_0(A)$-modules. If $n>0$, then Lemma \ref{looping1} shows that the map $\sym^{\bullet}(V^{\vee}[n])\rightarrow\Sigma^n\cinfty(V))$ is an equivalence, and we have an isomorphism $\pi_*(\sym^{\bullet}(V^{\vee}[n]))\cong \sym^{\bullet}(V^{\vee}[n])$ of graded $\R$-algebras. Let $(i)$ denote the shift by $i$ in the $p$-degree, then the factorization \[\sym^{\bullet}(V^{\vee}[n])\longrightarrow \sym^{\bullet}(V^{\vee}[n])\otimes\sym^{\bullet}(V^{\vee}[n](1))\longrightarrow\R\]
exhibits a graded free resolution of $\R$ as a graded $\sym^{\bullet}(V^{\vee}[n])$-module, where $\sym^{\bullet}(V^{\vee}[n])\otimes\sym^{\bullet}(V^{\vee}[n](1))$ is the quasi-free graded differential graded $\R$-algebra on $\mathrm{dim}(V)$ generators in bidegree $(n,0)$ and another $\mathrm{dim}(V)$ generators in bidegree $(n,1)$ whose differential is induced by the identity on $V$.
If $n=0$, Proposition \ref{projresolutiontransverse} implies that the factorization \[\cinfty(V)\longrightarrow\cinfty(V)\otimes\sym^{\bullet}(V^{\vee}(1))\longrightarrow\R\] exhibits a graded free resolution of $\R$ as a $\cinfty(V)$-module. It follows that for all $n\geq 0$ we have isomorphisms 
\[ E_2^{p,q} \cong H_p(\pi_*(A)\otimes \sym^{\bullet}(V^{\vee}[n](1)))_q.\]
Since $V$ is finite dimensional, each term of the bigraded $\pi_0(A)$-module $\pi_{*}(A)\otimes \sym^{\bullet}(V^{\vee}[n](1)))$ is a finite sum of the modules $\{\pi_k(A)\}_{k\geq0}$ and therefore geometric. Since the subcategory spanned by geometric modules is an abelian subcategory, we conclude that the spectral sequence $\{E_r^{p,q}\}_{r\geq 2}$ is a spectral sequence of geometric $\pi_0(A)$-modules. Since the spectral sequence is associated to a connective spectrum, it is concentrated in the first quadrant and the induced filtration on $\pi_q(B)$ is finite and satisfies $F^{k}\pi_q(B)=0$ for $k\leq -1$. We proceed by induction on the length of the filtration. We have $F^0\pi_q(B)=E_{\infty}^{0,q}$, which is geometric. Now suppose that $F^j\pi_q(B)$ is geometric for all $j<p$ and all $q$. We have a short exact sequence
\[ 0\longrightarrow F^{p-1}\pi_{q+p}(B)\longrightarrow F^p\pi_{q+p}(B) \longrightarrow E_{\infty}^{p,q}\longrightarrow 0 \]
for all $q$. Since the full subcategory spanned by geometric modules is stable under extensions, we conclude that $F^p\pi_{q+p}(B)$ is geometric as well. 
\end{proof}

We have given several criteria for the detection of geometric modules. As it turns out, we can characterize precisely the geometric simplicial $\cinfty$-rings that have the property that \emph{all} of their modules are geometric.
\begin{prop}
Let $A$ be a geometric simplicial $\cinfty$-ring such that $\specr\,A$ is locally compact, then the following are equivalent.
\begin{enumerate}[$(1)$]
    \item $\specr\,A$ is compact.
    \item Every $A$-module is geometric.
    \item $\Of_{\specr^{\wedge}\,A}$ is a compact object in $\Mod_{\Of_{\specr^{\wedge}\,A}}$.
\end{enumerate}
\end{prop}
\begin{proof}
For $\icat$ an \infcat admitting small limits and $X$ a locally compact Hausdorff space, we let $\shv^{\mathcal{K}}_{\icat}(X)\subset \pshv_{\icat}(\mathcal{K}(X))$ denote the \infcat of \emph{$\mathcal{K}$-sheaves} on $X$, as defined in \cite{HTT}, Section 7.3.4, where $\mathcal{K}(X)$ is the set of compact subsets of $X$, partially ordered by inclusion. Since the \infcat $\Mod$ is compactly generated, filtered colimits are left exact in $\Mod$, so using \cite{HTT}, Proposition 7.3.4.9, we deduce the existence of a canonical equivalence 
\[ \shv_{\Mod}(\specr A) \simeq \shv^{\mathcal{K}}_{\Mod}(\specr A) \]
such that in case $\specr\,A$ is compact, the global sections functor $\Gamma:\shv_{\Mod}(\specr A)\rightarrow \Mod$ is equivalent to the functor evaluating at $\specr A$. The full subcategory $\shv^{\mathcal{K}}_{\Mod}(\specr A)\subset \pshv_{\Mod}(\mathcal{K}(\specr A))$ is stable under filtered colimits, so we deduce that evaluation of $\mathcal{K}$-sheaves at $\specr A$ preserves filtered colimits. Since the inclusions of the fibres $\Mod_{\specr A}\subset \shv_{\Mod}(\specr A)$ and $\Mod_A\subset\Mod$ preserve and reflect colimits indexed by weakly contractible simplicial sets, we see that the global sections functor $\Gamma:\Mod_{\Of_{\specr\,
 A}}\rightarrow \Mod_{A}$ also preserves filtered colimits. For each $n\in \Z$, the full subcategories $\Mod^{\leq n}_A\subset\Mod_A$ and $\Mod_{\Of_{\specr\,A}}^{\leq n}\subset \Mod_{\Of_{\specr\,A}}$ are stable under filtered colimits and we have an equivalence $\Mod_{\Of_{\specr^{\wedge}\,A}}^{\leq n}\simeq\Mod_{\Of_{\specr\,A}}^{\leq n}$, so we conclude that the full subcategory of $\Mod_A$ spanned by $n$-truncated geometric $A$-module is stable under filtered colimits. Since $\Mod_A^{\leq n}$ is generated under filtered colimits by $n$-truncated $A$-modules perfect to order $n+1$, it follows from Proposition \ref{prop:perfectisgeom} and the geometricity of $A$ that every $n$-truncated $A$-module is geometric. Now $(1)\Rightarrow (2)$ follows from the fact that an $A$-module is geometric if and only if all its truncations are geometric. Now $(2)\Rightarrow (3)$ follows from the fact that $A$ is a compact generator of $\Mod_A$ and $(2)$ of Lemma \ref{lem:geometrization0} below, which asserts that we can identify $\Of_{\specr^{\wedge}\,A}$ with $\mspec^{\wedge}_A\,A$. For the reverse implication, we note that $(3)$ implies that $\mspec^{\wedge}_A$ carries perfect objects to compact objects, which implies that the geometrization is an $\omega$-accessible localization of $\Mod_A$ so that the inclusion $\Mod_A^{\gmt}\subset\Mod_A$ preserves filtered colimits. Since every perfect $A$-module is geometric, we conclude that $\Mod_A^{\gmt}=\Mod_A$. It remains to be shown that $(2)\Rightarrow (1)$. Suppose that $\specr\,A$ is not compact, then we can consider the ideal $I\subset\pi_0(A)$ of those $a\in \pi_0(A)$ that have the property that there exists a compact set $K\subset\specr\,A$ such that $x(a)=0$ for all $x\notin K$. It is not hard to see that $I$ is not a geometric $\pi_0(A)$-module.
\end{proof}
As an immediate consequence we have the following observation.
\begin{cor}
Let $M$ be a manifold, then $M$ is compact (as a topological space) if and only if $\cinfty_M$ is compact (as an object in the category $\mathrm{Mod}_{\cinfty_M}$ of sheaves of $\cinfty_M$-modules).   
\end{cor}
\subsubsection{The proof of Theorem \ref{thm:spectrumglobalsections}}
We turn to the proof of Theorem \ref{thm:spectrumglobalsections}. The main result to be established is that the composition $\Gamma\circ \spec$ is a reflection when restricted to Lindel\"{o}f simplicial $\cinfty$-rings. To prove this, it suffices to show that for all integers $n\geq0$, the map $\pi_n(A)\rightarrow \pi_n(\Gamma(\Of_{\specr^{\wedge}\,A}))$ exhibits a geometrization of $\pi_0(A)$-modules, in case $A$ is Lindel\"{o}f. Making this precise will require a few more auxiliary results.
\begin{cons}
The obvious $\infty$-operad map $\mathrm{MComm}^{\otimes}\rightarrow\mathrm{Comm}^{\otimes}$ carrying $(\langle n\rangle,T)$ to $\langle n\rangle$ determines a section $\Delta_{\icat}:\mathsf{CAlg}(\icat)\rightarrow \alg_{\mathrm{MComm}}(\icat)$ for any symmetric monoidal \infcat $\icat^{\otimes}$, and therefore also a section $\Delta:\sring\rightarrow\Mod^{\geq 0}$ carrying $A$ to the pair $(A,A$), which preserves limits and sifted colimits. Abusing notation, we also denote $\Delta$ the section functor $\ltop(\sring)\rightarrow\ltop(\Mod^{\geq 0})$ defined by unstraightening the natural transformation $\fun^{\mathrm{R}}((\_)^{op},\sring)\rightarrow \fun^{\mathrm{R}}((\_)^{op},\Mod^{\geq0})$ of functors $\ltop^{op}\rightarrow\catinfh$ induced by $\Delta$. In terms of exponentiated fibrations, we can identify $\Delta$ with the composition \[\ltop(\sring)\times_{\ltop}\overline{\rtop}^{op}\longrightarrow \sring\overset{\Delta}{\longrightarrow}\
\Mod^{\geq 0}\]
where the first functor is adjoint to the identity on $\ltop(\sring)$.
\end{cons}
\begin{lem}\label{lem:geometrization0}
The following hold true.
\begin{enumerate}[$(1)$]
\item The commuting diagram 
\[
\begin{tikzcd}
\Mod^{\geq0}\ar[r,hook] & \ltop(\Mod^{\geq0})\\
\sring\ar[u,"\Delta"]\ar[r,hook] & \ltop(\sring),\ar[u,"\Delta"]
\end{tikzcd}
\]
where the horizontal maps are induced by the inclusion $\{\spa\}\hookrightarrow\ltop$ of the initial object, is horizontally right adjointable.
\item Let $\ofxtop$ be a simplicial $\cinfty$-ringed \inftop and let $\ofxtop\rightarrow\ofytop$ be a map in $\ltop(\sring)$ exhibiting a relative spectrum for the geometry $\geodiffder$. Then the map $(\xtop,\Of_{\xtop},\Of_{\xtop})\rightarrow (\ytop,\Of_{\ytop},\Of_{\ytop})$ in $\ltop(\Mod^{\geq0})$ obtained by applying the diagonal functor $\Delta$ to $f$ exhibits a relative spectrum for the transformation $(\geodiffmod)_{\mathrm{disc}}\rightarrow\geodiffmod$. 
\end{enumerate}
\end{lem}
\begin{proof}
To prove $(1)$, we note that the horizontal maps are fully faithful, so it suffices to show that the vertical functors carry counits to counits. It follows from Proposition \ref{prop:fibreinclusionadjoint} that a map $f:\ofxtop\rightarrow\ofytop$ exhibits a counit transformation if and only if $\xtop\simeq\spa$ and $f$ is $p_{\sring}$-Cartesian, and similarly for $\Mod^{\geq0}$-valued sheaves. To prove $(2)$, we note that it follows from Proposition \ref{prop:pullbackcocartleftadj} that a morphism $(\xtop,\Of_{\xtop},\F)\rightarrow (\ytop,\Of_{\ytop},\mathcal{G})$ exhibits a unit transformation for the relative spectrum associated to the transformation $(\geodiffmod)_{\mathrm{disc}}\rightarrow\geodiffmod$ if and only if the underlying map $\ofxtop\rightarrow\ofytop$ exhibits a unit transformation for the relative spectrum $(\geodiffder)_{\mathrm{disc}}\rightarrow\geodiffder$ and $f$ is a coCartesian morphism for the fibration $\ltop(\Mod^{\geq0})\rightarrow\ltop(\sring)$. Thus $(1)$ and $(2)$ will follow from the following assertions: consider the diagram
\[
\begin{tikzcd}
\ltop(\Mod^{\geq 0})\ar[rr,"p_{\mathsf{ModTop}}"] \ar[dr,"q_{\Mod}"'] && \ltop(\sring)\ar[dl,"q_{\sring}"]\\ & \ltop
\end{tikzcd}
\]  
of biCartesian fibrations over $\ltop$. 
\begin{enumerate}[$(1')$]
\item The diagonal $\Delta$ carries $q_{\sring}$-Cartesian edges to $q_{\Mod}$-Cartesian edges.
\item The diagonal $\Delta$ carries any edge of $\ltop(\sring)$ to a $p_{\mathsf{ModTop}}$-coCartesian edge, that is $\Delta$ is a coCartesian section of $p_{\mathsf{ModTop}}$.
\end{enumerate}
The assertion $(1')$ is clear from the definition of $\Delta$. For $(2')$, recall that we can factor any map $f:\ofxtop\rightarrow\ofytop$ in $\ltop(\sring)$ as $\ofxtop\rightarrow (\ytop,f^*\Of_{\ytop})\rightarrow \ofytop$ where the first map is $q_{\sring}$-coCartesian and the second lies in $\shv_{\sring}(\ytop)$. The functor $\Delta|_{\shv_{\sring}(\ytop)}$ coincides with the pullback to $\shv_{\sring}(\ytop)$ of the functor $\Delta_{\shv_{\Mod_{\R}^{\geq0}}(\ytop)}:\mathsf{CAlg}(\shv_{\Mod_{\R}^{\geq 0}}(\ytop))\rightarrow\Mod(\shv_{\Mod_{\R}^{\geq 0}}(\ytop))$; it is straightforward to show from the description of the relative tensor product that $\Delta_{\shv_{\Mod_{\R}^{\geq0}}(\ytop)}$ is a coCartesian section of the fibration $\Mod(\shv_{\Mod_{\R}^{\geq 0}}(\ytop))\rightarrow \mathsf{CAlg}(\shv_{\Mod_{\R}^{\geq 0}}(\ytop))$. We are thus reduced to proving the following.
\begin{enumerate}
\item[$(2'')$] The section $\Delta$ carries $q_{\sring}$-coCartesian edges to $q_{\Mod}$-coCartesian edges.
\end{enumerate}
Let $q^2:\overline{\ltop}^{2}\rightarrow\ltop$ be the presentable fibration defined by the universal property that for any map of simplicial sets $K\rightarrow \ltop$, there is a canonical bijection between the set $\Hom_{(\sset)_{/\ltop}}(K,\overline{\ltop}^{2})$ and the set of maps $K\times_{\ltop}\overline{\rtop}^{op}\rightarrow \spa\times\spa$ that have the property that for any $k\in K$, the induced map $\{k\}\times_{\ltop}\overline{\rtop}^{op}\rightarrow \spa\times\spa$ preserves small limits. Note that the \infcat $\overline{\ltop}^{2}$ is equivalent to a pullback $\overline{\ltop}^{2}\simeq\overline{\ltop}\times_{\ltop}\overline{\ltop}$, with $\overline{\ltop}$ the universal topos fibration, so that $(q^2)^{-1}(\xtop)\simeq \xtop\times\xtop$. Evaluation at the two compact projective generators $(\cinfty(\R),0)$ and $(\R,\R)$ of $\Mod^{\geq 0}$ determines a conservative functor $\Mod^{\geq 0}\rightarrow\spa\times\spa$ which induces a fibrewise (over $\ltop$) conservative functor $\theta:\ltop(\Mod^{\geq0})\rightarrow\overline{\ltop}^2$. Similarly, we have a fibrewise conservative functor $\theta':\ltop(\sring)\rightarrow \overline{\ltop}$ fitting into a commutative diagram 
\[
\begin{tikzcd}
\ltop(\Mod^{\geq0}) \ar[d,"p_{\mathsf{ModTop}}"]\ar[r,"\theta"] &  \overline{\ltop}^2\ar[d,"p"]\\
\ltop(\sring) \ar[r,"\theta'"]&   \overline{\ltop}.
\end{tikzcd}
\]
where $p$ is induced by the first projection $\spa\times\spa\rightarrow\spa$. By construction, $\theta$ carries $q_{\Mod}$-Cartesian edges to $q$-Cartesian edges. The functor $\theta$ also preserves coCartesian edges; this amounts to the assertion that the diagram 
\[
\begin{tikzcd}
\shv_{\Mod^{\geq 0}}(\ytop) \ar[d] \ar[r]& \shv_{\Mod^{\geq 0}}(\xtop)\ar[d] \\
\ytop\times\ytop \ar[r] & \xtop\times\xtop
\end{tikzcd}
\]
is horizontally left adjointable, which is easy to verify. Since the functor $\theta$ is fibrewise conservative, it also detects coCartesian edges. We are reduced to proving that the functor $\Delta:\overline{\ltop}\rightarrow\overline{\ltop}^2$ induced by the diagonal $\spa\rightarrow\spa\times\spa$ carries $q$-coCartesian edges to $q^2$-coCartesian edges. This amounts to the assertion that the diagram
\[
\begin{tikzcd}
\ytop\times \ytop \ar[r,"f_*"] & \xtop\times\xtop\\
\xtop\ar[r,"f_*"]\ar[u,"\Delta"] & \xtop,\ar[u,"\Delta"]
\end{tikzcd}
\] 
is horizontally left adjointable. It suffices to show that the diagram is vertically right adjointable. Since the vertical functors are fully faithful and admit right adjoints given by product functors, this follows from the fact that $f_*$ preserves limits. 
\end{proof}
\begin{cor}\label{cor:geometrization2.5}
Let $(X,\Of_X)$ be a Postnikov complete simplicial $\cinfty$-ringed space such that $\Gamma(\Of_{X})$ is Lindel\"{o}f. Then $\Gamma(\Of_{X})$ is geometric.
\end{cor}
\begin{proof}
Combine Theorem \ref{thm:geometricmodulesderived}, Lemma \ref{lem:globalsecgeometric} and $(1)$ of Lemma \ref{lem:geometrization0}.
\end{proof}
\begin{prop}\label{prop:geometrization1}
Let $\daff_{\gmt_0}\subset\mathsf{Top}^{\mathrm{Pc},\loc}(\sring)$ be the full subcategory spanned by pairs $(X,\Of_{\xtop})$ such that $(X,\pi_0(\Of_{\xtop}))$ is an affine geometric $\cinfty$-scheme. Then the commuting diagram 
\[
\begin{tikzcd}
\sring_{\mathrm{Lin}}\ar[r,"\spec"]\ar[d,"\pi_0"] & \daff^{op}_{\gmt_0}\ar[d,"\pi_0"] \\
\cinfty\mathsf{ring}_{\mathrm{Lin}}\ar[r,"\spec"] & \cinfty\mathsf{Aff}^{op}_{\gmt} 
\end{tikzcd}
\]
is horizontally right adjointable. In particular, if $A$ is a lindel\"{o}f simplicial $\cinfty$-ring, the truncation of the unit map $\tau_{\leq 0}A\rightarrow\tau_{\leq 0}\Gamma(\Of_{\specr^{\wedge}\,A})$ exhibits a geometrization in the sense of Definition \ref{defn:completerings}. 
\end{prop}
\begin{proof}
We first show that $\Gamma:\daff_{\gmt_0}\rightarrow\sring$ takes values in $\sring_{\mathrm{Lin}}$. Let $\ofxtop\in \daff_{\gmt_0}$, then viewing $\ofxtop$ as a $\diff^{\mathrm{open}}$-structure, the truncation $\ofxtop\rightarrow (\xtop,\tau_{\leq 0}\Of_{\xtop})$ with $\tau_{\leq 0}\Of_{\xtop}=\tau_{\leq 0}^{\xtop}\circ \Of_{\xtop}$ exhibits a unit transformation for the adjunction between $\mathsf{Top}^{\loc}(\cinfty\mathsf{ring})$ and $\mathsf{Top}^{\mathrm{Pc},\loc}(\sring)$. Applying the diagonal functor $\Delta$, we obtain a map $\Of_{\xtop}\rightarrow\tau_{\leq 0}\Of_{\xtop}$ of $\Of_{\xtop}$-modules, which exhibits a $0$'th-truncation since $n$-truncations in $\shv_{\xtop}(\Mod^{\geq 0})$ are detected by the functor $\shv_{\xtop}(\Mod^{\geq0})\simeq\fun^{\pi}(\mathsf{VBCartSp},\xtop)\rightarrow\xtop\times\xtop$. Invoking Theorem \ref{thm:geometricmodulesderived}, we deduce that $\Gamma(\Of_{\xtop})\rightarrow\Gamma(\tau_{\leq 0}\Of_{\xtop})$ exhibits a 0'th truncation of $\Gamma(\Of_{\xtop})$-modules. We conclude by $(1)$ of Lemma \ref{lem:geometrization0} that $\Gamma(\Of_{\xtop})\rightarrow\Gamma(\tau_{\leq 0}\Of_{\xtop})$ is also a $0$'th truncation of simplicial $\cinfty$-rings. In particular, $\Gamma(\Of_{\xtop})$ is Lindel\"{o}f.  
Passing to right adjoints in the diagram above, we are required to show that the commuting diagram 
\[
\begin{tikzcd}
\sring_{\mathrm{Lin}} & \daff^{op}_{\gmt_0}\ar[l,"\Gamma"] \\
\cinfty\mathsf{ring}_{\mathrm{Lin}}\ar[u] & \cinfty\mathsf{Aff}_{\gmt}^{op}\ar[l,"\Gamma"] \ar[u]
\end{tikzcd}
\]
is vertically left adjointable, but this follows immediately from the argument we just provided.
\end{proof}

\begin{prop}\label{prop:geometrization2}
Let $A$ be a Lindel\"{o}f simplicial $\cinfty$-ring. Then for each integer $n\geq0$, the map $\pi_n(A)\rightarrow\pi_n(\Gamma(\Of_{\specr^{\wedge}\,A}))$ exhibits a geometrization of $\pi_0(A)$-modules. 
\end{prop}
\begin{proof}
We note that it follows from $(1)$ of Lemma \ref{lem:geometrization0} that taking global sections of the map 
\[(\spa,A,A)\longrightarrow  (\widehat{\shv}(\specr\,A),\Of_{\specr^{\wedge}\,A},\Of_{\specr^{\wedge}\,A})\]
obtained by applying $\Delta$ to the unit map 
\[(\spa,A)\longrightarrow (\widehat{\shv}(\specr\,A),\Of_{\specr^{\wedge}\,A})\]
 yields the map $(A,A)\rightarrow (\Gamma(\Of_{\specr^{\wedge}\,A}),\Gamma(\Of_{\specr^{\wedge}\,A}))$ obtained by applying $\Delta$ to the unit transformation $A\rightarrow\Gamma(\specr^{\wedge}\,A)$. Choose a factorization 
\[ (\spa,A,A)\longrightarrow (\spa,A,M)\longrightarrow (\widehat{\shv}(\specr\,A),\Of_{\specr^{\wedge}\,A},\Of_{\specr^{\wedge}\,A}) \]
where the second map is Cartesian for the fibration $\ltop(\Mod^{\geq0})\rightarrow\ltop(\sring)$, then we may identify $A\rightarrow M$ with the unit transformation $A\rightarrow \mathrm{M}\speco^{\wedge}_A\,A$, in view of $(2)$ of Lemma \ref{lem:geometrization0}. Applying $\Gamma$ to this factorization and using that $\Gamma$ preserves Cartesian edges, we deduce that the map $\Gamma(M)\rightarrow \Gamma(\Of_{\specr^{\wedge}\,A})$ is an equivalence. Now it follows from $(3)$ of Theorem \ref{thm:geometricmodulesderived} that for each $n\geq 0$, the map $\pi_n(A)\rightarrow \pi_n(\Gamma(\mspec^{\wedge}_{A}\,A))$ exhibits a geometrization of $\pi_0(A)$-modules. 
\end{proof}

\begin{prop}\label{prop:geometrization3}
Let $f:\spec\,A\rightarrow \ofxtop$ be a morphism of simplicial $\cinfty$-ringed \inftopoit. Suppose that the algebraic morphism $\widehat{\shv}(\specr\,A)\rightarrow\xtop$ is an equivalence, then $f$ is an equivalence if and only if the map $\Gamma(\Of_{\specr^{\wedge}\,A})\rightarrow \Gamma(\Of_{\xtop})$ is an equivalence.
\end{prop}
\begin{proof}
Let $f_*:\xtop\rightarrow \widehat{\shv}(\specr\,A)$ be the geometric morphism underlying $A$. It suffices to show that the map $\Of_{\specr^{\wedge}\,A}\rightarrow f_*(\Of_{\xtop})$ is an equivalence of $\Of_{\specr^{\wedge}\,A}$-modules. It follows from Theorem \ref{thm:geometricmodulesderived} that the global sections functor $\Gamma:\Mod_{\Of_{\specr^{\wedge}\,A}}\rightarrow \Mod_A$ is fully faithful, so it suffices to show that $\Gamma(\Of_{\specr^{\wedge}\,A})\rightarrow\Gamma(f_*(\Of_{\xtop}))$ is an equivalence. Using $(1)$ of Lemma \ref{lem:geometrization0} and the fact that $\Gamma:\ltop(\Mod)\rightarrow\Mod$ preserves Cartesian edges, we conclude that the map $\Gamma(\Of_{\specr^{\wedge}\,A})\rightarrow\Gamma(f_*(\Of_{\xtop}))$ is the result of applying $\Delta$ to the map $\Gamma(\Of_{\specr^{\wedge}\,A})\rightarrow\Gamma(\Of_{\xtop})$ of simplicial $\cinfty$-rings. 
\end{proof}
\begin{proof}[Proof of Theorem \ref{thm:spectrumglobalsections}]
To prove $(1)$, we first show that the image of $\Gamma\circ \spec$ is the full subcategory spanned by geometric simplicial $\cinfty$-rings. It follows from Proposition \ref{prop:geometrization1} that if $A$ is Lindel\"{o}f, then $\pi_0(\Gamma(\Of_{\specr^{\wedge}\,A}))$ is geometric and thus Lindel\"{o}f, so it follows from Corollary \ref{cor:geometrization2.5} that $\Gamma(\Of_{\specr^{\wedge}\,A})$ is geometric. If $A$ is already geometric, then Proposition \ref{prop:geometrization2} implies that the unit map $A\rightarrow\Gamma(\Of_{\specr^{\wedge}\,A})$ is an equivalence. To show that $\Gamma\circ\spec$ is a localization, it suffices to show that the indicated maps $f$ and $g$ in the commuting diagram 
\[
\begin{tikzcd}
A\ar[d]\ar[r] & \Gamma(\Of_{\specr^{\wedge}\,A})\ar[d,"f"] \\
\Gamma(\Of_{\specr^{\wedge}\,A})\ar[r,"g"] & \Gamma(\Of_{\specr^{\wedge}\,\Gamma(\Of_{\specr^{\wedge}\,A})})
\end{tikzcd}
\]
are equivalences by \cite{HTT}, Proposition 5.2.7.4, where the vertical arrows are unit transformations. Since $\Gamma(\Of_{\specr^{\wedge}\,A})$ is geometric as verified above, applying Proposition \ref{prop:geometrization2} again guarantees that $f$ is an equivalence. It follows that in the composition 
\[  \pi_n(A)\longrightarrow \pi_n(\Gamma(\Of_{\specr^{\wedge}\,A}))\overset{\pi_n(g)}{\longrightarrow} \pi_n(\Gamma(\Of_{\specr^{\wedge}\,\Gamma(\Of_{\specr\,A}^{\wedge})})) \]
of $\pi_0(A)$-modules, the first map and the composition exhibit a geometrization of $\pi_0(A)$-modules for every integer $n\geq 0$. It follows that the second map is an isomorphism so that $g$ is an equivalence as well.\\
We now shows $(2)$. In view of Proposition \ref{prop:geometrization3} and $(1)$, it suffices to show that the map $\widehat{\shv}(\specr\,A)\rightarrow \widehat{\shv}(\Gamma(\Of_{\specr^{\wedge}\,A}))$ is an equivalence of \inftopoit. This follows from Proposition \ref{prop:geometrization1} and Proposition \ref{prop:joyce1}. For $(3)$, we show that the counit map $\spec\,\Gamma(\Of_{\xtop})\rightarrow \ofxtop$ is an equivalence in case $\ofxtop$ satisfies $(a)$ and $(b)$. In view of Proposition \ref{prop:geometrization1} and Theorem \ref{thm:joyce}, the underlying map of \inftopoi is an equivalence. Invoking Propostion \ref{prop:geometrization3}, it suffices to show that the map $\Gamma(\Of_{\specr^{\wedge}\,\Gamma(\Of_{\xtop})})\rightarrow \Gamma(\Of_{\xtop})$ is an equivalence. By the unit-counit identities this is equivalent to proving that the unit map $\Gamma(\Of_{\xtop})\rightarrow\Gamma(\Of_{\specr^{\wedge}\,\Gamma(\Of_{\xtop})})$ is an equivalence. By virtue of $(1)$, it suffices to show that $\Gamma(\Of_{\xtop})$ is geometric, which is asserted by Corollary \ref{cor:geometrization2.5}.
\end{proof}

\subsection{Differential graded models: $C^{\infty}$dgas and dg-manifolds}
The goal of this subsection is to provide homological algebraic models for simplicial $\cinfty$-rings.
\begin{defn}
Let $\mathbf{cdga}_{\R}^{\geq 0}$ be the category of nonnegatively graded differentially graded algebras over $\R$ (grading conventions are homological). It comes with a canonical projection $\mathbf{cdga}_{\R}^{\geq 0}\rightarrow \mathsf{CAlg}^0_{\R}$ by restricting to degree 0, which is right adjoint to the obvious inclusion $\mathsf{CAlg}^0_{\R}\rightarrow \mathbf{cdga}^{\geq 0}_{\R}$. The category of \emph{$C^{\infty}$dgas}, denoted $C^{\infty}\mathbf{dga}^{\geq0}$, is the pullback $\mathbf{cdga}_{\R}^{\geq 0}\times_{\mathsf{CAlg}^0_{\R}}C^{\infty}\mathsf{ring}$. Concretely, a \emph{$C^{\infty}$dga} is a nonnegatively graded differentially graded algebra $A_{\bullet}$, such that $A_0$ has the structure of a $C^{\infty}$-ring compatible with its $\R$-algebra structure. A morphism of $C^{\infty}$dgas is a homomorphism of nonnegatively graded dg algebras that restricts to a morphism of $C^{\infty}$-rings in degree 0.
\end{defn}

\begin{rmk}
Occasionally, we will have to work with $C^{\infty}$dgas whose underlying chain complex is not connective. The obvious inclusion $\mathsf{CAlg}^0_{\R}\rightarrow\mathbf{cdga}_{\R}$ has a right adjoint $\mathbf{cdga}_{\R}\rightarrow \mathsf{CAlg}^0_{\R}$ which takes $A_{\bullet}$ to $\ker(\del_0)\subset A_0$. The category of \emph{unbounded $C^{\infty}$dgas}, denoted $C^{\infty}\mathbf{dga}$, is the pullback $\mathbf{cdga}_{\R}\times_{\mathsf{CAlg}_{\R}^0}C^{\infty}\mathsf{ring}$ whose objects are cdgas $A_{\bullet}$ for which $\ker(\del_0)$ is a $C^{\infty}$-ring.
\end{rmk}
\begin{rmk}\label{rmk:cinftycompactgen}
Consider the commuting diagram 
\[
\begin{tikzcd}
\mathbf{cdga}_{\R}^{\geq 0} \ar[d]\ar[r] & \mathbf{Mod}_{\R}^{\geq 0} \ar[d] \\
\calg^0_{\R} \ar[r] & \mathrm{Vect}_{\R}
\end{tikzcd}
\]
of 1-categories. By \cite{HA}, Corollary 3.2.3.2 the horizontal maps preserve and detect limits and sifted colimits and right vertical map preserves all colimits so the left vertical maps preserves sifted colimits as well. It follows from \cite{HTT}, Proposition 5.4.5.5 that the forgetful functors $\cinfty\mathbf{dga}^{\geq 0}\rightarrow \mathbf{cdga}_{\R}^{\geq 0}$ and $\cinfty\mathbf{dga}^{\geq 0}\rightarrow \calg_{\R}^0$ preserve sifted colimits. For unbounded $\cinfty$dgas, the same observation holds for filtered colimits, since $\ker(\del_0)$ preserves filtered colimits. It follows that $\cinfty\mathbf{dga}$ and $\cinfty\mathbf{dga}^{\geq 0}$ are compactly generated presentable 1-categories.
\end{rmk}
The proof of the following proposition is a straightforward verification.
\begin{prop}\label{prop:degzeroadj}
Let $A$ be a commutative $\R$-algebra. Let $A^0_{\bullet}$ be the nonnegatively graded cdga given by
\[  \ldots\longleftarrow 0\longleftarrow A \longleftarrow 0\longleftarrow \ldots, \]
where $A$ sits in degree $0$. Let $A_{\bullet}^{0,1}$ be the cdga given by 
\[  \ldots\longleftarrow 0\longleftarrow A \overset{\mathrm{id}}{\longleftarrow} A\longleftarrow 0\longleftarrow \ldots, \] 
where $A$ sits in degrees $0$ and $1$ and the underlying graded $\R$-algebra structure is that of the square zero extension of $A$ by the graded $A$-module $A[1]$, that is, $A$ shifted to degree $1$ (one immediately verifies that with this graded algebra structure, the identity $A\rightarrow A$ is a graded derivation). Then the following hold true.
\begin{enumerate}[$(1)$]
\item The pair $(A^0_{\bullet},A\overset{\mathrm{id}}{\rightarrow}A)$ exhibits a unit transformation at $A$ for the functors $(\_)_0$ and $\ker(\del_0)$.
\item The pair $(A^{0,1}_{\bullet},A\overset{\mathrm{id}}{\rightarrow}A)$ exhibits a counit transformation at $A$ for the functor $(\_)_0$.  
\end{enumerate}
Consequently, the functors $(\_)_0$ and $\ker(\del_0)$ preserve all limits and $(\_)_0$ preserves all colimits.
\end{prop}
We see that the category $\cinfty\mathbf{dga}^{\geq 0}$ arises as a pullback along a functor $(\_)_0:\mathbf{cdga}_{\R}^{\geq 0}\rightarrow\mathsf{CAlg}^0_{\R}$, which preserves limits and colimits. As it turns out, this property of the functor $(\_)_0$ gives us efficient means to compute colimits and free algebras in $\cinfty\mathbf{dga}^{\geq 0}$ in terms of the underlying cdga and $\cinfty$-ring in degree $0$. We will encounter this situation -that is, `algebraic' objects that are endowed with some sort of $\cinfty$-enhancement in low degrees- more often in this work and its successors, so we devote some time to formal considerations.
\begin{prop}\label{prop:leftadjseccocart}
Let $p:\icat\rightarrow\icatd$ be a categorical fibration and suppose that $p$ admits a left adjoint $f$. Let $\overline{e}:C\rightarrow C'$ be a map in $\icat$, then $\overline{e}$ is $p$-coCartesian if and only if the diagram 
\[
\begin{tikzcd}
fp(C)\ar[r] \ar[d,"\epsilon"] & fp(C') \ar[d,"\epsilon"]\\ 
C\ar[r,"\overline{e}"] & C'
\end{tikzcd}
\]
where both vertical arrows are counit transformations at $C$ and $C'$ respectively is a pushout. Moreover, for $C\in \icat$ and $e:p(C)\rightarrow D$ a map in $\icatd$, the following are equivalent.
\begin{enumerate}[$(1)$]
\item There exists a $p$-coCartesian lift $\overline{e}:C\rightarrow\widetilde{D}$ of $e$ starting at $C$.
\item There exists pushout diagram
\[
\begin{tikzcd}
fp(C)\ar[d,"\epsilon"]\ar[r,"f(e)"] & f(D)\ar[d] \\
C \ar[r] & \widetilde{D}.
\end{tikzcd}
\] 
where $\epsilon$ is the counit of $(f\adj p)$ at $C$ and the canonical map $D\overset{\eta}{\rightarrow} pf(D)\rightarrow p(\widetilde{D})$ is an equivalence, where the first map is the unit of $(f\adj p)$ at $D$.
\end{enumerate}
\end{prop}
\begin{proof}
The map $\overline{e}:C\rightarrow C'$ is $p$-coCartesian precisely if for every $C''\in\icat$ composition with $\overline{e}$ determines a pullback diagram 
\[
\begin{tikzcd}
\Hom_{\icat}(C',C'') \ar[d] \ar[r] &\Hom_{\icat}(C,C'')\ar[d]\\
\Hom_{\icatd}(p(C'),p(C''))  \ar[r] &\Hom_{\icatd}(p(C),p(C''))
\end{tikzcd}
\]
of spaces, but we have a commuting diagram 
\[
\begin{tikzcd}
\Hom_{\icat}(C',C'') \ar[d] \ar[r] &\Hom_{\icat}(C,C'')\ar[d]\\
\Hom_{\icat}(fp(C'),C'')  \ar[r]\ar[d,"\simeq"] &\Hom_{\icat}(fp(C),C'')\ar[d,"\simeq"]\\
\Hom_{\icatd}(p(C'),p(C''))  \ar[r] &\Hom_{\icatd}(p(C),p(C'')) 
\end{tikzcd}
\]
where the upper vertical maps are induced by the counit transformation. Thus $\overline{e}$ being $p$-coCartesian is equivalent to the upper square being a pullback, that is, the diagram 
\[
\begin{tikzcd}
fp(C)\ar[r] \ar[d] & fp(C') \ar[d]\\ 
C\ar[r,"\overline{e}"] & C'
\end{tikzcd}
\]
being a pushout. The implication $(1)\Rightarrow (2)$ follows easily from the first part of the proof and the unit-counit identities. Suppose that $(2)$ is satisfied, then we have a map $C\rightarrow\widetilde{D}$ lying over $p(C)\rightarrow p(\widetilde{D})$ fitting into the pushout diagram 
\[
\begin{tikzcd}
fp(C)\ar[d,"\epsilon"]\ar[r,"f(e)"] & f(D)\ar[d,"\zeta"] \\
C \ar[r] & \widetilde{D}.
\end{tikzcd}
\] 
The right vertical map $\zeta:f(D)\rightarrow\widetilde{D}$ factors as $f(D)\rightarrow fp(\widetilde{D}) \overset{\epsilon}{\longrightarrow} \widetilde{D}$
where $\epsilon$ is the counit of $(f\adj p)$ at $\widetilde{D}$ and the first map is the image under $f$ of the composition $D\overset{\eta}{\longrightarrow} pf(D)\overset{p(\zeta)}{\longrightarrow}p(\widetilde{D})$, where $\eta$ is the unit at $D$. By assumption, this composition is an equivalence so the diagram 
\[
\begin{tikzcd}
fp(C)\ar[d,"\epsilon"]\ar[r,"f(e)"] & f(D)\ar[r]& fp(\widetilde{D})  \ar[d,"\epsilon"] \\
C \ar[rr]& & \widetilde{D}.
\end{tikzcd}
\] 
is a pushout. By the first part of the proof, the map $C\rightarrow \widetilde{D}$ is $p$-coCartesian. Now consider the diagram
\[
\begin{tikzcd}
p(C)\ar[r,"e"]\ar[d,"\eta"] & D \ar[d,"\eta"] \\
pfp(C)\ar[d,"p(\epsilon)"]\ar[r,"pf(e)"] & pf(D)\ar[r]& pfp(\widetilde{D})  \ar[d,"p(\epsilon)"] \\
p(C) \ar[rr]& & p(\widetilde{D}).
\end{tikzcd}
\] 
The composition $p(C)\overset{\eta}{\rightarrow}pfp(C)\overset{p(\epsilon)}{\rightarrow} p(C)$ is an equivalence and by assumption, the composition $D\overset{\eta}{\rightarrow}pf(D)\overset{p(\zeta)}{\rightarrow}p(\widetilde{D})$ is an equivalence. We conclude that $e$ is equivalent to the lower horizontal map $p(C)\rightarrow p(\widetilde{D})$. Since this map has a $p$-coCartesian lift starting at $C$, so too does $e$, since $p$ is a categorical fibration.
\end{proof}
\begin{prop}\label{prop:colimleftadjcocart}
Let $p:\icat\rightarrow\icatd$ be a coCartesian fibration of \infcats and $g:\icatd'\rightarrow\icatd$ a functor admitting a left adjoint $f$. Let $\icat'$ be defined as the cone in the pullback diagram 
\[
\begin{tikzcd}
\icat'\ar[d,"p'"]\ar[r,"h"] & \icat\ar[d,"p"] \\
\icatd'\ar[r,"g"] & \icatd
\end{tikzcd}
\]
of simplicial sets. Let $K$ be a simplicial set and assume that $\icat$ and $\icatd'$ admit $K$-indexed colimits and that $\icatd'$ admits pushouts. Then $\icat'$ admits $K$-indexed colimits. Moreover, a diagram $\overline{\theta}:K^{\rhd}\rightarrow\icat'$ is a colimit diagram if and only if the following hold.
\begin{enumerate}[$(1)$]
\item Let $\theta=\overline{\theta}|_K$ and let $\overline{h\theta}:K^{\rhd}\rightarrow\icat'$ be a colimit diagram extending $h\theta$. Then the canonical comparison map $\overline{h\theta}(\infty)\rightarrow h\overline{\theta}(\infty)$ (defined up to contractible ambiguity) is $p$-coCartesian.
\item Since the object $ph\overline{\theta}(\infty)=gp'\overline{\theta}(\infty)$ lies in the image of $g$, the map $p\overline{h\theta}(\infty)\rightarrow ph\overline{\theta}(\infty)$ factors as
\[ p\overline{h\theta}(\infty) \longrightarrow gfp\overline{h\theta}(\infty)\longrightarrow p\overline{h\theta}(\infty)  \]
where the first map is the unit of $(f\adj g)$ and the second is the image of a map $e:fp\overline{h\theta}(\infty)\rightarrow p'\overline{\theta}(\infty)$ in $\icatd'$ under the diagram $fgp'\theta$. Let $\overline{p'\theta},\overline{fgp'\theta}:K^{\rhd}\rightarrow \icatd'$ be colimit diagrams extending $p'\theta$ and $fgp'\theta$ respectively, then $e$ fits into the canonical commuting diagram
\[
\begin{tikzcd}
\overline{fgp'\theta}(\infty) \ar[d]\ar[r] &fp\overline{h\theta}(\infty) \ar[d,"e"] \\ 
\overline{p'\theta}(\infty) \ar[r] & p'\overline{\theta}(\infty).
\end{tikzcd}
\]
Then this diagram is a pushout.
\end{enumerate}
Moreover, if $p$ preserves $K$-indexed colimits, then $(2)$ can be replaced by the following condition.
\begin{enumerate}
\item[$(2')$] The lower horizontal map $\overline{p'\theta}(\infty)\rightarrow p'\overline{\theta}(\infty)$ in the diagram above is an equivalence.
\end{enumerate}
\end{prop}
\begin{proof}
We show that for a functor $\theta:K\rightarrow \icat'$, there is, given our assumptions on $\icat$, $
\icatd'$ and $p$ and $g$, an extension $\overline{\theta}:K^{\rhd}\rightarrow \icat'$ satisfying $(1)$ and $(2)$ (or $(2')$) that is a colimit diagram. This will prove the proposition: since conditions $(1)$ and $(2)$ (or $(2')$) are stable under equivalence in the \infcat $\icat'_{K/}$, any other extension of $\theta$ that is a colimit diagram will also satisfy $(1)$ and $(2)$ (or $(2')$); conversely, for another extension $\widehat{\theta}$ satisfying $(1)$ and $(2)$ (or $(2')$) there is an essentially unique map $\overline{\theta}(\infty)\rightarrow \widehat{\theta}(\infty)$ because $\widehat{\theta}$ is a colimit diagram which is an equivalence, since $(1)$ and $(2)$ (or $(2')$) determine the images of $\overline{\theta}(\infty)$ and $\widehat{\theta}(\infty)$ in $\icat$, $\icatd$ and $\icatd'$ up to equivalence and $p$ is a categorical fibration. We have a pullback diagram of simplicial sets
\[
\begin{tikzcd}
\icat'_{\theta/} \ar[r] \ar[d] & \icat_{h\theta/}\ar[d] \\
\icatd'_{p'\theta/} \ar[r] & \icatd_{ph\theta/}.
\end{tikzcd}
\]   
According to \cite{HTT}, Proposition 2.4.3.1, the right vertical map in this diagram is a coCartesian fibration. By our assumptions, the conditions of \cite{HTT}, Lemma 5.2.5.2 are satisfied which guarantees that the lower horizontal map in the diagram above admits a left adjoint. It follows from Proposition \ref{prop:pullbackcocartleftadj} that the diagram is horizontally left adjointable. By assumption, the \infcat $\icat_{h\theta/}$ admits an initial object so $\icat'_{\theta/}$ admits an initial object as well, that is, $\theta$ admits a colimit. It follows from Proposition \ref{prop:pullbackcocartleftadj} that the canonical comparison map is $p$-coCartesian, so that $(1)$ is satisfied. Condition $(2)$ follows immediately from the description of the lower horizontal left adjoint provided by \cite{HTT}, Lemma 5.2.5.2. In case $p$ preserves $K$-indexed colimits, the map $\overline{fgp'\theta}(\infty) \rightarrow fp\overline{h\theta}(\infty)$ is an equivalence, so that the map $\overline{p'\theta}(\infty)\rightarrow p'\overline{\theta}(\infty)$ is an equivalence as well.   
\end{proof}

\subsubsection{The model structure on $\cinfty$dgas}
We deduce some useful tools for dealing with $\cinfty$dgas based on some results from the next section, in particular the flatness of the unit map $A\rightarrow F^{\cinfty}(A)^{\rmalg}$ (Proposition \ref{prop:freecinftycommutestruncation}) and apply them to construct model structures on nonnegatively graded and unbounded $\cinfty$dgas.
\begin{prop}\label{prop:freecinftydgacolim}
Let $(\_)^{\rmalg}_{\mathbf{dg}}:\cinfty\mathbf{dga}\rightarrow \mathbf{cdga}_{\R}$ be the functor forgetting the $\cinfty$-ring structure on the $0$-cycles. We also denote $(\_)^{\rmalg}_{\mathbf{dg}}$ for the functor $\cinfty\mathbf{dga}^{\geq 0}\rightarrow \mathbf{cdga}^{\geq 0}_{\R}$ forgetting the $\cinfty$-ring structure in degree $0$. The following hold true.
\begin{enumerate}[$(1)$]
\item The functor $(\_)^{\rmalg}_{\mathbf{dg}}$ admits a left adjoint $F^{\cinfty}_{\mathbf{dg}}$ and the pullback diagram 
\[  
\begin{tikzcd}
\cinfty\mathbf{dga} \ar[d] \ar[r,"(\_)^{\rmalg}_{\mathbf{dg}}"] & \mathbf{cdga}_{\R}\ar[d] \\
\cinfty\mathsf{ring}\ar[r,"(\_)^{\rmalg}"] & \calg^0_{\R}
\end{tikzcd}
\]
is horizontally left adjointable. Moroever, an unbounded $\cinfty$dga $B$ together with a morphism $A_{\bullet}\rightarrow (B_{\bullet})^{\rmalg}_{\mathbf{dg}}$ of cdgas exhibits a unit transformation for the upper horizontal adjunction if and only if the map $\ker(\del_{A_0})\rightarrow \ker(\del_{B_0})^{\rmalg}$ is isomorphic to the unit map $\ker(\del_{A_0})\rightarrow F^{\cinfty}(\ker(\del_{A_0}))^{\rmalg}$ and the commuting diagram 
\[
\begin{tikzcd}
\ker(\del_{A_0}) \ar[d] \ar[r] & \ker(\del_{B_0})^{\rmalg}\ar[d] \\
A_{\bullet} \ar[r] &(B_{\bullet})^{\rmalg}_{\mathbf{dg}}
\end{tikzcd}
\]
of cdgas is a pushout. The same assertion holds for nonnegatively graded $\cinfty$dgas.
\item Let $K$ be a small simplicial set, then a diagram $f:K^{\rhd}\rightarrow \cinfty\mathbf{dga}^{\geq 0}$ is a colimit diagram if and only if $f_0:K^{\rhd}\rightarrow \cinfty\mathsf{ring}$ is a colimit diagram and the square 
\[
\begin{tikzcd}
\colim_K f_0^{\rmalg} \ar[d]\ar[r] & \colim_K f^{\rmalg}_{\mathbf{dg}} \ar[d] \\
f(\infty)_0^{\rmalg}\ar[r] & f(\infty)^{\rmalg}
\end{tikzcd}
\]
is a pushout of cdgas.  
\end{enumerate}
\end{prop}
\begin{proof}
In the nonnegatively graded case, the proof of $(1)$ follows immediately from Propositions \ref{prop:degzeroadj}, \ref{prop:leftadjseccocart} and \ref{prop:pullbackcocartleftadj}. In the unbounded case, in order to invoke Propositions \ref{prop:leftadjseccocart} and \ref{prop:pullbackcocartleftadj}, we have to verify that condition $(*)$ of Proposition \ref{prop:pullbackcocartleftadj} is satisfied; invoking Propositions \ref{prop:degzeroadj} and \ref{prop:leftadjseccocart}, we must show that for any unbounded cdga $A_{\bullet}$ determining a unit map $\ker(\del_{A_0})\rightarrow F^{\cinfty}(\ker(\del_{A_0}))^{\rmalg}$ and a pushout diagram
\[
\begin{tikzcd}
\ker(\del_{A_0}) \ar[d]\ar[r] & F^{\cinfty}(\ker(\del_{A_0}))^{\rmalg}\ar[d] \\
A_{\bullet}\ar[r] & B_{\bullet}
\end{tikzcd}
\]
of cdgas, the right vertical map becomes an isomorphism upon applying $\ker(\del_0)$. For brevity, denote the commutative algebra $F^{\cinfty}(\ker(\del_{A_0}))^{\rmalg}$ by $F$. We can identify the underlying real vector space of $B_0$ with $F\oplus A_0/\ker(\del_{A_0})\otimes_{\ker(\del_{A_0})}F$ and the underlying real vector space of $B_{-1}$ with $A_{-1}\otimes_{\ker(\del_{A_0})}F$ so that the differential in degree $0$ is given by the composition
\[F\oplus A_0/\ker(\del_{A_0})\otimes_{\ker(\del_{A_0})}F  \longrightarrow A_0/\ker(\del_{A_0})\otimes_{\ker(\del_{A_0})}F \longrightarrow  A_{-1}\otimes_{\ker(\del_{A_0})}F\]
where the first map projects away the factor $F$ and the second map is induced by the differential $\del_{A_0}$. To identify the kernel of this composition with $F$, it suffices to show that the second map is injective. The differential $A_0/\ker(\del_{A_0})\rightarrow A_{-1}$ is a map of $\ker(\del_{A_0})$-modules, so we have an exact sequence 
\[ \tor^{\ker(\del_{A_0})}_{1}(A_{-1}/\mathrm{Im}(\del_{A_0}),F) \longrightarrow  A_0/\ker(\del_{A_0})\otimes_{\ker(\del_{A_0})}F\longrightarrow  A_{-1}\otimes_{\ker(\del_{A_0})}F.   \]
It suffices to argue the vanishing of the $F$-module $\tor^{\ker(\del_{A_0})}_{1}(A_{-1}/\mathrm{Im}(\del_{A_0}),F)$. This follows from Proposition \ref{prop:freecinftycommutestruncation}. Note that $(2)$ follows from Proposition \ref{prop:colimleftadjcocart}. 
\end{proof}
\begin{rmk}
We do not, at this time, have a similar description of colimits in unbounded $\cinfty$dgas. To apply the argument of Proposition \ref{prop:colimleftadjcocart} to the category $\cinfty\mathbf{dga}$, we would have to know that the canonical comparison map $\colim_Kf^{\rmalg}\rightarrow (\colim_kf)^{\rmalg}$ for diagrams $f:K\rightarrow \cinfty\mathsf{ring}$ were flat, which we do not (see Proposition \ref{prop:projectionflatconsequences}).
\end{rmk}
\begin{prop}\label{prop:freecinftydgahomotopical}
The functor $(F^{\cinfty}_{\mathbf{dg}})^{\rmalg}_{\mathbf{dg}}:\mathbf{cdga}_{\R}\rightarrow\mathbf{cdga}_{\R}$ is homotopical, that is, it preserves weak equivalences. Similarly, the functor $(F^{\cinfty}_{\mathbf{dg}})^{\rmalg}_{\mathbf{dg}}:\mathbf{cdga}^{\geq0}_{\R}\rightarrow\mathbf{cdga}^{\geq 0}_{\R}$ on nonnegatively graded objects is homotopical.
\end{prop}
\begin{proof}
We give the proof in the unbounded case, the one for the nonnegatively graded case being the same. Let $A_{\bullet}$ be an unbounded cdga, then it follows from Proposition \ref{prop:freecinftydgacolim} $(1)$ that the underlying cdga of $F^{\cinfty}_{\mathbf{dg}}(A_{\bullet})$ fits into a pushout diagram 
\[
\begin{tikzcd}
\ker(\del_{A_0}) \ar[d] \ar[r] & F^{\cinfty}(\ker(\del_{A_0}))^{\rmalg}\ar[d] \\
A_{\bullet} \ar[r] & F^{\cinfty}_{\mathbf{dg}}(A_{\bullet})^{\rmalg}.
\end{tikzcd}
\]  
of cdgas. Let $f:A_{\bullet}\rightarrow B_{\bullet}$ be a map of cdgas. It follows from Proposition \ref{prop:freecinftypushout} and the pushout diagram above that in case $f$ induces a surjection $\ker(\del_{A_0})\rightarrow\ker(\del_{B_0})$, then the diagram
\[
\begin{tikzcd}
A_{\bullet} \ar[r,"f"] \ar[d] & B_{\bullet}\ar[d] \\
F^{\cinfty}_{\mathbf{dg}}(A_{\bullet})^{\rmalg} \ar[r] & F^{\cinfty}_{\mathbf{dg}}(B_{\bullet})^{\rmalg}
\end{tikzcd}
\] 
is a pushout. It follows from Proposition \ref{prop:freecinftycommutestruncation} that the vertical maps are flat, so the diagram is also homotopy pushout in the model category $\mathbf{cdga}_{\R}$. It follows that in case $f:A_{\bullet}\rightarrow B_{\bullet}$ is a quasi-isomorphism of unbounded cdgas that induces a surjection on 0-cycles, the induced map $F^{\cinfty}_{\mathbf{dg}}(A_{\bullet})_{\mathbf{dg}}^{\rmalg}\rightarrow F^{\cinfty}_{\mathbf{dg}}(B_{\bullet})_{\mathbf{dg}}^{\rmalg}$ is a weak equivalence of cdgas. Suppose that $f:A_{\bullet}\rightarrow B_{\bullet}$ is a trivial fibration of cdgas and consider the morphism
\[
\begin{tikzcd}
A_1\ar[d,"f_1"] \ar[r,"\del_{A_1}"]& \ker(\del_{A_0})\ar[d,"f_0"] \ar[r]& H_0(A) \ar[d,"H_0(f)"] \ar[r] & 0\ar[d]    \\
B_1 \ar[r,"\del_{B_1}"]& \ker(\del_{B_0}) \ar[r]& H_0(B) \ar[r] & 0.  
\end{tikzcd}
\]
of exact sequences. By assumption, $f_1$ is a surjection and $H_0(f)$ is an isomorphism; it follows from the four lemma that $f_0$ is also a surjection (on $0$-cycles). We conclude that $(F^{\cinfty}_{\mathbf{dg}})^{\rmalg}_{\mathbf{dg}}$ carries trivial fibrations of cdgas to weak equivalences. Now we complete the proof by invoking Ken Brown's lemma and the fact that all objects of $\mathbf{cdga}_{\R}$ are fibrant.
\end{proof}

We have the following model structures on nonnegatively graded and unbounded $\cinfty$dgas, first observed in \cite{CR2,CR1}.

\begin{prop}[Carchedi-Roytenberg]\label{prop:modelcinftydga}
There are combinatorial model structures on $C^{\infty}\mathbf{dga}$ and $C^{\infty}\mathbf{dga}^{\geq0}$ that are transferred along the adjunction $(F^{\cinfty}_{\mathbf{dg}}\adj (\_)^{\rmalg}_{\mathbf{dg}})$ from $\mathbf{cdga}_{\R}$ and $\mathbf{cdga}^{\geq 0}_{\R}$. Specifically, a map $f$ of (nonnegatively graded or unbounded) $\cinfty$dgas is a fibration respectively a weak equivalence if and only if $f^{\rmalg}_{\mathbf{dg}}$ is a fibration respectively a weak equivalence, and the set of generating (trivial) cofibrations is the image under $F^{\cinfty}_{\mathbf{dg}}$ of the set of generating (trivial) cofibrations in $\mathbf{cdga}_{\R}$ and $\mathbf{cdga}^{\geq 0}_{\R}$. Explicitly, the generating (trivial) cofibrations of $\cinfty\mathbf{dga}$ are the following.
\begin{enumerate}[$(i)$]
\item The set of generating cofibrations of $\cinfty\mathbf{dga}$ contains the maps 
\begin{align}\label{gencof1}
 C^{\infty}(\R)\longrightarrow    C^{\infty}(\R)[\epsilon_1],\qquad \R[\epsilon_i]\longrightarrow \R[\epsilon_i,\epsilon_{i+1}],\,i\in \Z\setminus\{0\} .
\end{align} 
where $|\epsilon_i|=i$, and the differentials are given by $\del \epsilon_{i+1}=\epsilon_{i}$ for $i\neq 1$ and $\del\epsilon_1=x$ for $x$ the identity function on $\R$. Note that for $i=-1$, the kernel of $\del_{\R[\epsilon_{-1},\epsilon_{0}]_0}$ is isomorphic to $\R$, so that $\R[\epsilon_{-1},\epsilon_{0}]$ is canonically an unbounded $\cinfty$dga, even though it is a polynomial algebra in degree $0$.
\item  The set of generating trivial cofibrations $\cinfty\mathbf{dga}$ contains the maps
\begin{align}\label{gentrivcof1} \R\longrightarrow    C^{\infty}(\R)[\epsilon_1],\qquad  \R\longrightarrow   \R[\epsilon_i,\epsilon_{i+1}],\,i\in\Z\setminus \{0\}. 
\end{align} 
\end{enumerate}
The generating (trivial) cofibrations of $\cinfty\mathbf{dga}^{\geq 0}$ are the following.
\begin{enumerate}[$(i')$]
\item the set of generating cofibrations of $\cinfty\mathbf{dga}^{\geq 0}$ contains the maps 
\begin{align}\label{gencof}
\R\longrightarrow C^{\infty}(\R),\quad C^{\infty}(\R)\longrightarrow    C^{\infty}(\R)[\epsilon_1],\qquad \R[\epsilon_i]\longrightarrow \R[\epsilon_i,\epsilon_{i+1}],\,i\geq 1.
\end{align} 
\item The set of generating trivial cofibrations $\cinfty\mathbf{dga}^{\geq 0}$ contains the maps
\begin{align}\label{gentrivcof} \R\longrightarrow    C^{\infty}(\R)[\epsilon_1],\qquad  \R\longrightarrow   \R[\epsilon_i,\epsilon_{i+1}],\,i\geq 1. 
\end{align}
\end{enumerate} 
\end{prop}
\begin{proof}
We do the proof for unbounded $\cinfty$dgas; the proof for nonegatively graded $\cinfty$dgas is the same (in the remark below we give an alternative proof for the nonnegatively graded case). It follows from Remark \ref{rmk:cinftycompactgen} that $C^{\infty}\mathbf{dga}$ is compactly generated. Let $I$ and $J$ be the classes of generating cofibrations and generating trivial cofibrations in $(i)$ and $(ii)$. Using adjointness, we see that a map $g$ in $C^{\infty}\mathbf{dga}$ is a trivial fibration (i.e. a weak equivalence and a fibration) if and only if $g$ satisfies the right lifting property with respect to the morphisms in the class $I$, and that $g$ is a fibration if and only if it satisfies the right lifting property against the morphisms in the class $J$. Now it follows from the small object argument that all morphisms can be factored by a morphism in the class $\overline{J}$, the weak saturation of $J$, followed by a fibration. Similarly, all morphisms can be factored by a morphism in the class $\overline{I}$, the weak saturation of $I$, followed by a trivial fibration. By a standard argument, it is enough to show that every morphism in $\overline{J}$ is a weak equivalence as a morphism in $\mathbf{cdga}_{\R}$. This is an immediate consequence of Proposition \ref{prop:freecinftydgahomotopical}. 
\end{proof}
\begin{rmk}
In the nonnegatively graded setting, we can provide an argument for the existence of the transferred model structure that applies more generally (to any Fermat theory $\mathsf{Poly}_k\rightarrow \mathrm{T}$ over a commutative ring $k$ containing $\Q$). We wish to show that every morphism in the class $\overline{J}'$, the weak saturation of the class $(ii')$ of generating trivial cofibrations above, is a weak equivalence. Since the forgetful functor $(\_)_{\mathbf{dg}}^{\rmalg}:\cinfty\mathbf{dga}^{\geq0}\rightarrow \mathbf{cdga}^{\geq0}_{\R}$ preserves filtered colimits, weak equivalences in $\cinfty\mathbf{dga}^{\geq0}$ are stable under retracts and transfinite compositions so we are reduced to showing that a pushout of any map in $C^{\infty}\mathbf{dga}^{\geq 0}$ along any of the generating trivial cofibrations of \eqref{gentrivcof} is a weak equivalence. For the map $\R\rightarrow \R[\epsilon_i,\epsilon_{i+1}]$, this follows immediately from Proposition \ref{prop:freecinftydgacolim} $(2)$: a pushout along this map is simply a coproduct with $\R[\epsilon_i,\epsilon_{i+1}]$ in $\mathbf{cdga}_{\R}^{\geq 0}$. For the case of $\R\rightarrow C^{\infty}(\R)[\epsilon]$, we have to show that for any $C^{\infty}$dga $A_{\bullet}$, the canonical map \[f:A_{\bullet}\longrightarrow A_{\bullet}\otimes^{\infty}C^{\infty}(\R)[\epsilon]\]
is a quasi-isomorphism, where $\otimes^{\infty}$ now denotes the coproduct in $C^{\infty}\mathbf{dga}^{\geq 0}$. The map $f$ admits a retraction that fits into a pushout diagram 
\[\begin{tikzcd}
\cinfty(\R)[\epsilon] \ar[r] \ar[d] & \R  \ar[d] \\
A_{\bullet}\otimes^{\infty}C^{\infty}(\R)[\epsilon] \ar[r] & A_{\bullet}
\end{tikzcd} \]
of $\cinfty$dgas. To show that $f$ is a weak equivalence, it suffices to show the lower horizontal map in the diagram is one, but since the upper horizontal map is a weak equivalence and a surjection in degree 0, this is guaranteed by Proposition \ref{prop:freecinftydgacolim} $(2)$ and the fact that the model category $\mathbf{cdga}^{\geq 0}_{\R}$ is left proper.
\end{rmk}

We denote the $\infty$-category of (fibrant)-cofibrant $C^{\infty}$dgas localized at the weak equivalences by $C^{\infty}\mathsf{Alg}^{\geq0}$. Note that there is an obvious fully faithful and coproduct preserving functor $\cartsp^{op}\hookrightarrow C^{\infty}\mathsf{Alg}^{\geq0}$. This functor left Kan extends to yield a colimit preserving functor $\varphi:sC^{\infty}\mathsf{ring}\rightarrow C^{\infty}\mathsf{Alg}^{\geq 0}$. Similarly, we denote the \infcat of (fibrant)-cofibrant unbounded $\cinfty$dgas localized at the weak equivalences by $\cinfty\mathsf{Alg}$. The obvious left Quillen functor $\cinfty\mathbf{dga}^{\geq0}\rightarrow \cinfty\mathbf{dga}$ determines a functor $\cinfty\mathsf{Alg}^{\geq 0}\rightarrow\cinfty\alg$.  Using that fibrations and trivial fibrations in $C^{\infty}\mathbf{dga}^{\geq0}$ are stable under filtered colimits, we deduce that taking filtered colimits in $C^{\infty}\mathbf{dga}^{\geq0}$ preserves trivial fibrations and thus (by Ken Brown's Lemma and the fact that all objects are fibrant) all weak equivalences, so that filtered colimits are also homotopy colimits. 
\begin{thm}[$\cinfty$-Dold-Kan Correspondence]\label{smoothdoldkan}
The functor $\varphi:sC^{\infty}\mathsf{ring}\rightarrow C^{\infty}\mathsf{Alg}^{\geq 0}$ induced by the fully faithful inclusion $\cartsp^{op}\hookrightarrow C^{\infty}\mathsf{Alg}^{\geq0}$ is an equivalence of $\infty$-categories.
\end{thm}
\begin{proof}
We have a commuting diagram 
\[
\begin{tikzcd}
\mathsf{Poly}_{\R}^{op} \ar[d] \ar[r] & \cartsp^{op}\ar[d]\\
\calg^{\geq0}_{\R}\ar[r,"\mathbf{L}F^{\cinfty}_{\mathbf{dg}}"] & C^{\infty}\mathsf{Alg}^{\geq0}
\end{tikzcd}
\]
of \infcatst, where $\mathbf{L}F^{\cinfty}_{\mathbf{dg}}$ is the left derived functor of the free $C^{\infty}$dga functor. Passing to the sifted colimit completion (\cite{HTT}, Corollary 5.3.6.10), we obtain a commuting diagram  
\[
\begin{tikzcd}
\scring_{\R} \ar[d,"\simeq"] \ar[r,"F^{\cinfty}"] & \sring\ar[d,"\varphi"]\\
 \calg^{\geq 0}_{\R}\ar[r,"\mathbf{L}F^{\cinfty}_{\mathbf{dg}}"]& C^{\infty}\mathsf{Alg}^{\geq0}
\end{tikzcd}
\]
of presentable \infcats and functors admitting right adjoints between them. Let $U:C^{\infty}\mathsf{Alg}^{\geq0}\rightarrow \sring$ be a right adjoint to $\varphi$, and let $\icatd\subset \sring$ be the full subcategory spanned by those objects $C$ for which the unit map $C\rightarrow U(\varphi(C))$ is an equivalence. It suffices to show that $\icatd= \sring$, and that $U$ is conservative. Since $\varphi$ is a left Kan extension along the functor $\cartsp^{op}\hookrightarrow C^{\infty}\mathsf{Alg}^{\geq0}$, the full subcategory $\cartsp^{op}\subset \sring$ lies in $\icatd$. Passing to right adjoints in the diagram above, we have a diagram 
\[
\begin{tikzcd}
C^{\infty}\mathsf{Alg}^{\geq0} \ar[dr,"G"] \ar[rr,"U"] && \sring \ar[dl,"(\_)^{\rmalg}"] \\
& \calg^{\geq 0}_{\R}
\end{tikzcd}
\]
where $G=\mathbf{R}(\_)^{\rmalg}_{\mathbf{dg}}$ is the right derived forgetful functor. As $G$ and $(\_)^{\rmalg}$ are both conservative, $U$ is also conservative. Let $\mathcal{K}:=\{f_i:K_i\rightarrow C^{\infty}\mathsf{Alg}^{\geq0}\}$ be the collection of small diagrams in $C^{\infty}\mathsf{Alg}^{\geq0}$ such that
\begin{enumerate}[$(a)$]
    \item $G$ preserves colimits of diagrams in $\mathcal{K}$.
    \item $(\_)^{\rmalg}$ preserves colimits of diagrams in $\mathcal{K}$ after applying $U$.
\end{enumerate}
Now note that, as $(\_)^{\rmalg}$ is conservative, $U$ also preserves the colimits of the diagrams in $\mathcal{K}$. We observe the following:
\begin{enumerate}[$(1)$]
    \item All filtered diagrams are in $\mathcal{K}$, since the underived functor $G:C^{\infty}\mathbf{dga}^{\geq0}\rightarrow\mathbf{cdga}^{\geq0}_{\R}$ preserves ordinary filtered colimits, which are also homotopy colimits as the model structure on $\cinfty\mathbf{dga}^{\geq0}$ is combinatorial.
    \item Pushouts diagrams along the map \[\mathbf{L}F^{\cinfty}_{\mathbf{dg}}(\R[x]\rightarrow \R)= \varphi(C^{\infty}(\R)\rightarrow \R) \]
    are in $\mathcal{K}$. To see that $(a)$ holds, note that this map is modelled by the generating cofibration $C^{\infty}(\R)\rightarrow C^{\infty}(\R)[\epsilon]$ with $|\epsilon| =1$, which is also a cofibration in $\mathbf{cdga}^{\geq0}_{\R}$. As $\mathbf{cdga}^{\geq0}_{\R}$ is left proper, this suffices. Observe that applying $U$ to this map yields an effective epimorphism (because this can be checked by applying $(\_)^{\rmalg}$), so $(b)$ follows by Corollary \ref{cor:algeffepi}.
    \item Pushouts diagrams along the map \[\mathbf{L}F^{\cinfty}_{\mathbf{dg}}(\R[\epsilon_n]\rightarrow \R)= \varphi(\Sigma^nC^{\infty}(\R)\rightarrow \R) \]
    where $n\geq 1$ and $|\epsilon_n| =n$ are in $\mathcal{K}$. Again, $(a)$ holds because this map is modelled by the generating cofibration $\R[\epsilon_n]\rightarrow \R[\epsilon_n,\epsilon_{n+1}]$ and $(b)$ holds because applying $U$ yields an effective epimorphism.
\end{enumerate}
It follows that $U$ preserves the colimits described above, so $\icatd\subset\sring$ is stable under filtered colimits and pushouts along the maps $\Sigma^nC^{\infty}(\R)\rightarrow \R$ for $n\geq 0$. All good $\R$-cell objects in $\sring$ are constructed out of such colimits from the subcategory $\cartsp^{op}$ so Proposition \ref{prop:afpcell} shows that we indeed have $\icatd=\sring$.
\end{proof}
\begin{prop}
The functor $\sring\rightarrow\cinfty\mathsf{Alg}$ induced by the left Quillen functor $\cinfty\mathbf{dga}^{\geq 0}\rightarrow \cinfty\mathbf{dga}$ is fully faithful. Moreover, there is a commuting diagram 
\[
\begin{tikzcd}
\cinfty\alg \ar[d,"\tau_{\geq0}"] \ar[r,"(\_)^{\rmalg}"] & \calg_{\R} \ar[d,"\tau_{\geq0}"] \\
\sring \ar[r,"(\_)^{\rmalg}"] & \calg^{\geq0}_{\R}
\end{tikzcd}
\]
in $\mathsf{Pr}^{\mathrm{R}}_{\omega}$ (compactly generated presentable \infcats and right adjoint continuous functors between them) which is vertically left adjointable. The horizontal functors of this diagram are conservative.
\end{prop}
\begin{proof}
We have a strictly commuting diagram of left derived functors
\[
\begin{tikzcd}
\cinfty\mathbf{dga}^{\mathrm{fc}} & \mathbf{cdga}^{\mathrm{fc}}_{\R} \ar[l] \\
\cinfty\mathbf{dga}^{\geq 0,\mathrm{fc}} \ar[u] &\mathbf{cdga}^{\geq 0,\mathrm{fc}}_{\R} \ar[l] \ar[u]
\end{tikzcd}
\]
between fibrant-cofibrant objects, using the fact that all objects in these model categories are fibrant. All functors in this diagram preserve weak equivalences, so taking the coherent nerve applying the fibrant replacement functor in the model category of marked simplicial sets, we obtain a (strictly) commuting diagram of \infcatst. Each of these functors admits a right adjoint, and we obtain the desired diagram of \infcats commuting up to homotopy by passing the right adjoints. The unit of the adjunction on the left is an equivalence, so in order for the square to be $\tau_{\geq 0}$-left adjointable, it suffices to show that the unit of the adjunction on the right is also an equivalence, which is clear. 
\end{proof}
We do not know whether the model category structure we have constructed on $\cinfty\mathbf{dga}^{\geq0}$ is left proper (but we suspect it is not). We offer the following criterion for the recognition of homotopy pushouts in $\cinfty\mathbf{dga}^{\geq0}$.
\begin{lem}\label{lem:hpushoutcrit}
Let 
\[
\begin{tikzcd}
A_{\bullet}\ar[r] \ar[d] & B_{\bullet} \ar[d]\\
C_{\bullet}\ar[r] & D_{\bullet}
\end{tikzcd}
\]
be a diagram of $\cinfty$dgas. Suppose that
\begin{enumerate}[$(i)$]
    \item The map $H_0(A_{\bullet})\rightarrow H_0(C_{\bullet})$ is a surjection.
    \item The underlying diagram of cdgas is a homotopy pushout.
\end{enumerate}
Then the diagram is a homotopy pushout in $\cinfty\mathbf{dga}^{\geq 0}$.
\end{lem}
\begin{proof}
Since the functor $(\_)^{\rmalg}_{\mathbf{dg}}$ preserves weak equivalences, we have a commuting diagram 
\[
\begin{tikzcd}
\cinfty\mathbf{dga}^{\geq0} \ar[d] \ar[r,"(\_)^{\rmalg}_{\mathbf{dg}}"] & \mathbf{cdga}_{\R}^{\geq 0} \ar[d] \\
\sring\ar[r,"(\_)^{\rmalg}"] &  \scring_{\R}
\end{tikzcd}
\]
where the vertical functors implement localization at the weak equivalences. The left vertical functor carries the map $A_{\bullet}\rightarrow B_{\bullet}$ to an effective epimorphism, so the lemma follows from Corollary \ref{cor:algeffepi} and the commutativity of the diagram of \infcatst.
\end{proof}
\begin{cor}\label{cor:hpushoutcrit}
Let 
\[
\begin{tikzcd}
A_{\bullet}\ar[r] \ar[d] & B_{\bullet} \ar[d]\\
C_{\bullet}\ar[r] & D_{\bullet}
\end{tikzcd}
\]
be a pushout diagram of $\cinfty$dgas. Suppose that
\begin{enumerate}[$(i)$]
    \item Either the map $A_0\rightarrow B_0$ or the map $A_0\rightarrow C_0$ is a surjection.
    \item Either the map $H_0(A_{\bullet})\rightarrow H_0(B_{\bullet})$ or the map $H_0(A_{\bullet})\rightarrow H_0(C_{\bullet})$ is a surjection.
    \item Either the map $A_0\rightarrow B_0$ or the map $A_0\rightarrow C_0$ is a cofibration of cdgas.
\end{enumerate}
Then the diagram is a homotopy pushout in $\cinfty\mathbf{dga}^{\geq0}$.
\end{cor}
\begin{proof}
Condition $(i)$ implies that the diagram is also a pushout of cdga’s, by Proposition \ref{prop:freecinftydgacolim} $(2)$. Then condition $(iii)$ and the left properness of $\mathbf{cdga}_{\R}^{\geq 0}$ imply that the diagram is a homotopy pushout of cdga’s so we conclude using $(ii)$ and the previous lemma.
\end{proof}

\subsubsection{Dg-manifolds}
The remainder of this subsection is concerned with the provision of explicit differential graded models of affine derived $\cinfty$-schemes of finite presentation of finite presentation. We will import \emph{differential graded manifolds} into our theory, which are ubiquitous in mathematical physics and more physics adjacent literature, where they are also known as (homologically nonnegatively graded) \emph{$NQ$-supermanifold} and we will see that these dg-manifolds provide natural examples of affine derived $\cinfty$-schemes of finite presentation. Our eventual goal is to prove that the set of examples provided by dg-manifolds is exhaustive. In fact, in part II of this series we will prove something more precise due to Carchedi \cite{carchedidgman}: the category of dg-manifolds admits a natural homotopy theory as sketched in the introduction which Behrend-Liao-Xu in \cite{BLX} enhance to the structure of a category of fibrant objects, and we wish to show that the underlying \infcat coincides with $\dafffp$. In the last part of this subsection, we take a first step towards this result by showing that pullback of dg-manifolds along fibrations are carried to pullbacks of affine derived $\cinfty$-schemes.

\begin{ex}[Koszul $C^{\infty}$dgas and derived zero loci of smooth functions]\label{koszulcinftydga}
Let $f=(f_1,\ldots,f_m):\R^n\rightarrow \R^m$ be a smooth function. The derived zero locus $\mathrm{d}Z(f)$ of this function can be represented by a finitely presented $C^{\infty}$dga, the homotopy pushout of the diagram
\begin{equation*}
\begin{tikzcd}
C^{\infty}(\R^m)\ar[d,"\mathrm{ev}_0"'] \ar[r,"f^*"] & C^{\infty}(\R^n)\ar[d]\\
\R \ar[r]& \mathrm{d}Z(f) 
\end{tikzcd}
\end{equation*}
We can get a nice model for $\mathrm{d}Z(f)$ if we replace the map $C^{\infty}(\R^m)\rightarrow \R$ with the cofibration $C^{\infty}(\R^m)\rightarrow C^{\infty}(\R^m)[e_1,\ldots,e_m]$, with $\del e_i=x^i$ for $x^i$, $1\leq i\leq m$ the coordinate functions on $\R^m$. Since all objects in the diagram are cofibrant, the derived zero locus is modelled by the ordinary pushout of $C^{\infty}$dgas $C^{\infty}(\R^n)\otimes^{\infty}_{C^{\infty}(\R^m)}C^{\infty}(\R^m)[e_1,\ldots,e_m]$. We note that $C^{\infty}(\R^m)\rightarrow C^{\infty}(\R^m)[e_1,\ldots,e_m]$ is surjective in degree 0, so Proposition \ref{prop:freecinftydgacolim} $(2)$ asserts that the derived zero locus is given by the tensor product $C^{\infty}(\R^n)\otimes_{C^{\infty}(\R^m)}C^{\infty}(\R^m)[e_1,\ldots,e_m]$ of cdgas. This $C^{\infty}$dga is isomorphic to the \emph{Koszul $C^{\infty}$dga} $C^{\infty}(\R^n)[e_1,\ldots,e_m]$ with differential $\del e_i=f_i$, through the map 
\[C^{\infty}(\R^n)\otimes_{C^{\infty}(\R^m)}C^{\infty}(\R^m)[e_1,\ldots,e_m]\rightarrow C^{\infty}(\R^n)[e_1,\ldots,e_m],\quad h\otimes \left(g_0+\sum_{i=1}^m g_ie_i\right)\mapsto hf^*(g_0)+\sum_{i=1}^m hf^*(g_i)e_i.\]
\end{ex}

\begin{ex}[Kuranishi $C^{\infty}$dgas and derived critical loci of smooth functions]\label{ex:kurcdga} This example is largely a translation to the smooth setting of Vezzosi's notes on derived critical loci \cite{vez1}. Let $E\rightarrow M$ be a finite rank vector bundle over a manifold $M$. Generalizing the example above, we would like to find a convenient $C^{\infty}$dga model for the derived zero locus of some smooth section $s:M\rightarrow E$. We start by taking a suitable cofibrant replacement of the map $0^*:C^{\infty}(E)\rightarrow C^{\infty}(M)$ given by pulling back along the zero section: consider the $C^{\infty}$dga
\[ C^{\infty}(E)\otimes_{C^{\infty}(M)}\Gamma(\Lambda^{\bullet}E^{\vee}),\quad \del(f\otimes t)(x,v_x)=f(x,v_x)t|_x(v_x),\,x\in M,\,v\in E_x\text{ and }t\in \Gamma(\Lambda^{\bullet}E^{\vee}),\,|t| =1  \]
(as we explained in the previous example, it doesn't matter whether we take a pushout of cdgas or $C^{\infty}$dgas here because the map $C^{\infty}(M)\rightarrow \Gamma(\Lambda^{\bullet}E^{\vee})$ is an isomorphism in degree 0). We claim that the factorization
\[C^{\infty}(E)\longrightarrow C^{\infty}(E)\otimes_{C^{\infty}(M)}\Gamma(\Lambda^{\bullet}E^{\vee})\longrightarrow C^{\infty}(M) \]
is a cofibration followed by a trivial fibration. Indeed, we note that this factorization is functorial in $E$ so we only have to check the claim for $E$ a trivial bundle by stability of cofibrations and trivial fibrations under retracts and the fact that any vector bundle is a retract of a trivial one. In the case of the trivial bundle $\R^n\times M\rightarrow M$, the factorization is simply
\[C^{\infty}(M\times \R^n)\longrightarrow C^{\infty}(M\times \R^n)[e_1,\ldots,e_n]\longrightarrow C^{\infty}(M),\quad|e_i| =1,\,\del e_i=x^i,\,1\leq i\leq n. \]
Note that the first map is a pushout of a coproduct of generating cofibrations (and thus a cofibration), and the second map is a quasi-isomorphism by Lemma \ref{projresolutiontransverse} and degreewise surjective (and thus a trivial fibration). Since $s^*$, the pullback along the zero section, is an effective epimorphism so Corollary \ref{cor:hpushoutcrit} implies that the derived zero locus of the smooth section $s:M\rightarrow E$ is computed by the ordinary pushout of connective cdgas
\begin{equation*}
\begin{tikzcd}
C^{\infty}(E)\ar[d]\ar[r,"s^*"]& C^{\infty}(M)\ar[d]\\
C^{\infty}(E)\otimes_{C^{\infty}(M)}\Gamma(\Lambda^{\bullet}E^{\vee}) \ar[r]& \mathrm{d}Z(s)
\end{tikzcd}    
\end{equation*}
so the derived zero locus is simply the tensor product
\[ \mathrm{d}Z(s)=C^{\infty}(M)\otimes_{C^{\infty}(E)} C^{\infty}(E)\otimes_{C^{\infty}(M)}\Gamma(\Lambda^{\bullet}E^{\vee})\simeq \Gamma(\Lambda^{\bullet}E^{\vee}),\]
with its obvious structure of a $C^{\infty}$-ring in degree 0. One readily verifies that under this isomorphism, the differential maps to $\del t=t(s)$, $t \in \Gamma(E^{\vee})$. We call the $C^{\infty}$dga $(\Gamma(\Lambda^{\bullet}E^{\vee}),\del t=t(s))$ a \emph{Kuranishi $C^{\infty}$dga} for $\mathrm{d}Z(s)$; these are the global sections of the affine Kuranishi models we considered in the introduction. We will see in part II that any affine derived $\cinfty$-scheme of finite presentation $X$ for which the \emph{cotangent complex} $\cotan_X$ has Tor-amplitude in $[0,1]$ can be realized as a Kuranishi $C^{\infty}$dga for some finite rank vector bundle $E\rightarrow M$.
\end{ex}
\begin{rmk}
While Koszul $C^{\infty}$dgas are always cofibrant, Kuranishi $C^{\infty}$dgas are usually not. For instance, $\cinfty(\R\setminus\{0\})$ is not cofibrant as a $\cinfty$dga, since a lift of an invertible element along a surjection need not be invertible. However, there is an alternative model for $\sring$ in which this $\cinfty$dga is cofibrant: there is a localization functor on $\cinfty\mathbf{dga}^{\geq0}$ that carries an object $A_{\bullet}$ to the pushout $A_{\bullet}\oinfty_{A_0}\tilde{A}_0$, where $\tilde{A}_0$ is the germ of the zero locus of the zero'th differential on $A$. The essential image of this functor admits a model structure right transferred from $\cinfty\mathbf{dga}^{\geq0}$ which is Quillen equivalent to $\cinfty\mathbf{dga}^{\geq0}$ \cite{Pridham}.
\end{rmk} 
Note that it follows from Proposition \ref{prop:afpcell} that every $A\in \sring_{\fp}$ admits a presentation as a \emph{retract} of a cofibrant $\cinfty$dga of the form \[\cinfty(\R^n)[e^1_1,\ldots,e^1_{n_1},e^2_1,\ldots,e^2_{n_2},\ldots,\ldots,e^m_1,\ldots,e^m_{n_m}]\]
with $|e^i_{j}|=i$ and some differential. In the previous example, we claimed that when $\cotan_A$ has Tor-amplitude $[0,1]$, $A$ can by modelled by a Kuranishi $\cinfty$dga. More generally, using manifolds and nontrivial vector bundles in place Cartesian spaces and trivial vector bundles, it is no longer necessary to take retracts to obtain all finitely presented simplicial $\cinfty$-rings: suppose that the cotangent complex $\cotan_A$ has Tor-amplitude in $[0,n]$, then $A$ is modelled by a \emph{dg-manifold of amplitude $n$}, in the sense defined below.
\begin{defn}
Let $\mathrm{RingTop}_{\mathbf{dg}_{\R}}$ be the category defined as follows.
\begin{enumerate}[$(1)$]
    \item Objects are pairs $(X,\Of_{X\bullet})$ where $X$ is a topological space and $\Of_{X\bullet}\in \shv_{\mathbf{cdga}_{\R}^{\geq 0}}(X)$ is a sheaf of nonnegatively graded differential graded $\R$-algebras.
    \item Morphisms between pairs $(X,\Of_{X\bullet})\rightarrow (Y,\Of_{Y\bullet})$ are pairs $(f,\alpha)$ of continuous map $f:X\rightarrow Y$ together with a map $\alpha:f^*\Of_{Y\bullet}\rightarrow \Of_{X\bullet}$ of sheaves of nonnegatively graded cdgas over $\R$.
\end{enumerate}
Let $\mathrm{RingTop}^{\mathrm{loc}}_{\mathbf{dg}_{\R}}\subset \mathrm{RingTop}_{\mathbf{dg}_{\R}}$ denote the subcategory whose objects are pairs $(X,\Of_{X\bullet})$ for which $\Of_{X0}$, the sheaf of ordinary $\R$-algebras of degree 0 elements of $\Of_{X\bullet}$, is a local sheaf of $\R$-algebras and whose morphisms are pairs $(f,\alpha)$ for which $\alpha$ induces in degree 0 a local morphism of local sheaves of $\R$-algebras. By design, the obvious functor
\[ \mathrm{RingTop}_{\mathbf{dg}_{\R}}\longrightarrow\mathrm{Top},\quad\quad (X,\Of_{X\bullet})\longmapsto X \]
is a biCartesian fibration which restricts to a Cartesian fibration $\mathrm{RingTop}^{\mathrm{loc}}_{\mathbf{dg}_{\R}}\rightarrow \mathrm{Top}$.\\
Recall that the assignment $M\mapsto (M,\cinfty_M)$ determines a fully faithful functor from the category of manifolds into the category of locally ringed $\R$-algebras. We define a category $\mathrm{RingMan}_{\mathbf{dg}_{\R}}$ as the pullback $\mathrm{RingTop}^{\mathrm{loc}}_{\mathbf{dg}_{\R}}\times_{\mathrm{RingTop}^{\mathrm{loc}}_{\R}}\mathsf{Man}$, which may be identified with the full subcategory of $\mathrm{RingTop}^{\mathrm{loc}}_{\mathbf{dg}_{\R}}$ spanned by objects $(X,\Of_{X\bullet})$ for which $(X,\Of_{X0})$ is a manifold with its sheaf of $\cinfty$ functions.
\end{defn}
\begin{rmk}
Since the stalk of the sheaf of $\cinfty$ functions at each point on a manifold is a local ring with residue field $\R$, the category $\mathrm{RingMan}_{\mathbf{dg}_{\R}}$ is a full subcategory of $\mathrm{RingTop}_{\mathbf{dg}_{\R}}$.
\end{rmk}
The category of dg-manifolds admits a natural Grothendieck topology inherited from the larger category of $\mathrm{RingTop}_{\mathbf{dg}_{\R}}^{\mathrm{loc}}$ of locally dg-ringed spaces: a family $\{(Y_i,\Of_{Y_i\bullet})\rightarrow (X,\Of_{X\bullet})\}$ is a covering family if and only if each map $Y_i\rightarrow X$ is an \emph{open immersion}: a diffeomorphism onto an open subset $U_i$ that identifies $\Of_{Y_i\bullet}$ with $\Of_{X\bullet}|_{U_i}$. This Grothendieck topology is subcanonical. As a special case, we have the following lemma.
\begin{lem}\label{lem:opencoproddg}
Let $\coprod_{i} (U_i,\Of_{X\bullet}|_{U_i})\rightarrow (X,\Of_{X\bullet})$ be a countable coproduct of open immersions of dg-manifolds, then for each map $(Y,\Of_{\bullet})\rightarrow (X,\Of_{X\bullet})$ of dg-manifolds, the canonical map 
\[ \coprod_{i}(U_i\times_X Y,\Of_{Y\bullet}|_{U_i\times_X Y})\longrightarrow \left(\coprod_i(U_i,\Of_{X\bullet}|_{U_i})\right)\times_{(X,\Of_{X\bullet})}(Y,\Of_{\bullet}) \]
of dg-manifolds is an isomorphism.
\end{lem}

\begin{cons}[dg-spectrum]
Let $\icat\subset \cinfty\mathbf{dga}^{\geq0}$ denote the full subcategory spanned by objects $A_{\bullet}$ for which $A_0$ is the ring of $\cinfty$ functions on a manifold. Taking global sections defines a functor $\Gamma:\mathrm{RingMan}^{op}_{\mathbf{dg}_{\R}}\rightarrow\icat$. We define an adjoint to this functor: let $A_{\bullet}$ be a $\cinfty$dga with $A_0=\cinfty(M)$ for some manifold $M$. We define an object $\spec_{\mathbf{dg}}A_{\bullet}\in \mathrm{RingMan}_{\mathbf{dg}_{\R}}$ as follows.
\begin{enumerate}[$(1)$]
    \item The underlying manifold of $\spec_{\mathbf{dg}}A_{\bullet}$ is $M$.
    \item The structure sheaf $\Of_{\spec_{\mathbf{dg}}A \bullet}$ of cdgas over $M$ of $\spec_{\mathbf{dg}}A_{\bullet}$ is the sheafification of the presheaf
    \[ \widetilde{\Of_{\spec_{\mathbf{dg}}A \bullet}}:\mathrm{Open}(M)\longrightarrow \mathbf{cdga}_{\R}^{\geq 0},\quad \quad U\longmapsto \cinfty(U)\otimes_{\cinfty(M)}A_{\bullet}, \]
    the pushout of nonnegatively graded cdgas.
\end{enumerate}
As sheafification of nonnegatively graded cdgas commutes with taking degree $0$ elements, the sheaf $\Of_{\spec_{\mathbf{dg}}A 0}$ coincides with the sheaf of $\cinfty$ functions on $M$ so it follows that $\spec_{\mathbf{dg}}A_{\bullet}$ is indeed an object of $\mathrm{RingMan}_{\mathbf{dg}_{\R}}$. Taking global sections of the sheafification map $\widetilde{\Of_{\spec_{\mathbf{dg}}A \bullet}}\rightarrow \Of_{\spec_{\mathbf{dg}}A \bullet}$ determines a map $\phi:A_{\bullet}\rightarrow \Gamma(\Of_{\spec_{\mathbf{dg}}A \bullet})$. 
\end{cons}
\begin{prop}
Let $A_{\bullet}\in \icat$, then composing with $\phi$ induces for each $(Y,\Of_{Y\bullet})$ a bijection
\[ \Hom_{\mathrm{RingMan}^{op}_{\mathbf{dg}_{\R}}}(\spec_{\mathbf{dg}}A_{\bullet},(Y,\Of_{Y\bullet})) \longrightarrow \Hom_{\icat}(\Gamma(\Of_{\spec_{\mathbf{dg}}A\bullet}),\Gamma(\Of_{Y\bullet})) \longrightarrow \Hom_{\icat}(A_{\bullet},\Gamma(\Of_{Y\bullet})).  \]
\end{prop}
\begin{proof}
Since the map $\mathrm{RingTop}_{\mathbf{dg}_{\R}}\rightarrow\mathrm{Top}$ is a coCartesian fibration, we have an identification
\[\Hom_{\mathrm{RingMan}^{op}_{\mathbf{dg}_{\R}}}(\spec_{\mathbf{dg}}A_{\bullet},(Y,\Of_{Y\bullet}))\cong\coprod_{f\in \Hom_{\mathrm{Man}}(Y,M)}\Hom_{\shv_{\mathbf{cdga}_{\geq 0}}}(\Of_{\spec_{\mathbf{dg}}A \bullet},f_*\Of_{Y\bullet})\]
so that the map in the statement of the lemma is induced by the various global sections morphisms
\[ \Hom_{\shv_{\mathbf{cdga}_{\geq 0}}}(\Of_{\spec_{\mathbf{dg}}A \bullet},f_*\Of_{Y\bullet}) \longrightarrow \Hom_{\icat}(\Gamma(\Of_{\spec_{\mathbf{dg}}A}),\Gamma(\Of_{Y\bullet})) \]
as $f$ ranges over the smooth maps $Y\rightarrow M$. Since $f_*\Of_{Y\bullet}$ is a sheaf, the map above on hom spaces fits into a commuting diagram 
\[
\begin{tikzcd}
\Hom_{\shv_{\mathbf{cdga}_{\geq 0}}}(\Of_{\spec_{\mathbf{dg}}A \bullet},f_*\Of_{Y\bullet}) \ar[d]\ar[r]& \Hom_{\icat}(\Gamma(\Of_{\spec_{\mathbf{dg}}A}),\Gamma(\Of_{Y\bullet})) \ar[d]\\
\Hom_{\pshv_{\mathbf{cdga}_{\geq 0}}}(\widetilde{\Of_{\spec_{\mathbf{dg}}A \bullet}},f_*\Of_{Y\bullet}) \ar[r] & \Hom_{\icat}(A_{\bullet},\Gamma(\Of_{Y\bullet}))
\end{tikzcd}
\]
where the left vertical map is a bijection, so it suffices to show that the map induced by all various lower horizontal global sections morphisms is a bijection. Let $A_{\bullet}\rightarrow \Of_{Y_{\bullet}}(Y)$ be a morphism of $\cinfty$dgas and let $f:Y\rightarrow M$ be the map on manifolds corresponding to the map $\cinfty(M)\rightarrow \Of_{Y0}(Y)$ in degree $0$, then we have a composition 
\[ \cinfty_M\longrightarrow f_*\cinfty_Y\longrightarrow f_*\Of_{Y\bullet}  \]
of sheaves of cdgas, and upon evaluating at $M$, we have an extension of this composition to a commuting square
\[
\begin{tikzcd}
\cinfty(M) \ar[d] \ar[r] & \cinfty(Y) \ar[d] \\
A_{\bullet} \ar[r] & \Of_{Y_{\bullet}}(Y)
\end{tikzcd}
\]
which together determine the top horizontal map in the diagram 
\[
\begin{tikzcd}
\mathrm{Open}(M)^{op} \times \Delta^2\coprod_{\{M\}\times\Delta^2}\{M\}\times\Delta^1\times\Delta^1 \ar[r] \ar[d,hook]& \mathbf{cdga}_{\R}^{\geq 0} \\
\mathrm{Open}(M)^{op} \times\Delta^1\times\Delta^1 \ar[ur,dotted]
\end{tikzcd}
\]
where the diagonal map is a left Kan extension. Restricting the diagonal map to $\mathrm{Open}(M)^{op}\times \{1\}\times\Delta^1$ results in a map of presheaves $\widetilde{\Of_{\spec_{\mathbf{dg}}A \bullet}}\rightarrow f_*\Of_{Y\bullet}$ and this assignment is clearly a right inverse to the global sections map. Conversely, every map of presheaves $\widetilde{\Of_{\spec_{\mathbf{dg}}A \bullet}}\rightarrow f_*\Of_{Y\bullet}$ for some smooth map $f:Y\rightarrow N$ determines a diagram as above, which is a left Kan extension by definition of $\widetilde{\Of_{\spec_{\mathbf{dg}}A\bullet}}$, so this assignment is also a left inverse.  
\end{proof}
We let $\spec_{\mathbf{dg}}$ denote the resulting left adjoint to $\Gamma$. The following result is classical.
\begin{prop}[Batchelor's Theorem]\label{prop:batchelor}
Let $(X,\Of_{X\bullet})$ be an object in $\mathrm{RingMan}_{\mathbf{dg}_{\R}}$, then the following are equivalent.
\begin{enumerate}[$(1)$]
    \item There is a degreewise finite dimensional graded vector space $V$ concentrated in degrees $[1,n]$ and an open cover $\{U_i\subset X\}_i$ such that for each $i$, there is an isomorphism $\Of_{X\bullet}|_U \cong \sym_{\cinfty_{U_i}}(\cinfty_{U_i}\otimes V)$ of sheaves of graded algebras.
    \item There is a graded vector bundle $L$ on $X$ concentrated in degrees $[1,n]$ of degreewise constant finite rank, together with a map $L\rightarrow \Of_{X\bullet}$ of graded $\cinfty_X$-modules inducing an isomorphism $\sym_{\Of_{X0}}^{\bullet}L\simeq \Of_{X\bullet}$ of sheaves of graded algebras.
\end{enumerate}
\end{prop}
\begin{defn}
The category of \emph{dg-manifolds of amplitude $n$} is the full subcategory of $\mathrm{RingMan}_{\mathbf{dg}_{\R}}$ spanned by objects satisfying any of the equivalent conditions of Proposition \ref{prop:batchelor}. Dg-manifolds of amplitude 0 are simply manifolds. We let $\mathbf{dgMan}$ denote the category of dg-manifolds of variable amplitude.
\end{defn}
We will compare dg-manifolds with derived $\cinfty$-schemes via the global sections functor $\Gamma:\mathbf{dgMan}^{op}\rightarrow \cinfty\mathbf{dga}^{\geq0}$. 
\begin{prop}\label{prop:gammadgff}
The functor $\Gamma:\mathbf{dgMan}^{op}\rightarrow \cinfty\mathbf{dga}^{\geq0}$ is fully faithful. Let $A_{\bullet}\in \icat$, that is, $A_0\cong \cinfty(M)$ for a manifold $M$. Then $A_{\bullet}$ lies in the essential image of $\Gamma|_{\mathbf{dgMan}^{op}}$ if and only if there exists a graded vector bundle $L$ of degreewise constant finite rank on $M$ concentrated in degrees $[1,n]$ for some $n\in \Z_{\geq1}$ and map $\Gamma(L)\rightarrow A_{\bullet}$ of graded $A_0$-modules which induces an isomorphism $\sym^{\bullet}_{A_0}\Gamma(L)\rightarrow A_{\bullet}$ of graded commutative $A_0$-algebras.       
\end{prop}
\begin{proof}
Let $(X,\Of_{X\bullet})$ be a dg-manifold and choose a generating graded vector bundle $L$ on $X$ concentrated in degrees $[1,n]$ so that we have an isomorphism $\sym^{\bullet}_{\cinfty_X}L\rightarrow \Of_{X\bullet}$ of sheaves graded $\cinfty_X$-algebras. The object $\sym^{\bullet}_{\cinfty_X}L$ is the sheafification of the presheaf of graded $\cinfty_X$-algebras carrying an open $V\subset X$ to the graded free $\cinfty(V)$-algebra $\bigoplus_{i=0}^{\infty} \sym^i_{\cinfty(V)} \Gamma(L|_V)$. Note that in order to show that $\Gamma(\Of_{X\bullet})\in \icat$ satisfies the stated condition, it suffices to show that the presheaf $V\mapsto \bigoplus_{i=0}^{\infty} \sym^i_{\cinfty(V)} \Gamma(L|_V)$ is already a sheaf. Since $\Gamma(L)$ is a graded projective and finitely presented $\cinfty(X)$-module, the graded locally free sheaf $L$ can be identified with the graded presheaf $V\mapsto \cinfty(V)\otimes_{\cinfty(X)}L$, by Proposition \ref{prop:modulealreadysheaf}. It follows that for any $i\geq 0$, the tensor module sheaf $ L^{\otimes i}$ is the sheafification of the presheaf carrying $V$ to $\cinfty(V)\otimes_{\cinfty(X)}\Gamma(L)^{\otimes i}$, where $\Gamma(L)^{\otimes i}$ is the $i$-fold tensor product in the category of graded $\cinfty(X)$-modules. The module $\Gamma(L)^{\otimes i}$ is again a graded projective and finitely presented $\cinfty(X)$-module which implies that the presheaf $V\mapsto \cinfty(V)\otimes_{\cinfty(X)}\Gamma(L)^{\otimes i}$ is a sheaf. It is easy to see that the symmetric tensors form a subsheaf. In summary, the sheaf $\sym^i_{\cinfty_X}L$ coincides with the presheaf $V\mapsto \sym^i_{\cinfty(V)} \Gamma(L|_V)$. To see that the presheaf $V\mapsto \bigoplus_{i=0}^{\infty} \sym^i_{\cinfty(V)} \Gamma(L|_V)$ is a sheaf, it suffices to observe that for degree reasons, the map $\bigoplus_{i=0}^{\infty}\sym^{i}_{\cinfty_X}L\rightarrow \prod_{i=0}^{\infty}\sym_{\cinfty_X}^{i}L$ of presheaves of graded $\cinfty_X$-modules is an isomorphism.\\
We show that for $A_{\bullet}\subset\icat$ satisfying the condition above, the object $\spec_{\mathbf{dg}}A_{\bullet}$ is a dg-manifold and the unit map $A\rightarrow\Gamma(\Of_{\spec_{\mathbf{dg}}A_{\bullet}})$ is an isomorphism. Choose an isomorphism $\cinfty(M)\cong A_{0}$ and a graded vector bundle $L$ on $M$ inducing an isomorphism $\sym_{\cinfty(M)}^{\bullet}\Gamma(L)\rightarrow A_{\bullet}$. The vector bundle $L$, viewed as a sheaf of graded $\cinfty_M$-modules, can be identified with the presheaf $U\mapsto \cinfty(U)\otimes_{\cinfty(M)}\Gamma(L)$, so we have a map 
\[ L\longrightarrow \widetilde{\Of_{\spec_{\mathbf{dg}}A \bullet}}  \]
of presheaves of $\cinfty_X$-modules which induces an isomorphism $\sym_{\cinfty_U}L|_U\cong \cinfty(U)\otimes_{\cinfty(M)}\Of_{\spec_{\mathbf{dg}}A \bullet}$ over each open $U\subset M$. The same argument as the one employed above shows that $\widetilde{\Of_{\spec_{\mathbf{dg}}A \bullet}}$ is already a sheaf, so $\spec_{\mathbf{dg}}A_{\bullet}$ is a dg-manifold. The fact that the presheaf $\widetilde{\Of_{\spec_{\mathbf{dg}}A \bullet}}$ is a sheaf also immediately implies that the unit map $A_{\bullet}\rightarrow \Gamma(\Of_{\spec_{\mathbf{dg}}A_{\bullet}})$ is an isomorphism. We have shown that $\Gamma(\mathbf{dgMan}^{op})$ can be characterized as the subcategory of $\icat$ described in the statement of the proposition, and that the functor $\spec_{\mathbf{dg}}|_{\Gamma(\mathbf{dgMan}^{op})}$ embeds $\Gamma(\mathbf{dgMan}^{op})$ fully faithfully into $\mathbf{dgMan}^{op}$. We will be done once we argue that the spectrum functor is essentially surjective. Let $(X,\Of_{X\bullet})$ be a dg-manifold, then the preceding arguments readily imply that the canonical map 
\[  \widetilde{\Of_{\spec_{\mathbf{dg}}\Gamma(\Of_{X\bullet})}}\longrightarrow \Of_{X\bullet} \]
of presheaves of differential graded $\cinfty_X$-algebras is an isomorphism.
\end{proof}
Our next goal is to show that the global sections of a dg-manifold determine a finitely presented object in $\sring$. 
\begin{lem}\label{lem:dgmanfinpres}
Let $(X,\Of_{X\bullet})$ be a dg-manifold, then $\Gamma(\Of_{X\bullet})$ is homotopically finitely presented.
\end{lem}
\begin{proof}
We proceed by induction on the amplitude. Dg-manifolds of amplitude $0$ are manifolds whose algebras of functions are finitely presented in $\sring$. Let $(X,\Of_{X\bullet})$ be a dg-manifold of dimension $1$, then $\Gamma(\Of_{X\bullet})$ is isomorphic to $\Gamma(\Lambda^{\bullet}E)$ for some finite rank vector bundle $E$ over $X$ by Proposition \ref{prop:gammadgff}. The fact that this object is homotopically finitely presented is the content of Example \ref{ex:kurcdga}. Now let $(X,\Of_{X\bullet})$ be a dg-manifold of amplitude $n>1$ so that we can find a graded vector bundle $L$ concentrated in degrees $[1,n]$ and degreewise of constant finite rank (viewed as a graded locally free sheaf of $\cinfty_X$-modules) and a map of graded $\cinfty_X$-modules $L\rightarrow \Of_{X\bullet}$ inducing an isomorphism of graded $\cinfty_X$-algebras $\sym^{\bullet}_{\cinfty_X}L\overset{\simeq}{\rightarrow} \Of_{X\bullet}$, which yields a differential on $\sym^{\bullet}_{\cinfty_X}L$. Write $L=\oplus_{i=1}^n L_i[i]$, where $L_i$ denotes the locally free sheaf of degree $i$ local sections of $L$, and write $L_{<j}$ for the graded vector bundle $\oplus_{1\leq i<j}L_i[i]$, then the differential restricted to $L_n[n]$ determines a map $L_n[n-1]\rightarrow \sym^{\bullet}_{\cinfty_X}L_{<n}$ of sheaves of graded $\cinfty_X$-modules, which in turn induces a map $\rho:\sym^{\bullet}_{\cinfty_X}L_n[n-1]\rightarrow \sym^{\bullet}_{\cinfty_X}L_{<n}$ of sheaves of graded $\cinfty_X$-algebras. Endowing $\sym^{\bullet}_{\cinfty_X}L_n[n-1]$ with the zero differential and $\sym^{\bullet}_{\cinfty_X}L_{<n}$ with the differential making it a differential graded sheaf of subalgebras of $\sym^{\bullet}_{\cinfty_X}L$, $\rho$ becomes a map of sheaves of differential graded $\cinfty_X$-algebras. By construction, we have a pushout diagram 
\[
\begin{tikzcd}
\sym^{\bullet}_{\cinfty_X}L_n[n-1] \ar[d] \ar[r,"\rho"] &  \sym^{\bullet}_{\cinfty_X}L_{<n}\ar[d] \\
\sym^{\bullet}_{\cinfty_X}(L_n[n-1]\oplus L_n[n]) \ar[r] & \Of_{X\bullet}
\end{tikzcd}
\]
of sheaves of dg $\cinfty_X$-algebras, where the left vertical map is induced by the inclusion $L_n[n-1]\rightarrow L_n[n-1]\oplus L_n[n]$ and the differential on $\sym^{\bullet}_{\cinfty_X}(L_n[n-1]\oplus L_n[n])$ is the zero map on $L_n[n-1]$ and the identity $L_n\rightarrow L_n$ on $L_n[n]$. Taking global sections of the diagram above yields the pushout diagram
\[
\begin{tikzcd}
\sym^{\bullet}_{\cinfty(X)}\Gamma(L_n[n-1]) \ar[d] \ar[r] &  \sym^{\bullet}_{\cinfty(X)}\Gamma(L_{<n})\ar[d] \\
\sym^{\bullet}_{\cinfty(X)}\Gamma(L_n[n-1]\oplus L_n[n]) \ar[r] & \sym^{\bullet}_{\cinfty(X)}\Gamma(L),
\end{tikzcd}
\]
of $\cinfty$dgas. The objects $(X,\sym^{\bullet}_{\cinfty_X}L[n-1])$ and $(X,\sym^{\bullet}_{\cinfty_X}L_{<n})$ are dg-manifolds of amplitude $n-1$, and the homology of the object $\sym^{\bullet}_{\cinfty(X)}\Gamma(L[n-1]\oplus L[n])$ is $\cinfty(X)$ in degree 0 and is thus finitely presented in $\sring$, so invoking the inductive hypothesis, it suffices to show that the diagram above is a homotopy pushout. Since the left vertical map induces an isomorphism in degree 0 and surjection on connected components we may invoke Corollary \ref{cor:hpushoutcrit} to reduce to the assertion that the left vertical map is a cofibration of cdgas. The assignment 
\[  L_n\longmapsto  \left(\sym^{\bullet}_{\cinfty(X)}\Gamma(L_n[n-1])\longrightarrow\sym^{\bullet}_{\cinfty(X)}\Gamma(L_n[n-1]\oplus L_n[n])\right)   \]
determines a functor $\mathsf{Vect}^{op}_{X}\rightarrow \fun(\Delta^1,\cinfty\mathbf{dga}^{\geq0})$ from the opposite of the category of finite rank vector bundles on $X$, so we may suppose that $L_n$ is a trivial bundle, in which case the map under consideration is obviously a cofibration of cdgas. 
\end{proof}
\begin{lem}\label{lem:dgproducts}
Let $\{(X_i,\Of_{X_i\bullet})\}_{i\in I}$ be a countable collection of dg-manifolds whose coproduct in locally dg-ringed spaces lies in $\mathbf{dgMan}$ (that is, the dimensions of the underlying manifolds, the amplitudes and the ranks of all bundles are all bounded as functions $I\rightarrow\Z_{\geq 0}$). Then $\Gamma$ carries the coproduct of the collection $\{(X_i,\Of_{X_i\bullet})\}_{i\in I}$ to a coproduct of affine derived $\cinfty$-schemes of finite presentation.
\end{lem}
\begin{proof}
The coproduct $\coprod_i (X_i,\Of_{X_i\bullet})$ is taken in the category $\mathrm{RingMan}_{\mathbf{dg}_{\R}}$ which the functor $\Gamma:\mathrm{RingMan}_{\mathbf{dg}_{\R}}^{op}\rightarrow \icat\subset \cinfty\mathbf{dga}^{\geq0}$, as a right adjoint, carries to a countable product. Since $\coprod_i (X_i,\Of_{X_i\bullet})$ is a dg-manifold, Lemma \ref{lem:dgmanfinpres} guarantees that $\Gamma(\coprod_i (X_i,\Of_{X_i\bullet}))$ presents an affine derived $\cinfty$-scheme of finite presentation, so it suffices to show that $\Gamma(\coprod_i (X_i,\Of_{X_i\bullet}))\cong \prod_i \Gamma(X_i,\Of_{X_i\bullet})$ is a homotopy product in $\cinfty\mathbf{dga}^{\geq0}$ of the collection $\{\Gamma(X_i,\Of_{X_i\bullet})\}_{i\in I}$. This collection determines an injectively fibrant diagram $I\rightarrow \cinfty\mathbf{dga}^{\geq0}$ as all objects of $\cinfty\mathbf{dga}^{\geq0}$ are fibrant.
\end{proof}
In the introduction, we introduced a notion of \emph{weak equivalence} on the category of dg-manifolds. The resulting homotopy theory is made more accessible by the following result of Behrend-Liao-Xu.
\begin{thm}[\cite{BLX}]
Let $(Y,\Of_{Y\bullet})$ and $(X,\Of_{X\bullet})$ be dg-manifolds with associated graded vector bundles $K$ and $L$. Let us say that a morphism $(f,\alpha):(Y,\Of_{Y\bullet})\rightarrow (X,\Of_{X\bullet})$ of dg-manifolds is a \emph{fibration} if the $\cinfty$ map $f:Y\rightarrow X$ is a submersion and the induced map of vector bundles $f^*L\rightarrow K\subset \oplus_{i=0}^{\infty}\sym^i_{\cinfty_Y}K$ is a degreewise monomorphism of vector bundles on $Y$. The category $\mathbf{dgMan}$ with the weak equivalences of Definition \ref{def:kurweakeq} and the fibrations just described is a \emph{category of fibrant objects} (see for instance \cite{cis1}, section 7.5).        
\end{thm}
\begin{rmk}
Let $f:(Y,\Of_{Y\bullet})\rightarrow (X,\Of_{X\bullet})$ be a fibration of dg-manifolds, and suppose that the underlying map of manifolds is a projection $X\times\R^n\rightarrow X$. Then the induced map $\Gamma(\Of_{X\bullet})\rightarrow\Gamma(\Of_{Y\bullet})$ is a cofibration of $\cinfty$dgas.    
\end{rmk}
\begin{prop}
Let 
\[
\begin{tikzcd}
(D,\Of_{D\bullet}) \ar[d,"g"] \ar[r] & (C,\Of_{C\bullet}) \ar[d,"f"]\\
(B,\Of_{B\bullet}) \ar[r] & (A,\Of_{A\bullet})
\end{tikzcd}
\]
be a pullback diagram of dg-manifolds where $f$ and $g$ are fibrations, then the diagram 
\[
\begin{tikzcd}
\Gamma(\Of_{A\bullet}) \ar[d] \ar[r] & \Gamma(\Of_{B\bullet})\ar[d] \\
\Gamma(\Of_{C\bullet}) \ar[r] & \Gamma(\Of_{D\bullet})
\end{tikzcd}
\]
is a homotopy pushout diagram of $\cinfty$dgas.
\end{prop}
\begin{rmk}
In the proof below, we reference the \inftop $\dstack_{\lfp}:=\shv(\daff_{\fp})$ of \emph{derived $\cinfty$-stacks locally of finite presentation}, and the fact that the Yoneda embedding $j:\daff_{\fp}\hookrightarrow\dstack_{\lfp}$ is fully faithful, that is, the \'{e}tale topology on $\daff_{\fp}$ is subcanonical.
\end{rmk}
\begin{proof}
We prove this by showing that the property of being a homotopy pushout is local on $A$, $B$ and $C$ so that we may suppose that all three are Cartesian spaces. First, we show the following special case.
\begin{enumerate}
    \item[$(*)$] Let $\coprod_i f_i:\coprod_i(C_i,\Of_{C_i\bullet})\rightarrow (A,\Of_{A\bullet})$ be a countable coproduct of open immersion of dg-manifolds, then $\Gamma$ carries every pullback diagram along $f=\coprod_i f_i$ to a homotopy pushout diagram in $\cinfty\mathbf{dga}^{\geq0}$.
\end{enumerate}
To prove this, we note that Lemma \ref{lem:opencoproddg} shows that for any map $(B,\Of_{B\bullet})\rightarrow(A,\Of_{A\bullet})$ of dg-manifolds, the diagram
\[
\begin{tikzcd}
\coprod_i(C_i\times_AB,\Of_{B_i\bullet}|_{C_i\times_AB})\ar[d]\ar[r] & \coprod_i(C_i\Of_{C_i\bullet})\ar[d] \\
(B,\Of_{B\bullet}) \ar[r] & (A,\Of_{A\bullet})
\end{tikzcd}
\]
is a pullback. Since $\Gamma$ carries countable coproducts of dg-manifolds to coproducts of derived $\cinfty$-schemes by Lemma \ref{lem:dgproducts}, it suffices to show that the diagram 
\[
\begin{tikzcd}
\coprod_i\spec\Gamma(\Of_{B_i\bullet}|_{C_i\times_AB})\ar[d]\ar[r] & \coprod_i\spec\Gamma(\Of_{C_i\bullet})\ar[d] \\
\spec\Gamma(\Of_{C\bullet}) \ar[r] & \spec\Gamma(\Of_{D\bullet})
\end{tikzcd}
\]
of affine derived $\cinfty$-schemes of finite presentation is a pullback. Since the Yoneda embedding $j:\daff_{\fp}\rightarrow \dstack_{\lfp}$ preserves colimits of those \'{e}tale diagrams in  $\daff_{\fp}$ that admit a colimit in $\daff_{\fp}$, it suffices to prove that the diagram 
\[
\begin{tikzcd}
\coprod_ij(\spec\Gamma(\Of_{B_i\bullet}|_{C_i\times_AB}))\ar[d]\ar[r] & \coprod_ij(\spec\Gamma(\Of_{C_i\bullet}))\ar[d] \\
j(\spec\Gamma(\Of_{C\bullet})) \ar[r] & j(\spec\Gamma(\Of_{D\bullet}))
\end{tikzcd}
\]
is a pullback diagram of derived $\cinfty$-stacks locally of finite presentation. Since colimits are universal in $\dstack_{\lfp}$, it suffices to show that for each $i$, the diagram is a pullback, that is, we reduce to the case of a single open immersion. In this case, it is easy to see that the assertion to be verified follows from the following one.
\begin{enumerate}
    \item[$(*')$] Let $f:(C,\Of_{C\bullet})\rightarrow (A,\Of_{A\bullet})$ is an open immersion of dg-manifolds, then the induced map $\Gamma(\Of_{A\bullet})\rightarrow\Gamma(\Of_{C_{\bullet}})$ in the \infcat $\sring$ exhibits $\Gamma(\Of_{C_{\bullet}})$ as a localization with respect to any characteristic function of the open subset $C\subset A$.
\end{enumerate} 
To prove the latter assertion, we note that we have a pushout diagram of $\cinfty$dgas 
\[
\begin{tikzcd}
\cinfty(A) \ar[d] \ar[r] & \cinfty(C) \ar[d] \\
\Gamma(\Of_{A\bullet}) \ar[r] & \Gamma(\Of_{C\bullet}).
\end{tikzcd}
\]
The upper horizontal map exhibits a localization, so it suffices to show that this diagram is a homotopy pushout. The left vertical map is the identity in degree 0 and induces a surjection on connected components, so by Lemma \ref{lem:hpushoutcrit}, it suffices to observe that this diagram is a homotopy pushout of cdgas as the upper horizontal map is flat; we have thus established $(*)$. Since the map $C\rightarrow A$ is a submersion, its image $f(C)$ is open in $A$, inducing an open immersion $(f(C),\Of_{A\bullet}|_{f(C)})\rightarrow (A,\Of_{A\bullet})$ and a diagram 
\[
\begin{tikzcd}
(D,\Of_{D\bullet}) \ar[d] \ar[r] & (C,\Of_{C\bullet}) \ar[d]\\ 
(f(C)\times_AB,\Of_{f(C)\times_AB\bullet}) \ar[r] \ar[d]& (f(C),\Of_{A\bullet}|_{f(C)})\ar[d] \\
(B,\Of_{B\bullet}) \ar[r] & (A,\Of_{A\bullet})
\end{tikzcd}
\]
where both squares are pullbacks of dg-manifolds. The lower square determines a homotopy pushout upon taking global sections by $(*)$, so it suffices to show that the same holds for the upper one, that is, we may suppose that $f$ is a surjective submersion. As such, we may choose a countable open cover $\{U_i\rightarrow A\}_{i\in I}$ such that 
\begin{enumerate}[$(1)$]
    \item For each finite collection $i_1,\ldots,i_n\in I$, the intersection $U_{i_1\ldots i_n}:=U_{i_1}\cap\ldots\cap U_{i_n}$ is either empty or diffeomorphic to $\R^n$ (that is, $\{U_i\rightarrow A\}_i$ is a good open cover).
    \item For each finite collection $i_1,\ldots,i_n\in I$, the pullback $U_{i_1\ldots i_n}\times_A C$ admits a cover $V_j\rightarrow U_{i_1\ldots i_n}\times_A C$ and a diffeomorphism $V_j \cong \R^m\times  U_{i_1\ldots i_n}$.
\end{enumerate}
Since the cover $\{U_i\rightarrow A\}_i$ is countable, the locally dg-ringed space $\coprod_{\overline{i}\in I^n} (U_{\overline{i}},\Of_{A\bullet|_{U_{\overline{i}}}})$ is a dg-manifold. In other words, the map $h_A:\coprod_i(U_i,\Of_{A\bullet}|_{U_i})\rightarrow (A,\Of_{A\bullet})$ admits a \v{C}ech nerve in the category of dg-manifolds and Lemma \ref{lem:opencoproddg} shows that in each simplicial degree, this \v{C}ech nerve is a coproduct of dg-manifolds of the form $(U_{i_1,\ldots,i_n},\Of_{A\bullet}|_{U_{i_1,\ldots,i_n}})$. Similarly, we may pull back the cover $\{U_i\rightarrow A\}_i$ to $B$ and $C$ to obtain a diagram of dg-manifolds \[\begin{tikzcd}\coprod_i(U_i\times_A B,\Of_{B\bullet}|_{U_i\times_A B})\ar[d] \ar[r] &\coprod_i(U_i,\Of_{A\bullet}|_{U_i}) \ar[d] & \coprod_i(U_i\times_A C,\Of_{C\bullet}|_{U_i\times_A C})\ar[d]\ar[l] \\ (B,\Of_{B\bullet})  \ar[r]& (A,\Of_{A\bullet})& (C,\Of_{C\bullet})\ar[l]\end{tikzcd}\]
which admits a \v{C}ech nerve $\simpopplus\times\Lambda^{2}_0\rightarrow \mathbf{dgMan}$. Taking pullbacks in each simplicial degree determines a diagram $\simpopplus\times\Delta^1\times\Delta^1\rightarrow \mathbf{dgMan}$ which is a \v{C}ech nerve of its restriction to $\{\infty\}\times\Delta^1\times\Delta^1$. Consider the composition 
\[\Phi: \simpopplus\times\Delta^1\times\Delta^1\longrightarrow \mathbf{dgMan}\overset{\Gamma}{\longrightarrow} \sring_{\fp}^{op}\overset{\spec}{\longrightarrow} \dafffp,\]
then we wish to show that $\Phi|_{\{\infty\}\times\Delta^1\times\Delta^1}$ is a pullback diagram. First, we note that $\Phi$ is also a \v{C}ech nerve; this follows from $(*)$ and the fact that the composition $j\circ \spec\circ \Gamma$ preserves countable coproducts of dg-manifolds (that exist in dg-manifolds). In particular, the adjoint map $\simpopplus\rightarrow\fun(\Delta^1\times\Delta^1,\dstack_{\lfp})$ is a colimit diagram. Using $(*)$ again, we see that both transformations $\Phi|_{\simpopplus\times\{\Delta^1\}\times\{1\}}$ and $\Phi|_{\simpopplus\times\{1\}\times\{\Delta^1\}}$ are Cartesian transformations and therefore realization fibrations of simplicial objects, so it suffices to show that for each $[n]\in \simp$, the diagram $\Phi|_{\{[n]\}\times \Delta^1\times\Delta^1\}}$ is a pullback. This is a square diagram each object of which is a coproduct of representables. Using that coproducts are disjoint in \inftopoit, we reduce to the case where $A$ is a Cartesian space (recall that the cover $\{U_i\rightarrow X\}_i$ is good) and $C$ admits a cover $\{V_j\rightarrow C\}$ such that $V_j\rightarrow A$ is equivalent to a projection $\R^n\times A\rightarrow A$. Consider the \v{C}ech nerve of the diagram 
\[\begin{tikzcd}(B,\Of_{B\bullet})\ar[d] \ar[r] &(A,\Of_{A\bullet}) \ar[d] & \coprod_j(V_j,\Of_{C\bullet}|_{V_j})\ar[d]\ar[l] \\ (B,\Of_{B\bullet})  \ar[r]& (A,\Of_{A\bullet})& (C,\Of_{C\bullet}).\ar[l]\end{tikzcd}\]
Using the same argument we just employed we reduce to the case where the map $C\rightarrow A$ is diffeomorphic to a surjective projection $V \rightarrow \R^m$ from some open $V\subset\R^{n+m}$. We may choose a cover of $V$ by open balls in $\R^{n+m}$, then we apply the preceding argument again to reduce to the case of a projection $\R^{n+m}\rightarrow\R^{m}$ followed by an open embedding $\R^m\hookrightarrow \R^m$. Using $(*)$ again, we may assume that $V\rightarrow A$ is a projection $\R^{n+m}\rightarrow\R^m$. Now we apply this argument yet again to a good open cover of $B$ to reduce to the case of $B$ also being some Cartesian space. It 
 now suffices to argue that the diagram
\[
\begin{tikzcd}
\Gamma(\Of_{A\bullet}) \ar[d] \ar[r] & \Gamma(\Of_{B\bullet})\ar[d] \\
\Gamma(\Of_{C\bullet}) \ar[r] & \Gamma(\Of_{D\bullet})
\end{tikzcd}
\]
is a homotopy pushout, but in this diagram, all objects are cofibrant and the left vertical map is a cofibration of $\cinfty$dgas.
\end{proof}
We have now shown that the functor $\mathbf{dgMan}\rightarrow\dafffp$ carries homotopy pullback diagrams to pullback diagrams. To prove that this functor induces an equivalence after localizing at the weak equivalences of dg-manifolds, it suffices to argue that it detects weak equivalences and that for any dg-manifold $(X,\Of_{X\bullet})$ and any map $\Gamma(\Of_{X\bullet})\rightarrow A$ in $\sring_{\fp}$, there is a dg-manifold $(Y,\Of_{Y\bullet})$, a map $f:(Y,\Of_{Y\bullet})\rightarrow (X,\Of_{X\bullet})$ and a commuting diagram 
\[
\begin{tikzcd}
& A \ar[dr]\\
\Gamma(\Of_{Y\bullet}) \ar[rr,"\Gamma(f)"] \ar[ur,"\simeq"] && \Gamma(\Of_{X\bullet})
\end{tikzcd}
\]
in $\sring_{\fp}$. We will prove both assertions in part II of this series.
\begin{rmk}
It is important to keep in mind that dg-manifolds present \emph{affine} derived $\cinfty$-schemes of finite presentation. Our representability results in part III for differential geometric moduli problems \emph{do not in general} produce affine derived $\cinfty$-schemes of finite presentation and therefore cannot be modelled by dg-manifolds. The formalism of $\infty$-categorical structured spaces straightforwardly provides a good theory of non-affine derived $\cinfty$-schemes obtained by gluing affine ones in a fully coherent way, and we show that elliptic moduli problems are representable by such objects (if the moduli space in question is compact, then it is representable by an affine derived $\cinfty$-scheme of finite presentation, but to prove this we will need a good theory of non affine derived $\cinfty$-schemes in the first place). If one however insists on using the category of fibrant objects structure on dg-manifolds only, one has to resort to gluing in the \infcat of sheaves on dg-manifolds as we sketched in the introduction. \end{rmk}
\newpage
\section{Digression: flatness results for $\cinfty$-completions and flat ideals}
We have seen that resolving effective epimorphisms by morphisms dual to embeddings of graphs and some elementary properties of smooth functions lead to concrete ways of computing $\cinfty$-tensor products $A\oinfty_B C$, especially if one of the maps involved is an effective epimorphism (Corollary \ref{cor:algeffepi}). This result relates an operation induced by the extra $\cinfty$-structure on our derived rings to the underlying homotopical algebra. In this technical section, we take up several other such problems which arise naturally in derived $\cinfty$-geometry and are central to many constructions that follow in this work and its successors. The ideas in this section come from a variety of classical results on ideals of $\cinfty$-functions due to Whitney, Łojasiewicz, Malgrange and Tougeron \cite{Loja,Malgrange,Tougeron}.\\
Given a simplicial commutative $\R$-algebra, we may ask for a prescription for computing the homotopy groups of the free simplicial $\cinfty$-ring $F^{\cinfty}(A)$ on $A$. For instance, if $A$ is discrete, it is not a priori clear that $F^{\cinfty}(A)$ will be a discrete simplicial $\cinfty$-ring as well. We will show the following result.
\begin{prop}\label{prop:freecinftycommutestruncation}
Let $A$ be a simplicial commutative $\R$-algebra, then the unit map $A\rightarrow F^{\cinfty}(A)^{\rmalg}$ is flat (see Definition \ref{def:flatmorphism}) and thus induces for all $n\geq 0$ an equivalence
\[  \pi_n(A)\otimes_{\pi_0(A)}\pi_0(F^{\cinfty}(A)^{\rmalg})\simeq \pi_n(F^{\cinfty}(A)^{\rmalg}). \]
In particular, the Beck-Chevalley transformation $F^{\cinfty} \circ i_n\rightarrow i_n\circ F_n^{\cinfty}$ is an equivalence, where $i_n$ denotes the inclusions $\tau_{\leq n}\sring\subset\sring$ and $\tau_{\leq n}\scring_{\R}\subset\scring_{\R}$.
\end{prop}
We have already used this result to construct a model structure on unbounded $\cinfty$dgas; it will prove useful again in part II, where we transfer model structures from algebras with some extra structure (like a derivation, or a mixed graded structure in the sense of \cite{PTVV1,CPTVV1}) to $\cinfty$-algebras with extra structure. We will prove Proposition \ref{prop:freecinftycommutestruncation} momentarily. Another obvious question stems from the observation that Corollary \ref{cor:algeffepi} tells us very little about the homotopy groups of the coproduct $A\oinfty B$ of two simplicial $\cinfty$-rings, since the map $\R\rightarrow A$ is an effective epimorphism if and only if $\specr\,A$ is a single point. Proposition \ref{prop:freecinftycommutestruncation} gives a description of the homotopy groups of $A\oinfty B$ when $A$ and $B$ lie in the essential image of the functor $F^{\cinfty}$; indeed, in that case, we have isomorphisms 
\[\pi_n(A\oinfty B)\cong \pi_n(A\otimes B)\otimes_{\pi_0(A)\otimes\pi_0(B)}\pi_0(A\oinfty B)\] 
for all $n\geq 0$. It may then seem reasonable to expect that for any pair of simplicial $\cinfty$-rings $A$ and $B$, the canonical map
\begin{equation}\label{eq:coproductstrongmorphism} \pi_n(A\otimes B)\otimes_{\pi_0(A)\otimes\pi_0(B)}\pi_0(A\oinfty B)\longrightarrow \pi_n(A\oinfty B)\end{equation}
is an isomorphism. This assertion however is equivalent to an open problem in differential geometry.

\begin{prop}\label{prop:projectionflatconsequences}
The following are equivalent. 
    \begin{enumerate}[$(a)$]
        \item For any $m,n\in\Z_{\geq 1}$, the map $\cinfty(\R^{n})\otimes \cinfty(\R^m)\rightarrow\cinfty(\R^{n+m})$ induced by the projections $\R^{n+m}\rightarrow\R^n$ and $\R^{n+m}\rightarrow\R^m$ onto the first $n$ coordinates and the last $m$ coordinates respectively is a flat map of commutative $\R$-algebras.
        \item For any small diagram $\mathcal{J}:K\rightarrow \sring$, the canonical comparison map 
         \[ \underset{K}{\colim} \mathcal{J}^{\rmalg}\longrightarrow  (\underset{K}{\colim} \mathcal{J})^{\rmalg}   \]
        is a flat map of simplicial commutative $\R$-algebras. 
    \end{enumerate}
Furthermore, the following are equivalent.
     \begin{enumerate}[$(a)$]
         \item For any $m,n\in\Z_{\geq 1}$, the map $\cinfty(\R^{n})\rightarrow\cinfty(\R^{n+m})$ induced by the projection $\R^{n+m}\rightarrow\R^n$ onto the first $n$ coordinates is a flat map of commutative $\R$-algebras.
         \item For any simplicial $\cinfty$-ring $A$ and any $m\in \Z_{\geq 1}$, the map $A^{\rmalg}\rightarrow (\cinfty(\R^m)\oinfty A)^{\rmalg}$
         is a flat map of simplicial commutative $\R$-algebras.
         \item For any simplicial $\cinfty$-ring $A$ and any $m\in \Z_{\geq 1}$, the map $A^{\rmalg}\rightarrow (\cinfty(\R^m)\oinfty A)^{\rmalg}$ is strong, that is, the canonical map 
         \[ \pi_n(A)\otimes_{\pi_0(A)}\pi_0(A\oinfty \cinfty(\R^m))\longrightarrow \pi_n(A\oinfty\cinfty(\R^n))  \]
         is an isomorphism for all $n\geq 0$.
         \item For any $m\in\Z_{\geq 1}$ and any finitely generated ideal $I\subset \cinfty(\R^n)$, the first homotopy group of the coproduct $\cinfty(\R^n)/I\oinfty \cinfty(\R^m)$ taken in $\sring$ vanishes.
     \end{enumerate} 
\end{prop}

Clearly, the equivalent conditions of $(1)$ imply those of $(2)$. If we could establish the veracity of the conditions in Proposition \ref{prop:projectionflatconsequences}, we would also decide the question of \emph{left properness} of the model category structure on $\cinfty\mathbf{dga}^{\geq0}$ in the positive. Unfortunately, we haven't so far been able to prove that either of the conditions in this proposition are true, or provide a counterexample. There are, of course, specific classes of $0$-truncated simplicial $\cinfty$-rings the objects of which remain $0$-truncated under the coproduct in $\sring$, such as $\cinfty$-rings of functions on manifolds. Our next result asserts that the much larger class of simplicial $\cinfty$-rings of Whitney functions also has this property, generalizing the theorem of Reyes-Van Qu\^{e} to the derived setting.
\begin{thm}\label{thm:reyesvanquederived}
The class of discrete simplicial $\cinfty$-rings of Whitney functions is closed under finite coproducts in $\sring$: let $X\subset\R^n$ and $Y\subset\R^m$ be closed subsets, and let $C^{\infty}(X;\R^n)$ and $C^{\infty}(Y;\R^m)$ be the discrete simplicial $C^{\infty}$-rings of Whitney functions on $X$ and $Y$ respectively, then the canonical map
\[C^{\infty}(X;\R^n)\otimes^{\infty}C^{\infty}(Y;\R^m) \longrightarrow C^{\infty}(X\times Y;\R^{n+m}) \]
induced by the projections $X\times Y\rightarrow X$ and $X\times Y\rightarrow Y$ is an equivalence, where the tensor product is the coproduct of simplicial $C^{\infty}$-rings.
\end{thm}
We record a few consequences of Theorem \ref{thm:reyesvanquederived} and give its proof at the end of this subsection. As we will see, this result is essential and stands at the basis of all further conceptual and computational results involving Whitney functions in this work and its successors. In particular, Theorem \ref{thm:reyesvanquederived} makes possible the development of derived logarithmic $\cinfty$-geometry and derived $\cinfty$-geometry with corners, as it guarantees that the notion of the `subspace of positive elements' of a simplicial $\cinfty$-ring is well behaved (it is canonically endowed with the structure of homotopy coherent commutative monoid, or a \emph{$\Gamma$-object} in the sense of Segal, as we will show later). Recall that the category $\cartsp_{c}$ of Cartesian spaces with corners has as objects the Cartesian spaces with corners $\R^n\times\R^k_{\geq 0}$ and as morphisms the interior $b$-maps. Let $X\subset \R^n$ and $Y\subset \R^m$ correspond to closed quadrants of the form $\R^k\times \R^{n-k}_{\geq 0}$ and $\R^l\times \R^{m-l}_{\geq 0}$, then Theorem \ref{thm:reyesvanquederived} implies that the composition
\[ \cartsp_c\overset{j}{\longrightarrow} \sring_{c} \overset{\iota_c^*}{\longrightarrow}\sring \]
preserves coproducts, where $\iota_c$ is the fully faithful morphism of Lawvere theories $\cartsp\rightarrow\cartsp_c$. As $\iota_c^*$ is a sifted colimit completion of $\iota_c$, we find that $\iota_c^*$ preserves all small colimits which then implies the following result. 
\begin{cor}\label{cor:cornerforgetcolimitpreserve}
Consider the \infcat $\sring_{pc}$ of simplicial $C^{\infty}$-ring with pre-corners and the adjunction
\[ \begin{tikzcd} \sring \ar[r,shift left,"\iota_{c!}"] & \sring_{pc}. \ar[l,shift left,"\iota^*_c"]\end{tikzcd} \]
The functor $\iota^*_c$ carrying a simplicial $C^{\infty}$-ring with pre-corners to the underlying simplicial $C^{\infty}$-ring is a left adjoint. The right adjoint is given by the functor $\iota_{c*}:\sring\rightarrow\sring_{pc}$ obtained by adjunction from the functor
\begin{equation}\label{eq:cornerfunctor} \begin{tikzcd}
\cartsp_c^{op}\times \sring \ar[r,"j^{op}\times\mathrm{id}"]&[1,2em] \sring_{pc}^{op}\times\sring  \ar[r,"\iota_c^*\times\mathrm{id}"] &[1,2em] \sring^{op}\times\sring \ar[r,"\Hom_{\sring}"]&[1,2em] \spa
\end{tikzcd}\end{equation}
which is on objects given by the formula
\[ \iota_{c*}(A)(\R^k\times \R^{n-k}_{\geq 0}) =\Hom_{\sring}(\cinfty(\R^k\times \R^{n-k}_{\geq0}),A).\]
Moreover, both $\iota_{c!}$ and $\iota_{c*}$ are fully faithful.
\end{cor}
\begin{proof}
The existence of a right adjoint to $\iota_c^*$ is a consequence of the adjoint functor theorem. The composition $\pshv(\cartsp_c)\overset{L}{\rightarrow} \sring_{pc}\overset{\iota_c^*}{\rightarrow}\sring$ where $L$ is a left adjoint to the inclusion preserves colimits and is therefore a left Kan extension of $i^*_c\circ j^{op}:\cartsp^{op}_c\rightarrow\sring$ along the Yoneda embedding, and we can identify the functor obtained via adjunction from \eqref{eq:cornerfunctor} as a right adjoint to $\iota_c^*\circ L$. Since this right adjoint factors through $\sring_{pc}$ by Theorem \ref{thm:reyesvanquederived}, it is also right adjoint to $\iota^*_c$. Note that $\iota_{c!}$ is the sifted colimit completion of a fully faithful functor and is thus itself fully faithful. It follows from formal nonsense that $\iota_{c*}$ is also fully faithful.  
\end{proof}
This last corollary implies that the assignment carrying a simplicial $\cinfty$-ring $A$ to the \infcat of pairs $(A,A_c)$ is functorial.
\begin{prop}\label{prop:cornerpresfib}
The functor $\iota_c^*:\sring_{pc}\rightarrow \sring$ is a presentable fibration.
\end{prop}
\begin{proof}
Clearly, $\iota_c^*$ is a categorical fibration, so it suffices to show that $\iota_c^*$ is a Cartesian and coCartesian fibration with presentable fibres. Invoking Corollary \ref{cor:cornerforgetcolimitpreserve} and applying Proposition \ref{prop:leftadjseccocart} and its dual to $\iota^*_c$, we deduce that $\iota^*_c$ is a Cartesian and coCartesian fibration. For the assertion regarding presentability, we first note that the fibres of $\iota^*_c$ are accessible as the \infcat of accessible \infcats and accessible functors between them is stable under pullbacks in $\catinfh$ (note that the functor $\Delta^0\hookrightarrow \icat$ classifying some object $C\in\icat$ preserves colimits of weakly contractible diagrams for any \infcat $\icat$; in particular, this functor is $\kappa$-continuous for any regular cardinal $\kappa$). The presentability of the fibres now follows from the following formal argument.
\begin{enumerate}
    \item[$(*)$] Let $p:\icat\rightarrow\icatd$ be a coCartesian fibration among \infcats and $K$ a simplicial set. Let $f:K\rightarrow\icat_D$ be a diagram in the fibre over some object $D\in\icatd$. Let $i_D:\icat_D\subset\icat$ denote the inclusion, and suppose that the induced diagram $i_Df:K\rightarrow \icat$ admits a colimit and that $p$ preserves the colimit of $i_Df$. Then the diagram $f$ admits a colimit.
\end{enumerate}
We prove $(*)$. Let $C$ denote a colimit of $i_Df$ and denote $D'=p(C)$ so that we have a map $D\rightarrow D'$ for each $k\in K$. Pick one such map $e:D\rightarrow D'$. We have a diagram
\[
\begin{tikzcd}
K\ar[d,hook] \ar[r,"i_Df"] & \icat\ar[d] \\
K^{\rhd} \ar[r] \ar[ur]& \icatd
\end{tikzcd}
\]
wherein the diagonal carries the cone point to $C$. Since the lower horizontal functor $K^{\rhd}\rightarrow\icatd$ is a colimit diagram, the square is also a $p$-colimit diagram. It follows from \cite{HTT}, Proposition 4.3.1.9 that the object $C$ is a $p$-colimit of the diagram $e_!f:K\rightarrow \icat_{D'}\subset\icat$. Since $D'$ is a colimit of the constant diagram with domain $K$ on $D$, there is a map $e':D'\rightarrow D$ such that $e'\circ e\simeq \mathrm{id}_D$, so using \cite{HTT}, Proposition 4.3.1.10, we deduce that $e'_!(C)$ is a colimit of the diagram $e'_! e_!f\simeq f:K\rightarrow\icat_D$.
\end{proof}
We give one final application of Theorem \ref{thm:reyesvanquederived}, answering another question about the interaction of $\cinfty$-geometry and the categorical structure of $\sring$. Theorem \ref{thm:reyesvanquederived} asserts that the class of $\cinfty$-rings of Whitney functions is closed under finite coproducts; we may also ask whether the class of $\cinfty$-rings of Whitney functions on closed sets in a given $\R^n$ is closed under intersections. We verify this is the case for a class of self-intersections.
\begin{prop}\label{prop:whitneyfunctionpushout}
Let $\{X_i\subset\R\}_{i\in I}$ be a finite collection of closed subsets and let $X=\prod_i X_i \subset\R^n$ be their product as a closed subset of $\R^n$, with $n=|I|$. Let $p:\cinfty(\R^n)\rightarrow\cinfty(X;\R^n)$ be the quotient map onto the discrete simplicial $\cinfty$-ring of Whitney functions on $X$. Then the commuting diagram
\[
\begin{tikzcd}
\cinfty(\R^n) \ar[d,"p"] \ar[r,"p"] &\cinfty(X;\R^n) \ar[d,"\mathrm{id}"] \\
\cinfty(X;\R^n) \ar[r,"\mathrm{id}"] & \cinfty(X;\R^n)
\end{tikzcd}
\]
is a pushout in the \infcat $\sring$.
\end{prop}
Using more sophisticated methods from continuous Hochschild homology, we will prove in part II that the diagram above is a pushout for any $X\subset\R^n$ closed. In the interest of keeping this work somewhat self-contained, we have opted to give a more elementary proof here which applies to the special cases we will be primarily interested in, such as $X=\{0\}\times\R^k\subset\R^{n+k}$ or $X=\R^n\times\R^k_{\geq0}\subset\R^{n+k}$.
\begin{rmk}
It follows from Corollary \ref{cor:algeffepi} that the diagram above is also a pushout in $\scring_{\R}$. 
\end{rmk}
For the proof, we recall the following notion.
\begin{defn}
Let $X,Y\subset\R^n$ be closed subsets. The sets $X$ and $Y$ are \emph{regularly situated} if either $X\cap Y=\emptyset$ or for each $x_0\subset X\cap Y$, there is a neighbourhood $x_0\in V$ in $\R^n$ for which there are constants $C\in \R_{>0}$ and $\lambda \in\R_{\geq0}$ such that for each $x\in V\cap X$, we have the inequality 
\[  Cd(x,X\cap Y)^{\lambda}\leq d(x,Y), \]
where $d(\_,\_)$ denotes the Euclidean distance on $\R^n$.
\end{defn}
\begin{ex}
If $X\subset Y$, then $X$ and $Y$ are regularly situated. In particular, two copies of the same set $X$ are regularly situated.
\end{ex}
\begin{ex}
Let $X,Y\subset\R^n$ be subanalytic closed sets, then $X$ and $Y$ are regularly situated. This is proven by Bierstone-Milman \cite{BierMil}.
\end{ex}
\begin{rmk}
In part II, we will give the following characterization of the condition of being regularly situated for $X,Y\subset\R^n$ in terms of the derived intersection of the affine finitely generated $\cinfty$-schemes $(X,\cinfty_{(X;\R^n)})$ and $(Y,\cinfty_{(Y;\R^n)})$: $X$ and $Y$ are regularly situated if and only if the commuting square
\[
\begin{tikzcd}
\cinfty(X\cup Y;\R^n) \ar[d] \ar[r] &\cinfty(X;\R^n) \ar[d] \\
\cinfty(Y;\R^n) \ar[r] & \cinfty(X\cap Y;\R^n)
\end{tikzcd}
\]
of simplicial $\cinfty$-rings is a pullback if and only if it is a pushout.
\end{rmk}
\begin{defn}
Given a closed set $X\subset\R^n$, the space $\mathcal{M}(X;\R^n)$ of smooth functions $f:\R^n\setminus X\rightarrow \R$ that have the property that for any compact $K\subset \R^n$ and any multi-index $k\in \Z_{\geq 0}^{n}$, there exist constants $C,\alpha\in \R_{>0}$ such that for each $x\in K\setminus K\cap X$ the inequality 
\[ |D^{k}(f)(x)| \leq Cd(x,X)^{-\alpha}  \]
is satisfied, is the space of \emph{multipliers for the ideal $\liem^{\infty}_X$}: for any $\varphi\in \mathcal{M}(X;\R^n)$ and any $f\in \liem^{\infty}_X$, the function $f\varphi$ defined on $\R^n\setminus X$ uniquely extends to a $\cinfty$-function (still denoted $f\varphi$) on $\R^n$ that is flat on $X$.  
\end{defn}
We will require the following result.
\begin{lem}[Tougeron's Multiplier Lemma]\label{lem:tougeronmultiplier}
Let $X,Y\subset\R^n$ be closed and regularly situated, then there exists a multiplier $\varphi$ for the ideal $\liem^{\infty}_{X\cap Y}$ that equals $0$ in a neighbourhood of $X\setminus X\cap Y$ and equals $1$ in a neighbourhood of $Y\setminus X\cap Y$. 
\end{lem}
\begin{proof}
Lemme 4.5 of \cite{Tougeron}.
\end{proof}
\begin{proof}[Proof of Proposition \ref{prop:whitneyfunctionpushout}]
Applying Theorem \ref{thm:reyesvanquederived}, we may assume that $n=1$. It is obvious that the diagram in the statement of the proposition is a pushout after applying the 0'th truncation functor $\tau_{\leq 0}$, so it suffices to argue that the higher homotopy groups vanish. Since $\sring$ is a coCartesian symmetric monoidal \infcatt, the pushout $\cinfty(X;\R )\oinfty_{\cinfty(\R)}\cinfty(X;\R)$ is a colimit of the two sided Bar construction $\mathsf{Bar}_{\cinfty(\R)}(\cinfty(X;\R),\cinfty(X;\R))_{\bullet}$, the simplicial object 
\[\begin{tikzcd} \ldots\ar[r,shift left=3] \ar[r,shift left] \ar[r,shift right]\ar[r,shift right=3] & \cinfty(X)\oinfty \cinfty(\R)^{\oinfty 2} \oinfty \cinfty(X) \ar[r,shift left=2] \ar[r]\ar[r,shift right=2] & \cinfty(X)\oinfty \cinfty(\R) \oinfty \cinfty(X) \ar[r,shift left]\ar[r,shift right] & \cinfty(X)\oinfty \cinfty(X). \end{tikzcd}  \]
It follows from Theorem \ref{thm:reyesvanquederived} that $\mathsf{Bar}_{\cinfty(\R)}(\cinfty(X;\R),\cinfty(X;\R))_{k}\simeq  \cinfty(X\times\R^{k}\times X)$ and the face maps are induced by the various inclusions of small diagonals $X\times\R^{m}\times X\hookrightarrow X\times\R^{k}\times X$ for $m<k$. As geometric realizations are sifted, the colimit $|\mathsf{Bar}_{\cinfty(\R)}(\cinfty(X;\R),\cinfty(X;\R))_{\bullet}|$ may be computed in the \infcat $\Mod_{\R}$, where it becomes a geometric realization of a simplicial object in the heart. By the stable Dold-Kan correspondence, the homotopy groups of $|\mathsf{Bar}_{\cinfty(\R)}(\cinfty(X;\R),\cinfty(X;\R))_{\bullet}|$ as $\R$-vector spaces are computed by the spectral sequence associated to the filtered object determined by $\mathsf{Bar}_{\cinfty(\R)}(\cinfty(X;\R),\cinfty(X;\R))_{\bullet}$ now viewed as a simplicial object in $\R$-vector spaces, which collapses at the second page to the homology of the unnormalized chain complex $C(\mathsf{Bar}_{\cinfty(\R)}(\cinfty(X;\R),\cinfty(X;\R)))_{\bullet}$. It thus suffices to show that the higher homology groups of the \emph{normalized} chain complex $N(\mathsf{Bar}_{\cinfty(\R)}(\cinfty(X;\R),\cinfty(X;\R)))_{\bullet}$ vanish. This will be accomplished by constructing for each cycle in degrees $\geq 1$ an explicit boundary. Unraveling the definitions, we need to show the following.
\begin{enumerate}
\item[$(*)$] Let $k\geq 1$ and let $F(x,z_1,\ldots,z_k,y)$ be a Whitney function on $X\times\R^{k} \times X$ such that for all $1< j \leq k$, the Whitney function $F(x,z_1,\ldots,z_j,z_j,
\ldots,z_{k-1},y)$ on $X\times\R^{k-1}\times X$ vanishes (if $j=k$, then we have $F(x,z_1,\ldots,z_{k-1},y,y)$ i.e. we restrict the penultimate coordinate to $X$). If $F(x,x,z_{2},\ldots,z_k,z)$ also vanishes, then there exists a Whitney function $\widehat{F}(x,z_1,\ldots,z_{k+1},y)$ on $X\times\R^{k+1}\times X$ such that for all $1\leq j\leq k+1$, the Whitney function $\widehat{F}(x,z_1,\ldots,z_j,z_j,
\ldots,z_k,y)$ vanishes and $\widehat{F}(x,x,z_2,\ldots,z_{k+1},y)=F$. \end{enumerate}
Let $F$ be a cycle of degree $k$ and let $\Delta_i\subset\R\times\R^{k}\times \R$ for $0\leq i\leq k$ denote the small diagonal
\[ \Delta_i :=\{(x_0,z_1,\ldots,z_k,z_{k+1})\in \R\times\R^{k}\times \R;\,z_i=z_{i+1}\}.\]
Let $f$ be any representative of $F$. Since $F$ lies in the joint kernel of all face maps $d_i$ for $0<i\leq k$, the restriction $f|_{\Delta_i}$ is flat on $X\times\R^{k-1}\times X$ for all $0<i\leq k$, and since the differential of the normalized chain is the 0'th face maps, we see that $f|_{\Delta_0}$ is also flat on $X\times\R^{k-1}\times X$. The inclusion $\Delta_0\subset\R\times\R^{k}\times\R$ admits a smooth deformation retraction $r$ defined by 
\[ (z_0,z_1,\ldots,z_k,z_{k+1}) \longmapsto  (1/2(z_0+z_{1}),z_{2},\ldots,z_{k},z_{k+1}). \]
Pulling back functions along the composition 
\[ \R\times\R^{k}\times\R \overset{r}{\longrightarrow} \Delta_0 \hooklongrightarrow \R\times\R^{k}\times\R \]
yields an operator $r^*(\_)|_{\Delta_0}:\cinfty(\R\times\R^{k}\times\R)\rightarrow \cinfty(\R\times\R^{k}\times\R)$ that we denote $(\_)_{\Delta_0}$. By the vanishing properties of $F$, the function $f_{\Delta_0}$ is flat on the closed subset 
\[ (X\times\R^{k}\times X) \cap \Delta_0 \cong X\times\R^{k-1}\times X. \]
Now we claim that the sets $\Delta_0$ and $X\times\R^{k}\times X$ are regularly situated: indeed, for any $p=(z_0,z_1,\ldots,z_k,z_{k+1})$, the distance $d(p,\Delta_0)$ is $1/\sqrt{2}d(z_0,z_1)$ but if $p\in X\times\R^{k}\times X$, then $d(p,(X\times\R^{k}\times X) \cap \Delta_0)$ is also $1/\sqrt{2}d(z_0,z_1)$ which immediately implies that two sets in question are regularly situated. Tougeron's multiplier lemma provides a function $\varphi$ on $\R\times\R^{k}\times\R\setminus (X\times\R^{k}\times X) \cap \Delta_0$ that is $1$ in a neighbourhood of $\Delta_0\setminus (X\times\R^{k}\times X) \cap \Delta_0$ and $0$ in a neighbourhood of $X\times\R^{k}\times X\setminus (X\times\R^{k}\times X) \cap \Delta_0$. The function $\varphi f_{\Delta_0}$ is then a $\cinfty$ function on $\R\times\R^{k}\times\R$ that is flat on $X\times\R^{k}\times X$ and equals $f$ on $\Delta_0$. Now consider the function $\tilde{f}:=f-\varphi f_{\Delta_0}$, then the Whitney jet of $\tilde{f}$ is $F$ and $\tilde{f}$ vanishes along $\Delta_0$. It follows from Hadamard's Lemma that $\tilde{f}$ may be written as 
\[ \tilde{f} = (z_0-z_1)g(z_0,z_1\ldots,z_k,z_{k+1}).\]
It follows from the construction of the Hadamard quotient $g$ that $g|_{\Delta_i}$ is flat on $X\times\R^{k-1}\times X$ for $i\geq 1$. Now define the function $\widehat{f}:\R\times\R^{k+1}\times\R\rightarrow\R$ via the formula 
\[ (z_0,z_1,\ldots,z_{k+1},z_{k+2}) \longmapsto (z_1-z_2)g(z_0,z_2,z_3,\ldots,z_{k+1},z_{k+2}), \]
then one readily verifies that the Whitney jet of $\widehat{f}$ at $X\times\R^{k+1}\times X$ is a cycle of degree $k+1$ the boundary of which is $F$.
\end{proof}
\begin{rmk}\label{rmk:truncatepositivelts}
It follows rather trivially from Proposition \ref{prop:whitneyfunctionpushout} that the $\cinfty$ theoretical cotangent complex $\cotan_{\cinfty(X;\R^n)}$ for $X=\prod_i X_i$ with $X_i\subset\R$ closed (like rings of smooth functions on closed quadrants of the form $\cinfty(\R^n\times\R^k_{\geq0})$) is \emph{free on $n$ generators}. Indeed, for a pushout diagram 
\[
\begin{tikzcd}
A\ar[d] \ar[r] & B\ar[d] \\
C\ar[r] & D
\end{tikzcd}
\]
of simplicial $\cinfty$-rings, we have a canonical equivalence $\cotan_{B/A}\otimes_{B}D\simeq \cotan_{D/C}$. Since for any $A\in\sring$ the $A$-module $\cotan_{A/A}$ vanishes and $\cotan_{\cinfty(\R^n)}$ is free on $n$ generators, Proposition \ref{prop:whitneyfunctionpushout} guarantees that $\cotan_{\cinfty(X;\R^n)}$ is also free on $n$ generators. This should be viewed as an articulation of the idea that seen through the lens of deformation theory, the simplicial $\cinfty$-rings $\cinfty(X;\R^n)$ and $\cinfty(\R^n)$ are equivalent. We will give a number of consequences of this result in part II (where we also show that the preceding remarks are valid for any closed set $X\subset\R^n$); for instance, we will deduce that the map $\Hom(\cinfty(X;\R),A) \rightarrow \Hom(\cinfty(\R),A)$ of spaces is an inclusion of connected components, and if a map $A\rightarrow B$ of simplicial $\cinfty$-rings exhibits $B$ as an $m$-truncation of $A$, then the induced map
\[ \Hom_{\sring}(\cinfty(X;\R^n),A)\longrightarrow \Hom_{\sring}(\cinfty(X;\R^n),B) \]
exhibits an $m$-truncation of spaces. In fact, there unfortunately are several instances in this work where we use that $\cinfty(X;\R^n)$ has a free cotangent complex. Since a detailed discussion of the cotangent complex is not in order at this point, we ask the reader to tolerate a small amount of nonlinear logical interdependency and recognize that no circular reasoning occurs.
\end{rmk}
We proceed with the proof of Proposition \ref{prop:freecinftycommutestruncation}.
We have need of the following classical result of Malgrange.
\begin{prop}[Malgrange \cite{Malgrange}]
Let $\{f_1,\ldots,f_n\}$ be a collection of real analytic functions on $\R^n$, then the finitely generated ideal $(f_1,\ldots,f_n)$ is closed.  
\end{prop}
As an immediate corollary, we have the following.
\begin{cor}\label{cor:malgrangeflat}
Let $x\in \R^n$ and let $\Of_x^{\mathrm{an}}$ be the local ring of real analytic functions on $R^n$ at $x$, then the local morphism of local $\R$-algebras $\Of^{\mathrm{an}}_x\rightarrow \cinfty(\R^n)_x$ is faithfully flat.
\end{cor}

\begin{prop}\label{prop:flatnessofcompletion}
For every integer $n\geq 0$, the map of $\R$-algebras
\[  \R[x_1,\ldots,x_n]\longrightarrow \cinfty(\R^n) \]
determined by the $n$ coordinate functions is flat.
\end{prop}
\begin{proof} Clearly, we may suppose that $n\geq 1$. Consider the composition 
\[ \varphi:\R[x_1,\ldots,x_n]\longrightarrow \cinfty(\R^n) \longrightarrow \prod_{x\in\R^n}  \cinfty(\R^n)_x \]
where the second map is induced by the quotient maps $\cinfty(\R^n)\rightarrow \cinfty(\R^n)/\liem^g_x\simeq \cinfty(\R^n)_x$ sending a smooth function to its germ at $x$ as $x$ ranges over $\R^n$. We now prove the proposition under the assumption that the map $\varphi$ is flat. Let $I$ be an ideal of $\R[x_1,\ldots,x_n]$, then we should show that the top horizontal map in the commuting diagram 
\[
\begin{tikzcd}
\tor_0^{\R[x_1,\ldots,x_n]}(I, \cinfty(\R^n)) \ar[d] \ar[r] & \cinfty(\R^n)\ar[d] \\
\tor_0^{\R[x_1,\ldots,x_n]}(I,  \prod_{x\in\R^n} \cinfty(\R^n)_x) \ar[r] &  \prod_{x\in\R^n}  \cinfty(\R^n)_x
\end{tikzcd}
\]
is injective. By assumption, the lower horizontal map is injective, so it suffices to show that the left vertical map is injective. Unwinding the definitions, this amounts to the following assertion.
\begin{enumerate}
    \item[$(*)$] Let $\{P_i\}_{i\in J}$ and $\{f_i\}_{i\in J}$ be finite collections of real polynomials in $n$ variables and smooth functions in $n$ variables respectively. Suppose that there exist a finite index set $K$ together with a $K$-indexed collection $\{Q_{ik}\}$ of real polynomials in $n$ variables for each $i\in J$ such that the following hold.
    \begin{enumerate}[$(a)$]
        \item $\sum_i Q_{ik} P_i=0$ for each $k\in K$.
        \item For each $x\in \R^n$ there exists an open neighbourhood $x\in U_x$ and a $K$-indexed collection $\{g^x_k\}$ of smooth functions on $U_x$ such that $f_i=\sum_k g^x_kQ_{ik}$ on $U_x$.
    \end{enumerate}
    Then there exists a $K$-indexed collection $\{g_k\}$ of smooth functions on $\R^n$ such that $f_i=\sum_k g_kQ_{ik}$.
\end{enumerate}
To prove $(*)$, let $\{\psi_x\}_{x\in\R^n}$ be a partition of unity subordinate to the cover $\{U_x\rightarrow \R^n\}_{x\in \R^n}$, and define $g_k:=\sum_{x \in \R^n}\psi_xg_k^x$, then it is easy to see that $f_i=\sum_k g_kQ_{ik}$ holds for all $i\in J$. We are left to prove that $\varphi$ is a flat morphism. Since $\R[x_1,\ldots,x_n]$ is Noetherian hence coherent, $(6)$ of Proposition \ref{prop:flatproperties} guarantees that it suffices to show that for each $x\in \R^n$, the map $\R[x_1,\ldots,x_1]\rightarrow \cinfty(\R^n)_x$ is flat. As localizations are flat, we only have to show that the local homomorphism $\Of^{\mathrm{reg}}(\R^n)_x\rightarrow \cinfty(\R^n)_x$ is flat, where $\Of^{\mathrm{reg}}(\R^n)_x$ is the local ring of regular functions at $x$. We have a factorization 
\[ \Of^{\mathrm{reg}}(\R^n)_x \longrightarrow  \Of^{\mathrm{an}}(\R^n)_x \longrightarrow \cinfty(\R^n)_x \]
of local ring homomorphisms where $\Of^{\mathrm{an}}(\R^n)_x$ is the local ring of real analytic functions on $\R^n$ at $x$. The first map is a local morphism between Noetherian local rings that becomes an equivalence after formal completion at the maximal ideals and is thus faithfully flat, and the second map is faithfully flat by Corollary \ref{cor:malgrangeflat}.
\end{proof}
\begin{rmk}
Corollary \ref{cor:malgrangeflat} can also be used to prove the more powerful result that the unit of the relative spectrum functor sending a \emph{derived real analytic space} to the corresponding derived $\cinfty$-scheme is faithfully flat, but we will not develop such a theory in this work.
\end{rmk}
\begin{rmk}
In the proof of Proposition \ref{prop:flatnessofcompletion} we use that for every $x\in\R^n$ the map $\Of^{\mathrm{reg}}_x\rightarrow\cinfty(\R^n)_x$ is faithfully flat; beware however that the map $\R[x_1,\ldots,x_n]\rightarrow \cinfty(\R^n)$ is \emph{not} faithfully flat because $\R$ is not algebraically closed. The maximal ideals of $\R[x_1,\ldots,x_n]$ with residue field $\C$, such as $(x_1^2+1,x_2,\ldots,x_n)$, have the property that multiplying such an ideal with the module $\cinfty(\R^n)$ recovers all of $\cinfty(\R^n)$ and therefore do not lie in the image of the induced map on maximal ideal spectra.
\end{rmk}
The following proposition is one of many consequences of the resolution theorem for effective epimorphism. 
\begin{prop}\label{prop:freecinftypushout}
Let $f:A\rightarrow B$ be an effective epimorphism of simplicial commutative $\R$-algebras, then the natural diagram
\[
\begin{tikzcd}
A \ar[d] \ar[r] & B\ar[d] \\
F^{\cinfty}(A)^{\rmalg} \ar[r] & F^{\cinfty}(B)^{\rmalg}
\end{tikzcd}
\]
is a pushout in $\scring_{\R}$.
\end{prop}
\begin{proof}
Applying the unit transformation $\scring_{\R}\times\Delta^1\rightarrow\scring_{\R}$ of the adjunction $(F^{\cinfty}\adj (\_)^{\rmalg})$ yields a functor
\[\fun(\Delta^1,\scring_{\R}) \longrightarrow \fun(\Delta^1\times \Delta^1,\scring_{\R})  \]
carrying a map $A\rightarrow B$ to the diagram 
\[
\begin{tikzcd}
A \ar[d] \ar[r] & B\ar[d] \\
F^{\cinfty}(A)^{\rmalg} \ar[r] & F^{\cinfty}(B)^{\rmalg}.
\end{tikzcd}
\]
As the comonad $F^{\cinfty}(\_)^{\rmalg}$ preserves sifted colimits, this functor preserves sifted colimits. Since the full subcategory of $\fun(\Delta^1\times \Delta^1,\scring_{\R})$ spanned by pushout diagrams is stable under colimits, Proposition \ref{prop:resolveelementary} asserts that we may suppose that $A\rightarrow B$ is a graph inclusion. We wish to show that the diagram 
\[ 
\begin{tikzcd}
\R[x_1,\ldots,x_n,x_{n+1},\ldots,x_{n+m}] \ar[r] \ar[d] & \R[x_1,\ldots,x_n] \ar[d] \\
\cinfty(\R^{n+m}) \ar[r] & \cinfty(\R^{n})
\end{tikzcd}
\]
is a pushout in $\scring_{\R}$. It is not hard to see that as in the proof of Lemma \ref{diffdiscunramified}, an inductive argument reduces us to the case $m=1$. Then the horizontal maps in the diagram above are induced by the inclusion of the graph of a polynomial $P(\mathbf{x}):\R^n\rightarrow \R$ in $n$ variables. The $\R$-algebra $\R[x_1,\ldots,x_n]$ has a projective resolution $\R[x_1,\ldots,x_n,x_{n+1},\epsilon]$ with $\del \epsilon=x_{n+1}-P(\mathbf{x})$ as a $\R[x_1,\ldots,x_n,x_{n+1}]$-module, so using that the torsion spectral sequence of the pushout collapses at the second page, we see that the homotopy groups of the pushout are given by the homology of the complex $\cinfty(\R^{n+1})[\epsilon]$. Lemma \ref{projresolutiontransverse} asserts that the homology is indeed $\cinfty(\R^{n})$, concentrated in degree zero.
\end{proof}
\begin{proof}[Proof of Proposition \ref{prop:freecinftycommutestruncation}]
The class of simplicial commutative $\R$-algebras $A$ which have the property that the unit map $A\rightarrow F^{\cinfty}(A)$ is flat is stable under filtered colimits, so we may suppose that $A$ is of finite type over $\R$. Invoking Proposition \ref{fgcpt}, we can find an effective epimorphism $\R[x_1,\ldots,x_n]\rightarrow A$, so that Proposition \ref{prop:freecinftypushout} provides a pushout diagram 
\[
\begin{tikzcd}
\R[x_1,\ldots,x_n]\ar[d] \ar[r] & A\ar[d] \\
\cinfty(\R^n) \ar[r] & F^{\cinfty}(A)^{\rmalg}
\end{tikzcd}
\]
of simplicial commutative $\R$-algebras. Since flatness is stable under base change, we are done by Proposition \ref{prop:flatnessofcompletion}.
\end{proof}
\begin{cor}\label{cor:freecinftyflat1}
For every simplicial commutative $\R$-algebra $A$ and for every $\R$-point $p:\pi_0(A)\rightarrow \R$ corresponding to a maximal ideal $\liem_{p}$, the localization of the unit map $A_{\liem_p}\rightarrow F^{\cinfty}(A)^{\rmalg}_{\liem_p}$ is faithfully flat. 
\end{cor}
We proceed with the proof of Proposition \ref{prop:projectionflatconsequences}. 
\begin{proof}[Proof of Proposition {\ref{prop:projectionflatconsequences}}]
We prove the equivalence of the first statements. Note that $(b)\Rightarrow (a)$ is obvious. For the other direction, we choose for each $k\in K_0$ a small set $I_k$, an $I_k$-collection of positive integers $\{n^k_i\}_{i\in I_k}$ and an effective epimorphism $\bigotimes_{i\in I_k}^{\infty}\cinfty(\R^{n^k_i})\rightarrow \mathcal{J}(k)$, then Theorem \ref{thm:unramified2} shows that there is some small set $I$ and an $I$-indexed collection $\{n_i\}_{i\in I}$ of positive integers such that the comparison map $\colim_K\mathcal{J}^{\rmalg}\rightarrow(\colim_K\mathcal{J})^{\rmalg}$ is a pushout of the comparison map 
\[  \bigotimes_{i\in I}\cinfty(\R^{n_i})\longrightarrow \bigotimes^{\infty}_{i\in I}\cinfty(\R^{n_i}), \]
so it suffices to show that this map is flat. Since flat maps are stable under filtered colimits, we may assume that $I$ is finite. By induction, we may assume that $I$ is a set with two elements and we conclude. We prove the equivalence of the second list of statements. For $(a)\Rightarrow (b)$, we repeat the proof of $(1)$. $(b)\Rightarrow (c)$ and $(c)\Rightarrow (d)$ are obvious. For $(d)\Rightarrow (a)$, we note that if $R$ is a commutative ring and $M$ is a (discrete) $R$-module, then $M$ is flat if and only if $\tor_1^R(M,R/I)\cong 0$ for every finitely generated ideal $I\subset R$. Indeed, for any ideal $I\subset R$, we have a fibre sequence $I\rightarrow R\rightarrow R/I$ of discrete $\R$-modules, so using the long exact sequence associated to the fibre sequence 
\[ I\otimes_{R}M\longrightarrow M\longrightarrow R/I\otimes_{R}M  \]
yields an exact sequence
\[ 0\longrightarrow \tor_1^R(M,R/I) \longrightarrow  \tor^R_0(I,M)\longrightarrow  M \longrightarrow \tor_0^R(R/I,M)\]
and the vanishing of $\tor_1^R(M,R/I)$ is equivalent to the injectivity of the map $\tor^R_0(I,M)\rightarrow  M$. Applying this to the $\cinfty(\R^n)$-module $\cinfty(\R^{n+m})$ for $n,m\geq 1$, we deduce that in order to establish $(1)$, it suffices to show that the derived tensor product $\cinfty(\R^{n+m})\otimes_{\cinfty(\R^n)}\cinfty(\R^n)/I$ has vanishing first homotopy group for every finitely generated ideal $I$. We can identify this derived tensor product with the pushout
\[
\begin{tikzcd}
\cinfty(\R^n) \ar[d] \ar[r] & \cinfty(\R^n)/I\ar[d] \\
\cinfty(\R^{n+m}) \ar[r] & \cinfty(\R^{n+m})\otimes_{\cinfty(\R^n)}\cinfty(\R^n)/I
\end{tikzcd}
\]
of simplicial commutative $\R$-algebras. Using Theorem \ref{thm:unramified2}, this pushout is equivalent to the underlying simplicial commutative $\R$-algebra of the coproduct $\cinfty(\R^{n})/I\oinfty\cinfty(\R^m)$ of simplicial $\cinfty$-rings, which we assume has vanishing first homotopy group. 
\end{proof}
Contemplating the problem of flatness of the map $\colim_k\mathcal{J}^{\rmalg}\rightarrow(\colim_{k}\mathcal{J})^{\rmalg}$ is one way to make precise the following question: is the homotopy theory of simplicial $\cinfty$-rings determined by the underlying homotopy theory of simplicial $\R$-algebras, or is it more complicated? The point of the exercise above is to show that the answer to this question turns on the flatness of the family of maps $\{\cinfty(\R^n)\otimes\cinfty(\R^m)\rightarrow\cinfty(\R^{n+m})\}$. We believe that the answer to this question is negative. Recall the following notion.
\begin{defn}
A finitely generated ideal $I=(f_1,\ldots,f_k)\subset \cinfty(\R^n)$ is a \emph{Łojasiewicz ideal} if either of the following equivalent conditions are satisfied.
\begin{enumerate}[$(1)$]
    \item $\liem^{\infty}_{Z(I)}\subset I$.
    \item The function $f_1^2+\ldots+f_k^2$ satisfies \emph{Łojasiewicz' inequality}: for all compact $K\subset \R^n$ there exists a constant $C\in \R_{>0}$ and a constant $\alpha\in \R_{\geq 0}$ such that
    \[ f_1^2(x)+\ldots+f_k^2(x) \geq C d(x,Z(I))^{\alpha},\quad \forall x\in K \]
    where $d(x,Z(I))$ is the Euclidean distance between $x$ and $Z(I)$.
\end{enumerate}
\end{defn}
It follows immediately from characterization $(1)$ and Whitney's spectral theorem that finitely generated closed ideals are Łojasiewicz. In fact a finitely generated ideal of $\cinfty(\R)$ is closed if and only if it is Łojasiewicz, and in this case it can be shown that the ideal in question is necessarily principal. It is possible to show that if $I\subset\cinfty(\R^n)$ is Łojasiewicz, then the map 
 \[ \tor_0^{\cinfty(\R^n)}(I,\cinfty(\R^{n+m}))\longrightarrow\cinfty(\R^{n+m}) \]
is injective.
\begin{conj}
Let $f,g\in\cinfty(\R)$ by functions with the following property: the ideal $(f,g)$ in $\cinfty(\R)$ is \emph{not} a Łojasiewicz ideal. Then $\cinfty(\R)/(h)\oinfty\cinfty(\R)/(g)$ is not $0$-truncated.    
\end{conj}

We now turn to the proof of Theorem \ref{thm:reyesvanquederived}, for which we need to establish a few prelimenary results. In the absence of a general flatness result, we proceed by constructing, for certain ideals of $\cinfty(\R^n)$, small explicit resolutions that are acyclic for the base change induced by a projection $\R^{n+m}\rightarrow \R^n$ onto the first $n$ coordinates. The remainder of the results in this subsection hinges on the following lemma, due to Tougeron \cite{Tougeron} (for a single manifold), and extended to the form below by Reyes-van Qu\^{e} \cite{reyesvanque} (who attribute this generalization to Calder\'{o}n).
\begin{lem}[Tougeron's Flat Function Lemma]\label{lem:flatfunctionLemma}
Let $X\subset M$ and $Y\subset N$ be closed subsets of manifolds $M$ and $N$. Let $I$ be a countable set and let $\{\phi_i\}_{i\in I}$ be a set of functions on $M\times N$ that lie in $\liem^{\infty}_{X\times Y}$, that is, functions that are flat on $X\times Y$. Then there exists a characteristic function $\varphi_X$ for $M\setminus X$ and a characteristic function $\varphi_Y$ for $N\setminus Y$ that are flat on $X$ and $Y$ respectively such that the functions $\{\phi_i\}_{i\in I}$ are divisible by $\varphi_X+\varphi_Y$ and $\frac{\phi_i}{\varphi_X+\varphi_Y}$ is flat on $X\times Y$ for all $i\in I$. 
\end{lem}
\begin{proof}
See Theorem 4.11 of \cite{MR}.
\end{proof}
The following lemma shows that flat ideals in $C^{\infty}$-rings of smooth functions on manifolds behave for many purposes just as principal ideals.
\begin{lem}\label{lem:cofinalprincipalideal}
Let $M$ and $N$ be manifolds and let $X\subset M$ be a closed subset in a manifold. Denote by $I:=\liem_{X\times N}^{\infty}$ the closed ideal of functions flat on $X\times N$ viewed as a $C^{\infty}(M\times N)$-module, and let $K\subset \mathrm{Sub}_{\mathrm{fg}}(I)$ be the full subcategory of the filtered poset of finitely generated ideals contained in $I$ spanned by principal ideals contained in $I$ generated by functions depending only on coordinates in $M$. Then $K$ is filtered and the inclusion $K\subset\mathrm{Sub}_{\mathrm{fg}}(I)$ is left cofinal. 
\end{lem}
\begin{rmk}
As the $C^{\infty}(M\times N)$-module $\liem_{X\times N}^{\infty}$ is a colimit of the diagram $\mathrm{Sub}_{\mathrm{fg}}(I)$, it follows that for every closed $X\subset M$, the diagram $K^{\rhd}\rightarrow \mathrm{Mod}_{C^{\infty}(M\times N)}$ sending the cone to $\liem_{X\times N}^{\infty}$ is a colimit diagram.
\end{rmk}
\begin{proof}
We prove that the inclusion is left cofinal. According to \cite{HTT}, Theorem 4.1.3.1, we need to show that the poset $K_{J/}:=K\times_{\mathrm{Sub}_{\mathrm{fg}}(I)}\mathrm{Sub}_{\mathrm{fg}}(I)_{J/}$ is weakly contractible for every finitely generated ideal $J\subset I$. It suffices to show that $K_{J/}$ is filtered. Let $\{f_1,\ldots,f_n\}$ be a collection of functions that are flat on $X$ generating the ideal $J$, then it follows from Tougeron's flat function lemma that there exists a function $\varphi_X$ flat on $X$ and strictly positive outside $X$ on $M$ that divides each $f_i$ as a function on $M\times N$, so that $(f_1,\ldots,f_n)\subset (\varphi_X)\subset I$, that is, $K_{J/}$ is nonempty. Similarly, if we have a finite collection of functions $\{\varphi_X^j\}$ such that $J\subset (\varphi_X^j)$ for all $j$, then we apply the flat function lemma to the collection $\{\varphi_X^j\}$ to find a function $\varphi'_X$ such that $(\varphi_X^j)\subset (\varphi_X')$ for all $j$, so $K_{J/}$ is indeed filtered. The same argument shows that $K$ itself is filtered, which concludes the proof. \end{proof}
Let $f:M\rightarrow N$ be a submersion and let $X\subset N$ be a closed set, which determines the closed ideal sheaf $\liem^{\infty}_X\subset\cinfty_N$. In the following lemma, we employ a bootstrapping argument based on the previous lemma and the fact that for a principal ideal $(h)\subset\cinfty(N)$, the kernel of the generator $\cinfty(\R^n)\rightarrow (h)$ is again an ideal of the form $\liem^{\infty}_Y$ to establish the local vanishing of the higher homotopy groups $\pi_n(\cinfty_M\otimes_{\cinfty_N}\liem^{\infty}_X)$.
\begin{lem}\label{lem:acyclicideals}
Let $I$ belong to either of the following classes of ideals of $C^{\infty}(\R^n)$.
\begin{enumerate}[$(1)$]
    \item Principal ideals.
    \item Ideals of the form $\liem_{X}^{\infty}$ for $X\subset\R^n$ closed.
\end{enumerate}
Let $m\in \Z_{\geq 0}$ and let $Y\subset\R^m$ be any closed set, then as a $\cinfty(\R^n)$-module, $I$ is $\tor_0^{C^{\infty}(\R^n)}(\_,C^{\infty}(\R^{n+m})/\liem^{\infty}_{\R^n\times Y})$-acyclic and the map 
\[ \tor_0^{C^{\infty}(\R^n)}(I,C^{\infty}(\R^{n+m})/\liem^{\infty}_{\R^n\times Y})\longrightarrow \cinfty(\R^{n+m})/\liem^{\infty}_{\R^n\times Y}  \]
is a monomorphism, where the base change is along the composition
\[ \cinfty(\R^n)\longrightarrow\cinfty(\R^{n+m})\longrightarrow\cinfty(\R^{n+m})/\liem^{\infty}_{\R^n\times Y}, \]
the first map being induced by the projection $\R^{n+m}\rightarrow \R^n$.
\end{lem}
\begin{proof}
\begin{enumerate}[$(1)$]
    \item Let $g\in C^{\infty}(\R^n)$ be nonzero, then $hg=0$ if and only if $h\in \liem^0_{\mathrm{Supp}(g)}$. Since $\mathrm{Supp}(g)\subset \overline{\mathrm{Supp}(g)^{\circ}}$, we have the equality $\liem^0_{\mathrm{Supp}(g)}=\liem^{\infty}_{\mathrm{Supp}(g)}$, which establishes the fibre sequence
    \[ \liem_{\mathrm{Supp}(g)}^{\infty}\longrightarrow C^{\infty}(\R^n) \overset{1\mapsto g}{\longrightarrow} (g) \]
of discrete $C^{\infty}(\R^n)$-modules. From this fibre sequence, we obtain a fibre sequence 
  \[ \liem_{\mathrm{Supp}(g)}^{\infty}\otimes_{C^{\infty}(\R^n)}\cinfty(\R^{n+m})/\liem^{\infty}_{\R^n\times Y}\longrightarrow \cinfty(\R^{n+m})/\liem^{\infty}_{\R^n\times Y} \longrightarrow (g)\otimes_{C^{\infty}(\R^n)}\cinfty(\R^{n+m})/\liem^{\infty}_{\R^n\times Y} \]
of connective $\cinfty(\R^{n+m})/\liem^{\infty}_{\R^n\times Y}$-modules, because base change is right t-exact. The long exact sequence associated to this last fibre sequence yields equivalences
\[ \tor^{\cinfty(\R^n)}_k(\liem_{\mathrm{Supp}(g)}^{\infty},\cinfty(\R^{n+m})/\liem^{\infty}_{\R^n\times Y}) \cong \tor^{\cinfty(\R^n)}_{k+1}((g),\cinfty(\R^{n+m})/\liem^{\infty}_{\R^n\times Y})   \]
for all $k\geq 1$. Let $k>1$ and suppose for the sake of induction that we have proven that for all $1\leq j<k$ and all $g'\in C^{\infty}(\R^n)$, the torsion group $\tor^{\cinfty(\R^n)}_{j}((g'),\cinfty(\R^{n+m})/\liem^{\infty}_{\R^n\times Y})$ vanishes, then in view of Lemma \ref{lem:cofinalprincipalideal}, the torsion groups $\tor^{\cinfty(\R^n)}_{j}(\liem^{\infty}_X,\cinfty(\R^{n+m})/\liem^{\infty}_{\R^n\times Y})$ also vanish for any closed $X\subset \R^n$ since the Tor functors commute with filtered colimits. Taking $X=\mathrm{Supp}(g)$, the isomorphism above shows that also for $j=k$, the torsion group $\tor^{\cinfty(\R^n)}_{j}((g),\cinfty(\R^{n+m})/\liem^{\infty}_{\R^n\times Y})$ vanishes for all $g\in C^{\infty}(\R^n)$. It remains to prove the base case $k=1$: we have an exact sequence
\[\begin{tikzcd} 0 \ar[r] & \tor^{\cinfty(\R^n)}_{1}((g),\cinfty(\R^{n+m})/\liem^{\infty}_{\R^n\times Y}) \ar[r]\ar[draw=none]{d}[name=X, anchor=center]{} & \tor_0^{C^{\infty}(\R^n)}(\liem_{\mathrm{Supp}(g)}^{\infty},C^{\infty}(\R^{n+m})/\liem^{\infty}_{\R^n\times Y})  \ar[rounded corners,
            to path={ -- ([xshift=2ex]\tikztostart.east)
                      |- (X.center) \tikztonodes
                      -| ([xshift=-2ex]\tikztotarget.west)
                      -- (\tikztotarget)}]{dll}\\ \cinfty(\R^{n+m})/\liem^{\infty}_{\R^n\times Y} \ar[r] & \tor_0^{C^{\infty}(\R^n)}((g),C^{\infty}(\R^{n+m})/\liem^{\infty}_{\R^n\times Y}) \ar[r]&  0 \end{tikzcd}  \]
so it suffices to show that the map $\tor_0^{C^{\infty}(\R^n)}(\liem_{\mathrm{Supp}(g)}^{\infty},C^{\infty}(\R^{n+m})/\liem^{\infty}_{\R^n\times Y})\rightarrow \cinfty(\R^{n+m})/\liem^{\infty}_{\R^n\times Y}$ is a monomorphism, but Lemma \ref{lem:cofinalprincipalideal} implies that this map is a filtered colimit of maps of the form 
\begin{align}\label{eq:principalflat}
\tor_0^{C^{\infty}(\R^n)}((h),C^{\infty}(\R^{n+m})/\liem^{\infty}_{\R^n\times Y})\longrightarrow \cinfty(\R^{n+m})/\liem^{\infty}_{\R^n\times Y} 
\end{align}
for $(h)$ a principal ideal of the commutative ring $\cinfty(\R^n)$. Since the collection of monomorphisms is stable under filtered colimits in a Grothendieck abelian category, we are reduced to showing that each of the maps \eqref{eq:principalflat} is a monomorphism. To show injectivity, we note that as $(h)$ has a generator, all objects in the tensor product over $\cinfty(\R^{n})$ are pure tensors of the form $h(\mathbf{x})\otimes F(\mathbf{x},\mathbf{y})$. Let $f(\mathbf{x},\mathbf{y})$ be a function on $\R^{n+m}$ representing the Whitney jet $F$, then we should show that if $h(\mathbf{x})f(\mathbf{x},\mathbf{y})\in \liem^{\infty}_{\R^n\times Y}$, then $h(\mathbf{x})\otimes f(\mathbf{x},\mathbf{y})$ is carried to 0 by the projection $\tor_0^{C^{\infty}(\R^n)}((h),C^{\infty}(\R^{n+m}))\rightarrow \tor_0^{C^{\infty}(\R^n)}((h),C^{\infty}(\R^{n+m})/\liem^{\infty}_{\R^n\times Y})$. We have the following two claims.
\begin{enumerate}[$(a)$]
    \item If $h(\mathbf{x})f(\mathbf{x},\mathbf{y})\in \liem^{\infty}_{\R^n\times Y}$, then the function $D_{\mathbf{y}}^{\alpha}f$ also has the property that $h(\mathbf{x})D_{\mathbf{y}}^{\alpha}f(\mathbf{x},\mathbf{y})\in \liem^{\infty}_{\R^n\times Y}$  for any multi-index $\alpha\in\Z_{\geq 0}^{m}$
    \item If $h(\mathbf{x})f(\mathbf{x},\mathbf{y})\in \liem^{\infty}_{\R^n\times Y}$, then the function $D^{\beta}_{\mathbf{x}}f$ vanishes on $\mathrm{Supp}(h)\times Y$ for any multi-index $\beta\in \Z_{\geq0}^{n}$,
\end{enumerate}
The first claim is obvious. For $(b)$ we proceed by induction on the order $|\beta|$. If $|\beta|=0$, the claim follows because $h(\mathbf{x})f(\mathbf{x},\mathbf{y})$ vanishes on $\R^n\times Y$. Now the inductive hypothesis and the product rule guarantee that on $\mathrm{Supp}(h)\times Y$, we have $D^{\beta}_{\mathbf{x}}(h(\mathbf{x})f(\mathbf{x},\mathbf{y}))=h(\mathbf{x})D^{\beta}_{\mathbf{x}}f(\mathbf{x},\mathbf{y})$ and by assumption we have $D^{\beta}_{\mathbf{x}}(h(\mathbf{x})f(\mathbf{x},\mathbf{y}))=0$ on $\R^n\times Y$, which establishes the claim. Applying for each multi-index $\alpha\in \Z_{\geq0}^m$ claim $(b)$ to the function $D_{\mathbf{y}}^{\alpha}f(\mathbf{x},\mathbf{y})$, we conclude that $f(\mathbf{x},\mathbf{y})\in \liem^{\infty}_{\mathrm{Supp}(h)\times Y}$ which implies by the flat function lemma that $f(\mathbf{x},\mathbf{y})$ can be written as $g(\mathbf{x},\mathbf{y})(\varphi_n(\mathbf{x})+\varphi_m(\mathbf{y}))$ with $\varphi_n(\mathbf{x})\in \liem^{\infty}_{\mathrm{Supp}(h)}$ and $\varphi_m(\mathbf{y})\in \liem^{\infty}_{Y}$. Then in $\tor_0^{C^{\infty}(\R^n)}((h),C^{\infty}(\R^{n+m}))$ we have the identities 
\[h(\mathbf{x})\otimes f(\mathbf{x},\mathbf{y})= h(\mathbf{x})\varphi_n(\mathbf{x})\otimes g(\mathbf{x},\mathbf{y})+ h(\mathbf{x})\otimes g(\mathbf{x},\mathbf{y})\varphi_m(\mathbf{y})= h(\mathbf{x})\otimes g(\mathbf{x},\mathbf{y})\varphi_m(\mathbf{y}),\] which the projection $\tor_0^{C^{\infty}(\R^n)}((h),C^{\infty}(\R^{n+m}))\rightarrow \tor_0^{C^{\infty}(\R^n)}((h),C^{\infty}(\R^{n+m})/\liem^{\infty}_{\R^n\times Y})$ clearly maps to 0.
\item Combine $(1)$, Lemma \ref{lem:cofinalprincipalideal} and the fact that acyclic modules are stable under filtered colimits, as are monomorphism in a Grothendieck abelian category. 
\end{enumerate}
\end{proof}
\begin{prop}\label{prop:cinftycoproductflat}
Let $X\subset\R^n$ be a closed subset. Let $I\subset\cinfty(\R^m)$ be an ideal that is either principal or of the form $\liem^{\infty}_Y$ for some closed subset $Y\subset\R^m$. Then the unit map of the $0$'th truncation functor
\[ \cinfty(\R^n)/\liem^{\infty}_X\otimes^{\infty} \cinfty(\R^m)/I \longrightarrow \tau_{\leq 0}(\cinfty(\R^n)/\liem^{\infty}_X\otimes^{\infty} \cinfty(\R^m)/I)\]
is an equivalence.
\end{prop}
\begin{proof}
Let $I$ be an ideal of $\cinfty(\R^n)$ of the form given in the statement of the proposition and denote by $\tilde{I}$ the ideal generated by $I$ under the map $\cinfty(\R^n)\rightarrow\cinfty(\R^{n+m})$. If $I=(h)$ is principal, then $\tilde{I}=(h)$ where $h$ is now viewed as a function on $\R^{n+m}$ constant in the $m$-variables and if $I=\liem^{\infty}_X$, then the flat function lemma implies immediately that $\tilde{I}=\liem^{\infty}_{X\times\R^m}$. Consider the diagram 
\begin{equation}\label{eq:bigsquare}
\begin{tikzcd}
\R \ar[d]\ar[r] & \cinfty(\R^n) \ar[d]\ar[r] &  \cinfty(\R^n)/I \ar[d] \\
\cinfty(\R^m) \ar[r] \ar[d] &\cinfty(\R^{n+m}) \ar[r]\ar[d] &  \cinfty(\R^{n+m})/\tilde{I}  \ar[d] \\
\cinfty(\R^{m})/\liem^{\infty}_Y \ar[r] & \cinfty(\R^{n+m})/\liem^{\infty}_{\R^n\times Y} \ar[r]& \cinfty(\R^{n+m})/\liem^{\infty}_{\R^n\times Y}\oinfty_{\cinfty(\R^{n+m})}\cinfty(\R^{n+m})/\tilde{I}
\end{tikzcd}
\end{equation}
in which the upper left square and the lower right square are pushouts in $\sring$. We claim that the lower left square and upper right square are also pushouts of simplicial $\cinfty$-rings. It suffices to check the upper right one. Since the map $\cinfty(\R^n)\rightarrow\cinfty(\R^n)/I$ is an effective epimorphism, Corollary \ref{cor:algeffepi} asserts that it suffices to show that the associated diagram of simplicial $\R$-algebras is a pushout. Clearly, the diagram becomes a pushout after taking the 0'th truncation, so it suffices to show that the higher homotopy groups of the pushout $\cinfty(\R^n)/I\otimes_{\cinfty(\R^{n})}\cinfty(\R^{n+m})$ vanish. We have a fibre sequence 
 \[ I\longrightarrow C^{\infty}(\R^n) \longrightarrow \cinfty(\R^n)/I\]
 of discrete $\cinfty(\R^n)$-modules, so we get a fibre sequence 
 \[ I\otimes_{\cinfty(\R^n)}\cinfty(\R^{n+m})\longrightarrow \cinfty(\R^{n+m})\longrightarrow \cinfty(\R^n)/I\otimes_{\cinfty(\R^{n})}\cinfty(\R^{n+m})\]
 of connective $\cinfty(\R^{n+m})$-modules. By Lemma \ref{lem:acyclicideals} (for $Y=\R^m$), $\tor_n^{\cinfty(\R^n)}(I,\cinfty(\R^{n+m}))$ vanishes for all $n\geq 1$ and the map $\tor_0^{\cinfty(\R^n)}(I,\cinfty(\R^{n+m}))\rightarrow \cinfty(\R^{n+m})$ is a monomorphism, so the long exact sequence associated to the fibre sequence above guarantees the vanishing of $\tor_n^{\cinfty(\R^n)}(\cinfty(\R^n)/I,\cinfty(\R^{n+m}))$ for all $n\geq 1$, which establishes the claim. It follows that all squares in the diagram \eqref{eq:bigsquare} are pushouts, including the four rectangles and the big square. Therefore, we can identify the coproduct $\cinfty(\R^n)/I\oinfty\cinfty(\R^m)/\liem^{\infty}_Y$ with the pushout $\cinfty(\R^{n+m})/\liem^{\infty}_{\R^n\times Y}\oinfty_{\cinfty(\R^n)}\cinfty(\R^n)/I$, so in order to prove the proposition, it suffices to argue that this pushout is 0-truncated. Using Corollary \ref{cor:algeffepi}, we are reduced to showing that the higher homotopy groups of the pushout $\cinfty(\R^{n+m})/\liem^{\infty}_{\R^n\times Y}\otimes_{\cinfty(\R^n)}\cinfty(\R^n)/I$ of simplicial commutative $\R$-algebras vanish. We have a fibre sequence 
 \[ I\otimes_{\cinfty(\R^n)}\cinfty(\R^{n+m})/\liem^{\infty}_{\R^n\times Y}\longrightarrow \cinfty(\R^{n+m})/\liem^{\infty}_{\R^n\times Y}\longrightarrow \cinfty(\R^n)/I\otimes_{\cinfty(\R^{n})}\cinfty(\R^{n+m})/\liem^{\infty}_{\R^n\times Y}\]
 of connective $\cinfty(\R^{n+m})/\liem^{\infty}_{\R^n\times Y}$-modules. By Lemma \ref{lem:acyclicideals}, $\tor_n^{\cinfty(\R^n)}(I,\cinfty(\R^{n+m})/\liem^{\infty}_{\R^n\times Y})$ vanishes for all $n\geq 1$ and the map $\tor_0^{\cinfty(\R^n)}(I,\cinfty(\R^{n+m})/\liem^{\infty}_{\R^n\times Y})\rightarrow \cinfty(\R^{n+m})/\liem^{\infty}_{\R^n\times Y}$ is a monomorphism, so the long exact sequence associated to the fibre sequence above guarantees the vanishing of $\tor_n^{\cinfty(\R^n)}(\cinfty(\R^n)/I,\cinfty(\R^{n+m})/\liem^{\infty}_{\R^n\times Y})$ for all $n\geq 1$, which concludes the proof.
\end{proof}
\begin{proof}[Proof of Theorem \ref{thm:reyesvanquederived}]
By Proposition \ref{prop:cinftycoproductflat}, it suffices to show that the maps $\cinfty(X;\R^n)\rightarrow\cinfty(X\times Y;\R^{n+ m})$ and $\cinfty(Y;\R^m)\rightarrow\cinfty(X\times Y;\R^{n+ m})$ determine an isomorphism $\tau_{\leq 0}(C^{\infty}(X;\R^n)\otimes^{\infty}C^{\infty}(Y;\R^m))\cong \cinfty(X\times Y;\R^{n+ m})$ of ordinary $\cinfty$-rings. Unwinding the definitions, this amounts to the assertion that the ideal $\liem^{\infty}_{X\times Y}\subset\cinfty(\R^{n+m})$ is generated by the images of the ideals $\liem^{\infty}_X$ and $\liem^{\infty}_Y$, but this follows immediately from Tougeron's flat function lemma.
\end{proof}
\newpage
\section{Derived $\cinfty$-geometry with corners}\label{sec:cornerlog}
In this section we use the results of the previous sections to define \emph{derived $\cinfty$-schemes with corners} and study their basic properties. Our goal will be to introduce a tractable geometry $\geodiffderc$ controlling derived $\cinfty$-geometry with corners. To construct such a geometry, it seems natural to consider the \infcat $\sring_{pc}$ of algebras for the Lawvere theory of Cartesian spaces with corners. In the 1-categorical setting, such a theory has been developed by Joyce and Francis-Staite \cite{Joyfra}.
\begin{rmk}
For technical reasons, the geometry we will construct does not come equipped with a functor $\diffc\rightarrow\geodiffderc$ exhibiting a geometric envelope. We will instead introduce another pregeometry $\diffc'$ together with functors
\[\diffc \longleftarrow \diffc' \hooklongrightarrow \geodiffderc,  \]
the left one of which is a Morita equivalence of pregeometries and the fully faithful arrow on the right exhibits a geometric envelope.
\end{rmk}
\subsection{Simplicial $\cinfty$-rings with corners}
At first glance, one might be tempted to define $\geodiffderc$ as the compact objects of $\sring_{pc}$, exactly analogous to how $\geodiffder$ was introduced, but this turns out to be too naive: the functor \[(\cinfty(\_),\cinfty_b(\_)):\diffc\longrightarrow\sring_{pc}^{op}\]
does not preserve pullbacks along open inclusions. Consider, for instance, the following pullback
\[ 
\begin{tikzcd}
\R\ar[d,hook,"\exp"] \ar[r,equal] & \R\ar[d,hook,"\exp"] \\
\R_{\geq 0} \ar[r] & \R
\end{tikzcd}
\]
in $\diffc$. The vertical maps are admissible, so the diagram 
\[
\begin{tikzcd}
(\cinfty(\R),\cinfty_b(\R)) \ar[d,"\exp^*"] \ar[r] & (\cinfty(\R_{\geq 0}),\cinfty_b(\R_{\geq 0})) \ar[d,"\exp^*"] \\
(\cinfty(\R),\cinfty_b(\R)) \ar[r,equal] & (\cinfty(\R),\cinfty_b(\R))
\end{tikzcd}
\]
should be a pushout diagram in $\sring_{pc}$. Let $C$ denote the pushout of the diagram above, and let $\varphi:C\rightarrow (\cinfty(\R),\cinfty_b(\R))$ be the canonical morphism; we should verify whether or not this morphism is an equivalence. Corollary \ref{cor:cornerforgetcolimitpreserve} shows that the map $\cinfty(\R)\rightarrow \ev_{\R}(C)$ is an equivalence of spaces. Since evaluation at $\R_{\geq 0}$ preserves sifted colimits, the space $\ev_{\R_{\geq0}}(C)$ may be computed as the colimit of the simplicial object $\mathsf{Bar}_{\cinfty_b(\R)}(\cinfty_b(\R),\cinfty_b(\R_{\geq 0}))_{\bullet}$ which takes the form 
\[\begin{tikzcd} \ldots\ar[r,shift left=3] \ar[r,shift left] \ar[r,shift right]\ar[r,shift right=3] & \cinfty_b(\R\times \R^2 \times\R_{\geq 0}) \ar[r,shift left=2] \ar[r]\ar[r,shift right=2] & \cinfty_b(\R\times \R^1 \times\R_{\geq 0}) \ar[r,shift left]\ar[r,shift right] & \cinfty_b(\R \times\R_{\geq 0}). \end{tikzcd}  \]
In each simplicial level, we have the space of interior $b$-maps on a manifold of the form $\R\times\R^n\times \R_{\geq 0}$, which has one connected boundary component whose defining function is the last coordinate. All face maps preserve this boundary defining function, so it follows from Lemma \ref{lem:intblogification} that this simplicial object may be written as a product
\[ \Z_{\geq0}\times\left(\begin{tikzcd}\ldots\ar[r,shift left=3] \ar[r,shift left] \ar[r,shift right]\ar[r,shift right=3] & \cinfty_{>0}(\R\times \R^2 \times\R_{\geq 0}) \ar[r,shift left=2] \ar[r]\ar[r,shift right=2] & \cinfty_{>0}(\R\times \R^1 \times\R_{\geq 0}) \ar[r,shift left]\ar[r,shift right] & \cinfty_{>0}(\R \times\R_{\geq 0}) \end{tikzcd}\right).\]
where $\Z_{\geq0}$ is a constant simplicial object. The simplicial object in parentheses is equivalent to \[\begin{tikzcd} \ldots\ar[r,shift left=3] \ar[r,shift left] \ar[r,shift right]\ar[r,shift right=3] & \cinfty(\R\times \R^2 \times\R_{\geq 0}) \ar[r,shift left=2] \ar[r]\ar[r,shift right=2] & \cinfty(\R\times \R^1 \times\R_{\geq 0}) \ar[r,shift left]\ar[r,shift right] & \cinfty(\R \times\R_{\geq 0}), \end{tikzcd}  \] 
whose colimit can be identified with the pushout $\cinfty(\R)\oinfty_{\cinfty(\R)}\cinfty(\R_{\geq 0})\simeq\cinfty(\R)$. Since sifted colimits commute with products, we find that $\ev_{\R_{\geq0}}(C)\simeq \Z_{\geq0}\times \cinfty_{>0}(\R)$, and the map $\ev_{\R_{\geq 0}}(\varphi)$ is identified with the projection $\Z_{\geq0}\times \cinfty_{>0}(\R)\rightarrow \cinfty_{>0}(\R)= \cinfty_b(\R)$, which is not an equivalence. We could correct for the fact that the pullback diagram above is not preserved simply by localizing at the morphism $\varphi:(\cinfty(\R),\cinfty_b(\R))\rightarrow C$: we may consider the strongly reflective subcategory $\sring_{pc}[\{\varphi\}]$ so that the functor
\[ \diffc\hooklongrightarrow \sring_{pc}[\{\varphi\}]^{op} \]
does preserve our pullback. For the purposes of defining a transformation of pregeometries $\diffc\rightarrow\geodiffderc$ however, this appears too myopic; it seems we have to localize at \emph{all} comparison maps arising from applying $(\cinfty(\_),\cinfty_b(\_))$ to admissible pullbacks $\diffc$. Fortunately, it turns out that, as we will show in this section, the diagram
\[ 
\begin{tikzcd}
\R\ar[d,hook,"\exp"] \ar[r,equal] & \R\ar[d,hook,"\exp"] \\
\R_{\geq 0} \ar[r] & \R
\end{tikzcd}
\]
was (with the benefit of hindsight) fortuitously chosen: localizing \emph{only} at $\varphi$ already yields the correct ambient \infcatt. For the sequel, it will be convenient to introduce a different morphism which will yield the same localization.
\begin{defn}\label{defn:cinftycorners}
Consider the image of the map 
\[ \R_{>0}\hooklongrightarrow \R_{\geq 0} \]
under the functor $(\cinfty(\_),\cinfty_{b}(\_)):\diffc^{op}\hookrightarrow\sring_{pc}$. Denote by $\epsilon$ the counit $\epsilon:\iota_{c!}\iota_c^*\rightarrow\mathrm{id}$ and define an object $\mathcal{A}\in \sring_{pc}$ together with a map $\phi:\iota_{c!}\iota_c^* (\cinfty(\R_{> 0}),\cinfty_{b}(\R_{> 0})) \rightarrow\mathcal{A}$ via the pushout diagram
\[
\begin{tikzcd}
\iota_{c!}\iota_c^* (\cinfty(\R_{\geq 0}),\cinfty_{b}(\R_{\geq 0})) \ar[d] \ar[r,"\epsilon"] & (\cinfty(\R_{\geq 0}),\cinfty_{b}(\R_{\geq 0}))  \ar[d]\\ \iota_{c!}\iota_c^* (\cinfty(\R_{> 0}),\cinfty_{b}(\R_{> 0})) 
\ar[r,"\phi"]& \mathcal{A}.  
\end{tikzcd}
\]
We let $S=\{\phi\}$, the one element set containing the morphism $\phi$. The \infcat of \emph{simplicial $\cinfty$-rings with corners}, denoted $\sring_c$, is the presentable \infcat of $S$-local objects of $\sring_{pc}$. 
\end{defn}
\begin{rmk}\label{rmk:cinftycorners}
Unraveling the definition, a simplicial $\cinfty$-ring with pre-corners $(A,A_c)$ is $S$-local just in case the upper horizontal map in the pullback diagram 
\[
\begin{tikzcd}
A_c\times_{A_{\geq 0}}A_{>0} \ar[r]\ar[d] & A_{>0} \ar[d] \\
A_c\ar[r] & A_{\geq 0}
\end{tikzcd}
\]
of spaces is an equivalence, where we use the notation \[A_{\geq 0}:=\Hom_{\sring}(\cinfty(\R_{\geq0}),A),\quad \text{and}\quad A_{>0}:=\Hom_{\sring}(\cinfty(\R_{>0}),A)\]
for $A$ a simplicial $\cinfty$-ring. This subsection will be concerned with the simplicial commutative monoid structure on the space $A_{\geq0}$ induced by the homotopy coherent $\cinfty$-operations. To this end, it turns out to be crucial to establish, as we will in a moment, that the right vertical map $A_{>0}\rightarrow A_{\geq 0}$ -an inclusion of connected components- coincides with the largest subgroup contained in the simplicial commutative monoid $A_{\geq 0}$. Contrary to the 1-categorical case, this is not immediate and depends on a computation of the cotangent complex of $\cinfty(\R_{\geq0})$ which is deferred to part II (see Remark \ref{rmk:truncatepositivelts}).
\end{rmk}
\begin{rmk}\label{rmk:altset}
On $(\cinfty(\R),\cinfty_b(\R))$ and $(\cinfty(\R_{>0}),\cinfty_b(\R_{>0}))$, the counit $\iota_{c!}\iota_c^*\rightarrow\mathrm{id}$ is an equivalence, so for each $(A,A_c)\in \sring_{pc}$, there is a diagram 
\[
\begin{tikzcd}
A_c\times_{A_{\geq 0}}A_{>0} \ar[r]\ar[d] & A_{>0}\ar[r,equal] \ar[d] &   A_{>0} \ar[d]\\
A_c\ar[r] & A_{\geq 0} \ar[r] & A
\end{tikzcd}
\]
Here, the right square is a pullback. Let $S'=\{\varphi\}$, the one element set containing the map $\varphi:(\cinfty(\R),\cinfty_b(\R))\rightarrow C$ from the discussion above. Unwinding the definitions, we see that $(A,A_c)$ is $S$-local if and only if it is $S'$-local.
\end{rmk}
Since the forgetful functor $\iota_c^*:\sring_{pc}\rightarrow\sring$ preserves colimits and carries the counit $\epsilon:\iota_{c!}\iota_c^*\rightarrow\mathrm{id}$ to the identity, $\iota_c^*$ carries the map $\phi$ of Definition \ref{defn:cinftycorners} to an equivalence. From the universal property of cocontinuous localizations, we deduce that $\iota_c^*$ factors via a left adjoint $\sring_c\rightarrow\sring$. This functor coincides with the composition $\sring_c\hookrightarrow\sring_{pc}\rightarrow\sring$, which is a right adjoint. Note that both adjoints of this functor are fully faithful, so the argument of Proposition \ref{prop:cornerpresfib} grants the following result.
\begin{prop}
The functor $\sring_c\rightarrow \sring$ is a presentable fibration. Moreover, the inclusion $\sring_c\hookrightarrow \sring_{pc}$ preserves Cartesian edges, and the localization $L:\sring_{pc}\rightarrow\sring_{c}$ preserves coCartesian edges.
\end{prop}
\begin{rmk}
As Corollary \ref{cor:accessiblecorner} asserts, the localization $\sring_c\subset\sring_{pc}$ is $\omega$-accessible, that is, $\sring_c\subset\sring_{pc}$ is stable under filtered colimits. In particular, every compact object in $\sring_c$ is a retract of an object in the image of $L$. Since $\iota_{c!}L$ is equivalent to $L$ and idempotents may be lifted along coCartesian fibrations, we deduce that for any compact object $(A,A_c)$ in $\sring_c$, there is some $A'_c$ such that $(A,A_c')$ is compact in $\sring_{pc}$. 
\end{rmk}
We will use the \infcat $\sring_c$ to define a derived geometry generated by manifolds with corners. To this end, we first make an observation concerning finitely generated and compact objects in $\sring_{pc}$.
\begin{prop}\label{prop:compactcornersfair}
The following hold true.
\begin{enumerate}[$(1)$]
    \item The functor $\iota_c^*$ carries finitely generated objects of $\sring_{pc}$ to finitely generated objects of $\sring$.
    \item The functor $\iota_c^*$ carries finitely presented objects of $\sring_{pc}$ into the full subcategory $\sring_{\gmt}\subset\sring$.
\end{enumerate}
\end{prop}
\begin{rmk}
In the proof below, we use that a finitely presented simplicial $\cinfty$-ring with pre-corners has a perfect cotangent complex, which follows formally from the fact that simplicial $\cinfty$-rings of the form $\cinfty(\R^n\times\R^k_{\geq0})$ have a perfect cotangent complex. We also use that if a map $f:A\rightarrow B$ of simplicial $\cinfty$-rings has a perfect \emph{relative cotangent complex} $\cotan_f$ and $\pi_0(f)$ exhibits $\pi_0(B)$ as finitely presented over $\pi_0(B)$, then $f$ exhibits $B$ as finitely presented over $B$. 
\end{rmk}

\begin{proof}
For $(1)$, we need to show that the right adjoint $\iota_{c*}$ preserves colimits of filtered diagrams consisting only of monomorphisms. From the general theory of algebraic theories it is enough to check that the functors $\ev_{\R^n\times\R^k_{\geq 0}}\iota_{c*}:\sring\rightarrow\spa$ have this property, but these functors are corepresented by the finitely generated objects $\cinfty(\R^n\times\R_{\geq 0}^{k})$. For $(2)$, we suppose that $(A,A_c)$ is finitely presented in $\sring_{pc}$. Consider a finite presentation of $(\pi_0(A),\pi_0(A_c))$, that is, a coequalizer diagram 
\[ \begin{tikzcd}(\cinfty(\R^p\times\R^q_{\geq0}),\cinfty_b(\R^p\times\R^q_{\geq0})) \ar[r,shift left]\ar[r,shift right]& (\cinfty(\R^n\times\R^m_{\geq0}),\cinfty_b(\R^n\times\R^m_{\geq0})) \ar[r] & (\pi_0(A),\pi_0(A_c)). \end{tikzcd}  \]
As $\iota_c^*:\cinfty\mathsf{ring}_c\rightarrow \cinfty\mathsf{ring}$ preserves colimits, we have a coequalizer diagram 
\[ \begin{tikzcd} \cinfty(\R^p\times\R^q_{\geq0}) \ar[r,shift left]\ar[r,shift right]& \cinfty(\R^n\times\R^m_{\geq0})\ar[r] & \pi_0(A). \end{tikzcd}  \]
Since we have an epimorphism of $\cinfty$-rings $\cinfty(\R^{p+q})\rightarrow \cinfty(\R^p\times\R^q_{\geq0}) $, we also have a coequalizer diagram 
\[ \begin{tikzcd} \cinfty(\R^{p+q}) \ar[r,shift left]\ar[r,shift right]& \cinfty(\R^n\times\R^m_{\geq0})\ar[r] & \pi_0(A) \end{tikzcd}  \]
which shows that $\pi_0(A)$ is finitely presented over $\cinfty(\R^n\times\R^m_{\geq0})$ \emph{as a $\cinfty$-ring}. Since the latter object is free in $\sring_{pc}$, the map $\cinfty(\R^n\times\R^m_{\geq0})\rightarrow \pi_0(A)$ lifts to a map $f:\cinfty(\R^n\times\R^m_{\geq0})\rightarrow A$ in $\sring$, and as we have just verified, $\pi_0(f)$ is finitely presented. The cotangent complexes of both $\cinfty(\R^n\times\R^m_{\geq0})$ and $A$ are perfect, so the relative cotangent complex $\cotan_f$ is also perfect, which implies that $A$ is finitely presented over $\cinfty(\R^n\times\R^m_{\geq0})$. Invoking Proposition \ref{afpisfair}, we deduce that $A$ is geometric.
\end{proof}
Recall that geometricity implies that if $(A,A_c)$ is a compact object of $\sring_{pc}$, then the underlying simplicial $\cinfty$-ring $A$ has the property that $\pi_0(A)$ is finitely generated and germ determined, and for each $n\geq 1$, the object $\pi_n(A)$ has the property that module of global sections of the sheafification of the presheaf
\[  U_a\longmapsto  \pi_n(A)\otimes_{\pi_0(A)}\pi_0(A)[a^{-1}]  \]
coincides with $\pi_n(A)$. More briefly, for $(A,A_c)$ compact, the unit map
\[ A\longrightarrow \Gamma\spec\,A \]
is an equivalence.
\begin{defn}\label{defn:cornergeometry}
Let $\geodiffderc$ be the opposite of the full subcategory of $\sring_c^{op}$ spanned by compact objects. We define the notions of an admissible morphism and admissible covering in $\geodiffderc$ as follows.
\begin{enumerate}[$(1)$]
    \item A morphism $f:\mathrm{Spec}\,(A,A_c)\rightarrow \mathrm{Spec}\, (B,B_c)$ is admissible if and only if there exists some $b\in \pi_0(B)$ such that the underlying map $B\rightarrow A$ of simplicial $\cinfty$-rings exhibits $A$ as a localization of $B$ by $b$ and $f$ is a coCartesian morphism for the fibration $\sring_c
    \rightarrow \sring$.
    \item A collection of morphisms $\{\mathrm{Spec}\,(B_i,B_{ic})\rightarrow \mathrm{Spec}\, (B,B_c)\}$ generates a covering sieve if and only if the underlying collection $\{\mathrm{Spec}\,B_i\rightarrow \mathrm{Spec}\, B\}$ of morphisms among geometric (cf. Proposition \ref{prop:compactcornersfair} and the preceding remark) simplicial $\cinfty$-rings generates a covering sieve for the \'{e}tale topology on $\sring_{\gmt}^{op}\simeq\daff_{\gmt}$.
\end{enumerate}
Let $\diffc'\subset \geodiffderc$ be the smallest full subcategory of $\geodiffderc$ that contains the objects $(\cinfty(\R^n),\cinfty_b(\R^n))$ for all $n$ and satisfies the following condition: should $f:\mathrm{Spec}\,(A,A_c)\rightarrow \mathrm{Spec}\, (B,B_c)$ be admissible and $\mathrm{Spec}\, (B,B_c)\in \diffc'$, then also $\mathrm{Spec}\,(A,A_c)\in \diffc'$ (hence we require that the inclusion $\diffc'\subset\geodiffderc$ is a categorical fibration).
\end{defn}
\begin{rmk}\label{rmk:cocartesiancorners}
Suppose that $(A,A_c)\rightarrow (B,B_c)$ is (the opposite of) an admissible morphism in $\sring_c$, then it might not be a priori clear that $(B,B_c)$ is a compact object of $\sring_c$. To see this is the case, we note that the assumption that $A\rightarrow B$ is a localization of simplicial $\cinfty$-rings provides a map $\cinfty(\R)\rightarrow A$ and an equivalence $B\simeq \cinfty(\R\setminus\{0\})\oinfty_{\cinfty(\R)}A$. We have an adjoint map $L\iota_{c!}(\cinfty(\R))\rightarrow (A,A_c)$ and we can form the pushout
\[
\begin{tikzcd}
L\iota_{c!}(\cinfty(\R)) \ar[r] \ar[d] & L\iota_{c!}(\cinfty(\R\setminus\{0\})) \ar[d] \\
(A,A_c)  \ar[r]& (C,C_c).
\end{tikzcd}
\]
Since $\sring_c\rightarrow\sring$ preserves colimits, the map $A\rightarrow C$ coincides with $A\rightarrow B$. Because the upper horizontal map is coCartesian by the description of coCartesian edges in Proposition \ref{prop:cornerpresfib}, the lower horizontal map is coCartesian as well, as all colimits in $\sring_c$ are relative colimits. It follows that $(C,C_c)\simeq (B,B_c)$, that is, being admissible in $\sring_c$ is equivalent to fitting into a pushout diagram as above, which shows that $(B,B_c)$ is compact if $(A,A_c)$ is, since $L\iota_{c!}$ preserves compact objects.
\end{rmk}
The main results of this subsection are summarized in the following theorem.
\begin{thm}\label{thm:geoenvcorners}
The following hold true.
\begin{enumerate}[$(i)$]
    \item The \infcat $\sring_c$ is compactly generated, that is, the canonical functor $\mathrm{Pro}(\geodiffderc)\rightarrow \sring_c^{op}$ is an equivalence.
    \item Definition \ref{defn:cornergeometry} furnishes the structure of a geometry on $\geodiffderc$ and the structure of a pregeometry on $\diffc'$ such that the inclusion $\diffc'\subset\geodiffderc$ is a transformation of pregeometries.
    \item The functor $\diffc'\subset \geodiffderc$ exhibits a geometric envelope.
    \item The functor $(\cinfty_{\_},\cinfty_{b\_}):\diffc\rightarrow\rtop(\geodiffderc)$ is fully faithful and preserves pullbacks along admissible maps.
    \item Denote by $\spec_c$ the functor $\spec^{\geodiffderc}$, then $\spec_c:\diffc'\rightarrow \rtop(\geodiffderc)$ takes values in the essential image of $(\cinfty_{\_},\cinfty_{b\_})$ and determines a Morita equivalence of pregeometries
    \[ \diffc'\longrightarrow \diffc. \]
\end{enumerate}

\end{thm}
\subsubsection{Positive log structures}
The proof of theorem \ref{thm:geoenvcorners} will require a number of prelimenaries. First observe that the geometry structure on $\geodiffderc$ makes reference to the coCartesian morphisms of $\iota_c^*$, which involve the formation of certain pushouts in $\sring_c$ and are more difficult to characterize explicitly than its Cartesian morphism, which are obtained by taking certain pullbacks in $\sring_c$ and are therefore detectable on the underlying spaces. To improve our understanding of the fibres of $\sring_c\rightarrow\sring$ and its coCartesian morphisms, we will establish a structural result of independent interest which relates simplicial $\cinfty$-rings with corners to an algebraic model for $\cinfty$-geometry with corners and more general singularities. The latter theory is a derived and differential geometric version of \emph{logarithmic geometry} in the sense of Fontaine-Illusie, Kato and Ogus \cite{Kato,Ogus}. 
\begin{rmk}
While we make no use of this perspective, the theory of positive logarithmic $\cinfty$-geometry we expose in this subsection could have been developed entirely in a model categorical setting, as is done by Sagave, Sch\"{u}rg and Vezzosi and Bhatt \cite{SSV,Bhatt}, at the cost of rendering many arguments significantly more cumbersome. In particular, it is not hard to see that the equivalence of Theorem \ref{thm:logcinftycorners} is induced by a Quillen equivalence between combinatorial model categories. We leave it as an exercise for the sufficiently industrious reader to make the necessary translations. 
\end{rmk}
\begin{rmk}
Apart from the works of Sagave-Sch\"{u}rg-Vezzosi and Bhatt, the derived antecedents of this section include the work on logarithmic structures for $\einfty$-ring spectra and applications to THH of Rognes, Sagave and Schlichtkrull \cite{Rognes,RSS}. In differential geometry, the origins of logarithmic ideas trace back to the \emph{$b$-geometry} of Melrose \cite{Melrose1}, made explicit in the work of Kottke-Melrose \cite{KottkeMelrose}, and especially that of Gillam-Molcho \cite{GilMol}. 
\end{rmk}
\begin{nota}
As in our notation, the set $\N$ does not contain 0, we write $\Z_{\geq 0}$ for the free commutative monoid on one generator. The commutative product in a generic commutative monoid is written additively $(\_+\_)$, while the product in a commutative monoid coming from a commutative algebra is written multiplicatively (by juxtaposing the elements being multiplied).
\end{nota}
\begin{cons}\label{cons:logcinftyrings}
Let $\mathsf{CartSp}_{\geq 0}\subset\cartsp_c$ be the full subcategory spanned by the objects of the form $\R^n_{\geq 0}$; this determines a ($1$-sorted) Lawvere theory. Recall the notation $\mathsf{FCMon}$ for the category of finitely generated and free commutative monoids. We define a functor $\theta:\mathsf{FCMon}^{op}\rightarrow \mathsf{CartSp}_{\geq 0}$ as follows.
\begin{enumerate}[$(1)$]
    \item $\theta$ carries the free commutative monoid $\Z_{\geq 0}^n$ to $\R^n_{\geq 0}$.
    \item $\theta$ carries a morphism $f:\Z_{\geq 0}^{n}\leftarrow \Z_{\geq0}^m$ determined by an $m$-tuple $\{(k^i_1,\ldots,k^i_n)\in \Z^n_{\geq0}\}_{1\leq i\leq m}$ to the smooth map $\R^n_{\geq 0}\rightarrow\R^m_{\geq 0}$ given by 
    \[ (x_1,\ldots,x_n) \longmapsto \left(\prod_{1\leq j\leq n}x^{k^1_j}_{j},\ldots,\prod_{1\leq j\leq n}x^{k^m_j}_{j}\right).\]
\end{enumerate}
Restricting along the product preserving functor $\theta$ induces a functor 
\[ \theta^*: \fun^{\pi}(\cartsp_{\geq 0},\spa)\longrightarrow s\mathsf{CMon}\]
in $(\prl_{\mathrm{Proj}})^{op}$ preserving limits and sifted colimits. Composing $\theta^*$ with the functor induced by the product preserving full subcategory inclusion $\iota_{\geq 0}:\cartsp_{\geq 0}\hookrightarrow\cartsp_c$ yields a limit and sifted colimit preserving functor $\theta^*\iota_{\geq0}^*:\sring_{pc}\rightarrow s\mathsf{CMon}$. Corollary \ref{cor:cornerforgetcolimitpreserve} provides a right adjoint $\iota_{c*}$ to the functor $\iota_c^*:\sring_{pc}\rightarrow \sring$ induced by the inclusion $\iota_c:\cartsp\hookrightarrow\cartsp_c$. The composite functor $\theta^*\iota_{\geq0}^*\iota_{c*}$ carries simplicial $\cinfty$-rings to simplicial commutative monoids, and we will denote this functor by $(\_)_{\geq 0}:\sring\rightarrow s\mathsf{CMon}$. We define the presentable \infcat of \emph{positive prelog simplicial $\cinfty$-rings} as the cone in the pullback diagram
\[
\begin{tikzcd}
\splring \ar[d] \ar[r] & \fun(\Delta^1,s\mathsf{CMon})\ar[d,"\ev_{1}"] \\
\sring \ar[r,"(\_)_{\geq0}"] & s\mathsf{CMon}
\end{tikzcd}
\]
among presentable \infcats and functors that admit left adjoints between them. An object of $\splring$ consists of a pair $(A,M\rightarrow A_{\geq 0})$ where $A$ is a simplicial $\cinfty$-ring and $M\rightarrow A_{\geq 0}$ is a map of simplicial commutative monoids. \\
We define a functor $\sring_{pc}\rightarrow \splring$ as follows. Composing the unit transformation $\mathrm{id}\rightarrow \iota_{c*}\iota_c^*$ with $\theta^*\iota_{\geq 0}^*$ yields a functor
\[ \sring_{pc}\longrightarrow \fun(\Delta^1,\sring_{pc})\longrightarrow\fun(\Delta^1,s\mathsf{CMon}) \]
which participates as the top horizontal map in the strictly commuting diagram 
\[
\begin{tikzcd}
\sring_{pc} \ar[d,"\iota_c^*"] \ar[r] & \fun(\Delta^1,s\mathsf{CMon})\ar[d,"\ev_{1}"] \\
\sring \ar[r,"(\_)_{\geq0}"] & s\mathsf{CMon}
\end{tikzcd}
\]
among \infcatst; hence we obtain an induced functor 
\[
\begin{tikzcd}
\sring_{pc} \ar[dr,"\iota_c^*"'] \ar[rr,"\Xi"] && \splring \ar[dl,"p"] \\
& \sring
\end{tikzcd}
\]
which is given on objects by the assignment $(A,A_c)\mapsto (A,A_c\rightarrow A_{\geq 0})$. From the description of $\iota_c^*$-Cartesian edges in Proposition \ref{prop:cornerpresfib} and the fact that $\sring_{pc}\rightarrow s\mathsf{CMon}$ preserves limits we immediately deduce that $\Xi$ carries $\iota_c^*$-Cartesian edges to $p$-Cartesian edges.
\end{cons}
\begin{rmk}
The functor $\Xi$ of Construction \ref{cons:logcinftyrings} does not take $\iota^*_c$-coCartesian edges to $p$-coCartesian edges and therefore merely induces a lax natural transformation between straightened functors $\mathrm{St}^{+,\mathrm{co}}(\iota^*_c)\Rightarrow \mathrm{St}^{+,\mathrm{co}}(p)$. By the results of \cite{HHLN}, straightening/unstraightening yields equivalences
\[ \fun(\icat,\catinf)_{\mathrm{lax}} \simeq \mathsf{biCart}^{(\mathrm{co})lax}_{\icat} \simeq \fun(\icat^{op},\catinf)_{\mathrm{colax}}. \]
For each map $A\rightarrow B$ of simplicial $\cinfty$-rings, requisite the 2-morphism is given by the Beck-Chevalley transformation. \end{rmk}
To aid our analysis, we recall some facts about simplicial abelian groups and simplicial commutative monoids.
\begin{lem}
Consider $s\mathsf{CMon}$ with its coCartesian symmetric monoidal structure and $\spa$ with its Cartesian symmetric monoidal structure, then the forgetful functor $s\mathsf{CMon}\rightarrow \spa$ induced by evaluation at $\Z_{\geq 0}$ has a canonical symmetric monoidal structure.
\end{lem}
\begin{proof}
We sketch two proofs. According to \cite{HA}, Theorem 2.3.4.18, the functor $f:s\mathsf{CMon}\rightarrow \mathsf{Mon}_{\einfty}$ classifies an \infop map $s\mathsf{CMon}^{\coprod}\rightarrow \spa^{\times}$ lifting the functor evaluating at $\Z_{\geq 0}$. Unwinding the definitions, this functor is symmetric monoidal if and only if $f$ preserves finite coproducts, which is the case. \\
For another argument, it is not hard to see that the functor of $1$-categories $\mathsf{FCMon}\rightarrow \mathsf{Set}$ has a canonical symmetric monoidal structure, and the relevant symmetric monoidal functor can be obtained by symmetric monoidal left Kan extension.
\end{proof}
\begin{prop}\label{prop:scomsab}
A simplicial commutative monoid $A$ is \emph{grouplike} if the commutative monoid $\pi_0(A)$ is a (necessarily abelian) group. Let $s\mathsf{CMon}^{\mathrm{gp}}$ denote the full subcategory spanned by the grouplike commutative monoids.
\begin{enumerate}[$(1)$]
    \item The full subcategory inclusion $s\mathsf{CMon}^{\mathrm{gp}}\subset s\mathsf{CMon}$ admits a right adjoint (that we will denote $(\_)^{\times}$, the \emph{$\infty$-group of units}). 
    \item The full subcategory inclusion $s\mathsf{CMon}^{\mathrm{gp}}\subset s\mathsf{CMon}$ admits a left adjoint (that we will denote $(\_)^{\mathrm{gp}}$, the \emph{group completion}).
    \item Let $\mathsf{FAb}\subset s\mathsf{CMon}^{\mathrm{gp}}$ denote the full subcategory spanned by finitely generated free abelian groups, which is an idempotent complete Lawvere theory. Let $s\mathsf{Ab}$ be the \infcat of algebras for this theory, then the inclusion $\mathsf{FAb}\subset s\mathsf{CMon}^{\mathrm{gp}}$ induces an equivalence of \infcats $s\mathsf{Ab}\simeq s\mathsf{CMon}^{\mathrm{gp}}$.
    \item There is a functor $\spect^{\geq 0}\rightarrow s\mathsf{Ab}$ in $\prl_{\mathrm{Proj}}$ fitting into a pushout diagram 
    \[
    \begin{tikzcd}
    \mathsf{Mon}_{\einfty} \ar[d,"(\_)^{\mathrm{gp}}"] \ar[r] & s\mathsf{CMon}\ar[d,"(\_)^{\mathrm{gp}}"] \\
    \spect^{\geq0} \ar[r] & s\mathsf{Ab}
    \end{tikzcd}
    \]
    in $\prl_{\mathrm{Proj}}$.
\end{enumerate}
\end{prop}
\begin{proof}
\begin{enumerate}[$(1)$]
   \item To see that the inclusion $s\mathsf{CMon}^{\mathrm{gp}}\subset s\mathsf{CMon}$ admits a right adjoint, let $\pi_0(A)^{\times}\subset\pi_0(A)$ be the submonoid on the invertible elements of $\pi_0(A)$, that is, the largest subgroup contained in $\pi_0(A)$, and consider the pullback diagram 
   \[
   \begin{tikzcd}
   A^{\times} \ar[d] \ar[r] & A\ar[d] \\
   \pi_0(A)^{\times} \ar[r] & \pi_0(A),
   \end{tikzcd}
   \]
   in $s\mathsf{CMon}$, then $A^{\times}$ is clearly grouplike and for each grouplike simplicial commutative monoid $B$, the map of spaces
   \[ \Hom_{s\mathsf{CMon}}(B,A^{\times}) \longrightarrow \Hom_{s\mathsf{CMon}}(B,A) \]
   is a pullback of the map of sets 
   \[ \Hom_{\mathsf{CMon}}(\pi_0(B),\pi_0(A)^{\times}) \longrightarrow \Hom_{\mathsf{CMon}}(\pi_0(B),\pi_0(A)), \]
   which is a bijection as the operation $(\_)^{\times}:\mathsf{CMon}\rightarrow\mathsf{Ab}$ is right adjoint to the inclusion of abelian groups into commutative monoids. Thus, the inclusion of connected components $A^{\times}\rightarrow A$ exhibits a colocalization.
    \item Consider the full subcategory $s\mathsf{CMon}^{\geq 1}\subset s\mathsf{CMon}$ spanned by objects $A$ for which the underlying space is 1-connective, which is stable under colimits. It follows from the previous lemma that the underlying space functor $s\mathsf{CMon}\rightarrow\spa$ carries the initial object to a final object, but as the underlying space functor reflects limits, we conclude that $s\mathsf{CMon}$ is pointed, so we have a suspension/looping adjunction
    \[ \begin{tikzcd}  s\mathsf{CMon} \ar[r,shift left,"\Sigma"] &  s\mathsf{CMon}. \ar[l,shift left,"\Omega"]  \end{tikzcd}  \]
    It follows from the previous lemma that the underlying space of $\Sigma A$ is the colimit of the Bar construction $|\mathsf{Bar}_A(*,*)_{\bullet}|$ which is $1$-connective. Unwinding the definitions, we can identify the functor 
    \[ s\mathsf{CMon}\longrightarrow\fun(\simpop,\spa),\quad A\longmapsto \mathsf{Bar}_A(*,*)_{\bullet} \]
    with the composition
    \[  s\mathsf{CMon}\longrightarrow\mathsf{Mon}_{\einfty}\longrightarrow \mathsf{Mon}_{\mathbb{E}_1}\subset\fun(\simpop,\spa). \]
    We conclude that $A$ is a grouplike simplicial commutative monoid if and only if the simplicial object $\mathsf{Bar}_A(*,*)_{\bullet}$ is a group object. Moreover, since all groups are effective in $\spa$ and in $s\mathsf{CMon}$, the functor $A\mapsto |\mathsf{Bar}_A(*,*)_{\bullet}|$ carries grouplike simplicial commutative monoids to \v{C}ech nerves. The functor $\Omega$ factor as
    \[  s\mathsf{CMon}\overset{\simeq}{\longrightarrow} \fun(\simpopplus,s\mathsf{CMon})'\subset \fun(\simpopplus,s\mathsf{CMon}) \overset{\ev_{[1]}}{\longrightarrow} s\mathsf{CMon},  \]
    where $ \fun(\simpopplus,s\mathsf{CMon})'$ denotes the full subcategory spanned by \v{C}ech nerves $U_{\bullet}$ such that $U_0\simeq *$. The first equivalence restricts to one $s\mathsf{CMon}^{\geq 1}\simeq \mathsf{Grp}^{+}(s\mathsf{CMon})$ between $1$-connective objects and \v{C}ech nerves $U_{\bullet}$ with $U_0\simeq *$ that are colimit diagrams. Let $U_{\bullet}$ be a \v{C}ech nerve with $U_0\simeq *$, then $U_{\bullet}|_{\simpop}$ is a group object in $\mathsf{sCMon}$; then $\pi_0(U_1)$ is a group so that $U_{1}$ is grouplike, since group object in commutative monoids are abelian groups by the classical Eckmann-Hilton argument. It follows that the adjunction $(\Sigma\adj\Omega)$ restricts to give an adjunction 
     \[ \begin{tikzcd}  s\mathsf{CMon}^{\mathrm{gp}} \ar[r,shift left,"\Sigma"] &  s\mathsf{CMon}^{\geq 1} \ar[l,shift left,"\Omega"]  \end{tikzcd}  \]
     which is an equivalence: if $A$ is grouplike, then the Bar construction $U_{\bullet}:=\mathsf{Bar}_A(*,*)_{\bullet}$ is a \v{C}ech nerve so the canonical map $A=U_1\rightarrow *\times_{U_{-1}}*$ is an equivalence. Conversely, let $B$ be a $1$-connective object and $V_{\bullet}$ the \v{C}ech nerve of $*\rightarrow B$, then we should show that the canonical map $\Sigma V_1\rightarrow B$ is an equivalence. Let $V_{\bullet}'$ denote a right Kan extension of the diagram $W:(\simpopplus)_{\leq 1}\rightarrow\mathsf{sCMon}$ given by 
     \[\begin{tikzcd} V_1\ar[r,shift left]\ar[r,shift right] & * \ar[r] & \Sigma V_1  \end{tikzcd}\]
     along the inclusion $(\simpopplus)_{\leq 1}\subset\simpopplus$, then we have an induced map $\alpha:V'_{\bullet}\rightarrow V_{\bullet}$ which restricts to the identity on $\simpop_{\leq 1}$. Because $V_1$ is grouplike, the diagram $W$ is a right Kan extension of $W|_{(\simpopplus)_{\leq 0}}$, which implies by \cite{HTT}, Proposition 4.3.2.8 that $V_{\bullet}'$ is a \v{C}ech nerve. We conclude that $\alpha|_{\simpop}$ is a morphism of group objects such that $\alpha_1$ is the identity, but this implies that $\alpha$ is an equivalence. It follows that the composition
     \[ s\mathsf{CMon}\overset{\Sigma}{\longrightarrow}  s\mathsf{CMon}^{\geq 1}\underset{\simeq}{\overset{\Omega}{\longrightarrow}}   s\mathsf{CMon}^{\mathrm{gp}}  \]
     is a left adjoint to the inclusion. 
    \item It follows from $(1)$, $(2)$ and \cite{HA}, Proposition 7.1.4.12 that it suffices to argue that the essential image of $(\_)^{\mathrm{gp}}$ on $\mathsf{FCMon}$ consists of finitely generated free abelian groups. On underlying spaces, we can identify the map $\Sigma\Z_{\geq 0}\rightarrow \Sigma \Z$ induced by the inclusion $\Z_{\geq0}\rightarrow \Z$ with the map $\beta:B\icat \rightarrow B\icatd$ of classifying spaces, where $\icat$ and $\icatd$ are the single object categories with space of morphisms $\Z_{\geq0}$ and $\Z$ respectively. We can identify the fibre product $\icat\times_{\icatd}\icatd_{*/}$ with the poset category $\Z$, which has contractible classifying space, so we deduce that $\beta$ is an equivalence by Quillen's Theorem A. It follows that the unit of the group completion $\Z_{\geq0}\rightarrow (\Z_{\geq0})^{\mathrm{gp}}$ is equivalent to $\Z_{\geq0}\rightarrow \Z$. 
    \item Consider the functor $s\mathsf{CMon}\rightarrow \mathsf{Mon}_{\einfty}$ induced by the transformation of algebraic theories $\F\rightarrow\mathsf{FCMon}$, then $A\in s\mathsf{CMon}$ is grouplike if and only if the associated $\einfty$-space is grouplike, but we have an equivalence $\mathsf{Mon}_{\einfty}^{\mathrm{gp}}\simeq\spect^{\geq 0}$ (by \cite{HA}. Remark 5.2.6.26 for instance), so we have a pullback diagram 
    \[
    \begin{tikzcd}
    s\mathsf{Ab} \ar[d] \ar[r] & s\mathsf{CMon} \ar[d] \\
    \spect^{\geq 0} \ar[r]& \mathsf{Mon}_{\einfty}
    \end{tikzcd}
    \]
    of \infcats and conservative functors preserving limits and colimits.
\end{enumerate}
\end{proof}
\begin{rmk}\label{rmk:units}
We have seen that for any simplicial $\cinfty$-ring $A$, the space 
\[A_{\geq 0}:=\Hom_{\sring}(\cinfty(\R_{\geq 0}),A)\]
admits a natural structure of a simplicial commutative monoid. The simplicial $\cinfty$-ring $A$ evidently also admits the structure of a simplicial commutative monoid; the requisite forgetful functor $(\_)^{\mathsf{Mon}}:\sring\rightarrow s\mathsf{CMon}$ can be defined in two (naturally equivalent) ways: we can define a functor $\mathsf{FCMon}^{op}\rightarrow\cartsp$ via the same formulae that appear in Construction \ref{cons:logcinftyrings}, or we can take the composition $\sring\overset{(\_)^{\rmalg}}{\rightarrow}\calg_{\R}^{\geq0}\rightarrow s\mathsf{CMon}$ where the second functor is induced by the lax monoidal functor $\mathsf{Mod}^{\geq0}_{\R}\rightarrow \spa$. The functors $(\_)_{\geq0}$ and $(\_)^{\mathsf{Mon}}$ may be combined by defining a functor $\Theta:\mathsf{FCMon}^{op}\times\Delta^1\rightarrow \cartsp_c$ as follows.
\begin{enumerate}[$(1)$]
    \item $\Theta$ carries the object $(\Z_{\geq0}^n,0)$ to $\R^n_{\geq 0}$ and the object $(\Z_{\geq0}^m,1)$ to $\R^m$.
    \item $\Theta$ carries a morphism $f:(\Z_{\geq0}^n,0)\leftarrow (\Z_{\geq0}^m,0)$ determined by an $m$-tuple $\{(k^i_1,\ldots,k^i_n)\in \Z^n_{\geq 0}\}_{1\leq i\leq m}$ to the map $\R^n_{\geq 0}\rightarrow\R^m_{\geq 0}$ given by 
     \[ (x_1,\ldots,x_n) \longmapsto \left(\prod_{1\leq j\leq n}x^{k^1_j}_{j},\ldots,\prod_{1\leq j\leq n}x^{k^m_j}_{j}\right).\]
    The morphisms $(\Z_{\geq0}^n,1)\leftarrow(\Z_{\geq0}^m,1)$ and $(\Z_{\geq0}^n,0)\leftarrow(\Z_{\geq0}^m,1)$ are carried to morphisms $\R^n\rightarrow\R^m$ and $\R^n_{\geq0}\rightarrow\R^m$ respectively, defined by the same formula.
\end{enumerate}
Composing $\iota_{c*}$ with $\Theta^*:\fun(\cartsp_c,\spa)\rightarrow\fun(\mathsf{FCMon}^{op}\times\Delta^1,\spa)$ yields a natural transformation $\sring\rightarrow\fun(\Delta^1,s\mathsf{CMon})$ that lifts, for each $A\in \sring$, the map of spaces $A_{\geq 0}\rightarrow A$ induced by the map $\cinfty(\R)\rightarrow \cinfty(\R_{\geq 0})$ to a map of simplicial commutative monoids $A_{\geq 0}\rightarrow A^{\mathsf{Mon}}$. In Remark \ref{rmk:truncatepositivelts} we argued that the natural map 
\[\pi_0(A_{\geq0})\longrightarrow\pi_0(A)_{\geq 0}=\Hom_{\cinfty\mathsf{ring}}(\cinfty(\R_{\geq 0}),\pi_0(A))\]
is an equivalence. Since the map $\cinfty(\R)\rightarrow \cinfty(\R_{\geq 0})$ is a regular epimorphism of $\cinfty$-rings, we have an injection $\pi_0(A)_{\geq 0}\hookrightarrow\pi_0(A)$, which is obtained by applying the functor taking connected components to the map $A_{\geq0}\rightarrow A$ and the isomorphism $\pi_0(A_{\geq 0})\cong \pi_0(A)_{\geq 0}$. We conclude that the commutative monoid structure on $\pi_0(A)$ restricts to one on the subset $\pi_0(A_{\geq 0})$, and this latter structure then coincides with the one coming from the simplicial commutative monoid structure on $A_{\geq 0}$ defined in Construction \ref{cons:logcinftyrings}. We use this observation to identify the group of units of $\pi_0(A_{\geq 0})$: the group $\pi_0(A)^{\times}\hookrightarrow \pi_0(A)$ of invertible elements coincides with the map $\Hom_{\cinfty\mathsf{ring}}(\cinfty(\R\setminus\{0\}),A)\rightarrow \Hom_{\cinfty\mathsf{ring}}(\cinfty(\R),A)$ by definition of the localization. Thus, if $x\in \pi_0(A_{\geq 0})$ is invertible as an element in $\pi_0(A)$ we have a commuting diagram 
\[
\begin{tikzcd}
\cinfty(\R) \ar[d] \ar[r] & \cinfty(\R_{\geq 0}) \ar[d,"x"] \\
\cinfty(\R\setminus\{0\}) \ar[r]& \pi_0(A)
\end{tikzcd}
\]
of $\cinfty$-rings, so the map classifying $x$ factors through the pushout $\cinfty(\R_{>0})\rightarrow \pi_0(A)$, which shows that the inverse of $x$ lies in the submonoid $\pi_0(A_{\geq 0})$. It follows that the group of units $\pi_0(A_{\geq 0})^{\times}$ is given by a pullback $\pi_0(A_{\geq 0})\times_{\pi_0(A)}\pi_0(A)^{\times}$, so the monomorphism $A_{\geq 0}^{\times}\rightarrow A_{\geq0}$ fits as the left vertical map into a pullback diagram
\[
\begin{tikzcd}
A_{\geq 0}^{\times} \ar[d] \ar[r] & \Hom_{\sring}(\cinfty(\R\setminus\{0\}),A) \ar[d] \\
A_{\geq 0} \ar[r] & \Hom_{\sring}(\cinfty(\R),A).
\end{tikzcd}
\]
As a result, this map coincides with the map $\Hom_{\sring}(\cinfty(\R_{>0}),A)\rightarrow \Hom_{\sring}(\cinfty(\R_{\geq0}),A)$. Thus, the map $A_{>0}\rightarrow A_{\geq 0}$ of Remark \ref{rmk:cinftycorners} coincides with the inclusion of the $\infty$-group of units $A_{\geq 0}^{\times}\hookrightarrow A_{\geq 0}$.
\end{rmk}
\begin{rmk}\label{rmk:0truncatedplrings}
We give one more application of Remark \ref{rmk:truncatepositivelts}. Abusing notation, we denote $(\_)_{\geq0}:\cinfty\mathsf{ring}\rightarrow \mathsf{CMon}$ for the functor given by $A\mapsto \Hom_{\cinfty\mathsf{ring}}(\cinfty(\R_{\geq 0}),A)$, and define a category $\cinfty\mathsf{PLog}$ as the pullback $\cinfty\mathsf{ring}\times_{\mathsf{CMon}}\fun(\Delta^1,\mathsf{CMon})$. We have a diagram 
\[
\begin{tikzcd}
\cinfty\mathsf{ring} \ar[d,hook] \ar[r,"(\_)_{\geq0}"] & \mathsf{CMon} \ar[d,hook] \\
\sring\ar[r,"(\_)_{\geq0}"] & s\mathsf{CMon}
\end{tikzcd}
\]
which commutes up to canonical homotopy, determining a fully faithful functor $g:\cinfty\mathsf{PLog}\rightarrow \splring$. The vertical maps admit left adjoint functors denoted by $\pi_0$ and the associated Beck-Chevalley transformation at an object $A\in\sring$ is obtained by applying $\pi_0$ to the map of simplicial commutative monoids $f:A_{\geq 0}\rightarrow \pi_0(A)_{\geq 0}$ induced by the unit map $A\rightarrow A_{\geq 0}$. It follows from Remark \ref{rmk:truncatepositivelts} that $f$ exhibits a 0-truncation, so we have an equivalence $\pi_0((\_)_{\geq 0})\simeq \pi_0(\_)_{\geq0}$ which provides a left adjoint $\pi_0:\splring\rightarrow\cinfty\mathsf{PLog}$ to $g$. This adjunction is equivalent to the $0$'th truncation $\tau_{\leq 0}$. To see this, it suffices to show that an object $(A,M\rightarrow A_{\geq0})$ is $0$-truncated if and only if it lies in the essential image of $g$ which is easily seen to consist of those objects $(B,N\rightarrow B_{\geq0})$ where $B$ is a 0-truncated simplicial $\cinfty$-ring and $N$ is a 0-truncated simplicial commutative monoid. The `only if' direction follows immediately from the fact that both $p:\splring\rightarrow \sring$ and $\splring\rightarrow\fun(\Delta^1,s\mathsf{CMon})$ preserve limits. For the `if' direction, we suppose that $A$ and $M$ are $0$-truncated, then we have for any $(B,N\rightarrow B_{\geq0})\in \splring$ and any map $f:B\rightarrow A$ a fibre sequence 
\[\Hom_{(s\mathsf{CMon})/_{B_{\geq 0}}}(N,M\times_{A_{\geq 0}}B_{\geq 0}) \longrightarrow \Hom_{\splring}((B,N\rightarrow B_{\geq0}),(A,M\rightarrow A_{\geq0})) \longrightarrow\Hom_{\sring}(B,A)  \]
since $p:\splring\rightarrow\sring$ is a Cartesian fibration. To conclude that $\Hom_{\splring}((B,N\rightarrow B_{\geq0}),(A,M\rightarrow A_{\geq0}))$ is $0$-truncated, it suffices to argue that the base and the fibre spaces are $0$-truncated. As $A$ is 0-truncated, the base space is also 0-truncated and as $M\rightarrow A_{\geq 0}$ is a 0-truncated morphism, the map $M\times_{A_{\geq 0}}B_{\geq 0}\rightarrow B_{\geq 0}$ is too so the fibre is also 0-truncated. \\
Using an analogous argument, it can be shown that $\Xi:\sring_{pc}\rightarrow\splring$ takes $n$-truncations to $n$-truncations for all $n\geq 0$, that is, the relevant Beck-Chevalley map provides an equivalence $\tau_{\leq n}\circ \Xi\simeq \Xi^n\circ \tau_{\leq n}$, where $\Xi^n$ is the functor $\tau_{\leq n}\sring_{pc}\rightarrow\tau_{\leq n}\splring$ induced by $\Xi$.
\end{rmk}

\begin{defn}
Let $A$ be a simplicial commutative monoid and let $M\in (s\mathsf{CMon})_{/A}$ be a prelog structure on $A$, then $M$ is a \emph{log structure} on $A$ if the upper horizontal map in the pullback diagram
\[
\begin{tikzcd}
M\times_{A}A^{\times} \ar[d] \ar[r] & A^{\times} \ar[d] \\
M\ar[r] & A
\end{tikzcd}
\]
is an equivalence, where the right vertical map is the counit of the coreflective embedding $s\mathsf{Ab}\subset s\mathsf{CMon}$, that is, the inclusion of connected components determined by the invertible elements in the commutative monoid $\pi_0(A)$. We denote by $\mathsf{Log}_{A}\subset (s\mathsf{CMon})_{/A}$ the full subcategory spanned by log structures and $\slring\subset\splring$ the full subcategory spanned by objects $(A,M\rightarrow A_{\geq0})$ such that the prelog structure $M$ is a log structure on $A_{\geq 0}$.
\end{defn}
\begin{rmk}
A prelog structure $M\rightarrow A$ is a log structure if and only if the canonical maps $M^{\times}\rightarrow A^{\times}$ and $M^{\times}\rightarrow M\times_{A}A^{\times}$ are both equivalences.
\end{rmk}
The following proposition is an immediate consequence of Remarks \ref{rmk:cinftycorners} and \ref{rmk:units}.
\begin{prop}\label{prop:logcinfty}
The functor $\Xi:\sring_{pc}\rightarrow\splring$ restricted to $\sring_c$ takes values in $\slring$. Denoting the resulting functor $\sring_c\rightarrow\splring$ by $\Xi_{\mathsf{Log}}$, the commuting diagram 
\[
\begin{tikzcd}
\sring_c\ar[d,hook] \ar[r,"\Xi_{\mathsf{Log}}"] & \slring \ar[d,hook] \\
\sring_{pc}\ar[r,"\Xi"] & \splring
\end{tikzcd}
\]
is a homotopy pullback diagram of \infcatst.
\end{prop}
The construction $(A,A_c)\mapsto (A,A_c\rightarrow A_{\geq 0})$ implemented by the functors $\Xi$ and $\Xi_{\mathsf{Log}}$ is obviously conservative. $\Xi$ and $\Xi_{\mathsf{Log}}$ also preserve limits and sifted colimits (as we will show shortly) so we might like to interpret them as forgetful functors. The notion of a simplicial $\cinfty$-ring with corners appears prima facie strictly more structured than a positive prelog simplicial $\cinfty$-ring, as $\Xi$ forgets the $\cinfty$ information contained in $A_c$. When we restrict to logarithmic structures however, we see that there is no loss of information at all.
\begin{thm}\label{thm:logcinftycorners}
The functor $\Xi_{\mathsf{Log}}:\sring_c\rightarrow \slring$ is an equivalence of \infcatst.
\end{thm}
As we will see, this result grants us control over the coCartesian morphisms of $\iota_c^*$, which reduces the computation of limits and colimits in $\sring_c$ to limits and colimits in $\sring$ and in \infcats of log structures. The proof of Theorem \ref{thm:logcinftycorners} requires a few prelimenaries. Our first order of business is to understand the relative left adjoint to the inclusion $\slring\subset\splring$. The following result is familiar from the usual theory of log structures on monoids, albeit that the proof is somewhat more involved since we do not take recourse to point-set arguments.
\begin{prop}\label{prop:logaccessible}
Denote by $p_{\mathsf{Log}}$ the composition $\slring\subset\splring\overset{p}{\rightarrow}\slring$.
\begin{enumerate}[$(1)$]
    \item The functor $p_{\mathsf{Log}}$ is a Cartesian fibration and the inclusion $\slring\hookrightarrow\splring$ carries Cartesian edges to Cartesian edges.
    \item For each simplicial commutative monoid $A$, the fully faithful inclusion $\mathsf{Log}_A\subset(s\mathsf{CMon})_{/A}$ preserves sifted colimits and admits a left adjoint. 
    \item The inclusion $\slring\subset \splring$ admits a left adjoint relative to $\sring$.
\end{enumerate}
\end{prop}
\begin{proof}
The proof of $(1)$ amounts to the assertion that if $M\rightarrow A$ is a log structure on a simplicial commutative monoid $A$ and $B\rightarrow A$ is any morphism of simplicial commutative monoids, then $B\times_{A}M$ is a log structure on $B$, which is a straightforward check. For $(2)$, let $M\rightarrow A$ be a log structure, and consider the pushout diagram of prelog structures over $A$:
\[
\begin{tikzcd}
M\times_AA^{\times} \ar[d] \ar[r] & A^{\times} \ar[d] \\
M \ar[r,"f"] & N.
\end{tikzcd}
\]
It suffices to show that $N$ is a log structure over $A$ and that restriction along the morphism $f$ induces, for each log structure $M'\rightarrow A$, an equivalence
\[\Hom_{(s\mathsf{CMon})_{/A}}(N,M') \overset{\simeq}{\longrightarrow}  \Hom_{(s\mathsf{CMon})_{/A}}(M,M').\]
The following assertion will enjoy verification at the end of the proof.
\begin{enumerate}
    \item[$(*)$] The diagram 
\[
\begin{tikzcd}
0\ar[r]\ar[d] & 0\ar[d]  \\
N \ar[r] & A.
\end{tikzcd}
\]
is a pullback square of simplicial commutative monoids.
\end{enumerate}
We have a diagram
\[
\begin{tikzcd}
 N^{\times} \ar[d]\ar[r,"g"] & N\times_AA^{\times} \ar[r,"h"]\ar[d] & A^{\times}\ar[d] \\
 N\ar[r,equal] & N\ar[r] & A.
\end{tikzcd}
\]
Both maps $N^{\times}\rightarrow N$ and $N\times_AA^{\times}\rightarrow N$ are inclusions of connected components, so the map $N^{\times}\rightarrow N\times_AA^{\times}$ is one as well. We first show that the map $h\circ g:N^{\times}\rightarrow A^{\times}$ is an equivalence. Consider the diagram
\[
\begin{tikzcd}
0\ar[r]\ar[d] & 0\ar[r]\ar[d] & 0\ar[d]  \\
N^{\times}\ar[r] & N\times_AA^{\times} \ar[d]\ar[r] & A^{\times}\ar[d] \\
& N \ar[r] & A.
\end{tikzcd}
\]
The right upper square is a pullback diagram since the right outer rectangle is one, by $(*)$. Because the map $N^{\times}\rightarrow N\times_AA^{\times}$ is an inclusion of connected components, the upper rectangle is also a pullback diagram of simplicial abelian groups, and therefore also a pullback diagram of connective spectra. Since the map $N^{\times}\rightarrow A^{\times}$ is an effective epimorphism (i.e. 0-connective), the upper rectangle is also a pullback diagram of spectra. Then it is a pushout diagram, so the map $N^{\times}\rightarrow A^{\times}$ is an equivalence. It follows that $A^{\times}$ is a retract of $N\times_AA^{\times}$. Choose an element $x\in \pi_0(N\times_AA^{\times})$, then $h(x)$ is invertible in $\pi_0(A^{\times})$ so admits an inverse $y$. Consider the element $z:=g((h\circ g)^{-1}(y))\in \pi_0(A)$, then $h(x+z)=h(x)+h(g((h\circ g)^{-1}(y)))=h(x)+y$, which is the unit. By $(*)$, we have $h^{-1}(0)=0$, so $z$ is an inverse of $x$. It follows that $g$ and therefore also $h$ is an equivalence. We now proceed by showing that the map on morphism spaces induced by restricting along $f$ induces an equivalence for each log structure $M'\rightarrow A$. The relevant map is the the upper horizontal one in a pullback diagram 
\[
\begin{tikzcd}
\Hom_{(s\mathsf{CMon})_{/A}}(N,M') \ar[d] \ar[r] &   \Hom_{(s\mathsf{CMon})_{/A}}(M,M') \ar[d] \\
\Hom_{(s\mathsf{CMon})_{/A}}(A^{\times},M') \ar[r] & \Hom_{(s\mathsf{CMon})_{/A}}(M\times_AA^{\times},M')
\end{tikzcd}
\]
of spaces, so it suffices to argue that the lower horizontal map is an equivalence. In fact, we claim that both the domain and codomain of this map are weakly contractible. Note that both $A^{\times}$ and $M\times_AA^{\times}$ lie in the image of the functor $(s\mathsf{CMon})_{/A^{\times}}\rightarrow (s\mathsf{CMon})_{/A}$. From the adjunction \[\begin{tikzcd}(s\mathsf{CMon})_{/A^{\times}}\ar[r,shift left]& (s\mathsf{CMon})_{/A}\ar[l,shift left]\end{tikzcd}\] 
we have for any $L\rightarrow A^{\times}$ an equivalence $\Hom_{(s\mathsf{CMon})_{/A}}(L,M') \simeq 
\Hom_{(s\mathsf{CMon})_{/A^{\times}}}(L,M'\times_AA^{\times})$. As $M'$ is a log structure, the map $M'\times_AA^{\times}\rightarrow A^{\times}$ is an equivalence, so $M'\times_AA^{\times}$ is a final object in the \infcat of prelog structures over $A^{\times}$, which proves our claim. We have constructed a left adjoint to the inclusion $\mathsf{Log}_A\subset (s\mathsf{CMon})_{/A}$, so $(2)$ follows from the observation that this inclusion is stable under sifted colimits, as sifted colimits are universal in $s\mathsf{CMon}$.\\
It is an immediate consequence of $(1)$, $(2)$ and \cite{HA}, Proposition 7.3.2.6 that the inclusion $\slring\subset\splring$ admits a left adjoint relative to $\sring$. \\
We are left to prove assertion $(*)$. The diagram 
\[
\begin{tikzcd}
M\times_{A}A^{\times} \ar[d] \ar[r] & A^{\times} \ar[d] \\
M \ar[r] & A
\end{tikzcd}
\]
of simplicial commutative monoids induces an $M\times_{A}A^{\times}$-bilinear (in the sense of \cite{HA}, section 4.4.4) map $M\times A^{\times}\rightarrow A$ which is encoded by the simplicial object $\mathsf{Bar}_{M\times_AA^{\times}}(M\times_AA^{\times},A^\times)_{\bullet}$ being equipped with an augmentation to $A$. The map $N\rightarrow A$ can be identified with the canonical map $|\mathsf{Bar}_{M\times_AA^{\times}}(M\times_AA^{\times},A^\times)_{\bullet}| \rightarrow A$. We have morphisms of simplicial objects \[\mathsf{Bar}_{M\times_AA^{\times}}(M,A^\times)_{\bullet}\overset{\alpha}{\longleftarrow} \mathsf{Bar}_{M\times_AA^{\times}}(M\times_AA^{\times},A^\times)_{\bullet} \overset{\beta}{\longrightarrow} \mathsf{Bar}_{M\times_AA^{\times}}(M\times_AA^{\times},0)_{\bullet}  \]
induced by the $(M\times_AA^{\times})$-module morphisms $M\times_AA^{\times}\rightarrow M$ and $A^{\times}\rightarrow 0$; in particular, for each $[n]\in \simp$, we have maps of spaces
\[ M\times (M\times_AA^{\times})^{\times n}\times A^{\times }\longleftarrow   M\times_AA^{\times} \times (M\times_AA^{\times})^{\times n}\times A^{\times} \longrightarrow M\times_AA^{\times} \times (M\times_AA^{\times})^{\times n}\times *, \]
where the left map is an inclusion of connected components and the right map projects away the factor $A^{\times}$. The map $\alpha:\mathsf{Bar}_{M\times_AA^{\times}}(M\times_AA^{\times},A^\times)_{\bullet}\rightarrow \mathsf{Bar}_{M\times_AA^{\times}}(M,A^\times)_{\bullet}$ fits as the left vertical map into a diagram 
\[
\begin{tikzcd}
\mathsf{Bar}_{M\times_AA^{\times}}(M\times_AA^{\times},A^\times)_{\bullet} \ar[d,"\alpha"] \ar[r] & A^{\times} \ar[d] \\
\mathsf{Bar}_{M\times_AA^{\times}}(M,A^\times)_{\bullet} \ar[r] & A.
\end{tikzcd}
\]
Since both vertical maps are inclusions of connected components in each simplicial level, it follows from an easy check on connected components that this diagram is a pullback diagram of simplicial objects. Since colimits are universal in spaces, it suffices to show that the colimit of the simplicial object defined as the cone in the pullback diagram
\[
\begin{tikzcd}
\mathsf{Bar}_{M\times_AA^{\times}}(M\times_AA^{\times},A^\times)_{\bullet}\times_{A^{\times}} 0 \ar[r] \ar[d] & 0\ar[d]\\
\mathsf{Bar}_{M\times_AA^{\times}}(M\times_AA^{\times},A^\times)_{\bullet}  \ar[r] & A^{\times}
\end{tikzcd}
\]
of simplicial objects is contractible. The map $A^{\times}\rightarrow 0$ induces a commuting diagram 
\[
\begin{tikzcd}
\mathsf{Bar}_{M\times_AA^{\times}}(M\times_AA^{\times},A^\times)_{\bullet}  \ar[r]\ar[d,"\beta"] & A^{\times} \ar[d]\\
\mathsf{Bar}_{M\times_AA^{\times}}(M\times_AA^{\times},0)_{\bullet} \ar[r] & 0.
\end{tikzcd}
\]
Since the left vertical map projects away the factor $A^{\times}$ in each simplicial level, this diagram is a pullback diagram. It follows that the composite map 
\[ \mathsf{Bar}_{M\times_AA^{\times}}(M\times_AA^{\times},A^\times)_{\bullet}\times_{A^{\times}}0\longrightarrow \mathsf{Bar}_{M\times_AA^{\times}}(M\times_AA^{\times},0)_{\bullet} \]
is an equivalence, as it is a pullback along the map $0\rightarrow 0$. We conclude by observing that the augmented simplicial object $\mathsf{Bar}_{M\times_AA^{\times}}(M\times_AA^{\times},0)_{\bullet}\rightarrow 0$ is a colimit diagram.
\end{proof}
\begin{rmk}\label{rmk:logificationfunctor}
The proof above gives an explicit description of the value of the left adjoint $\splring\rightarrow\slring$ on a prelog structure $M\rightarrow A_{\geq0}$ as the pushout
\[
\begin{tikzcd}
M\times_{A_{\geq 0}}A_{\geq0}^{\times} \ar[d] \ar[r] & A_{\geq 0}^{\times} \ar[d] \\
M \ar[r,"f"] & N.
\end{tikzcd}
\]
If a log structure $N\rightarrow A_{\geq0}$ fits into a pushout diagram as above, we say that $f$ \emph{exhibits $N$ as a logification of $M$} (with respect to some simplicial commutative monoid $A_{\geq0}$). We denote the resulting left adjoint, the \emph{logification functor} by $L_{\mathsf{Log}}:\splring\rightarrow\slring$. Also note that by virtue of Theorem \ref{thm:logcinftycorners}, a map $(A,A_c)\rightarrow (B,B_c)$ is a coCartesian arrow of $\sring_c$ precisely if $A_c\rightarrow B_c$ exhibits a logification in the \infcat of prelog structures over $B_{\geq 0}$.
\end{rmk}
At this point, we give criteria for the recognition of limits and colimits in $\splring$. 

\begin{lem}\label{lem:ploglimcolim}
The following hold true.
\begin{enumerate}[$(1)$]
    \item The functor $p:\splring\rightarrow\sring$ preserves all limits and colimits.
    \item The functor \[\ev_{\mathsf{PLog}}:\splring\longrightarrow\fun(\Delta^1,s\mathsf{CMon})\overset{\ev_0}{\longrightarrow}s\mathsf{CMon}\]
    preserves all limits and colimits.
\end{enumerate}
\end{lem}
\begin{proof}
Only the statements involving colimits are not immediate. Since $p$ is a presentable fibration over a presentable base, $p$ preserves colimits. Now choose a small diagram $f:K\rightarrow \splring$, then we wish to show that the map $\colim \ev_{\mathsf{PLog}}\circ f\rightarrow \ev_{\mathsf{PLog}}(\colim f)$ is an equivalence. Let $G$ denote the functor $\splring\rightarrow \fun(\Delta^1,s\mathsf{CMon})$. Since an edge $\Delta^1\rightarrow\fun(\Delta^1,s\mathsf{CMon})$ is $\ev_1$-coCartesian if and only if the composition $\Delta^1\rightarrow\fun(\Delta^1,s\mathsf{CMon})\overset{\ev_0}{\rightarrow} s\mathsf{CMon}$ is an equivalence, we are required to show that the map $\colim Gf\rightarrow G(\colim f)$ is a coCartesian edge of $\fun(\Delta^1,s\mathsf{CMon})$. This follows from Proposition \ref{prop:colimleftadjcocart}.
\end{proof}
\begin{cor}\label{cor:ploglimcolimcons}
The functor $p\times \ev_{\mathsf{PLog}}:\splring\rightarrow \sring\times s\mathsf{CMon}$ is conservative and preserves all limits and colimits.
\end{cor}
\begin{cor}\label{cor:slringcompactlygen}
The inclusion $\slring\subset\splring$ preserves filtered colimits; in other words, the localization $L_{\mathsf{Log}}$ is $\omega$-accessible.
\end{cor}
\begin{proof}
Let $\mathcal{J}:K\rightarrow \slring$ be filtered diagram and denote by $(A,M\rightarrow A_{\geq0})$ a colimit of $\mathcal{J}$. According to the proof of Lemma \ref{lem:ploglimcolim}, we have a commuting diagram 
\[
\begin{tikzcd}
\colim_{i\in\mathcal{J}} M_i \ar[d] \ar[r,"\simeq"] & M\ar[d]\\
\colim_{i\in\mathcal{J}} (A_i)_{\geq 0} \ar[r] & A_{\geq 0} 
\end{tikzcd}
\]
in $s\mathsf{CMon}$. Since the functor $s\mathsf{CMon}_{/(A_{\geq 0}^{\times}\rightarrow A_{\geq 0})}\rightarrow s\mathsf{CMon}_{/A_{\geq 0}}$ is fully faithful and for each $i\in K$, the composition $(A_i)^{\times}_{\geq 0}\rightarrow A_{\geq 0}$ factors through $A_{\geq 0}^{\times}$, there is a map $\colim_{i\in\mathcal{J}}(A_i)_{\geq 0}^{\times}\rightarrow A_{\geq 0}^{\times}$ fitting into a commuting diagram 
\[
\begin{tikzcd}
\colim_iM_i\times_{(A_i)_{\geq 0}} (A_i)^{\times}_{\geq 0} \ar[d] \ar[r,"\simeq"]&  \colim_i (A_i)^{\times}_{\geq 0} \ar[d,"\gamma"] \ar[r] &  A^{\times}_{\geq 0}\ar[d] \\
\colim_i M_i \ar[r] & \colim_i (A_i)_{\geq 0}\ar[r] & A_{\geq 0} 
\end{tikzcd}
\]
where the indicated map $\gamma$ is an inclusion of connected components. Here, the upper left map is a colimit of equivalences and therefore an equivalence. As filtered colimits commute with finite limits in $s\mathsf{CMon}$, the left square is a pullback. It suffices to show that the right square is a pullback and that the upper right map is an equivalence. Using the natural transformation $A_{\geq 0}\rightarrow A$ of Remark \ref{rmk:units}, we deduce the existence of a commuting diagram 
\[
\begin{tikzcd}
\colim_i (A_i)^{\times}_{\geq 0} \ar[d,"\gamma"] \ar[r] &  A^{\times}_{\geq 0} \ar[d] \\
\colim_i (A_i)_{\geq 0}\ar[r]\ar[d,"\delta"] & A_{\geq 0} \ar[d] \\
\colim_iA_i\ar[r,equal] & A
\end{tikzcd}
\]
As the composition $\delta\circ\gamma$ is also an inclusion of components, it is clear that the outer square is a pullback so that the upper horizontal map is an equivalence. We will be done once we show that $x \in \colim_i (A_i)_{\geq 0}$ lies in the image of $\gamma$ if and only if $x$ is invertible in $A$. Suppose $x$ factors through some $(A_i)_{\geq 0}$. The `only if' direction is obvious, and in the other direction we see that there must be some $j$ such that $x^{-1}\in A_j$ and some upper bound $k$ for $\{i,j\}\subset K$ such that $x$ is invertible in $A_k$. Since we have $(A_k)^{\times}_{\geq 0}\simeq (A_k)_{\geq0}\times_{A_k}A_k^{\times}$, $x$ is also invertible in $(A_k)_{\geq 0}$, so $x$ lies in the image of $\gamma$ as required. 
\end{proof}
\begin{prop}\label{prop:xiproperties}
The functor $\Xi$ of Construction \ref{cons:logcinftyrings} has the following properties.
\begin{enumerate}[$(1)$]
    \item $\Xi$ is conservative.
    \item $\Xi$ preserves limits and sifted colimits. 
    \item $\Xi$ is monadic.
    \item Let $\Upsilon$ be a left adjoint to $\Xi$, then for each $(A,M\rightarrow A_{\geq 0})$, the unit map $(A,M\rightarrow A_{\geq 0})\rightarrow \Xi\Upsilon(A,M\rightarrow A_{\geq 0})$ maps to an equivalence under $p$.
    \item $\Xi$ is a left Kan extension and a $p$-left Kan extension of its restriction to the image of $j:\cartsp_c^{op}\hookrightarrow\sring_c$.
\end{enumerate}
\end{prop}
\begin{proof}
Consider the functor $\rho:\splring\rightarrow \spa\times\spa$ obtained by taking the product of the functors \[\splring\overset{p}{\longrightarrow} \sring\overset{\ev_{\{\R\}}}{\longrightarrow}\spa\]
and 
\[\splring\overset{\ev_{\mathsf{PLog}}}{\longrightarrow} s\mathsf{CMon}\overset{\ev_{\Z_{\geq 0}}}{\longrightarrow}\spa.\]
Then $\rho$ is conservative and preserves limits and sifted colimits by Lemma \ref{lem:ploglimcolim}. We have a commuting diagram 
\[
\begin{tikzcd}
\sring_{pc} \ar[rr,"\Xi"] \ar[dr,"\ev_{\R}\times\ev_{\R_{\geq 0}}"'] && \splring \ar[dl,"\rho"] \\
& \spa\times\spa
\end{tikzcd}
\]
of \infcatst. Since the left diagonal map is conservative and preserves limits and sifted colimits, we deduce $(1)$ and $(2)$. Note that $(3)$ is an immediate consequence $(1)$ and $(2)$, the presentability of both $\sring_{pc}$ and $\splring$ and Lurie's Barr-Beck theorem. To prove $(4)$, \cite{HA}, Proposition 7.3.2.6 guarantees that it suffices to show that for each $A\in \sring$, the functor $\Xi_A$ between the fibres at $A$ admits a left adjoint since $\Xi$ preserves Cartesian edges. Suppose $q:\icat\rightarrow\icatd$ is a presentable fibration and $D\in \icatd$ an object, then a diagram $K^{\rhd}\rightarrow \icat_D$ where $K$ is weakly contractible is a colimit diagram if and only if it is a $q$-colimit diagram if and only if it is a colimit diagram in $\icat$. As $\Xi$ preserves sifted colimits, it follows that $\Xi_A$ also preserves sifted colimits for each $A\in\sring$. To conclude that $\Xi_A$ admits a left adjoint, it suffices to prove that $\Xi_A$ preserves limits, by the adjoint functor theorem and the presentability of the fibres. This follows from the following relative version of assertion $(*)$ of Proposition \ref{prop:cornerpresfib}, the proof of which uses the same techniques and is left to the reader. \begin{enumerate}
    \item[$(**)$] Let $p:\icat\rightarrow\icatd$ and $q:\icat'\rightarrow\icatd$ be coCartesian fibrations among \infcats and let $f:\icat\rightarrow\icat'$ be a morphism in $\mathsf{coCart}_{\icatd}$. Let $K$ be a simplicial set and let $g:K\rightarrow\icat_D$ be a diagram in the fibre over some object $D\in\icatd$. Let $i_D:\icat_D\subset\icat$ denote the inclusion, and suppose that the induced diagram $i_Dg:K\rightarrow \icat$ admits a colimit and that $p$ and $f$ preserve the colimit of $i_Dg$. Then the diagram $g$ admits a colimit and the functor $f_{D}:\icat_D\rightarrow\icat_D'$ preserves this colimit.
\end{enumerate}
Note that $(5)$ follows immediately from \cite{HTT}, Proposition 5.5.8.15.
\end{proof}
\begin{cor}\label{cor:accessiblecorner}
The localization $\sring_c\subset\sring_{pc}$ is an $\omega$-accessible localization. In particular, $\sring_c$ is compactly generated.
\end{cor}
\begin{proof}
It suffices to show that the inclusion $\sring_c\subset\sring_{pc}$ preserves filtered colimits. To see this, combine Propositions \ref{prop:logcinfty} and \ref{prop:xiproperties} and Corollary \ref{cor:slringcompactlygen}. 
\end{proof}
The functor $(\_)_{\geq 0}:\sring\rightarrow s\mathsf{CMon}$ does not preserve filtered colimits (it only preserves $\kappa$-filtered colimits for regular cardinals $\kappa$ for which $\cinfty(\R_{\geq 0})$ is $\kappa$-compact in $\sring$; such a cardinal is necessarily uncountable by Tougeron's flat function Lemma), so we cannot conclude that $\splring$ is compactly generated solely from the knowledge that it arises as a pullback of compactly generated presentable \infcatst. Nevertheless, we have the following result.
\begin{prop}\label{prop:plogprojectivelygenerated}
The \infcat $\splring$ is the \infcat of algebras for a 2-sorted Lawvere theory (in particular, $\splring$ is compactly generated). More precisely, consider the wide subcategory $\cartsp_c^{\rhd}\subset \cartsp$ whose morphisms are interior $b$-maps $f:\R^n\times\R^k_{\geq 0}\rightarrow \R^m\times\R^j_{\geq 0}$ that satisfy the following condition.
\begin{enumerate}
    \item[$(*)$] $f$ pulls back every  boundary defining function of $\R^m\times\R^j_{\geq 0}$ to a product of boundary defining functions on $\R^n\times\R^k_{\geq 0}$.
\end{enumerate}
We may repeat Construction \ref{cons:logcinftyrings} for $\cartsp_c^{\rhd}$, which results in a functor $\Xi^{\rhd}$. Then the functor $\Xi^{\rhd}$ induces an equivalence 
\[ \sring_{pc}^{\rhd}\overset{\simeq}{\longrightarrow} \splring. \]
\end{prop}
\begin{proof}
It follows from Proposition \ref{prop:lawvereprojgen} that the \infcat $\sring\times s\mathsf{CMon}$ is the \infcat of algebras for the 2-sorted Lawvere theory $\cartsp\times \mathsf{FCMon}^{op}$. It follows from Corollary \ref{cor:ploglimcolimcons} and \cite{HA}, Proposition 7.1.4.12 that $\splring$ is generated under sifted colimits by the essential image of the map 
\[ \cartsp^{op}\times\mathsf{FCMon} \overset{j}{\hooklongrightarrow} \sring\times s\mathsf{CMon} \overset{F}{\longrightarrow} \splring, \]
which consists of compact projective objects, where $F$ is a left adjoint to $p\times \ev_{\mathsf{PLog}}$. Let $\mathrm{T}^{op}\subset \splring$ denote this essential image which is equivalent to its full subcategory spanned by objects of the form $(\cinfty(\R^n_{\geq0}\times \R^k),\Z_{\geq 0}^n\rightarrow \cinfty_{\geq0}(\R^n_{\geq0}\times \R^k))$, then $\mathrm{T}$ is a 2-sorted Lawvere theory and the full subcategory inclusion $\mathrm{T}^{op}\subset \splring$ induces an equivalence $s\mathrm{T}\mathsf{Alg}\simeq \splring$. We are left to show that the functor $\Xi^{\rhd}$ is an equivalence. Since $\Xi^{\rhd}$ is a right adjoint that preserves sifted colimits, its left adjoint $V:\splring\rightarrow\sring_{pc}^{\rhd}$ carries $\mathrm{T}^{op}$ into the full subcategory $\icat_0\subset \sring^{\rhd}_{pc}$ spanned by compact projective objects, which contains $\cartsp^{\rhd}$. It suffices to show that the resulting functor $\mathrm{T}\rightarrow \icat_0^{op}$ factors through $\cartsp^{\rhd}$ as an equivalence. To see it is essentially surjective, note that the diagram 
\[
\begin{tikzcd}
\sring_{pc}^{\rhd}\ar[dr,"\iota^*_{c}\times\ev_{\R_{\geq0}}"']\ar[rr,"\Xi^{\rhd}"] && \splring\ar[dl,"p\times\ev_{\mathsf{PLog}}"] \\
& \sring\times s\mathsf{CMon}.
\end{tikzcd}
\]
induces a diagram 
\[
\begin{tikzcd}
\icat_0^{op} && \mathrm{T}\ar[ll] \\
&  \cartsp\times\mathsf{FCMon}^{op} \ar[ur] \ar[ul].
\end{tikzcd}
\]
so we conclude using that $\cartsp\times\mathsf{FCMon}^{op}\rightarrow \icat_0^{op}$ factors through $\cartsp^{\rhd}_c$ as an essentially surjective functor. For fully faithfulness, we note that Proposition \ref{prop:xiproperties} establishes that $\Xi^{\rhd}$ is a right adjoint relative to $\sring$, so we have a natural equivalence $\iota_c^*\circ V\simeq p$ which yields for each pair of objects $A:=(\cinfty(\R^{n}\times\R^k_{\geq 0}),\Z^k_{\geq0}\rightarrow\cinfty_{\geq0}(\R^n\times \R_{\geq0}^k))$  and $B:=(\cinfty(\R^{m}\times\R^j_{\geq 0}),\Z^j_{\geq0}\rightarrow\cinfty_{\geq0}(\R^m\times \R_{\geq0}^j))$ a commuting diagram 
\[
\begin{tikzcd}
\Hom_{\sring^{\rhd}_c}(j(\R^n\times\R_{\geq0}^k),j(\R^m\times\R^j_{\geq0})) \ar[dr] && \Hom_{\splring}(A,B) \ar[dl] \ar[ll] \\
& \Hom_{\sring}(\cinfty(\R^n\times\R_{\geq0}^k),\cinfty(\R^m\times\R^j_{\geq0}))
\end{tikzcd}
\]
of $0$-truncated spaces (the upper right space is 0-truncated by Remark \ref{rmk:0truncatedplrings}). The result will thus be established if we can argue that on connected components, both diagonal maps are injective and have the same image, which is a direct inspection. 
\end{proof}
\begin{rmk}
In view of Proposition \ref{prop:plogprojectivelygenerated}, the functor $\Xi:\sring_{pc}\rightarrow\splring$ can be identified with the functor $\sring_{pc}\rightarrow\sring_{pc}^{\rhd}$ given by the subcategory inclusion $\cartsp^{\rhd}_{c}\rightarrow\cartsp_c$.
\end{rmk}
\begin{rmk}
The reasoning applied in this subsection is valid in algebraic (derived) logarithmic geometry as well, showing that the \infcat of simplicial prelog rings over some commutative ring $k$ is also projectively generated. 
\end{rmk}
The main observation not of formal nature underlying Theorem \ref{thm:logcinftycorners} is contained in the following lemma.
\begin{lem}\label{lem:intblogification}
Let $M$ be a manifold with faces and let $H_1(M)$ be the set of connected boundary components. Consider the map $e_M:\Z_{\geq 0}^{H_1(M)}\rightarrow \cinfty_{\geq 0}(M)$ of commutative monoids induced by the map of sets $H^1(M)\rightarrow \cinfty_{\geq 0}(M)$ carrying the boundary component $S$ to a function defining $S$. Then $e_M$ takes values in the submonoid of interior $b$-maps and the commuting triangle
\[
\begin{tikzcd}
\Z_{\geq 0}^{H_1(M)}\ar[rr]\ar[dr,"e_M"'] && \cinfty_{b}(M)\ar[dl] \\
& \cinfty_{\geq 0}(M)
\end{tikzcd}
\]
exhibits $\cinfty_{b}(M)\subset \cinfty_{\geq 0}(M)$ as the logification of $e_M:\Z_{\geq 0}^{H_1(M)}\rightarrow\cinfty_{\geq 0}(M)$.
\end{lem}
\begin{proof}
Clearly, functions defining boundary components on $M$ are interior $b$-maps. Consider the diagram 
\[
\begin{tikzcd}
0\ar[d] \ar[r] & \cinfty_{> 0}(M) \ar[d] \ar[r,equal] & \cinfty_{> 0}(M) \ar[d] \\
\Z^{H_1(M)}_{\geq0} \ar[r] & \cinfty_{b}(M)\ar[r] &\cinfty_{\geq 0}(M)
\end{tikzcd}
\]
of commutative monoids. It is easy to see that both squares are pullbacks, so it suffices to show that the left square is also a pushout of simplicial commutative monoids; that is, the maps $Z_{\geq0}^{H_1(M)}\rightarrow \cinfty_{b}(M)$ and $ \cinfty_{> 0}(M)\rightarrow \cinfty_{b}(M)$ exhibit $\cinfty_{b}(M)$ as a coproduct of $\Z_{\geq 0}^{H_1(M)}$ and $ \cinfty_{> 0}(M)$. The symmetric monoidal structure on $s\mathsf{CMon}$ is coCartesian and the symmetric monoidal structure on $\spa$ is Cartesian, so after unwinding definitions, we are reduced to producing an equivalence of spaces $\cinfty_{b}(M)\simeq \Z_{\geq 0}^{H_1(M)}\times \cinfty_{> 0}(M)$ (which is just a bijection of sets in this case) such that the induced maps $\cinfty_{> 0}(M)\rightarrow \cinfty_{> 0}(M)$ and $\Z_{\geq 0}^{H_1(M)}\rightarrow\Z_{\geq 0}^{H_1(M)}$ are equivalent to the identity, and the maps $\cinfty_{> 0}(M)\rightarrow\Z_{\geq 0}^{H_1(M)}$ and $\Z_{\geq 0}^{H_1(M)}\rightarrow \cinfty_{> 0}(M)$ are equivalent to the zero morphism. We get the desired bijection of sets $\cinfty_{b}(M)\cong \Z_{\geq 0}^{H_1(M)}\times \cinfty_{> 0}(M)$ from the observation that every interior $b$-map $f:M\rightarrow \R_{\geq 0}$ can be written as $h_{S_1}^{m_1}\ldots h_{S_n}^{m_n}g$ with a unique $g\in \cinfty_{>0}(M)$ and a unique tuple $(h_S)_{{H_1(M)}}\in\Z_{\geq 0}^{H_1(M)}$, the indicated coefficients associated to the $\{S_j\}$ being the only ones that are nonzero.
\end{proof}
\begin{cor}\label{cor:chilcolimit}
The composition $\sring_{pc}\overset{\Xi}{\rightarrow}\splring\overset{L_{\mathsf{Log}}}{\rightarrow}\slring$ preserves all colimits.
\end{cor}
\begin{proof}
As $L_{\mathsf{Log}}$ preserves colimits and $\Xi$ preserves sifted colimits, the composition $L_{\mathsf{Log}}\Xi$ is a left Kan extension of its restriction to the essential image of the Yoneda embedding $j:\cartsp_c^{op}\hookrightarrow\sring_{pc}$, so it suffices to show that the composition $L_{\mathsf{Log}}\Xi j$ preserves coproducts. Contemplate the commuting diagrams 
\[
\begin{tikzcd}
(\cinfty(\R^n_{\geq 0}\times\R^k),\Z^{n}_{\geq 0}\rightarrow \cinfty_{\geq 0}(\R^n_{\geq 0}\times\R^k)) \ar[d] \ar[r,"\alpha"] & (\cinfty(\R^{n+m}_{\geq 0}\times\R^{k+l}),\Z^{n+m}_{\geq 0}\rightarrow \cinfty_{\geq 0}(\R^{n+m}_{\geq 0}\times\R^{k+l})) \ar[d,"\gamma"] \\
(\cinfty(\R^n_{\geq 0}\times\R^k),\cinfty_{b}(\R^n_{\geq 0}\times\R^k)\rightarrow \cinfty_{\geq 0}(\R^n_{\geq 0}\times\R^k)) \ar[r,"\alpha'"] & (\cinfty(\R^{n+m}_{\geq 0}\times\R^{k+l}),\cinfty_{b}(\R^{n+m}_{\geq 0}\times\R^{k+l})\rightarrow \cinfty_{\geq 0}(\R^{n+m}_{\geq 0}\times\R^{k+l})) 
\end{tikzcd}
\]
and
\[
\begin{tikzcd}
(\cinfty(\R^{n+m}_{\geq 0}\times\R^{k+l}),\Z^{n+m}_{\geq 0}\rightarrow \cinfty_{\geq 0}(\R^{n+m}_{\geq 0}\times\R^{k+l})) \ar[d,"\gamma"] & (\cinfty(\R^m_{\geq 0}\times\R^l),\Z^{m}_{\geq 0}\rightarrow \cinfty_{\geq 0}(\R^m_{\geq 0}\times\R^l)) \ar[l,"\beta"'] \ar[d] \\
(\cinfty(\R^{n+m}_{\geq 0}\times\R^{k+l}),\cinfty_{b}(\R^{n+m}_{\geq 0}\times\R^{k+l})\rightarrow \cinfty_{\geq 0}(\R^{n+m}_{\geq 0}\times\R^{k+l})) & (\cinfty(\R^m_{\geq 0}\times\R^l),\cinfty_{b}(\R^m_{\geq 0}\times\R^l)\rightarrow \cinfty_{\geq 0}(\R^m_{\geq 0}\times\R^l)). \ar[l,"\beta'"]
\end{tikzcd}
\]
We wish to show that the maps $\alpha'$ and $\beta'$ exhibit a coproduct in $\slring$. Using Corollary \ref{cor:ploglimcolimcons}, we deduce that the maps $\alpha$ and $\beta$ exhibit a coproduct in $\splring$ as the underlying diagram of $\cinfty$-rings and the underlying diagram of (finitely generated free) simplicial commutative monoids exhibit a coproduct. By virtue of Lemma \ref{lem:intblogification}, the vertical maps in the diagrams above exhibit logifications, so we conclude by observing that logification, as a left adjoint, preserves coproducts.
\end{proof}
\begin{lem}\label{lem:sectionlog}
The composition
\[ \sring \overset{\iota_{c!}}{\longrightarrow} \sring_{pc} \overset{\Xi}{\longrightarrow} \splring \overset{L_{\mathsf{Log}}}{\longrightarrow}\slring \]
is equivalent to the composition
\[  \sring \overset{s}{\longrightarrow} \splring \overset{L_{\mathsf{Log}}}{\longrightarrow} \slring, \]
where $s$ is a left adjoint to $p:\splring\rightarrow\sring$.
\end{lem}
\begin{proof}
Consider the full subcategory $\icat\subset\fun_{\sring}(\sring,s\cinfty\mathsf{Log})$ spanned by sections $F$ satisfying the following conditions.
\begin{enumerate}[$(1)$]
    \item $F$ preserves sifted colimits.
    \item For each $n\geq 0$, $F$ carries the object $\cinfty(\R^n)$ to an initial object in the fibre over $\cinfty(\R^n)$.
\end{enumerate}
Sections satisfying $(1)$ are precisely left Kan extensions of their restriction along the full subcategory inclusion $\cartsp^{op}\subset \sring$ so this restriction induces an equivalence between $\icat$ and the full subcategory of 
\[\fun_{\cartsp^{op}}(\cartsp^{op},\cartsp^{op}\times_{\sring}s\cinfty\mathsf{Log})\]
spanned by sections $f$ that carry each object of $\cartsp^{op}$ to an initial object in the fibre. The projection $q:\cartsp^{op}\times_{\sring}s\cinfty\mathsf{Log}\rightarrow \cartsp^{op}$ is a Cartesian fibration, so each such functor is a left adjoint to $q$. It follows that the set of equivalence classes of objects of $\icat$ consists of a single element. We conclude by observing that both functors in the statement of the lemma satisfy $(1)$ and $(2)$.
\end{proof}
\begin{lem}\label{lem:xilogpreservescolimits}
The functor $\Xi_{\mathsf{Log}}:\sring_c\rightarrow\slring$ preserves all colimits.
\end{lem}
\begin{proof}
It suffices to argue that $L_{\mathsf{Log}}\Xi$ carries the set $S=\{\phi\}$ of Definition \ref{defn:cinftycorners} into the set of equivalences of $\slring$, as it then follows from the universal property of cocontinuous localizations that the functor $L_{\mathsf{Log}}\Xi$ factors through $\sring_c$ as a colimit preserving functor. Since $\Xi$ restricted to $\sring_c$ takes values in $\slring$, we consequently deduce that $\Xi_{\mathsf{Log}}$ is equivalent to $L_{\mathsf{Log}}\Xi$ and therefore preserves colimits. \\
The functor $L_{\mathsf{Log}}\Xi:\scring_{pc}\rightarrow \slring$ preserves colimits by Corollary \ref{cor:chilcolimit} so it carries the pushout diagram
\[
\begin{tikzcd}
\iota_{c!}\iota_c^* (\cinfty(\R_{\geq 0}),\cinfty_{b}(\R_{\geq 0})) \ar[d] \ar[r,"\epsilon"] & (\cinfty(\R_{\geq 0}),\cinfty_{b}(\R_{\geq 0}))  \ar[d]\\ \iota_{c!}\iota_c^* (\cinfty(\R_{> 0}),\cinfty_{b}(\R_{> 0})) 
\ar[r]& \mathcal{A}  
\end{tikzcd}
\]
of Definition \ref{defn:cinftycorners} to a pushout diagram in $\slring$. It follows from Lemma \ref{lem:sectionlog} that the pushout diagram above is carried to a pushout diagram
\[
\begin{tikzcd}
(\cinfty(\R_{\geq 0}),\cinfty_{\geq0}(\R_{\geq 0})^{\times}\rightarrow \cinfty_{\geq0}(\R_{\geq 0})) \ar[d] \ar[r] & (\cinfty(\R_{\geq 0}),\cinfty_{b}(\R_{\geq 0})\rightarrow \cinfty_{\geq0}(\R_{\geq 0}))  \ar[d]\\ (\cinfty(\R_{> 0}),\cinfty_{\geq0}(\R_{> 0})^{\times}\rightarrow \cinfty_{\geq0}(\R_{> 0})) 
\ar[r]& L_{\mathsf{Log}}\Xi(\mathcal{A}),
\end{tikzcd}
\]
where the left vertical map is a coCartesian morphism between initial log structures. Since the functor $\slring\rightarrow \sring$ preserves colimits, the map on underlying simplicial $\cinfty$-rings of the lower horizontal map in the diagram above is an equivalence. Since the left vertical map is a $p_{\mathsf{Log}}$-coCartesian edge and the diagram is a $p_{\mathsf{Log}}$-pushout, the right vertical map is also $p_{\mathsf{Log}}$-coCartesian. Therefore, we are reduced to verifying that the logification of $(\cinfty(\R_{> 0}),\cinfty_{b}(\R_{\geq 0})\rightarrow \cinfty_{\geq0}(\R_{> 0}))$ is the initial log structure. Consider the pullback diagram 
\[
\begin{tikzcd}
M \ar[r] \ar[d]& \cinfty_{> 0}(\R_{> 0}) \ar[d]\\
\cinfty_{b}(\R_{\geq 0}) \ar[r] & \cinfty_{\geq 0}(\R_{> 0}).
\end{tikzcd}
\]
Recalling the description of the logification functor, we wish to show that the map $\cinfty_{> 0}(\R_{> 0})\rightarrow \cinfty_{>0}(\R_{>0})\coprod_{M}\cinfty_{b}(\R_{\geq 0})$ is an equivalence. It is sufficient to argue that the left vertical map in the diagram above is an equivalence, which is equivalent to the assertion that if $f:\R_{\geq0}\rightarrow\R_{\geq 0}$ is an interior $b$-map, then the restriction $f|_{\R_{>0}}$ factors through $\R_{>0}\hookrightarrow\R_{\geq 0}$, but this holds by definition of interior $b$-maps.
\end{proof}
\begin{proof}[Proof of Theorem \ref{thm:logcinftycorners}]
Let $F$ denote a left adjoint to $\Xi_{\mathsf{Log}}$. Since $\Xi_{\mathsf{Log}}$ is conservative, it suffices to argue that the unit transformation $\mathrm{id}\rightarrow\Xi_{\mathsf{Log}}F$ is an equivalence. Since both $\Xi_{\mathsf{Log}}$ and $F$ preserve colimits and the objects $L_{\mathsf{Log}}(\cinfty(\R^n_{\geq 0}\times\R^k),\Z^n_{\geq0}\rightarrow \cinfty_{\geq 0}(\R^n_{\geq 0}\times\R^k))$ generate $\slring$ under sifted colimits, we need only check that the unit is an equivalence on this collection of objects. It follows from (the proof of) Lemma \ref{lem:xilogpreservescolimits} that $\Xi$ carries the strong saturation $\overline{S}$ of the set $S=\{\phi\}$ to the set of maps in $\splring$ that become an equivalence after applying $L_{\mathsf{Log}}$. In particular, for any localization $X\rightarrow L(X)$ in $\sring_{pc}$, the map $\Xi (X)\rightarrow \Xi L(X)\simeq \Xi_{\mathsf{Log}}L(X)$ in $\splring$, whose codomain lies in $\slring$, becomes an equivalence upon logifying and is therefore also a localization, that is, the diagram 
\[
\begin{tikzcd}
\sring_c\ar[d,hook] \ar[r,"\Xi_{\mathsf{Log}}"] & \slring \ar[d,hook] \\
\sring_{pc}\ar[r,"\Xi"] & \splring
\end{tikzcd}
\]
is vertically left adjointable. Then the resulting commuting diagram 
\[
\begin{tikzcd}
\sring_c \ar[r,"\Xi_{\mathsf{Log}}"] & \slring  \\
\sring_{pc}\ar[r,"\Xi"] \ar[u,"L"]& \splring \ar[u,"L_{\mathsf{Log}}"]
\end{tikzcd}
\]
is tautologically vertically right adjointable, and therefore also horizontally left adjointable, that is, the logification functor carries unit transformations of the lower adjunction to unit transformations of the upper one. It follows from Proposition \ref{prop:plogprojectivelygenerated} that the object $(\cinfty(\R_{\geq 0}^{n}\times\R^k),\cinfty_{b}(\R_{\geq 0}^{n}\times\R^k))$ together with the triangle
\[
\begin{tikzcd}
\Z_{\geq 0}^n\ar[rr]\ar[dr,"e_n"'] && \cinfty_{b}(\R^n_{\geq0}\times\R^k) \ar[dl] \\
& \cinfty_{\geq 0}(\R^n_{\geq0}\times\R^k)
\end{tikzcd}
\]
is a unit transformation at $(\cinfty(\R_{\geq 0}^{n}\times\R^k),\Z_{\geq 0}^n\rightarrow\cinfty_{\geq 0}(\R^n\times\R^k))$. This map exhibits a logification by Lemma \ref{lem:intblogification} and is therefore carried to an equivalence by $L_{\mathsf{Log}}$.
\end{proof}
We now turn to the proof of Theorem \ref{thm:geoenvcorners}. \begin{lem}\label{lem:localizealgtheory}
Let $\mathrm{T}$ be a Lawvere theory and let $s\mathrm{T}\alg$ be the associated \infcat of algebras. Let $S$ be small set of morphisms in $s\mathrm{T}\alg$ and denote by $s\mathrm{T}\alg[S^{-1}]\subset s\mathrm{T}\alg$ the strongly reflective full subcategory spanned by $S$-local objects. Let $\icat$ be an idempotent complete \infcat that admits finite limits and denote by $\fun^{\pi}(\mathrm{T},\icat)[S^{-1}]\subset\fun^{\pi}(\mathrm{T},\icat)$ the full subcategory spanned by functors $F:\mathrm{T}\rightarrow\icat$ for which the following condition is satisfied.
\begin{enumerate}[$(*)$]
    \item For each object $C\in \icat$, the composition
    \[ \mathrm{T}\overset{F}{\longrightarrow} \icat\overset{\Hom_{\icat}(C,\_)}{\longrightarrow}\spa   \]
    is $S$-local in $s\mathrm{T}\alg$.
\end{enumerate}
Suppose that the inclusion $s\mathrm{T}\alg[S^{-1}]\subset s\mathrm{T}\alg$ preserves filtered colimits, then restriction along the functor $\mathrm{T}^{op}\overset{j}{\hookrightarrow} s\mathrm{T}\alg\overset{L}{\rightarrow} s\mathrm{T}\alg[S^{-1}]$ induces an equivalence
\[ \fun^{\mathrm{lex}}(s\mathrm{T}\alg[S^{-1}]^{op}_{\fp},\icat) \overset{\simeq}{\longrightarrow}  \fun^{\pi}(\mathrm{T},\icat)[S^{-1}]\]
with inverse given by a functor taking right Kan extensions along $\mathrm{T}\hookrightarrow s\mathrm{T}\alg_{\fp}\rightarrow s\mathrm{T}\alg[S^{-1}]_{\fp}$.
\end{lem}
\begin{proof}
The Yoneda embedding $j:\icat\hookrightarrow\pshv(\icat)$ induces a commuting diagram 
\[
\begin{tikzcd}
 \fun'(s\mathrm{T}\alg[S^{-1}]^{op},\icat) \ar[d,hook]\ar[r] &\fun^{\mathrm{lex}}(s\mathrm{T}\alg[S^{-1}]^{op}_{\fp},\icat) \ar[d,hook] \ar[r] & \fun^{\pi}(\mathrm{T},\icat)  \ar[d,hook] \\
 \fun'(s\mathrm{T}\alg[S^{-1}]^{op},\pshv(\icat)) \ar[r]&
\fun^{\mathrm{lex}}(s\mathrm{T}\alg[S^{-1}]^{op}_{\fp},\pshv(\icat)) \ar[r,"(Lj)^*"] & \fun^{\pi}(\mathrm{T},\pshv(\icat)) 
\end{tikzcd}
\]
where $\fun'(s\mathrm{T}\alg[S^{-1}]^{op},\icat)$ and $\fun'(s\mathrm{T}\alg[S^{-1}]^{op},\pshv(\icat))$ denote full subcategories of functors preserving small limits. As $\pshv(\icat)$ admits small limits and the \infcat $s\mathrm{T}\alg[S^{-1}]$ is compactly generated by virtue of the assumption that the inclusion $s\mathrm{T}\alg[S^{-1}]\subset s\mathrm{T}\alg$ preserves filtered colimits, the lower left horizontal restriction 
\[ \fun'(s\mathrm{T}\alg[S^{-1}]^{op},\pshv(\icat)) \overset{\simeq}{\longrightarrow} \fun^{\mathrm{lex}}(s\mathrm{T}\alg[S^{-1}]^{op}_{\fp},\pshv(\icat))\]
is an equivalence after Remark \ref{rmk:lexfunctors}. The composition $\fun'(s\mathrm{T}\alg[S^{-1}]^{op},\pshv(\icat))\rightarrow \fun^{\pi}(\mathrm{T},\pshv(\icat))$ factors via the restriction \[r:\fun_S'(s\mathrm{T}\alg,\pshv(\icat))\longrightarrow \fun^{\pi}(\mathrm{T},\pshv(\icat)),\]
where $\fun_S'(s\mathrm{T}\alg,\pshv(\icat))$ is the full subcategory spanned by limit preserving functors $F:s\mathrm{T}\alg^{op}\rightarrow \pshv(\icat)$ carrying the set $S$ to into the set of equivalences of $\pshv(\icat)$. This is the case for such a functor $F$ if and only if for each $C\in \icat$, the functor $\ev_{C}\circ F:s\mathrm{T}\alg^{op}\rightarrow\spa$ carries the set $S$ into the set of equivalences in $\spa$, but since $\ev_{C}\circ F$ preserves limits and is therefore representable, this corresponds to the associated representing object $A\in s\mathrm{T}\alg$ being $S$-local. Let $\fun^{\pi}(\mathrm{T},\pshv(\icat))[S^{-1}]$ be the full subcategory spanned by limit preserving functors $F$ such that $\ev_{C}\circ F:s\mathrm{T}\alg^{op}\rightarrow\spa$ carries the set $S$ into the set of equivalences in $\spa$. Since the representing object $A$ of $\ev_{C}\circ F$ may be identified with the functor $\ev_{C}\circ F \circ j$, we conclude that the restriction $r$ takes values in $\fun^{\pi}(\mathrm{T},\pshv(\icat))[S^{-1}]$ and determines an equivalence $\fun_S'(s\mathrm{T}\alg,\pshv(\icat))\simeq \fun^{\pi}(\mathrm{T},\pshv(\icat))[S^{-1}]$. It follows that the restriction $(Lj)^*$ factors up to homotopy through $\fun^{\pi}(\mathrm{T},\pshv(\icat))[S^{-1}]$ and is an equivalence onto its essential image, but as $\fun^{\pi}(\mathrm{T},\pshv(\icat))[S^{-1}]\subset \fun^{\pi}(\mathrm{T},\pshv(\icat))$ is a replete full subcategory, $(Lj)^*$ itself factors through $\fun^{\pi}(\mathrm{T},\pshv(\icat))[S^{-1}]$ and determines an equivalence $\fun^{\mathrm{lex}}(s\mathrm{T}\alg[S^{-1}]^{op}_{\fp},\pshv(\icat))\simeq \fun^{\pi}(\mathrm{T},\pshv(\icat))[S^{-1}]$. Since we have an isomorphisms of simplicial sets $\fun^{\pi}(\mathrm{T},\pshv(\icat))[S^{-1}]\times_{\fun^{\pi}(\mathrm{T},\pshv(\icat))}\fun^{\pi}(\mathrm{T},\icat)\cong \fun^{\pi}(\mathrm{T},\icat)[S^{-1}]$, we deduce that restriction along $jL$ induces the top horizontal map in the commuting diagram 
\[
\begin{tikzcd}
\fun^{\mathrm{lex}}(s\mathrm{T}\alg[S^{-1}]^{op}_{\fp},\icat) \ar[d,hook] \ar[r] & \fun^{\pi}(\mathrm{T},\icat)[S^{-1}] \ar[d,hook] \\
\fun^{\mathrm{lex}}(s\mathrm{T}\alg[S^{-1}]^{op}_{\fp},\pshv(\icat)) \ar[r,"\simeq"] & \fun^{\pi}(\mathrm{T},\pshv(\icat))[S^{-1}].
\end{tikzcd}
\]
By assumption on $\icat$, the essential image of the Yoneda embedding is stable under finite limits and retracts in $\pshv(\icat)$, so using that every object of $s\mathrm{T}\alg[S^{-1}]^{op}_{\fp}$ is a retract of a finite limit of objects in the essential image of $\mathrm{T}\rightarrow s\mathrm{T}\alg\overset{L}{\rightarrow} s\mathrm{T}\alg[S^{-1}]$ we conclude that the top horizontal map is an equivalence. 
\end{proof}
\begin{proof}[Proof of Theorem \ref{thm:geoenvcorners} $(i)$, $(ii)$, $(iii)$]
We verify the claims made in the statement of the theorem.
\begin{enumerate}[$(i)$]
    \item The \infcat $\sring_c$ is compactly generated. This was checked in Corollary \ref{cor:accessiblecorner}.
    \item Definition \ref{defn:cornergeometry} determines the structure of a geometry on $\geodiffderc$. We need to check that admissible morphisms are stable under pullbacks, retracts and that, if $g$ is admissible and $h$ another map with codomain the domain of $g$, then $h$ is admissible if and only if $g\circ h$ is admissible. Since localizations are stable under pushouts of simplicial $\cinfty$-rings and the functor $\iota_c^*$ preserves colimits, it suffices to show that a pushout in $\sring_c$ along a $\iota_c^*$-coCartesian morphism is again $\iota_c^*$-coCartesian. This is the case since all colimits in $\sring_c$ are $\iota_c^*$-colimits. Similarly, we know that localizations of morphisms of simplicial $\cinfty$-rings are stable under retracts, so we conclude that admissible morphisms in $\geodiffderc$ are stable under retracts from the observation that coCartesian morphisms are (which in turn follows from the fact that pullback squares are stable under retracts). Repeating this line of argument once more, we obtain the last verification from the corresponding verification for localizations, together with \cite{HTT}, Proposition 2.4.1.7.
    \item The inclusion $\diffc'\rightarrow \geodiffderc$ exhibits a geometric envelope. Choose an idempotent complete \infcat admitting finite limits, then we have a commuting diagram 
    \[
    \begin{tikzcd}
    & \fun^{\pi\mathrm{ad}}(\diffc',\icat) \ar[dr,"\theta''"] \\ \fun^{\mathrm{lex}}(\geodiffderc,\icat) \ar[ur,"\theta"] \ar[rr,"\theta'"] && \fun^{\pi}(\cartsp_c,\icat)[S^{-1}].
    \end{tikzcd}
    \]
    Let us first argue that the restriction functor $\theta''$ indeed takes values in $\fun^{\pi}(\cartsp_c,\icat)[S^{-1}]$: composing with the functor $\Hom_{\icat}(C,\_):\icat\rightarrow \spa$ for some $C\in \icat$, we may replace $\icat$ by $\spa$ and $\fun^{\pi}(\cartsp_c,\icat)[S^{-1}]$ by $\sring_c$. We note that $\fun^{\pi\mathrm{ad}}(\diffc',\spa)$ is an $\omega$-accessible localization of $\pshv(\diffc'^{op})$ and that restriction along $\cartsp_c\rightarrow \diffc'$ induces the functor $\theta'':\fun^{\pi\mathrm{ad}}(\diffc',\spa)\rightarrow\sring_{pc}$ which preserves limits and filtered colimits. To conclude that $\theta''$ factors through $\sring_c$, we need to show that its left adjoint $F$ carries the set $S'$ of Remark \ref{rmk:altset} into the set of equivalences of $\fun^{\pi\mathrm{ad}}(\diffc',\spa)$. We have a commuting diagram 
    \[ \begin{tikzcd}
    \cartsp^{op} \ar[d,"j"] \ar[r] & \diffc^{\prime op}\ar[d,"j"] \\
    \pshv(\cartsp^{op}) \ar[r] \ar[d]& \pshv(\diffc^{\prime op}) \ar[d] \\
    \sring_{pc}\ar[r,"F"] & \fun^{\pi\mathrm{ad}}(\diffc',\spa)
    \end{tikzcd} \]
    where the lower square is obtained by passing to left adjoints in the square 
    \[
    \begin{tikzcd}
    \pshv(\cartsp^{op})& \pshv(\diffc^{\prime op}) \ar[l]\\
    \sring_{pc} \ar[u,hook] & \fun^{\pi\mathrm{ad}}(\diffc',\spa) \ar[u,hook] \ar[l]
    \end{tikzcd}
    \]
    where the horizontal functors are induced by pulling back along $\cartsp\rightarrow\diffc'$. It follows that $F$ carries the map in $S'$ to the lower horizontal one in the pushout diagram 
    \[
    \begin{tikzcd}
    j(\R) \ar[d,"j(\exp)"'] \ar[r] & j(\R_{\geq 0}) \ar[d] \\
    j(\R) \ar[r,"\phi"] & \mathcal{A},
    \end{tikzcd}
    \]
    but the Yoneda embedding $j:\diffc'^{op}\rightarrow\fun^{\pi\mathrm{ad}}(\diffc',\spa)$ preserves pushouts along admissible maps, so $\phi$ is indeed an equivalence. It follows from Lemma \ref{lem:localizealgtheory} that the functor $\theta'$ is an equivalence, so it suffices to show that the functor $\theta''$ is an equivalence. The functor taking right Kan extensions along the inclusion $\cartsp_c\hookrightarrow\diffc'$ is a right adjoint to $\theta''$ which can be identified with the inverse of $\theta'$ composed with $\theta$. The counit map of this adjunction is an equivalence, and it is obvious from the definition of $\diffc'$ that the unit map is also an equivalence.   
\end{enumerate}
\end{proof}
We will complete the proof of Theorem \ref{thm:geoenvcorners} at the end of this section; first, we remark on the discrepancy between $\diffc$ and $\diffc'$. The pregeometry $\diffc'$ is \emph{not} equivalent to $\diffc$, nor does the functor $(\cinfty(\_),\cinfty_b(\_)):\diffc\rightarrow\sring_c$ take values in $(\geodiffderc)^{op}$. Indeed, we have the following alternative.
\begin{lem}
Let $M$ be a manifold with faces, then $(\cinfty(M),\cinfty_b(M))$ is a compact object in $\sring_c$ if and only if $M$ has finitely many connected boundary components.
\end{lem}
\begin{proof}[Proof sketch]
First, consider $M$ a manifold with faces with infinitely many boundary components. Lemma \ref{lem:intblogification} shows that there is an equivalence $\cinfty_b(M)\cong \Z^{M_1(M)}\coprod \cinfty_{\geq0}(M)^{\times}$, where $M_1(M)$ is the set of connected boundary components of $M$. It follows that the sharpening of $\cinfty_b(M)$ is infinitely generated so, as the simplicial commutative monoid associated to any finitely generated simplicial $\cinfty$-rings with corners has finitely generated sharpening, the object $(\cinfty(M),\cinfty_b(M))$ cannot be compact in $\sring_c$. The converse follows from the following assertions.
\begin{enumerate}
    \item[$(*)$] For every manifold with faces $M$, there exists an interior $b$-map $M\hookrightarrow \R^{n}\times\R^k_{\geq 0}$ which is a \emph{p-embedding} of manifolds with corners (see \cite{Melrose1} for an explanation of the terminology). 
    \item[$(**)$] Let $S\subset M$ be a $p$-embedded submanifold, then $S$ admits a tubular neighbourhood.
\end{enumerate}
To prove $(*)$, we use the boundary flowout map $M\hookrightarrow M^{\circ}$ to embed $M$ into its interior, and then apply the Whitney embedding theorem to embed $M^{\circ}$ into $\R^n$ for some $n>>1$ resulting in a closed embedding $f:M\hookrightarrow\R^n$. Choose a finite complete set of boundary defining functions $\{\rho_H\}_{H\in M_1(M)}$, then the map $f\prod_{H\in M_1(M)}\rho_H:M\rightarrow \R^n\times\R^k_{\geq 0}$ with $k=|M_1(M)|$ is a embedding. The fact that every $p$-embedded submanifold admits a tubular neighbourhood is proven verbatim as in the case without corners.
\end{proof}
\begin{rmk}
A similar argument as the one presented in the previous lemma yields that every object $(\cinfty(U),M\rightarrow\cinfty_{\geq0}(U))$ in $\diffc'$ must have finitely generated sharpening, but it is certainly possible that as an open subset $U\subset\R_{\geq 0}^n\times\R^{k}$ has infinitely many boundary components, so admissible morphisms in $\geodiffderc$ may not `create' sufficiently many boundary defining functions. Both these issues disappear when we apply the spectrum functor $\spec_c$, since every manifold with faces may always be covered by opens that admit an embedding $U\rightarrow\R_{\geq0}^n\times\R^{k}$ onto an open subset with \emph{connected} boundary. 
\end{rmk}
\subsection{Sheaves of log structures}
Let $s_c:\sring\rightarrow \slring$ be a left adjoint to the projection $p_{\mathsf{Log}}:\slring\rightarrow \sring$. This functor carries compact objects to compact objects and determines a transformation of geometries $s_c:\geodiffder\rightarrow\geodiffderc$. This transformation is coCartesian for the functor $\mathsf{Forget}:\mathsf{Geo}\rightarrow\catinf^{\mathrm{lex},\mathrm{Idem}}$, that is $\geodiffderc$ is endowed with the coarsest admissibility structure and compatible topology that makes $s_c$ a transformation of geometries. Let $\xtop$ be an \inftopt, then composition with $s_c$ determines a functor $\fun^{\lex}(\geodiffderc,\xtop)\rightarrow\fun^{\lex}(\geodiffder,\xtop)$ which is equivalent to the functor $p^{\xtop}_{\mathsf{Log}}:\shv_{\sring_c}(\xtop)\rightarrow\shv_{\sring}(\xtop)$ that composes with $p_{\mathsf{Log}}:\sring_c\rightarrow\sring$. Since the functor $(\_)_{\geq 0}$ admits a left adjoint, the \infcat $\shv_{\splring}$ fits into a pullback diagram 
\[
\begin{tikzcd}
\shv_{\splring}\ar[r]\ar[d,"p^{\xtop}_{\mathsf{PLog}}"] & \fun(\Delta^1,\shv_{s\mathsf{CMon}}(\xtop))\ar[d,"\ev_0"]\\
\shv_{\sring}(\xtop) \ar[r] &\shv_{s\mathsf{CMon}}(\xtop)
\end{tikzcd}
\]
so the functor $p^{\xtop}_{\mathsf{PLog}}:\shv_{\splring}(\xtop)\rightarrow\shv_{\sring}(\xtop)$ given by composing with $p_{\mathsf{PLog}}$ is a presentable fibration. We will denote objects in $\shv_{\splring}(\xtop)$ by pairs $(\Of,\mathcal{M}\rightarrow(\Of_{\xtop})_{\geq 0})$, and we say that $\mathcal{M}$ is a \emph{sheaf of prelog structures over $(\Of_{\xtop})_{\geq 0}$ on $\xtop$} (or just a \emph{prelog structure over $(\Of_{\xtop})_{\geq 0}$} if it is clear over which \inftop we work). Here we let $(\_)_{\geq 0}$ denote the functor $\shv_{\sring}(\xtop)\rightarrow \shv_{s\mathsf{CMon}}(\xtop)$ given by composing with $(\_)_{\geq 0}$. Since the $\infty$-group of units functor $(\_)^{\times}:s\mathsf{CMon}\rightarrow s\mathsf{Ab}$ also a admits a left adjoint, the \infcat $\shv_{\slring}(\xtop)$ is equivalent to the full subcategory $\shv_{\slring}(\xtop)\subset\shv_{\xtop}(\splring)$ spanned by the pairs $(\Of_{\xtop},\mathcal{M}\rightarrow (\Of_{\xtop})_{\geq 0})$ for which the upper horizontal map in the pullback diagram 
\[
\begin{tikzcd}
\mathcal{M}\times_{(\Of_{\xtop})_{\geq 0}}(\Of_{\xtop})^{\times}_{\geq 0} \ar[d] \ar[r]&(\Of_{\xtop})^{\times}_{\geq 0}\ar[d] \\
\mathcal{M}\ar[r] & (\Of_{\xtop})_{\geq 0}
\end{tikzcd}
\]
is an equivalence. Here, the right vertical map in the diagram above is the counit of the right adjoint to the inclusion $\shv_{s\mathsf{Ab}}(\xtop)\rightarrow \shv_{s\mathsf{CMon}}(\xtop)$, which is given by composing with the right adjoint $(\_)^{\times}$. Thus, we will say that the fibre of $p^{\xtop}_{\mathsf{Log}}$ over $\Of_{\xtop}$ is the \emph{\infcat of log structures over $(\Of_{\xtop})_{\geq 0}$} and denote it $\mathsf{Log}((\Of_{\xtop})_{\geq 0})$. Since pullbacks of log structures are log structures, the functor $p^{\xtop}_{\mathsf{Log}}:\shv_{\slring}(\xtop)\rightarrow \shv_{\sring}(\xtop)$ is a Cartesian fibration.
\begin{prop}\label{prop:sheaveslogification}
Let $\xtop$ be an \inftopt, then the following hold true.
\begin{enumerate}[$(1)$]
\item The full subcategory inclusion $\shv_{\slring}(\xtop)\subset\shv_{\splring}(\xtop)$ admits a left adjoint relative to $\shv_{\slring}(\xtop)$ so that for each $\Of_{\xtop}\in\shv_{\sring}(\xtop)$, the inclusion $\mathsf{Log}((\Of_{\xtop})_{\geq 0})\subset \shv_{s\mathsf{CMon}}(\xtop)_{/(\Of_{\xtop})_{\geq 0}}$ admits a localization functorial in $\Of_{\xtop}$. 
\item Let $\Of_{\xtop}$ be a sheaf of simplicial $\cinfty$-rings on $\xtop$, then a map $f:\mathcal{M}\rightarrow \mathcal{N}$ of prelog structures over $(\Of_{\xtop})_{\geq 0}$ exhibits a localization if and only if $\mathcal{N}$ is a log structure and $f$ fits into a pushout diagram
\[
\begin{tikzcd}
\mathcal{M}\times_{(\Of_{\xtop})_{\geq 0}}(\Of_{\xtop})_{\geq 0}^{\times}\ar[d]\ar[r] & (\Of_{\xtop})_{\geq 0}^{\times}\ar[d] \\
\mathcal{M} \ar[r,"f"] & \mathcal{N}
\end{tikzcd}
\] 
of sheaves of simplicial commutative monoids on $\xtop$.
\item The functor $p^{\xtop}_{\mathsf{Log}}$ is a presentable fibration. 
\end{enumerate}
\end{prop}
\begin{proof}
The functor $\shv_{\slring}(\xtop)\subset \shv_{\splring}(\xtop)$ admits a left adjoint $L^{\xtop}_{\mathsf{Log}}$ given by the composition 
\[\shv_{\splring}(\xtop)\subset\fun(\xtop^{op},\splring)\overset{L_{\mathsf{Log}}\circ\_}{\longrightarrow}\fun(\xtop^{op},\slring)\overset{L_{\slring}}{\longrightarrow} \shv_{\slring}(\xtop)   \] 
where $L_{\slring}$ is the sheafification functor for $\slring$-valued sheaves. It follows that a map $f:X\rightarrow Y$ in $\shv_{\splring}(\xtop)$ ehxibits a localization for this adjunction if and only if $f$ factorizes as $X\overset{g}{\rightarrow} Y\overset{h}{\rightarrow} Z$ where $g$ is exhibits a localization for the functor $L_{\mathsf{Log}}\circ \_$ and $h$ exhibits a sheafification in the \infcat $\fun(\xtop^{op},\slring)$. Since $L_{\mathsf{Log}}$ is a left adjoint relative to $\sring$, the functor $p^{\xtop}_{\mathsf{Log}}$ carries $g$ to an equivalence in $\fun(\xtop^{op},\sring)$. Since the underlying presheaf of simplicial $\cinfty$-rings of the object $Y$ is already a sheaf and $h$ exhibits a sheafification, the verification that $L^{\xtop}_{\mathsf{Log}}$ is a left adjoint relative to $\shv_{\sring}(\xtop)$ will be complete once we show that the functor $\fun(\xtop^{op},\slring)\rightarrow\fun(\xtop^{op},\sring)$ carries sheafifications to sheafifications, that is, the diagram 
\[
\begin{tikzcd}
\shv_{\slring}\ar[d]\ar[r] &\fun(\xtop^{op},\slring) \ar[d] \\
\shv_{\sring}\ar[r] & \fun(\xtop^{op},\sring)
\end{tikzcd}
\]
is horizontally left adjointable. It suffices to show the diagram is vertically right adjointable, which amounts to the assertion that for $\Of_{\xtop}$ a sheaf of simplicial $\cinfty$-rings on $\xtop$, the functor $\xtop^{op}\overset{\Of_{\xtop}}{\rightarrow}\sring\overset{\iota_{c*}}{\rightarrow}\sring_c$ is already a sheaf, which is the case as $\iota_{c*}$ is a right adjoint. This completes the proof of $(1)$. Note that $(2)$ follows immediately from the analysis just performed and Remark \ref{rmk:logificationfunctor}. It follows from \cite{HA}, Lemma 2.2.1.11 that $p^{\xtop}_{\mathsf{Log}}$ is a coCartesian fibration. Since $p^{\xtop}_{\mathsf{Log}}$ is also a Cartesian fibration, we complete the proof of $(3)$ by noting that it has presentable fibres because the inclusion $\mathsf{Log}((\Of_{\xtop})_{\geq 0})\subset \shv_{s\mathsf{CMon}}(\xtop)_{/(\Of_{\xtop})_{\geq 0}}$ preserves sifted colimits, as colimits are universal in $\xtop$. 
\end{proof}
Let $\Of_{\xtop}$ be a sheaf of simplicial $\cinfty$-rings on an \inftop $\xtop$, then we say that a morphism $f:\mathcal{M}\rightarrow\mathcal{N}$ of sheaves of simplicial commutative monoids over $(\Of_{\xtop})_{\geq 0}$ \emph{exhibits $\mathcal{N}$ as a logification of $\mathcal{M}$ (with respect to $(\Of_{\xtop})_{\geq 0}$)} if $\mathcal{N}$ is a log structure over $(\Of_{\xtop})_{\geq 0}$ and $f$ is a localization morphism for the adjunction constructed in $(1)$ of Proposition \ref{prop:sheaveslogification}.
\begin{prop}\label{prop:localringcornerproperties}
The following hold true.
\begin{enumerate}[$(1)$]
	\item Let $\xtop$ be an \inftopt, then a left exact functor $\geodiffderc\rightarrow\xtop$ is a $\geodiffderc$-structure if and only if the composition $\geodiffder\rightarrow \geodiffderc \rightarrow \xtop$ is a $\geodiffder$-structure, and a morphism $\Of_{\xtop}\rightarrow \Of'_{\xtop}$ of left exact functors $\geodiffderc\rightarrow\xtop$ is local if and only if the underlying morphism $s_c^*\Of_{\xtop}\rightarrow s_c^*\Of_{\xtop}'$ is local. In other words, the two squares in the diagram 
	\[
	\begin{tikzcd}
	\fun^{\lex,\loc}(\geodiffderc,\xtop) \ar[d,"s_c^*"]\ar[r] & \fun^{\lex}(\geodiffderc,\xtop)\ar[d,"s_c^*"]  & \str_{\geodiffderc}(\xtop) \ar[d,"s_c^*"]\ar[l] \\
	\fun^{\lex,\loc}(\geodiffder,\xtop) \ar[r] & \fun^{\lex}(\geodiffder,\xtop)  & \str_{\geodiffder}(\xtop) \ar[l] 
	\end{tikzcd}
	\]
	are pullbacks.  
    \item The functor $s_c^*:\strloc_{\geodiffderc}(\xtop)\rightarrow\strloc_{\geodiffder}(\xtop)$ induced by composition with $s_c$ is a presentable fibration. The fibre over a $\geodiffder$-structure $\Of_{\xtop}$ can be identified with the \infcat of log structures on $(\Of_{\xtop})_{\geq 0}$.
    \item The functor $s_c^*$ admits a left adjoint that carries each $\Of_{\xtop}$ to an initial object in the fibre over $\Of_{\xtop}$.
    \item The functor $\strloc_{\geodiffderc}(\spa)\rightarrow \sring_c$ is fully faithful and preserves all colimits.
    \item The object $(\R,\R_{\geq0})\in\strloc_{\geodiffderc}(\spa)$ is final. Let $L_{(\R,\R_{\geq0})}:(\sring_{c})_{/(\R,\R_{\geq0})}\rightarrow \strloc_{\geodiffder}(\spa)$ be the left adjoint to the fully faithful functor $\strloc_{\geodiffderc}(\spa)\simeq \strloc_{\geodiffderc}(\spa)_{/(\R,\R_{\geq0})}\rightarrow (\sring_{c})_{/(\R,\R_{\geq0})}$ provided by Proposition \ref{prop:localfactorization}, then a map $f:(A,A_c)\rightarrow (B,B_c)$ is a unit transformation at $(A,A_c)\rightarrow (\R,\R_{\geq0})$ if and only if the map $A\rightarrow B$ is a unit transformation at $A\rightarrow \R$ for the functor $L_{\R}$ and $f$ is $p_{\mathsf{Log}}$-coCartesian, that is $A_c\rightarrow B_c$ exhibits a logification of prelog structures over $B_{\geq 0}$.
    \end{enumerate}
    \end{prop}
\begin{proof}
The assertion $(1)$ follows immediately from the fact that $\geodiffderc$ is endowed with the coarsest geometry structure for which $s_c$ is a transformation of geometries. It follows from $(1)$ that the functor $s_c^*:\strloc_{\geodiffderc}(\xtop)\rightarrow\strloc_{\geodiffder}(\xtop)$ is a pullback of the functor $\shv_{\sring_c}(\xtop)\rightarrow \shv_{\sring}(\xtop)$ induced by composition with $p_{\mathsf{Log}}$, which together with Proposition \ref{prop:sheaveslogification} implies $(2)$. The assertion $(3)$ follows from $(2)$ and \cite{HTT}, Proposition 5.2.4.3. Applying $(1)$ in the case $\xtop=\spa$ yields a pullback diagram
\[
\begin{tikzcd}
\strloc_{\geodiffderc}(\spa)\ar[r] \ar[d] & \sring_c\ar[d,"p_{\mathsf{Log}}"] \\
\strloc_{\geodiffder}(\spa) \ar[r] & \sring
\end{tikzcd}
\]  
so $(4)$ follows from Proposition \ref{prop:localringproperties} and \cite{HTT}, Lemma 5.4.5.5. To prove $(5)$, we first observe that $(\R,\R_{\geq0})$ is a final object in the fibre over a final object in $\strloc_{\geodiffder}(\spa)$ so this object is final because $p_{\mathsf{Log}}$ is a Cartesian fibration. We have pullback diagrams
\[
\begin{tikzcd}
\strloc_{\geodiffderc}(\spa)_{/(\R,\R_{\geq0})}\ar[r] \ar[d] & (\sring_c)_{/(\R,\R_{\geq0})}\ar[d,"p_{\mathsf{Log}}"] \\
\strloc_{\geodiffder}(\spa)_{/(\R,\R_{\geq0})} \ar[r] & \sring_{/(\R,\R_{\geq0})}
\end{tikzcd}
\]  
of simplicial sets, so we conclude by invoking Proposition \ref{prop:pullbackcocartleftadj}.
\end{proof}    
\begin{cons}
Recall the biCartesian fibration $\overline{\rtop}^{op}\rightarrow\ltop$. Define a simplicial set $\widetilde{\ltop(\sring_c)}$ over $\ltop$ by the universal property that for any map of simplicial sets $K\rightarrow\ltop$, there is a canonical bijection
\[ \Hom_{(\sset)_{/\ltop}}(K,\widetilde{\ltop(\sring_c)})\cong \Hom_{\sset}(K\times_{\ltop}\overline{\rtop}^{op},\sring_c)  ,\]
then it follows from \cite{HTT}, Corollary 3.2.2.12 that in the commuting diagram of simplicial sets
\[ \begin{tikzcd}
\widetilde{\ltop(\sring_c)}\ar[dr,"q_{\sring_c}"']\ar[rr,"p_{\mathsf{LogTop}}"] && \widetilde{\ltop(\sring)}\ar[dl,"q_{\sring}"] \\ &\ltop
\end{tikzcd} \]
all three maps are biCartesian fibrations and $p_{\mathsf{LogTop}}$ carries $q_{\sring_c}$-Cartesian and -coCartesian edges to $q_{\sring}$-Cartesian and -coCartesian edges. Let $\ltop(\sring_c)\subset \widetilde{\ltop(\sring_c)}$ be the full subcategory spanned by functors $\xtop^{op}\rightarrow\sring_c$ that preserve small limits, then $p_{\mathsf{LogTop}}$ restricts to a Cartesian fibration $\ltop(\sring_c)\rightarrow\ltop(\sring)$ that we also denote by $p_{\mathsf{LogTop}}$. We will call objects of $\ltop(\sring_c$) \emph{simplicial $\cinfty$-ringed \inftopoi with corners} with some disregard for good grammatical practice (\inftopoi do not have corners), and write them as triples $(\xtop,\Of_{\xtop},\mathcal{M})$.
\end{cons}
The preceding construction provides a Cartesian fibration $p_{\mathsf{LogTop}}:\ltop(\sring_c)\rightarrow\ltop(\sring)$ the fibre of which over some $\ofxtop$ is the \infcat of log structures over $(\Of_{\xtop})_{\geq 0}$. For $f^*:\xtop\rightarrow\ytop$ an algebraic morphism, the induced functor $\shv_{\sring_c}(\ytop)\rightarrow\shv_{\sring_c}(\xtop)$ can be identified with the functor $\fun^{\lex}(\geodiffderc,\ytop)\rightarrow\fun^{\lex}(\geodiffderc,\xtop)$ that composes with $f_*$, the geometric morphism right adjoint to $f^*$, so it has a left exact left adjoint given by composing with $f^*$. Thus, $q_{\sring_c}$ is also a coCartesian fibration. A sheaf of simplicial $\cinfty$-rings with corners is the data of a pair $(\Of_{\xtop},\mathcal{M}\rightarrow (\Of_{\xtop})_{\geq 0})$ of a sheaf of simplicial $\cinfty$-rings and a sheaf of simplicial commutative monoids over $(\Of_{\xtop})_{\geq 0}$, which the functor $f^*$ carries to a pair $f^*(\Of_{\xtop},\mathcal{M}\rightarrow (\Of_{\xtop})_{\geq 0})=(\Of_{\ytop},\mathcal{N}\rightarrow (\Of_{\ytop})_{\geq 0})$. We now show that $\Of_{\ytop}$ and $\mathcal{N}$ can be identified with $f^*\Of_{\xtop}$ and $f^*\mathcal{M}$, the pullback of sheaves of simplicial $\cinfty$-rings and sheaves of simplicial commutative monoids respectively, and we will describe the map $f^*\mathcal{M}\rightarrow (f^*\Of_{\xtop})_{\geq 0}$ explicitly.
\begin{prop}\label{prop:pullbacklogcomponents}
Let $f^*:\xtop\rightarrow\ytop$ be an algebraic morphism of \inftopoit, then the following hold.
\begin{enumerate}[$(1)$]
\item The strictly commuting diagram 
\[
\begin{tikzcd}
\shv_{\sring_c}(\ytop)\ar[r]\ar[d,"p^{\ytop}_{\mathsf{Log}}"] & \shv_{\sring_c}(\xtop)\ar[d,"p^{\xtop}_{\mathsf{Log}}"] \\
\shv_{\sring}(\ytop)\ar[r] & \shv_{\sring}(\xtop)
\end{tikzcd}
\]
of \infcats where the vertical maps are given by composition with the limit preserving functor $p_{\mathsf{Log}}$ is horizontally left adjointable.
\item The strictly commuting diagram 
\[
\begin{tikzcd}
\shv_{\sring_c}(\ytop)\ar[r]\ar[d] & \shv_{\sring_c}(\xtop)\ar[d] \\
\shv_{s\mathsf{CMon}}(\ytop)\ar[r] & \shv_{s\mathsf{CMon}}(\xtop)
\end{tikzcd}
\]
of \infcats where the vertical functors are given by composition with the limit preserving functor 
\[\sring_c\hooklongrightarrow\sring_{pc}\overset{\ev_{\mathsf{PLog}}}{\longrightarrow} s\mathsf{CMon}\]
 is horizontally left adjointable. 
\end{enumerate}
\end{prop}
\begin{proof}
Since the functors $\sring_c\rightarrow \sring$ and $\sring_c\rightarrow s\mathsf{CMon}$ factors through $\sring_c\subset\sring_{pc}$, it suffices to show that the diagram 
\begin{equation}\label{eq:la1}
\begin{tikzcd}
\shv_{\sring_c}(\ytop)\ar[r]\ar[d] & \shv_{\sring_c}(\xtop)\ar[d] \\
\shv_{\splring}(\ytop)\ar[r] & \shv_{\splring}(\xtop)
\end{tikzcd}
\end{equation}
is horizontally left adjointable and that the diagrams
\begin{equation}\label{eq:la2}
\begin{tikzcd}
\shv_{\splring}(\ytop)\ar[r]\ar[d,"p^{\ytop}_{\mathsf{PLog}}"] & \shv_{\splring}(\xtop)\ar[d,"p^{\xtop}_{\mathsf{PLog}}"] \\
\shv_{\sring}(\ytop)\ar[r] & \shv_{\sring}(\xtop)
\end{tikzcd} \quad\quad
\begin{tikzcd}
\shv_{\splring}(\ytop)\ar[r]\ar[d] & \shv_{\splring}(\xtop)\ar[d] \\
\shv_{s\mathsf{CMon}}(\ytop)\ar[r] & \shv_{s\mathsf{CMon}}(\xtop)
\end{tikzcd}
\end{equation}
are horizontally left adjointable. The vertical functors in the diagram \eqref{eq:la1} admit left adjoints, so it suffices that the associated diagram of left adjoints 
\[
\begin{tikzcd}
\shv_{\sring_c}(\ytop)& \shv_{\sring_c}(\xtop)\ar[l] \\
\shv_{\splring}(\ytop)\ar[u] & \shv_{\splring}(\xtop)\ar[u]\ar[l]
\end{tikzcd}
\]
is vertically right adjointable. Since the right adjoints are fully faithful, this amounts to the following  assertion.
\begin{enumerate}
\item[$(*)$] The pullback functor $\shv_{\splring}(\xtop)\rightarrow \shv_{\splring}(\ytop)$ carries the full subcategory $\shv_{\sring_c}(\xtop)$ into $\shv_{\sring_c}(\ytop)$. 
\end{enumerate}
We prove $(*)$. It follows from Proposition \ref{prop:plogprojectivelygenerated} that $\splring$ is compactly generated. Let $\icat\subset \splring$ denote the full subcategory spanned by compact objects, then the pullback of $\splring$-valued sheaves coincides with the functor
\[ \fun^{\lex}(\icat,\xtop)\longrightarrow \fun^{\lex}(\icat,\ytop),\quad \quad \Of_{\xtop}\longmapsto f^*\circ\Of_{\xtop}.  \]
It follows from Corollary \ref{cor:accessiblecorner} that $L_{\mathsf{Log}}:\splring\rightarrow\slring$ carries $\icat$ into $\geodiffder$, so that the functor $\shv_{\sring_c}(\xtop)\subset\shv_{\splring}(\xtop)$ is given by composing left exact functors $\geodiffderc\rightarrow \xtop$ with $L_{\mathsf{Log}}$. It follows that we can identify the full subcategory of $\shv_{\splring}(\xtop)$-valued sheaves spanned by sheaves taking values in $\sring_c$ with the full subcategory of those left exact functors $\icat\rightarrow\xtop$ that factor (up to homotopy) through $L_{\mathsf{Log}}:\icat\rightarrow\geodiffderc$. Obviously, the functor $f^*\circ\_:\fun^{\lex}(\icat,\xtop)\rightarrow\fun^{\lex}(\icat,\ytop)$ carries this full subcategory to the corresponding full subcategory of $\shv_{\sring_c}$-valued sheaves.\\
We now show that the left diagram of \eqref{eq:la2} is horizontally left adjointable. It suffices to show that the diagram is vertically right adjointable. This follows from the fact that the functors $p_{\mathsf{PLog}}^{\xtop}$ and $p_{\mathsf{PLog}}^{\ytop}$ admit fully faithful right adjoints given by composing with the fully faithful right adjoint to $p_{\mathsf{PLog}}$, which carries $A\in\sring$ to $(A,A_{\geq 0}\overset{\mathrm{id}}{\rightarrow}A_{\geq 0})$. Since the functor $\ev_{\mathsf{PLog}}:\splring\rightarrow s\mathsf{CMon}$ also admits a fully faithful right adjoint, the same argument shows that the right diagram of \eqref{eq:la2} is horizontally left adjointable.
\end{proof}
\begin{rmk}
Let $f^*:\xtop\rightarrow \ytop$, the diagram 
\[
\begin{tikzcd}
\shv_{\sring_c}(\ytop)\ar[r]\ar[d] & \fun(\Delta^1,\shv_{\sring_c}(\xtop))\ar[d] \\
\shv_{s\mathsf{CMon}}(\ytop)\ar[r] & 
\fun(\Delta^1,\shv_{s\mathsf{CMon}}(\xtop))
\end{tikzcd}
\]
where the vertical functors compose with the assignment $(A,A_c\rightarrow A_{\geq 0})\mapsto (A_c\rightarrow A_{\geq 0})$ is \emph{not} horizontally left adjointable (because composing sheaves of simplicial $\cinfty$-rings with the functor $A\mapsto A_{\geq 0}$ does not furnish a left adjointable square). Thus, taking horizontal left adjoints determines for each $(\xtop,\Of_{\xtop},\mathcal{M})$ a diagram 
\[
\begin{tikzcd}
f^*\mathcal{M}\ar[d]\ar[r] & f^*\mathcal{M}\ar[d]\\
f^*((\Of_{\xtop})_{\geq 0})\ar[r] & (f^*\Of_{\xtop})_{\geq 0}
\end{tikzcd}
\]
where the horizontal morphisms are the Beck-Chevalley maps associated to the square above. From Proposition \ref{prop:pullbacklogcomponents} we see that the upper horizontal map is an equivalence.
\end{rmk}
In fact, Proposition \ref{prop:pullbacklogcomponents} implies that the functor $p_{\mathsf{LogTop}}$ is a coCartesian fibration as well. 

\begin{prop}\label{prop:logtopcocart}
The functor $p_{\mathsf{LogTop}}:\ltop(\sring_c)\rightarrow\ltop(\sring)$ is a presentable fibration. A morphism $f:(\xtop,\Of_{\xtop},\mathcal{M})\rightarrow (\ytop,\Of_{\ytop},\mathcal{N})$ of simplicial $\cinfty$-ringed \inftopoi with corners is $p_{\mathsf{LogTop}}$-coCartesian if and only if it factors as a composition
\[(\xtop,\Of_{\xtop},\mathcal{M})\longrightarrow (\xtop,f^*\Of_{\xtop},f^*\mathcal{M}) \longrightarrow(\ytop,\Of_{\ytop},\mathcal{N}) \]
where $f^*:\xtop\rightarrow\ytop$ is the algebraic morphism underlying $f$, the first map is $q_{\sring_c}$-coCartesian and the second map exhibits a logification with respect to $(\Of_{\ytop})_{\geq0}$.   
\end{prop}
\begin{proof}
It is clear that $p_{\mathsf{LogTop}}$ is a Cartesian fibration with presentable fibres and Proposition \ref{prop:sheaveslogification} asserts that for each \inftop $\xtop$, the fibre $p_{\mathsf{LogTop}}^{-1}(\xtop)\rightarrow q^{-1}_{\sring}(\xtop)$ is a presentable fibration. To show that $p_{\mathsf{LogTop}}$ is a presentable fibration, it suffices to show that for each algebraic morphism $f^*:\xtop\rightarrow \ytop$, the diagram 
\[
\begin{tikzcd}
\shv_{\ytop}(\sring_c)\ar[r,"f_*"]\ar[d,"p^{\xtop}_{\mathsf{Log}}"] & \shv_{\xtop}(\sring_c)\ar[d,"p^{\ytop}_{\mathsf{Log}}"] \\
\shv_{\ytop}(\sring)\ar[r,"f_*"]& \shv_{\xtop}(\sring)
\end{tikzcd}
\]
is horizontally left adjointable, in view of Lemma \ref{lem:rcocart}. This follows immediately from Proposition \ref{prop:pullbacklogcomponents}. The description of $p_{\mathsf{LogTop}}$ is also a consequence of Proposition \ref{prop:pullbacklogcomponents}.
\end{proof}
 \begin{prop}
 The following hold true.
 \begin{enumerate}[$(1)$]  
    \item The diagram 
    \[
    \begin{tikzcd}
    \ltop^{\loc}(\sring_c) \ar[d,"p_{\mathsf{LogTop}}"] \ar[r] & \ltop(\sring_c)\ar[d,"p_{\mathsf{LogTop}}"] \\
    \ltop^{\loc}(\sring) \ar[r]& \ltop(\sring)
    \end{tikzcd}
    \]  
 is a pullback diagram of \infcatst.
 \item The functor $\ltop^{\loc}(\sring_c)\rightarrow \ltop^{\loc}(\sring_c)$ is a presentable fibration. The fibre over a locally simplicial $\cinfty$-ringed \inftop $\ofxtop$ can be identified with the \infcat of log structures on $(\Of_{\xtop})_{\geq 0}$.
 \item Consider the relative spectrum 
 \[ \ltop(\sring_c)\longrightarrow \ltop^{\loc}(\sring_c), \]
 then a map $f:(\xtop,\Of_{\xtop},\mathcal{M})\rightarrow (\ytop,\Of_{\ytop},\mathcal{N})$ of simplicial $\cinfty$-ringed \inftopoi with corners exhibits a unit transformation if and only if the map $\ofxtop\rightarrow\ofytop$ exhibits a relative spectrum for the geometry $\geodiffder$ and $f$ is $p_{\mathsf{LogTop}}$-coCartesian, that is, we have a pushout diagram 
 \[ 
 \begin{tikzcd}
f^*\mathcal{M}\times_{(\Of_{\ytop})_{\geq 0}}(\Of_{\ytop})_{\geq 0}^{\times} \ar[r]\ar[d] & (\Of_{\ytop})_{\geq 0}^{\times} \ar[d]\\
 f^*\mathcal{M} \ar[r] & \mathcal{N}
 \end{tikzcd}
 \] 
 of sheaves of simplicial commutative monoids on $\ytop$, where $f^*:\xtop\rightarrow \ytop$ is the algebraic morphism underlying $f$.
    \item The relative spectrum $\spec^{\geodiffderc}_{\geodiffder}$ may be identified with the section of $\ltop^{\loc}(\sring_c)\rightarrow \ltop^{\loc}(\sring)$ that carries each $(\xtop,\Of_{\xtop})$ to an initial object in the fibre over $(\xtop,\Of_{\xtop})$; in particular, $\spec^{\geodiffderc}_{\geodiffder}$ is fully faithful. 
\end{enumerate}
\end{prop}
\begin{proof}
Since $\ltop^{\loc}(\sring_c)\subset\ltop(\sring)$ is the subcategory whose objects are $\geodiffderc$-structures and whose morphisms are local morphisms, $(1)$ follows immediately from Proposition \ref{prop:localringcornerproperties}. Now $(2)$ is a consequence of $(1)$ and Proposition \ref{prop:logtopcocart}. The assertion $(3)$ follows from $(1)$ and Proposition \ref{prop:pullbackcocartleftadj}. The assertion $(4)$ follows from $(1)$ and \cite{HTT}, Proposition 5.2.4.3.
\end{proof}
Let $\mathsf{Top}^{\mathrm{Pc}}(\sring_c)$ denote the full subcategory of $\rtop(\sring_c)$ spanned by triples $(\xtop,\Of_{\xtop},\mathcal{N})$ such that $\xtop$ is Postnikov complete, weakly $(-1)$-complicial and spatial. Identifying spatial weakly $(-1)$-complicial Postnikov complete \inftopoi with sober topological spaces, we have a biCartesian fibration $\mathsf{Top}^{\mathrm{Pc}}(\sring_c)\rightarrow\mathsf{Top}$ that we also denote $q_{\sring_c}$. The fibre over $X$ is the \infcat $\widehat{\shv}_{\sring_c}(X)^{op}$, the opposite of the Postnikov complete sheaves of simplicial $\cinfty$-rings with corners over $X$. 
\begin{rmk} Consider the full subcategory of $\mathsf{Top}^{\mathrm{Pc}}(\sring_c)$ spanned by pairs $(X,\Of_X)$ where $\Of_X:\widehat{\shv}(X)\rightarrow \sring_c$ takes values in $\cinfty\mathsf{ring}_c$. This is a 1-category that can be identified with the nerve of the category defined as follows.
\begin{enumerate}
\item[$(O)$] Objects are triples $(X,\Of_X,\mathcal{M})$ where $X$ is a sober topological space, $\Of_X$ is a sheaf of $\cinfty$-rings on $X$ and $\mathcal{M}$ is a sheaf of commutative monoids on $X$ that is also a log structure on $(\Of_X)_{\geq 0}$.
\item[$(M)$] Morphisms are triples $(f,\alpha,\beta):(X,\Of_X,\mathcal{M})\rightarrow (Y,\Of_Y,\mathcal{N})$ where $f:X\rightarrow Y$ is a continuous map, $\alpha:f^*\Of_Y\rightarrow \Of_X$ is a morphism of sheaves of $\cinfty$-rings and $\beta:\mathcal{N}\rightarrow\mathcal{M}$ is a morphism of prelog structures over $(\Of_{X})_{\geq 0}$, where we view $f^*\mathcal{N}$ as a prelog structure over $(\Of_{X})_{\geq 0}$ via the map $f^*\mathcal{N}\rightarrow (f^*\Of_{Y})_{\geq 0}\overset{\alpha}{\rightarrow} (\Of_X)_{\geq 0}$.
\end{enumerate}  
We will denote this category by $\mathsf{Top}(\cinfty\mathsf{ring}_c)$. 
\end{rmk}
Let $\mathsf{Top}^{\mathrm{Pc},\loc}(\sring_c)$ be the intersection $\ltop^{\loc}(\sring_c)\cap\mathsf{Top}^{\mathrm{Pc}}(\sring_c)$. As for locally simplicial $\cinfty$-ringed spaces, we have the following.
\begin{prop}
The subcategory inclusion $\mathsf{Top}^{\mathrm{Pc},\loc}(\sring_c)\subset \mathsf{Top}^{\mathrm{Pc}}(\sring_c)$ is full and is stable under small limits. 
\end{prop} 
\begin{proof}
This follows from $(4)$ of Proposition \ref{prop:localringcornerproperties}.
\end{proof}
It follows from the preceding results that the diagram
\[
\begin{tikzcd}
\mathsf{Top}^{\mathrm{Pc},\mathrm{loc}}(\sring_c)^{op}\ar[d]\ar[r,"\Gamma"] & \sring_c\ar[d,"p_{\mathsf{Log}}"] \\
\mathsf{Top}^{\mathrm{Pc},\loc}(\sring)^{op} \ar[r,"\Gamma"] & \sring
\end{tikzcd}
\]
is horizontally left adjointable. We let $\spec_c$ denote a left adjoint to the upper global sections functor; it carries a pair $(A,M\rightarrow A_{\geq 0})$ to the triple $(\specr\,A,\Of_{\specr^{\wedge}\,A},\mathcal{M})$, where $\mathcal{M}$ is the logification with respect to $(\Of_{\specr^{\wedge}\,A})_{\geq0}$ of the sheafification of the constant presheaf of simplicial commutative monoids with value $M$. The diagram above induces for each simplicial $\cinfty$-ring $A$ a functor 
\[ \mathrm{LSpec}_A:\mathsf{Log}(A_{\geq0})\longrightarrow \mathsf{Log}((\Of_{\specr^{\wedge}\,A})_{\geq 0}). \]
\begin{prop}\label{prop:logbasechange}
The following hold true. 
\begin{enumerate}[$(1)$]
\item Let $A$ be a simplicial $\cinfty$-ring, then the functor $\mathrm{LSpec}_A$ factors as 
\[ \mathsf{Log}(A_{\geq0})\longrightarrow \mathsf{Log}_{\underline{A}_{\geq 0}} \longrightarrow \mathsf{Log}_{(\widetilde{\Of_{\speco\,A}})_{\geq0}} \longrightarrow \mathsf{Log}_{(\Of_{\specr^{\wedge}\,A})_{\geq 0}},  \]
where the first map is induced by the algebraic morphism $i^*:\spa\rightarrow \pshv((\sring^{\mathrm{ad}}_{A/})^{op})$ from the initial \inftop and $\underline{A}=i^*A$ is the constant presheaf on $A$, the second functor takes the logification with respect to the presheaf $(\widetilde{\Of_{\speco\,A}})_{\geq0}$ and the third functor is induced by the geometric morphism $\pshv((\sring^{\mathrm{ad}}_{/A})^{op})\rightarrow\widehat{\shv}(\specr\,A)$ that sheafifies, takes the spatial reflection and the Postnikov completion.  
\item The functor $\spec_c$ carries $p_{\mathsf{Log}}$-coCartesian edges to $p_{\mathsf{LogTop}}$-coCartesian edges. In particular, for any map $A\rightarrow B$ of simplicial $\cinfty$-rings, there is a canonical homotopy rendering the diagram 
\[
\begin{tikzcd}
\mathsf{Log}(A_{\geq 0})\ar[d,"\mathrm{LSpec}_A"]\ar[r] &\mathsf{Log}(B_{\geq 0}) \ar[d,"\mathrm{LSpec}_B"] \\
\mathsf{Log}((\Of_{\specr^{\wedge}\,A})_{\geq 0}) \ar[r] & \mathsf{Log}((\Of_{\specr^{\wedge}\,B})_{\geq 0})  
\end{tikzcd}
\]
commutative.
\end{enumerate}
\end{prop}
\begin{proof}
The first assertion is an immediate consequence of the definition of $\spec_c$ and the second follows from Lemma \ref{lem:ajointablepreserve}.
\end{proof}
Let $A\rightarrow B$ an effective epimorphism of Lindel\"{o}f simplicial $\cinfty$-rings. We have seen that the for modules instead of log structures, the analogous diagram of $(2)$ of Proposition \ref{prop:logbasechange} of modules is vertically right adjointable. Applying this to $\R$-points $\R\rightarrow A$ shows that for any sheaf $\F$ of $\Of_{\specr^{\wedge}\,A}$-modules, the counit map $\mspec^{\wedge}_A\,\Gamma(\F)\rightarrow \F$ is an equivalence, so that $\Gamma$ is fully faithful. For logarithmic structures, this is not the case. To see this, we consider the category of manifolds with corners and interior $b$-maps among them, which we can identify with the full subcategory of $\mathsf{Top}^{\mathrm{Pc},\loc}(\sring_c)$ spanned by triples $(X,\Of_X,\mathcal{M})$ for which the following conditions are satisfied.
\begin{enumerate}[$(1)$]
\item $X$ is paracompact Hausdorff.
\item There exists an open cover $\{U_i\rightarrow X\}_i$ such that $(U,\Of_X|_U,\mathcal{M}|_U)$ is equivalent to $\spec_c(\cinfty(V),\cinfty_b(V))$ for some open $V\subset\R^{n}\times\R_{\geq0}^k$. 
\end{enumerate}
Then $(X,\Of_X)$ is an affine derived geometric $\cinfty$-scheme by Theorem \ref{thm:spectrumglobalsections}, so the counit map $\spec_c\Gamma(\Of_X,\mathcal{M})\rightarrow  (X,\Of_X,\mathcal{M})$ is an equivalence if and only if for each $x\in X$, the diagram in $(2)$ of Proposition \ref{prop:logbasechange} is vertically right adjointable for the map $x:\Gamma(\Of_X)\rightarrow \R$. Unwinding the definitions, the Beck-Chevalley transformation at $(X,\Of_X,\mathcal{M})$ is the map on associated log structures of the map $\Gamma(\mathcal{M})\rightarrow \mathcal{M}_x$, where $\mathcal{M}_x$ is the stalk of $\mathcal{M}$ at $x$. We conclude that a manifold with corners $X$ is an \emph{affine} derived $\cinfty$-scheme with corners if and only if for each $x\in X$, the map $\cinfty_b(X)\rightarrow \cinfty_{bX,x}$ exhibits a logification with respect to $(\cinfty(X)_x)_{\geq 0}$, where $\cinfty_{bX,x}$ is the stalk at $x$ of the sheaf $\cinfty_{bX}$ of interior $b$-maps on $X$. This is the case precisely if $X$ is a manifold with \emph{faces}.  
\begin{rmk}
Let $\ofxtop$ be a derived $\cinfty$-scheme and let $\mathcal{M}$ be a sheaf of log structures over $(\Of_{\xtop})_{\geq 0}$. For any $U\in\xtop$, we let $(\xtop_{/U},\Of_{\xtop}|_U,\mathcal{M}|_U)$ denote the codomain of a $q_{\sring_c}$-coCartesian lift of the algebraic morphism $\xtop\rightarrow\xtop_{/U}$ starting at $(\xtop,\Of_{\xtop},\mathcal{M})$. A log structure over $(\Of_{\xtop})_{\geq 0}$ is \emph{quasi-coherent} if there is an effective epimorphism $\coprod_iU_i\rightarrow 1_{\xtop}$ such that $(\xtop_{/U_i},\Of_{\xtop}|_{U_i},\mathcal{M}|_{U_i})$ is an affine derived $\cinfty$-scheme with corners. Likewise, we say that a morphism $\mathcal{M}\rightarrow\mathcal{N}$ of log structures over $(\Of_{\xtop})_{\geq 0}$ is \emph{quasi-coherent} if there exists an effective epimorphism $\coprod_i U_i\rightarrow 1_{\xtop}$ such that $(\xtop_{/{U_i}},\Of_{\xtop}|_{U_i},\mathcal{M}|_{U_i})\rightarrow (\xtop_{/{U_i}},\Of_{\xtop}|_{U_i},\mathcal{N}|_{U_i})$ lies in the image of the functor $\spec_c$. Let $\mathsf{QLog}((\Of_{\xtop})_{\geq 0})\subset \mathsf{Log}((\Of_{\xtop})_{\geq 0})$ be the subcategory whose objects are quasi-coherent log structures and whose morphisms are quasi-coherent morphisms. The assignment 
\[  \ofxtop\longmapsto \mathsf{QLog}((\Of_{\xtop})_{\geq 0})  \]
determines a sheaf of \infcats on $\mathsf{d}\cinfty\sch$. In future work, we will consider full subsheaves of $\mathsf{QLog}$ well adapted to the problem of blowing up derived $\cinfty$-schemes and resolving singularities. 
\end{rmk}
\begin{rmk}
The spectrum functor is compatible with truncations: let $\geodiffc$ be the 0-stub of $\geodiffderc$, whose underlying \infcat is the opposite of the 1-category of compact objects of $\cinfty\mathsf{ring}_c$, then the canonical transformation of geometries $f:\geodiffderc\rightarrow\geodiffc$ identifies $\ltop(\cinfty\mathsf{ring}_c)$ with the full subcategory of $\ltop(\sring_c)$ spanned by triples $(\xtop,\Of_{\xtop},\mathcal{M})$ for which $\Of_{\xtop}$ and $\mathcal{M}$ are $0$-truncated. Since $\geodiffderc$ is a geometric envelope of a $0$-truncated pregeometry, it follows from Proposition \ref{prop:pgeotruncfunctor} that the inclusion $\ltop(\cinfty\mathsf{ring}_c)\subset\ltop(\sring_c)$ has a left adjoint $\tau_{\leq 0}$ which carries a \emph{$\diffc'$-structure} $(\xtop,\Of)$ to $(\xtop,\tau_{\leq 0}\Of)$, where $\tau_{\leq 0}\Of_{\xtop}$ is the composition 
\[ \diffc'\overset{\Of}{\longrightarrow} \xtop\overset{\tau_{\leq 0}}{\longrightarrow}\xtop. \] 
We have a diagram
\[
\begin{tikzcd}
\sring_c\ar[d,"\tau_{\leq 0}"] \ar[r,"\spec_c"] & \mathsf{Top}^{\mathrm{Pc,\loc}}(\sring_c)^{op}\ar[d,"\tau_{\leq 0}"] \\
\cinfty\mathsf{ring}_c \ar[r,"\spec_c"] & \mathsf{Top}(\sring_c)^{op}
\end{tikzcd}
\]
which commutes up to canonical homotopy. This diagram is also vertically right adjointable, since truncatedness is preserved by algebraic morphisms and logification. For $(A,A_c)$ a $\cinfty$-ring with corners, the spectrum $(X,\Of_X,\mathcal{M})$ has been constructed before by Joyce and Francis-Staite in \cite{Joyfra}.
\end{rmk}
We conclude this section by finishing the proof of Theorem \ref{thm:geoenvcorners}.
\begin{proof}[Proof of Theorem \ref{thm:geoenvcorners} $(iv)$, $(v)$]
We have identified the 1-category $\diffc$ of manifolds with faces as the full subcategory of $\mathsf{Top}(\sring_c)$ spanned by triples $(X,\Of_X,\mathcal{M})$ that satisfy the following conditions.
\begin{enumerate}[$(1)$]
\item $X$ is paracompact Hausdorff.
\item $X$ may be covered by opens $\{U_i\subset X\}$ such that $(U,\Of_X|_U,\mathcal{M}|_U)$ is equivalent to $\spec_c\,(\cinfty(V),\cinfty_b(V))$ for some open $V\subset \R^n\times\R^k_{\geq 0}$.
\item The counit $\spec_c\,\Gamma(\Of_X,\mathcal{M})\rightarrow (X,\Of_X,\mathcal{M})$ is an equivalence.
\end{enumerate} 
For any open $V\subset \R^n\times\R^k_{\geq 0}$, the spectrum of $(\cinfty(V),\cinfty_b(V))$ lies in $\diffc$, so we have a functor $\spec_c:\diffc'\rightarrow \diffc$. To see that this is a transformation of pregeometries, it suffices to show that it preserves pullbacks along admissible maps. Since this is a local question, it suffices to consider admissible maps of the form $U\subset\R^n\times\R^k_{\geq0}$, where $U$ is an open subset with \emph{connected} boundary. Unwinding the definitions, it is enough to show that for such open subsets, the map 
\[ (\cinfty(\R^n\times\R^k_{\geq 0}),\cinfty_b(\R^n\times\R^k_{\geq 0})) \longrightarrow (\cinfty(U),\cinfty_b(U))  \]
is a coCartesian morphism in $\sring_c$. It follows from Theorem \ref{thm:logcinftycorners} that it suffices to show that the morphism is a logification. To see this, we note that the composite map 
\[\Z^k_{\geq 0}\longrightarrow \cinfty_b(\R^k)\longrightarrow\cinfty_b(U) \subset\cinfty_{\geq 0}(U) \]
where the first functor exhibits a logification and specifies the boundary defining functions, fits into a diagram 
\[
\begin{tikzcd}
0\ar[d]\ar[r] &  \Z^{k-|S|}_{\geq 0} \ar[d] \ar[r]& \cinfty_{>0}(U) \ar[d] \\
\Z^S_{\geq 0} \ar[r] & \Z^k_{\geq 0} \ar[r] & \cinfty_{\geq 0}(U)
\end{tikzcd}
\]
where $S\subset \{1,\ldots,k\}$ is the subset determining the boundary defining functions on $U$ (which indeed form a subset as $U$ has connected boundary). The left square is a pushout and the right square is a pullback, so $\Z^S_{\geq 0}\rightarrow \cinfty_b(U)$ exhibits a logification. \\
Let $\pregeo$ be the categorical mapping cylinder of the functor $\spec_c$ defined as follows.
\begin{enumerate}
    \item[$(O)$] An object of $\pregeo$ is either an object of $\diffc'$ or an object of $\diffc$.
    \item[$(M)$] Morphism sets are given by 
    \begin{align*}
        \Hom_{\pregeo}(M,N) = \begin{cases}
        \Hom_{\diffc'}(M,N) \quad  &M,N\in \diffc' \\
        \Hom_{\diffc}(M,N)  \quad &M,N\in \diffc\\
        \Hom_{\diffc}(\spec_c(M),N)  \quad &M\in \diffc', \,N\in \diffc\\
        \emptyset \quad &M\in \diffc, \,N\in \diffc'.
        \end{cases}
    \end{align*}
\end{enumerate} 
There are obvious full subcategory inclusions $i:\diffc'\subset \pregeo$ and $j:\diffc\subset\pregeo$. Note that the latter admits a retraction $r:\pregeo\rightarrow \diffc$ such that $r\circ i=\spec_c$ defined on objects by $r(M)=M$ if $M\in \diffc$ and $r(M)=\spec_c(M)$ if $M\in \diffc'$. We first argue the following.
\begin{enumerate}
    \item[$(*)$] Let $U\subset N$ be an open submanifold of a manifold with faces. Let $M \in \diffc'$ corresponding to some simplicial $\cinfty$-ring with corners $(A,A_c)$ and suppose that there is a map $M\rightarrow N$ in $\pregeo$. Then there is a pullback diagram 
    \[
    \begin{tikzcd}
    V\ar[d,"g"] \ar[r] & U\ar[d] \\
    M\ar[r] & N
    \end{tikzcd}
    \]
    with $V$ (necessarily) in $\diffc'$ and $g$ an admissible map.
\end{enumerate}
We note that we have some pullback diagram 
\[
\begin{tikzcd}
U\ar[d] \ar[r] & N\ar[d] \\
\R\setminus\{0\} \ar[r] & \R 
\end{tikzcd}
\]
of manifolds with faces, and it is not hard to see that the admissible map opposite to $(A,A_c)\rightarrow(A[a^{-1}],A[a^{-1}]_c)$ corresponding to the map $(\cinfty(\R),\cinfty_b(\R))\rightarrow \Gamma(N) \rightarrow (A,A_c)$ fits into the desired pullback diagram. It follows from \cite{dagv}, Lemma 1.2.14 that taking left Kan extensions along $\pregeo$ carries $\diffc'$-structures to $\diffc$-structures. Let $\fun'(\pregeo,\xtop)\subset\fun(\pregeo,\xtop)$ be the full subcategory spanned by functors $F:\pregeo\rightarrow\xtop$ such that
\begin{enumerate}[$(a)$]
    \item For each $M\in \diffc'$, the map $F(M)\rightarrow F(\spec_c(M))$ is an equivalence in $\xtop$.
    \item The restriction $F|_{\diffc}$ is a $\diffc$-structure on $\xtop$. 
\end{enumerate}
We note that $(a)$ is equivalent to the assertion that $F$ is a right Kan extension of $F|_{\diffc}$, so it follows from \cite{HTT}, Proposition 4.3.2.15 that the restriction \[\fun(\pregeo,\xtop)\longrightarrow\fun(\diffc,\xtop)\] induces a trivial Kan fibration $\fun'(\pregeo,\xtop)\rightarrow\str_{\diffc}(\xtop)$. Choose a section $s$ of this trivial fibration, then the map $\str_{\diff}(\xtop)\rightarrow\fun(\diffc',\xtop)$ factors as 
\[ \str_{\diff}(\xtop)\overset{s}{\longrightarrow}  \fun'(\pregeo,\xtop) \longrightarrow \fun(\diffc',\xtop) \]
where the second functor is induced by restriction along $i$. Thus, our work will be done once we show that the restriction \[\fun(\pregeo,\xtop)\longrightarrow\fun(\diffc',\xtop)\]
induces a trivial fibration $\fun'(\pregeo,\xtop)\rightarrow\str_{\diffc'}(\xtop)$. In view of \cite{HTT}, Proposition 4.3.2.15, it suffices to show that 
\begin{enumerate}[$(i)$]
    \item A functor $F:\pregeo\rightarrow \xtop$ lies in $\fun'(\pregeo,\xtop)$ if and only if its restriction $F|_{\diffc'}$ is a $\diffc'$-structure and $F$ is a left Kan extension of $F|_{\diffc'}$.
    \item Every functor $F_0\in \str_{\diffc'}(\xtop)$ admits a left Kan extension along $i:\diffc'\hookrightarrow\pregeo$. 
\end{enumerate}
We note that the (essential) smallness of $\pregeo$ and the presentability of $\xtop$ guarantee that $(ii)$ is satisfied, so we show $(i)$. We first show the `if' direction. Let $F:\pregeo\rightarrow\xtop$ be a left Kan extension of its restriction to $\diffc'$, which we assume is a $\diffc'$-structure. Then we should show that $F|_{\diffc}$ is a $\diffc$-structure and that for all $M\in \diffc'$, the map $F(M)\rightarrow F(\spec_c\,M)$ is an equivalence. It follows from \cite{dagv}, Lemma 1.2.14 that $F|_{\diffc}$ is a $\diffc$-structure. For the second assertion to be verified, we note that fact that $F|_{\diffc}$ and $F|_{\diffc'}$ are $\diffc$- and $\diffc'$-structures respectively guarantees that we may assume that subspace of $\R^{n}\times\R^k_{\geq0}$ associated to $M\in \diffc'$ has connected boundary, since any $M$ admits an admissible cover $\{U_i\rightarrow M\}_i$ such that each iterated intersection $U_{i_1}\times_M\ldots\times_M U_{i_n}$ has this property (in fact, we may assume that $M$ is equivalent to $(\cinfty(\R^m\times\R^{l}_{\geq 0}),\cinfty_b(\R^m\times\R^{l}_{\geq 0})$). This assumption implies in particular that the unit transformation $M\rightarrow \Gamma(\spec_c\,M)$ is an equivalence. In turn, this implies that the map $M\rightarrow \spec_c\,M$ in $\pregeo$ is an initial object in the \infcat $\diffc'\times_{\diffc'}\pregeo_{/\spec_c\,M}$. By definition of left Kan extension, we thus have $F(M)\simeq F(\spec_c\,M)$. Conversely, we assume that $F\in \fun'(\pregeo,\xtop)$, then $F|_{\diffc'}$ is a $\diffc'$-structure since we can identify $F|_{\diffc'}$ with the pullback of $F|_{\diffc}$ along $\spec_c:\diffc'\rightarrow\diffc$, which is a transformation of pregeometries. We need only show that $F$ is a left Kan extension of $F|_{\diffc'}$. Let $F'$ be a left Kan extension of $F_{\diffc'}$, then we invoke \cite{dagv}, Lemma 1.2.14 again to deduce that the left Kan extension $F'$ of $F_{\diffc'}$ is a $\diffc$-structure. The natural transformation $F'|_{\diffc}\rightarrow F|_{\diffc}$ is an equivalence on the full subcategory spanned by manifolds with faces that are open submanifolds of some $\R^{n}\times\R^{k}_{\geq 0}$ with connected boundary, so we conclude that $F'\simeq F$ as both functors are $\diffc$-structures.\\
It remains to be shown that the equivalence $\str_{\diffc}(\xtop)\simeq\str_{\diffc'}(\xtop)$ identifies the subcategories of local morphisms. To see this, we note that a morphism $\alpha:\Of\rightarrow\Of'$ of $\diffc$-structures is local if and only if the morphism of $\diff$-structures $\alpha|_{\diff}$ is local, where we view $\diff$ as a full subcategory of $\diffc$. Since the same observation holds for locality of morphisms of $\diffc'$-structures, we conclude. \end{proof}
\begin{rmk}
It is an immediate consequence of Theorem \ref{thm:geoenvcorners} (and the general theory of geometries and pregeometries) that there are preferred equivalences between
\begin{enumerate}[$(i)$]
    \item The \infcat $\rtop(\diffc)$ of \inftopoi equipped with $\diffc$-structures.
    \item The \infcat of locally simplicial $\cinfty$-ringed \inftopoi with corners.
    \item The \infcat of locally simplicial $\cinfty$-ringed \inftopoi equipped with positive log structures.
\end{enumerate}
These equivalence restrict to one between the \infcat of derived $\cinfty$-schemes with corners locally of finite presentation and the \infcat of $0$-localic $\geodiffderc$-schemes locally of finite presentation. A 1-categorical version of this result was obtained by Francis-Staite in her recent thesis \cite{FrSt}; she compared the positive log differentiable spaces of Gillam-Molcho with interior $\cinfty$-schemes with corners; both classes of objects form full subcategories of all of the equivalent four \infcats described above. 
\end{rmk}
\begin{rmk}
In applications to moduli theory, such as the construction of representing stacks for elliptic moduli problems in future work, derived $\cinfty$-schemes with corners will usually be locally given by a retract of the zero set of a section of a vector bundle over a manifold with faces. Such derived $\cinfty$-schemes with corners will have the simplest possible nontrivial corners/log structures: they have \emph{free sharpening}. \end{rmk}

\newpage
\appendix

\section{Colimits in Lawvere theories and resolutions of diagrams}
Let $A$ be a $\cinfty$-ring, and let $I$ be a finitely generated ideal in $A$. A choice of generators for the ideal $I$ determines a map $g_I:\cinfty(\R^k)\rightarrow A$ and we can identify the relative tensor product $A^{\rmalg}\otimes_{\cinfty(\R^k)}\R$ of commutative rings determined by $g_I$ and the map $\ev_0:\cinfty(\R^k)\rightarrow\R$ evaluating at 0 with the quotient ring $A^{\rmalg}/I$. Since $g_I$ and $\ev_0$ are morphisms of $\cinfty$-rings, we can also consider the pushout $A\oinfty_{\cinfty(\R^k)}\R$ of $\cinfty$-rings and we have a canonical comparison map $(A\oinfty_{\cinfty(\R^k)}\R)^{\rmalg}\rightarrow A^{\rmalg}\otimes_{\cinfty(\R^k)}\R$. It is a consequence of \emph{Hadamard's lemma} that this map is actually an isomorphism. Our goal in this appendix will be to generalize the situation just described in several directions.
\begin{enumerate}[$(1)$]
    \item We will replace the ordinary categories of $\cinfty$-rings and commutative $\R$-algebras with the \infcats of $\spa$-valued algebras.
    \item In place of the functor $(\_)^{\rmalg}$, we will consider the functor $f^*:s\mathrm{T}'\alg\rightarrow s\mathrm{T}\alg$ induced by an arbitrary transformation of Lawvere theories $f:\mathrm{T}\rightarrow\mathrm{T}'$.
    \item We will consider not just pushouts along free $\mathrm{T}'$-algebras, but arbitrary colimit diagrams in $s\mathrm{T}'\alg$.
\end{enumerate}
The main result proven in this appendix formulates a condition that guarantees that we can express colimits in $s\mathrm{T}'\alg$ in terms of those in $s\mathrm{T}\alg$ in a very precise sense. We highlight two important consequences: given a transformation $f$ satisfying the pushout axiom $(P)$ below, the canonical comparison map $\colim_K f^*\mathcal{J} \rightarrow f^*(\colim_K\mathcal{J})$ for \emph{any} diagram $\mathcal{J}:K\rightarrow s\mathrm{T}'\alg$ is a pushout of the comparison map for coproducts of free $\mathrm{T}'$-algebras, which are a priori the only colimits we understand. Moreover, $f$ \emph{preserves} pushouts along effective epimorphisms. These results are important computational tools, and we will invoke them frequently in this work and its successors.\\
We first single out a class of effective epimorphism of finitely generated free $\mathrm{T}$ algebras.  
\begin{defn}\label{defn:graphinclusion}
Let $\mathrm{T}$ be a Lawvere theory. A morphism in $\mathrm{T}$ is a \emph{graph inclusion} if it is equivalent to a morphism of the form $\mathrm{id}_X\times g:X\rightarrow X\times Y$ for some $g:X\rightarrow Y$.
\end{defn}

\begin{thm}\label{thm:unramified2}
Let $f:\mathrm{T}\rightarrow\mathrm{T}'$ be a transformation of Lawvere theories. Suppose that $f$ satisfies the following \emph{pushout axiom}.
\begin{enumerate}
    \item[$(P)$] For each pair of morphisms $g:X\rightarrow Y$ and $h:W\rightarrow Z$ in $\mathrm{T}'$, the natural diagram 
    \begin{equation}\label{eq:precolimcompcare}
    \begin{tikzcd}
    j(f(X)\times f(Y)\times f(W)\times f(Z)) \ar[r] \ar[d] &  j(f(X)\times f(W)) \ar[d] \\
    jf(X\times Y\times W\times Z) \ar[r] & j(f(X\times W))
    \end{tikzcd}
    \end{equation}
    where the upper horizontal map is induced by the two graph inclusions $X\rightarrow X\times Y$ and $W\rightarrow W\times Z$ and the lower horizontal map is induced by the graph inclusion $X\times W\rightarrow X\times Y\times W\times Z$, is a pushout in $s\mathrm{T}\alg$. 
\end{enumerate}
Let $L \subset K$ be an inclusion of small simplicial sets that is bijective on vertices and suppose we are given a diagram 
\[ \mathcal{J}:L\times\Delta^1\coprod_{L\times\{1\}}K\times\{1\} \longrightarrow  s\mathrm{T}'\alg\]
such that for each $k\in K$, the edge $\mathcal{J}|_{\{k\}\times\Delta^1}$ is an effective epimorphism. Denote $\mathcal{J}_0=\mathcal{J}|_{L\times\{0\}}$ and $\mathcal{J}_1=\mathcal{J}|_{K\times\{1\}}$ then the natural diagram 
\begin{equation}\label{eq:colimcompcare}
\begin{tikzcd}
\underset{L}{\colim} f^*(\mathcal{J}_0) \ar[r] \ar[d] & \underset{K}{\colim} f^*(\mathcal{J}_1)\ar[d] \\
f^*(\underset{L}{\colim} \mathcal{J}_0) \ar[r]  & f^*(\underset{K}{\colim} \mathcal{J}_1)
\end{tikzcd}
\end{equation}
in $s\mathrm{T}\alg$ is a pushout.
\end{thm}
\begin{rmk}
It is not difficult to see that a transformation $f:\mathrm{T}\rightarrow\mathrm{T}'$ satisfying the conclusion of the theorem above also satisfies the pushout axiom $(P)$, so the conclusion of the theorem is equivalent to the pushout axiom $(P)$.
\end{rmk}
\begin{cor}\label{cor:unramifiedness2}
Let $f:\mathrm{T}\rightarrow\mathrm{T}'$ be a transformation of Lawvere theories satisfying the pushout axiom $(P)$. Let $h:K\rightarrow s\mathrm{T}'\alg$ be a finite diagram of finitely generated $\mathrm{T}'$-algebras. Then there is a collection of free $\mathrm{T}'$-algebras $\{t_k\}_{k\in K}$ and a pushout diagram 
\[
\begin{tikzcd}
j(\underset{k\in k}{\prod}f(t_k)) \ar[d] \ar[r] & \underset{K}{\colim}f^*\mathcal{J} \ar[d] \\
jf(\underset{k\in K}{\prod}t_k) \ar[r] & f^*(\underset{K}{\colim}\mathcal{J} ).
\end{tikzcd}
\]
in $s\mathrm{T}\alg$.
\end{cor}
\begin{proof}
Choose for each $k\in K$ an effective epimorphism $j(t_k)\rightarrow h(k)$ and apply Theorem \ref{thm:unramified2} to the resulting diagram $\mathrm{sk}_0(K)\times\Delta^1\coprod_{\mathrm{sk}_0(K)\times\{1\}}K\times\{1\}$, where $\mathrm{sk}_0(K)$ is the 0-skeleton of $K$.
\end{proof}
\begin{cor}\label{cor:pushouteffepi}
Let $f:\mathrm{T}\rightarrow\mathrm{T}'$ be a transformation of Lawvere theories satisfying the pushout axiom $(P)$. Then $f^*:s\mathrm{T}'\mathsf{Alg}\rightarrow s\mathrm{T}\mathsf{Alg}$ preserves pushouts along effective epimorphisms. 
\end{cor}
\begin{proof}
Consider a diagram
\[
\begin{tikzcd}
A\ar[d] \ar[r,"g"] & B\\
C
\end{tikzcd}
\]
for which $g$ is an effective epimorphism, and let $\tau:\Lambda^{2}_{0}\times\Delta^1\rightarrow s\mathrm{T}'\mathsf{Alg}$ be the diagram such that $\tau|_{\Lambda^{2}_{0}\times\{1\}}=\sigma$ and $\tau|_{\Lambda^{2}_{0}\times\{0\}}$ is the diagram 
\[
\begin{tikzcd}
A\ar[d] \ar[r,"\mathrm{id}"] & A\\
C,
\end{tikzcd}
\]
and the maps $A\rightarrow A$, $C\rightarrow C$ and $A\rightarrow B$ are the maps $\mathrm{id}_A$, $\mathrm{id}_C$ and $g$ respectively. Then $\tau$ satisfies the conditions of Theorem \ref{thm:unramified2} so that the comparison map $f^*(B)\coprod_{f^*(A)} f^*(C)\rightarrow f^*(B\coprod_A C)$ is a pushout of the map $f^*(A)\coprod_{f^*(A)} f^*(C)\rightarrow f^*(A\coprod_A C)$, which is an equivalence.
\end{proof}
\subsection{Graph inclusions}
The proof of Theorem \ref{thm:unramified2} requires some preparation: we need to resolve an arbitrary effective epimorphism by graph inclusions. Such a resolution is constructed for an arbitrary pregeometry $\pregeo$ in \cite{dagix}, sections 2 and 3, but we have no need of that generality, so we only treat the case of Lawvere theories, using somewhat different methods\footnote{The construction in section 3 of \cite{dagix} is a generalization of the one below if we view Lawvere theories as discrete pregeometries. The reason we offer an alternative proof is that the construction in \emph{loc. cit.} does not appear to be completely correct: in construction 3.7 and remark 3.8, from the data of an \inftop $\xtop$ and a map $\alpha:\Of\rightarrow \Of'$ of local $\pregeo$-structures on $\xtop$, a certain Cartesian fibration $p:\overline{\icate}\rightarrow\xtop$ over an \inftop is constructed, which resolves the map $\alpha$ in a suitable sense. Informally, the \infcat $p^{-1}(U)$ is given by pairs $(\Of^F,S)$ of a `free' $\pregeo$-structure on $U$ and a finite set $S$, a map $\Of^F\rightarrow \Of'|_U$ of local $\pregeo$-structures on $\xtop_{/U}$ and a commuting diagram 
\[
\begin{tikzcd}[ampersand replacement=\&]
 \coprod_S \Of^F \ar[d]\ar[r] \& \Of^F \ar[d] \\
\Of|_U \ar[r,"\alpha|_U"] \& \Of'|_U
\end{tikzcd}
\]
where the upper horizontal morphism is the fold map. It is claimed that this fibre admits coproducts and thus is sifted. However, this does not appear to hold in general, since the only reasonable choice for a coproduct of a pair of data $(\Of^F,S)$ and $(\Of^{F'},S')$ as above in this \infcat is the pair $(\Of^F\coprod \Of^{F'},S\coprod S')$, but there is no reason for the existence of a unique map $\coprod_{S\coprod S'} \Of^F\coprod \Of^{F'}\rightarrow \Of|_U$ that makes the requisite diagram commute, since we are not in general given maps $\coprod_{S'}\Of^F\rightarrow \Of|_U$ and $\coprod_{S}\Of^{F'}\rightarrow \Of|_U$. For the argument to go through, it seems one needs to consider the \infcat of finite tuples $\{(\Of_1^F,S_1),(\Of_2^F,S_2),\ldots\}$ equipped with data as above, but we will not attempt a formal construction at this point}. These resolutions turn out to be a remarkably powerful technical device in themselves, for which we will find many uses. We formalize it in the following proposition.
\begin{prop}[Free resolutions of effective epimorphisms]\label{prop:resolveelementary}
Let $\mathrm{T}$ be a Lawvere theory. Let $\icat$ be an \infcat that admits sifted colimits and let $\icat_0\subset\icat$ be a full subcategory stable under sifted colimits. Suppose we are given a functor $F:\fun(\Delta^1,s\mathrm{T}\alg)\rightarrow \icat$ such that \begin{enumerate}[$(1)$]
    \item $F$ preserves sifted colimits.
    \item For every graph inclusion $g:X\rightarrow X\times Y$ of free simplicial $\mathrm{T}$-algebras, the object $F(j(g))$ lies in $\icat_0$.
\end{enumerate}
Then $F$ carries every effective epimorphism of $s\mathrm{T}\alg$ into $\icat_0$.
\end{prop}
\begin{proof}
Let $\mathrm{T}$ be a Lawvere theory, and $I$ a small index set together with a functor $t_{\_}:I\rightarrow \mathrm{T}$ whose image minimally generates $\mathrm{T}$ under products. The functor 
\[ \ev_{I} : s\mathrm{T}\alg\longrightarrow \fun(I,\spa) \]
adjoint to the functor
\[ s\mathrm{T}\alg \times I =\coprod_{i\in I}s\mathrm{T}\alg\overset{\coprod_{i\in I}\ev_{t_i}}{\longrightarrow} \spa \]
is conservative and preserves limits and sifted colimits and is thus monadic. Let $\mathrm{Free}_{\mathrm{T}}:\fun(I,\spa)\rightarrow s\mathrm{T}\alg$ be a left adjoint to $\ev_{I}$, determined up to equivalence by $\mathrm{Free}_{\mathrm{T}}(*_i)=j(t_i)$ for $i\in I$, where $*_i:I\rightarrow \spa$ carries $i$ to the final space $*$ and all other indices to the initial empty space. Let $K$ be a simplicial set, then the induced adjunction 
\[  \begin{tikzcd}
 \fun(K,s\mathrm{T}\alg) \ar[r,shift left]  &  \fun(K\times I,\spa) \ar[l,shift left]
\end{tikzcd}  \]
is again monadic; letting $K=\Delta^1$, we deduce that for each map $\alpha:A\rightarrow B$ of simplicial $\mathrm{T}$-algebras, there exists an ($\ev_{I}$-split) augmented simplicial object $\alpha_{\bullet}:\simpopplus \times\Delta^1\rightarrow s\mathrm{T}\alg$, the Bar resolution $\mathsf{Bar}_{\ev_I\circ \mathrm{Free}_{\mathrm{T}}}(\ev_I\circ \mathrm{Free}_{\mathrm{T}},\alpha)$, such that $\alpha_{-1}=f$, $\alpha_{\bullet}:\simpopplus\rightarrow \fun(\Delta^1,s\mathrm{T}\alg)$ is a colimit diagram and each $\alpha_n$ is the image of $\ev_I(\alpha_{n-1})$ under $\mathrm{Free}_{\mathrm{T}}$. Recall that we are given a sifted colimit preserving functor
\[ F:\fun(\Delta^1,s\mathrm{T}\alg) \longrightarrow\icat  \]
and a full subcategory $\icat_0\subset\icat$ stable under sifted colimits such that $F(j(g))$ lies in $\icat_0$ for every graph inclusion. We now show that the proposition follows from the following two claims.
\begin{enumerate}
    \item[$(*)$] If $\alpha$ is an effective epimorphism, the map $\mathrm{Free}_{\mathrm{T}}(\ev_I(\alpha))$ is also an effective epimorphism.
    \item[$(**)$] If $\alpha$ is an effective epimorphism, then the object $F(\mathrm{Free}_{\mathrm{T}}(\ev_I(\alpha)))$ lies in $\icat_0$. 
\end{enumerate}
Indeed, if $\alpha$ is an effective epimorphism, then $(*)$ and the construction of $\alpha_{\bullet}$ guarantee that for every $n\geq 0$, the map $\alpha_n$ is an effective epimorphism. It follows from $(**)$ that $F(\alpha_n)$ lies in $\icat_0$ for each $n\geq 0$. Since $F$ preserves sifted colimits and $\icat_0$ is stable under sifted colimits, we conclude. \\
We prove $(*)$. We have a strictly commuting diagram of right adjoints
\[
\begin{tikzcd}
s\mathrm{T}\alg\ar[r,"\ev_{I}"] & \spa^I\\
\tau_{\leq 0}s\mathrm{T}\alg \ar[r,"\ev_I"] \ar[u,hook] & \mathsf{Set}^I \ar[u,hook] 
\end{tikzcd}
\]
hence a commuting diagram of left adjoints
\[
\begin{tikzcd}
s\mathrm{T}\alg\ar[d,"\pi_0"]& \spa^I\ar[d,"\pi_0"]\ar[l,"\mathrm{Free}_{\mathrm{T}}"']\\
\tau_{\leq 0}s\mathrm{T}\alg  & \mathsf{Set}^I \ar[l,"\mathrm{Free}^0_{\mathrm{T}}"'].
\end{tikzcd}
\]
The fact that the truncation functor $\tau_{\leq 0}:s\mathrm{T}\alg\rightarrow \tau_{\leq 0}s\mathrm{T}\alg$ is given by composing product preserving functors with $\pi_0:\spa\rightarrow\mathsf{Set}$ means precisely that the diagram of left adjoints above is horizontally right adjointable. Given a map $\alpha:A\rightarrow B$, it follows that the map $\pi_0(\mathrm{Free}_{\mathrm{T}}(\ev_I(\alpha)))$ is given by applying the free functor $\mathrm{Free}^0_{\mathrm{T}}$ to the map $\pi_0(\ev_I(\alpha))$. By assumption, this latter map is a surjection, so it suffices to show that $\mathrm{Free}^0_{\mathrm{T}}$ carries surjections of $I$-indexed sets to effective epimorphisms of discrete simplicial $\mathrm{T}$-algebras. We conclude by observing that every surjection of $I$-sets is a regular epimorphism, that the left adjoint $\mathrm{Free}^0_{\mathrm{T}}$ preserves regular epimorphisms and that all regular epimorphism in $\tau_{\leq 0}s\mathrm{T}\alg$ are effective as the latter category has kernel pairs.\\
We prove $(**)$. The \infcat $\fun(I,\spa)$ is isomorphic to the nerve of the Kan-enriched category $\fun(I,\mathsf{Kan})=\prod_{i\in I}\mathsf{Kan}$. Let $\alpha:A\rightarrow B$ be a map in $s\mathrm{T}\alg$, then after applying a factorization as a trivial cofibration followed by a fibration, we may assume that the morphism $\ev_I(\alpha)$ in $\fun(I,\spa)$ is an $I$-indexed collection $\{\ev_{t_i}(A)\rightarrow \ev_{t_i}(B)\}_{i\in I}$ of Kan fibrations between Kan complexes. If $\alpha$ is an effective epimorphism, then for each $i$, the map $\ev_{t_i}(A)\rightarrow \ev_{t_i}(B)$ is a surjection on connected components and thus a surjection in each simplicial degree since it is a Kan fibration. We may view $\ev_{t_i}(A)$ and $\ev_{t_i}(B)$ as constant bisimplicial objects, so that we can think of the collection $\{\ev_{t_i}(A)\rightarrow \ev_{t_i}(B)\}_{i\in I}$ as a morphism in the category $\fun(\simpop,\mathsf{Set}_{\simp}^{I})$. Applying the diagonal functor
\[\fun(\simpop,\mathsf{Set}_{\simp}^{I}) \longrightarrow \mathsf{Set}_{\simp}^{I} \]
returns the morphism $\ev_I(\alpha)$ so we deduce that $\ev_I(\alpha)$ is a colimit of the diagram 
\[ \{\ev_{t_i}(A)\rightarrow \ev_{t_i}(B)\}_{i\in I}:\simpop\times\Delta^1 \longrightarrow \mathsf{Set}^I \longrightarrow  \mathsf{Set}_{\simp}^I \longrightarrow \spa^I,  \]
where the last morphism implements localization at the weak equivalences. The functor $\mathsf{Set} \rightarrow  \mathsf{Set}_{\simp} \rightarrow \spa$ can be identified with the inclusion of $0$-truncated spaces, so we conclude that the diagram above is in each simplicial degree $n$ and for each $i\in I$ given by a surjective map
\[ \coprod_{\ev_{t_i}(A)_n}*\longrightarrow \coprod_{\ev_{t_i}(B)_n}* \]
of discrete spaces. Since $\mathrm{Free}_{\mathrm{T}}$ preserves colimits, the map $\mathrm{Free}_{\mathrm{T}}(\ev_I(\alpha))$ arises as the geometric realization of a simplicial diagram that is in each simplicial degree $n$ given by
\begin{equation}\label{eq:resol}
      \coprod_i\coprod_{\ev_{t_i}(A)_n} j(t_i)\longrightarrow \coprod_i\coprod_{\ev_{t_i}(B)_n}j(t_i).
\end{equation}  
Using that $F$ preserves sifted colimits and that $\icat_0\subset\icat$ is stable under sifted colimits again, we are reduced to proving that $F$ carries every morphism of the form \eqref{eq:resol} into $\icat_0$. Writing $I$ as a filtered colimit of its finite subsets, we may assume that $I$ is finite (using that $F$ preserves, and $\icat_0\subset \icat$ is stable under, filtered colimits). Writing for each $i\in I$, the set $\ev_{t_i}(B)_n$ as a filtered colimit of its finite subsets, we may assume that $\ev_{t_i}(B)_n$ is finite. Writing for each $x\in \ev_{t_i}(B)_n$ the set $\ev_{t_i}(\alpha)^{-1}_n(x)\subset \ev_{t_i}(A)_n$ as a filtered colimit of its finite subsets, we may assume that $\ev_{t_i}(A)_n$ is finite. Now we observe that with all sets indexing the coproducts finite, the map \eqref{eq:resol} is a graph inclusion.
\end{proof}

\subsection{The proof of Theorem \ref{thm:unramified2}}
We turn to the proof of Theorem \ref{thm:unramified2} which proceeds by simultaneously decomposing $K$ and the simplicial subset $L$ into simple pieces for which the theorem can be shown to hold using Proposition \ref{prop:resolveelementary}. Our implementation of this idea is unfortunately somewhat technical and extends the material of Section 4.2.3 of \cite{HTT} in a rather ad-hoc manner (certainly much more general results are imaginable). We first need a few preparatory lemmas.

\begin{lem}\label{lem:LKcolim}
Let $L\rightarrow K$ be a map of simplicial sets and let 
\[ f:L\times\Delta^1\coprod_{L\times\{1\}}K\times\{1\}\longrightarrow \icat  \]
be a diagram valued in an \infcat $\icat$. Then the collection of dotted maps
\[ \begin{tikzcd}
 L\times\Delta^1\coprod_{L\times\{1\}}K\times\{1\} \ar[d,hook] \ar[r,"f"] \ar[d]& \icat \\
L^{\rhd}\times\Delta^1\coprod_{L^{\rhd}\times\{1\}}K^{\rhd}\times\{1\}\ar[ur,"\overline{f}"',dotted]
\end{tikzcd} \]
rendering the diagram commutative and having the property that $\overline{f}|_{L^{\rhd}\times\{0\}}:L^{\rhd}\rightarrow\icat$ and $\overline{f}|_{K^{\rhd}\times\{1\}}:K^{\rhd}\rightarrow\icat$ are colimit diagrams, is parametrized by a contractible Kan complex if both $f|_{L\times\{0\}}$ and $f|_{K\times\{1\}}$ admit a colimit and is empty otherwise.
\end{lem}
\begin{proof}
Clearly, the space of extensions is empty if either $f|_{L\times\{0\}}$ or $f|_{K\times\{1\}}$ does not admit a colimit. Let 
\[ \icatd\subset \fun(L^{\rhd}\times\Delta^1\coprod_{L^{\rhd}\times\{1\}}K^{\rhd}\times\{1\},\icat)\cong \fun(L^{\rhd}\times\Delta^1,\icat)\times_{\fun(L^{\rhd},\icat)}\fun(K^{\rhd},\icat)\]
be the full subcategory spanned by functors $\overline{f}$ for which $\overline{f}|_{L^{\rhd}\times\{0\}}$ and $\overline{f}|_{K^{\rhd}\times\{1\}}$ are colimit diagrams, and let 
\[ \icatd'\subset \fun(L\times\Delta^1\coprod_{L\times\{1\}}K\times\{1\},\icat)\cong \fun(L\times\Delta^1,\icat)\times_{\fun(L,\icat)}\fun(K,\icat)\]
be the full subcategory spanned by functors $f$ for which $f|_{L\times\{0\}}$ and $f|_{K\times\{1\}}$ admit colimit diagrams. We show that the natural projection $p:\icatd\rightarrow\icatd'$ is a trivial Kan fibration.
We have a factorization \[L\times\Delta^1\coprod_{L\times\{1\}}K\times\{1\}\overset{i}{\hooklongrightarrow}  L\times\Delta^1\coprod_{L\times\{1\}}K^{\rhd}\times\{1\} \overset{i'}{\hooklongrightarrow} L^{\rhd}\times\Delta^1\coprod_{L^{\rhd}\times\{1\}}K^{\rhd}\times\{1\}. \]
The map $i$ is a pushout of the inclusion $K\hookrightarrow K^{\rhd}$ and the map $i'$ is a pushout of the inclusion 
\[ L\times\Delta^1\coprod_{L\times\{1\}}L^{\rhd}\times\{1\}\hooklongrightarrow L^{\rhd}\times\Delta^1, \] 
so $p$ is a trivial fibration if the functors $p_K:\icatd_K\rightarrow\icatd_K'$ and $p_L:\icatd_L\rightarrow\icatd_L'$ are trivial fibrations, where $\icatd_K\subset\fun(K^{\rhd},\icat)$ is the full subcategory spanned by colimit diagrams, $\icatd_L\subset\fun(L^{\rhd}\times \Delta^{1},\icat)$ is the full subcategory spanned by diagrams $\overline{g}$ such that $\overline{g}|_{L^{\rhd}\times\{1\}}$ is a colimit diagram, and $\icatd_K'$ and $\icatd_L'$ are defined similarly. We give the proof that $p_L$ is a trivial fibration, the other case being similar and easier. It is easy to see that $p_L$ is a categorical fibration (\cite{HTT} Proposition 3.1.2.3), so it suffices to show that this functor is a categorical equivalence. Note that the \infcat $\icatd_L'$ consists of solid diagrams
\[
\begin{tikzcd}
L\ar[d,hook]\ar[r] & \fun(\Delta^1,\icat)\ar[d,"\ev_1"] \\
L^{\rhd} \ar[r] \ar[ur,dotted]& \icat
\end{tikzcd}
\]
for which the composition $L\rightarrow \fun(\Delta^1,\icat)\overset{\ev_0}{\rightarrow}\icat$ admits a colimit, and the \infcat $\icatd_L$ consists of dotted lifts for which the composition $L^{\rhd}\rightarrow \fun(\Delta^1,\icat)\overset{\ev_0}{\rightarrow}\icat$ is a colimit. Combining the fact that $\ev_1$ is a coCartesian fibration associated to the functor carrying a map $h:X\rightarrow Y$ to the projection $\icat_{/Y}\simeq\icat_{/f}\rightarrow\icat_{/X}$ and \cite{HTT}, Proposition 1.2.13.8, 4.3.1.9 and 4.3.1.10, we see that the first condition is tantamount to the condition that the solid diagram admits an $\ev_1$-colimit and the second states that the dotted lift is an $\ev_1$-colimit. In view of \cite{HTT}, Proposition 4.1.3.8, we may assume that $L$ is an \infcatt, in which case the desired result follows form \cite{HTT}, Proposition 4.3.2.15.
\end{proof}
\begin{lem}\label{lem:effepicolimstable}
Let $\mathrm{T}$ be a Lawvere theory, $L\subset K$ an inclusion of simplicial sets that is bijective on vertices and suppose we are given a diagram
\[ f:L\times\Delta^1\coprod_{L\times\{1\}}K\times\{1\}\longrightarrow s\mathrm{T}\alg. \]
If $f|_{\{k\}\times\Delta^1}$ is an effective epimorphism for each $k\in K$, then the canonical map $\colim_L f|_{L\times\{0\}}\rightarrow\colim_K f|_{K\times\{1\}}$ provided by Lemma \ref{lem:LKcolim} is an effective epimorphism as well. 
\end{lem}
\begin{proof}
The map 
\[\mathrm{sk}_0(L)\times\Delta^1\coprod_{\mathrm{sk}_0(L)\times\{1\}}\mathrm{sk}_0(K)\times\{1\}\longrightarrow L\times\Delta^1\coprod_{L\times\{1\}}K\times\{1\}\]
induces a commuting diagram 
\[ 
\begin{tikzcd}
\coprod_{k\in K_0} f|_{L\times\{0\}}(k) \ar[d] \ar[r,"\alpha"] & \coprod_{k\in K_0} f|_{K\times\{1\}}(k) \ar[d,"\beta"] \\
\underset{L}{\colim} f|_{L\times\{0\}}  \ar[r] & \underset{K}{\colim} f|_{K\times\{1\}}.
\end{tikzcd}
\]
It suffices to show that $\alpha$ and $\beta$ are effective epimorphisms. Note that the diagram above is the result of applying the localization functor $L:\pshv(\mathrm{T}^{op})\rightarrow s\mathrm{T}\alg$ to the same diagram with colimits and coproducts taken in $\pshv(\mathrm{T}^{op})$ and $\alpha$ and $\beta$ can be written as $L(\alpha')$ and $L(\beta')$ for $\alpha'$ the upper horizontal and $\beta'$ the right vertical map respectively of that diagram in $\pshv(\mathrm{T}^{op})$. Since effective epimorphisms in $s\mathrm{T}\alg$ are defined via the inclusion $s\mathrm{T}\alg\subset \pshv(\mathrm{T}^{op})$, it follows from \cite{HTT}, Corollary 6.2.3.11 and Lemma 6.2.3.13 that $\alpha'$ and $\beta'$ are effective epimorphisms. We will be done once we show that $L$ preserves effective epimorphisms as an endofunctor on $\pshv(\mathrm{T}^{op})$. Let $h:X_0\rightarrow |\check{C}_{\bullet}(h)|$ be an effective epimorphism, then $L(X_0)\rightarrow L(|\check{C}_{\bullet}(h)|)$ is equivalent to the map $L(X_{\bullet})_0\rightarrow  |L(\check{C}_{\bullet}(h)))|$ since the inclusion $s\mathrm{T}\alg\subset \pshv(\mathrm{T}^{op})$ preserves sifted colimits. Now we conclude by observing that for a simplicial diagram $X_{\bullet}$ in an \inftopt, the map $X_0\rightarrow|X_{\bullet}|$ is an effective epimorphism, which is an immediate consequence of \cite{HTT}, Lemma 6.2.3.13.
\end{proof}
We now give a proof of a very special case of Theorem \ref{thm:unramified2}.
\begin{lem}\label{lem:c3}
Let $f:\mathrm{T}\rightarrow\mathrm{T}'$ be a transformation of Lawvere theories satisfying the pushout axiom $(P)$. Suppose we are given a functor
\[
\mathcal{J}:\Delta^n\longrightarrow s\mathrm{T}'\alg.
\]
Let $v$ denote the final vertex of $\Delta^n$, then the natural diagram 
\[
\begin{tikzcd}
\coprod_{v'\in \Delta^n} f^*\mathcal{J}(v') \ar[d] \ar[r] &  f^*\mathcal{J}(v) \ar[d,equal] \\
f^*(\coprod_{v'\in \Delta^n}\mathcal{J}(v'))   \ar[r] & f^*\mathcal{J}(v)
\end{tikzcd}
\]
is a pushout.
\end{lem}
\begin{proof}
With an easy inductive argument, we may reduce to the case $n=1$. Note that the full subcategory of $\fun(\Delta^1,s\mathrm{T}'\alg)$ consisting of diagrams $\mathcal{J}$ for which the conclusion of the lemma holds is closed under sifted colimits, so (using the Bar resolution and the argument of Proposition \ref{prop:resolveelementary}, for instance) we may suppose that the arrow is of the form $j(t)\rightarrow j(s)$, in which case the conclusion of the lemma is the pushout axiom $(P)$ for the special case of $W=Z=*$.
\end{proof}
The next lemma is a device by which we may replace $L$ and $K$ by simpler simplicial sets. Before we proceed, we recall some notation. Let $K$ be a simplicial set and  $F:\icate\rightarrow (\sset)_{/K}$ be a functor from an ordinary category $\icate$ determining for each $I\in \icate$ a map $\pi_I:F(I)\rightarrow K$. The \emph{$F$-parametrized join of $K$ and $\icate$}, denoted $K_F$ is the simplicial set whose $n$-simplices are given by 
\begin{itemize}
    \item A map $p:\Delta^n\rightarrow \Delta^1$ corresponding to a decomposition $[n]=[0,k]\cup [k+1,n]$
    \item A map $e_0:\Delta^{\{p^{-1}(0)\}}\rightarrow K$.
    \item A map $e_1:\Delta^{\{p^{-1}(1)\}}\rightarrow \icate$ corresponding to a sequence 
    \[ I_{k+1}\longrightarrow I_{k+2}\longrightarrow \ldots \longrightarrow I_n. \]
    \item A map $e_1':\Delta^{\{p^{-1}(0)\}}\rightarrow F(I_{k+1})$ fitting into a commuting diagram 
    \[
    \begin{tikzcd}
    & F(I_{k+1}) \ar[d,"\pi_{I_{k+1}}"] \\
    \Delta^{\{p^{-1}(0)\}} \ar[ur,"e_1'"] \ar[r,"e_1"]& K.
    \end{tikzcd}
    \]
\end{itemize}
For $L\subset K$ a simplicial subset of $K$, we let $L_F$ denote the simplicial subset of $K_F$ whose $n$-simplices correspond to data as above such that the maps $\Delta^{\{p^{-1}(0)\}}\rightarrow K$ and $\Delta^{\{p^{-1}(0)\}}\rightarrow F(I_{k+1})$ factor through $L$ and $F(I_{k+1})\times_K L$ respectively; we may identify this simplicial set with the $F_L$-parametrized join of $L$ and $\icate$ where $F_L$ is the functor $\icate\overset{F}{\rightarrow}(\sset)_{/K}\overset{\_\times_KL}{\rightarrow}(\sset)_{/L}$. Suppose we are given a diagram
\[  \mathcal{J}:L\times \Delta^1\coprod_{L\times\{1\}}K\times\{1\} \longrightarrow s\mathrm{T}'\alg.  \]
so that we have for each $I\in \icate$ a diagram
\[ F(I)\times_KL\times \Delta^1\coprod_{F(I)\times_KL\times\{1\}}F(I)\times\{1\} \longrightarrow s\mathrm{T}'\alg. \]
For $\sigma$ an $n$-simplex of $K$, consider the colimit of the diagram \[\icate\overset{F}{\longrightarrow} (\sset)_{/K}{\longrightarrow} \set \hooklongrightarrow\spa\]
where the second functor carries $S\rightarrow K$ to the set of $n$-simplices of $\Delta^n\times_KS$, and the last functor is the inclusion of 0-truncated spaces. The colimit above is naturally isomorphic in the homotopy category to the nerve of the category of pairs $(I,\sigma')$ with $I\in \icate$ and $\sigma'\in F(I)_n$ such that $\pi_I(\sigma')=\sigma$. Denote this category by $\icate_{\sigma}$ and let $\icate_{\sigma}'\subset\icate_{\sigma}$ be the full subcategory spanned by pairs $(I,\sigma')$ for which $\sigma'$ is degenerate. 
\begin{rmk}\label{rmk:decomp}
Note that if $\sigma$ lies in $L\subset K$, then defining $\icate_{\sigma}$ and $\icate_{\sigma}'$ using the functor $F_L$ instead of $F$ yields the same categories.
\end{rmk}
\begin{lem}\label{lem:decomposetool}
Suppose that the following conditions are satisfied.
\begin{enumerate}[$(a)$]
    \item For each nondegenerate simplex $\sigma$ of $K$, the simplicial set $\icate_{\sigma}$ is weakly contractible.
    \item For each degenerate simplex $\sigma$ of $K$, the inclusion $\icate_{\sigma}'\subset\icate_{\sigma}$ is a weak homotopy equivalence.
    \end{enumerate}
Then there exists a functor $G:\icate^{\rhd}\rightarrow \fun(\Delta^1\times\Delta^1,s\mathrm{T}\alg)$ with the following properties. 
\begin{enumerate}[$(\bullet_1)$]
    \item For each $I\in \icate$, the diagram $G(I)$ is the natural diagram
    \[
    \begin{tikzcd}
    \underset{F(I)\times_KL}{\colim} f^*(\mathcal{J}_0|_{F(I)\times_K L}) \ar[r] \ar[d] & \underset{F(I)}{\colim}     f^*(\mathcal{J}_1|_{F(I)})\ar[d] \\
    f^*(\underset{F(I)\times_KL}{\colim} \mathcal{J}_0|_{F(I)\times_K L}) \ar[r] &     f^*(\underset{F(I)}{\colim} \mathcal{J}_1|_{F(I)})
    \end{tikzcd}
    \]
    \item The diagram $G(\infty)$ is equivalent to the diagram \eqref{eq:colimcompcare}. 
    \item The diagram $\ev_{\{1\}\times\Delta^1}\circ G$ is the image under $f^*$ of a colimit diagram in $\fun(\Delta^1,s\mathrm{T}'\alg)$.
    \item The diagram $\ev_{\{0\}\times\Delta^1}\circ G$ is a colimit diagram.
    \end{enumerate}
    \end{lem}
\begin{proof}
 Let $p:\mathcal{M}\rightarrow\Delta^1$ be a Cartesian fibration associated to $f^*:s\mathrm{T}'\alg\rightarrow s\mathrm{T}\alg$ and consider the \infcat of sections $\mathcal{Z}:=\fun_{\Delta^1}(\Delta^1,\mathcal{M})$, that is, the oplax limit of the arrow $\Delta^1\rightarrow\catinfh$ classifying $f^*$. Let $S$ denote the simplicial set $L\times \Delta^1\coprod_{L\times\{1\}}K\times\{1\}$; we first show that we can extend the diagram $S\rightarrow s\mathrm{T}'\alg$ to $\mathcal{Z}$ such that evaluation at $0$ coincides with the diagram $f^*\circ \mathcal{J}$. This amounts to finding a dotted lift in the diagram 
\[
\begin{tikzcd}
S \ar[r] \ar[d,hook] & \mathcal{M} \ar[d,"p"] \\
S\times\Delta^1 \ar[r,"\pi_{\Delta^1}"']\ar[ur,dotted] & \Delta^1 
\end{tikzcd}
\]
such that for each $s\in S$, the arrow $\{s\}\times \Delta^1\rightarrow\mathcal{M}$ is $p$-Cartesian, where the upper horizontal arrow is the composition $S\rightarrow s\mathrm{T}'\alg\simeq p^{-1}(1)\subset \mathcal{M}$ and the lower horizontal map is the projection onto the second factor. Since $p$ is a Cartesian fibration, we can find such a lift. Denote the corresponding diagram $S\rightarrow\mathcal{Z}$ by $\mathcal{J}_f$. Note that the diagram \eqref{eq:colimcompcare} arises as follows: \cite{HTT}, Proposition 5.1.2.2 shows that the diagrams $L\times\{0\}\rightarrow \mathcal{Z}$ and $K\times\{1\}\rightarrow \mathcal{Z}$ admit colimits, so we can apply Lemma \ref{lem:LKcolim} to $\mathcal{J}_f$ to obtain a diagram 
\[  \overline{\mathcal{J}_f}:L^{\rhd}\times\Delta^1\coprod_{L^{\rhd}\times\Delta^1}K^{\rhd}\times\{1\}\longrightarrow \mathcal{Z}. \]
Restricting to $\{\infty\}\times\Delta^1$, we obtain a commuting diagram 
\[
\begin{tikzcd}
\Delta^1\times\Delta^1 \ar[rr,"s"]\ar[dr,"\pi_2"'] &&  \mathcal{M} \ar[dl,"p"] \\
& \Delta^1
\end{tikzcd}
\]
which we may extend to a diagram $\Delta^{1}\times \Delta^2\rightarrow \mathcal{M}$ such that $\{i\}\times\Delta^{\{1,2\}}$ is a Cartesian edge for $i=0,1$. The restricted diagram $\Delta^1\times \Delta^{\{0,1\}}\subset \Delta^1\times \Delta^2\rightarrow \mathcal{M}$ lies in $p^{-1}(0)\simeq s\mathrm{T}\alg$ and corresponds to the diagram \eqref{eq:colimcompcare}, by \cite{HTT}, Proposition 5.1.2.2. In analogy with \cite{HTT}, Proposition 4.2.3.4 we show the following.
\begin{enumerate}
    \item[$(*)$] There exists an extension $q:L_F\times\Delta^1\coprod_{L_F\times\{1\}}K_F\times\{1\}\rightarrow\mathcal{Z}$ of $\mathcal{J}_f$ such that the composition $\icate\times\{0\}\subset L_F\times\{0\}\rightarrow\mathcal{Z}$ carries each $I \in \icate$ to a colimit of the composition $F(I)\times_K L\rightarrow L{\rightarrow}\mathcal{Z}$, and the composition $\icate\times\{1\}\subset K_F\times\{1\}\rightarrow\mathcal{Z}$ carries each $I \in \icate$ to a colimit of the composition $F(I)\rightarrow K{\rightarrow}\mathcal{Z}$
\end{enumerate}
Choose a well-ordering on the set of $\icate_{\mathrm{nd}}$ of nondegenerate simplices of $\icate$ such that $d(\sigma)< d(\tau)$ implies $\sigma<\tau$ where $d(\_)$ denotes the dimension. Let $\gamma$ be the order type of $\icate_{\mathrm{nd}}$ so that we have a bijection $\alpha\mapsto \sigma_{\alpha}$ for $\alpha<\gamma$. Let $Y_{\alpha}\subset \icate$ be the simplicial subset spanned by the nondegenerate simplices $\{\sigma_{\beta}\}_{\beta<\alpha}$. Let $LK_{\alpha}$ be the simplicial subset of $L_F\times\Delta^1\coprod_{L_F\times\{1\}}K_F\times\{1\}$ consisting of simplices $\tau$ for which $\tau\cap \icate\times\Delta^1$ lies in $Y_{\alpha}\times\Delta^1$ and note that $LK_0=L\times\Delta^1\coprod_{L\times\{1\}}K\times\{1\}$. We construct a compatible sequence of maps $\{q_{\alpha}\}_{\alpha}:LK_{\alpha}\rightarrow\mathcal{Z}$ with $q_0=\mathcal{J}_f$ and set $q=\colim_{\alpha<\gamma}q_{\alpha}$. By a standard argument, we may reduce to the case of successor ordinals, that is, given $q_{\alpha}$, we construct $q_{\alpha+1}$. Suppose that $Y_{\alpha+1}$ is obtained from $Y$ by adding a vertex $X_I$ corresponding to some object $I\in \icate$, then the inclusion $LK_{\alpha}\subset LK_{\alpha+1}$ is a pushout of the inclusion
\[ F(I)\times_L K\times\Delta^1\coprod_{F(I)\times_L K\times\{1\}} F(I)\times\{1\}\hooklongrightarrow (F(I)\times_L K)^{\rhd}\times\Delta^1\coprod_{(F(I)\times_L K)^{\rhd}\times\{1\}} F(I)^{\rhd}\times\{1\}  \]
and we use Lemma \ref{lem:LKcolim} to find an extension $q_{\alpha+1}$ satisfying $(*)$. If $Y_{\alpha+1}$ is obtained from $Y_{\alpha}$ by adding some nondegenerate $n$-simplex $\sigma$ for $n>0$ corresponding to some sequence of maps $I_0\rightarrow\ldots\rightarrow I_n$ in $\icate$, then the inclusion $LK_{\alpha}\subset LK_{\alpha+1}$ is a pushout of the inclusion
\[ F(I_0)\times_K L\star \del\sigma\times\Delta^1\coprod _{F(I_0)\times_K L\star \del\sigma\times\{1\}} F(I_0)\star \del \sigma\times\{1\}\hooklongrightarrow F(I_0)\times_K L\star \sigma\times\Delta^1\coprod _{F(I_0)\times_K L\star \sigma\times\{1\}} F(I_0)\star \sigma\times\{1\}  \]
This map factorizes as a pushout of the inclusion
\begin{equation}\label{eq:join1} F(I_0)\star\del\sigma \hooklongrightarrow F(I_0)\star \sigma \end{equation}
followed by a pushout of the inclusion
\begin{equation}\label{eq:join2} (F(I_0)\times_K L\star \del\sigma)\times\Delta^1\coprod _{(F(I_0)\times_K L\star \del\sigma)\times\{1\}} (F(I_0)\times_K L\star \sigma)\times\{1\}\hooklongrightarrow (F(I_0)\times_K L\star \sigma)\times\Delta^1.  \end{equation}
To find an extension along the map \eqref{eq:join1} is to find a solution to the lifting problem
\[
\begin{tikzcd}
F(I_0)\star\del\sigma \ar[d,hook] \ar[r] & \mathcal{Z} \\
F(I_0)\star\sigma \ar[ur,dotted] 
\end{tikzcd}
\]
but since the the restriction of the upper horizontal map to the join of $F(I_0)$ with the initial vertex of $
\del\sigma$ is a colimit diagram by construction, the desired extension may be found. To find an extension along the map \eqref{eq:join2} is to find a solution to the lifting problem
\[ 
\begin{tikzcd}
F(I_0)\times_K L\star \del\sigma \ar[d,hook] \ar[r] & \fun(\Delta^1,\mathcal{Z}) \ar[d,"\ev_1"] \\
F(I_0)\times_K L\star \sigma \ar[r] \ar[ur,dotted] & \mathcal{Z},
\end{tikzcd}
\]
but since the restriction of the upper horizontal map to the join of $F(I_0)\times_K L$ with the initial vertex of $\del\sigma$ is an $\ev_1$-colimit diagram by the proof of Lemma \ref{lem:LKcolim}, the desired extension may be found, which concludes the construction of $q$. We note that our construction and \cite{HTT}, Proposition 4.2.3.4 guarantee that the maps 
\[ \mathcal{Z}_{q|_{L_F\times\{0\}}/}\longrightarrow \mathcal{Z}_{\mathcal{J}_f|_{L\times\{0\}}/},\quad\quad \mathcal{Z}_{q|_{K_F\times\{1\}}/}\longrightarrow \mathcal{Z}_{\mathcal{J}_f|_{K\times\{1\}}/} \]
are both trivial fibrations. Thus, upon applying Lemma \ref{lem:LKcolim} to $q$, we obtain a diagram 
\[ \overline{q}:L_F^{\rhd}\times\Delta^1\coprod_{L_F^{\rhd}\times\{1\}}K_F^{\rhd}\times\{1\} \longrightarrow \mathcal{Z}  \]
such that the composition 
\[  L^{\rhd}\times\Delta^1\coprod_{L^{\rhd}\times\{1\}}K^{\rhd}\times\{1\} \hooklongrightarrow L_F^{\rhd}\times\Delta^1\coprod_{L_F^{\rhd}\times\{1\}}K_F^{\rhd}\times\{1\} \overset{\overline{q}}{\longrightarrow} \mathcal{Z} \]
is equivalent to $\overline{\mathcal{J}_f}$ (here we use the contractibility in the statement of Lemma \ref{lem:LKcolim}). It follows from $(a)$, $(b)$, Remark \ref{rmk:decomp} and \cite{HTT}, Proposition 4.2.3.8 that the maps
\[  \mathcal{Z}_{q|_{L_F\times\{0\}}/} \longrightarrow \mathcal{Z}_{q|_{\icate\times\{0\}}/},\quad\quad\mathcal{Z}_{q|_{K_F\times\{1\}}/}\longrightarrow\mathcal{Z}_{q|_{\icate\times\{1\}}/} \]
are both trivial fibrations, which implies that the composition
\[  \icate^{\rhd}\times\Delta^1\hooklongrightarrow L_F^{\rhd}\times\Delta^1\coprod_{L_F^{\rhd}\times\{1\}}K_F^{\rhd}\times\{1\} \overset{\overline{q}}{\longrightarrow} \mathcal{Z}   \]
corresponds to a colimit diagram $\icate^{\rhd}\rightarrow\fun(\Delta^1,\mathcal{Z})$. We obtain a diagram $\icate^{\rhd}\times\Delta^1\times\Delta^1\rightarrow\mathcal{M}$ which we may extend to a diagram on $\icate^{\rhd}\times\Delta^1\times\Delta^2$ such that for each $X\in \icate^{\rhd}\times\Delta^1$, the edge $\{X\}\times\Delta^{\{1,2\}}$ is $p$-Cartesian. The restriction to $\icate^{\rhd}\times\Delta^1\times\Delta^{\{0,1\}}$ lies in $p^{-1}(0)\simeq  s\mathrm{T}\alg$ and corresponds to a functor
\[  G:\icate^{\rhd}\times\Delta^1\times \Delta^1\longrightarrow  s\mathrm{T}\alg\]
satisfying $(\bullet_1)$ through $(\bullet_4)$.
\end{proof}
\begin{proof}[Proof of Theorem \ref{thm:unramified2}]
Let $P$ be the set of finite simplicial subsets of $K$ which is partially ordered and filtered, determining a functor $F:P\rightarrow (\sset)_{/K}$ satisfying the conditions of Lemma \ref{lem:decomposetool}. We deduce the existence of a diagram $G_P:P^{\rhd}\rightarrow \fun(\Delta^1\times\Delta^1,s\mathrm{T}\alg)$ satisfying $(\bullet_1)$ through $(\bullet_4)$. Using $(\bullet_4)$, $(\bullet_3)$ and the fact that $f^*$ commutes with sifted colimits, we see that $G_P$ is a colimit diagram. Using $(\bullet_1)$ and $(\bullet_2)$, we are reduced to proving the theorem for finite simplicial sets. We now proceed by induction on the dimension, that is, we show that the theorem holds for simplicial sets $K$ of dimension $n$ for an arbitrary inclusion $L\subset K$ that is bijective on vertices provided the same statement holds for simplicial sets of lower dimension. For the base case, we assume that $K$ is a finite collection of points, so that the inclusion $L\subset K$ is an equality. It suffices to prove that diagrams of the form 
\[
\begin{tikzcd}
f^*(X_1)\coprod f^*(X_2)\coprod \ldots\coprod f^*(X_n)\ar[r]\ar[d] & f^*(Y_1)\coprod f^*(Y_2)\coprod\ldots\coprod f^*(Y_n)\ar[d] \\
f^*(X_1\coprod X_2\coprod \ldots\coprod X_n) \ar[r] & f^*(Y_1\coprod Y_2\coprod \ldots\coprod Y_n)
\end{tikzcd}
\]
are pushouts, where each $X_i\rightarrow Y_i$ is an effective epimorphism. With an easy inductive argument we may reduce to the case where the coproducts are binary. Denote the effective epimorphisms by $\alpha:A\rightarrow C$ and $\beta:B\rightarrow D$. Consider the composition 
\[ \phi:\fun(\Delta^1,s\mathrm{T}'\alg)\times \{\beta\}\longrightarrow \fun(\del\Delta^1,\fun(\Delta^1,s\mathrm{T}'\alg))\longrightarrow \fun((\del\Delta^1)^{\rhd},\fun(\Delta^1,s\mathrm{T}'\alg)),  \]
where the second functor is a functor taking colimits, a section of the trivial fibration provided by \cite{HTT}, Proposition 4.3.2.15. Now consider the restriction functor
\[ \fun(\del\Delta^1\star \Delta^1,\fun(\Delta^1,s\mathrm{T}\alg)) \longrightarrow \fun((\del\Delta^1)^{\rhd},\fun(\Delta^1,s\mathrm{T}\alg))  \]
induced by the full subcategory inclusion $i:(\del\Delta^1)^{\rhd}= \del\Delta^1\star \Delta^{\{1\}}\subset\del\Delta^1\star \Delta^1$. Let $\fun'(\del\Delta^1\star \Delta^1,\fun(\Delta^1,s\mathrm{T}\alg))\subset \fun(\del\Delta^1\star \Delta^1,\fun(\Delta^1,s\mathrm{T}\alg))$ be the full subcategory spanned by functors which are left Kan extensions along $i$. The restriction map $\fun'(\del\Delta^1\star \Delta^1,\fun(\Delta^1,s\mathrm{T}\alg))\rightarrow \fun((\del\Delta^1)^{\rhd},\fun(\Delta^1,s\mathrm{T}\alg))$ is a trivial fibration by \cite{HTT}, Proposition 4.3.2.15. Choosing a section of this fibration and composing with the restriction $\Delta^1\subset \del\Delta^1\star \Delta^1$ yields a functor \[ \chi:\fun((\del\Delta^1)^{\rhd},\fun(\Delta^1,s\mathrm{T}\alg)) \longrightarrow\fun(\Delta^1,\fun(\Delta^1,s\mathrm{T}\alg))\cong\fun(\Delta^1\times\Delta^1,s\mathrm{T}\alg). \]
Composing $f^*\circ \phi$ with $\chi$ then gives a functor
\[\fun(\Delta^1,s\mathrm{T}'\alg) \longrightarrow \fun(\Delta^1\times\Delta^1,s\mathrm{T}\alg) \]
which carries a map $g:R\rightarrow S$ of simplicial $\mathrm{T}'$-algebras to the diagram
\[
\begin{tikzcd}
 f^*(R)\coprod f^*(B)\ar[d] \ar[r,"f^*(g)\coprod f^*(\beta)"] &[2em]  f^*(S)\coprod f^*(D) \ar[d] \\
f^*(R\coprod B) \ar[r,"f^*(g\coprod \beta)"] &[2em] f^*(S\coprod D) 
\end{tikzcd}
\]
Note that this functor preserves sifted colimits and that the collection of pushout diagrams is closed under colimits in $\fun(\Delta^1\times\Delta^1,s\mathrm{T}\alg)$. Invoking Proposition \ref{prop:resolveelementary}, we see that it is sufficient to argue that the natural commuting diagram 
\[
    \begin{tikzcd}
    f^*j(X\times Y)\coprod f^*(B)\ar[r]\ar[d]& f^*j(X)\coprod f^*(D)\ar[d]\\
     f^*(j(X\times Y)\coprod B) \ar[r] & f^*(j(X)\coprod D)
    \end{tikzcd}
    \]
is a pushout for every graph inclusion $X\rightarrow X\times Y$. Now we apply this argument again to the effective epimorphism $B\rightarrow D$ to reduce to the diagram 
\[
    \begin{tikzcd}
    f^*j(X\times Y)\coprod f^*j(W\times Z)\ar[r]\ar[d]& f^*j(X)\coprod f^*j(W)\ar[d]\\
     f^*(j(X\times Y)\coprod j(W\times Z)) \ar[r] & f^*(j(X)\coprod j(W)).
    \end{tikzcd}
\]
We extend this diagram like so
\[
    \begin{tikzcd}
     f^*j(X)\coprod f^*j(Y)\coprod f^*j(W)\coprod f^*j(Z)\ar[r]\ar[d]& f^*j(X)\coprod f^*j(W)\ar[d,equal]\\
    f^*j(X\times Y)\coprod f^*j(W\times Z)\ar[r]\ar[d]& f^*j(X)\coprod f^*j(W)\ar[d]\\
     f^*(j(X\times Y)\coprod j(W\times Z)) \ar[r] & f^*(j(X)\coprod j(W)).
    \end{tikzcd}
\]
Now we observe that the top square and the outer rectangle are pushouts by virtue of the pushout axiom $(P)$, so the bottom square is a pushout as well.\\
Next we assume that $K$ has a finite number of nondegenerate simplices in dimension $n>0$ and no nondegenerate simplices of higher dimension, and we assume that the theorem holds for simplicial sets of dimension $<n$. Order the set of nondegenerate $n$-dimensional simplices via a bijection to some $\{1,\ldots,p\}$, and let $K_n^i$ be the simplicial set spanned by the nondegenerate simplices of dimension $<n$ and the nondegenerate $n$-simplices smaller in the order than $i\leq p$. Let $K_n^0:=K_{n-1}$ be the $(n-1)$-skeleton of $K$, then for $i\geq 1$, $K_n^i$ is obtained as a pushout $K^{i-1}_{n}\coprod_{\del\Delta^{n}}\Delta^n$ determined by a functor $\Lambda^2_0\rightarrow (\sset)_{/K_n^i}$. Let $v$ denote the final vertex of $\Delta^n\subset K_n^i$, then $\colim_{\Delta^n}f^*(\mathcal{J}_1|_{\Delta^n})\simeq f^*(\colim_{\Delta^n}\mathcal{J}_1|_{\Delta^n})\simeq f^*\mathcal{J}_1(v)$. For brevity, we will write
\[ \Delta^n_L:=\del\Delta^n\times_KL,\quad\quad \del\Delta^n_L:=\del\Delta^n\times_KL,\quad \quad L^i_n:=K^{i}_n\times_KL.\]
Apply Lemma \ref{lem:decomposetool} to the functor $\simpop\rightarrow (\sset)_{/K_n^i}$ carrying $[m]$ to $K^{i-1}_n\coprod \del\Delta^{n}\coprod \ldots \coprod \del\Delta^{n} \coprod \Delta^{n}$, then using $(\bullet_1)$ through $(\bullet_4)$ and the fact that $f^*$ preserve sifted colimits, we are reduced to proving that as $m$ ranges over the positive integers, the diagrams
\begin{equation*}
\begin{tikzcd}
\underset{\Delta^n_L}{\colim}f^*(\mathcal{J}_0|_{\Delta^n_L})  \coprod\underset{\del\Delta^n_L}{\colim}f^*(\mathcal{J}_0|_{\del\Delta^n_L})^{\coprod_m}\coprod \underset{L^{i-1}_n}{\colim}f^*(\mathcal{J}_0|_{L^{i-1}_n}) \ar[r] \ar[d] & f^*\mathcal{J}_1(v) \coprod\underset{\del\Delta^n}{\colim}f^*(\mathcal{J}_1|_{\del\Delta_n})^{\coprod_m}\coprod\underset{K^{i-1}_n}{\colim}f^*(\mathcal{J}_1|_{K^{i-1}_n})\ar[d] \\ f^*(\underset{\Delta^n_L}{\colim}\mathcal{J}_0|_{\Delta^n_L} \coprod(\underset{\del\Delta^n_L}{\colim}\mathcal{J}_0|_{\del\Delta^n_L})^{\coprod_m}\coprod\underset{L^{i-1}_n}{\colim}\mathcal{J}_0|_{L^{i-1}_n})\ar[r]  &f^*(\mathcal{J}_1(v) \coprod(\underset{\del\Delta^n}{\colim}\mathcal{J}_1|_{\del\Delta_n})^{\coprod_m}\coprod\underset{K^{i-1}_n}{\colim}\mathcal{J}_1|_{K^{i-1}_n})
\end{tikzcd}
\end{equation*}
with $m$ coproducts in the middle are pushouts for all $1\leq i\leq p$. We proceed by induction on $i$ and prove the inductive step and the base case simultaneously. We claim that the diagrams
\[
\begin{tikzcd}
 \underset{L^{i-1}_n}{\colim}f^*(\mathcal{J}_0|_{L^{i-1}_n}) \ar[r] \ar[d] & \underset{K^{i-1}_n}{\colim}f^*(\mathcal{J}_1|_{K^{i-1}_n})\ar[d]  \\
 f^*(\underset{L^{i-1}_n}{\colim}\mathcal{J}_0|_{L^{i-1}_n}) \ar[r] &  f^*(\underset{K^{i-1}_n}{\colim}\mathcal{J}_1|_{K^{i-1}_n})
\end{tikzcd} \quad \quad
\begin{tikzcd}
\underset{\del\Delta^n_L}{\colim}f^*(\mathcal{J}_0|_{\del\Delta^n_L}) \ar[r] \ar[d] &\underset{\del\Delta^n}{\colim}f^*(\mathcal{J}_1|_{\del\Delta_n})\ar[d] \\ f^*(\underset{\del\Delta^n_L}{\colim}\mathcal{J}_0|_{\del\Delta^n_L})\ar[r]  &f^*(\underset{\del\Delta^n}{\colim}\mathcal{J}_1|_{\del\Delta_n})
\end{tikzcd}
\]
and the diagram 
\[
\begin{tikzcd}
\underset{\Delta^n_L}{\colim}f^*(\mathcal{J}_0|_{\Delta^n_L}) \ar[r] \ar[d] &f^*\mathcal{J}_1(v)\ar[d] \\ f^*(\underset{\Delta^n_L}{\colim}\mathcal{J}_0|_{\Delta^n_L})\ar[r]  &f^*\mathcal{J}_1(v)
\end{tikzcd}
\]
are pushouts. For the two upper diagrams, we apply the inductive hypotheses on $i$ and on the dimension. For the lower one, we reason as follows. If $\Delta^n_L\subset\Delta^n$ contains the nondegenerate simplex in dimension $n$, then $\Delta^n_L=\Delta^n$ and the diagram is a pushout as both vertical maps are equivalences. If $\Delta^n_L\neq\Delta^n$, then $\Delta^n_L$ is of dimension $<n$ and we extend the diagram as follows
\[
\begin{tikzcd}
\coprod_{v'\in \Delta^n} f^*\mathcal{J}_0(v')\ar[r] \ar[d] &
\underset{\Delta^n_L}{\colim}f^*(\mathcal{J}_0|_{\Delta^n_L}) \ar[r] \ar[d] &f^*\mathcal{J}_1(v)\ar[d] \\ f^*(\coprod_{v'\in \Delta^n}\mathcal{J}_0(v')) \ar[r]&f^*(\underset{\Delta^n_L}{\colim}\mathcal{J}_0|_{\Delta^n_L})\ar[r]  &f^*\mathcal{J}_1(v).
\end{tikzcd}
\]
The inductive hypothesis applied to the $<n$-dimensional simplicial set $\Delta^n_L$ guarantees that the left square is a pushout, so it suffices to show that the outer rectangle is one as well. This rectangle is also the outer rectangle in the diagram 
\[
\begin{tikzcd}
\coprod_{v'\in \Delta^n} f^*\mathcal{J}_0(v')\ar[r] \ar[d] &
\coprod_{v'\in \Delta^n} f^*\mathcal{J}_1(v') \ar[r] \ar[d] &f^*\mathcal{J}_1(v)\ar[d] \\ f^*(\coprod_{v'\in \Delta^n}\mathcal{J}_0(v')) \ar[r]&f^*(\coprod_{v'\in \Delta^n}\mathcal{J}_1(v'))\ar[r]  &f^*\mathcal{J}_1(v).
\end{tikzcd}
\]
Here, the left square is a pushout by the inductive hypothesis, and the right square is a pushout by Lemma \ref{lem:c3}. We return to the square diagram each object of which is composed of $m+2$ summands. First suppose that $m=0$, so that we have binary coproducts. The relevant diagram factorizes as the composition of the following three squares.
\[
\begin{tikzcd}
\underset{\Delta^n_L}{\colim}f^*(\mathcal{J}_0|_{\Delta^n_L}) \coprod \underset{L^{i-1}_n}{\colim}f^*(\mathcal{J}_0|_{L^{i-1}_n}) \ar[r] \ar[d] & f^*\mathcal{J}_1(v) \coprod\underset{K^{i-1}_n}{\colim}f^*(\mathcal{J}_1|_{K^{i-1}_n})\ar[d]  \\
 \underset{\Delta^n_L}{\colim}f^*(\mathcal{J}_0|_{\Delta^n_L}) \coprod f^*(\underset{L^{i-1}_n}{\colim}\mathcal{J}_0|_{L^{i-1}_n}) \ar[r] \ar[d] & f^*\mathcal{J}_1(v) \coprod f^*(\underset{K^{i-1}_n}{\colim}\mathcal{J}_1|_{K^{i-1}_n})\ar[d]  \\
 f^*(\underset{\Delta^n_L}{\colim}\mathcal{J}_0|_{\Delta^n_L}) \coprod f^*(\underset{L^{i-1}_n}{\colim}\mathcal{J}_0|_{L^{i-1}_n}) \ar[r] \ar[d] & f^*\mathcal{J}_1(v) \coprod f^*(\underset{K^{i-1}_n}{\colim}\mathcal{J}_1|_{K^{i-1}_n})\ar[d]  \\
f^*(\underset{\Delta^n_L}{\colim}\mathcal{J}_0|_{\Delta^n_L}\coprod \underset{L^{i-1}_n}{\colim}\mathcal{J}_0|_{L^{i-1}_n}) \ar[r] & f^*(\mathcal{J}_1(v) \coprod\underset{K^{i-1}_n}{\colim}\mathcal{J}_1|_{K^{i-1}_n}).
\end{tikzcd}
\]
We have just shown that the upper and middle square are coproducts of pushout diagrams. For the lower square, we note that Lemma \ref{lem:effepicolimstable} guarantees that the maps $\underset{\Delta^n_L}{\colim}\mathcal{J}_0|_{\Delta^n_L}\rightarrow \mathcal{J}_1(v)$ and $\underset{L^{i-1}_n}{\colim}\mathcal{J}_0|_{L^{i-1}_n}\rightarrow \underset{K^{i-1}_n}{\colim}\mathcal{J}_1|_{K^{i-1}_n}$ are effective epimorphisms. Now the assertion that the lower square is a pushout is the statement of the theorem for $K=\del\Delta^1$, the base case of our induction. The cases of $m>0$ coproducts in the middle are handled inductively using the same argument.
\end{proof}
\newpage

\printbibliography
\end{document}